\documentclass[letterpaper,12pt]{book}

\usepackage[utf8]{inputenc}
\usepackage[T1]{fontenc}
\usepackage{textcomp}
\usepackage{libertinus}

\usepackage[english]{babel}

\usepackage{amsmath,amssymb,amsfonts,amsthm}
\usepackage{mathrsfs}

\usepackage[
    lmargin=2.75cm,
    rmargin=2.75cm,
    top=4cm,
    bottom=3.5cm
]{geometry}

\usepackage[
    hidelinks,
    colorlinks=true,
    linkcolor=blue,
    citecolor=red,
    urlcolor=blue
]{hyperref}

\usepackage{fancyhdr}

\usepackage[pdftex]{graphicx}
\usepackage{epstopdf}

\usepackage{tikz}
\usepackage{tikz-cd}

\usetikzlibrary{positioning}
\usetikzlibrary{matrix}
\usetikzlibrary{shapes.geometric}

\usepackage{pgfplots}
\pgfplotsset{compat=newest}

\usepackage{float}
\usepackage{caption}

\usepackage{makeidx}

\theoremstyle{plain}

\newtheorem{theorem}{Theorem}[chapter]
\newtheorem{corollary}[theorem]{Corollary}
\newtheorem{lemma}[theorem]{Lemma}
\newtheorem{proposition}[theorem]{Proposition}
\newtheorem{definition}[theorem]{Definition}
\newtheorem{example}[theorem]{Example}

\newtheorem{problem}[theorem]{Problem}
\newtheorem{exercise}{Exercise}[chapter]
\newtheorem{remark}[theorem]{Remark}

\renewenvironment{proof}[1][Proof]{%
    \par\noindent\textit{#1. }\ignorespaces
    \pushQED{\qed}%
}{%
    \popQED\par
}

\title{\bf {Introducción a la teoría de la medida e integración}}
\author{María de los Ángeles Sandoval Romero\\
        Hugo Guadalupe Reyna Castañeda\\
        Luis Antonio Cedeño Pérez}
\date{Departamento de Matemáticas, Facultad de Ciencias\\
      Universidad Nacional Autónoma de México\\[0.5cm]
       Semestre 2026-2\\
       Febrero-Mayo 2026}

\usepackage{emptypage}
\let\olditemize\itemize
\def\itemize{\olditemize\itemsep=-1pt}
\usepackage{makeidx}
\makeindex

\begin{document}

\begin{titlepage}

\vspace*{3cm}

\begin{center}

{\Large\bfseries Introduction to Measure and Integration Theory}

\vspace{0.5cm}

{\large
María de los Ángeles Sandoval-Romero\\[0.1cm]
Hugo Guadalupe Reyna-Castañeda\\[0.1cm]
Luis Antonio Cedeño-Pérez}

\vspace{0.5cm}

{\large
Departamento de Matemáticas, Facultad de Ciencias\\
Universidad Nacional Autónoma de México\\
First Preliminary Version:
February–May 2026
}

\vspace{8cm}

{\small
2020 \textit{Mathematics Subject Classification:} 28A20, 28A25, 46E30, 26A42

\textit{Keywords:} measurable functions, measures, Lebesgue integral, Lebesgue spaces, modes of convergence.
}

\end{center}

\end{titlepage}

\thispagestyle{empty} 

\chapter*{Preface}

These notes originated from several courses on Mathematical Analysis II taught at the Faculty of Sciences of the National Autonomous University of Mexico (UNAM). Their primary purpose is to provide a rigorous and accessible introduction to measure and integration theory, aimed principally at undergraduate students in mathematics and related areas who are encountering these topics for the first time.

Measure theory occupies a central position in modern analysis. In addition to providing the natural foundation for the Lebesgue integral, its ideas and methods play an essential role in areas as diverse as probability theory, differential equations, functional analysis, ergodic theory, harmonic analysis, and numerous contemporary mathematical models. For this reason, we consider it important to present its fundamental concepts from a perspective that combines mathematical rigor, conceptual motivation, and intuitive development.

In preparing this text, we have sought to emphasize not only the formal constructions, but also the ideas that historically led to the development of the theory. In particular, special attention is devoted to the problem of limiting processes, the limitations of the Riemann integral, and the fundamental role played by the notions of measure, convergence, and approximation in modern analysis.

The approach adopted throughout the text seeks to maintain a balance between generality and pedagogical clarity. Whenever possible, we have favored detailed proofs and illustrative examples, with the intention of allowing the reader to develop a gradual and solid understanding of the concepts. In addition, exercises and small projects are included in order to complement the theory and encourage a more active participation on the part of the student.

This version constitutes the first edition of the manuscript and, as with any work in progress, it may still contain typographical errors, inaccuracies, or aspects susceptible to improvement, both in the exposition and in some of the proofs. We would sincerely appreciate any comments, corrections, or suggestions that may contribute to improving future versions of these notes.

We hope that these notes will be useful to those beginning the study of measure theory and that they may serve as an invitation to explore more deeply some of the most beautiful and fundamental ideas in mathematical analysis. This work was carried out within the framework of the PAPIME project PE100726, supported by DGAPA-UNAM.

\begin{flushright}
Mexico, 2026
\end{flushright}

\lhead[{\scriptsize \thepage}]{ {\scriptsize \rightmark }}
\rhead[{\scriptsize \leftmark}]{ {\scriptsize \thepage}}

\markboth
{{\scriptsize CONTENTS}}
{{\scriptsize CONTENTS}}

{
\hypersetup{linkcolor=black}
\tableofcontents
}

\lhead[{\scriptsize \thepage}]{ {\scriptsize \rightmark }}
\rhead[{\scriptsize \leftmark}]{ {\scriptsize \thepage}}

\markboth
{{\scriptsize CONTENTS}}
{{\scriptsize CONTENTS}}

\thispagestyle{empty} 

\part{The Lebesgue Integral}
\chapter{The problem of integration}\label{Capitulo1}

\markboth{{\scriptsize 1. THE PROBLEM OF INTEGRATION}}
{{\scriptsize 1. THE PROBLEM OF INTEGRATION}}

To motivate the emergence of measure theory, we begin with a review of the construction of the Riemann integral\footnote{Georg Friedrich Bernhard Riemann (1826--1866) was a German mathematician who made fundamental contributions to analysis and differential geometry, some of which paved the way for the development of general relativity. His name is associated with, among other concepts, the zeta function, the Riemann hypothesis, the Riemann integral, Riemann's lemma, Riemann manifolds and surfaces, as well as Riemannian geometry.} on $\mathbb{R}$, studied in the differential and integral calculus courses typically taken during the first semesters of an undergraduate degree in mathematics.

\begin{figure}[ht!]
\centering
\includegraphics[scale=0.235]{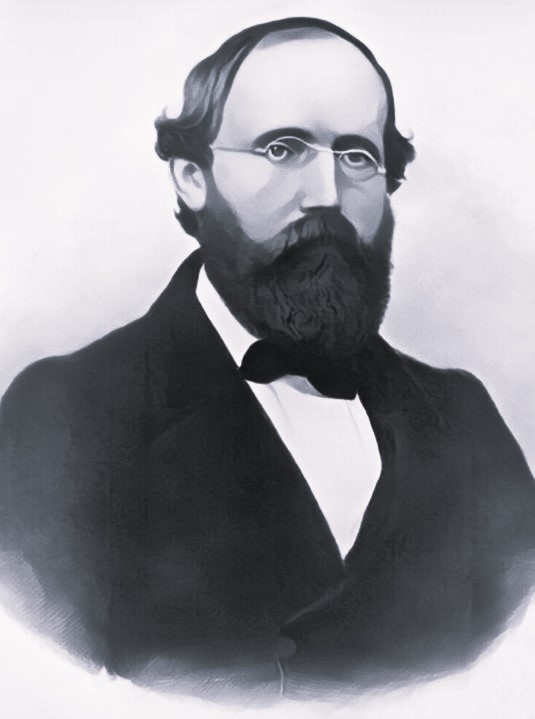} 
\begin{center}
G. Riemann (1826-1866)
\end{center}
\end{figure}

The fundamental idea behind the Riemann integral is to approximate the area under the graph of a function by means of sums of areas of rectangles. In order to formalize this procedure, we introduce the concept of a partition of an interval.

A \textit{partition} \index{partition} $P$ of the interval $[a,b]$ is a finite sequence $\{x_0,x_1,\ldots,x_n\}$ of real numbers such that
$$
a=x_{0} < x_{1} < \cdots < x_{n}=b,
$$
and its \textit{norm} is defined as the real number
$$
\Vert P \Vert:=\max_{1\leq i \leq n}\{x_i-x_{i-1}\}.
$$

We denote by $\mathfrak{P}[a,b]$ the collection of all partitions of the interval $[a,b]$.

Let $f:[a,b]\to\mathbb{R}$ be a nonnegative bounded function, and let $P\in\mathfrak{P}[a,b]$ be a partition. For each $i\in\{1,\ldots,n\}$, define
$$
m_i:=\inf\{f(x):x\in[x_{i-1},x_i]\}
\quad\text{and}\quad
M_i:=\sup\{f(x):x\in[x_{i-1},x_i]\}.
$$

For example, if $f$ is increasing on $[a,b]$, then for every $i\in\{1,\ldots,n\}$ we have
$$
m_i=f(x_{i-1})
\quad\text{and}\quad
M_i=f(x_i).
$$

The information provided by the values $m_i$ and $M_i$ makes it possible to construct approximations to the area under the graph of $f$ through sums of areas of rectangles associated with each subinterval of the partition.

\begin{figure}[ht!]
\begin{minipage}[r]{0.5\textwidth}
\begin{center}
\begin{tikzpicture}[yscale=1.6]
\draw [thick, white, fill=lime!15]  (0.5,0)--(0.5,0.5)--(1,0.5)--(1,0);
\draw [thick, white, fill=lime!15]  (1,0)--(1,0.85)--(1.5,0.85)--(1.5,0);
\draw [thick, white, fill=lime!15]  (1.5,0)--(1.5,0.73)--(2,0.73)--(2,0);
\draw [thick, white, fill=lime!15]  (2,0)--(2,0.55)--(2.5,0.55)--(2.5,0);
\draw [thick, white, fill=lime!15]  (3.75,0)--(3.75,1.2)--(4.35,1.2)--(4.35,0);
\draw [thick, white, fill=lime!15]  (4.35,0)--(4.35,1.5)--(5,1.5)--(5,0);
\draw[->, gray] (-0.5,0) -- (6,0);
\draw [->, gray] (0,-0.1) -- (0,2);
\draw[thick] plot[smooth] coordinates
{(0.5,0.5)(1,1)(3,0.51)(4,1.5)(5,1.5)};

\draw[dotted] (0.5,0)--(0.5,0.5); \draw (0.5,0)node{$_{|}$}; \draw (0.5,-0.05) node[below]{$_{a}$};
\draw[dotted] (5,0)--(5,1.5); \draw (5,0)node{$_{|}$}; \draw (5,-0.05) node[below]{$_{b}$};

\draw[dotted] (1,0)--(1,0.85); \draw (1,0)node{$_{|}$}; \draw (1,-0.05) node[below]{$_{x_{1}}$};
\draw[dotted] (1.5,0)--(1.5,0.85); \draw (1.5,0)node{$_{|}$}; \draw (1.5,-0.05) node[below]{$_{x_{2}}$};

\draw[dotted] (2,0)--(2,0.75); \draw (2,0)node{$_{|}$}; \draw (2,-0.05) node[below]{$_{x_{3}}$};
\draw[dotted] (2.5,0)--(2.5,0.55); \draw (2.5,0)node{$_{|}$}; \draw (2.5,-0.05) node[below]{$_{x_{4}}$};
\draw (3.1,-0.05) node[below]{$\ldots$};
\draw[dotted] (3.75,0)--(3.75,1.2); \draw (3.75,0)node{$_{|}$}; \draw (3.73,-0.05) node[below]{$_{x_{n-2}}$};
\draw[dotted] (4.35,0)--(4.35,1.5); \draw (4.35,0)node{$_{|}$}; \draw (4.5,-0.05) node[below]{$_{x_{n-1}}$};

\draw[dotted] (0.5,0.5)--(1,0.5);
\draw[dotted] (1.5,0.85)--(1,0.85);
\draw[dotted] (2,0.73)--(1.5,0.73);
\draw[dotted] (2.5,0.55)--(2,0.55);
\draw[dotted] (4.35,1.2)--(3.75,1.2);
\draw[dotted] (5,1.5)--(4.35,1.5);
\end{tikzpicture}
\end{center}
\end{minipage} \hfill 
 \begin{minipage}[l]{0.5\textwidth}
\begin{center}
\begin{tikzpicture}[yscale=1.6]
\draw [thick, white, fill=lime!15]  (0.5,0)--(0.5,1)--(1,1)--(1,0);
\draw [thick, white, fill=lime!15]  (1,0)--(1,1)--(1.5,1)--(1.5,0);
\draw [thick, white, fill=lime!15]  (1.5,0)--(1.5,0.9)--(2,0.9)--(2,0);
\draw [thick, white, fill=lime!15]  (2,0)--(2,0.75)--(2.5,0.75)--(2.5,0);
\draw [thick, white, fill=lime!15]  (3.75,0)--(3.75,1.57)--(4.35,1.57)--(4.35,0);
\draw [thick, white, fill=lime!15]  (4.35,0)--(4.35,1.57)--(5,1.57)--(5,0);
\draw[->, gray] (-0.5,0) -- (6,0);
\draw [->, gray] (0,-0.1) -- (0,2);
\draw[thick] plot[smooth] coordinates
{(0.5,0.5)(1,1)(3,0.51)(4,1.5)(5,1.5)};

\draw[dotted] (0.5,0)--(0.5,1); \draw (0.5,0)node{$_{|}$}; \draw (0.5,-0.05) node[below]{$_{a}$};
\draw[dotted] (5,0)--(5,1.5); \draw (5,0)node{$_{|}$}; \draw (5,-0.05) node[below]{$_{b}$};

\draw[dotted] (1,0)--(1,1); \draw (1,0)node{$_{|}$}; \draw (1,-0.05) node[below]{$_{x_{1}}$};
\draw[dotted] (1.5,0)--(1.5,1); \draw (1.5,0)node{$_{|}$}; \draw (1.5,-0.05) node[below]{$_{x_{2}}$};

\draw[dotted] (2,0)--(2,0.9); \draw (2,0)node{$_{|}$}; \draw (2,-0.05) node[below]{$_{x_{3}}$};
\draw[dotted] (2.5,0)--(2.5,0.75); \draw (2.5,0)node{$_{|}$}; \draw (2.5,-0.05) node[below]{$_{x_{4}}$};
\draw (3.1,-0.05) node[below]{$\ldots$};
\draw[dotted] (3.75,0)--(3.75,1.57); \draw (3.75,0)node{$_{|}$}; \draw (3.73,-0.05) node[below]{$_{x_{n-2}}$};
\draw[dotted] (4.35,0)--(4.35,1.57); \draw (4.35,0)node{$_{|}$}; \draw (4.5,-0.05) node[below]{$_{x_{n-1}}$};

\draw[dotted] (0.5,1)--(1,1);
\draw[dotted] (1.5,1)--(1,1);
\draw[dotted] (2,0.9)--(1.5,0.9);
\draw[dotted] (2.5,0.75)--(2,0.75);
\draw[dotted] (4.35,1.57)--(3.75,1.57);
\draw[dotted] (5,1.57)--(4.35,1.57);
\end{tikzpicture}
\end{center}
\end{minipage}
\end{figure}

We define the \textbf{lower Riemann--Darboux sum} of $f$ associated with the partition $P$ by \index{sum!lower Riemann--Darboux}
$$
\underline{\mathcal{S}}(f,P)
:= \sum_{i=1}^n m_i\,(x_i-x_{i-1}).
$$

Similarly, we define the \textbf{upper Riemann--Darboux sum} of $f$ associated with the partition $P$ by \index{sum!upper Riemann--Darboux}
$$
\overline{\mathcal{S}}(f,P)
:= \sum_{i=1}^n M_i\,(x_i-x_{i-1}).
$$

It is not difficult to verify that the lower and upper sums satisfy the following properties\footnote{See, for example, \cite{Bartle1}.}:

\begin{proposition}\label{11}
Let $f:[a,b]\to\mathbb{R}$ be a nonnegative bounded function.
\begin{itemize}
\item[\textit{(a)}] $\underline{\mathcal{S}}(f,P)\leq \overline{\mathcal{S}}(f,P)$ for every $P\in\mathfrak{P}[a,b]$.

\item[\textit{(b)}] If $P,Q\in\mathfrak{P}[a,b]$ are such that $Q$ is a refinement of $P$ (that is, $P\subseteq Q$), then
$$
\underline{\mathcal{S}}(f,P)\leq \underline{\mathcal{S}}(f,Q)
\quad\text{and}\quad
\overline{\mathcal{S}}(f,Q)\leq \overline{\mathcal{S}}(f,P).
$$

\item[\textit{(c)}] $\underline{\mathcal{S}}(f,P)\leq \overline{\mathcal{S}}(f,Q)$.
\end{itemize}
\end{proposition}

We denote by $\underline{\mathcal{S}}(f)$ the set of all lower sums and by $\overline{\mathcal{S}}(f)$ the set of all upper sums of a nonnegative bounded function
$f:[a,b]\to\mathbb{R}$. That is,
$$
\begin{aligned}
\underline{\mathcal{S}}(f)
&:=\{\underline{\mathcal{S}}(f,P):P\in\mathfrak{P}[a,b]\},\\
\overline{\mathcal{S}}(f)
&:=\{\overline{\mathcal{S}}(f,P):P\in\mathfrak{P}[a,b]\}.
\end{aligned}
$$

Observe that both $\underline{\mathcal{S}}(f)$ and $\overline{\mathcal{S}}(f)$ are nonempty. This allows us to give the following definition.

\begin{definition}\label{12}
Let $f:[a,b]\to\mathbb{R}$ be a nonnegative bounded function. We define the \textbf{lower Riemann integral} of $f$ on $[a,b]$ as the supremum of the set $\underline{\mathcal{S}}(f)$, and we denote it by
$$
\underline{R}\int_a^b f(x)\,dx.
$$
\index{integral!lower Riemann}

The \textbf{upper Riemann integral} of $f$ on $[a,b]$ is defined as the infimum of the set $\overline{\mathcal{S}}(f)$, and we denote it by
$$
\overline{R}\int_a^b f(x)\,dx.
$$
\index{integral!upper Riemann}
\end{definition}

It can be shown that the lower and upper integrals of a nonnegative bounded function $f:[a,b]\to\mathbb{R}$ always exist and satisfy
$$
\underline{R}\int_a^b f(x)\,dx
\leq
\overline{R}\int_a^b f(x)\,dx.
$$

In these terms, we define the Riemann integral of a nonnegative bounded function on $[a,b]$.

\begin{definition}\label{13}
Let $f:[a,b]\to\mathbb{R}$ be a nonnegative bounded function. We say that $f$ is \textit{Riemann integrable} on $[a,b]$ \index{function!Riemann integrable} if and only if the lower and upper Riemann integrals of $f$ on $[a,b]$ coincide.

In this case, we denote the common value by
$$
R\int_a^b f(x)\,dx,
$$
and call it the \textbf{Riemann integral of $f$} on $[a,b]$. \index{integral!Riemann}
\end{definition}

The following result plays a fundamental role in the theory of the Riemann integral, since it shows that the integrability of a function can be characterized without explicitly computing the integral, but rather through control of the difference between the lower and upper sums associated with suitable partitions. We assume that the reader is familiar with this result, and therefore we omit its proof\footnote{See, for example, \cite{Bartle1, Rudin2}.}.

\begin{theorem}\label{14}
Let $f:[a,b]\to\mathbb{R}$ be a nonnegative bounded function. Then $f$ is Riemann integrable on $[a,b]$ if and only if for every $\varepsilon>0$ there exists a partition $P\in\mathfrak{P}[a,b]$ (depending on $\varepsilon$) such that
$$
\overline{\mathcal{S}}(f,P)-\underline{\mathcal{S}}(f,P)<\varepsilon.
$$
\end{theorem}

An immediate consequence of the previous criterion is given by the following result, which provides an abundance of examples of Riemann integrable functions.

\begin{theorem}\label{15}
Let $f:[a,b]\to\mathbb{R}$ be a nonnegative continuous function on $[a,b]$. Then $f$ is Riemann integrable on $[a,b]$.
\end{theorem}

\begin{proof}
Since $f$ is continuous on $[a,b]$ and $[a,b]$ is a compact subset of $\mathbb{R}$, it follows that $f$ is uniformly continuous on $[a,b]$. Thus, given $\varepsilon>0$, there exists $\delta>0$ such that for every $x,y\in[a,b]$,
$$
|f(x)-f(y)|<\tfrac{\varepsilon}{2(b-a)}
\qquad \text{whenever}\quad
|x-y|<\delta.
$$

It is possible to construct a partition $P\in\mathfrak{P}[a,b]$ such that $\Vert P\Vert<\delta$. Hence, $|x_i-x_{i-1}|\leq \Vert P\Vert < \delta$ for every $i=1,\ldots,n$.

For each $i=1,\ldots,n$, the function $f$ is continuous on the interval $[x_{i-1},x_i]$, and therefore there exist points $x_i^{\ast},y_i^{\ast}\in[x_{i-1},x_i]$ such that $f(x_i^{\ast})\leq f(x)\leq f(y_i^{\ast})$ for every $x\in[x_{i-1},x_i]$. Consequently,
$$
\begin{aligned}
\overline{\mathcal{S}}(f,P)
&=
\sum_{i=1}^{n} M_i(x_i-x_{i-1})
=
\sum_{i=1}^{n} f(y_i^{\ast})(x_i-x_{i-1}),\\
\underline{\mathcal{S}}(f,P)
&=
\sum_{i=1}^{n} m_i(x_i-x_{i-1})
=
\sum_{i=1}^{n} f(x_i^{\ast})(x_i-x_{i-1}).
\end{aligned}
$$

Therefore,
$$
\begin{aligned}
\overline{\mathcal{S}}(f,P)-\underline{\mathcal{S}}(f,P)
&=
\sum_{i=1}^{n}
\bigl(f(y_i^{\ast})-f(x_i^{\ast})\bigr)
(x_i-x_{i-1})\\
&<
\sum_{i=1}^{n}
\frac{\varepsilon}{2(b-a)}
(x_i-x_{i-1})\\
&=
\frac{\varepsilon}{2(b-a)}
\sum_{i=1}^{n}(x_i-x_{i-1})\\
&=
\frac{\varepsilon}{2}
<
\varepsilon.
\end{aligned}
$$

This proves the result.
\end{proof}

A natural question to ask at this point is whether every nonnegative bounded function is Riemann integrable. The answer is negative, and we shall present the classical counterexample: \emph{the Dirichlet function\footnote{Johann Peter Gustav Lejeune Dirichlet (1805--1859) was a German mathematician to whom the modern definition of a function is commonly attributed. He was educated in Germany and later in France, where he studied with several of the most renowned mathematicians of the time, including Fourier. His methods introduced an entirely new perspective, and his results are among the most important in mathematics.}}. \index{function!Dirichlet}

\begin{proposition}\label{16}
There exists a nonnegative bounded function on $[0,1]$ that is not Riemann integrable on $[0,1]$.
\end{proposition}

\begin{figure}[ht!]
\centering
\includegraphics[scale=0.25]{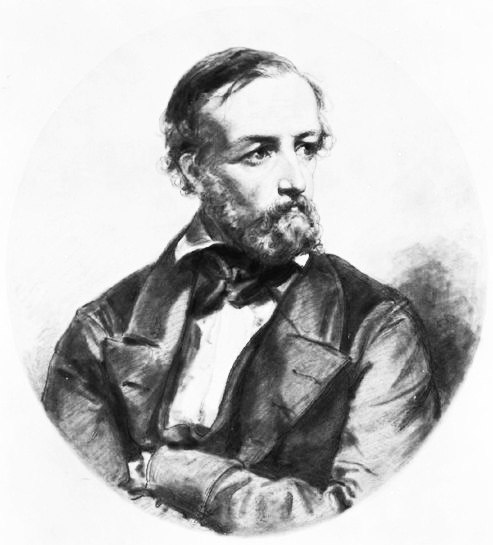} 
\begin{center}
J. Dirichlet (1805-1859)
\end{center}
\end{figure}

\begin{proof}
The Dirichlet function $D:[0,1]\to\mathbb{R}$ is defined by
$$
D(x):=
\left\{
\begin{array}{lcl}
1 && \text{if}\quad x\in\widetilde{\mathbb{Q}},\\
0 && \text{if}\quad x\notin\widetilde{\mathbb{Q}},
\end{array}
\right.
$$
where $\widetilde{\mathbb{Q}}:=[0,1]\cap\mathbb{Q}$.

Clearly, $D$ is nonnegative and bounded on $[0,1]$.

We claim that $D$ is not Riemann integrable on $[0,1]$. Let $P=\{x_0,x_1,\ldots,x_n\}$ be a partition of $[0,1]$. Since both $\widetilde{\mathbb{Q}}$ and its complement $[0,1]\smallsetminus\widetilde{\mathbb{Q}}$ are dense in $[0,1]$, it follows that for every $i=1,\ldots,n$,
$$
[x_{i-1},x_i]\cap\widetilde{\mathbb{Q}}\neq\varnothing
\quad\text{and}\quad
[x_{i-1},x_i]\cap([0,1]\smallsetminus\widetilde{\mathbb{Q}})
\neq\varnothing.
$$

Consequently,
$$
\begin{aligned}
m_i
&=
\inf\{D(x):x\in[x_{i-1},x_i]\}
=
\inf\{0,1\}
=
0,\\
M_i
&=
\sup\{D(x):x\in[x_{i-1},x_i]\}
=
\sup\{0,1\}
=
1.
\end{aligned}
$$

Hence,
$$
\begin{aligned}
\overline{\mathcal{S}}(D,P)
&=
\sum_{i=1}^{n}M_i(x_i-x_{i-1})
=
\sum_{i=1}^{n}1\cdot(x_i-x_{i-1})
=
1,\\
\underline{\mathcal{S}}(D,P)
&=
\sum_{i=1}^{n}m_i(x_i-x_{i-1})
=
\sum_{i=1}^{n}0\cdot(x_i-x_{i-1})
=
0.
\end{aligned}
$$

Since $P$ is an arbitrary partition of $[0,1]$, we conclude that
$$
\underline{\mathcal{S}}(D)=0
\quad\text{and}\quad
\overline{\mathcal{S}}(D)=1.
$$

Therefore,
$$
\underline{R}\int_{0}^{1}D(x)\,dx
=
0
\neq
1
=
\overline{R}\int_{0}^{1}D(x)\,dx,
$$
which implies that $D$ is not Riemann integrable on $[0,1]$.
\end{proof}

The question we must now ask is the following: \textit{what is it that fails in the construction of the Riemann integral that prevents us from computing the integral of a function such as the Dirichlet function?}

The first thing to observe is that, because of the way this function is defined, no matter how fine a partition of the interval $[0,1]$ may be, the function attains both values of its image on each of its subintervals. Consequently, the values of the infimum and supremum are always equal to $0$ and $1$, respectively.

\begin{figure}[ht!]
\centering
\begin{tikzpicture}[scale=3.5]

\draw[->, gray] (0,0) -- (1.1,0);

\draw[->,gray] (0,-0.1) -- (0,1.2);

\draw (0,1) node{${-}$};
\draw (1,0) node{$_{|}$};
\draw (0,1) node[left]{$_{1}$};
\draw (0,0) node[left]{$_{0}$};
\draw (1,-0.01) node[below]{$_{1}$};
\foreach \x in {0,0.2,0.4,0.6,0.8,1} {
    \draw[dashed,gray] (\x,0) -- (\x,1);
}

\foreach \x in {0.05,0.15,0.25,0.35,0.45,0.55,0.65,0.75,0.85,0.95} {
    \fill[blue] (\x,0) circle (0.5pt);
}
\foreach \x in {0.1,0.2,0.3,0.4,0.5,0.6,0.7,0.8,0.9} {
    \fill[red] (\x,1) circle (0.5pt);
}

\end{tikzpicture}
\end{figure}

This behavior is not due to the values taken by the function themselves, but rather to the way in which those values are distributed throughout the domain. The value $1$ is attained on a countable subset of the interval, whereas the value $0$ is attained on its complement, which is intuitively ``much larger.'' The Riemann construction is unable to distinguish between these sets, since every partition of the domain contains both rational and irrational points.

This observation suggests that, in order to integrate functions with extreme discontinuities, it is not sufficient merely to refine the partition of the domain; rather, it is necessary to incorporate information about the ``size'' of the sets on which the function takes certain values. In other words, the problem is not which values the function takes, but how many points in the domain correspond to each value.

From this perspective, it becomes natural to consider an alternative approach to integration: instead of subdividing the domain into subintervals, one partitions the set of values taken by the function and analyzes the subsets of the domain associated with each of them. That is, one studies the inverse images of the function and systematically assigns a notion of ``size'' to such sets.

It should be emphasized that this change of perspective does not depend essentially on the Euclidean structure of the interval $[a,b]$. Indeed, the idea of grouping points in the domain according to the values taken by a function and assigning an appropriate notion of ``size'' to the corresponding subsets can be formulated for real-valued functions defined on an arbitrary nonempty set $X$.

\begin{figure}[ht!]
\centering
\begin{tikzpicture}[xscale=1.5,yscale=1.5]
\draw[-, gray] (-3.5,0) -- (3.5,0); 
\draw (3.5,0) node[right]{$_{X}$};
\draw[->, gray] (0,-0.15) -- (0,2);
\draw (0,2) node[right]{$_{\mathbb{R}^{\geq0}}$};
\draw[ultra thick] (3,0) arc[start angle=0, end angle=180, x radius=3, y radius=1.5];

\draw[-,dashed, gray] (-3.5,0.25)--(3.5,0.25);

\draw[-,dashed,gray] (-3.5,0.5)--(3.5,0.5);

\draw[-,dashed,gray] (-3.5,0.75)--(3.5,0.75);

\draw[-,dashed,gray] (-3.5,1)--(3.5,1);

\draw[-,dashed,gray] (-3.5,1.25)--(3.5,1.25);

\draw[-,dashed,gray] (-3.5,1.6)--(3.5,1.6);

\draw[-,red,thick] (-1.7,1.25)--(1.65,1.25);

\draw[-,red,thick] (-2.2,1)--(-1.7,1);
\draw[-,red,thick] (1.65,1)--(2.2,1);

\draw[-,red,thick] (-2.6,0.75)--(-2.2,0.75);
\draw[-,red,thick] (2.2,0.75)--(2.6,0.75);

\draw[-,red,thick] (-2.8,0.5)--(-2.6,0.5);
\draw[-,red,thick] (2.6,0.5)--(2.8,0.5);

\draw[-,red,thick] (-3,0.25)--(-2.8,0.25);
\draw[-,red,thick] (2.8,0.25)--(3,0.25);

\draw[-,red,thick] (-3.5,0)--(-3,0);
\draw[-,red,thick] (3,0)--(3.5,0);

\draw[-,dashed] (-1.7,1.25)--(-1.7,1);
\draw[-,dashed] (1.65,1.25)--(1.65,1);

\draw[-,dashed] (-2.2,1)--(-2.2,0.75);
\draw[-,dashed] (2.2,1)--(2.2,0.75);

\draw[-,dashed] (-2.6,0.75)--(-2.6,0.5);
\draw[-,dashed] (2.6,0.75)--(2.6,0.5);

\draw[-,dashed] (-2.8,0.5)--(-2.8,0.25);
\draw[-,dashed] (2.8,0.5)--(2.8,0.25);
\end{tikzpicture}
\end{figure}

This suggests that the fundamental problem is not, in the first instance, the definition of a new notion of integration, but rather the construction of a coherent concept of ``size'' for certain subsets of an arbitrary nonempty set $X$. It is precisely \textit{measure theory} that provides the appropriate formal framework for addressing this problem.

In this context, the first objective will be to study those classes of subsets of $X$ for which the intuitive notion of ``size'' can be rigorously formalized. These subsets will constitute the so-called \textit{measurable sets}, and the analysis of their properties will serve as the starting point for the development of the theory.

Once this structure has been established, we shall consider real-valued functions defined on $X$ that are compatible with it, namely, those whose inverse images of measurable sets are themselves measurable sets. These functions, called \textit{measurable functions}, will form the natural class on which it becomes possible to define a new notion of integration. Within this framework, the \textit{Lebesgue integral} will be constructed. It will constitute the main object of study in this first part of the text and represents the most natural extension of the Riemann integral to considerably more general spaces and functions.

\section{The integral and limiting processes}

The fact that the Dirichlet function is not Riemann integrable conceals a deeper difficulty related to the problem of integration. Beyond constituting an isolated pathological example, this phenomenon reveals a central issue in the historical development of analysis: the behavior of the integral with respect to limiting processes.

The study of the convergence of sequences of functions occupied a central position in the development of mathematical analysis during the nineteenth century. The need to rigorously justify the interchange between the limit and integral symbols led to the introduction of different notions of convergence, such as pointwise convergence and uniform convergence, each having different implications for continuity, integrability, and the global behavior of the functions involved.

These notions are not only fundamental from a theoretical point of view, but also play an essential role in various areas of analysis and its applications, such as the study of Fourier series, the solution of differential equations, and the analysis of approximation methods. For these reasons, before proceeding toward more general notions of integration, it is convenient to examine in some detail the different types of convergence of sequences of functions and their interaction with the Riemann integral. In particular, we are interested in addressing the following fundamental problem of integration:

\begin{problem}\label{17}
\it Let $[a,b]$ be an interval in $\mathbb{R}$. Is the limit of a sequence of integrable functions $f_k:[a,b]\to\mathbb{R}$ itself integrable on $[a,b]$? If so, does the following equality hold
$$
\lim_{k\to\infty} R\int_a^b f_k(x)\,dx
=
R\int_a^b \lim_{k\to\infty} f_k(x)\,dx ?
$$
\end{problem}

To address the previous problem within the framework of the Riemann integral, it is necessary to carefully specify what is meant by the convergence of a sequence of functions. Indeed, unlike the case of sequences of real numbers, there exist several notions of convergence for functions, each of which captures different aspects of the behavior of the sequence and has different consequences for continuity, integrability, and the possibility of interchanging limits with analytic operations.

In what follows, we introduce two fundamental notions of convergence for sequences of real-valued functions: pointwise convergence and uniform convergence.

\begin{definition}\label{18}
Let $X$ be a nonempty set. We say that a sequence of functions $f_k:X\to\mathbb{R}$ \textbf{converges pointwise} \index{convergence!pointwise} on $X$ to a function $f:X\to\mathbb{R}$ if $f_k(x)\to f(x)$ for every $x\in X$. That is, if for every $\varepsilon>0$ and every $x\in X$, there exists $k_0\in\mathbb{N}$ (depending on $\varepsilon$ and on $x$) such that
$$
|f_k(x)-f(x)|<\varepsilon
\qquad \forall\, k\geq k_0.
$$

The function $f$ is called the \textbf{pointwise limit} of the sequence $(f_k)$ and is usually denoted by $\lim_{k\to\infty}f_k(x)$.
\end{definition}

\begin{definition}\label{19}
Let $X$ be a nonempty set. We say that a sequence of functions $f_k:X\to\mathbb{R}$ \textbf{converges uniformly} \index{convergence!uniform} on $X$ to a function $f:X\to\mathbb{R}$ if, given
$\varepsilon>0$, there exists $k_0\in\mathbb{N}$ (depending on $\varepsilon$) such that
$$
|f_k(x)-f(x)|<\varepsilon
\qquad \forall\, k\geq k_0,\ \forall\, x\in X.
$$

The function $f$ is called the \textbf{uniform limit} of the sequence $(f_k)$.
\end{definition}

Some basic properties and illustrative examples of these notions of convergence are presented in the exercises section.

The following example shows that, in the context of the Riemann integral, the first question raised in Problem~\ref{17} does not, in general, have an affirmative answer.

\begin{example}\label{110}
There exists a sequence of Riemann integrable functions on $[0,1]$ whose pointwise limit is not a Riemann integrable function on $[0,1]$.
\end{example}

\begin{proof}
Let $\widetilde{\mathbb{Q}}$ denote the set of all rational numbers contained in the interval $[0,1]$. It is well known that $\widetilde{\mathbb{Q}}$ is a countable set. Fix $q_1:=0$ and consider an enumeration $\{q_1,q_2,\ldots\}$ of the elements of $\widetilde{\mathbb{Q}}$, all distinct from one another.

Define $D_1:[0,1]\to\mathbb{R}$ by
$$
D_1(x):=
\left\{
\begin{array}{lcl}
1 && \text{if}\quad x=q_1,\\
0 && \text{if}\quad x\neq q_1.
\end{array}
\right.
$$

\begin{figure}[ht!]
\centering
\begin{tikzpicture}[xscale=1.5,yscale=1]
	\draw[->, gray] (0,0) -- (4.25,0); \draw [->, gray] (0,0) -- (0,3.45); 
	\draw[-, gray] (-0.23,0)--(0,0); \draw[-,gray] (0,0)--(0,-0.23);
	
	\draw[-, ultra thick] (0.02,0)--(4,0);
	
\draw (0,3) node{$_{-}$}; \draw (0,3) node[left]{$_{1}$};
\draw (4,0) node{$_{|}$}; \draw (4,-0.1) node[below]{$_{1}$};
\draw (0,-0.1) node[below]{$_{0}$};
\draw (0,3) node{$_{\bullet}$};
\draw (0,0) node{$_{\circ}$};

\draw (1,0) node{$_{|}$};
\draw (1,-0.1) node[below]{$_{\delta}$};

\draw (2,0) node{$_{|}$};
\draw (2,-0.1) node[below]{$_{\frac{\varepsilon}{2}}$};
\end{tikzpicture}
\begin{center}
$D_1(x)$
\end{center}
\end{figure}

Let $\varepsilon\in(0,1)$. Choose $\delta\in(0,1)$ such that
$\delta<\frac{\varepsilon}{2}$.

Defining the partition
$P_{\delta}:=\{0<\delta<1\}$
of the interval $[0,1]$, we have
$$
\overline{\mathcal{S}}(D_1,P_{\delta})=\delta
\quad\text{and}\quad
\underline{\mathcal{S}}(D_1,P_{\delta})=0,
$$
and therefore
$$
\overline{\mathcal{S}}(D_1,P_{\delta})
-
\underline{\mathcal{S}}(D_1,P_{\delta})
=
\delta
<
\varepsilon.
$$

That is, $D_1$ is Riemann integrable on $[0,1]$.

Now define $D_2:[0,1]\to\mathbb{R}$ by
$$
D_2(x):=
\left\{
\begin{array}{lcl}
1 && \text{if}\quad x\in\{q_1,q_2\},\\
0 && \text{if}\quad x\notin\{q_1,q_2\}.
\end{array}
\right.
$$

Let $\varepsilon\in(0,1)$. Choose $\delta_1,\delta_2\in(0,1)$ such that
$\delta_1+\delta_2<\frac{\varepsilon}{2}$
and
$\delta_2<\min\left\{\frac{q_2}{2},\frac{1-q_2}{2}\right\}$.

Define the partitions
$P_{\delta_1}:=\{0<\delta_1<1\}$
and
$P_{\delta_2}:=\{0<q_2-\delta_2<q_2+\delta_2<1\}$
of $[0,1]$.

Let
$P:=P_{\delta_1}\cup P_{\delta_2}\in\mathfrak{P}[0,1]$.

Then
$$
\overline{\mathcal{S}}(D_2,P)
-
\underline{\mathcal{S}}(D_2,P)
\leq
\overline{\mathcal{S}}(D_2,P_{\delta_1})
-
\underline{\mathcal{S}}(D_2,P_{\delta_1})
=
\delta_1
<
\varepsilon.
$$

Consequently, $D_2$ is Riemann integrable on $[0,1]$.

\begin{figure}[ht!]
    \centering
    \begin{tikzpicture}[xscale=1.5,yscale=1]
	\draw[->, gray] (0,0) -- (4.25,0); \draw [->, gray] (0,0) -- (0,3.45); 
	\draw[-, gray] (-0.23,0)--(0,0); \draw[-,gray] (0,0)--(0,-0.23);
	
	\draw[-, ultra thick] (0.02,0)--(4,0);
	
\draw (0,3) node{$_{-}$}; \draw (0,3) node[left]{$_{1}$};
\draw (4,0) node{$_{|}$}; \draw (4,-0.1) node[below]{$_{1}$};
\draw (0,-0.1) node[below]{$_{0}$};
\draw (0,3) node{$_{\bullet}$};
\draw (0,0) node{$_{\circ}$};

\draw (0.25,0) node{$_{|}$}; \draw (0.25,-0.1) node[below]{$_{q_2}$};
\draw (0.25,3) node{$_{\bullet}$};
\draw (0.25,0) node{$_{\circ}$};
\end{tikzpicture}
\begin{center}
$D_2(x)$
\end{center}
\end{figure}

Let $N\in\mathbb{N}$ with $N>1$ be arbitrary but fixed. Define $D_N:[0,1]\to\mathbb{R}$ by
$$
D_N(x):=
\left\{
\begin{array}{lcl}
1 && \text{if}\quad x\in\{q_1,q_2,\ldots,q_N\},\\
0 && \text{if}\quad x\notin\{q_1,q_2,\ldots,q_N\}.
\end{array}
\right.
$$

\begin{figure}[ht!]
    \centering
    \begin{tikzpicture}[xscale=1.5,yscale=1]
	\draw[->, gray] (0,0) -- (4.25,0); \draw [->, gray] (0,0) -- (0,3.45); 
	\draw[-, gray] (-0.23,0)--(0,0); \draw[-,gray] (0,0)--(0,-0.23);
	
	\draw[-, ultra thick] (0.02,0)--(4,0);
	
\draw (0,3) node{$_{-}$}; \draw (0,3) node[left]{$_{1}$};
\draw (4,0) node{$_{|}$}; \draw (4,-0.1) node[below]{$_{1}$};
\draw (0,-0.1) node[below]{$_{0}$};
\draw (0,3) node{$_{\bullet}$};
\draw (0,0) node{$_{\circ}$};

\draw (0.25,0) node{$_{|}$}; \draw (0.25,-0.1) node[below]{$_{q_2}$};
\draw (0.25,3) node{$_{\bullet}$};
\draw (0.25,0) node{$_{\circ}$};

\draw (0.5,0) node{$_{|}$}; \draw (0.5,-0.1) node[below]{$_{q_3}$};
\draw (0.5,3) node{$_{\bullet}$};
\draw (0.5,0) node{$_{\circ}$};

\draw (0.75,0) node{$_{|}$}; \draw (0.75,-0.1) node[below]{$_{q_4}$};
\draw (0.75,3) node{$_{\bullet}$};
\draw (0.75,0) node{$_{\circ}$};

\draw (1,0) node{$_{|}$}; \draw (1,-0.1) node[below]{$_{q_5}$};
\draw (1,3) node{$_{\bullet}$};
\draw (1,0) node{$_{\circ}$};

\draw (1.5,0) node{$_{\ldots}$}; \draw (1.5,-0.1) node[below]{$_{\ldots}$};
\draw (1.5,3) node{$_{\ldots}$};

\draw (2,0) node{$_{|}$}; 
\draw (2,3) node{$_{\bullet}$};
\draw (2,0) node{$_{\circ}$};

\draw (2.25,0) node{$_{|}$}; 
\draw (2.25,3) node{$_{\bullet}$};
\draw (2.25,0) node{$_{\circ}$};

\draw (2.5,0) node{$_{|}$}; \draw (2.5,-0.1) node[below]{$_{q_{N}}$};
\draw (2.5,3) node{$_{\bullet}$};
\draw (2.5,0) node{$_{\circ}$};

\end{tikzpicture}
\begin{center}
$D_N(x)$
\end{center}
\end{figure}

Let $\varepsilon\in(0,1)$. For each $j=1,\ldots,N$, choose $\delta_j\in(0,1)$ such that
$\delta_j\leq\frac{\varepsilon}{2^{N+1}}$
and
$\delta_j<\min\left\{\frac{q_j}{2},\frac{1-q_j}{2}\right\}$. It is then clear that
$\delta_1+\cdots+\delta_N<\frac{\varepsilon}{2}$. Now define, for each $j=1,\ldots,N$, the partition
$P_j:=\{0<q_j-\delta_j<q_j+\delta_j<1\}$
of $[0,1]$. Defining
$P:=P_1\cup\cdots\cup P_N\in\mathfrak{P}[0,1]$,
we obtain
$$
\overline{\mathcal{S}}(D_N,P)
-
\underline{\mathcal{S}}(D_N,P)
\leq
\overline{\mathcal{S}}(D_N,P_1)
-
\underline{\mathcal{S}}(D_N,P_1)
=
\delta_1
<
\varepsilon.
$$

That is, $D_N$ is Riemann integrable on $[0,1]$.

Therefore, $(D_k)$ is a sequence of Riemann integrable functions on $[0,1]$ that converges pointwise to the Dirichlet function $D$ on $[0,1]$. Indeed, given $x\in[0,1]$, if $x$ is rational, then there exists $k_0$ such that $x=q_{k_0}$, and hence $D_k(x)=1$ for every $k\geq k_0$, whereas if $x$ is irrational, then $D_k(x)=0$ for every $k\geq1$. In both cases,
$D_k(x)\to D(x)$ as $k\to\infty$. However, $D$ is not Riemann integrable on $[0,1]$.
\end{proof}

The previous example shows that pointwise convergence of Riemann integrable functions does not guarantee the integrability of the limit. Concerning the second question raised in Problem~\ref{17}, it is possible to construct a sequence of Riemann integrable functions that converges pointwise to a function that is also integrable, but for which the limit and the integral cannot be interchanged. We propose this as an exercise [Exercise~\ref{Ej15}].

The difficulties encountered in the study of the Riemann integral show that its limitations are not restricted to the description of the sets on which the function takes certain values, but also affect its behavior with respect to the pointwise convergence of sequences of functions.

Since the introduction of the Lebesgue integral was announced in the previous section, it is natural to ask to what extent this new notion makes it possible to overcome these limitations. The second objective of this first part will be, within the framework of the Lebesgue integral, to guarantee the integrability of the pointwise limit of a sequence of integrable functions and to establish the conditions under which it is possible to interchange the limit and integration symbols.

\section{Exercises}

\begin{exercise}\label{Ej11}
Let $X$ be a nonempty set, let $f_k:X\to\mathbb{R}$ be a sequence of functions, and let $f:X\to\mathbb{R}$ be a function.
\begin{itemize}
\item[(a)] Prove that if $(f_k)$ converges uniformly to $f$ on $X$, then $(f_k)$ converges pointwise to $f$ on $X$.

\item[(b)] Give an example of a sequence of functions that converges pointwise to a function on $X$ but does not converge uniformly on $X$.
\end{itemize}
\end{exercise}

\begin{exercise}\label{Ej12}
Let $X$ be a nonempty subset of $\mathbb{R}$, let $f_k:X\to\mathbb{R}$ be a sequence of functions, and let $f:X\to\mathbb{R}$ be a function.
\begin{itemize}
\item[(a)] Prove that if $(f_k)$ converges uniformly to $f$ on $X$ and each $f_k$ is continuous on $X$, then $f$ is continuous on $X$.

\item[(b)] Give an example of a sequence of continuous functions that converges pointwise to a discontinuous function.
\end{itemize}
\end{exercise}

Let $X$ be a nonempty set and let $f_k:X\to\mathbb{R}$ be a sequence of functions. We say that the sequence $(f_k)$ is \textbf{uniformly Cauchy} on $X$ if, for every $\varepsilon>0$, there exists $k_0\in\mathbb{N}$ such that \index{sequence!uniformly Cauchy}
$$
|f_j(x)-f_k(x)|<\varepsilon
\qquad
\forall\, k,j\geq k_0,\quad \forall\, x\in X.
$$

\begin{exercise}[Uniform Cauchy Convergence Criterion]\index{criterion!uniform Cauchy convergence}\label{Ej13}
Let $X$ be a nonempty set and let $f_k:X\to\mathbb{R}$ be a sequence of functions. Prove that $(f_k)$ converges uniformly on $X$ if and only if $(f_k)$ is uniformly Cauchy on $X$.
\end{exercise}

{\setlength{\parindent}{0pt}

\begin{exercise}\label{Ej14}
Let $f_k:[a,b]\to\mathbb{R}$ be a sequence of Riemann integrable functions such that $(f_k)$ converges uniformly to a function $f:[a,b]\to\mathbb{R}$. Prove that $f$ is Riemann integrable and that the following equality holds:
$$
\lim_{k\to\infty}R\int_a^b f_k(x)\,dx
=
R\int_a^b \lim_{k\to\infty}f_k(x)\,dx.
$$
\end{exercise}

(Hint: Use the Uniform Cauchy Convergence Criterion.)}

\begin{exercise}\label{Ej15}
On $[0,1]$, consider the sequence of functions $f_k:[0,1]\to\mathbb{R}$ defined by
$$
f_k(x):=
\left\{
\begin{array}{lcl}
2k^2x && \text{if}\quad 0\leq x\leq \frac{1}{2k},\\
2k(1-kx) && \text{if}\quad \frac{1}{2k}\leq x\leq \frac{1}{k},\\
0 && \text{if}\quad \frac{1}{k}\leq x\leq 1.
\end{array}
\right.
$$

\begin{itemize}
\item[(a)] Prove that $f_k$ is Riemann integrable and that
$$
R\int_0^1 f_k(x)\,dx=\frac{1}{2}
$$
for every $k\in\mathbb{N}$.

\item[(b)] Prove that $f_k\to0$ pointwise on $[0,1]$.

\item[(c)] Conclude that
$$
\lim_{k\to\infty}R\int_0^1 f_k(x)\,dx
\neq
R\int_0^1 \lim_{k\to\infty}f_k(x)\,dx.
$$
\end{itemize}
\end{exercise}
\chapter{Classes of subsets}\label{Chapter2}

\markboth{{\scriptsize 2. CLASSES OF SUBSETS}}
{{\scriptsize 2. CLASSES OF SUBSETS}}

Set theory constitutes one of the fundamental pillars of mathematics, since it underlies a large part of its structures and developments. In particular, it is essential to study more sophisticated families of sets for which a meaningful notion of measurability can be established. We begin by studying certain classes of subsets of a given nonempty set, upon which we shall introduce the fundamental concepts of set algebra.

Just as, in the theory of metric spaces, a topology determines the class of sets on which notions such as continuity and compactness are developed, in measure theory a $\sigma$-algebra specifies the collection of subsets for which the concept of measure can be defined coherently. In this chapter we shall show that every nonempty set can be associated with a $\sigma$-algebra. When our underlying set is $\mathbb{R}$, we shall work with the Borel $\sigma$-algebra. The study of these structures has important applications, for example, in probability theory.

We shall also extend the concept of sequence to the setting of subsets of a given set, emphasizing its close relationship with the usual notion of a sequence. In particular situations, the use of such sequences will allow us to analyze the behavior of these structures under limiting processes.

Finally, we shall see that Dynkin classes provide, in certain situations, an effective method for verifying that a family of sets is a $\sigma$-algebra, a fact that is crucial in establishing a variety of important properties.

\section{Definitions and basic properties}

From this point onward, we shall consider a fixed nonempty set $X$ and denote by $\mathcal{P}(X)=\{A:A\subset X\}$ the class of all subsets of $X$.

\begin{definition}\label{21}
A nonempty class $\mathcal{R}\subset\mathcal{P}(X)$ is called a \textbf{ring} \index{ring} of subsets of $X$ if it satisfies the following conditions:
\begin{itemize}
\item[\rm(R1)] $A\smallsetminus B\in\mathcal{R}$ for every $A,B\in\mathcal{R}$.
\item[\rm(R2)] $A\cup B\in\mathcal{R}$ for every $A,B\in\mathcal{R}$.
\end{itemize}
\end{definition}

Let $\mathcal{R}$ be a ring of subsets of $X$. Since $\mathcal{R}$ is nonempty, there exists a subset $A$ of $X$ such that $A\in\mathcal{R}$. Using property (R1) from the previous definition, we obtain that $A\smallsetminus A=\varnothing\in\mathcal{R}$. That is, the empty set is always an element of a ring $\mathcal{R}$.

A characterization of the definition of a ring of subsets of $X$ is given in the following proposition.

\begin{proposition}\label{22}
A nonempty class $\mathcal{R}\subset\mathcal{P}(X)$ is a \textbf{ring} if and only if $\mathcal{R}$ is closed under intersections and symmetric differences.
\end{proposition}

\begin{proof}
Suppose that $\mathcal{R}\subset\mathcal{P}(X)$ is a ring. By Definition~\ref{21}, we have that $A\smallsetminus B$ and $A\cup B\in\mathcal{R}$ for every $A,B\in\mathcal{R}$. Thus, from the identities
$$
A\bigtriangleup B=(A\smallsetminus B)\cup(B\smallsetminus A)
\quad\text{and}\quad
A\cap B=(A\cup B)\smallsetminus(A\bigtriangleup B),
$$
it follows that $A\bigtriangleup B,A\cap B\in\mathcal{R}$ for every $A,B\in\mathcal{R}$.

Conversely, suppose that $\mathcal{R}\subset\mathcal{P}(X)$ is a class closed under intersections and symmetric differences. From the following set identities
$$
A\smallsetminus B=A\cap(A\bigtriangleup B)
\quad\text{and}\quad
A\cup B=A\bigtriangleup(B\smallsetminus A),
$$
we conclude that $\mathcal{R}$ is a ring of subsets of $X$.
\end{proof}

\begin{definition}\label{23}
A nonempty class $\mathcal{A}\subset\mathcal{P}(X)$ is called an \textbf{algebra} \index{algebra} of subsets of $X$ if it satisfies the following conditions:
\begin{itemize}
\item[\rm(A1)] $X\in\mathcal{A}$.
\item[\rm(A2)] $A\smallsetminus B\in\mathcal{A}$ for every $A,B\in\mathcal{A}$.
\item[\rm(A3)] $A\cup B\in\mathcal{A}$ for every $A,B\in\mathcal{A}$.
\end{itemize}
\end{definition}

It follows directly from properties (A1) and (A2) that, if $\mathcal{A}$ is an algebra of subsets of $X$, then $X\smallsetminus A\in\mathcal{A}$ for every $A\in\mathcal{A}$. Consequently, an algebra $\mathcal{A}$ of subsets of $X$ is closed under complements.

\begin{definition}\label{24}
A nonempty class $\mathcal{SR}\subset\mathcal{P}(X)$ is called a \textbf{$\sigma$-ring} \index{sigma@$\sigma$!ring} of subsets of $X$ if it satisfies the following conditions:
\begin{itemize}
\item[\,\,\rm(SR1)] $A\smallsetminus B\in\mathcal{SR}$ for every $A,B\in\mathcal{SR}$.
\item[\,\,\rm(SR2)] For every sequence $(A_k)$ of elements of $\mathcal{SR}$, we have $\bigcup_{k=1}^{\infty}A_k\in\mathcal{SR}$.
\end{itemize}
\end{definition}

We have the following result.

\begin{proposition}\label{25}
Every $\sigma$-ring $\mathcal{SR}$ of subsets of $X$ is a ring of subsets of $X$.
\end{proposition}

\begin{proof}
It suffices to verify property (R2). Since $\mathcal{SR}$ is a nonempty class, there exists $A\subset X$ such that $A\in\mathcal{SR}$, and from property (SR1) we conclude that
$$
A\smallsetminus A=\varnothing\in\mathcal{SR}.
$$

Let $A,B\in\mathcal{SR}$ be arbitrary and define the following sequence $(A_k)$ in $\mathcal{SR}$:
$$
\begin{aligned}
A_1&:=A,\\
A_2&:=B,\\
A_k&:=\varnothing \qquad \forall\, k\geq3.
\end{aligned}
$$

Applying property (SR2), we obtain
$$
\bigcup_{k=1}^{\infty}A_k=A\cup B\in\mathcal{SR},
$$
and therefore $\mathcal{SR}$ is a ring of subsets of $X$.
\end{proof}

It can be proved [Exercise~\ref{E21}] that a $\sigma$-ring of subsets of $X$ is also closed under countable intersections. However, contrary to Proposition~\ref{22}, a class of subsets of $X$ that is closed under countable intersections and symmetric differences is not necessarily a $\sigma$-ring [Exercise~\ref{E24}].

\begin{definition}\label{26} \index{set!measurable}
A nonempty class $\mathsf{S}\subset\mathcal{P}(X)$ is called a \textbf{$\sigma$-algebra} \index{sigma@$\sigma$!algebra} of subsets of $X$ if it satisfies the following conditions:
\begin{itemize}
\item[\rm(S1)] $X\in\mathsf{S}$.
\item[\rm(S2)] $A\smallsetminus B\in\mathsf{S}$ for every $A,B\in\mathsf{S}$.
\item[\rm(S3)] For every sequence $(A_k)$ of elements of $\mathsf{S}$, we have $\bigcup_{k=1}^{\infty}A_k\in\mathsf{S}$.
\end{itemize}

A \textbf{measurable space} is a nonempty set $X$ equipped with a $\sigma$-algebra $\mathsf{S}$. We denote it by $(X,\mathsf{S})$. \index{space!measurable}
\end{definition}

\begin{example}\label{27}
If $X$ is an arbitrary set containing more than one element, then $\mathcal{P}(X)$ contains at least two $\sigma$-algebras, namely,
\begin{itemize}
\item[(a)] $\mathsf{S}_1=\{\varnothing,X\}$.
\item[(b)] $\mathsf{S}_2=\mathcal{P}(X)$.
\end{itemize}

The $\sigma$-algebra in part (a) is the smallest $\sigma$-algebra that can be associated with the set $X$ and is known as the \textbf{trivial $\sigma$-algebra on $X$}, whereas the $\sigma$-algebra in part (b) is the largest one and is known as the \textbf{discrete $\sigma$-algebra on $X$}.
\end{example}

The following characterization of a $\sigma$-algebra of subsets of $X$ is, in some situations, more convenient for proofs.

\begin{proposition}\label{28}
A nonempty class $\mathsf{S}$ is a $\sigma$-algebra of subsets of $X$ if and only if it satisfies the following properties:
\begin{itemize}
\item[(a)] $X\in\mathsf{S}$.
\item[(b)] $X\smallsetminus A\in\mathsf{S}$ for every $A\in\mathsf{S}$.
\item[(c)] For every sequence $(A_k)$ of elements of $\mathsf{S}$, we have $\bigcup_{k=1}^{\infty}A_k\in\mathsf{S}$.
\end{itemize}
\end{proposition}

\begin{proof}
Suppose that $\mathsf{S}$ is a $\sigma$-algebra of subsets of $X$. By Definition~\ref{26}, it suffices to verify that $\mathsf{S}$ is closed under complements. Let $A\in\mathsf{S}$. From (S1) and (S2), we conclude that $X\smallsetminus A\in\mathsf{S}$.

Conversely, suppose that $\mathsf{S}\subset\mathcal{P}(X)$ satisfies \textit{(a)}, \textit{(b)}, and \textit{(c)}. It is immediate that $\varnothing\in\mathsf{S}$ and, consequently, that $\mathsf{S}$ is closed under finite unions.

We now show that $\mathsf{S}$ is closed under differences. Let $A,B\in\mathsf{S}$ be arbitrary. From the identity
$$
X\smallsetminus(A\smallsetminus B)=(X\smallsetminus A)\cup B,
$$
we conclude that $X\smallsetminus(A\smallsetminus B)\in\mathsf{S}$ and therefore $A\smallsetminus B\in\mathsf{S}$ by part \textit{(b)}.
\end{proof}

The following example can be proved using the previous proposition.

\begin{example}\label{29}
If $X$ is a set with exactly three elements, namely $X=\{a,b,c\}$, then
$$
\mathsf{S}=\{\varnothing,\{a,b\},\{c\},X\}
$$
is a $\sigma$-algebra of subsets of $X$.
\end{example}

By following the proof of Proposition~\ref{25}, one may conclude directly that every $\sigma$-algebra of subsets of $X$ is an algebra of subsets of $X$.

In general, it can be shown that a ring of subsets of $X$ is not necessarily a $\sigma$-ring and, similarly, that an algebra is not always a $\sigma$-algebra [Exercise~\ref{E218}]. However, one way to relate these classes of subsets of $X$ is the following: every ring $\mathcal{R}$ with finitely many elements is a $\sigma$-ring, and every algebra $\mathcal{A}$ with finitely many elements is a $\sigma$-algebra. The proofs of these statements are straightforward and are left as exercises [Exercise~\ref{E28}].

We now examine some properties relating different $\sigma$-algebras on the same set $X$.

\begin{proposition}\label{210}
Let $\mathcal{I}$ be an arbitrary nonempty set and let $\{\mathsf{S}_i:i\in\mathcal{I}\}$ be a family of $\sigma$-algebras of subsets of $X$. Then
$$
\mathsf{S}=\bigcap_{i\in\mathcal{I}}\mathsf{S}_i
$$
is a $\sigma$-algebra of subsets of $X$.
\end{proposition}

\begin{proof}
Since $X\in\mathsf{S}_i$ for every $i\in\mathcal{I}$, it follows that $X\in\mathsf{S}$. Let $A,B\in\mathsf{S}$. Then $A,B\in\mathsf{S}_i$ for every $i\in\mathcal{I}$ and, since each $\mathsf{S}_i$ is a $\sigma$-algebra, we have $A\smallsetminus B\in\mathsf{S}_i$ for every $i\in\mathcal{I}$. Consequently, $A\smallsetminus B\in\mathsf{S}$.

Finally, if $(A_k)$ is a sequence of elements of $\mathsf{S}$, then $(A_k)$ is a sequence of elements of $\mathsf{S}_i$ for every $i\in\mathcal{I}$. Hence,
$$
\bigcup_{k=1}^{\infty}A_k\in\mathsf{S}_i
$$
for every $i\in\mathcal{I}$ and therefore
$$
\bigcup_{k=1}^{\infty}A_k\in\mathsf{S}.
$$

This completes the proof.
\end{proof}

By the nature of the previous proposition, we may ask whether the union of $\sigma$-algebras is again a $\sigma$-algebra. In general, the answer is negative, as shown in the following example.

\begin{example}\label{211}
The union of two $\sigma$-algebras is not necessarily a $\sigma$-algebra.
\end{example}

\begin{proof}
Let $X=\{a,b,c\}$ and consider the $\sigma$-algebras
$$
\mathsf{S}_1=\{\varnothing,\{a\},\{b,c\},X\}
$$
and
$$
\mathsf{S}_2=\{\varnothing,\{a,b\},\{c\},X\}.
$$

Then
$$
\mathsf{S}_1\cup\mathsf{S}_2
=
\{\varnothing,\{a\},\{b,c\},\{a,b\},\{c\},X\}
$$
is not a $\sigma$-algebra of subsets of $X$, since it does not satisfy property (S3) of Definition~\ref{26}. Indeed, if we define the sequence of elements of $\mathsf{S}_1\cup\mathsf{S}_2$ by
$$
\begin{aligned}
A_1&:=\{a\},\\
A_2&:=\{c\},\\
A_k&:=\varnothing,\qquad \forall\, k\geq3,
\end{aligned}
$$
then
$$
\bigcup_{k=1}^{\infty}A_k
=
\{a\}\cup\{c\}
=
\{a,c\}
\notin
\mathsf{S}_1\cup\mathsf{S}_2.
$$
\end{proof}

The following proposition establishes a sufficient condition for the union of two $\sigma$-algebras to be a $\sigma$-algebra.

\begin{proposition}\label{212}
Let $\mathsf{S}_1$ and $\mathsf{S}_2$ be two $\sigma$-algebras of subsets of $X$ such that $\mathsf{S}_1\subset\mathsf{S}_2$. Then $\mathsf{S}_1\cup\mathsf{S}_2$ is a $\sigma$-algebra of subsets of $X$.
\end{proposition}

Its proof is quite straightforward and is therefore left as an exercise [Exercise~\ref{E212}].

\section{\texorpdfstring{$\sigma$}{}-algebra generated}

Example~\ref{27} provides at least two $\sigma$-algebras for an arbitrary nonempty set $X$ containing at least two elements. However, there is a way to construct different $\sigma$-algebras from an arbitrary class of subsets of $X$, from which we shall obtain interesting properties.

\begin{theorem}[$\sigma$-algebra generated]\label{213}
Let $\mathcal{C}\subset\mathcal{P}(X)$ be arbitrary. Then there exists a unique $\sigma$-algebra \index{sigma@$\sigma$!generated algebra!by a class $\sigma(\mathcal{C})$} $\sigma(\mathcal{C})\subset\mathcal{P}(X)$ satisfying the following two properties:
\begin{itemize}
\item[(a)] $\mathcal{C}\subset\sigma(\mathcal{C})$.
\item[(b)] If $\mathsf{S}$ is a $\sigma$-algebra of subsets of $X$ such that $\mathcal{C}\subset\mathsf{S}$, then $\sigma(\mathcal{C})\subset\mathsf{S}$.
\end{itemize}
\end{theorem}

\begin{proof}
Let
$$
\mathscr{F}
=
\left\{
\mathsf{S}\subset\mathcal{P}(X):
\mathsf{S}\ \text{is a}\ \sigma\text{-algebra and}\ \mathcal{C}\subset\mathsf{S}
\right\}.
$$

Observe that $\mathscr{F}$ is nonempty since $\mathcal{P}(X)\in\mathscr{F}$.

Define
$$
\sigma(\mathcal{C})
:=
\bigcap
\left\{
\mathsf{S}\subset\mathcal{P}(X):
\mathsf{S}\in\mathscr{F}
\right\}.
$$

\textit{(a):} Proposition~\ref{210} guarantees that $\sigma(\mathcal{C})$ is a $\sigma$-algebra of subsets of $X$, and clearly
$\mathcal{C}\subset\sigma(\mathcal{C})$
since
$\mathcal{C}\subset\mathsf{S}$
for every $\mathsf{S}\in\mathscr{F}$.

\textit{(b):} If $\mathsf{S}\subset\mathcal{P}(X)$ is a $\sigma$-algebra such that $\mathcal{C}\subset\mathsf{S}$, then $\mathsf{S}\in\mathscr{F}$ and therefore, by the property that the intersection of sets is always contained in each of the sets being intersected,
$\sigma(\mathcal{C})\subset\mathsf{S}$.

Finally, if $\widetilde{\mathsf{S}}$ is a $\sigma$-algebra of subsets of $X$ satisfying \textit{(a)} and \textit{(b)}, then
$\widetilde{\mathsf{S}}\subset\sigma(\mathcal{C})$
and
$\sigma(\mathcal{C})\subset\widetilde{\mathsf{S}}$.
This proves uniqueness.
\end{proof}

Property \textit{(b)} allows us to interpret $\sigma(\mathcal{C})$ as the smallest $\sigma$-algebra of subsets of $X$ containing the class $\mathcal{C}$. A simple example is the following.

{\setlength{\parindent}{0pt}
\begin{example}\label{214}
Let $A,B\subset X$ be such that $A\cap B=\varnothing$. Define the class $\mathcal{C}\subset\mathcal{P}(X)$ by
$$
\mathcal{C}:=\{A,B\}.
$$

In general, this collection is not a $\sigma$-algebra (why not?). However, we can enlarge it by adding certain subsets of $X$ in order to obtain the $\sigma$-algebra generated by $\mathcal{C}$.

To construct a $\sigma$-algebra containing $\mathcal{C}$, we use Proposition~\ref{28} and Theorem~\ref{213}. Let $\mathsf{S}$ be a $\sigma$-algebra of subsets of $X$ such that $\mathcal{C}\subset\mathsf{S}$. Then $A,B\in\mathsf{S}$ and, by definition of $\sigma$-algebra, also $X\in\mathsf{S}$. By Proposition~\ref{28}, the union of elements of $\mathsf{S}$ belongs to $\mathsf{S}$, and the complement of every element of $\mathsf{S}$ also belongs to $\mathsf{S}$. In particular,
$$
A\cup B,\varnothing,X\smallsetminus A,X\smallsetminus B,
X\smallsetminus(A\cup B)\in\mathsf{S}.
$$

Therefore, the collection
$$
\mathcal{F}
=
\{X,A,B,A\cup B,X\smallsetminus A,
X\smallsetminus B,
X\smallsetminus(A\cup B),\varnothing\}
$$
is the smallest $\sigma$-algebra of subsets of $X$ containing $\mathcal{C}$. We note that the assumption $A\cap B=\varnothing$ is essential throughout the argument (why?; see Exercise~\ref{E222}).
\end{example}}

We now examine some important properties of the $\sigma$-algebra generated by an arbitrary class.

\begin{theorem}\label{215}
Let $\mathcal{C},\mathcal{C}'\subset\mathcal{P}(X)$ be arbitrary. The following properties hold:
\begin{itemize}
\item[(a)] If $\mathcal{C}\subset\mathcal{C}'$, then $\sigma(\mathcal{C})\subset\sigma(\mathcal{C}')$.
\item[(b)] If $\mathcal{C}$ is a $\sigma$-algebra, then $\sigma(\mathcal{C})=\mathcal{C}$.
\item[(c)] $\sigma(\sigma(\mathcal{C}))=\sigma(\mathcal{C})$.
\item[(d)] If $\mathcal{C}\subset\sigma(\mathcal{C}')$ and $\mathcal{C}'\subset\sigma(\mathcal{C})$, then $\sigma(\mathcal{C})=\sigma(\mathcal{C}')$.
\end{itemize}
\end{theorem}

\begin{proof}
\textit{(a):} Theorem~\ref{213} guarantees that $\mathcal{C}\subset\mathcal{C}'\subset\sigma(\mathcal{C}')$. Thus, $\sigma(\mathcal{C}')$ is a $\sigma$-algebra containing $\mathcal{C}$, and therefore $\sigma(\mathcal{C})\subset\sigma(\mathcal{C}')$ by Theorem~\ref{213}.

\textit{(b):} Since $\mathcal{C}\subset\mathcal{C}$ and $\mathcal{C}$ is a $\sigma$-algebra, it follows from Theorem~\ref{213} that $\sigma(\mathcal{C})\subset\mathcal{C}$. On the other hand, the same theorem gives $\mathcal{C}\subset\sigma(\mathcal{C})$.

\textit{(c):} This is a direct consequence of the previous part.

\textit{(d):} Since $\mathcal{C}\subset\sigma(\mathcal{C}')$, parts \textit{(a)} and \textit{(c)} imply that $\sigma(\mathcal{C})\subset\sigma(\sigma(\mathcal{C}'))=\sigma(\mathcal{C}')$. Similarly, by parts \textit{(a)} and \textit{(c)}, if $\mathcal{C}'\subset\sigma(\mathcal{C})$, then $\sigma(\mathcal{C}')\subset\sigma(\mathcal{C})$.

Therefore, $\sigma(\mathcal{C})=\sigma(\mathcal{C}')$.
\end{proof}

\begin{remark}\label{216}
The previous definition, together with the analogous properties corresponding to $\mathcal{R}(\mathcal{C})$ (the ring generated by $\mathcal{C}$), $\mathcal{A}(\mathcal{C})$ (the algebra generated by $\mathcal{C}$), and $\mathcal{SR}(\mathcal{C})$ (the $\sigma$-ring generated by $\mathcal{C}$), remain valid when $\mathcal{C}$ is an arbitrary class of subsets of $X$. In fact, $\mathcal{R}(\mathcal{C})\subset\mathcal{SR}(\mathcal{C})\subset\sigma(\mathcal{C})$ and $\mathcal{R}(\mathcal{C})\subset\mathcal{A}(\mathcal{C})\subset\sigma(\mathcal{C})$. \index{generated!ring} \index{sigma@$\sigma$!generated ring!by a class $\mathcal{SR}(\mathcal{C})$} \index{generated!algebra} \index{algebra!generated!by a class $\mathcal{A}(\mathcal{C})$} \index{ring!generated!by a class $\mathcal{R}(\mathcal{C})$}
\end{remark}

From the previous remark, we may conclude the following: if $\mathcal{C}$ is a finite class of subsets of $X$, then $\mathcal{A}(\mathcal{C})$ is a finite algebra of subsets of $X$ [Exercise~\ref{E29}] and, therefore, $\mathcal{A}(\mathcal{C})$ is a $\sigma$-algebra of subsets of $X$ such that $\mathcal{C}\subset\mathcal{A}(\mathcal{C})$ [Exercise~\ref{E28}]. Consequently,
$\sigma(\mathcal{C})=\mathcal{A}(\mathcal{C})$.

Let us now consider some examples.

\begin{example}\label{217}
If $\mathcal{C}=\varnothing$ or $\mathcal{C}=\{\varnothing\}$, then
$$
\mathcal{R}(\mathcal{C})=\mathcal{SR}(\mathcal{C})=\{\varnothing\}
\quad\text{and}\quad
\mathcal{A}(\mathcal{C})=\sigma(\mathcal{C})=\{\varnothing,X\}.
$$
\end{example}

\begin{proof}
It suffices to prove the result for $\mathcal{R}(\mathcal{C})$ and $\mathcal{A}(\mathcal{C})$, respectively [Exercise~\ref{E29}].

Clearly, $\{\varnothing\}$ is a ring of subsets of $X$ and satisfies $\mathcal{C}\subset\{\varnothing\}$. Therefore, by Theorem~\ref{213}, it follows that $\mathcal{R}(\mathcal{C})\subset\{\varnothing\}$.

On the other hand, since $\varnothing\in\mathcal{R}$ for every ring of subsets of $X$, we have $\{\varnothing\}\subset\mathcal{R}(\mathcal{C})$,
and hence equality holds.

By Example~\ref{27}, the class $\{\varnothing,X\}$ is an algebra of subsets of $X$ such that $\mathcal{C}\subset\{\varnothing,X\}$. Consequently,
$\mathcal{A}(\mathcal{C})\subset\{\varnothing,X\}$.

Conversely, $\varnothing,X\in\mathcal{A}$ for every algebra of subsets of $X$ by properties (A1) and (A2). Therefore,
$\{\varnothing,X\}\subset\mathcal{A}(\mathcal{C})$.
\end{proof}

\begin{example}\label{218}
If $\mathcal{C}\subset\mathcal{P}(X)$ is a ring, then
$$
\mathcal{A}(\mathcal{C})=\{A\subset X:A\in\mathcal{C}\quad\text{or}\quad X\smallsetminus A\in\mathcal{C}\}.
$$
\end{example}

\begin{proof}
For simplicity, let us denote
$$
\mathfrak{A}
:= \{A\subset X : A\in\mathcal{C} \,\,\text{or}\,\,
X\smallsetminus A\in\mathcal{C}\},
$$
and prove that $\mathfrak{A}$ is an algebra of subsets of $X$ containing the class $\mathcal{C}$.

Clearly, $\mathfrak{A}$ is nonempty, since $X\in\mathfrak{A}$, because $\mathcal{C}$ is a ring of subsets of $X$ and therefore $\varnothing\in\mathcal{C}$.

We now proceed to verify the following claims.

{\scshape Claim 1:} $\mathfrak{A}$ is closed under complements.

Let $A\in\mathfrak{A}$. If $A\in\mathcal{C}$, then $X\smallsetminus A\in\mathfrak{A}$, since
$$
X\smallsetminus(X\smallsetminus A)=A\in\mathcal{C}.
$$

On the other hand, if $X\smallsetminus A\in\mathcal{C}$, then by definition $X\smallsetminus A\in\mathfrak{A}$.

Consequently, $X\smallsetminus A\in\mathfrak{A}$ for every $A\in\mathfrak{A}$.

{\scshape Claim 2:} $\mathfrak{A}$ is closed under finite unions.

Let $A,B\in\mathfrak{A}$ be arbitrary. Consider the following cases.

{\scshape Case 1:} If $A,B\in\mathcal{C}$, then $A\cup B\in\mathcal{C}$ by property (R2), and therefore $A\cup B\in\mathfrak{A}$.

{\scshape Case 2:} If $X\smallsetminus A,X\smallsetminus B\in\mathcal{C}$, then by Proposition~\ref{22},
$$
(X\smallsetminus A)\cap(X\smallsetminus B)
=
X\smallsetminus(A\cup B)\in\mathcal{C}.
$$

Consequently, $A\cup B\in\mathfrak{A}$.

{\scshape Case 3:} If $X\smallsetminus A\in\mathcal{C}$ and $B\in\mathcal{C}$, then by property (R1),
$$
(X\smallsetminus A)\smallsetminus B
=
(X\smallsetminus A)\cap(X\smallsetminus B)
=
X\smallsetminus(A\cup B)\in\mathcal{C},
$$
and consequently $A\cup B\in\mathfrak{A}$.

{\scshape Case 4:} The case in which $A\in\mathcal{C}$ and $X\smallsetminus B\in\mathcal{C}$ is proved analogously to the previous case.

From the previous cases, we conclude that $\mathfrak{A}$ is an algebra of subsets of $X$ containing the class $\mathcal{C}$ [Exercise~\ref{E215}]. Consequently,
$$
\mathcal{A}(\mathcal{C})\subset\mathfrak{A}.
$$

The reverse inclusion is immediate. Indeed, if $A\in\mathcal{C}$, then $A\in\mathcal{A}(\mathcal{C})$. On the other hand, if $X\smallsetminus A\in\mathcal{C}$, then $X\smallsetminus A\in\mathcal{A}(\mathcal{C})$ and, since $\mathcal{A}(\mathcal{C})$ is an algebra, it follows that $A\in\mathcal{A}(\mathcal{C})$. Therefore,
$\mathfrak{A}\subset\mathcal{A}(\mathcal{C})$,
which completes the proof.
\end{proof}

Given a class $\mathcal{C}$ of subsets of $X$ and a subset $A$ of $X$, we define
$$
A\cap\mathcal{C}:=\{A\cap B:B\in\mathcal{C}\},
$$
which is clearly a class of subsets of $A$ and, therefore, $\sigma(A\cap\mathcal{C})$ is a $\sigma$-algebra of subsets of $A$.

We now present some important properties of this construction.

\begin{proposition}\label{219}
Let $\mathcal{C}\subset\mathcal{P}(X)$ and let $A$ be a fixed nonempty subset of $X$. Then the collection $A\cap\sigma(\mathcal{C})$ is a $\sigma$-algebra of subsets of $A$.
\end{proposition}

\begin{proof}
Clearly, $A\in A\cap\sigma(\mathcal{C})$ since $A=A\cap X$ and $X\in\sigma(\mathcal{C})$. Let now $B_1,B_2\in\sigma(\mathcal{C})$ be arbitrary. From the identity between subsets of $X$
$$
(A\cap B_1)\smallsetminus(A\cap B_2)=A\cap(B_1\smallsetminus B_2)
$$
it follows directly that $(A\cap B_1)\smallsetminus(A\cap B_2)\in A\cap\sigma(\mathcal{C})$ since $B_1\smallsetminus B_2\in\sigma(\mathcal{C})$. This proves that $A\cap\sigma(\mathcal{C})$ satisfies properties (S1) and (S2).

Finally, let $(A_k)$ be a sequence of elements of $A\cap\sigma(\mathcal{C})$. Then there exists a sequence $(B_k)$ of elements of $\sigma(\mathcal{C})$ such that $A_k=A\cap B_k$ for every $k\in\mathbb{N}$. Consequently,
$$
\bigcup_{k=1}^{\infty}A_k=\bigcup_{k=1}^{\infty}(A\cap B_k)=A\cap\left(\bigcup_{k=1}^{\infty}B_k\right)\in A\cap\sigma(\mathcal{C})
$$
since $\bigcup_{k=1}^{\infty}B_k\in\sigma(\mathcal{C})$. Thus, $A\cap\sigma(\mathcal{C})$ satisfies property (S3).

Therefore, $A\cap\sigma(\mathcal{C})$ is a $\sigma$-algebra of subsets of $A$.
\end{proof}

\begin{proposition}\label{220}
Let $\mathcal{C}\subset\mathcal{P}(X)$ and let $A$ be a fixed nonempty subset of $X$. The class defined by
$$
\mathcal{F}:=\{B\in\sigma(\mathcal{C}):A\cap B\in\sigma(A\cap\mathcal{C})\}
$$
is a $\sigma$-algebra of subsets of $X$ such that $\mathcal{F}=\sigma(\mathcal{C})$.
\end{proposition}

\begin{proof}
It is immediate that $X\in\mathcal{F}$ since $X\in\sigma(\mathcal{C})$ and $A=A\cap X\in\sigma(A\cap\mathcal{C})$.

Let $B_1,B_2\in\mathcal{F}$ be arbitrary. Then $B_1,B_2\in\sigma(\mathcal{C})$ and are such that $A\cap B_1,A\cap B_2\in\sigma(A\cap\mathcal{C})$. Since $\sigma(A\cap\mathcal{C})$ is a $\sigma$-algebra of subsets of $A$, it follows that
$$
(A\cap B_1)\smallsetminus(A\cap B_2)=A\cap(B_1\smallsetminus B_2)\in\sigma(A\cap\mathcal{C})
$$
and therefore $B_1\smallsetminus B_2\in\mathcal{F}$ since $B_1\smallsetminus B_2\in\sigma(\mathcal{C})$.

Let $(B_k)$ be a sequence of elements of $\mathcal{F}$. Then $B_k\in\sigma(\mathcal{C})$ and $A\cap B_k\in\sigma(A\cap\mathcal{C})$ for every $k\in\mathbb{N}$, so that $(B_k)$ and $(A\cap B_k)$ are sequences of elements of $\sigma(\mathcal{C})$ and $\sigma(A\cap\mathcal{C})$, respectively. Thus,
$$
\bigcup_{k=1}^{\infty}B_k\in\sigma(\mathcal{C})
\quad\text{and}\quad
\bigcup_{k=1}^{\infty}(A\cap B_k)=A\cap\left(\bigcup_{k=1}^{\infty}B_k\right)\in\sigma(A\cap\mathcal{C})
$$
which proves that $\bigcup_{k=1}^{\infty}B_k\in\mathcal{F}$.

Consequently, $\mathcal{F}$ is a $\sigma$-algebra of subsets of $X$ such that $\mathcal{F}\subset\sigma(\mathcal{C})$.

On the other hand, if $B\in\mathcal{C}$, then $B\in\sigma(\mathcal{C})$ and
$A\cap B\in A\cap\mathcal{C}\subset\sigma(A\cap\mathcal{C})$. That is, $\mathcal{C}\subset\mathcal{F}$ and therefore $\sigma(\mathcal{C})\subset\mathcal{F}$. This proves that $\mathcal{F}=\sigma(\mathcal{C})$.
\end{proof}

\begin{theorem}\label{221}
Let $\mathcal{C}\subset\mathcal{P}(X)$ and let $A$ be a fixed nonempty subset of $X$. Then,
$$
\sigma(A\cap\mathcal{C})=A\cap\sigma(\mathcal{C}).
$$
\end{theorem}

\begin{proof}
By Proposition~\ref{219}, $A\cap\sigma(\mathcal{C})$ is a $\sigma$-algebra of subsets of $A$ such that $A\cap\mathcal{C}\subset A\cap\sigma(\mathcal{C})$ since $\mathcal{C}\subset\sigma(\mathcal{C})$. Then, by part \textit{(b)} of Theorem~\ref{213}, we have $\sigma(A\cap\mathcal{C})\subset A\cap\sigma(\mathcal{C})$.

Conversely, by Proposition~\ref{220}, the class $\mathcal{F}:=\{B\in\sigma(\mathcal{C}):A\cap B\in\sigma(A\cap\mathcal{C})\}$ is a $\sigma$-algebra of subsets of $X$ that coincides with $\sigma(\mathcal{C})$. Consequently,
$A\cap\sigma(\mathcal{C})\subset\sigma(A\cap\mathcal{C})$.
\end{proof}

\begin{remark}\label{222}
The previous theorem also remains valid if one considers the ring, algebra, and $\sigma$-ring generated by the class $A\cap\mathcal{C}$.
\end{remark}

\begin{lemma}\label{223}
Let $\mathcal{C}\subset\mathcal{P}(X)$ be nonempty.
\begin{itemize}
    \item[(a)] The set $\mathcal{R}$ defined by
$$
\mathcal{R}=\left\{ A \subset X\,:\, \text{there exist } N \in \mathbb{N}\, \text{ and }\, C_{1},\ldots,C_{N} \in \mathcal{C}\, \text{ such that }\, A \subset \bigcup_{j=1}^{N} C_{j} \right\}
$$
is a ring of subsets of $X$.
\item[(b)] The set $\mathcal{SR}$ defined by
$$
\mathcal{SR}=\left\{ A \subset X\,:\, \text{there exists a sequence } \,(C_{k})\, \text{ in } \,\mathcal{C}\, \text{ such that }\, A \subset \bigcup_{k=1}^{\infty} C_{k} \right\}
$$
is a $\sigma$-ring of subsets of $X$.
\end{itemize}
\end{lemma}

\begin{proof}
\textit{(a):} It is clear that $\varnothing \in \mathcal{R}$ since $\varnothing \subset C$ for every $C \in \mathcal{C}$. Thus, $\mathcal{R}$ is a nonempty class of subsets of $X$.

Let $A,B \in \mathcal{R}$. There exist $C_{1},\ldots,C_{N} \in \mathcal{C}$ and $D_{1},\ldots,D_{M} \in \mathcal{C}$ such that
$$
A \subset \bigcup_{i=1}^{N} C_{i} \quad \mbox{and}\quad B \subset \bigcup_{j=1}^{M} D_{j}.
$$

Since
$$
A \smallsetminus B \subset A \subset \bigcup_{i=1}^{N} C_{i} \quad \mbox{and}\quad A \cup B \subset \left(\bigcup_{i=1}^{N} C_{i}\right) \cup \left( \bigcup_{j=1}^{M} D_{j} \right),
$$
we conclude that $A\smallsetminus B$ and $A \cup B$ are elements of $\mathcal{R}$ and, therefore, that $\mathcal{R}$ is a ring of subsets of $X$.

\textit{(b):} Again, it is clear that $\mathcal{SR}$ is nonempty since $\varnothing \in \mathcal{SR}$. Let $A,B \in \mathcal{SR}$. There exist sequences $(C_k)$ and $(D_k)$ of elements of $\mathcal{C}$ such that
$$
A \subset \bigcup_{k=1}^{\infty} C_{k} \quad \mbox{and}\quad B \subset \bigcup_{k=1}^{\infty} D_{k}.
$$

Since $A \smallsetminus B \subset A$, we have
$$
A \smallsetminus B \subset \bigcup_{k=1}^{\infty} C_{k}
$$
which implies that $A \smallsetminus B \in \mathcal{SR}$.

Let now $(A_{j})$ be a sequence of elements of $\mathcal{SR}$. For each $j \in \mathbb{N}$, there exists a sequence $(G_{k,j})$ of elements of $\mathcal{C}$ such that
$$
A_{j} \subset \bigcup_{k=1}^{\infty} G_{k,j}.
$$

Consequently,
$$
\bigcup_{j=1}^{\infty} A_{j} \subset \bigcup_{j=1}^{\infty} \left( \bigcup_{k=1}^{\infty} G_{k,j}  \right) = \bigcup_{(j,k) \in \mathbb{N}^2} G_{k,j}.
$$

Let $\phi:\mathbb{N}\to \mathbb{N}^{2}$ be a bijection and define $C_{\ell}:=G_{\phi(\ell)}$ for every $\ell \in \mathbb{N}$. Then $(C_{\ell})$ is a sequence of elements of $\mathcal{C}$ and satisfies
$$
\bigcup_{(j,k) \in \mathbb{N}^2} G_{k,j}=\bigcup_{\ell=1}^{\infty} C_{\ell}.
$$

Therefore,
$$
\bigcup_{j=1}^{\infty} A_{j} \subset \bigcup_{\ell=1}^{\infty} C_{\ell}
$$
and we conclude that $\bigcup_{j=1}^{\infty} A_{j} \in \mathcal{SR}$.
\end{proof}

\begin{theorem}\label{224}
Let $\mathcal{C}\subset\mathcal{P}(X)$ be nonempty.
\begin{itemize}
\item[(1)] Every $A\in\mathcal{R}(\mathcal{C})$ is contained in the union of finitely many elements of $\mathcal{C}$.
\item[(2)] Every $A\in\mathcal{SR}(\mathcal{C})$ is contained in the union of some sequence of elements of $\mathcal{C}$.
\end{itemize}
\end{theorem}

\begin{proof}
\textit{(1)} and \textit{(2):} Consider the classes $\mathcal{R}$ and $\mathcal{SR}$ defined in Lemma~\ref{223}. The statement follows immediately from the fact that $\mathcal{C}\subset\mathcal{R}$ and $\mathcal{C}\subset\mathcal{SR}$, since every family containing $\mathcal{C}$ must, in particular, contain the smallest class with the property under consideration (see Remark~\ref{216}).
\end{proof}

Let $X$ and $Y$ be two nonempty sets and let $f:X\to Y$ be a function. Given a class $\mathcal{B}$ of subsets of $Y$, we denote by
$$
f^{-1}(\mathcal{B}):=\left\{ f^{-1}(B)\subset X:B\in\mathcal{B} \right\}
$$
the class of inverse images of elements of $\mathcal{B}$.

\begin{proposition}\label{225}
If $f:X \to Y$ is a function and $\mathcal{B}$ is a class of subsets of $Y$, then
$$
\sigma\left( f^{-1}(\mathcal{B}) \right) = f^{-1}(\sigma(\mathcal{B})).
$$
\end{proposition}

\begin{proof}
Since $\sigma(\mathcal{B})$ is a $\sigma$-algebra of subsets of $Y$ such that $\mathcal{B}\subset\sigma(\mathcal{B})$, it follows that $f^{-1}(\sigma(\mathcal{B}))$ is a $\sigma$-algebra of subsets of $X$ such that $f^{-1}(\mathcal{B})\subset f^{-1}(\sigma(\mathcal{B}))$ [Exercise~\ref{E210}]. Theorem~\ref{213} guarantees that
$$
\sigma (f^{-1}(\mathcal{B})) \subset f^{-1}(\sigma(\mathcal{B})).
$$

Conversely: Define the class
$$
\mathfrak{F}:=\{B \subset Y\,:\,f^{-1}(B) \in \sigma(f^{-1}(\mathcal{B}))\},
$$
which is a $\sigma$-algebra of subsets of $Y$ [Exercise~\ref{E211}]. Now, if $B \in \mathcal{B}$, then
$f^{-1}(B) \in f^{-1}(\mathcal{B}) \subset \sigma(f^{-1}(\mathcal{B}))$,
so that $\mathcal{B}\subset\mathfrak{F}$. Therefore,
$\sigma(\mathcal{B})\subset\mathfrak{F}$,
and from the definition of $\mathfrak{F}$ it follows that
$f^{-1}(\sigma(\mathcal{B}))\subset\sigma (f^{-1}(\mathcal{B}))$.
\end{proof}

\section{\texorpdfstring{$\sigma$}{}-algebra of Borel}

In a basic course on real analysis, the notion of an open set constitutes the foundation upon which essential concepts such as convergence, continuity, and compactness are formulated and generalized. The abstract structure that allows these notions to be formalized and systematically studied is that of a topological space. A \textbf{topological space} $(X,\tau)$ is a nonempty set $X$ endowed with a \textbf{topology} $\tau$, that is, a nonempty class $\tau$ of subsets of $X$ satisfying the following properties:
\begin{itemize} \index{space!topological} \index{topology}
    \item[(a)] $X,\varnothing \in \tau$.
    \item[(b)] If $A,B \in \tau$, then $A \cap B \in \tau$.
    \item[(c)] If $\{A_j \,:\, j \in \mathcal{J}\}$ is an arbitrary family of elements of $\tau$, then $\bigcup_{j \in \mathcal{J}} A_j \in \tau$.
\end{itemize}

The elements of $\tau$ are called the \textbf{open subsets} of $X$ relative to $\tau$.

Let $(X,\tau)$ be a topological space. The class $\tau$ not only describes the topological structure of $X$, but also allows one to construct canonically an associated measurable structure. Indeed, the Borel $\sigma$-algebra of $X$ is defined by $\mathcal{B}(X):=\sigma(\tau)$, that is, the smallest $\sigma$-algebra containing all open subsets of $X$ relative to $\tau$.

The complement of an open subset of $X$ relative to $\tau$ is called a \textbf{closed subset} of $X$ relative to $\tau$. It can be shown that the class of closed subsets of $X$ relative to $\tau$ also generates the Borel $\sigma$-algebra of $X$ [Exercise \ref{E217}].

In this way, the topology naturally induces a measurable structure compatible with it. The Borel $\sigma$-algebra constitutes the canonical domain on which measures respecting the topological structure of the space are defined, thereby establishing an essential connection between topology and measure theory. In the particular case of the real numbers $\mathbb{R}$, this construction gives rise to the Borel $\sigma$-algebra, introduced in 1898 by the French mathematician Émile Borel\footnote{Félix Édouard Justin Émile Borel (1871--1956) was a French mathematician and politician. Together with René-Louis Baire and Henri Lebesgue, he was one of the pioneers of measure theory and its applications to probability theory. One of his probability books introduced the entertaining thought experiment that entered popular culture under the name of the infinite monkey theorem. In addition, he published research on game theory.}. \index{sigma@$\sigma$!algebra!of Borel on $X$}

\begin{figure}[ht!]
\centering
\includegraphics[scale=0.275]{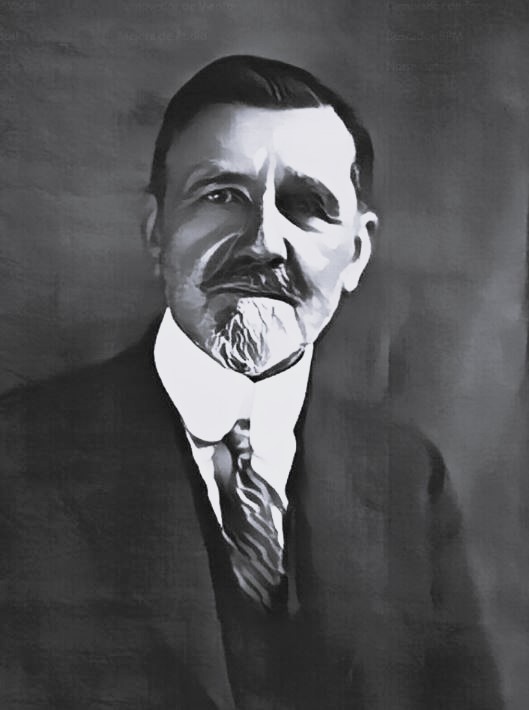} 
\begin{center}
Émile Borel (1871-1956)
\end{center}
\end{figure} 

\begin{definition}\label{226}
Let $\mathbb{R}$ and $\tau_{\mathbb{R}} \subset \mathcal{P}(\mathbb{R})$ be the usual topology of $\mathbb{R}$. We define the \textbf{Borel $\sigma$-algebra on $\mathbb{R}$}\index{sigma@$\sigma$!algebra!of Borel on $\mathbb{R}$}, which we denote by $\mathcal{B}(\mathbb{R}),$ as the $\sigma$-algebra generated by $\tau_{\mathbb{R}}$, that is,
$$
\mathcal{B}(\mathbb{R}):=\sigma(\tau_{\mathbb{R}}).
$$
\end{definition}

Note the following: every element of $\tau_{\mathbb{R}}$, that is, every open set in $\mathbb{R}$, can be written as the union of an at most countable family of open intervals with endpoints in $\mathbb{R}$. Thus, if $\mathcal{I}$ is the class of open intervals with endpoints in $\mathbb{R}$, namely,
$$
\mathcal{I}=\left\{ (a,b)\,:\, a < b\,\,\text{with}\,\,a,b \in \mathbb{R} \right\},
$$
then $\tau_{\mathbb{R}}  \subset \sigma(\mathcal{I})$. Hence, since $\mathcal{I} \subset \tau_{\mathbb{R}} \subset  \mathcal{B}(\mathbb{R})$, part \textit{(d)} of Theorem \ref{215} guarantees that
$$
\mathcal{B}(\mathbb{R})=\sigma(\mathcal{I})=\sigma\left(\left\{ (a,b)\,:\, a < b\,\,\text{with}\,\,a,b \in \mathbb{R} \right\} \right).
$$

The elements of $\mathcal{B}(\mathbb{R})$ are called \textit{Borel sets}, \textit{Borelians}, or \textit{Borel measurable sets}. \index{set!Borel measurable}

What kinds of elements belong to the Borel $\sigma$-algebra? Since $\mathcal{I} \subset \sigma(\mathcal{I})=\mathcal{B}(\mathbb{R})$, every open interval with endpoints in $\mathbb{R}$ is a Borel measurable set.

Closed intervals with endpoints in $\mathbb{R}$, $[a,b]$ $(a < b)$, can be written as the countable intersection of open intervals in the following way:
$$
[a,b]=\bigcap_{k=1}^{\infty} \left(a-\frac{1}{k}, b + \frac{1}{k} \right).
$$

Since $\mathcal{B}(\mathbb{R})$ is a $\sigma$-algebra, the countable intersection of open intervals is also an element of $\mathcal{B}(\mathbb{R})$. Thus, since $\left(a-\frac{1}{k} , b + \frac{1}{k} \right) \in \mathcal{B}(\mathbb{R})$ for every $k \in \mathbb{N}$, we conclude that every closed interval is a Borel measurable set.

Likewise, for $a,b \in \mathbb{R}$
$$
(a,+\infty)=\bigcup_{k=1}^{\infty}(a,a+k) \in \mathcal{B}(\mathbb{R})
\quad\mbox{and}\quad (-\infty,b)=\bigcup_{k=1}^{\infty}(b-k,b) \in \mathcal{B}(\mathbb{R}).
$$

Consequently,
$$
[a,+\infty) = \bigcap_{k=1}^{\infty} \left(a-\frac{1}{k}, + \infty\right) \in \mathcal{B}(\mathbb{R})
\quad\mbox{and}\quad(-\infty,b]=\bigcap_{k=1}^{\infty}\left(-\infty, b + \frac{1}{k}\right) \in \mathcal{B}(\mathbb{R}).
$$

Analogously, we obtain that
$$
[a,b)=\bigcap_{k=1}^{\infty} \left(a-\frac{1}{k},b\right) \in\mathcal{B}(\mathbb{R})
\quad\mbox{and}\quad (a,b]=\bigcap_{k=1}^{\infty} \left(a,b+\frac{1}{k}\right)  \in \mathcal{B}(\mathbb{R}).
$$

Sets consisting of a single number are also Borel measurable sets since
$$
\{a \}=\bigcap_{k=1}^{\infty} \left(a-\frac{1}{k}, a + \frac{1}{k}\right)\in \mathcal{B}(\mathbb{R}).
$$

Complements, countable intersections, and countable unions of these sets are all Borel measurable.

\begin{proposition}\label{227}
For any real numbers $a,b$ such that $a<b$, the intervals $[a,b]$, $(a,+\infty)$, $(-\infty,b)$, $[a,b)$, $(a,b]$, and $\{a\}$ are all elements of $\mathcal{B}(\mathbb{R})$.
\end{proposition}

\begin{proposition}\label{228}
The sets $\mathbb{N}, \mathbb{Z}, \mathbb{Q}$, and $\mathbb{I}$ are Borel measurable.
\end{proposition}

\begin{proof}
The proof follows directly from the following equalities between Borel sets:
$$
\mathbb{N}=\bigcup_{k=1}^{\infty} \{k\}\in \mathcal{B}(\mathbb{R}),
$$
$$
\mathbb{Z}=\bigcup_{k=1}^{\infty} \{k-1 \} \cup \bigcup_{k=1}^{\infty} \{-k\}\in \mathcal{B}(\mathbb{R}),
$$
$$
\mathbb{Q}=\bigcup_{\substack{n,m\in\mathbb{Z} \\m\neq 0} } \left\{ \frac{n}{m} \right\}\in \mathcal{B}(\mathbb{R}),
$$
$$
\mathbb{I}=\mathbb{R}\smallsetminus \mathbb{Q} \in \mathcal{B}(\mathbb{R}).
$$

This completes the proof.
\end{proof}

\begin{proposition}\label{229}
The following classes of subsets of $\mathbb{R}$ generate the Borel $\sigma$-algebra:
\begin{itemize}
\item[(1)] $\mathcal{C}=\{ [a,b]\,:\, a< b \,\,\,(a,b \in \mathbb{R} ) \}$.
\item[(2)] $\mathcal{C}=\{ (a,b]\,:\, a< b \,\,\,(a,b \in \mathbb{R} ) \}$.
\item[(3)] $\mathcal{C}=\{ [a,b)\,:\, a< b \,\,\,(a,b \in \mathbb{R} ) \}$.
\item[(4)] $\mathcal{C}=\{ (a,+\infty)\,:\, a \in \mathbb{R}  \}$.
\item[(5)] $\mathcal{C}=\{ (-\infty,b]\,:\, b \in \mathbb{R}  \}$.
\end{itemize}
\end{proposition}

\begin{proof}
\textit{(1)} By part \textit{(d)} of Theorem \ref{215}, it suffices to show that $\mathcal{C} \subset \sigma(\mathcal{I})=\mathcal{B}(\mathbb{R})$ and $\mathcal{I}\subset \sigma(\mathcal{C})$.

From the above, it is clear that every closed interval $[a,b] \in \mathcal{C}$ belongs to $\mathcal{B}(\mathbb{R})$ since
$$
[a,b]=\bigcap_{k=1}^{\infty} \left(a-\frac{1}{k}, b + \frac{1}{k} \right).
$$

Therefore $\mathcal{C} \subset \mathcal{B}(\mathbb{R})$.

Conversely, every open interval $(a,b) \in \mathcal{I}$ can be written as
$$
(a,b)=\bigcup_{k=1}^{\infty} \left[a+\frac{1}{k}, b - \frac{1}{k} \right]
$$
which establishes that $\mathcal{I} \subset \sigma(\mathcal{C})$.

The remaining parts are proved in a completely analogous way [Exercise \ref{E227}].
\end{proof}

\begin{definition}\label{230} \index{set!type G delta@type $\mathcal{G}_{\delta}$} \index{set!type F sigma@type $\mathcal{F}_{\sigma}$}
A set is of \textbf{type $\mathcal{G}_{\delta}$} if it is the intersection of a countable family of open sets, and a set is of \textbf{type $\mathcal{F}_{\sigma}$} if it is the union of a countable family of closed sets.
\end{definition}

The sets belonging to the classes of Definition \ref{230} are all Borel measurable. We now define the class of sets of \textbf{type $\mathcal{G}_{\sigma\,\delta}$} as those that can be written as the countable union of sets of type $\mathcal{G}_{\delta}$, and the class of sets of \textbf{type $\mathcal{F}_{\delta\,\sigma}$} as those that can be written as the countable intersection of sets of type $\mathcal{F}_{\sigma}$. Continuing this process, one can show [Exercise \ref{E228}] that the following two chains of classes of Borel measurable sets are obtained:
$$
\begin{array}{c}
\mathcal{G}_{\delta} \subset \mathcal{G}_{\sigma\delta} \subset \mathcal{G}_{\delta\sigma\delta} \subset \mathcal{G}_{\sigma\delta\sigma\delta} \subset \cdots\\
\mathcal{F}_{\sigma} \subset \mathcal{F}_{\delta\sigma} \subset \mathcal{F}_{\sigma\delta\sigma} \subset \mathcal{F}_{\delta\sigma\delta\sigma} \subset \cdots \\
\end{array}
$$

It is natural to ask whether the collection $\mathcal{B}(\mathbb{R})$ contains all subsets of $\mathbb{R}$. The answer is no; that is, one can prove that there exists a subset $\mathsf{V}$ (Vitali) of $\mathbb{R}$ that is not Borel measurable. The construction of this set will be presented in a later chapter.

Given any subset $A \subset \mathbb{R}$, by endowing it with the subspace topology
$$
\tau_A := \{ A \cap \mathcal{I} : \mathcal{I} \in \tau_{\mathbb R} \},
$$
the pair $(A,\tau_A)$ becomes a topological space. In particular, this space comes naturally equipped with a Borel $\sigma$-algebra $\mathcal{B}(A)$. The purpose of what follows is to show that this $\sigma$-algebra is not a new or different structure, but rather coincides exactly with the restriction of the Borel $\sigma$-algebra of the real numbers, $\mathcal{B}(\mathbb{R})$, to the subset $A$.

\begin{proposition}\label{231}
Let $A \subset \mathbb{R}$. The Borel $\sigma$-algebra of $A$, considered as a topological subspace of $\mathbb{R}$, coincides with the $\sigma$-algebra $A\cap \mathcal{B}(\mathbb{R})$.
\end{proposition}

\begin{proof}
It is immediate to note that $\tau_{A}=A\cap \tau_{\mathbb{R}}$ and to apply Theorem \ref{221}.
\end{proof}

\begin{example}\label{232}
Consider $A=\mathbb{N}$. For each $j \in \mathbb{N}$, the open interval
$$
\mathcal{I}_{j}:=\left(j-\frac{1}{4},j+\frac{1}{4} \right) \subset \mathbb{R}
$$
satisfies $\mathcal{I}_{j}\cap \mathbb{N}=\{j\}$ and, therefore, $\{j\} \in \tau_{\mathbb{N}}$.

Since every subset $B$ of $\mathbb{N}$ can be written as the at most countable union
$$
B=\bigcup_{j \in B}\{j\},
$$
we conclude that $B \in \tau_{\mathbb{N}}$. Consequently,
$$
\tau_{\mathbb{N}} = \mathcal{P}(\mathbb{N})
\quad\text{and}\quad
\sigma(\tau_{\mathbb{N}})=\mathcal{B}(\mathbb{N})=\mathcal{P}(\mathbb{N}).
$$
\end{example}

Proposition \ref{231} remains valid in any topological space $(X,\tau)$ by endowing each subset $A$ of $X$ with the subspace topology defined by $\tau_{A}:=A \cap \tau$.

\subsection{The space \texorpdfstring{$\mathbb{R}^{N}$}{}}

Given $\bar{x} \in \mathbb{R}^{N}$ for $N>1$, define
$$
\begin{aligned}
    \Vert \bar{x}\Vert_{p}&:=\left(\sum_{j=1}^{N} |x_{j}|^{p} \right)^{1/p} \quad \text{if }\,p\in [1,\infty),\\
    \Vert \bar{x}\Vert_{\infty}&:=\max_{1\leq j \leq N}|x_{j}|.
\end{aligned}
$$

In a basic course on real analysis it is shown that $\Vert \cdot \Vert_{p}$ is a norm on $\mathbb{R}^{N}$ for every $p \in [1,\infty]$\footnote{See, for example, \cite[Chapter 2]{Clapp}.}.

The topology on $\mathbb{R}^{N}$ induced by the norm $\Vert \cdot \Vert_{p}$ for every $p \in [1,\infty]$ is described as follows: Given $\bar{x} \in \mathbb{R}^{N}$ and $\varepsilon>0$, we define the \textbf{open ball} centered at $\bar{x}$ and radius $\varepsilon$ with respect to the norm $\Vert \cdot \Vert_{p}$ for $p \in [1,\infty]$ as the set
$$
B_{p}(\bar{x},\varepsilon):=\left\{ \bar{y}\in \mathbb{R}^{N}\,:\,\Vert \bar{x}-\bar{y} \Vert_{p} <\varepsilon \right\}.
$$

Therefore, a subset $\mathcal{O}$ of $\mathbb{R}^{N}$ is open with respect to the norm $\Vert \cdot \Vert_{p}$ for $p \in [1,\infty]$ if, for every $\bar{x} \in \mathcal{O}$, there exists $\varepsilon>0$ such that $B_{p}(x,\varepsilon) \subset \mathcal{O}$.

\begin{figure}[h!]
\centering

\begin{minipage}{0.32\textwidth}
\centering
\begin{tikzpicture}[scale=1.2]
    \draw[dotted, fill=blue!10] (1,0) -- (0,1) -- (-1,0) -- (0,-1) -- cycle;

    \draw[->,gray] (-1.6,0) -- (1.6,0);
    \draw[->,gray] (0,-1.6) -- (0,1.6);

    \node at (0,-1.9) {$B_{1}(\bar{0},1)$};
\end{tikzpicture}
\end{minipage}
\hfill
\begin{minipage}{0.32\textwidth}
\centering
\begin{tikzpicture}[scale=1.2]
    \draw[dotted, fill=red!10] (0,0) circle (1);

    \draw[->,gray] (-1.6,0) -- (1.6,0);
    \draw[->,gray] (0,-1.6) -- (0,1.6);

    \node at (0,-1.9) {$B_{2}(\bar{0},1)$};
\end{tikzpicture}
\end{minipage}
\hfill
\begin{minipage}{0.32\textwidth}
\centering
\begin{tikzpicture}[scale=1.2]

    \draw[dotted, fill=lime!15] (-1,-1) rectangle (1,1);

    \draw[->,gray] (-1.6,0) -- (1.6,0);
    \draw[->,gray] (0,-1.6) -- (0,1.6);
    \node at (0,-1.9) {$B_{\infty}(\bar{0},1)$};
    
\end{tikzpicture}
\end{minipage}
\end{figure}

We denote by $\tau_{p}$ the class of open subsets of $\mathbb{R}^{N}$ with respect to the norm $\Vert \cdot \Vert_{p}$ for $p \in [1,\infty]$. Therefore, $(\mathbb{R}^{N},\tau_{p})$ is a topological space for every $p \in [1,\infty]$.

Since $(\mathbb{R}^{N},\tau_{p})$ is a topological space for every $p \in [1,\infty]$, each of these topologies generates its corresponding Borel $\sigma$-algebra, which we denote by $\mathcal{B}_{p}(\mathbb{R}^{N})=\sigma(\tau_{p})$.

One of the fundamental properties of $\mathbb{R}^{N}$ endowed with any of the norms $\Vert \cdot \Vert_{p}$ for $p \in [1,\infty]$ is that its open sets can be described by means of a countable family of open balls. Indeed, $\mathcal{O}$ is open in $\mathbb{R}^{N}$ with respect to the norm $\Vert \cdot \Vert_{p}$ for $p \in [1,\infty]$ if and only if there exists a sequence $(B_{k})$ of open balls with respect to the norm $\Vert \cdot \Vert_{p}$ for $p \in [1,\infty]$ such that $\mathcal{O}=\bigcup_{k=1}^{\infty} B_{k}$\footnote{See, for example, \cite{Clapp, Rudin2}.}. Thus, if $\mathcal{B}_{p}$ is the class of open balls in $\mathbb{R}^{N}$ with respect to the norm $\Vert \cdot \Vert_{p}$ for $p \in [1,\infty]$, then $\mathcal{B}_{p}(\mathbb{R}^{N})=\sigma(\mathcal{B}_{p})$.

We now show that the Borel $\sigma$-algebras of $\mathbb{R}^{N}$ associated with the topologies $\tau_{p}$ coincide for every $p \in [1,\infty]$. This fact is a consequence of the equivalence of the norms $\Vert \cdot \Vert_{p}$ on $\mathbb{R}^{N}$ and, therefore, of the fact that they generate the same collection of open sets in $\mathbb{R}^{N}$.

\begin{theorem}\label{233}
$\mathcal{B}_{p}(\mathbb{R}^{N})=\sigma(\tau_{p})=\sigma(\tau_{q})=\mathcal{B}_{q}(\mathbb{R}^{N})$ for any $p,q\in [1,\infty]$.
\end{theorem} \index{sigma@$\sigma$!algebra!of Borel on $\mathbb{R}^{N}$}

\begin{proof}
We will prove that $\tau_{p}=\tau_{q}$ for any $p,q\in [1,\infty]$.

{\scshape Case 1.} \quad $p=1$ and $q \in (1,\infty]$.

$\subset):$ Let $\bar{x}\in \mathbb{R}^{N}$ be arbitrary. Note that
$$
|x_{i}|^{q} \leq \sum_{j=1}^{N} |x_{j}|^{q} \quad \forall i=1,\ldots,N
$$
and, therefore,
$$
|x_{i}| \leq \left( \sum_{j=1}^{N} |x_{j}|^{q}\right)^{1/q} \quad \forall i=1,\ldots,N.
$$

Consequently, $\Vert \bar{x}\Vert_{1} \leq N\Vert \bar{x}\Vert_{q}$.

Let $\mathcal{O} \in \tau_{1}$ and $\bar{x} \in \mathcal{O}$ be given. There exists $\varepsilon>0$ such that $B_{1}(\bar{x},\varepsilon) \subset \mathcal{O}$. From the previous inequality we conclude that $B_{q}(\bar{x},\frac{\varepsilon}{N}) \subset B_{1}(\bar{x},\varepsilon) \subset \mathcal{O}$ and therefore that $\mathcal{O} \in \tau_{q}$.

$\supset):$ Conversely, for arbitrary $\bar{x} \in \mathbb{R}^{N}$, one has
$$
|x_{i}|\leq \sum_{j=1}^{N}|x_{j}|=\Vert \bar{x} \Vert_{1} \quad \forall i=1,\ldots,N,
$$
and therefore
$$
|x_{i}|^{q} \leq \Vert \bar{x}\Vert_{1}^{q} \quad \forall i=1,\ldots,N.
$$

Thus, $\Vert \bar{x}\Vert_{q} \leq N^{1/q}\Vert\bar{x}\Vert_{1}$.

Given $\mathcal{O} \in \tau_{q}$ and $\bar{x} \in \mathcal{O}$, there exists $\varepsilon>0$ such that $B_{q}(\bar{x},\varepsilon) \subset \mathcal{O}$. The previous inequality guarantees that $B_{1}\left(\bar{x},\frac{\varepsilon}{N^{1/q}}\right) \subset B_{q}(\bar{x},\varepsilon) \subset \mathcal{O}$. Consequently, $\mathcal{O} \in \tau_{1}$.

{\scshape Case 2.} \quad $p \in [1,\infty)$ and $q =\infty$.

Again, let $\bar{x} \in \mathbb{R}^{N}$ be arbitrary. Then,
$$
|x_{i}|^{p} \leq \max_{1\leq j\leq N}|x_{j}|^{p} \quad \forall i=1,\ldots,N
$$
and hence
$$
\sum_{i=1}^{N}|x_{i}|^{p} \leq N \max_{1\leq j\leq N}|x_{j}|^{p}.
$$

That is,
$$
\Vert \bar{x} \Vert_{p} \leq N^{1/p}\big(\max_{1\leq j \leq N}|x_{j}|^{p}\big)^{1/p}=N^{1/p}\Vert \bar{x}\Vert_{\infty}.
$$

Now,
$$
|x_{i}|^{p}\leq \sum_{j=1}^{N}|x_{j}|^{p} \quad \forall i=1,\ldots,N
$$
so that
$$
|x_{i}|\leq \left(\sum_{j=1}^{N}|x_{j}|^{p}\right)^{1/p} \quad \forall i=1,\ldots,N.
$$

Consequently, $\Vert \bar{x}\Vert_{\infty} \leq \Vert \bar{x} \Vert_{p}$.

The equality between $\tau_{p}$ and $\tau_{\infty}$ follows from the inclusions
$$
B_{p}(\bar{x},\varepsilon)\subset B_{\infty}(\bar{x},\varepsilon) \quad \text{and} \quad B_{\infty}\left(\bar{x},\frac{\varepsilon}{N^{1/p}}\right) \subset B_{p}(\bar{x},\varepsilon)
$$
valid for any $\bar{x} \in \mathbb{R}^{N}$ and $\varepsilon>0$.
\end{proof}

Consequently, we shall denote by $\mathcal{B}(\mathbb{R}^{N})$ the Borel $\sigma$-algebra of $\mathbb{R}^{N}$, without making reference to the specific norm under consideration.

Exercise \ref{E229} presents several classes of subsets of $\mathbb{R}^{N}$ that generate $\mathcal{B}(\mathbb{R}^{N})$, which are essentially analogous to those described in Proposition \ref{229} for $\mathbb{R}$.

\section{Sequences of sets}

In this section we make precise the concept of a sequence of subsets of a given set $X$ and establish an appropriate notion of convergence for this type of sequence. These ideas make it possible to describe the asymptotic behavior of families of sets and are fundamental for the subsequent development of the theory, in particular for the study of limits, continuity, and properties of measures.

\begin{definition}\label{234}
A \textbf{sequence of subsets} of a nonempty set $X$ is a function $f:\mathbb{N}\to\mathcal{P}(X)$. The value of this function at $k\in\mathbb{N}$ is called the $k$-th term of the sequence and is denoted by $A_k=f(k)$. A sequence is simply denoted by $(A_k)$.
\index{sequence!of subsets}.
\end{definition}

\begin{example}\label{235}
Let $(a_k)$ be a sequence of real numbers. Define a sequence $(A_k)$ of subsets of $\mathbb{R}$ by
$$
A_k:=\{x\in\mathbb{R}: x\le a_k\}=(-\infty,a_k].
$$

In this way, to each term of the numerical sequence $(a_k)$ there corresponds a subset of $\mathbb{R}$.
\end{example}

In many situations it is convenient to decompose a sequence of sets into a family of pairwise disjoint subsets that preserves the same accumulated information. The following proposition formalizes this idea under minimal algebraic hypotheses and will be useful in several later constructions.

\begin{proposition}\label{236}
Let $\mathfrak{R} \subset \mathcal{P}(X)$ be a class closed under finite unions and differences. Let $(A_k)$ be a sequence of elements of $\mathfrak{R}$. Then, there exists a sequence $(E_k)$ of elements of $\mathfrak{R}$ with the following properties:
\begin{itemize}
\item[(a)] $E_k \subset A_k$ for every $k \in \mathbb{N}$,
\item[(b)] $E_k \cap E_j = \varnothing$ for every $k \neq j$,
\item[(c)] $\displaystyle\bigcup_{k=1}^{j} E_{k} = \bigcup_{k=1}^{j} A_k$ for every $j \in \mathbb{N}$.
\end{itemize}
\end{proposition}

\begin{proof}
\textit{(a):} We define the sequence $(E_k)$ as follows:
$$
E_{1}:=A_1 \quad\mbox{and}\quad E_k :=A_k\smallsetminus \bigcup_{i=1}^{k-1}A_{i}\quad \forall k \geq 2.
$$

By definition, $E_k \subset A_k$ for every $k \in \mathbb{N}$ and $(E_k)$ is a sequence of elements of $\mathfrak{R}$ since this class is closed under finite unions and differences.

\begin{figure}[ht!]
\centering
\begin{tikzpicture}[scale=1]

\begin{scope}
  \clip (-3,-2) rectangle (3,2); 

  \draw[dotted] (-1.2,0) circle (1.4);
  \draw[dotted] ( 0.8,0) circle (1.4);
  \draw[dotted] (-0.2,1.1) circle (1.4);
\end{scope}

\begin{scope}
  \clip (-1.2,0) circle (1.4);
  \fill[lime!30] (-3,-3) rectangle (3,3);
\end{scope}

\begin{scope}
  \clip (0.8,0) circle (1.4);      
  \begin{scope}
    \clip (-1.2,0) circle (1.4);   
    \fill[white] (-3,-3) rectangle (3,3); 
  \end{scope}
  \fill[red!20] (-3,-3) rectangle (3,3);
\end{scope}

\begin{scope}
  \clip (-0.2,1.1) circle (1.4);   
  \begin{scope}
    \clip (-1.2,0) circle (1.4);
    \fill[white] (-3,-3) rectangle (3,3);
  \end{scope}
  \begin{scope}
    \clip (0.8,0) circle (1.4);
    \fill[white] (-3,-3) rectangle (3,3);
  \end{scope}
  \fill[blue!20] (-3,-3) rectangle (3,3);
\end{scope}

\draw[dotted] (-1.2,0) circle (1.4);
  \draw[dotted] ( 0.8,0) circle (1.4);
  \draw[dotted] (-0.2,1.1) circle (1.4);

  \node at (-0.2,2) {$_{E_{1}}$};
  \node at (1.5,0.2) {$_{E_{2}}$};
  \node at (-1.75,0.2) {$_{E_3}$};
\end{tikzpicture}
\end{figure}

\textit{(b):} Let $k \neq j$ and, without loss of generality, suppose that $k > j$. Then
$$
E_k :=A_{k}\smallsetminus \bigcup_{i=1}^{k-1}A_{i}=A_k\, \cap \left( X\smallsetminus  \left( \bigcup_{i=1}^{k-1} A_i\right)\right)  \subset \bigcap_{i=1}^{k-1} (X\smallsetminus A_i) \subset X\smallsetminus A_{j}
$$
since $j \in \{1,2,\ldots,k-1 \}$ and $E_{j} \subset A_{j}$ by \textit{(a)}. Thus,
$$
E_k \cap E_j \subset (X\smallsetminus A_j) \cap E_j \subset (X\smallsetminus A_j) \,\cap A_j = \varnothing
$$
which establishes that $E_k \cap E_j = \varnothing$.

\textit{(c):} Let $j \in \mathbb{N}$. Since $E_k \subset A_k$ for every $k \in \mathbb{N}$, then $\displaystyle\bigcup_{k=1}^{j} E_k \subset \bigcup_{k=1}^j A_k$.

Conversely, if $x \in \displaystyle\bigcup_{k=1}^{j} A_k$ then
$$
x \in A_{k}\quad \mbox{for some}\quad k \in \{1,\ldots,j \}.
$$

Let $k_{0}$ be the first natural number such that $x \in A_{k_{0}}$. Then,
$$
 x \in A_{k_{0}}\smallsetminus \bigcup_{k=1}^{k_{0}-1} A_k =E_{k_{0}}
$$
which shows that $x \in \displaystyle\bigcup_{k=1}^{j} E_k$.
\end{proof}

For a sequence of subsets $(A_k)$ of a set $X$, it is natural to introduce notions of limit that describe their long-term behavior. For this purpose, we define below the limit inferior and the limit superior.

\begin{definition}\label{237}
Let $(A_k)$ be a sequence of subsets of $X$. We define the \textbf{limit inferior} \index{limit!inferior of subsets} and the \textbf{limit superior} \index{limit!superior of subsets} of $(A_k)$ as the following subsets of $X$:
\begin{itemize}
\item[(1)] $\displaystyle\liminf_{k \to \infty} A_k = \displaystyle\bigcup_{k=1}^{\infty}\bigcap_{j=k}^{\infty} A_j$.
\item[(2)] $\displaystyle\limsup_{k \to \infty} A_k = \displaystyle\bigcap_{k=1}^{\infty}\bigcup_{j=k}^{\infty} A_j$.
\end{itemize}
\end{definition}

\begin{proposition}\label{238}
For every sequence $(A_k)$ of subsets of $X$, the following hold:
\begin{itemize}
\item[(a)] $\displaystyle\liminf_{k \to \infty} A_k \subset \lim\sup_{k \to \infty} A_k .$
\item[(b)] $X\smallsetminus \displaystyle\liminf_{k \to \infty} A_k = \limsup_{k \to \infty} (X\smallsetminus A_k)$.
\item[(c)] $X\smallsetminus \displaystyle\limsup_{k \to \infty} A_k = \liminf_{k \to \infty} (X\smallsetminus A_k)$.
\end{itemize}
\end{proposition}

\begin{proof}
\textit{(a):} The proof is straightforward and is left as an exercise [Exercise \ref{E230}].

\textit{(b) \text{ and } (c):} These follow by applying De Morgan's laws,
$$
\begin{aligned}
X\smallsetminus  \liminf_{k \to \infty} A_{k} &= X \smallsetminus  \bigcup_{k=1}^{\infty}\bigcap_{j=k}^{\infty } A_{j} \\
&= \bigcap_{k=1}^{\infty} \left(X \smallsetminus  \bigcap_{j=k}^{\infty} A_{j} \right)\\
&= \bigcap_{k=1}^{\infty} \bigcup_{j=k}^{\infty} (X\smallsetminus A_{j})\\
&= \limsup_{k \to \infty} (X\smallsetminus  A_{k}) ;
\end{aligned}
$$

$$
\begin{aligned}
X\smallsetminus  \limsup_{k \to \infty} A_{k} &= X \smallsetminus  \bigcap_{k=1}^{\infty}\bigcup_{j=k}^{\infty } A_{j} \\
&= \bigcup_{k=1}^{\infty} \left(X \smallsetminus  \bigcup_{j=k}^{\infty} A_{j} \right)\\
&= \bigcup_{k=1}^{\infty} \bigcap_{j=k}^{\infty} (X\smallsetminus A_{j})\\
&= \liminf_{k \to \infty} (X\smallsetminus  A_{k}).
\end{aligned}
$$

This completes the proof.
\end{proof}

\begin{definition}\label{239}
Let $(A_k)$ be a sequence of subsets of $X$. If there exists a subset $A$ \index{limit!of subsets} of $X$ such that
$$
\liminf_{k \to \infty} A_k=\limsup_{k \to \infty} A_k =A
$$
then we say that the sequence \textbf{converges} to the set $A$. The set $A$ is called the limit of the sequence $(A_k)$ and will be denoted by
$$
\lim_{k \to \infty} A_{k} = A, \quad \text{or}\quad A_k\to A.
$$
\end{definition}

Exercise \ref{E231} establishes a relationship between the convergence of a sequence of subsets and the pointwise convergence of a sequence of functions.

\begin{example}\label{240}
Let $A \subset X$ be arbitrary but fixed and nonempty. Define the sequence $(A_k)$ by
$$
A_{k}:=\left\{\begin{array}{lcl}
X\smallsetminus A & & if\,\, k \,\,is\,\, even,\\
A & & if\,\, k \,\,is\,\, odd.
\end{array}
\right.
$$

The sequence $(A_k)$ is not convergent.
\end{example}

\begin{proof}
Since
$$
\liminf_{k \to \infty} A_k = \displaystyle\bigcup_{k=1}^{\infty}\bigcap_{j=k}^{\infty} A_j=\bigcup_{k=1}^{\infty} \varnothing=\varnothing
$$
and
$$
\limsup_{k \to \infty} A_k = \displaystyle\bigcap_{k=1}^{\infty}\bigcup_{j=k}^{\infty} A_j =\bigcap_{k=1}^{\infty} X=X,
$$
then $\displaystyle\lim_{k \to \infty} A_k$ does not exist.
\end{proof}

\begin{example}\label{241}
Let $X=[-1,1]$ and define for each $k \in  \mathbb{N}$ the set $A_k:=[-1/k,0]$ if $k$ is odd and $A_k:=[0,1/k]$ if $k$ is even. Then $(A_k)$ is a convergent sequence.
\end{example}

\begin{proof}
Since
$$
\liminf_{k \to \infty} A_k = \displaystyle\bigcup_{k=1}^{\infty}\bigcap_{j=k}^{\infty} A_j=\bigcup_{k=1}^{\infty} [0,0]=\{0\}
$$
and
$$
\limsup_{k \to \infty} A_k = \displaystyle\bigcap_{k=1}^{\infty}\bigcup_{j=k}^{\infty} A_j =\bigcap_{k=1}^{\infty} ([-1/k,0]\cup[0,1/k])=[0,0]=\{0\},
$$
then
$$
\lim_{k \to \infty} A_{k}=\liminf_{k \to \infty} A_k = \limsup_{k \to \infty} A_k = \{ 0 \}.
$$

Therefore, the sequence $(A_k)$ converges to the set $\{ 0 \}$.
\end{proof}

In addition to the notions of limit, it is possible to introduce a concept of monotonicity for sequences of sets. In particular, a sequence may be increasing or decreasing according to the inclusion relation between its terms.

\begin{proposition}\label{242}
Let $(A_k)$ be a sequence \index{sequence!of monotone subsets @ of monotone subsets} of subsets of a given set $X$.
\begin{itemize}
\item[(a)] If $(A_k)$ is increasing, that is, $A_{k} \subset A_{k+1}$ for every $k \in \mathbb{N}$, then $\displaystyle\lim_{k \to \infty} A_k = \bigcup_{k=1}^{\infty} A_k$.
\item[(b)] If $(A_k)$ is decreasing, that is, $A_{k+1} \subset A_{k}$ for every $k \in \mathbb{N}$, then $\displaystyle\lim_{k\to\infty} A_k = \bigcap_{k=1}^{\infty} A_k$.
\end{itemize}
\end{proposition}

\begin{proof}
\textit{(a):} Since the sequence $(A_{k})$ is increasing, then
$$
\bigcup_{j=k}^{\infty} A_{j} = \bigcup_{j=1}^{\infty} A_{j} \quad \text{and}\quad \bigcap_{j=k}^{\infty} A_{j} = A_{k} \qquad \text{for every }k \in \mathbb{N}.
$$

Therefore,
$$
\begin{aligned}
\limsup_{k \to \infty} A_{k} &= \bigcap_{k=1}^{\infty}\bigcup_{j=k}^{\infty} A_{j} = \bigcap_{k=1}^{\infty} \bigcup_{j=1}^{\infty} A_{j} = \bigcup_{j=1}^{\infty} A_{j}\\
\mbox{and}\,\,\lim\inf_{k \to \infty} A_{k} &= \bigcup_{k=1}^{\infty}\bigcap_{j=k}^{\infty} A_{j} = \bigcup_{k=1}^{\infty} A_{k}.
\end{aligned}
$$

That is, $\lim_{k \to\infty}A_{k}=\bigcup_{k=1}^{\infty}A_{k}$.

\begin{figure}[ht!]
\begin{minipage}[r]{0.5\textwidth}
\begin{center}
\begin{tikzpicture}[xscale=0.65, yscale=0.65]
\draw[fill=red!20] (0,0) circle (3);
\draw[fill=red!15] (0,0) circle (2);
\draw[fill=red!10] (0,0) circle (1);

\draw (0,0) node{$_{A_{k}}$};
\draw (0,1.5) node{$_{A_{k+1}}$};
\draw (0,2.65) node{${\vdots}$};
\draw (0,3.3) node{$_{A_{k+j}}$};
\end{tikzpicture}
\begin{center}
Increasing sequence $(A_k)$
\end{center}
\end{center}
\end{minipage} \hfill 
 \begin{minipage}[l]{0.5\textwidth}
\begin{center}
\begin{tikzpicture}[xscale=0.65,yscale=0.65]
\draw[fill=blue!20] (0,0) circle (3);
\draw[fill=blue!15] (0,0) circle (2);
\draw[fill=blue!10] (0,0) circle (1);

\draw (0,0) node{$_{A_{k+j}}$};
\draw (0,1.65) node{${\vdots}$};
\draw (0,2.5) node{$_{A_{k+1}}$};
\draw (0,3.3) node{$_{A_{k}}$};
\end{tikzpicture}
\begin{center}
Decreasing sequence $(A_k)$
\end{center}
\end{center}
\end{minipage}
\end{figure}

\textit{(b):} Since $(A_{k})$ is a decreasing sequence, then $(X \smallsetminus A_{k})$ is an increasing sequence. The previous part guarantees that
$$
\liminf_{k \to \infty} (X \smallsetminus A_{k})=\limsup_{k \to \infty}(X\smallsetminus A_{k})=\bigcup_{k=1}^{\infty}(X \smallsetminus A_{k})=X \smallsetminus \bigcap_{k=1}^{\infty}A_{k}.
$$

From parts (\textit{b}) and (\textit{c}) of Proposition \ref{238} we conclude that
$$
\limsup_{k \to \infty} A_{k} = \liminf_{k \to \infty}A_{k} = \bigcap_{k=1}^{\infty}A_{k}.
$$

Therefore, $\lim_{k \to \infty}A_{k}=\bigcap_{k=1}^{\infty}A_{k}$.
\end{proof}

\begin{remark}\label{243}
Let $(A_k)$ be a sequence of subsets of $X$ such that $A_k \in \mathsf{S}$ for every $k \in \mathbb{N}$, where $\mathsf{S}$ is a $\sigma$-algebra of subsets of $X$. Then,
$$
\liminf_{k \to \infty} A_k \in \mathsf{S} \quad \text{and} \quad \limsup_{k \to \infty} A_{k} \in \mathsf{S}.
$$
\end{remark}

\section{Dynkin classes}

Let $\mathcal{C} \subset \mathcal{P}(X)$ be a given class such that all its elements satisfy a property $\wp$. One way to prove that all elements of $\sigma(\mathcal{C})$ also satisfy $\wp$ is to show that the set $\{A \subset X\,:\,A\,\,\mbox{satisfies}\,\,\wp \}$ is a $\sigma$-algebra; however, in specific cases this turns out to be complicated. In this section we establish two notions of classes of subsets of $X$ whose study will facilitate the problem described above. These classes were introduced in 1961 by the mathematician E. B. Dynkin\footnote{Eugene Borisovich Dynkin (1924-2014) was a Soviet and American mathematician. He made contributions to the fields of probability and algebra, especially to semisimple Lie groups, Lie algebras, and Markov processes. Dynkin diagrams, Dynkin systems, and Dynkin's lemma are named after him.}.

\begin{figure}[ht!]
\centering
\includegraphics[scale=0.25]{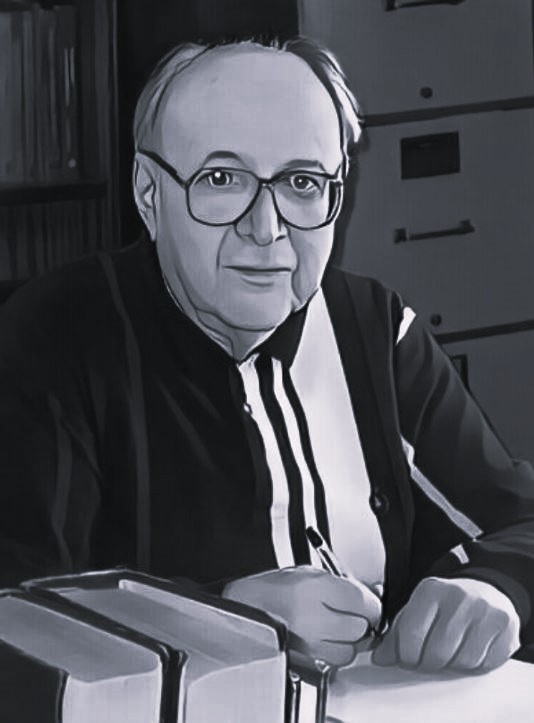} 
\begin{center}
Eugene Dynkin (1924-2014)
\end{center}
\end{figure}

\begin{definition}\label{244}
A nonempty class $\mathcal{C} \subset \mathcal{P}(X)$ is called a \textbf{$\pi$-system} of subsets of $X$ \index{pi@$\pi$!system} if $A_{1} \cap \cdots \cap A_{N} \in \mathcal{C}$ for any $A_{1},\ldots,A_{N} \in \mathcal{C}$.
\end{definition}

\begin{example}\label{245}
Let $X$ be an arbitrary set. The classes $\mathcal{C}_{1}=\{\varnothing,X \}$ and $\mathcal{C}_{2}=\mathcal{P}(X)$ are $\pi$-systems. However, unlike a $\sigma$-algebra, the class $\mathcal{C}_1$ is not the smallest $\pi$-system that can be associated with the set $X$ (see {\rm Example \ref{27}}), since $\mathcal{C}_{0}=\{X\}$ is also a $\pi$-system and $\mathcal{C}_{0} \subset \mathcal{C}_1$.
\end{example}

\begin{example}\label{246}
Let $(X,\tau)$ be a topological space. The classes $\tau$ and $X \smallsetminus \tau:=\{X \smallsetminus A\,:\,A \in \tau\}$ are $\pi$-systems.
\end{example}

\begin{example}\label{247}
Let $\mathbb{R}^{2}_{2}=(\mathbb{R}^{2},\Vert\cdot\Vert_{2})$. The class $\mathcal{C}=\{ B \subset \mathbb{R}^{2}_{2}\,:\,B\,\,\mbox{is an open ball} \}$ is not a $\pi$-system.
\end{example}

\begin{definition}\label{248}
 A nonempty class $\mathcal{L} \subset \mathcal{P}(X)$ is called a \textbf{$\lambda$-system} (or Dynkin system or Dynkin class) of subsets of $X$ \index{lambda@$\lambda$!system} \index{class!Dynkin} if it satisfies the following properties:
\begin{itemize}
\item[\rm (L1)] $X \in \mathcal{L}$.
\item[\rm (L2)] $A\smallsetminus B \in \mathcal{L}$ for any $A,B \in \mathcal{L}$ with $B \subset A$.
\item[\rm (L3)] If $(A_k)$ is an increasing sequence of elements of $\mathcal{L}$, then $\displaystyle\bigcup_{k=1}^{\infty} A_k \in \mathcal{L}.$
\end{itemize}
\end{definition}

Let $\mathcal{L}$ be a $\lambda$-system and let $A \in \mathcal{L}$ be arbitrary. From properties (L1) and (L2) we obtain that $X\smallsetminus A \in \mathcal{L}$ since $A \subset X$. Therefore, $\mathcal{L}$ is closed under complements.

\begin{proposition}\label{249}
Let $X$ be a set and $\mathcal{L} \subset \mathcal{P}(X)$ a $\lambda$-system.
\begin{itemize}
\item[(1)] If $(A_k)$ is a decreasing sequence of elements of $\mathcal{L}$ then $\displaystyle\bigcap_{k=1}^{\infty} A_k \in \mathcal{L}$.
\item[(2)] If $A,B \in \mathcal{L}$ are disjoint then $A \cup B \in \mathcal{L}$.
\end{itemize}
\end{proposition}

\begin{proof}
\textit{(1):} For each $k \in \mathbb{N}$, $A_{k+1} \subset A_{k}$ if and only if $X\smallsetminus A_{k} \subset X\smallsetminus A_{k+1}$. Then $(X\smallsetminus A_k)$ is an increasing sequence of elements of $\mathcal{L}$ and, from property (L3), we obtain that
$$
\bigcup_{k=1}^{\infty} (X\smallsetminus A_k) = X\smallsetminus\bigcap_{k=1}^{\infty} A_k \in \mathcal{L}.
$$

Consequently, $\displaystyle\bigcap_{k=1}^{\infty} A_k \in \mathcal{L}$ since $\mathcal{L}$ is closed under complementation.

\textit{(2):} Since $A \cap B = \varnothing$ if and only if $B \subset X\smallsetminus A$, then $(X\smallsetminus A)\smallsetminus B=X\smallsetminus (A\cup B) \in \mathcal{L}$ by property (L2). Therefore, $A \cup B \in \mathcal{L}$.
\end{proof}

\begin{example}\label{250}
Let $X$ be a nonempty set and $\mathsf{S}$ a $\sigma$-algebra of subsets of $X$. Clearly, $\mathsf{S}$ satisfies all the properties of a $\lambda$-system and a $\pi$-system. However, not every $\lambda$-system is a $\sigma$-algebra of subsets of $X$. Indeed, let $X=\{1,2,3,4 \}$ and
$$
\mathcal{L}=\left\{ \varnothing, \{1,2 \},\{1,3\},\{1,4\},\{2,3\},\{2,4\},\{3,4\},X \right\}.
$$

Clearly $X \in \mathcal{L}$. Let $A,B \in  \mathcal{L}$ with $B \subset A$. The cases when $B=\varnothing$ or $B=A$ are trivial since $A\smallsetminus \varnothing=A \in \mathcal{L}$ and $A\smallsetminus A=\varnothing \in \mathcal{L}$. Thus, we have the following: $A=X$ and $B \neq X \neq \varnothing$. If $B=\{1,2 \}$ then $A\smallsetminus B=\{3,4\} \in \mathcal{L}$. The remaining cases are proved similarly. Consequently, (L2) holds. Finally, for any increasing sequence $(A_k)$ of elements of $\mathcal{L}$ one has $\bigcup_{k=1}^{\infty} A_k =X \in \mathcal{L}$. Therefore, $\mathcal{L}$ is a $\lambda$-system but not a $\sigma$-algebra since property (S2) does not hold. Indeed, if $A=\{1,2\}$ and $B=\{2,3\}$ then $A\smallsetminus B=\{1\} \not \in \mathcal{L}$.
\end{example}

\begin{proposition}\label{251}
Every class $\mathcal{L} \subset \mathcal{P}(X)$ that is both a $\pi$-system and a $\lambda$-system is itself a $\sigma$-algebra.
\end{proposition}

\begin{proof}
Since $\mathcal{L}$ is closed under complementation, contains $X$, and is closed under increasing unions, it is enough to prove that $\mathcal{L}$ is closed under unions [Exercise \ref{E26}]. For any $A,B \in \mathcal{L}$ one has that $A \cup B= A \cup (B\smallsetminus (A \cap B))$ where $A \cap (B\smallsetminus(A \cap B))=\varnothing$ and $B\smallsetminus (A \cap B) \in \mathcal{L}$ by (L2). That is, $A \cup B$ is a disjoint union of elements of $\mathcal{L}$ and, consequently, $A \cup B \in \mathcal{L}$ by Proposition \ref{249}.
\end{proof}

The converse of the previous result is also valid since every $\sigma$-algebra satisfies the properties of both a $\pi$ and a $\lambda$-system.

\begin{proposition}\label{252}
The intersection $\bigcap_{i\,\in\,\mathcal{I}}\, \mathcal{C}_{i}$ of any family $\{\mathcal{C}_{i}\,:\,i\in\mathcal{I} \}$ of $\pi$-systems is a $\pi$-system.
\end{proposition}

\begin{proposition}\label{253}
The intersection $\bigcap_{\tau\,\in\,\mathcal{T}}\, \mathcal{L}_{\tau}$ of any family $\{\mathcal{L}_{\tau}\,:\,\tau\in\mathcal{T} \}$ of $\lambda$-systems is a $\lambda$-system.
\end{proposition}

The proofs of these statements are straightforward and are left as an exercise [Exercise \ref{E240}]. These results guarantee that the following definitions make sense.

\begin{definition}\label{254}
Let $\mathcal{E} \subset \mathcal{P}(X)$ be arbitrary. Consider the following nonempty classes in $\mathcal{P}(X)$
$$\mathscr{C}=\{\mathcal{C} \subset \mathcal{P}(X)\,:\,\mathcal{C}\,\,\mbox{is a}\,\,\pi-\mbox{system}\,\,\mbox{and}\,\,\mathcal{E} \subset \mathcal{C} \}, $$
$$
\mathscr{L}=\{\mathcal{L} \subset \mathcal{P}(X)\,:\,\mathcal{L}\,\,\mbox{is a}\,\,\lambda-\mbox{system}\,\,\mbox{and}\,\,\mathcal{E} \subset \mathcal{L} \}.
$$

Then, the \textbf{$\pi$-system generated by $\mathcal{E}$}, denoted by $\pi(\mathcal{E})$, is the class \index{pi@$\pi$! generated system! by a class $\pi(\mathcal{E})$} 
$$
\pi(\mathcal{E}):=\bigcap\{\mathcal{C} \subset \mathcal{P}(X)\,:\,\mathcal{C} \in \mathscr{C}\}.
$$

The \textbf{$\lambda$-system generated by $\mathcal{E}$}, $\mathcal{L}(\mathcal{E})$, is the class \index{lambda@$\lambda$!generated system!by a class $\mathcal{L}(\mathcal{E})$}
$$
\mathcal{L}(\mathcal{E}):=\bigcap\{\mathcal{L} \subset \mathcal{P}(X)\,:\,\mathcal{L} \in \mathscr{L} \}.
$$
\end{definition}

Analogously to the case of the $\sigma$-algebra generated by $\mathcal{E}$, it is possible to prove the uniqueness of the classes defined above. Likewise, the properties corresponding to Theorem \ref{215}, Proposition \ref{219}, Proposition \ref{220}, and Theorem \ref{221} remain valid for the classes $\pi(\mathcal{E})$ and $\mathcal{L}(\mathcal{E})$, as will be seen in the proof of the following theorem. This result gathers together the fundamental aspects needed to solve the problem posed at the beginning of this section.

\begin{theorem}[Dynkin]\label{255} \index{theorem!Dynkin}
Let $\mathcal{C} \subset \mathcal{P}(X)$ be a $\pi$-system and $\mathcal{L}_{0} \subset \mathcal{P}(X)$ a $\lambda$-system such that $\mathcal{C} \subset \mathcal{L}_{0}$. Then $\sigma(\mathcal{C}) \subset \mathcal{L}_{0}$.
\end{theorem}

To prove Dynkin's theorem we will use the following three lemmas under their general hypotheses. Let $\mathscr{L}$ be the class of all $\lambda$-systems of subsets of $X$ that contain the $\pi$-system $\mathcal{C}$. Then,
$$
\mathcal{L}(\mathcal{C})=\bigcap\{\mathcal{L} \subset \mathcal{P}(X)\,:\,\mathcal{L} \in \mathscr{L}\}.
$$

\begin{lemma}\label{256}
For every $A \in \mathcal{C}$, define the class $\mathcal{L}_{A}:=\{F \subset X\,:\, F \cap A \in \mathcal{L}(\mathcal{C}) \}$. Then $\mathcal{L}_{A}$ is a $\lambda$-system such that $\mathcal{C} \subset \mathcal{L}_{A}$ and $\mathcal{L}(\mathcal{C}) \subset \mathcal{L}_{A}$.
\end{lemma}

\begin{proof}
Since $A \cap X = A$ and $A\in \mathcal{L}(\mathcal{C})$, it follows that $X \in \mathcal{L}_{A}$. 

Let $E,F \in \mathcal{L}_{A}$ be such that $E \subset F$. Then, $E \cap A, F \cap A \in \mathcal{C} \subset \mathcal{L}(\mathcal{C})$ and satisfy $E \cap A \subset F \cap A$. By property (L2) we have $(F \cap A)\smallsetminus (E \cap A)=(F\smallsetminus E)\cap A \in \mathcal{L}(\mathcal{C})$ and therefore $F\smallsetminus E \in \mathcal{L}_{A}$. 

Let $(F_k)$ be an increasing sequence of elements of $\mathcal{L}_{A}$. Since $F_k \cap A \in \mathcal{C}\subset \mathcal{L}(\mathcal{C})$ for each $k \in \mathbb{N}$, the sequence $(F_k \cap A)$ is increasing in $\mathcal{L}(\mathcal{C})$ so that $\bigcup_{k=1}^{\infty}\,(F_k \cap A) = (\bigcup_{k=1}^{\infty} F_k) \cap A \in \mathcal{L}(\mathcal{C})$ by property (L3). Then, $\bigcup_{k=1}^{\infty} F_k \in \mathcal{L}_{A}$ and we conclude that $\mathcal{L}_{A}$ is a $\lambda$-system.

Finally, $\mathcal{C} \subset \mathcal{L}_{A}$. Indeed, let $B \in \mathcal{C}$ be arbitrary. Since $A \in \mathcal{C}$ and $\mathcal{C}$ is a $\pi$-system, it follows that $A \cap B \in \mathcal{C} \subset \mathcal{L}(\mathcal{C})$ and thus $B \in \mathcal{L}_{A}$. 

Consequently, $\mathcal{L}_{A}$ is a $\lambda$-system containing $\mathcal{C}$ and therefore $\mathcal{L}(\mathcal{C}) \subset \mathcal{L}_{A}$.
\end{proof}

\begin{lemma}\label{257}
For every $A \in \mathcal{L}(\mathcal{C})$ define the class $\mathcal{L}_{A}':=\{ F \subset X\,:\,  F\cap A \in \mathcal{L}(\mathcal{C}) \}$. Then $\mathcal{C} \subset \mathcal{L}_{A}'$ and $\mathcal{L}(\mathcal{C}) \subset \mathcal{L}_{A}'$.
\end{lemma}

\begin{proof}
Analogously to the previous lemma, it can be verified that the class $\mathcal{L}_{A}'$ is a $\lambda$-system. Let us now show that $\mathcal{C} \subset \mathcal{L}_{A}'$.
Let $D \in \mathcal{C}$ be arbitrary. By Lemma \ref{256} we have $\mathcal{L}(\mathcal{C}) \subset \mathcal{L}_{D}$. Since $A \in \mathcal{L}(\mathcal{C})$, it follows that $A \cap D \in \mathcal{L}(\mathcal{C})$, which implies that $D \in \mathcal{L}_{A}'$. Consequently,
$\mathcal{C} \subset \mathcal{L}_{A}'$.

Since $\mathcal{L}_{A}'$ is a $\lambda$-system containing $\mathcal{C}$, it follows that $\mathcal{L}(\mathcal{C}) \subset \mathcal{L}_{A}'$.
\end{proof}

\begin{lemma}\label{258}
$\mathcal{L}(\mathcal{C})$ is a $\pi$-system.
\end{lemma}

\begin{proof}
Let $A,B\in\mathcal{L}(\mathcal{C})$. By Lemma \ref{257}, one has $\mathcal{L}(\mathcal{C})\subset\mathcal{L}'_{A}$ and therefore $B\in\mathcal{L}'_{A}$. Consequently, $A\cap B\in\mathcal{L}(\mathcal{C})$, which proves that $\mathcal{L}(\mathcal{C})$ is a $\pi$-system.
\end{proof}

\vspace{0.15cm}

\begin{proof}[Proof of Theorem \ref{255}]
    Lemma \ref{258} guarantees that $\mathcal{L}(\mathcal{C})$ is a $\pi$-system. Since $\mathcal{L}(\mathcal{C})$ is a $\lambda$-system, then $\mathcal{L}(\mathcal{C})$ is a $\sigma$-algebra by Proposition \ref{251} and $\mathcal{C} \subset \mathcal{L}(\mathcal{C})$. Therefore $\sigma(\mathcal{C})\subset \mathcal{L}(\mathcal{C}) \subset \mathcal{L}_{0}$.
\end{proof}

\begin{corollary}\label{259}
If $\mathcal{C} \subset \mathcal{P}(X)$ is a $\pi$-system, then $\sigma(\mathcal{C})=\mathcal{L}(\mathcal{C})$.
\end{corollary}

\begin{proof}
Theorem \ref{255} guarantees that $\sigma(\mathcal{C}) \subset  \mathcal{L}(\mathcal{C})$. Conversely: by Proposition \ref{251}, $\sigma(\mathcal{C})$ is a $\lambda$-system containing $\mathcal{C}$ and therefore $\mathcal{L}(\mathcal{C}) \subset \sigma(\mathcal{C})$.
\end{proof}

\begin{remark}\label{260}
Let $\mathcal{C} \subset \mathcal{P}(X)$ be such that every $A \in \mathcal{C}$ satisfies a certain property $\wp$. According to {\rm Corollary \ref{259}}, if $\mathcal{C}$ is a $\pi$-system, then in order to prove that every element of $\sigma(\mathcal{C})$ satisfies the property $\wp$, it suffices to verify that the class $\mathcal{L}=\{ F \subset X\,:\,F\,\,\mbox{satisfies}\,\,\wp\}$ is a $\lambda$-system.
\end{remark}

The problem posed in this section can also be solved through the well-known \textbf{monotone class lemma}, a result that we leave here as an exercise. Some important applications of these results will be seen later.

\section{Exercises}
 
\begin{exercise} \label{E21}
Let $(A_k)$ be a sequence of subsets of an arbitrary set $X$. Prove that
$$\bigcap_{k=1}^\infty A_k= A_1\smallsetminus \bigcup_{k=1}^\infty (A_1\smallsetminus A_k).$$
\end{exercise}

{\setlength{\parindent}{0pt}
\begin{exercise} \label{E22}
Given a subset $A$ of $X$, define its \textbf{characteristic function} $\chi_{A}:X \to \{0,1\}$ by:
$$
\chi_{A}(x):=\left\{
\begin{array}{lcl}
1 & &\mbox{if}\,\,\,x \in A,\\
0 & &\mbox{if}\,\,\,x \not\in A.
\end{array}
\right.
$$

Let $\mathfrak{R} \subset \mathcal{P}(X)$ be a given nonempty class. Prove that $\mathfrak{R}$ is a ring of subsets of $X$ if and only if the set $\{ \chi_{A}\,:\,A \in \mathfrak{R} \}$ is an algebraic ring with the operations of addition and multiplication modulo 2.
\end{exercise}

(Hint: First prove that $\chi_{A} + \chi_{B} \equiv \chi_{A \bigtriangleup B}$ and $\chi_{A} \cdot \chi_{B} \equiv \chi_{A \cap B}$ (mod 2).).
}

\begin{exercise} \label{E23}
Let $X$ be an uncountable set. Prove that
\begin{itemize}
\item[(1)] The set $\mathcal{R}=\left\{ A \subset X\,:\, A\,\,\mbox{is finite or countable}\right\}$ is a $\sigma$-ring of subsets of $X$.
\item[(2)] The set $\mathsf{S}=\left\{ A \subset X\,:\, A\,\,\mbox{or}\,\,X\smallsetminus A\,\,\mbox{is finite or countable}\right\}$ is a $\sigma$-algebra of subsets of $X$.
\end{itemize}

Show moreover that $\mathcal{R} \neq \mathsf{S} \neq \mathcal{P}(X)$.
\end{exercise}

{\setlength{\parindent}{0pt}
\begin{exercise} \label{E24}
Prove, by means of an example, the following:
\begin{itemize}
    \item[(a)] A nonempty family $\mathcal{C}\subset \mathcal{P}(X)$ closed under countable intersections and symmetric difference is not necessarily a $\sigma$-ring.
    \item[(b)] A nonempty family $\mathcal{C} \subset \mathcal{P}(X)$ closed under countable unions and intersections and symmetric difference such that $X \in \mathcal{C}$ is not necessarily a $\sigma$-ring. 
\end{itemize}
\end{exercise}

(Hint: Consider $X=\mathbb{R}$ and $\mathcal{C}=\{(-\infty,a)\,:\,a\leq +\infty \}\cup \{ (-\infty,b)\,:\,b<+\infty \} $).}

{\setlength{\parindent}{0pt}
\begin{exercise} \label{E25}
Prove that if $\mathcal{A}$ is an algebra of subsets of $X$ with the property that for every sequence of pairwise disjoint subsets $(A_k)$ of $\mathcal{A}$ one has $\bigcup_{k=1}^{\infty} A_k \in \mathcal{A}$, then $\mathcal{A}$ is a $\sigma$-algebra of subsets of $X$.
\end{exercise}

(Hint: Use Proposition \ref{236}).}

\begin{exercise} \label{E26}
Let $\mathcal{A}$ be an algebra of subsets of $X$. Prove that $\mathcal{A}$ is a $\sigma$-algebra of subsets of $X$ if and only if $\mathcal{A}$ is closed under increasing unions.
\end{exercise}

\begin{exercise} \label{E27}
We say that a subset $A$ of $\mathbb{N}$ has asymptotic density if there exists $\theta \in \mathbb{R}$ such that
$$
\lim_{N \to\infty} \frac{|A\cap \{1,\ldots,N\}|}{N}=\theta,
$$
where $|A \cap \{1,\ldots,N\}|$ denotes the cardinality of the set $A \cap \{1,\ldots,N\}$.

Let $\mathfrak{A}=\{A \subset \mathbb{N}\,:\,A\,\text{ has asymptotic density}\}$. Is $\mathfrak{A}$ an algebra of subsets of $\mathbb{N}$? Is $\mathfrak{A}$ a $\sigma$-algebra of subsets of $\mathbb{N}$? Justify your answers in detail.
\end{exercise}

\begin{exercise} \label{E28}
Prove that every ring $\mathcal{R}$ of subsets of $X$ with finitely many elements is a $\sigma$-ring and that every algebra $\mathcal{A}$ of subsets of $X$ with finitely many elements is a $\sigma$-algebra.
\end{exercise}

\begin{exercise} \label{E29}
Let $X$ be a nonempty set. Given a subset $A$ of $X$, denote
$$
A^{\varepsilon}:=\left\{
\begin{array}{l c l}
A  &  & \text{if }\,\varepsilon=0,\\
X \smallsetminus A &  & \text{if }\,\varepsilon=1.\\
\end{array}
\right.
$$

Let $\mathcal{C}=\{A_{1},\ldots,A_{N}\}$ be a finite class of subsets of $X$. For each $\varepsilon:=(\varepsilon_{1},\ldots,\varepsilon_{N}) \in \{0,1\}^{N}$, define the set $B_{\varepsilon}:=\bigcap_{j=1}^{N} A^{\varepsilon_{j}}_{j}$.

Prove the following:
\begin{itemize}
    \item[(a)] $\mathcal{A}(\mathcal{C})=\left\{ \bigcup_{\varepsilon \in \mathcal{I}}B_{\varepsilon}\,:\, \mathcal{I} \subset \{0,1\}^{N} \,\,\text{ and }\,\, \mathcal{I} \neq \varnothing\right\} \cup \{\varnothing\}$.
    \item[(b)]  $|\mathcal{A}(\mathcal{C})| \leq 2^{2^{N}}$
\end{itemize}

What can you conclude about $\mathcal{R}(\mathcal{C})$, $\mathcal{SR}(\mathcal{C})$ and $\sigma(\mathcal{C})$? Justify your answer in detail.
\end{exercise}

\begin{exercise} \label{E210}
Let $X$ and $Y$ be nonempty sets. Prove that if $f:X \to Y$ is a function and $\mathsf{T}$ is a $\sigma$-algebra of subsets of $Y$, then
$$
f^{-1}(\mathsf{T})=\{ f^{-1}(A)  \,:\, A \in \mathsf{T}\}
$$
is a $\sigma$-algebra of subsets of $X$.
\end{exercise}

\begin{exercise} \label{E211}
Let $X$ and $Y$ be nonempty sets. Prove that if $f:X \to Y$ is a function and $\mathsf{S} \subset \mathcal{P}(X)$ is a $\sigma$-algebra of subsets of $X$, then
$$
\mathcal{F}=\{ F \subset Y\,:\, f^{-1}(F) \in \mathsf{S}\}
$$
is a $\sigma$-algebra of subsets of $Y$.
\end{exercise}

\begin{exercise} \label{E212}
Let $X$ be a nonempty set.
\begin{itemize}
    \item[(a)] Let $\mathsf{S}_{1}$ and $\mathsf{S}_{2}$ be two $\sigma$-algebras of subsets of $X$ such that $\mathsf{S}_{1} \subset \mathsf{S}_{2}$. Prove that $\mathsf{S}_{1} \cup \mathsf{S}_{2}$ is a $\sigma$-algebra of subsets of $X$.
    \item[(b)] Let $(\mathsf{S}_{k})$ be an increasing sequence of $\sigma$-algebras of subsets of $X$. Prove that $\mathcal{A}=\bigcup_{k=1}^{\infty} \mathsf{S}_{k}$ is an algebra of subsets of $X$.
    \item[(c)] Give an example of an increasing sequence $(\mathsf{S}_k)$ of $\sigma$-algebras of subsets of $X$ such that $\bigcup_{k=1}^{\infty} \mathsf{S}_{k}$ is not a $\sigma$-algebra of subsets of $X$.
\end{itemize}
\end{exercise}

\begin{exercise} \label{E213}
Let $\mathsf{S}_{1}$ and $\mathsf{S}_{2}$ be two $\sigma$-algebras on a set $X$. Prove or provide a counterexample to determine whether the difference $\mathsf{S}_{1} \smallsetminus \mathsf{S}_{2}:=\{A\smallsetminus B\,:\, A\in \mathsf{S}_{1}\,\,\text{and}\,\,B \in \mathsf{S}_{2}\}$ is a $\sigma$-algebra of subsets of $X$.
\end{exercise}

\begin{exercise} \label{E214}
Prove that a nonempty class $\mathsf{S}$ of subsets of $X$ is a $\sigma$-algebra if and only if $\mathsf{S}$ satisfies the following conditions
\begin{itemize}
\item[(1)] $\varnothing \in \mathsf{S}$.
\item[(2)] $X\smallsetminus A \in \mathsf{S}$ for every $A \in \mathsf{S}$.
\item[(3)] For every sequence $(A_k)$ of elements of $\mathsf{S}$ one has $\bigcap_{k=1}^{\infty} A_k \in \mathsf{S}$.
\end{itemize}
\end{exercise}

\begin{exercise} \label{E215}
Prove that a nonempty class $\mathcal{A}$ of subsets of $X$ is an algebra if and only if $\mathcal{A}$ satisfies the following conditions
\begin{itemize}
\item[(1)] $X \in \mathcal{A}$.
\item[(2)] $X\smallsetminus A  \in \mathcal{A}$ for every $A \in \mathcal{A}$.
\item[(3)] $A \cup B \in \mathcal{A}$ for every $A,B \in \mathcal{A}$.
\end{itemize}
\end{exercise}

\begin{exercise} \label{E216}
Prove that a nonempty class $\mathcal{A}$ of subsets of $X$ is an algebra if and only if $\mathcal{A}$ satisfies the following conditions
\begin{itemize}
\item[(1)] $X \in \mathcal{A}$.
\item[(2)] $X\smallsetminus A \in \mathcal{A}$ for every $A \in \mathcal{A}$.
\item[(3)] $A \cap B \in \mathcal{A}$ for every $A,B \in \mathcal{A}$. 
\end{itemize}
\end{exercise}

\begin{exercise} \label{E217}
Let $X$ be a nonempty set.
\begin{itemize}
    \item[(a)] Let $\mathsf{S}$ be a $\sigma$-algebra of subsets of $X$. Prove that the collection $X\smallsetminus \mathsf{S}:=\{X\smallsetminus A\,:\, A \in \mathsf{S} \}$ is a $\sigma$-algebra of subsets of $X$. Verify moreover that $\mathsf{S}$ and $X\smallsetminus \mathsf{S}$ coincide.
    \item[(b)] Let $\mathcal{C}$ be a class of subsets of $X$. Prove that if $X\smallsetminus \mathcal{C}:=\{X\smallsetminus A\,:\,A \in \mathcal{C}\}$, then $\sigma(\mathcal{C})=\sigma(X\smallsetminus \mathcal{C})$.
    \item[(c)] Let $(X,\tau)$ be a topological space. How do you interpret the previous item when $\mathcal{C}=\tau$?
\end{itemize}
\end{exercise}

{\setlength{\parindent}{0pt}
\begin{exercise} \label{E218}
Let $X=(0,1]$ and define the collection $\mathcal{A}$ consisting of the empty set and subsets of the form
$$
\bigcup_{i=1}^{k} (a_i,b_i]
$$
where $(a_i,b_i] \subset (0,1]$ with $(a_i,b_i] \cap (a_j,b_j] = \varnothing$ for $i \neq j$ and $k \in \mathbb{N}$. Prove that $\mathcal{A}$ is an algebra but not a $\sigma$-algebra.
\end{exercise}

(Hint: Use Exercise \ref{E216}).}

{\setlength{\parindent}{0pt}
\begin{exercise} \label{E219}
Let $X=\mathbb{R}$ and let $\mathcal{A}=\left\{\mbox{finite disjoint unions of intervals of the form}\,\,(a,b], \right.$ \\
$\left. (-\infty,b]\,\,\mbox{and}\,\, (a,+\infty) \right\}$. Prove that $\mathcal{A}$ is an algebra but not a $\sigma$-algebra. 
\end{exercise}

(Hint: Use Exercise \ref{E216}).}

\begin{exercise} \label{E220}
Let $X$ be nonempty and let $\{A_{i}\,:\,i \in \mathcal{I}\}$ be a partition of $X$. Prove that the class defined by
$$
\mathfrak{A}:=\left\{ A_{j}:=\cup_{j\in\mathcal{J}}\,A_{j}\,:\, \mathcal{J} \subset \mathcal{I} \,\,\text{ and }\,\, \mathcal{J}\neq \varnothing\right\}\cup\{\varnothing\}
$$
is a $\sigma$-algebra of subsets of $X$.
\end{exercise}

\begin{exercise} \label{E221}
Let $X$ be fixed nonempty and let $\mathcal{A}=\{A_{1},\ldots,A_{N}\}$ be a finite partition of $X$. Give a complete description of $\sigma(\mathcal{A})$.

What happens if the partition is countable? Justify your answer.
\end{exercise}

{\setlength{\parindent}{0pt}
\begin{exercise} \label{E222}
Let $X$ be fixed nonempty and let $A,B \subset X$ be arbitrary. Define the class $\mathcal{C} \subset \mathcal{P}(X)$ by $\mathcal{C}:=\{A,B\}$ and explicitly determine all the elements of $\sigma(\mathcal{C})$.
\end{exercise}

(Hint: In the most general case, the total number of elements of $\sigma(\mathcal{C})$ is 16).}

{\setlength{\parindent}{0pt}
\begin{exercise} \label{E223}
Let $\mathcal{C} \subset \mathcal{P}(X)$ be fixed. Prove that for every $A \in \sigma(\mathcal{C})$, there exists a countable subfamily $\mathcal{C}_{0} \subset \mathcal{C}$ such that $A \in \sigma(\mathcal{C}_{0})$.
\end{exercise}

(Hint: Let $\mathsf{S}=\cup \{\sigma(\mathcal{C}')\,:\,\mathcal{C}' \subset \mathcal{C}\,\,\mbox{is countable} \}$. Prove that $\mathsf{S}$ is a $\sigma$-algebra and that $\mathcal{C}\subset \mathsf{S}=\sigma(\mathcal{C})$.)}

{\setlength{\parindent}{0pt}
\begin{exercise} \label{E224}
Let $X$ be uncountable and $\mathcal{C}=\{A \subset X\,:\,A\,\,\mbox{is finite} \}$. Prove that
\begin{itemize}
\item[(a)] $\mathcal{R}(\mathcal{C})=\{C \subset X\,:\,C\,\,\mbox{is finite}\}$.
\item[(b)] $\mathcal{A}(\mathcal{C})=\{A \subset X\,:\,A\,\,\mbox{or}\,\,X\smallsetminus A\,\,\mbox{is finite}\}$.
\item[(c)] $\mathcal{SR}(\mathcal{C})=\{A \subset X\,:\,A\,\,\mbox{is finite or countable}\}$.
\item[(d)] $\sigma(\mathcal{C})=\{A \subset X\,:\,A\,\,\mbox{or}\,\,X\smallsetminus A \mbox{ is finite or countable}\}$
\end{itemize}
\end{exercise}

(Hint: Use Exercise \ref{E23}).}

\begin{exercise} \label{E225}
Let $X$ be arbitrary nonempty, $\mathsf{S} \subset \mathcal{P}(X)$ a $\sigma$-algebra and $F \in \mathcal{P}(X)\smallsetminus \mathsf{S}$. Prove that
$$
\sigma(\mathsf{S}\cup \{F\})=\left\{ (A \cap F) \cup (B\smallsetminus F)\,:\,A,B \in \mathsf{S} \right\}.
$$
\end{exercise}

\begin{exercise} \label{E226}
Let $f:X \to Y$ be a function and $\mathcal{B} \subset \mathcal{P}(Y)$ given. Prove that $\mathcal{A}(f^{-1}(\mathcal{B}))=f^{-1}(\mathcal{A}(\mathcal{B}))$ where $\mathcal{A}(\mathcal{B})$ is the algebra of subsets of $Y$ generated by the class $\mathcal{B}$ and $\mathcal{A}(f^{-1}(\mathcal{B}))$ is the algebra of subsets of $X$ generated by the class $f^{-1}(\mathcal{B})$.
\end{exercise}

\begin{exercise} \label{E227}
Prove that the following classes of subsets of $\mathbb{R}$ generate the Borel $\sigma$-algebra
\begin{itemize}
\item[(1)] $\mathcal{C}=\{ [a,b]\,:\, a< b \,\,\,(a,b \in \mathbb{R} ) \}$.
\item[(2)] $\mathcal{C}=\{ (a,b]\,:\, a< b \,\,\,(a,b \in \mathbb{R} ) \}$.
\item[(3)] $\mathcal{C}=\{ [a,b)\,:\, a< b \,\,\,(a,b \in \mathbb{R} ) \}$.
\item[(4)] $\mathcal{C}=\{ (a,+\infty)\,:\, a \in \mathbb{R}  \}$.
\item[(5)] $\mathcal{C}=\{ (-\infty,b]\,:\, b \in \mathbb{R}  \}$.
\end{itemize}
\end{exercise}

\begin{exercise} \label{E228}
Prove the following inclusions:
$$
\begin{array}{c}
\mathcal{G}_{\delta} \subset \mathcal{G}_{\sigma\delta} \subset \mathcal{G}_{\delta\sigma\delta} \subset \mathcal{G}_{\sigma\delta\sigma\delta} \subset \cdots\\
\mathcal{F}_{\sigma} \subset \mathcal{F}_{\delta\sigma} \subset \mathcal{F}_{\sigma\delta\sigma} \subset \mathcal{F}_{\delta\sigma\delta\sigma} \subset \cdots \\
\end{array}
$$
\end{exercise}

\begin{exercise} \label{E229}
Prove that the following classes of subsets of $\mathbb{R}^{N}$ generate the Borel $\sigma$-algebra $\mathcal{B}(\mathbb{R}^{N})$.
\begin{itemize}
\item[(1)] $\mathcal{C}=\left\{ (a_{1},b_{1}) \times (a_{2},b_{2}) \times \cdots \times (a_{N},b_{N})\,:\, a_{i},b_{i} \in \mathbb{R}\,\,\text{with}\,\,a_{i}< b_{i}\,\,\forall i =1,\ldots,N \right\}$. 
\item[(2)] $\mathcal{C}=\{ [a_{1},b_{1}] \times [a_{2},b_{2}]\times \cdots \times[a_{N},b_{N}]\,:\,  a_{i},b_{i} \in \mathbb{R}\,\,\text{with}\,\,a_{i}< b_{i}\,\,\forall i =1,\ldots,N \}$.
\item[(3)] $\mathcal{C}=\{ (-\infty,a_{1}) \times (-\infty,a_{2}) \times \cdots \times (-\infty,a_{N})\,:\, a_{i} \in \mathbb{R} \,\,\,\forall i =1,\ldots,N\}$.
\item[(4)] $\mathcal{C}=\{ (a_{1},+\infty) \times (a_{2},+\infty)\times \cdots \times (a_{N},+\infty)\,:\,  a_{i} \in \mathbb{R} \,\,\,\forall i =1,\ldots,N  \}$.
\end{itemize}
\end{exercise}

\begin{exercise} \label{E230}
Let $(A_k)$ be a sequence of subsets of an arbitrary set $X$. Prove the following:
\begin{itemize}
\item[(a)] $x \in \displaystyle\liminf_{k \to \infty} A_k  \quad \Leftrightarrow\quad  \mbox{there exists}\,\,k=k(x)\in \mathbb{N}\,\,\mbox{such that}\,\, x \in A_{j} \,\,\mbox{for every}\,\,j \geq k.$

\item[(b)] $x \in \displaystyle\limsup_{k \to \infty} A_k\quad \Leftrightarrow\quad x \in A_k\,\,\mbox{for infinitely many }\,k\mbox{'s}.$

\item[(c)] $\displaystyle\liminf_{k \to \infty} A_k \subset \lim\sup_{k \to \infty} A_k.$
\end{itemize}
\end{exercise}

Let $X$ be a nonempty set and $(f_k)$ a sequence of functions $f_k:X \to \mathbb{R}$ such that, for each $x \in X$, the sequence $(f_k(x))$ is bounded in $\mathbb{R}$. The lower limit and upper limit functions of the sequence $(f_k)$ are then defined as the functions $\liminf_{k \to \infty}f_k,\,\limsup_{k \to \infty}f_k:X \to \mathbb{R}$ given by
$$
\begin{aligned}
\left(\liminf_{k \to \infty}f_k\right)(x)&:=\liminf_{k \to \infty}f_k(x):=\sup_{k \geq 1}\inf_{j \geq k}\big(f_j(x) \big),\\
\left(\limsup_{k \to \infty}f_k\right)(x)&:=\limsup_{k \to \infty}f_k(x):=\inf_{k \geq 1}\sup_{j \geq k}\big(f_j(x) \big) .
\end{aligned}
$$

\begin{exercise} \label{E231}
Let $(A_k)$ be a sequence of subsets of an arbitrary set $X$. Prove the following:
\begin{itemize}
\item[(a)] $\,\, \displaystyle\lim_{k \to \infty} A_k = A \quad\Leftrightarrow\quad \lim_{k \to \infty} (X\smallsetminus A_k) = X\smallsetminus A$.
\item[(b)] $\,\, \displaystyle\liminf_{k \to \infty} \chi_{A_k} = \chi_{\liminf (A_k)} \quad \mbox{ and }\quad \limsup_{k \to \infty} \chi_{A_{k}} = \chi_{\limsup (A_k)}$ on $X$.
\item[(c)] $\,\, \displaystyle\lim_{k \to \infty} A_k = A \quad\Leftrightarrow\quad \lim_{k \to \infty} \chi_{A_{k}} = \chi_{A}$ pointwise on $X$.
\end{itemize}
\end{exercise}

\begin{exercise} \label{E232}
Let $(a_k)$ be a sequence of nonnegative real numbers converging to the number $a \geq 0$ in $\mathbb{R}$. Let $A_k = [0,a_k]$ for every $k \in \mathbb{N}$. Compute $\displaystyle\liminf_{k \to\infty} A_k$ and $\displaystyle\limsup_{k \to \infty} A_k$.
\end{exercise}

\begin{exercise} \label{E233}
Let $(A_k)$ and $(B_k)$ be two sequences of subsets of $X$. Prove the following:
\begin{itemize}
\item[(a)] $\displaystyle\limsup_{k \to \infty} (A_k \cap B_k) \subset \limsup_{k \to \infty} A_k \cap \limsup_{k \to \infty} B_k$. Give an example in which the inclusion is proper and prove that equality holds if $\cap$ is replaced by $\cup$.
\item[(b)] $\displaystyle\liminf_{k \to \infty} A_k \cup \liminf_{k \to \infty} B_k \subset\displaystyle\liminf_{k \to \infty} (A_k \cup B_k)$. Give an example in which the inclusion is proper and prove that equality holds if $\cup$ is replaced by $\cap$.
\item[(c)] $\displaystyle\limsup_{k \to \infty} A_k \bigtriangleup \limsup_{k \to \infty} B_k \subset \limsup_{k \to \infty} (A_k \bigtriangleup B_k)$.
\end{itemize}
\end{exercise}

{\setlength{\parindent}{0pt}
\begin{exercise} \label{E234}
Let $(A_k)$ be a sequence of subsets of $X$. Define $D_1:=A_1$, $D_2:=A_1 \bigtriangleup A_2$ and, in general, $D_{k+1}:=D_k \bigtriangleup A_{k+1}$ for $k=1,2,\ldots$. Prove that
$$
\liminf_{k \to \infty} D_k = \limsup_{k \to \infty} D_k \quad\Leftrightarrow\quad \liminf_{k \to \infty} A_k = \limsup_{k \to \infty} A_k = \varnothing.
$$
\end{exercise}

(Hint: Use Exercise \ref{E231}).}

\begin{exercise} \label{E235}
Determine whether the following sequences of subsets $(A_k)$ of a set $X$ converge.
\begin{itemize}
\item[(a)] $A_k=\varnothing$ if $k$ is odd and $A_k=X$ if $k$ is even.
\item[(b)] If $X=\mathbb{R}$ and $A_{k}:=\left(0 ,1 + (-1)^k \right)$.
\end{itemize}
\end{exercise}

\begin{exercise} \label{E236}
Let $(A_k)$ be a sequence of subsets of a set $X$ and let $A$ be a subset of $X$ such that $\displaystyle\lim_{k \to \infty} A_k =A$. Prove that, for every $B \subset X$,
\begin{itemize}
\item[(a)] $\displaystyle\lim_{k \to \infty} (A_k \cap B)=A\cap B$.
\item[(b)] $\displaystyle\lim_{k \to \infty} (A_k \cup B)=A \cup B$.
\item[(c)] $\displaystyle\lim_{k \to \infty} (A_k \smallsetminus B)=A \smallsetminus B$.
\item[(d)] $\displaystyle\lim_{k \to \infty} (A_k \bigtriangleup B)=A \bigtriangleup B$.
\end{itemize}
\end{exercise}

\begin{exercise} \label{E237}
Let $X=\mathbb{R}$. Prove that the following classes
$\mathcal{P}_{1}:=\{ (-\infty,a]\,:\, a \in \mathbb{R} \}$ and
$\mathcal{P}_{2}:=\{ (a,b]\,:\,a,b\in\mathbb{R}\,\,\mbox{with}\,\,a<b \} \cup \{ \varnothing\}$
are $\pi$-systems.
\end{exercise}

\begin{exercise} \label{E238}
Prove that a class $\mathcal{L}\subset \mathcal{P}(X)$ is a $\lambda$-system if and only if it satisfies the following properties:
\begin{itemize}
\item[(1)] $X \in \mathcal{L}$.
\item[(2)] $X\smallsetminus A \in \mathcal{L}$ for every $A \in \mathcal{L}$.
\item[(3)] For every disjoint sequence $(A_k)$ of elements of $\mathcal{L}$,
$\displaystyle\bigcup_{k=1}^{\infty} A_k \in \mathcal{L}$.
\end{itemize}
\end{exercise}

\begin{exercise} \label{E239}
Prove that the intersection $\bigcap_{i\,\in\,\mathcal{I}}\, \mathcal{C}_{i}$ of any family $\{\mathcal{C}_{i}\,:\,i\in\mathcal{I} \}$ of $\pi$-systems is a $\pi$-system.
\end{exercise}

\begin{exercise} \label{E240}
Prove that the intersection $\bigcap_{\tau\,\in\,\mathcal{T}}\, \mathcal{L}_{\tau}$ of any family $\{\mathcal{L}_{\tau}\,:\,\tau\in\mathcal{T} \}$ of $\lambda$-systems is a $\lambda$-system.
\end{exercise}

\section{Project: The Monotone Class Lemma}

\subsection*{2.7.1.\quad Objective}

Let $X$ be a fixed nonempty set.

\begin{definition} \label{261}
A nonempty collection $\mathcal{K} \subset \mathcal{P}(X)$ is called a \textbf{monotone class} if for every increasing sequence $(A_k)$ in $\mathcal{K}$ and every decreasing sequence $(B_k)$ in $\mathcal{K}$, the sets $\bigcup_{k=1}^{\infty}\,A_k$ and $\bigcap_{k=1}^{\infty}\,B_k$ belong to $\mathcal{K}$.
\end{definition}

For example, every $\sigma$-ring of subsets of $X$ is a monotone class.

The objective of this project is to prove the well-known monotone class lemma, a result essentially similar to Dynkin’s theorem (see Theorem \ref{255}).

\begin{theorem}[Monotone Class Lemma] \label{262}
Let $\mathcal{A} \subset \mathcal{P}(X)$ be an algebra and let $\mathcal{K}_{0} \subset \mathcal{P}(X)$ be a monotone class such that $\mathcal{A} \subset \mathcal{K}_{0}$. Then $\sigma(\mathcal{A}) \subset \mathcal{K}_{0}$.
\end{theorem}

\subsection*{2.7.2.\quad Procedure}

\begin{itemize}
    \item[1.] Prove that every $\sigma$-algebra of subsets of $X$ is a monotone class of subsets of $X$, but that not every monotone class of subsets of $X$ is necessarily a $\sigma$-algebra of subsets of $X$.

    \item[2.] We say that a nonempty class $\mathcal{KA}$ of subsets of $X$ is a monotone algebra if it is both a monotone class and an algebra.
    
    Prove that every monotone algebra $\mathcal{KA}$ of subsets of $X$ is a $\sigma$-algebra.

    \item[3.] Let $\{\mathcal{K}_{i}\,:\,i\in\mathcal{I} \}$ be an arbitrary family of monotone classes of subsets of $X$.
    
    Prove that the intersection $\bigcap_{i\,\in\,\mathcal{I}}\, \mathcal{K}_{i}$ is a monotone class of subsets of $X$.
 
    \item[4.] Let $\mathcal{C}$ be an arbitrary nonempty class of subsets of $X$. Prove the following statements:
    \begin{itemize}
        \item[a)] There exists a unique monotone class of subsets of $X$, denoted by $\mathcal{K}(\mathcal{C})$, containing the class $\mathcal{C}$ and contained in every monotone class containing $\mathcal{C}$.
        \item[b)] The monotone class generated by $\mathcal{C}$ is contained in the $\sigma$-algebra generated by $\mathcal{C}$.
    \end{itemize}

    \item[5.] Let $\mathcal{A}$ be an algebra of subsets of $X$ and define
    $$\mathcal{K}_{c}:=\left\{A \in \mathcal{K}(\mathcal{A})\,:\,X\smallsetminus A \in \mathcal{K}(\mathcal{A}) \right\}.$$
    Prove that $\mathcal{K}_{c}$ is a monotone class such that $\mathcal{A} \subset \mathcal{K}_{c}$ and therefore that $\mathcal{K}(\mathcal{A})=\mathcal{K}_{c}$.

    \item[6.] Let $\mathcal{A}$ be an algebra of subsets of $X$ and define
    $$\mathcal{K}_{1}:=\left\{A \in \mathcal{K}(\mathcal{A})\,:\,A \cap B \in \mathcal{K}(\mathcal{A})\,\,\mbox{ for every}\,\,B \in \mathcal{A} \right\}.$$
    Prove that $\mathcal{K}_{1}$ is a monotone class such that $\mathcal{A} \subset \mathcal{K}_{1}$ and therefore that $\mathcal{K}(\mathcal{A})=\mathcal{K}_{1}$.

    \item[7.] Let $\mathcal{A}$ be an algebra of subsets of $X$ and define
    $$\mathcal{K}_{2}:=\left\{A \in \mathcal{K}(\mathcal{A})\,:\,A \cap B \in \mathcal{K}(\mathcal{A})\,\,\mbox{ for every}\,\,B \in \mathcal{K}(\mathcal{A}) \right\}.$$
    Prove that $\mathcal{K}_{2}$ is a monotone class such that $\mathcal{A} \subset \mathcal{K}_{2}$ and therefore that $\mathcal{K}(\mathcal{A})=\mathcal{K}_{2}$.

    \item[8.] Let $\mathcal{A}$ be an algebra of subsets of $X$. Prove that $\mathcal{K}(\mathcal{A})$ is an algebra of subsets of $X$.

    \item[9.] Prove that if $\mathcal{A}$ is an algebra and $\mathcal{K}_{0}$ is a monotone class of subsets of $X$ such that $\mathcal{A} \subset \mathcal{K}_{0}$, then $\sigma(\mathcal{A}) \subset \mathcal{K}_{0}$.

    \item[10.] Finally, prove that if $\mathcal{A}$ is an algebra of subsets of $X$, then $\sigma(\mathcal{A}) = \mathcal{K}(\mathcal{A})$.
\end{itemize}

\subsection*{2.7.3.\quad Remarks}

Let $\mathcal{A} \subset \mathcal{P}(X)$ be a class such that every $A \in \mathcal{A}$ satisfies a certain property $\wp$. Analogously to Theorem \ref{255}, if $\mathcal{A}$ is an algebra, then in order to verify that all elements of $\sigma(\mathcal{A})$ also satisfy the property $\wp$, it suffices to prove that the class
$$
\mathcal{K}=\{A \subset X : A \text{ satisfies the property } \wp\}
$$
is a monotone class.

The relationship between Theorem \ref{255} and Theorem \ref{262} lies in the fact that every algebra is, in particular, a $\pi$-system, and every $\lambda$-system is a monotone class. These results play a fundamental role, for example, in the proof of uniqueness of measures and in the development of the theory of product measures.
\chapter{Measurable functions with respect to a \texorpdfstring{$\sigma$}{}-algebra} \label{Capitulo3}
\markboth{{\scriptsize 3. MEASURABLE FUNCTIONS WITH RESPECT TO A SIGMA ALGEBRA}}{ {\scriptsize 3. MEASURABLE FUNCTIONS WITH RESPECT TO A SIGMA ALGEBRA}}

Given a function $f : X \to Y$ between two sets endowed with a $\sigma$-algebra, we are interested in studying the properties it preserves in the structure of measurable sets. We say that a function is measurable if the inverse image of every measurable set is also measurable. In this chapter we will work mainly with functions $f : X \to \mathbb{R}$, considering on $\mathbb{R}$ the Borel $\sigma$-algebra and assuming that $X$ is a measurable space endowed with a $\sigma$-algebra $\mathsf{S}$.

We will define the set $\mathbb{M}(X, \mathsf{S})$ of real-valued measurable functions with respect to a $\sigma$-algebra $\mathsf{S}$ of subsets of $X$. This set forms a vector space under the usual operations of addition and scalar multiplication, and is moreover closed under other countable operations. We will prove that every real-valued measurable function can be obtained as the pointwise limit of a sequence of simple functions, that is, measurable functions taking only finitely many values on elements of $\mathsf{S}$.

It will be useful to introduce a ``larger'' set than $\mathbb{R}$ in which we may conveniently work with all kinds of limits and which, as we shall see in the next chapter, will serve as a foundation for the study of measures. This set is the extended real line $\overline{\mathbb{R}}$. Finally, we will generalize the concept of measurable function to functions taking values in the complex numbers and establish conditions that allow their measurability to be characterized more simply, making use of the results developed in the previous sections.

\section{Definition and examples}

\begin{definition} \label{31}
A \textbf{measurable space} is a pair $(X,\mathsf{S})$ in which $X$ is a nonempty set and $\mathsf{S}$ is a $\sigma$-algebra of subsets of $X$.
\end{definition}

\begin{definition} \label{32}
Let $(X,\mathsf{S})$ and $(Y, \mathsf{T})$ be measurable spaces. A function $f:X \to Y$ is said to be \textbf{measurable with respect to $\mathsf{S}$ and $\mathsf{T}$} if $f^{-1}(B) \in \mathsf{S}$ for every $B \in \mathsf{T}$, that is, if $f^{-1}(\mathsf{T}) \subset \mathsf{S}$.
\end{definition}

The measurability of a function $f$ depends on the $\sigma$-algebras being considered on $X$ and $Y$. To emphasize this, we will sometimes use the notation
$$
f:(X,\mathsf{S}) \to (Y,\mathsf{T})
$$
instead of $f:X \to Y$.

In the case where $Y=\mathbb{R}$ and $\mathsf{T}=\mathcal{B}(\mathbb{R})$, we will say that $f:X \to \mathbb{R}$ is \textbf{$\mathsf{S}$-measurable} instead of measurable with respect to $\mathsf{S}$ and $\mathcal{B}(\mathbb{R})$. Likewise, if $X=Y=\mathbb{R}$ and $\mathsf{S}=\mathsf{T}=\mathcal{B}(\mathbb{R})$, we will simply say that $f$ is \textbf{Borel measurable}.

We denote by $\mathbb{M}(X,\mathsf{S})$ the set of real-valued functions $f:X \to \mathbb{R}$ that are $\mathsf{S}$-measurable, and by $\mathbb{M}^{+}(X,\mathsf{S})$ the set of nonnegative real-valued functions that are $\mathsf{S}$-measurable.

The following result establishes a sufficient condition for the measurability of a function between measurable spaces, formulated in terms of a generating class.
 
\begin{theorem} \label{33}
Let $(X,\mathsf{S})$ and $(Y,\mathsf{T})$ be measurable spaces and let $f:X \to Y$ be a function. If there exists $\mathcal{C} \subset \mathcal{P}(Y)$ such that $\sigma(\mathcal{C})=\mathsf{T}$ and $f^{-1}(\mathcal{C}) \subset \mathsf{S}$, then $f$ is measurable with respect to $\mathsf{S}$ and $\mathsf{T}$.
\end{theorem}

\begin{proof}
Proposition \ref{225} ensures that $f^{-1}(\sigma(\mathcal{C}))=\sigma(f^{-1}(\mathcal{C}))$ and, since $\sigma(\mathcal{C})=\mathsf{T}$, it follows that $f^{-1}(\mathsf{T})=\sigma(f^{-1}(\mathcal{C}))$. From the inclusion $f^{-1}(\mathcal{C}) \subset \mathsf{S}$ we obtain that $\sigma(f^{-1}(\mathcal{C}))\subset  \mathsf{S}$ and therefore that $f^{-1}(\mathsf{T})\subset \mathsf{S}$.
\end{proof}

In the particular case of $\mathsf{S}$-measurable functions, we obtain the following result.

\begin{proposition} \label{34}
\it Let $(X,\mathsf{S})$ be a measurable space and let $f:X \to \mathbb{R}$ be a function. The following statements are equivalent:
\begin{itemize}
\item[(a)] $f$ is $\mathsf{S}$-measurable.
\item[(b)] $f^{-1}((c,+\infty))$ belongs to $\mathsf{S}$ for every $c \in \mathbb{R}$.
\item[(c)] $f^{-1}((-\infty,c])$ belongs to $\mathsf{S}$ for every $c \in \mathbb{R}$.
\item[(d)] $f^{-1}((-\infty,c))$ belongs to $\mathsf{S}$ for every $c \in \mathbb{R}$.
\item[(e)] $f^{-1}([c,+\infty))$ belongs to $\mathsf{S}$ for every $c \in \mathbb{R}$.
\end{itemize}
\end{proposition}

\begin{proof}
$(a)\Rightarrow (b): $ Since $f$ is $\mathsf{S}$-measurable, we have $f^{-1}(\mathcal{B}(\mathbb{R})) \subset \mathsf{S}$. Moreover, because $(c,+\infty) \in \mathcal{B}(\mathbb{R})$ for every $c \in \mathbb{R}$, it follows that $f^{-1}((c,+\infty)) \in \mathsf{S}$ for every $c \in \mathbb{R}$.

$(b)\Rightarrow (c):$ Since $(-\infty,c]=\mathbb{R}\smallsetminus(c,+\infty)$ for every $c \in \mathbb{R}$, we have
$$
f^{-1}((-\infty,c])=f^{-1}\left( \mathbb{R}\smallsetminus(c,+\infty) \right)=f^{-1}(\mathbb{R})\smallsetminus f^{-1}((c,+\infty)) =X\smallsetminus f^{-1}((c,+\infty)) \in \mathsf{S}
$$
for every $c \in \mathbb{R}$.

$(c) \Rightarrow (d):$ For every $c \in \mathbb{R}$, the interval $(-\infty,c)$ can be written as
$$
(-\infty,c)=\bigcup_{k=1}^{\infty} \left(-\infty, c -\frac{1}{k} \right]
$$
and therefore
$$
f^{-1}((-\infty,c))=f^{-1}\left( \bigcup_{k=1}^{\infty} \left(-\infty, c -\frac{1}{k} \right] \right) =\bigcup_{k=1}^{\infty} f^{-1}\left(\left(-\infty, c -\frac{1}{k} \right]\right) \in \mathsf{S}
$$
for every $c \in \mathbb{R}$.

$(d) \Rightarrow (e):$ For every $c \in \mathbb{R}$ we have
$$
[c,+\infty)=\mathbb{R}\smallsetminus(-\infty,c).
$$

Therefore
$$
f^{-1}([c,+\infty))=f^{-1}(\mathbb{R}\smallsetminus(-\infty,c))=f^{-1}(\mathbb{R})\smallsetminus f^{-1}((-\infty,c))=X\smallsetminus f^{-1}((-\infty,c)) \in \mathsf{S}
$$
for every $c \in \mathbb{R}$.

$(e)\Rightarrow (a):$ Consider the class $\mathcal{C}=\left\{[c,+\infty)\,:\, c\in \mathbb{R} \right\}$ of subsets of $\mathbb{R}$. Then $\sigma(\mathcal{C})=\mathcal{B}(\mathbb{R})$ [Exercise \ref{E227}] and, by hypothesis, $f^{-1}(\mathcal{C}) \subset \mathsf{S}$. Theorem \ref{33} ensures that $f$ is $\mathsf{S}$-measurable.
\end{proof}

Let us now look at some examples.

\begin{example} \label{35}
Let $(X,\mathsf{S})$ be a measurable space. The identity function
$id:(X,\mathsf{S}) \to (X,\mathsf{S})$, defined by $id(x)=x$, is measurable with respect to $\mathsf{S}$.
\end{example}

\begin{proof}
This is immediate since $id^{-1}(B)=B$ for every $B \subset X$.
\end{proof}

In general, the identity function need not be measurable when two different $\sigma$-algebras $\mathsf{S}$ and $\mathsf{T}$ are considered on the same set $X$. Indeed, the mapping
$$
id:(X,\mathsf{S}) \to (X,\mathsf{T})
$$
is measurable with respect to $\mathsf{S}$ and $\mathsf{T}$ if and only if $\mathsf{T}\subset \mathsf{S}$. The proof of this statement is a simple exercise [Exercise \ref{E36}].

\begin{example} \label{36}
Let $(X,\mathsf{S})$ and $(Y,\mathsf{T})$ be measurable spaces and let $y_{0} \in Y$. The constant function $f:X \to Y$ equal to $y_{0}$ is measurable with respect to $\mathsf{S}$ and $\mathsf{T}$.
\end{example}

\begin{proof}
For every $B \subset Y$ with $B \in \mathsf{T}$, we have
$$
f^{-1}(B)=\left\{
\begin{array}{cl}
X & \mbox{if}\,\,\,y_{0} \in B,\\
\varnothing & \mbox{if}\,\,\, y_{0} \notin B.
\end{array}
\right.
$$

In both cases, $f^{-1}(B) \in \mathsf{S}$ and therefore $f^{-1}(\mathsf{T}) \subset \mathsf{S}$.
\end{proof}

\begin{example} \label{37}
An upper semicontinuous (\textbf{u.s.c.}) function $f:\mathbb{R} \to \mathbb{R}$ is Borel measurable.
\end{example}

\begin{proof}
Since $f$ is u.s.c., for every $c \in \mathbb{R}$ the set
$\{x \in \mathbb{R}\,:\, f(x) < c \}=f^{-1}((-\infty,c))$
is open in $\mathbb{R}$ [Exercise \ref{E311}]. Proposition \ref{34} ensures that $f$ is Borel measurable.
\end{proof}

\begin{example} \label{38}
A lower semicontinuous (\textbf{l.s.c.}) function $f:\mathbb{R} \to \mathbb{R}$ is Borel measurable.
\end{example}

\begin{proof}
Since $f$ is l.s.c., for every $c \in \mathbb{R}$ the set
$\{x \in \mathbb{R}\,:\, f(x) > c \}=f^{-1}((c,+\infty))$
is open in $\mathbb{R}$ [Exercise \ref{E311}]. Proposition \ref{34} ensures that $f$ is Borel measurable.
\end{proof}

\begin{example} \label{39}
A continuous function $f:\mathbb{R} \to \mathbb{R}$ on $\mathbb{R}$ is Borel measurable.
\end{example}

\begin{proof}
Since every continuous function $f$ is both u.s.c. and l.s.c., Examples \ref{37} and \ref{38} imply that $f$ is Borel measurable.
\end{proof}

A generalization of the previous example may be found in Exercise \ref{E311}.

\begin{proposition} \label{310}
Let $(X,\mathsf{S})$ be a measurable space, let $f:X\to\mathbb{R}$ be an $\mathsf{S}$-measurable function, and let $\varphi:\mathbb{R}\to\mathbb{R}$ be a Borel measurable function. Then $\varphi \circ f : X \to \mathbb{R}$ is $\mathsf{S}$-measurable.
\end{proposition}

\begin{proof}
From the following equality
$$
(\varphi \circ f)^{-1}(\mathcal{B}(\mathbb{R}))=f^{-1}(\varphi^{-1}(\mathcal{B}(\mathbb{R}))),
$$
the result follows.
\end{proof}

\begin{example} \label{311}
Let $(\mathbb{R}^2, \mathcal{B}(\mathbb{R}^2))$ and $(\mathbb{R},\mathcal{B}(\mathbb{R}))$ be measurable spaces. Then
$f=(f_1,f_2):\mathbb{R}^2 \to \mathbb{R}^2$
is $\mathcal{B}(\mathbb{R}^2)$-measurable if and only if
$f_i:\mathbb{R}^2 \to \mathbb{R}$
is $\mathcal{B}(\mathbb{R}^2)$-measurable for each $i=1,2$.
\end{example}

\begin{figure}[ht!]
\centering
\begin{tikzpicture}[xscale=0.5,yscale=0.5]

\draw[fill=gray!21, dotted] (-4,-1.5)--(4,-1.5)--(4,1.5)--(-4,1.5)--(-4,1.5);
\draw[fill=gray!21, dotted] (1.5,-4)--(1.5,4)--(-1.5,4)--(-1.5,-4)--(1.5,-4);

\draw[fill=red!15, dotted] (-1.5,-1.5)--(-1.5,1.5)--(1.5,1.5)--(1.5,-1.5)--(-1.5,-1.5);
\draw[->, gray] (-4,0)--(4,0); 
\draw[->, gray] (0,-4)--(0,4); 


\draw (-1.5,0) node{$_{|}$};  \draw (-1.5,0) node[below]{$_{a_{1}}$};
\draw (1.5,0) node{$_{|}$};  \draw (1.5,0) node[below]{$_{b_{1}}$};

\draw (0,-1.5) node{$_{-}$}; \draw (0,-1.5) node[left]{$_{a_{2}}$};
\draw (0,1.5) node{$_{-}$}; \draw (0,1.5) node[left]{$_{b_{2}}$};
\end{tikzpicture}
\begin{center}
$(a_1,b_1)\times (a_{2},b_{2})$
\end{center}
\end{figure}

\begin{proof}
Suppose that $f$ is $\mathcal{B}(\mathbb{R}^2)$-measurable and let $\pi_{i} : \mathbb{R}^{2} \to \mathbb{R}$ denote the $i$-th projection of $\mathbb{R}^{2}$ onto $\mathbb{R}$ for $i=1,2$. Since $\pi_i$ is continuous on $\mathbb{R}^{2}$, it is therefore $\mathcal{B}(\mathbb{R}^2)$-measurable [Exercise \ref{E311}]. Thus, because $f_{i}= \pi_{i} \circ f$ for $i=1,2$, Proposition \ref{310} implies that $f_{i}$ is $\mathcal{B}(\mathbb{R}^2)$-measurable.

Conversely, suppose that $f_{i}$ is $\mathcal{B}(\mathbb{R}^2)$-measurable for $i=1,2$, and let
$R=(a_1,b_1) \times (a_2,b_2)$
be an open rectangle in $\mathbb{R}^{2}$. Then
$$
\begin{aligned}
f^{-1}(R)&=\{ (x,y) \in \mathbb{R}^{2}\,:\, f(x,y) \in R  \} \\
&= \{(x,y) \in \mathbb{R}^{2} \,:\, a_1 < f_1(x,y) < b_1 \} \cap \{(x,y) \in \mathbb{R}^{2} \,:\, a_2 < f_2(x,y) < b_2 \}  \\
&=f^{-1}_{1}(a_{1},b_{1}) \cap f_{2}^{-1}(a_{2},b_{2}) \in \mathcal{B}(\mathbb{R}^{2}).
\end{aligned}
$$

Therefore, if
$$
\mathcal{C}=\{(a,b)\times (c,d)\,:\,a,b,c,d\in \mathbb{R}\,\,\text{with}\,\,a<b \,\,\text{and}\,\,c<d\},
$$
then $f^{-1}(\mathcal{C}) \subset \mathcal{B}(\mathbb{R}^{2})$. Theorem \ref{33} ensures that $f$ is $\mathcal{B}(\mathbb{R}^2)$-measurable since $\sigma(\mathcal{C})=\mathcal{B}(\mathbb{R}^2)$ [Exercise \ref{E229}].
\end{proof}

\begin{definition} \label{312}
Given a subset $A$ of $X$, we define the \textbf{characteristic function} (\textbf{indicator function}) of $A$, denoted by $\chi_{A}:X \to \{0,1\}$, as follows
$$
\chi_{A}(x):=\left\{
\begin{array}{l c l}
1 & & \text{if}\,\,x \in A,\\
0 & & \text{if}\,\, x\not\in A.
\end{array}
\right.
$$
\end{definition}

In what follows, this function will play a fundamental role throughout the text; for this reason, it is convenient to keep in mind its properties, which are presented in [Exercise \ref{E33}].

The following proposition provides a characterization of its measurability.

\begin{proposition} \label{313}
Let $(X,\mathsf{S})$ be a measurable space. The function $\chi_{A}$ is $\mathsf{S}$-measurable if and only if $A \in \mathsf{S}$.
\end{proposition}

\begin{figure}[ht!]
    \centering
	\begin{tikzpicture}[xscale=0.8,yscale=0.8]
	\draw [fill=lime!15,  dotted]  (2.3,0)--(2.3,2.8)--(1,2.8)--(1,0);
	\draw[->, gray] (0,0) -- (4,0); \draw [->, gray] (-1,0) -- (-1,4); 
	\draw[-] (-1,0)--(0,0); 
\draw (-1,2.8) node{$-$}; \draw (-1,2.8) node[left] {$1$};
\draw (-1,0) node{$-$}; \draw (-1,0) node[left] {$0$};
\draw (1.6,-0.05) node[below] {$A$};
\draw [ultra thick] (1,2.8)--(2.3,2.8);
\draw [ultra thick] (1,0)--(-1,0); \draw [ultra thick] (2.3,0)--(3.8,0);
\draw (1,2.8) node{$_\bullet$}; \draw (2.3,2.8) node{$_\bullet$};
\end{tikzpicture}
\begin{center}
$\chi_{A}(x)$
\end{center}
\end{figure}

\begin{proof}
$\Rightarrow):$ Since $A = \chi_{A}^{-1}((0,+\infty)) = \{ x \in X\,:\, \chi_{A}(x) > 0 \}$ and $\chi_A$ is $\mathsf{S}$-measurable, then $A \in \mathsf{S}$ by Proposition \ref{34}. 

$\Leftarrow):$ For $c \in \mathbb{R}$ we have
$$
\chi_{A}^{-1}((c,+\infty))=\{x \in X \,:\, \chi_{A}(x)>c \}=\left\{
\begin{array}{c c l}
\varnothing & & \mbox{if}\,\, c \geq 1,\\
A & & \mbox{if}\,\, 0 \leq c <1,\\
X & & \mbox{if}\,\, c<0.
\end{array}
\right.
$$

Consequently, $\chi_{A}^{-1}((c,+\infty)) \in \mathsf{S}$ for every $c \in \mathbb{R}$. Applying Proposition \ref{34}, we conclude that $\chi_{A}$ is $\mathsf{S}$-measurable.
\end{proof}

Finite linear combinations of characteristic functions play a central role in measure and integration theory. In particular, these functions allow more general classes of functions to be approximated through limiting procedures, so it is natural to introduce and analyze them in detail.

\begin{definition} \label{314}
Let $(X,\mathsf{S})$ be a measurable space. We say that a function $s:X \to \mathbb{R}$ is \textbf{$\mathsf{S}$-simple} if $s$ takes only finitely many values in $\mathbb{R}$ and is $\mathsf{S}$-measurable.
\end{definition}

\begin{proposition} \label{315}
Let $(X,\mathsf{S})$ be a measurable space. A function $s:X \to \mathbb{R}$ is $\mathsf{S}$-simple if and only if $s$ is a linear combination of characteristic functions of pairwise disjoint sets in $\mathsf{S}$. That is, there exist distinct $\alpha_{1},\ldots,\alpha_{N} \in \mathbb{R}$ and pairwise disjoint sets $A_{1},\ldots,A_{N} \in \mathsf{S}$ whose union is equal to $X$ such that
$$
s=\sum_{j=1}^{N} \alpha_{j}\chi_{A_{j}}.
$$

The above expression is known as the \textbf{canonical representation} of $s$.
\end{proposition}

\begin{figure}[ht!]
    \centering
	\begin{tikzpicture}[xscale=0.8,yscale=0.8]
	\draw[->,gray] (0,0) -- (4,0); \draw [->,gray] (0,0) -- (0,4); 
	 \draw[-,gray] (0,0)--(0,-1.2);

\draw (0,2.8) node{$-$}; \draw (0,2.8) node[left] {$_{\alpha_{i}}$};
\draw [ thick] (1.3,2.8)--(2.22,2.8);
\draw (1.3,2.8) node{$\bullet$}; 
\draw (2.3,2.8) node{$\circ$};

\draw (0,2.3) node{$-$}; \draw (0,2.3) node[left] {$_{\alpha_{k}}$};
\draw [ thick] (2.3,2.3)--(2.9,2.3);
\draw (2.3,2.3) node{$\bullet$}; 
\draw (3,2.3) node{$\circ$};

\draw (0,1) node{$-$}; \draw (0,1) node[left] {$_{\alpha_{j}}$};
\draw [ thick] (0,1)--(1.12,1);
\draw (0,1) node{$\bullet$}; 
\draw (1.2,1) node{$\circ$};

\draw (0,-1) node{$-$}; \draw (0,-1) node[left] {$_{\alpha_{\ell}}$};
\draw [ thick] (3,-1)--(4,-1);
\draw (3,-1) node{$\bullet$}; 

\draw[densely dotted, gray] (1.25,2.8)--(1.25,1.1);
\draw[densely dotted, gray] (1.25,0.9)--(1.25,0);

\draw[densely dotted, gray] (2.3,2.7)--(2.3,0); 
\draw[densely dotted, gray] (3,-1)--(3,2.25); 
\end{tikzpicture}
\begin{center}
Example of an $\mathsf{S}$-simple function
\end{center}
\end{figure}

\begin{proof}
$\Rightarrow):$ Since $s$ is $\mathsf{S}$-simple, it takes only finitely many values. Let $\alpha_{1},\ldots,\alpha_{N} \in \mathbb{R}$ be all the distinct values taken by $s$.

Define, for each $j=1,\ldots,N$, the set
$$
A_{j}=s^{-1}(\{\alpha_{j} \})=\{ x \in X\,:\, s(x)=\alpha_{j} \}.
$$

Thus, $A_{1},\ldots,A_{N} \in \mathsf{S}$ since $s$ is $\mathsf{S}$-measurable, $X=\bigcup_{j=1}^{N} A_{j}$, and $A_{i} \cap A_{j} = \varnothing$ for every $i \neq j$. Consequently,
$$
s=\sum_{j=1}^{N} \alpha_{j}\chi_{A_{j}}.
$$

$\Leftarrow):$ If
$$
s=\sum_{j=1}^{N} \alpha_{j}\chi_{A_{j}},
$$
then for every $c \in \mathbb{R}$ we have
$$
\{ x\in X\,:\, s(x) > c\}=\left\{
\begin{array}{c c l}
\varnothing &  & \mbox{if}\quad \alpha_{j} \leq c\quad \forall \,j=1,\ldots,N,\\
&&\\
\bigcup \{ A_{j}\,:\, \alpha_{j} > c \} & & \mbox{otherwise}.
\end{array}
\right.
$$

Proposition \ref{34} ensures that $s$ is $\mathsf{S}$-measurable.
\end{proof}

Algebraic operations between $\mathsf{S}$-simple functions again produce $\mathsf{S}$-simple functions, as established in the following proposition.

\begin{proposition} \label{316}
Let $(X,\mathsf{S})$ be a measurable space. If $s,t:X \to \mathbb{R}$ are $\mathsf{S}$-simple functions, then $s+t$, $\gamma \,s$ with $\gamma \in \mathbb{R}$, $st$, and $\frac{1}{s}$ (if $s \neq 0$) are also $\mathsf{S}$-simple functions.
\end{proposition}

The proof is straightforward and is left as an exercise [Exercise \ref{E35}].

Given a function $f:X \to \mathbb{R}$, we denote by
$$
f^{+}(x):=\max\{f(x),0\}\quad\mbox{and}\quad f^{-}(x):=\max\{-f(x),0\}
$$
the positive part and the negative part of $f$, respectively.

It is immediate to verify that $f^{+}$ and $f^{-}$ are nonnegative functions such that $f=f^{+}-f^{-}$ and $|f|=f^{+}+f^{-}$ [Exercise \ref{E32}].

\begin{figure}[ht!]
\begin{minipage}[r]{0.5\textwidth}
\begin{center}
\begin{tikzpicture}[xscale=0.8, yscale=1]
\draw[->, gray] (0,0) -- (6.5,0);
\draw [->, gray] (0,-1.2) -- (0,1.2);
\draw[ultra thick] plot[smooth] coordinates
{(0.00,0.00)(0.79,0.71)(1.57,1.00)(2.36,0.71)
(3.14,0.00)(3.93,-0.71)(4.71,-1.00)(5.50,-0.71)(6.28,-0.00)};
\draw (2,1.2) node{$f$};
\end{tikzpicture}
\end{center}
\end{minipage} \hfill 
 \begin{minipage}[l]{0.5\textwidth}
\begin{center}
\begin{tikzpicture}[xscale=0.8,yscale=1]
\draw[->,gray] (0,0) -- (7,0);
\draw [->,gray] (0,-1.2) -- (0,1.2);
\draw[ultra thick, red] plot[smooth] coordinates
{(0,0)(0.79,0.71)(1.57,1.00)(2.36,0.71)
(3.14,0.00)};
\draw[-, ultra thick, red] (3.14,0)--(6.28,0);

\draw (2,1.25) node[red]{$f^{+}$};

\draw[ultra thick, blue] plot[smooth] coordinates
{(3.14,0.00)(3.93,0.71)(4.71,1.00)(5.50,0.71)(6.28,0.00)};
\draw[-, ultra thick, blue] (3.14,0)--(0,0);
\draw (4.5,1.25) node[blue]{$f^{-}$};
\end{tikzpicture}
\end{center}
\end{minipage}
\end{figure}

\begin{proposition} \label{317}
Let $(X,\mathsf{S})$ be a measurable space and let $f:X \to \mathbb{R}$ be an $\mathsf{S}$-measurable function. Then both $f^{+}$ and $f^{-}$ are $\mathsf{S}$-measurable functions.
\end{proposition}

\begin{proof}
We will prove only that $f^{+}$ is $\mathsf{S}$-measurable and leave as an exercise the proof that $f^{-}$ is $\mathsf{S}$-measurable [Exercise \ref{E32}]. 

Let $c \in \mathbb{R}$ be arbitrary.

{\scshape Case 1.}\quad If $c\geq 0$, then $(f^{+})^{-1}((c,+\infty))=f^{-1}(c,+\infty) \in \mathsf{S}$ by Proposition \ref{34}, since $f$ is $\mathsf{S}$-measurable.

{\scshape Case 2.}\quad If $c<0$, then $(f^{+})^{-1}((c,+\infty))=X \in \mathsf{S}$ by property (S1).

Therefore, $f^{+}$ is $\mathsf{S}$-measurable.
\end{proof}

\section{Approximation by simple functions}

The following lemmas provide a characterization of measurable functions as those that arise as pointwise limits of simple functions. \index{lemma!approximation by simple functions}

\begin{lemma} \label{318}
Let $(X,\mathsf{S})$ be a measurable space and let $f:X \to \mathbb{R}$ be a nonnegative $\mathsf{S}$-measurable function. Then there exists a sequence $(s_k)$ of $\mathsf{S}$-simple functions such that
\begin{itemize}
\item[(a)] $0 \leq s_{k} \leq s_{k+1} \leq f$ for each $k \in \mathbb{N}$.
\item[(b)] $s_k \to f$ pointwise on $X$.
\item[(c)] If $f$ is bounded, then $s_{k} \to f$ uniformly on $X$.
\end{itemize}
\end{lemma}

\begin{proof}
Fix $k \in \mathbb{N}$. For each $j \in \{0, 1, 2, \ldots, k2^{k}-1\}$ define the sets
$$
\begin{aligned}
A_{k}(j)&:=\left\{ x \in X\,:\, \frac{j}{2^{k}} \leq f(x) < \frac{j+1}{2^{k}}  \right\},\\
A_{k}(k2^{k})&:=\left\{ x \in X\,:\, k \leq f(x)  \right\}.
\end{aligned}
$$
which belong to $\mathsf{S}$ since $f$ is an $\mathsf{S}$-measurable function. Moreover, the family $\{A_{k}(j)\,:\, 0 \leq j \leq k2^{k} \}$ forms a partition of the set $X$ because $f$ is nonnegative.

Therefore, the function $s_k : X \to \mathbb{R}$ given by
$$
s_{k}(x):=\sum_{j=0}^{k2^k} \frac{j}{2^k}\,\chi_{A_{k}(j)}(x)
$$
is $\mathsf{S}$-simple (see Proposition \ref{315}).

\begin{figure}[ht!]
\centering
\begin{tikzpicture}[xscale=4.35, yscale=4.75]
\draw[domain= 0:1.5] plot(\x,{0.5*\x^2 -(1/3)*(\x)^3} );
\draw[domain= -1.5:0] plot(\x,{-0.5*\x^2 -(1/6)*(\x)^6} );

\draw[->, gray] (-1.5,-0.78)--(1.522,-0.78);
\draw[->, gray] (-1.5,-0.78)--(-1.5,0.5);

\draw (1.5,-0.78) node[right]{$_{(X,\mathsf{S})}$}; 
\draw (-1.5,0.5) node[left]{$_{(\mathbb{R}^{+},\ \mathcal{B}(\mathbb{R}))}$}; 

\draw[gray] (-1.5,-0.7) node{$_{-}$}; \draw[-, densely dotted,gray] (-1.5,-0.7)--(1.5,-0.7);
\draw[gray] (-1.5,-0.6) node{$_{-}$}; \draw[-, densely dotted,gray] (-1.5,-0.6)--(1.5,-0.6);
\draw[gray] (-1.5,-0.5) node{$_{-}$}; \draw[-, densely dotted,gray] (-1.5,-0.5)--(1.5,-0.5);
\draw[gray] (-1.5,-0.4) node{$_{-}$}; \draw[-, densely dotted,gray] (-1.5,-0.4)--(1.5,-0.4);
\draw[gray] (-1.5,-0.3) node{$_{-}$}; \draw[-, densely dotted,gray] (-1.5,-0.3)--(1.5,-0.3);
\draw[gray] (-1.5,-0.2) node{$_{-}$}; \draw[-, densely dotted,gray] (-1.5,-0.2)--(1.5,-0.2);
\draw[gray] (-1.5,-0.1) node{$_{-}$}; \draw[-, densely dotted,gray] (-1.5,-0.1)--(1.5,-0.1);
\draw[gray] (-1.5,0.0) node{$_{-}$}; \draw[-, densely dotted,gray] (-1.5,0.0)--(1.5,0.0);
\draw[gray] (-1.5,0.1) node{$_{-}$}; \draw[-, densely dotted,gray] (-1.5,0.1)--(1.5,0.1);
\draw[gray] (-1.5,0.2) node{$_{-}$}; \draw[-, densely dotted,gray] (-1.5,0.2)--(1.5,0.2);
\draw[gray] (-1.5,0.3) node{$_{-}$}; \draw[-, densely dotted,gray] (-1.5,0.3)--(1.5,0.3);

\draw[ultra thick] (-1.51,-0.78)--(-1.48,-0.78);  \draw[-, densely dotted] (-1.48,-0.7)--(-1.48,-0.78);
\draw[ultra thick] (-1.48,-0.7)--(-1.46,-0.7);  \draw[-, densely dotted] (-1.46,-0.6)--(-1.46,-0.7);
\draw[ultra thick] (-1.47,-0.6)--(-1.45,-0.6);  \draw[-, densely dotted] (-1.45,-0.5)--(-1.45,-0.6);
\draw[ultra thick] (-1.45,-0.5)--(-1.43,-0.5);  \draw[-, densely dotted] (-1.43,-0.4)--(-1.43,-0.5);
\draw[ultra thick] (-1.43,-0.4)--(-1.41,-0.4);  \draw[-, densely dotted] (-1.41,-0.3)--(-1.41,-0.4);
\draw[ultra thick] (-1.41,-0.3)--(-1.38,-0.3);  \draw[-, densely dotted] (-1.38,-0.2)--(-1.38,-0.3);
\draw[ultra thick] (-1.38,-0.2)--(-1.35,-0.2);  \draw[-, densely dotted] (-1.35,-0.1)--(-1.35,-0.2);
\draw[ultra thick] (-1.35,-0.1)--(-1.32,-0.1);  \draw[-, densely dotted] (-1.32,0)--(-1.32,-0.1);

\draw[ultra thick] (-1.32,0)--(-1.27,0);  \draw[-, densely dotted] (-1.27,0.1)--(-1.27,0);
\draw[ultra thick] (-0.45,0)--(0.57,0);  \draw[-, densely dotted] (-0.45,0.1)--(-0.45,0);
\draw[-, densely dotted] (0.57,0.1)--(0.57,0);
\draw[ultra thick] (1.32,0)--(1.5,0); \draw[-, densely dotted] (1.32,0.1)--(1.32,0);

\draw[ultra thick] (-1.27,0.1)--(-1.21,0.1);  \draw[-, densely dotted] (-1.21,0.2)--(-1.21,0.1);
\draw[ultra thick] (-0.65,0.1)--(-0.45,0.1);  \draw[-, densely dotted] (-0.65,0.2)--(-0.65,0.1);
\draw[ultra thick] (0.57,0.1)--(1.32,0.1);

\draw[ultra thick] (-1.21,0.2)--(-1.12,0.2);  \draw[-, densely dotted] (-1.12,0.3)--(-1.12,0.2);
\draw[ultra thick] (-0.85,0.2)--(-0.65,0.2);  \draw[-, densely dotted] (-0.85,0.3)--(-0.85,0.2);

\draw[ultra thick] (-0.85,0.3)--(-1.12,0.3);

\draw (-1.5,-0.78) node[left]{$_{0/2^k}$}; 
\draw (-1.5,-0.7) node[left]{$_{1/2^k}$}; 
\draw (-1.5,-0.6) node[left]{$_{2/2^k}$}; 
\draw (-1.5,-0.5) node[left]{$_{3/2^k}$};
\draw (-1.5,-0.37) node[left]{$_{\vdots}$};
\draw (-1.5,-0.27) node[left]{$_{\vdots}$}; 
\draw (-1.5,-0.17) node[left]{$_{\vdots}$}; 
\draw(-1.5,-0.1) node[left]{$_{j/2^k}$};
\draw (-1.5,0.0) node[left]{$_{(j+1)/2^k}$};
\draw (-1.5,0.13) node[left]{$_{\vdots}$}; 
\draw (-1.5,0.2) node[left]{$_{(k2^k-1)/2^k}$}; 
\draw (-1.5,0.3) node[left]{$_{k}$}; 
\end{tikzpicture}
\end{figure}

\textit{(a):} It is immediate that $0 \leq s_{k} \leq f$ for every $k \in \mathbb{N}$. Let $x \in X$ be arbitrary but fixed and let $k \in \mathbb{N}$. We will show that $s_{k}(x) \leq s_{k+1}(x)$, for which we consider three cases:\\[0.005cm]

{\scshape Case 1.} \quad $k < k+1 \leq f(x)$.

In this case, we have $x \in A_{k}(k2^k) \cap A_{k+1}((k+1)2^{k+1})$, and therefore
$$
s_{k}(x)=\frac{k2^k}{2^k}=k \quad \mbox{and}\quad s_{k+1}(x)=\frac{(k+1)2^{(k+1)}}{2^{k+1}}=k+1.
$$

Consequently, $s_{k}(x) \leq s_{k+1}(x)$.\\[0.005cm]

{\scshape Case 2.}\quad $k \leq f(x) < k+1$.

It is clear that $s_k(x)=k$. Defining $j_{\ast}:=k2^{k+1}$, we have that $0<j_{\ast} < (k+1)2^{k+1}-1$. Thus, since
$$
\frac{j_\ast}{2^{k+1}} \leq f(x) < k+1,
$$
it follows that $x \in A_{k+1}(j_0)$ for some $0<j_\ast \leq j_0 \leq (k+1)2^{k+1}-1$. Consequently,
$$
s_{k+1}(x)=\frac{j_0}{2^{k+1}} \geq \frac{j_\ast}{2^{k+1}} = k = s_k(x)
$$
which proves the result in this case.\\[0.005cm]

{\scshape Case 3.}\quad $f(x)<k$.

There exists $j_{\ast} \in \{0,1,\ldots,k2^k-1 \}$ such that
$$
f(x) \in \left[ \frac{j_{\ast}}{2^k}, \frac{j_{\ast}+1}{2^k} \right)
$$
and therefore $s_{k}(x)=\frac{j_{\ast}}{2^k}$.

Observe that $\frac{j_{\ast}}{2^k}=\frac{2j_{\ast}}{2^{k+1}}$ and $\frac{j_{\ast}+1}{2^k}=\frac{2(j_{\ast}+1)}{2^{k+1}}$, so that
$$
\left[ \frac{j_{\ast}}{2^k}, \frac{j_{\ast}+1}{2^k} \right)=\left[ \frac{2j_{\ast}}{2^{k+1}}, \frac{2j_{\ast}+1}{2^{k+1}} \right) \cup \left[ \frac{2j_{\ast}+1}{2^{k+1}}, \frac{2j_{\ast}+2}{2^{k+1} } \right).
$$

Thus,
$$
f(x) \in \left[ \frac{2j_{\ast}}{2^{k+1}}, \frac{2j_{\ast}+1}{2^{k+1}} \right)\quad \mbox{or}\quad f(x) \in \left[ \frac{2j_{\ast}+1}{2^{k+1}}, \frac{2j_{\ast}+2}{2^{k+1} } \right)
$$
and consequently
$$
s_{k+1}(x)=\left\{
\begin{array}{ccl}
\frac{2j_{\ast}}{2^{k+1}} & & \mbox{if}\,\,\, f(x) \in \left[ \frac{2j_{\ast}}{2^{k+1}}, \frac{2j_{\ast}+1}{2^{k+1}} \right),\\
&&\\
\frac{2j_\ast+1}{2^{k+1}} & &  \mbox{if}\,\,\, f(x) \in \left[ \frac{2j_{\ast}+1}{2^{k+1}}, \frac{2j_{\ast}+2}{2^{k+1} } \right).\\
\end{array}
\right.
$$

Consider the following subcases:

{\scshape Subcase 3.1.} 
$$
s_{k+1}(x)=\frac{2j_{\ast}}{2^{k+1}}=\frac{j_\ast}{2^k}\geq s_k(x).
$$

{\scshape Subcase 3.2.} 
$$
s_{k+1}(x)=\frac{2j_\ast +1}{2^{k+1}}=\frac{j_\ast}{2^k}+\frac{1}{2^{k+1}} \geq \frac{j_\ast}{2^k} =s_{k}(x).
$$

This proves that $(s_k)$ is nondecreasing.

\textit{(b):} Let $\varepsilon >0$. Choose $k(\varepsilon) \in \mathbb{N}$ such that $\frac{1}{2^{k(\varepsilon)}} < \varepsilon.$ Since $f$ is nonnegative, for each $x \in X$ there exists $k(x) \in \mathbb{N}$ such that $f(x) < k(x)$.

Let $x \in X$. Take $k_0:=\max\{k(\varepsilon),k(x) \} \in \mathbb{N}$. Then, for $k \geq k_0$, we have $f(x)< k$ and therefore there exists $j \in \{ 0,\ldots,k2^{k}-1 \}$ such that
$$
s_{k}(x),f(x) \in \left[ \frac{j}{2^k},\frac{j+1}{2^k} \right).
$$

Consequently,
$$
|s_k(x)-f(x)| < \frac{1}{2^k} < \varepsilon \qquad \forall k \geq k_0,
$$
that is, $(s_k)$ converges pointwise to $f$ on $X$.

\textit{(c):} Since $f$ is bounded, there exists $c >0$ such that $f(x) < c$ for every $x \in X$. Thus, it is possible to find $k(c) \in \mathbb{N}$ such that $c < k(c)$.

Let $\varepsilon >0$. Choose $k(\varepsilon) \in \mathbb{N}$ such that $\frac{1}{2^{k(\varepsilon)}} < \varepsilon.$ Defining $k_0:=\max\{k(\varepsilon),k(c) \} \in \mathbb{N}$ and following the idea of the previous part, we obtain
$$
|s_k(x)-f(x)| < \frac{1}{2^k} < \varepsilon \qquad \forall k \geq k_0,\quad \forall x \in X.
$$

That is, $(s_k)$ converges uniformly to $f$ on $X$.
\end{proof}

\begin{lemma} \label{319}
Let $(X,\mathsf{S})$ be a measurable space and let $f:X \to \mathbb{R}$ be an $\mathsf{S}$-measurable function. Then there exists a sequence $(r_k)$ of $\mathsf{S}$-simple functions such that
\begin{itemize}
\item[(a)] $|r_k| \leq |f|$ for every $k \in \mathbb{N}$.
\item[(b)] $r_k \to f$ pointwise on $X$.
\item[(c)] If $f$ is bounded, then $r_k \to f$ uniformly on $X$.
\end{itemize}
\end{lemma}

\begin{proof}
\textit{(a)} and \textit{(b):} For an $\mathsf{S}$-measurable function $f:X \to \mathbb{R}$, consider the nonnegative and $\mathsf{S}$-measurable functions $f^{+},f^{-}:X \to \mathbb{R}$ (see Proposition \ref{317}).

Applying Lemma \ref{318} to $f^{+}$ and $f^{-}$, we obtain two nondecreasing sequences $(s_k)$ and $(t_k)$ of nonnegative $\mathsf{S}$-simple functions such that $s_k(x) \to f^{+}(x)$ and $t_k(x) \to f^{-}(x)$ for each $x \in X$.

Define $r_{k}:=s_{k}-t_{k}$ for each $k \in \mathbb{N}$. Then $(r_k)$ is a sequence of $\mathsf{S}$-simple functions (see Proposition \ref{316}) such that
$$
|r_k|:=|s_k - t_k| \leq |s_k| + |t_k| \leq |f^+|+|f^-| =f^{+} + f^{-} = |f|\quad \forall k \in \mathbb{N}.
$$ 

Moreover, it is clear that $r_k(x) \to f(x)$ for every $x \in X$ since
$$
f(x)=f^{+}(x)-f^{-}(x)=\displaystyle\lim_{k \to \infty} s_k(x) - \lim_{k \to \infty} t_k(x) = \lim_{k \to \infty}(s_k- t_k)(x) = \lim_{k \to \infty} r_k(x)
,
$$
for every $x \in X$.

\textit{(c):} Since $f$ is bounded, there exists $c>0$ such that $|f(x)| \leq c$ for every $x \in X$. Using the fact that $f^{+},f^{-}$ are nonnegative and satisfy $f^{+}+f^{-} = |f|$, we conclude that $f^{+}$ and $f^{-}$ are bounded $\mathsf{S}$-measurable functions. Thus, by part \textit{(c)} of Lemma \ref{318}, $s_k \to f^{+}$ and $t_k \to f^-$ uniformly on $X$, and therefore $r_k \to f$ uniformly on $X$.
\end{proof}

The following result establishes the conditions under which the measurability criterion for a sequence of $\mathsf{S}$-measurable functions remains valid when passing to a limiting process.

\begin{theorem} \label{320}
Let $(X,\mathsf{S})$ be a measurable space and let $(f_k)$ be a sequence of $\mathsf{S}$-measurable functions such that $(f_k(x))$ is bounded for every $x \in X$. The following functions are $\mathsf{S}$-measurable.
\begin{itemize}
\item[(a)] $m_k(x):=\displaystyle\min_{1\leq j \leq k}f_j(x)$ for fixed $k \in \mathbb{N}$.
\item[(b)]  $M_{k}(x):=\displaystyle\max_{1\leq j \leq k} f_j(x)$ for fixed $k \in \mathbb{N}$.
\item[(c)] $\left(\displaystyle\inf_{k \in \mathbb{N}} f_k\right)(x):=\displaystyle\inf_{k \in \mathbb{N}} f_k(x)$.
\item[(d)] $\left(\displaystyle\sup_{k\in \mathbb{N}} f_k\right)(x):=\displaystyle\sup_{k\in \mathbb{N}} f_k(x)$.
\item[(e)] $\left(\displaystyle\liminf_{k \to \infty} f_k\right)(x) :=\displaystyle\sup_{k \geq 1}\inf_{ j \geq k} f_j(x)$.
\item[(f)] $\left(\displaystyle\limsup_{k \to \infty} f_k\right)(x):=\displaystyle\inf_{k \geq 1}\sup_{ j \geq k} f_j(x)$.
\end{itemize}
\end{theorem}

\begin{proof}
Let $c \in \mathbb{R}$ be arbitrary.

\textit{(a):} Let $x \in X$. It is clear that $m_k(x) > c$ if and only if $f_j(x) > c$ for every $j=1,\ldots,k$. Consequently,
$$
m_k^{-1}((c,+\infty))=\{x \in X\,:\, m_k(x)> c\} = \bigcap_{j=1}^{k}\{x \in X\,:\, f_j(x)> c\} = \bigcap_{j=1}^{k} f_{j}^{-1}((c,+\infty))
$$
where $f_{j}^{-1}((c,+\infty)) \in \mathsf{S}$ for every $j=1,\ldots,k$. Therefore, $m_k^{-1}((c,+\infty)) \in \mathsf{S}$, and this proves that $m_k$ is $\mathsf{S}$-measurable.

\textit{(b):} Let $x \in X$. Observe that $M_k(x) > c$ if and only if there exists $j_0\in \{1,\ldots,k\}$ such that $f_{j_0}(x) > c$. Hence,
$$
M_k^{-1}((c,+\infty))=\{x \in X\,:\, M_k(x)> c\} = \bigcup_{j=1}^{k}\{x \in X\,:\, f_j(x)> c\} = \bigcup_{j=1}^{k} f_{j}^{-1}((c,+\infty))
$$
from which it follows that $M_k^{-1}((c,+\infty)) \in \mathsf{S}$ since $f_j^{-1}((c,+\infty)) \in \mathsf{S}$ for every $j=1,\ldots,k$.

\begin{figure}[ht!]
\begin{minipage}[l]{0.5\textwidth}
\begin{center}
\begin{tikzpicture}[xscale=2.15, yscale=2.15]
\draw[domain= 0:2,thick] plot(\x,{{exp(-\x)}} );
\draw[domain= 0:2,thick] plot(\x,{{2*\x*exp(-\x^2)}} );
\draw[domain= 0:2,thick] plot(\x,{{3*\x^2*exp(-\x^3)}} );
\draw[domain= 0:2,thick] plot(\x,{{4*\x^3*exp(-\x^4)}} );

\draw[->,gray] (0,-0.01)--(2.1,-0.01); 
\draw[->,gray] (0,-0.01)--(0,1.55); 
\draw (1.16,-0.01) node[below]{$_{x_0}$};
\draw (0,0.32) node{$_{-}$};
\draw (0,0.32) node[left]{$_{m_k(x_{0})}$};
\draw[-, densely dotted,gray] (1.16,-0.01)--(1.16,0.32);
\draw[-, densely dotted,gray] (0,0.32)--(2.1,0.32);

\draw (0,0.15) node{$_{-}$};
\draw (0,0.15) node[left]{$_{c}$};
\draw[-, densely dotted,gray] (0,0.15)--(2.1,0.15);
\end{tikzpicture}
\end{center}
\end{minipage} \hfill 
 \begin{minipage}[r]{0.5\textwidth}
\begin{center}
\begin{tikzpicture}[xscale=2.15, yscale=2.15]
\draw[domain= 0:2,thick] plot(\x,{{exp(-\x)}} );
\draw[domain= 0:2,thick] plot(\x,{{2*\x*exp(-\x^2)}} );
\draw[domain= 0:2,thick] plot(\x,{{3*\x^2*exp(-\x^3)}} );
\draw[domain= 0:2,thick] plot(\x,{{4*\x^3*exp(-\x^4)}} );

\draw[->,gray] (0,-0.01)--(2.1,-0.01); 
\draw[->,gray] (0,-0.01)--(0,1.55); 
\draw (1.16,-0.01) node[below]{$_{x_0}$};
\draw (0,1) node{$_{-}$};
\draw (0,1) node[left]{$_{M_k(x_{0})}$};
\draw[-, densely dotted,gray] (1.16,-0.01)--(1.16,1);
\draw[-, densely dotted,gray] (0,1)--(2.1,1);

\draw (0,0.5) node{$_{-}$};
\draw (0,0.5) node[left]{$_{c}$};
\draw[-, densely dotted,gray] (0,0.5)--(2.1,0.5);
\end{tikzpicture}
\end{center}
\end{minipage}
\end{figure}

\textit{(c):} Let $x \in X$. If $\inf_{k \geq 1}f_k(x) \geq c$, then $c$ is a lower bound of the set $\{ f_k(x)\,:\, k \in \mathbb{N}\}$ and therefore $f_k(x) \geq c$ for every $k \in \mathbb{N}$. Conversely, if $f_k(x) \geq c$ for every $k \in \mathbb{N}$, then by the definition of infimum it follows that $\inf_{k \geq 1}f_k(x) \geq c$.

Therefore,
$$
\left\{x \in X\,:\, \inf_{k \geq 1}f_k(x) \geq c\right\}=\bigcap_{k=1}^{\infty}f_k^{-1}([c,+\infty)) \in \mathsf{S}
$$
since $f_k$ is an $\mathsf{S}$-measurable function for every $k \in \mathbb{N}$.

\textit{(d):} Let $x \in X$. If $\sup_{k \geq 1}f_k(x) \leq c$, then $c$ is an upper bound of the set $\{ f_k(x)\,:\, k \in \mathbb{N}\}$ and therefore $f_k(x) \leq c$ for every $k \in \mathbb{N}$. Now, if $f_k(x) \leq c$ for every $k \in \mathbb{N}$, then by the definition of supremum we conclude that $\sup_{k \geq 1}f_k(x) \leq c$.

Consequently,
$$
\left\{x \in X\,:\, \sup_{k \geq 1}f_k(x) \leq c\right\}=\bigcap_{k=1}^{\infty}f_k^{-1}((-\infty,c] ) \in \mathsf{S}
$$
since $f_k$ is an $\mathsf{S}$-measurable function for every $k \in \mathbb{N}$.

\textit{(e):} Define, for each $k \in \mathbb{N}$, the function $g_k:=\inf_{j \geq k} f_j$. Part \textit{(c)} ensures that $(g_k)$ is a sequence of $\mathsf{S}$-measurable functions. Then, part \textit{(d)} implies that $\sup_{k \geq 1} g_k$ is an $\mathsf{S}$-measurable function. That is, $\liminf_{k \to \infty} f_k$ is $\mathsf{S}$-measurable.

\textit{(f)} Defining $h_k:=\sup_{j \geq k} f_j$ for $k \in \mathbb{N}$ and applying parts \textit{(c)} and \textit{(d)}, we conclude that $\inf_{ k \geq 1} h_k=\limsup_{k \to \infty}f_k$ is an $\mathsf{S}$-measurable function.

This completes the proof.
\end{proof}

The following result establishes that the pointwise limit of a sequence of $\mathsf{S}$-measurable functions is again an $\mathsf{S}$-measurable function. This fact marks a fundamental difference with other classes of functions studied in previous courses; for example, the pointwise limit of continuous functions need not be continuous [Exercise \ref{Ej12}].

\begin{corollary} \label{321}
Let $(X,\mathsf{S})$ be a measurable space. If $f_k:X \to \mathbb{R}$ is $\mathsf{S}$-measurable for every $k \in \mathbb{N}$ and $(f_k)$ converges pointwise to $f$ on $X$, then $f:X \to \mathbb{R}$ is $\mathsf{S}$-measurable.
\end{corollary}

\begin{proof}
Since $(f_k(x))$ converges to $f(x)$ in $\mathbb{R}$ for each $x \in X$, then [Exercise \ref{E317}]
$$
f(x)=\liminf_{k \to \infty} f_k(x) = \limsup_{k \to \infty} f_k(x)\quad \mbox{for each}\,\,x \in X
$$
and therefore $f$ is $\mathsf{S}$-measurable by parts \textit{(e)} and \textit{(f)} of Theorem \ref{320}.
\end{proof}

We therefore obtain the following characterization of $\mathsf{S}$-measurable functions.

\begin{theorem} \label{322}
Let $(X,\mathsf{S})$ be a measurable space. Then $f:X \to \mathbb{R}$ is an $\mathsf{S}$-measurable function if and only if $f$ is the pointwise limit of a sequence of $\mathsf{S}$-simple functions.
\end{theorem}

\begin{proof}
$\Rightarrow):$ This follows from Lemma \ref{319}. $\Leftarrow):$ This is immediate from Corollary \ref{321}.
\end{proof}

Using this criterion for $\mathsf{S}$-measurability, it is possible to establish algebraic properties of $\mathsf{S}$-measurable functions, as shown in the following result.

\begin{theorem} \label{323}
Let $(X,\mathsf{S})$ be a measurable space and let $f,g:X \to \mathbb{R}$ be $\mathsf{S}$-measurable functions. Then $\gamma \,f$ with $\gamma \in \mathbb{R}$, $f+g$, and $f g$ are $\mathsf{S}$-measurable.
\end{theorem}

\begin{proof}
By Theorem \ref{322}, there exist two sequences $(s_k)$ and $(t_k)$ of $\mathsf{S}$-simple functions such that $s_k \to f$ and $t_k \to g$ pointwise on $X$. Thus, the sequences $(\gamma s_k)$, $(s_k+t_k)$, and $(s_k  t_k)$ of $\mathsf{S}$-simple functions (see Proposition \ref{316}) converge pointwise on $X$ to the functions $\gamma f$, $f+g$, and $f g$, respectively. Theorem \ref{322} ensures that $\gamma f$, $f+g$, and $fg$ are $\mathsf{S}$-measurable.
\end{proof}

We conclude this section with an interesting application of the approximation lemmas for simple functions.

\begin{proposition} \label{324}
Let $(X,\mathsf{S})$ be a measurable space. Given a function $f:X \to \mathbb{R}$, consider the $\sigma$-algebra of subsets of $X$ defined by
$$
\sigma(f):=\left\{ f^{-1}(B) \,:\, B \in \mathcal{B}(\mathbb{R}) \right\}.
$$

Then, for every $\sigma(f)$-measurable function $\eta:X \to \mathbb{R}$, there exists a Borel-measurable function $\phi:\mathbb{R} \to \mathbb{R}$ such that $\eta= \phi \circ f$ on $X$.
\end{proposition}

\begin{proof}
Consider the following cases:

{\scshape Case 1.}\quad $\eta$ is a $\sigma(f)$-simple function with canonical representation $\eta=\sum_{j=1}^{N} \alpha_{j}\chi_{A_{j}}$.

For each $j=1,\ldots,N$, there exists $B_{j} \in \mathcal{B}(\mathbb{R})$ such that $A_{j}=f^{-1}(B_{j})$. Define the function $\phi :\mathbb{R} \to \mathbb{R}$ by
$$
\phi(\xi):=\sum_{j=1}^{N} \alpha_{j}\chi_{B_{j}}(\xi).
$$

It is clear that the class $\{B_{1},\ldots,B_{N}\}$ forms a partition of $\mathbb{R}$ consisting of subsets in $\mathcal{B}(\mathbb{R})$. Proposition \ref{315} then ensures that $\phi$ is a $\mathcal{B}(\mathbb{R})$-simple function. 

Let $j=1,\ldots,N$ and let $x \in X$. If $x \in A_{j}=f^{-1}(B_{j})$, then $f(x) \in B_{j}$ and therefore $\chi_{B_{j}}(f(x))=1=\chi_{A_j}(x)$. Now, if $x \in X\smallsetminus f^{-1}(B_{j})$, then $f(x) \in \mathbb{R} \smallsetminus B_{j}$, from which we conclude that $\chi_{B_{j}}(f(x))=0=\chi_{A_j}(x)$. Hence,
$$
\phi(f(x))=\sum_{j=1}^{N} \alpha_{j} \chi_{B_{j}}(f(x)) =\sum_{j=1}^{N} \alpha_{j} \chi_{A_{j}}(x)=\eta(x)
$$
for every $x \in X$.

This proves the equality $\eta=\phi \circ f$ in this case.

{\scshape Case 2.}\quad $\eta$ is an arbitrary $\sigma(f)$-measurable function.

Lemma \ref{319} ensures that there exists a sequence $(\eta_k)$ of $\sigma(f)$-simple functions such that $\eta_k \to \eta$ pointwise on $X$. By {\scshape Case 1}, for each $k \in \mathbb{N}$, there exists a Borel-measurable function $\phi_{k}:\mathbb{R} \to \mathbb{R}$ such that $\eta_{k} = \phi_{k} \circ f$ on $X$. 

Let $A:=\{\xi \in \mathbb{R}\,:\,\big(\phi_{k}(\xi)\big)\,\text{ converges in }\mathbb{R}\}$. Since $(\phi_{k})$ is a sequence of Borel-measurable functions, Exercise \ref{315} ensures that $A \in \mathcal{B}(\mathbb{R})$. Define $\phi:\mathbb{R} \to \mathbb{R}$ by
$$
\phi(\xi):=\left\{
\begin{array}{lcl}
     \lim_{k \to \infty}\phi_{k}(\xi) & & \text{if }\,\xi \in A,  \\
     0 & & \text{if }\,\xi \not\in A.
\end{array}
\right.
$$

Exercise \ref{310} ensures that $\phi$ is a Borel-measurable function. 

Let $x \in X$. By Lemma \ref{319} and {\scshape Case 1}, we have $f(x) \in A$ and therefore
$$
\phi(f(x))=\lim_{k \to \infty} \phi_k(f(x)) =\lim_{k \to \infty} \eta_{k}(x)=\eta(x).
$$

This proves that $\eta=\phi \circ f$ on $X$, and the proof is complete.
\end{proof}

\section{The extended real line}

To study sequences of real-valued functions in a more complete way, particularly when uniform bounds are not available, it is convenient to enlarge the set of real numbers. For this purpose, the symbols $+\infty$ and $-\infty$ are introduced, allowing one to describe precisely phenomena of divergence and unbounded growth. This extension provides an appropriate framework for formulating notions of convergence that arise naturally in measure and integration theory.

\begin{definition} \label{325}
We define the set of extended real numbers (\textbf{extended real line}) \index{extended real numbers $\overline{\mathbb{R}}$}, denoted by $\overline{\mathbb{R}}$, as
$$
\overline{\mathbb{R}} := \{-\infty \} \cup \mathbb{R} \cup \{ + \infty \}.
$$
\end{definition}

We endow $\overline{\mathbb{R}}$ with the operations of addition (+) and multiplication ($\cdot$) between real numbers and these symbols as follows:
\begin{itemize}
\item[\it (1)] $x+(\pm \infty):=(\pm \infty)+x=(\pm \infty)+(\pm \infty)=\pm \infty$ for every $x \in \mathbb{R}$.
\item[\it (2)] $(\pm \infty) \cdot (\pm \infty) := +\infty$ and $(\pm \infty) \cdot (\mp\infty) := -\infty$.
\item[\it (3)] 
$$
x \cdot (\pm \infty) := (\pm \infty) \cdot x: = \left\{
\begin{array}{l c l}
\pm \infty & &\mbox{if}\,\, x>0,\\
0 & & \mbox{if}\,\,x=0,\\
\mp \infty & & \mbox{if}\,\,x<0.
\end{array}
\right.
$$
\end{itemize}

The following expressions will not be defined:
$$
(+\infty) + (-\infty),\quad (-\infty)+(+\infty),\quad \frac{x}{\pm \infty}\quad \text{and}\quad \frac{x}{0} \qquad \forall x \in \mathbb{R}.
$$

The total order considered on $\overline{\mathbb{R}}$ will be the usual one on $\mathbb{R}$ together with the condition that, for every $x \in \mathbb{R}$,
$$
-\infty < x < \infty.
$$

\begin{definition} \label{326}
Let $A\subset\overline{\mathbb{R}}$. If $A$ is not bounded above, we define $\sup (A):= + \infty$. Analogously, if $A$ is not bounded below, we define $\inf (A):=-\infty$.
\end{definition}

\begin{definition} \label{327}
Let $(x_k)$ be a sequence of elements of $\overline{\mathbb{R}}$. We define the \textbf{limit inferior} and \textbf{limit superior} of $(x_k)$, denoted respectively by $\displaystyle\liminf_{k \to \infty} x_k$ and $\displaystyle\limsup_{k \to \infty} x_k$, as follows: \index{limit!inferior of extended real numbers} \index{limit!superior of extended real numbers}
\begin{itemize}
\item[(1)] $\displaystyle\liminf_{k \to \infty} x_k := \sup_{k \geq 1} \inf_{j \geq k} x_j$
\item[(2)] $\displaystyle\limsup_{k \to \infty} x_k := \inf_{k \geq 1} \sup_{j \geq k} x_j.$
\end{itemize}

The limit inferior and limit superior of a sequence $(x_k)$ always exist in $\overline{\mathbb{R}}$ and satisfy
$\displaystyle\liminf_{k \to \infty} x_k \leq \displaystyle\limsup_{k \to \infty} x_k$.
\end{definition}

It is easy to verify [Exercise \ref{E317}] that a sequence of extended real numbers $(x_k)$ converges in $\overline{\mathbb{R}}$ if and only if
$$
\liminf_{k \to \infty} x_k = \limsup_{k \to \infty} x_k.
$$ 

These concepts will play a fundamental role in the subsequent development of the theory. For this reason, many of their properties are stated only as exercises [Exercise \ref{E317}], since they are concepts studied in the first courses on differential and integral calculus.

The next step will be to introduce a $\sigma$-algebra of subsets of $\overline{\mathbb{R}}$ in order to work with the extended real line as a measurable space.

\begin{definition} \label{328}
We define the \textbf{extended Borel $\sigma$-algebra}, \index{sigma@$\sigma$!algebra!extended Borel} as the $\sigma$-algebra of subsets of $\overline{\mathbb{R}}$ given by
$$
\mathcal{B}(\overline{\mathbb{R}}):=\sigma\left(\{-\infty \} \cup \mathcal{B}(\mathbb{R}) \cup \{+\infty \}\right).
$$
\end{definition}

If $\mathcal{C} \subset \mathcal{P}(\mathbb{R})$ is any class of subsets of $\mathbb{R}$ such that $\sigma(\mathcal{C}) = \mathcal{B}(\mathbb{R})$ [Exercise \ref{E227}], then
$$
\mathcal{B}(\overline{\mathbb{R}})=\sigma\left(\{-\infty \} \cup \mathcal{C} \cup \{+\infty \}\right).
$$

From all of the above, the extended real numbers $\overline{\mathbb{R}}$ endowed with the extended Borel $\sigma$-algebra, $(\overline{\mathbb{R}},\mathcal{B}(\overline{\mathbb{R}}))$, form a measurable space, and thus we have the following definition:

\begin{definition} \label{329}
Let $(X,\mathsf{S})$ be a measurable space. A function $f:X \to \overline{\mathbb{R}}$ is said to be \textbf{$\mathsf{S}$-measurable} if $f^{-1}(\mathcal{B}(\overline{\mathbb{R}})) \subset \mathsf{S}$. \index{function!$\mathsf{S}$-measurable}
\end{definition}

Following the idea of Definition \ref{328}, we can characterize in the next proposition an $\mathsf{S}$-measurable function with values in the extended real numbers in a way analogous to Proposition \ref{34}.

\begin{proposition} \label{330}
Let $(X,\mathsf{S})$ be a measurable space and let $f:X \to \overline{\mathbb{R}}$ be a function. The following statements are equivalent:
\begin{itemize}
\item[(a)] $f$ is $\mathsf{S}$-measurable.
\item[(b)] $f^{-1}((c,+\infty])$ belongs to $\mathsf{S}$ for every $c \in \mathbb{R}$.
\item[(c)] $f^{-1}([-\infty,c])$ belongs to $\mathsf{S}$ for every $c \in \mathbb{R}$.
\item[(d)] $f^{-1}([-\infty,c))$ belongs to $\mathsf{S}$ for every $c \in \mathbb{R}$.
\item[(e)] $f^{-1}([c,+\infty])$ belongs to $\mathsf{S}$ for every $c \in \mathbb{R}$.
\end{itemize}

Where $(c,+\infty]:=(c,+\infty) \cup \{+\infty \}$, $[c,+\infty]:=[c,+\infty)\cup \{+\infty \}$, $[-\infty, c):=\{-\infty \} \cup (-\infty,c)$ and $[-\infty ,c]:=\{-\infty \} \cup (-\infty,c]$ for every $c \in \mathbb{R}$.
\end{proposition}

Observe that in the previous proposition, $c$ must be a real number since, in general, it is not true that if $f^{-1}(\{c \}) \in \mathsf{S}$ for every $c \in \overline{\mathbb{R}}$, then $f:X \to \overline{\mathbb{R}}$ is $\mathsf{S}$-measurable [Exercise \ref{E319}], and therefore the equivalence between parts \textit{(e)} and \textit{(a)} would fail.

The following theorem characterizes an $\mathsf{S}$-measurable function taking values in the extended real numbers through the study of a function taking values in the real numbers, to which we may apply the results developed in the previous sections.

\begin{theorem} \label{331}
Let $(X,\mathsf{S})$ be a measurable space. A function $f:X \to \overline{\mathbb{R}}$ is $\mathsf{S}$-measurable if and only if
\begin{itemize}
\item[(a)] The sets $A_{+\infty}(f):=\{x \in X\,:\, f(x)=+\infty \}$ and $A_{-\infty}(f):=\{x \in X\,:\, f(x)=-\infty \}$ belong to $\mathsf{S}$.
\item[(b)] The real-valued function $f_{0}:X \to \mathbb{R}$ defined by
$$
f_{0}(x):=\left\{
\begin{array}{l c l}
f(x) & & \mbox{if}\,\,x \not \in A_{+\infty}(f) \cup A_{-\infty}(f),\\
0 & &\mbox{if}\,\,x  \in A_{+\infty}(f) \cup A_{-\infty}(f),
\end{array}
\right.
$$
is $\mathsf{S}$-measurable.
\end{itemize}
\end{theorem}

\begin{figure}[ht!]
\begin{minipage}[l]{0.5\textwidth}
\begin{center}
	\begin{tikzpicture}[xscale=0.7,yscale=0.65]
\draw[->, gray] (-4,0)--(6,0); 
\draw[->, gray] (0,-4)--(0,4); 
\draw[thick,blue, domain=-3:0.5] plot (\x, {1/(\x-1)^2});
\draw[thick,blue, domain=2:5.5] plot (\x, {-1/(\x-1.5)^2});
\draw[-, dotted] (0.6,-4)--(0.6,4);
\draw[-, dotted] (1.9,-4)--(1.9,4);
\draw (1.25,4.3) node{$_{f:X \to \overline{\mathbb{R}}}$};
\end{tikzpicture}
\end{center}
\end{minipage} \hfill 
 \begin{minipage}[r]{0.5\textwidth}
\begin{center}
\begin{tikzpicture}[xscale=0.7,yscale=0.65]
\draw[->, gray] (-4,0)--(6,0); 
\draw[->,gray] (0,-4)--(0,4); 
\draw[thick, domain=-3:0.5] plot (\x, {1/(\x-1)^2});
\draw[thick, domain=2:5.5] plot (\x, {-1/(\x-1.5)^2});
\draw[-, dotted] (0.6,-4)--(0.6,4);
\draw[-, dotted] (1.9,-4)--(1.9,4);
\draw (1.25,4.3) node{$_{f_{0}:X \to\mathbb{R}}$};
\draw (0.6,0) node{$\bullet$};
\draw (1.9,0) node{$\bullet$};
\draw (-0.2,0.2) node{$_{0}$};
\end{tikzpicture}
\end{center}
\end{minipage}
\end{figure}

\begin{proof}
$\Rightarrow):$ Since $\{ \pm \infty \} = [\pm \infty, \pm \infty]$ and $f$ is $\mathsf{S}$-measurable, Proposition \ref{330} ensures that the sets
$$
A_{\pm \infty} (f)  =f^{-1}(\{\pm \infty \}) = f^{-1}([\pm \infty, \pm \infty])
$$
belong to $\mathsf{S}$. This proves \textit{(a)}.

Now let us show that $f_{0}$ is $\mathsf{S}$-measurable. Let $c \in \mathbb{R}$ be arbitrary.

{\scshape Case 1.}\quad If $c \geq 0$, then
$$
\begin{aligned}
f^{-1}_{0}((c,+\infty))=\{x \in X\,:\, f_{0}(x)>c \}&=\{x \in X\,:\, f(x)>c \}\smallsetminus\{x \in X\,:\, f(x)=+\infty \}\\
&=f^{-1}((c,+\infty])\smallsetminus A_{+\infty}(f) \in \mathsf{S}
\end{aligned}
$$
by the previous argument and Proposition \ref{330}.

{\scshape Case 2.}\quad If $c<0$, then
$$
\begin{aligned}
f^{-1}_{0}((c,+\infty))=\{x \in X\,:\, f_{0}(x)>c \}&=\{x \in X\,:\, f(x)>c \}\cup \{x \in X\,:\, f(x)=-\infty \}\\
&=f^{-1}((c,+\infty])\cup A_{-\infty}(f) \in \mathsf{S}
\end{aligned}
$$
again by Proposition \ref{330} and the previous argument.

Therefore, $f_{0}$ is $\mathsf{S}$-measurable, and this proves \textit{(b)}.

\begin{figure}[ht!]
\begin{minipage}[l]{0.5\textwidth}
\begin{center}
	\begin{tikzpicture}[xscale=0.7,yscale=0.65]
\draw[->, gray] (-4,0)--(6,0); 
\draw[->, gray] (0,-4)--(0,4); 
\draw[thick, domain=-3:0.5, blue] plot (\x, {1/(\x-1)^2});
\draw[thick, domain=2:5.5, blue] plot (\x, {-1/(\x-1.5)^2});
\draw[-, dotted] (0.6,-4)--(0.6,4);
\draw[-, dotted] (1.9,-4)--(1.9,4);
\draw (1.25,4.3) node{$_{f:X \to\mathbb{R}}$};

\draw (0,1.5) node{$-$}; \draw (-0.45,1.65) node{$c$};
\draw[-,dotted] (-3.6,1.5) -- (5.5,1.5);

\end{tikzpicture}
\end{center}
\end{minipage} \hfill 
 \begin{minipage}[r]{0.5\textwidth}
\begin{center}
	\begin{tikzpicture}[xscale=0.7,yscale=0.65]
\draw[->, gray] (-4,0)--(6,0); 
\draw[->, gray] (0,-4)--(0,4); 
\draw[thick, domain=-3:0.5] plot (\x, {1/(\x-1)^2});
\draw[thick, domain=2:5.5] plot (\x, {-1/(\x-1.5)^2});
\draw[-, dotted] (0.6,-4)--(0.6,4);
\draw[-, dotted] (1.9,-4)--(1.9,4);
\draw (1.25,4.3) node{$_{f_{0}:X \to\mathbb{R}}$};
\draw (0.6,0) node{$\bullet$};
\draw (1.9,0) node{$\bullet$};
\draw (-0.2,0.2) node{$_{0}$};
\draw (0,1.5) node{$-$}; \draw (-0.45,1.65) node{$c$};
\draw[-,dotted] (-3.6,1.5) -- (5.5,1.5);

\end{tikzpicture}
\end{center}
\end{minipage}
\end{figure}

\begin{figure}[ht!]
\begin{minipage}[l]{0.5\textwidth}
\begin{center}
	\begin{tikzpicture}[xscale=0.7,yscale=0.65]
\draw[->, gray] (-4,0)--(6,0); 
\draw[->, gray] (0,-4)--(0,4); 
\draw[thick, domain=-3:0.5, blue] plot (\x, {1/(\x-1)^2});
\draw[thick, domain=2:5.5, blue] plot (\x, {-1/(\x-1.5)^2});
\draw[-, dotted] (0.6,-4)--(0.6,4);
\draw[-, dotted] (1.9,-4)--(1.9,4);
\draw (1.25,4.3) node{$_{f:X \to\mathbb{R}}$};

\draw (0,-1.5) node{$-$}; \draw (-0.45,-1.35) node{$c$};
\draw[-,dotted] (-3.6,-1.5) -- (5.5,-1.5);
\end{tikzpicture}
\end{center}
\end{minipage} \hfill 
 \begin{minipage}[r]{0.5\textwidth}
\begin{center}
	\begin{tikzpicture}[xscale=0.7,yscale=0.65]
\draw[->, gray] (-4,0)--(6,0); 
\draw[->, gray] (0,-4)--(0,4); 
\draw[thick, domain=-3:0.5] plot (\x, {1/(\x-1)^2});
\draw[thick, domain=2:5.5] plot (\x, {-1/(\x-1.5)^2});
\draw[-, dotted] (0.6,-4)--(0.6,4);
\draw[-, dotted] (1.9,-4)--(1.9,4);
\draw (1.25,4.3) node{$_{f_{0}:X \to\mathbb{R}}$};
\draw (0.6,0) node{$\bullet$};
\draw (1.9,0) node{$\bullet$};
\draw (-0.2,0.2) node{$_{0}$};

\draw (0,-1.5) node{$-$}; \draw (-0.45,-1.35) node{$c$};
\draw[-,dotted] (-3.6,-1.5) -- (5.5,-1.5);
\end{tikzpicture}
\end{center}
\end{minipage}
\end{figure}

$\Leftarrow):$ Let $c \in \mathbb{R}$ be arbitrary.

{\scshape Case 1.}\quad If $c \geq 0$, then
$$
\begin{aligned}
f^{-1}((c,+\infty])=\{x \in X\,:\,f(x)>c \}&=\{ x\in X\,:\, f_{0}(x)>c\} \cup A_{+\infty}(f)\\
&=f_{0}^{-1}((c,+\infty)) \cup A_{+\infty}(f) \in \mathsf{S}
\end{aligned}
$$
since $f_{0}$ is $\mathsf{S}$-measurable and $A_{\pm \infty} \in \mathsf{S}$.

{\scshape Case 2.} \quad If $c < 0$, then
$$
\begin{aligned}
f^{-1}((c,+\infty])=\{x \in X\,:\,f(x)>c \}&=\{ x\in X\,:\, f_{0}(x)>c\} \smallsetminus A_{-\infty}(f)\\
&=f_{0}^{-1}((c,+\infty)) \smallsetminus A_{-\infty}(f) \in \mathsf{S}.
\end{aligned}
$$
Thus, by Proposition \ref{330}, we conclude that $f$ is $\mathsf{S}$-measurable.
\end{proof}

Once again, we characterize measurable functions $f:X \to \overline{\mathbb{R}}$ as those that arise as the pointwise limit of a sequence of simple functions.

\begin{theorem} \label{332}
Let $(X,\mathsf{S})$ be a measurable space and let $f:X \to \overline{\mathbb{R}}$ be a function. Then $f$ is $\mathsf{S}$-measurable if and only if there exists a sequence $(s_k)$ of $\mathsf{S}$-simple functions ($s_k:X \to \mathbb{R}$) such that, for every $x\in X$,
$$
f(x)=\lim_{k \to \infty} s_{k}(x).
$$
\end{theorem}

\begin{proof}
$\Rightarrow):$ Suppose that $f$ is $\mathsf{S}$-measurable. Theorem \ref{331} ensures that the sets $A_{+\infty}(f),A_{-\infty}(f) \in \mathsf{S}$ and that the function $f_0:X \to \mathbb{R}$ given by
$$
f_{0}(x):=\left\{
\begin{array}{l c l}
f(x) & & \mbox{if}\,\,x \not \in A_{+\infty}(f) \cup A_{-\infty}(f),\\
0 & &\mbox{if}\,\,x  \in A_{+\infty}(f) \cup A_{-\infty}(f),
\end{array}
\right.
$$
is $\mathsf{S}$-measurable.

Observe that $f(x)=f_0(x)$ if $x \not \in A_{+\infty}(f) \cup A_{-\infty}(f)$. Consider the following cases:

{\scshape Case 1.}\quad $f$ is nonnegative.

Clearly, $f_0(x) \geq 0$ for every $x \in X$. Applying Lemma \ref{318} to $f_0$, we conclude that there exists a nondecreasing sequence $(s_{0,k})$ of $\mathsf{S}$-simple functions such that $s_{0,k}(x) \to f_0(x)$ for every $x\in X$.

For each $k \in \mathbb{N}$ define $s_k:X \to \mathbb{R}$ by
$$
s_k(x):=\left\{
\begin{array}{l c l}
s_{0,k}(x) & & \mbox{if}\,\,x \not \in A_{+\infty}(f),\\
k & &\mbox{if}\,\,x  \in A_{+\infty}(f),
\end{array}
\right.
$$
which is clearly an $\mathsf{S}$-simple function. Therefore, $(s_k)$ is a nondecreasing sequence of $\mathsf{S}$-simple functions converging pointwise on $X$, since
$$
\begin{aligned}
\lim_{k \to \infty} s_k(x)&=\lim_{k \to \infty}s_{0,k}(x)=f_0(x)=f(x)\quad \mbox{if}\,\,\,x \not\in A_{+\infty}(f),\\
\lim_{k \to \infty} s_k(x)&=\lim_{k \to \infty} k=+\infty=f(x)\quad\mbox{if}\,\,\,x \in A_{+\infty}(f).\\
\end{aligned}
$$

This proves the result.

{\scshape Case 2.}\quad $f$ is an arbitrary function.

Consider the nonnegative $\mathsf{S}$-measurable functions $f^{+}$ and $f^{-}$. Observe that if $f(x_{0})=+ \infty$, then $f^{+}(x_{0})=+\infty$ and $f^{-}(x_0) = 0$. Analogously, if $f(x_{1})=-\infty$, then $f^{+}(x_1)=0$ and $f^{-}(x_1)=+\infty$. Consequently, the expression $f^{+}-f^{-}$ is well defined in $\overline{\mathbb{R}}$ and equals $f$. Applying {\scshape Case 1} to the functions $f^{+}$ and $f^{-}$, we conclude that there exist sequences of $\mathsf{S}$-simple functions $(r_k)$ and $(t_k)$ such that, for every $x\in X$,
$$
\lim_{k \to \infty} r_k(x)=f^{+}(x)\quad \mbox{and}\quad \lim_{k \to \infty} t_k(x)=f^{-}(x).
$$

Thus, defining the sequence of $\mathsf{S}$-simple functions $(s_k):=(r_k-t_k)$, it follows that for every $x\in X$,
$$
\lim_{k \to \infty} s_k(x)= \lim_{k \to \infty} r_k(x) - \lim_{k \to \infty} t_k(x) = f^{+}(x)-f^{-}(x)=f(x),
$$
which proves the result.

$\Leftarrow):$ There exists a sequence $(s_k)$ of $\mathsf{S}$-simple functions such that $\displaystyle\lim_{k \to \infty} s_k(x)=f(x)$ for every $x \in X$. Then, for every $c \in \mathbb{R}$,
$$
\begin{aligned}
f(x) > c &\quad \Leftrightarrow\quad f(x) - c > 0\\
&\quad \Leftrightarrow\quad\mbox{there exists}\,\, r\in \mathbb{N}\,\,\mbox{such that}\,\, f(x)-c>\frac{2}{r}\\
&\quad \Leftrightarrow\quad\mbox{there exist}\,\, r\in \mathbb{N}\,\,\mbox{and}\,\,k=k(x,r)\in \mathbb{N}\,\,\mbox{such that}\,\, s_n(x)-c>\frac{1}{r}\,\,\mbox{for every}\,\,n \geq k\\
&\quad \Leftrightarrow\quad x \in  \bigcup_{r=1}^{\infty} \bigcup_{k=1}^{\infty} \bigcap_{n=k}^{\infty}\left\{x\in X\,:\, s_n(x) > c + \frac{1}{r} \right\}.
\end{aligned}
$$

Consequently,
$$
f^{-1}((c,+\infty])=\{x \in X\,:\, f(x)> c \} = \bigcup_{r=1}^{\infty} \bigcup_{k=1}^{\infty} \bigcap_{n=k}^{\infty}\left\{x\in X\,:\, s_n(x) > c + \frac{1}{r} \right\} \in \mathsf{S}
$$
by Proposition \ref{330}. That is, $f$ is $\mathsf{S}$-measurable.
\end{proof}

\begin{remark}\label{333}
The converse implication in the previous theorem could be proved in a way analogous to the real-valued case, by establishing the measurability of the functions appearing in Theorem \ref{320} while considering, in general, that the sequence $(f_k)$ of $\mathsf{S}$-measurable functions is not bounded and taking into account Definition \ref{326}. We propose this here as an exercise {\rm [Exercise \ref{E323}]}.
\end{remark}

Consider two functions $f,g:X \to \overline{\mathbb{R}}$ taking values in the extended real numbers. Suppose that there exist $x_{0},x_{1} \in X$ such that
$$
f(x_0)=+\infty \quad \mbox{and}\quad g(x_0)=-\infty,
$$
$$
f(x_1)=-\infty \quad \mbox{and}\quad g(x_1)=+\infty.
$$

In these cases the expression $(f+g)(x)$ is not well defined, so it becomes necessary to specify its definition in order to work with algebraic combinations of such functions. Thus, if
$$
x \in \left(A_{+\infty}(f) \cap A_{-\infty}(g) \right) \cup \left(A_{-\infty}(f) \cap A_{+\infty}(g)  \right)
$$
we define $(f+g)(x):=0$. In any other case we define
$$
(f+g)(x):=f(x)+g(x).
$$

With this convention, the reader should not conclude that we have now defined expressions that were previously undefined (see Definition \ref{325}).

\begin{theorem} \label{334}
Let $(X,\mathsf{S})$ be a measurable space, let $f,g :X \to \overline{\mathbb{R}}$ be $\mathsf{S}$-measurable functions, and let $\gamma \in \mathbb{R}$. Then $\gamma f$, $f+g$, and $fg$ are $\mathsf{S}$-measurable.
\end{theorem}

\begin{proof}
By Theorem \ref{332}, there exist two sequences $(s_k)$ and $(t_k)$ of $\mathsf{S}$-simple functions such that $s_k \to f$ and $t_k \to g$ pointwise on $X$. Defining the sequence $(\gamma s_k)$ of $\mathsf{S}$-simple functions, we have $\gamma s_k \to \gamma f$, and therefore $\gamma f$ is $\mathsf{S}$-measurable.

For $f+g$, let us write
$$
A:=\left(A_{+\infty}(f) \cap A_{-\infty}(g) \right) \cup \left(A_{-\infty}(f) \cap A_{+\infty}(g)  \right) \in \mathsf{S}.
$$

Then $\chi_{X\smallsetminus A}$ is $\mathsf{S}$-measurable by Proposition \ref{313}. Thus, the sequence of $\mathsf{S}$-simple functions $(s_k + t_k)\chi_{X\smallsetminus A}$ converges to $f+g$ by the previous remark. Therefore $f+g$ is $\mathsf{S}$-measurable.

Finally, for $fg$, consider $A_{0}(f)=\{ x \in X\,:\, f(x)=0\}$ and $A_{0}(g)=\{x \in X\,:\, g(x)=0\}$. Then the sequences of $\mathsf{S}$-simple functions $(s_k\,\chi_{X\smallsetminus A_0(f)})$ and $(t_k\,\chi_{X\smallsetminus A_0(g)})$ satisfy
$$
(s_k\,\chi_{X\smallsetminus A_0(f)})\cdot(t_k\,\chi_{X\smallsetminus A_0(g)}) \to fg \quad \text{pointwise on }X,
$$
and therefore $fg$ is $\mathsf{S}$-measurable.
\end{proof}

If we now consider a function $f:X \to \overline{\mathbb{R}}$, then in the expression $\frac{1}{f}$ we encounter a situation similar to the previous one whenever there exist $x \in X$ such that $f(x)=-\infty$, $0$, or $+\infty$. Following the observations made above, for such points we define $\left( \frac{1}{f} \right)(x):=0$, and in this way we may establish its $\mathsf{S}$-measurability whenever $f$ is an $\mathsf{S}$-measurable function [Exercise \ref{E318}].

\section{Measurable functions with complex values}

Recall that a complex number $z$ is of the form
$$
z = x + iy, \qquad x,y \in \mathbb{R},
$$
where $x$ is called the \textbf{real part} of $z$ and $y$ the \textbf{imaginary part}. We denote these quantities by $\Re(z)$ and $\Im(z)$, respectively. The set of complex numbers is denoted by $\mathbb{C}$.

Moreover, $\mathbb{C}$ may be regarded as a vector space over itself, with the operations of addition and scalar multiplication given by
$$
\begin{aligned}
z + \omega &:= (\Re(z)+\Re(\omega)) + i(\Im(z)+\Im(\omega)), \\
\gamma\, z &:= \big(\Re(\gamma)\Re(z)-\Im(\gamma)\Im(z)\big)
+ i\big(\Re(\gamma)\Im(z)+\Im(\gamma)\Re(z)\big),
\end{aligned}
\qquad \forall\, z,\omega,\gamma \in \mathbb{C}.
$$

The space $\mathbb{C}$ is also a normed space with the norm induced by the complex modulus, defined by
$$
|z| := \sqrt{|\Re(z)|^{2}+|\Im(z)|^{2}}.
$$

Thus, if we identify $\mathbb{C}$ with $\mathbb{R}^{2}$ through the
bijection $\psi:\mathbb{C}\to\mathbb{R}^{2}$ given by
$$
\psi(z):=(\Re(z),\Im(z)),
$$
it can be shown that $\mathbb{C}$ and $\mathbb{R}^{2}$ are equivalent as
metric spaces; indeed, $\psi$ is Lipschitz continuous and so is its inverse. Consequently, the Borel $\sigma$-algebra on $\mathbb{C}$,
denoted by $\mathcal{B}(\mathbb{C})$, may be defined analogously to the
Borel $\sigma$-algebra on $\mathbb{R}^{2}$ [Exercise \ref{E327}].\index{sigma@$\sigma$!algebra!Borel on $\mathbb{C}$}

\begin{definition} \label{335} \index{sigma@$\sigma$!algebra!Borel on $\mathbb{C}$}
Let $(X,\mathsf{S})$ be a measurable space and let $f:X \to \mathbb{C}$ be a function. We say that $f$ is \textbf{$\mathsf{S}$-measurable} if $f^{-1}(\mathcal{B}(\mathbb{C})) \subset \mathsf{S}$. \index{function!$\mathsf{S}$-measurable}
\end{definition}

Given a function $f:X \to \mathbb{C}$, we denote by
$$
(\Re\,f)(x):=\Re(f(x))\quad \mbox{and}\quad (\Im\,f)(x):=\Im(f(x))
$$
the real and imaginary parts of $f$, respectively.

The following result establishes a characterization of $\mathsf{S}$-measurable complex-valued functions in terms of their real and imaginary parts.

\begin{theorem} \label{336}
Let $(X,\mathsf{S})$ be a measurable space and let $f:X \to \mathbb{C}$ be a function. The following conditions are equivalent:
\begin{itemize}
\item[(a)] $f$ is $\mathsf{S}$-measurable.
\item[(b)] For every $c,d \in \mathbb{R}$, the set $\{ x\in X\,:\, (\Re f)(x)>c,\,\,(\Im f)(x)>d \}$ belongs to $\mathsf{S}$.
\item[(c)] $\Re f$ and $\Im f : X \to \mathbb{R}$ are $\mathsf{S}$-measurable.
\end{itemize}
\end{theorem}

\begin{proof}
$(a)\Rightarrow(b):$ Let $c,d \in \mathbb{R}$. Since $(c,+\infty) \times (d,+\infty) \in \mathcal{B}(\mathbb{C})$, then
$$
f^{-1}((c,+\infty) \times (d,+\infty))=\{ x\in X\,:\, (\Re f)(x)>c,\,\,(\Im f)(x)>d \} 
$$
belongs to $\mathsf{S}$.

\begin{figure}[ht!]
\begin{minipage}[l]{0.5\textwidth}
\begin{center}
	\begin{tikzpicture}[xscale=0.55,yscale=0.55]
\draw[fill=lime!25, dotted] (4,4)--(0.6,4)--(0.6,-4)--(4,-4);
\draw[->, gray] (-1,0)--(4,0); 
\draw[->, gray] (0,-4)--(0,4); 


\draw (4,0) node[right]{$_{\Re(f)}$};
\draw (0,4) node[left]{$_{\Im(f)}$};

\draw (0.6,0) node{$_{|}$}; \draw (0.45,0) node[below]{$_{c}$};

\end{tikzpicture}
\end{center}
\end{minipage} \hfill 
 \begin{minipage}[r]{0.5\textwidth}
\begin{center}
	\begin{tikzpicture}[xscale=0.55,yscale=0.55]
\draw[fill=lime!25, dotted] (4,4)--(0.6,4)--(0.6,-4)--(4,-4);
\draw[fill=blue!25, dotted] (4,4)--(0.6,4)--(0.6,-0.5)--(4,-0.5);
\draw[fill=blue!20, dotted] (4,4)--(0.6,4)--(0.6,-1)--(4,-1);
\draw[fill=blue!15, dotted] (4,4)--(0.6,4)--(0.6,-1.5)--(4,-1.5);
\draw[fill=blue!10, dotted] (4,4)--(0.6,4)--(0.6,-2)--(4,-2);

\draw[->, gray] (-1,0)--(4,0); 
\draw[->, gray] (0,-4)--(0,4); 


\draw (4,0) node[right]{$_{\Re(f)}$};
\draw (0,4) node[left]{$_{\Im(f)}$};

\draw (0.6,0) node{$_{|}$}; \draw (0.45,0) node[below]{$_{c}$};
\draw (0,-0.5) node{$_{-}$}; \draw (0,-0.5) node[left]{$_{-1}$};
\draw (0,-1) node{$_{-}$}; \draw (0,-1) node[left]{$_{-2}$};
\draw (0,-1.5) node{$_{-}$}; \draw (0,-1.5) node[left]{$_{-3}$};
\draw (0,-2) node{$_{-}$}; \draw (0,-2) node[left]{$_{-4}$};
\draw (0,-2.5) node[left]{$_{\vdots}$};
\draw (2.5,-2.5) node{$_{\vdots}$};

\draw[-,dotted] (0.6,-0.5)--(4,-0.5);
\draw[-,dotted] (0.6,-1)--(4,-1);
\draw[-,dotted] (0.6,-1.5)--(4,-1.5);
\draw[-,dotted] (0.6,-2)--(4,-2);
\end{tikzpicture}
\end{center}
\end{minipage}
\end{figure}

$(b)\Rightarrow(c):$ For every $c,d \in \mathbb{R}$ we have
$$
\begin{aligned}
\Re f^{-1}((c,+\infty))=\{x \in X\,:\, \Re f (x) > c \}=\bigcup_{k=1}^{\infty} \{x \in X\,:\, (\Re f)(x)>c,\,\,(\Im f)(x)>-k \}
\end{aligned}
$$
which belongs to $\mathsf{S}$, and
$$
\begin{aligned}
\Im f^{-1}((d,+\infty))=\{x \in X\,:\, \Im f (x) > d \}=\bigcup_{k=1}^{\infty} \{x \in X\,:\, (\Re f)(x)>-k,\,\,(\Im f)(x)>d \}
\end{aligned}
$$
also belongs to $\mathsf{S}$.

Therefore, $\Re f$ and $\Im f$ are $\mathsf{S}$-measurable.

$(c)\Rightarrow (a):$ Let $\mathcal{C}:=\{(c,+\infty) \times (d,+\infty)\,:\, c,d \in \mathbb{R} \}$. Then,
$$
\begin{aligned}
f^{-1}((c,+\infty) \times (d,+\infty))&=\{x \in X\,:\, (\Re f(x),\Im f(x)) \in(c,+\infty) \times (d,+\infty)\}\\
&=\{x \in X\,:\, \Re f (x)>c \} \,\cap\,\{x \in X\,:\, \Im f(x)>d\}\\
&=\Re f^{-1}((c,+\infty)) \cap \Im f^{-1} ((d,+\infty)) \in \mathsf{S}
\end{aligned}
$$
which implies that $f^{-1}(\mathcal{C}) \subset \mathsf{S}$. Theorem \ref{33} ensures that $f$ is $\mathsf{S}$-measurable since $\sigma(\mathcal{C})=\mathcal{B}(\mathbb{C})$ [Exercise \ref{E229}].
\end{proof}

\begin{figure}[ht!]
\centering
\begin{tikzpicture}[xscale=0.55,yscale=0.55]
\draw[fill=lime!20, dotted] (4,4)--(0.6,4)--(0.6,-4)--(4,-4);
\draw[fill=blue!15, dotted] (-4,4)--(-4,0.5)--(4,0.5)--(4,4);
\draw[fill=gray!21, dotted] (4,4)--(0.6,4)--(0.6,0.5)--(4,0.5);

\draw[->, gray] (-4,0)--(4,0); 
\draw[->, gray] (0,-4)--(0,4); 


\draw (4,0) node[right]{$_{\Re(f)}$};
\draw (0,4) node[left]{$_{\Im(f)}$};

\draw (0.6,0) node{$_{|}$};  \draw (0.45,0) node[below]{$_{c}$};
\draw (0,0.5) node{$_{-}$}; \draw (0,0.5) node[left]{$_{d}$};

\end{tikzpicture}
\end{figure}

The equivalence between parts \textit{(a)} and \textit{(c)} was proved, in a different way, in Example \ref{311} when identifying $\mathbb{C}$ with $\mathbb{R}^{2}$ endowed with the usual topology.

Algebraic operations and complex conjugation of $\mathsf{S}$-measurable functions produce $\mathsf{S}$-measurable functions. The statements of these results, as well as their proofs, are left as exercises [Exercise \ref{E328}].

Before concluding this chapter, let us observe the following. For every complex number $z \in \mathbb{C}$, the inequalities
$$
|\Re(z)|,\,|\Im(z)| \le |z| \le |\Re(z)|+|\Im(z)|
$$
hold.

These inequalities allow us to characterize the convergence of a sequence of complex numbers $(z_k)$ in terms of the convergence of the real sequences $\big(\Re(z_k)\big)$ and $\big(\Im(z_k)\big)$. Consequently, it is possible to characterize the $\mathsf{S}$-measurability of a function $f:X \to \mathbb{C}$ through the pointwise limit of a sequence of $\mathsf{S}$-measurable functions $f_k:X \to \mathbb{C}$ [Exercise \ref{E329}].
 
\section{Exercises}

\begin{exercise} \label{E31}
Prove that every monotone function $f:\mathbb{R} \to \mathbb{R}$ is Borel-measurable.
\end{exercise}

\begin{exercise} \label{E32}
Let $X$ be a fixed nonempty set.
\begin{itemize}
    \item[1.] Prove that if $f:X \to \mathbb{R}$ is a function, then the functions $f^{+}:X \to \mathbb{R}$ and $f^{-}:X \to \mathbb{R}$ defined by $f^{+}(x):=\max\{0,f(x)\}$ and $f^{-}(x)=\max\{-f(x),0\}$ are nonnegative and satisfy $f=f^{+}-f^{-}$ and $|f|=f^{+}+f^{-}$.
    \item[2.] Let $(X,\mathsf{S})$ be a measurable space. Prove that the function $f^{-}:X \to \mathbb{R}$ defined by $f^{-}:=\max\{-f,0\}$, associated with an $\mathsf{S}$-measurable function $f:X \to \mathbb{R}$, is also $\mathsf{S}$-measurable.
    \item[3.] Is it true that if $f^{+}$ and $f^{-}$ are $\mathsf{S}$-measurable, then $f$ is $\mathsf{S}$-measurable? Justify your answer in detail.
\end{itemize}
\end{exercise}

\begin{exercise} \label{E33}
Let $A,B \subset X$. Prove that:
\begin{itemize}
\item[(a)] $A \subset B$ if and only if $\chi_{A} \leq \chi_{B}$.
\item[(b)] $\chi_{A\cap B}=\chi_{A}\cdot \chi_{B}=\min(\chi_{A},\chi_{B})$.
\item[(c)] $\chi_{A \cup B}=\chi_{A}+\chi_{B}-\chi_{A}\cdot\chi_{B}=\max(\chi_{A},\chi_{B})=1-(1-\chi_{A})(1-\chi_{B})$.
\item[(d)] $\chi_{X\smallsetminus A}=1-\chi_{A}$.
\item[(e)] $\chi_{A \bigtriangleup B}=|\chi_{A}-\chi_{B}|=\chi_{A}+\chi_{B}-2\cdot\chi_{A}\cdot\chi_{B}$.
\end{itemize}
\end{exercise}

{\setlength{\parindent}{0pt}
\begin{exercise} \label{E34}
Provide an example of a function $f:\mathbb{R} \to \mathbb{R}$ that is not $\mathcal{B}(\mathbb{R})$-measurable.
\end{exercise}

(Hint: In Chapter 1 we mentioned that there exists a subset $V$ (Vitali) of $\mathbb{R}$ such that $V$ is not a Borel set).}

\begin{exercise} \label{E35}
Let $(X,\mathsf{S})$ be a measurable space. Prove that if $s,t:X \to \mathbb{R}$ are $\mathsf{S}$-simple functions, then $s+t$, $\gamma \,s$ with $\gamma \in \mathbb{R}$, $st$, and $\frac{1}{s}$ (if $s \neq 0$) are also $\mathsf{S}$-simple functions.
\end{exercise}

\begin{exercise} \label{E36}
\begin{itemize}
    \item[(1)] Prove that the identity function $id:(X,\mathsf{S}) \to (\mathbb{R}, \mathcal{B}(\mathbb{R}))$, $id(x)=x$, is not $\mathsf{S}$-measurable when $X=\{1,2,3\}$ and $\mathsf{S}=\{\varnothing, \{1\}, \{2,3\}, X \}$.
     \item[(2)] Let $X$ be a fixed nonempty set and let $\mathsf{S},\mathsf{T}$ be two $\sigma$-algebras of subsets of $X$. Prove that $id:(X,\mathsf{S}) \to (X,\mathsf{T})$ is measurable relative to $\mathsf{S}$ and $\mathsf{T}$ if and only if $\mathsf{T} \subset \mathsf{S}$.
\end{itemize}
\end{exercise}

\begin{exercise} \label{E37}
Consider the measurable space $(X,\mathsf{S})$ with $\mathsf{S}=\{ \varnothing , X \}$. Prove that a function $f:X \to \mathbb{R}$ is $\mathsf{S}$-measurable if and only if $f$ is constant.
\end{exercise}

\begin{exercise} \label{E38}
Consider the measurable space $(X,\mathsf{S})$ with $X=\{-1,0,1 \}$ and the $\sigma$-algebra $\mathsf{S}=\{\varnothing, \{0\}, \{-1,1\}, X \}$. Let $id:(X,\mathsf{S}) \to (\mathbb{R},\mathcal{B}(\mathbb{R}))$ be the identity function. Prove that $|id|$ and $id^{2}$ are $\mathsf{S}$-measurable functions, but that $id$ is not.
\end{exercise}

\begin{exercise} \label{E39}
Let $X$ be an arbitrary nonempty set, let $\mathcal{A}=\{A_1,\ldots,A_N\}$ be a finite partition of $X$, and let $f:X \to \mathbb{R}$ be a function. Prove that $f$ is $\sigma(\mathcal{A})$-measurable if and only if $f$ is constant on each element of the partition $\mathcal{A}$.

Therefore, every function measurable with respect to $\sigma(\mathcal{A})$ takes at most $N$ distinct values.
\end{exercise}

\begin{exercise} \label{E310}
Let $(X,\mathsf{S})$ be a measurable space, let $f,g : X \to \mathbb{R}$ be $\mathsf{S}$-measurable functions, and let $A \in \mathsf{S}$. Prove that the function $h:X \to \mathbb{R}$ defined by
$$
h(x):=\left\{
\begin{array}{ccl}
f(x) & & \mbox{if}\,\,\, x \in A,\\
g(x) & & \mbox{if}\,\,\, x \notin A,
\end{array}
\right.
$$
is $\mathsf{S}$-measurable.
\end{exercise}

\begin{exercise} \label{E311}
\begin{itemize}
    \item[(a)] Let $X=(X,d)$ be a metric space and let $f:X \to \overline{\mathbb{R}}$ be a function. Prove that $f$ is lower semicontinuous if and only if $f^{-1}(c,+\infty)$ is open in $X$ for every $c \in \mathbb{R}$. Conclude that every lower semicontinuous function is $\mathcal{B}(X)$-measurable.
    \item[(b)] Let $X=(X,d)$ be a metric space and let $f:X \to \overline{\mathbb{R}}$ be a function. Prove that $f$ is upper semicontinuous if and only if $f^{-1}(-\infty,c)$ is open in $X$ for every $c \in \mathbb{R}$. Conclude that every upper semicontinuous function is $\mathcal{B}(X)$-measurable.
    \item[(c)] Let $(X,\tau_{X})$ and $(Y,\tau_{Y})$ be two topological spaces. Prove that if $f:X \to Y$ is continuous, then $f$ is measurable relative to $\mathcal{B}(X)$ and $\mathcal{B}(Y)$.
\end{itemize}
\end{exercise}

{\setlength{\parindent}{0pt}
\begin{exercise} \label{E312}
Let $f:\mathbb{R} \to \mathbb{R}$ be differentiable. Prove that both $f$ and $f':\mathbb{R} \to \mathbb{R}$ are Borel-measurable functions.
\end{exercise}

(Hint: $f'(x)=\displaystyle\lim_{k \to \infty} k \left[ f(x + \frac{1}{k})-f(x)  \right]$ for $x \in \mathbb{R}$. Prove that $f_{a}(x):=f(x+a)$ is Borel-measurable for each fixed $a \in \mathbb{R}$).}

\begin{exercise} \label{E313}
Let $\mathfrak{F} \subset \mathcal{P}(\mathbb{R})$ be a $\sigma$-algebra of subsets of $\mathbb{R}$. Given the measurable spaces $(\mathbb{R},\mathfrak{F})$ and $(\mathbb{R},\mathcal{B}(\mathbb{R}))$, prove that $\mathcal{B}(\mathbb{R}) \subset \mathfrak{F}$ if and only if every continuous function $f: (\mathbb{R},\mathfrak{F}) \to (\mathbb{R},\mathcal{B}(\mathbb{R}))$ is $\mathfrak{F}$-measurable.

That is, the Borel $\sigma$-algebra is the smallest $\sigma$-algebra that makes continuous functions measurable.
\end{exercise}

\begin{exercise} \label{E314}
Let $f,g:X \to \mathbb{R}$ be functions. Prove that $f(x)+g(x)<c$ if and only if there exist rational numbers $r_{x}$ and $s_{x}$ such that $r_{x}+s_{x}<c$, $f(x)<r_{x}$, and $g(x)<s_{x}$. Let $(X,\mathsf{S})$ be a measurable space. Hence prove, without using simple functions, that $f+g$ is $\mathsf{S}$-measurable whenever $f$ and $g$ are $\mathsf{S}$-measurable.
\end{exercise}

\begin{exercise} \label{E315}
Let $X$ be a fixed nonempty set and let $(f_k)$ be a sequence of functions from $X$ into $\mathbb{R}$. Define
$$
A:=\{x \in X\,:\, (f_k(x))\,\,\mbox{converges in }\,\mathbb{R} \}.
$$

Prove that
$$
A= \bigcap_{n=1}^{\infty} \bigcup_{k=1}^{\infty} \bigcap_{p=1}^{\infty} \left\{ x \in X\,:\, |f_k(x)-f_{k+p}(x)| < \frac{1}{n}\right\}.
$$

Conclude that if $(f_k)$ is a sequence of $\mathsf{S}$-measurable functions, then $A \in \mathsf{S}$.
\end{exercise}

\begin{exercise} \label{E316}
Let $X$ be a fixed nonempty set and let $(f_k)$ be a sequence of functions from $X$ into $\mathbb{R}$. If $a \in \mathbb{R}$ is given and
$$
B=\{x \in X\,:\, f_k(x)\nrightarrow a \},
$$
prove that
$$
B=\bigcup_{n=1}^{\infty} \bigcap_{k=1}^{\infty} \bigcup_{\ell=k}^{\infty}\left\{ x\in X\,:\, |f_{\ell}(x)-a|>\frac{1}{n}\right\}.
$$

Conclude that if $(f_k)$ is a sequence of $\mathsf{S}$-measurable functions, then $B \in \mathsf{S}$.
\end{exercise}

\begin{exercise} \label{E317}
Let $(x_k)$ and $(y_k)$ be sequences of real numbers. Prove the following statements:
\begin{itemize}
    \item[(a)] $\liminf_{k \to \infty} x_{k}$ and $\limsup_{k \to \infty} x_{k}$ always exist in $\overline{\mathbb{R}}$ and satisfy
    $$\liminf_{k \to \infty} x_{k} \leq \limsup_{k \to \infty} x_{k}.$$
    
    \item[(b)] $(x_k)$ converges in $\overline{\mathbb{R}}$ if and only if
    $$\liminf_{k \to \infty} x_k = \limsup_{k \to \infty} x_k.$$
    
    \item[(c)] $\liminf_{k \to \infty}x_{k} = \ell_{\ast} \in \mathbb{R}$ if and only if for every $\varepsilon >0$, the set
    $$\{k \in \mathbb{N}\,:\,x_{k}<\ell_{\ast}+\varepsilon\}$$
    is infinite and the set
    $$\{k \in \mathbb{N}\,:\,x_{k}<\ell_{\ast}-\varepsilon\}$$
    is finite.
    
    \item[(d)] $\limsup_{k \to \infty}x_{k} = \ell^{\ast} \in \mathbb{R}$ if and only if for every $\varepsilon >0$, the set
    $$\{k \in \mathbb{N}\,:\,x_{k}>\ell^{\ast}-\varepsilon\}$$
    is infinite and the set
    $$\{k \in \mathbb{N}\,:\,x_{k}>\ell^{\ast}+\varepsilon\}$$
    is finite.
    
    \item[(e)] If $x_k \leq y_k$ for every $k \in \mathbb{N}$, then
    $$\liminf_{k \to \infty}x_k \leq \liminf_{k\to \infty}y_k \quad \mbox{and} \quad \limsup_{k \to \infty}x_k \leq \limsup_{k\to \infty}y_k.$$
    
    \item[(f)] The lower limit is superadditive, that is,
    $$\liminf_{k \to \infty}(x_{k} + y_{k}) \geq \liminf_{k \to \infty} x_{k} + \liminf_{k \to \infty} y_{k}.$$ 

    Give an example where the inequality is strict.
    
     \item[(g)] The upper limit is subadditive, that is,
    $$\limsup_{k \to \infty} (x_{k}+y_{k})\leq \limsup_{k \to \infty}x_{k} + \limsup_{k \to \infty} y_{k}.$$ 
    
    Give an example where the inequality is strict.
    
    \item[(h)] For every $\gamma >0$,
    $$
\liminf_{k \to \infty}(\gamma\,x_{k})=\gamma\liminf_{k \to \infty}x_{k} \quad \text{and}\quad\limsup_{k \to \infty} (\gamma\,x_{k})=\gamma\limsup_{k \to \infty}x_{k}.
    $$
    
    \item[(i)] For every $\gamma <0$,
    $$
\liminf_{k \to \infty}(\gamma\,x_{k})=\gamma\limsup_{k \to \infty}x_{k} \quad \text{and}\quad\limsup_{k \to \infty} (\gamma\,x_{k})=\gamma\liminf_{k \to \infty}x_{k}.
    $$
\end{itemize}
\end{exercise}

\begin{exercise} \label{E318}
Let $(X,\mathsf{S})$ be a measurable space and let $f:X \to \overline{\mathbb{R}}$ be an $\mathsf{S}$-measurable function. Define $\frac{1}{f}:X \to \mathbb{R}$ as follows:
$$
\left( \frac{1}{f} \right)(x):=\left\{
\begin{array}{l c l}
0 & & \mbox{if}\,\,f(x)=0,\pm\infty,\\
&&\\
\frac{1}{f(x)} & & \mbox{if}\,\,f(x) \neq 0,\pm\infty.
\end{array}
\right.
$$

Prove that $\frac{1}{f}$ is $\mathsf{S}$-measurable.
\end{exercise}

{\setlength{\parindent}{0pt}
\begin{exercise} \label{E319}
Let $(X,\mathsf{S})$ be a measurable space and let $f:X \to \overline{\mathbb{R}}$ be a function. Prove that if $f$ is $\mathsf{S}$-measurable, then
$$
f^{-1}(\{c\})=\{x \in X\,:\, f(x)=c \} \in \mathsf{S}
$$
for every $c \in \overline{\mathbb{R}}$.

Is the converse true? Justify your answer in detail.
\end{exercise}

(Hint: Consider $X=\overline{\mathbb{R}}$ endowed with its Borel $\sigma$-algebra. Define a function with respect to a non-Borel-measurable subset $V$ of $[0,1]$ such that $\{x \in X\,:\, f(x)=c \}$ contains at most one point.)}

\begin{exercise} \label{E320}
Let $(X,\mathsf{S})$ be a measurable space and let $A,B \in \mathsf{S}$ be such that $A \cup B =X$. Given the $\sigma$-algebras
$$
\mathsf{S}_{A}:=\{ E \cap A \,:\, E \in \mathsf{S}\}
$$
and
$$
\mathsf{S}_{B}:=\{ F \cap B\,:\, F \in \mathsf{S} \},
$$
prove that if $f_{A}:A \to \overline{\mathbb{R}}$ and $f_{B}:B \to \overline{\mathbb{R}}$ are $\mathsf{S}_{A}$-measurable and $\mathsf{S}_{B}$-measurable functions, respectively, such that
$$
f_{A}(x)=f_{B}(x)\quad \mbox{for every}\quad x \in A \cap B,
$$
then the function $f:X \to \overline{\mathbb{R}}$ defined by
$$
f(x):=\left\{
\begin{array}{lcl}
f_{A}(x) & &\mbox{if}\,\,\, x\in A,\\
f_{B}(x) & &\mbox{if}\,\,\, x\in B,
\end{array}
\right.
$$
is $\mathsf{S}$-measurable.
\end{exercise}

{\setlength{\parindent}{0pt}
\begin{exercise} \label{E321}
Let $(X,\mathsf{S})$ be a measurable space, let $f,h:X \to \overline{\mathbb{R}}$ be $\mathsf{S}$-measurable functions, and let $g:X \to \overline{\mathbb{R}}$ be given. Prove that if $g$ is $\mathsf{S}$-measurable, then $\mbox{\rm mid}(f,g,h)$ \it is $\mathsf{S}$-measurable. Conversely, if $\mbox{\rm mid}(f,g,h)$ is $\mathsf{S}$-measurable for every pair of $\mathsf{S}$-measurable functions $f$ and $h$, then $g$ is $\mathsf{S}$-measurable. Here $\mbox{\rm mid}(f,g,h)$ denotes the middle value among $f,g,$ and $h$.
\end{exercise}

(Hint: First prove that
$$
\mbox{\rm mid}(f,g,h)=\min(\max(f,g),\max(g,h),\max(f,h)) ).
$$
)}

\begin{exercise} \label{E322}
Let $(X,\mathsf{S})$ be a measurable space and let $f:X \to \overline{\mathbb{R}}$ be an $\mathsf{S}$-measurable function. For each positive real number $\beta$, define the $\beta$-truncation of $f$ as the function $f_{\beta}:X \to \mathbb{R}$ given by
$$
f_{\beta}(x):=\left\{
\begin{array}{lcl}
-\beta & &\mbox{if}\,\,\,f(x)<-\beta,\\
f(x) & &\mbox{if}\,\,\,|f(x)| \leq \beta,\\
\beta & &\mbox{if}\,\,\, f(x)>\beta.
\end{array}
\right.
$$

Prove that $f_{\beta}$ is $\mathsf{S}$-measurable.
\end{exercise} 

\begin{exercise} \label{E323}
Let $(X,\mathsf{S})$ be a measurable space and let $f_k:X \to \overline{\mathbb{R}}$, $k \in \mathbb{N}$, be a sequence of $\mathsf{S}$-measurable functions. Prove the results of Theorem \ref{320} for this sequence of functions; that is, prove that the following $\overline{\mathbb{R}}$-valued functions are $\mathsf{S}$-measurable:
\begin{itemize}
\item[(a)] $m_k(x):=\displaystyle\min_{1\leq j \leq k}f_j(x)$ for fixed $k \in \mathbb{N}$;
\item[(b)]  $M_{k}(x):=\displaystyle\max_{1\leq j \leq k} f_j(x)$ for fixed $k \in \mathbb{N}$;
\item[(c)] $\left(\displaystyle\inf_{k \in \mathbb{N}} f_k\right)(x):=\displaystyle\inf_{k \in \mathbb{N}} f_k(x)$;
\item[(d)] $\left(\displaystyle\sup_{k\in \mathbb{N}} f_k\right)(x):=\displaystyle\sup_{k\in \mathbb{N}} f_k(x)$;
\item[(e)] $\left(\displaystyle\liminf_{k \to \infty} f_k\right)(x) :=\displaystyle\sup_{k \geq 1}\inf_{ j \geq k} f_j(x)$;
\item[(f)] $\left(\displaystyle\limsup_{k \to \infty} f_k\right)(x):=\displaystyle\inf_{k \geq 1}\sup_{ j \geq k} f_j(x)$.
\end{itemize}
\end{exercise}

\begin{exercise} \label{E324}
Let $(X,\mathsf{S})$ be a measurable space and let $Y$ be a nonempty set. Given a function $f:X \to Y$, consider the $\sigma$-algebra of subsets of $Y$ defined by
$$
\mathcal{F}(f):=\{F \subset Y\,:\,f^{-1}(F) \in \mathsf{S} \}.
$$

Prove that $f:(X,\mathsf{S}) \to (Y,\mathcal{F}(f))$ is measurable. Moreover, prove that if $\mathsf{T}$ is a $\sigma$-algebra of subsets of $Y$ such that $f:(X,\mathsf{S}) \to (Y,\mathsf{T})$ is measurable, then $\mathsf{T} \subset \mathcal{F}(f)$.

That is, $\mathcal{F}(f)$ is the largest $\sigma$-algebra that makes $f$ measurable.
\end{exercise}

\begin{exercise} \label{E325}
Let $(X,\mathsf{S})$ be a measurable space and let $f:X \to \mathbb{R}$ be a function. Consider the $\sigma$-algebra of subsets of $X$ defined by
$$
\sigma(f):=\{f^{-1}(B)\,:\, B \in \mathcal{B}(\mathbb{R}) \}.
$$

Prove that $\sigma(f) \subset \mathsf{S}$ if and only if $f$ is an $\mathsf{S}$-measurable function.

That is, $\sigma(f)$ is the smallest $\sigma$-algebra that makes $f$ measurable. \index{sigma@$\sigma$!algebra generated!by a function $\sigma(f)$}
\end{exercise}

\begin{exercise} \label{E326}
Let $N >1$. Prove the following statements:
\begin{itemize}
    \item[(a)] The projection onto the $i$-th coordinate is continuous, that is, $\pi_{i}:\mathbb{R}^{N} \to \mathbb{R}$ given by $\pi(\bar{x}):=x_{i}$ is continuous on $\mathbb{R}^{N}$. 
    \item[(b)] $f=(f_1,\ldots,f_N):\mathbb{R}^{N} \to \mathbb{R}^{N}$ is $\mathcal{B}(\mathbb{R}^{N})$-measurable if and only if each function $f_{i}:\mathbb{R}^{N} \to \mathbb{R}$ is $\mathcal{B}(\mathbb{R}^{N})$-measurable for every $i=1,\ldots,N$.
    \item[(c)] Let $(X,\mathsf{S})$ be a measurable space and, for each $i=1,\ldots,N$, let $f_{i}:X \to \mathbb{R}$ be a function. The function $f:X \to \mathbb{R}^{N}$ defined by 
    $$
f(x):=\left(f_{1}(x),\ldots,f_{N}(x) \right)
    $$
    is measurable relative to $\mathsf{S}$ and $\mathcal{B}(\mathbb{R}^{N})$ if and only if each $f_{i}$ is $\mathsf{S}$-measurable for every $i=1,\ldots,N$.
\end{itemize}
\end{exercise}

\begin{exercise} \label{E327}
Prove the following statements:
\begin{itemize}
    \item[(a)] For every $z \in \mathbb{C}$, the following inequalities hold:
    $$
|\Re(z)|,\,|\Im(z)| \le |z| \le |\Re(z)|+|\Im(z)|.
    $$
    
    \item[(b)] The function $\psi:\mathbb{C} \to \mathbb{R}^{2}$ given by
    $$
\psi(z):=\left( \Re(z),\Im(z) \right)
    $$
    is bijective. Also determine its inverse $\psi^{-1}:\mathbb{R}^{2} \to \mathbb{C}$.
    
    \item[(c)] The functions $\psi:\mathbb{C} \to \mathbb{R}^{2}$ and $\psi^{-1}:\mathbb{R}^{2} \to \mathbb{C}$ from the previous item are Lipschitz continuous.
    
    \item[(d)] Conclude that
    $$
\mathcal{B}(\mathbb{C})=\mathcal{B}(\mathbb{R}^{2}).
    $$
\end{itemize}
\end{exercise}

\begin{exercise} \label{E328}
Let $(X,\mathsf{S})$ be a measurable space. For $\mathsf{S}$-measurable functions $f,g:X \to \mathbb{C}$, prove the corresponding results to Theorem \ref{333}, including complex conjugation; that is, prove that:
\begin{itemize}
    \item[(a)] $\gamma\,f:X \to \mathbb{C}$ is $\mathsf{S}$-measurable for every $\gamma \in \mathbb{C}$.
    
    \item[(b)] $f+g:X \to \mathbb{C}$ is $\mathsf{S}$-measurable.
    
    \item[(c)] $f\cdot g:X \to \mathbb{C}$ is $\mathsf{S}$-measurable.
    
    \item[(d)] $\overline{f}:X \to \mathbb{C}$ is $\mathsf{S}$-measurable, where
    $$
\overline{f}(x):=\overline{f(x)}=\Re(f(x))-i\Im(f(x)).
    $$
\end{itemize}
\end{exercise}

\begin{exercise} \label{E329}
Let $(X,\mathsf{S})$ be a measurable space, let $f_k:X \to \mathbb{C}$ be a sequence of functions, and let $f:X \to \mathbb{C}$ be a function. Prove the following statements:
\begin{itemize}
    \item[(a)] $(f_{k})$ converges pointwise to $f$ on $X$ if and only if the sequences of functions $(\Re(f_k))$ and $(\Im(f_k))$ converge pointwise to $\Re(f)$ and $\Im(f)$ on $X$, respectively.
    
    \item[(b)] If each $f_k$ is $\mathsf{S}$-measurable for every $k \in \mathbb{N}$ and $(f_k)$ converges pointwise to $f$ on $X$, then $f$ is $\mathsf{S}$-measurable.
    
    \item[(c)] The function $f:X \to \mathbb{C}$ is $\mathsf{S}$-measurable if and only if there exists a sequence of $\mathsf{S}$-measurable functions $(g_k)$, with $g_k:X\to\mathbb{C}$, such that $(g_k)$ converges pointwise to $f$ on $X$.
\end{itemize}
\end{exercise}

\begin{exercise}\label{E330}
 Let $(X,\mathsf{S})$ be a measurable space and let $f:X\to\mathbb{C}$ be a function. Prove that $f$ is $\mathsf{S}$-measurable if and only if the functions
$$
\varphi_1:X\to\mathbb{C},\qquad \varphi_1(x):=e^{i\,\Re(f(x))},
$$
$$
\varphi_2:X\to\mathbb{C},\qquad \varphi_2(x):=e^{i\,\Im(f(x))}
$$
are $\mathsf{S}$-measurable.   
\end{exercise}
\chapter{Measures on \texorpdfstring{$\sigma$}{}-algebras} \label{Capitulo4}
\markboth{{\scriptsize 4. MEASURES ON SIGMA ALGEBRAS}}{ {\scriptsize 4. MEASURES ON SIGMA ALGEBRAS}}

Measure theory provides a general framework for quantifying the ``size'' of sets, extending and unifying notions such as length, area, and volume. Beyond classical geometry, this theory is essential in probability and statistics, mathematical analysis, and numerous applications in science and engineering. At its core, it transforms the geometric intuition of measuring into a rigorous and flexible tool that can be applied even in abstract spaces.

In this first part, we shall work with measures defined on a $\sigma$-algebra. Formally, a measure is a function $\mu$ that assigns to each set in the $\sigma$-algebra a number in $[0,+\infty]$, satisfying certain additivity properties.

Throughout this chapter we shall study the elementary properties of measures and present relevant examples. We shall introduce the concept of a measure space, which includes the notion of null sets. These play an essential role in determining when a measure space is complete, a property that allows the extension of the measure to broader classes of subsets. Likewise, we shall develop the concept of statements that hold almost everywhere and prove the uniqueness theorems for measures.

\section{Definitions and examples}

\begin{definition} \label{41}
Let $(X,\mathsf{S})$ be a measurable space. A \textbf{measure} \index{measure!measure} on $(X,\mathsf{S})$ is a function $\mu:\mathsf{S} \to \overline{\mathbb{R}}$ satisfying the following properties:
\begin{itemize}
\item[\rm(M1)] $\mu(\varnothing)=0$.
\item[\rm(M2)] $\mu(A) \geq 0$ for every $A \in \mathsf{S}$.
\item[\rm(M3)] If $(A_k)$ is a sequence of pairwise disjoint elements of $\mathsf{S}$, then
$$
\mu \left(\bigcup_{k=1}^{\infty} A_k \right)=\sum_{k=1}^{\infty} \mu(A_k).
$$

This last equality is called the \textbf{$\sigma$--additivity} of $\mu$. \index{sigma@$\sigma$!additivity} 
\end{itemize}
\end{definition}

\begin{definition} \label{42} \index{measure!finite} \index{measure!$\sigma$-finite}
Let $(X,\mathsf{S})$ be a measurable space and let $\mu:\mathsf{S} \to \overline{\mathbb{R}}$ be a measure on $(X,\mathsf{S})$. We say that $\mu$ \textbf{is finite} if $\mu(X) < +\infty$.

We say that $\mu$ is \textbf{$\sigma$-finite} if there exists a sequence $(A_k)$ of elements of $\mathsf{S}$ such that $\mu(A_k)<+\infty$ for every $k \in \mathbb{N}$ and
$$
X=\bigcup_{k=1}^{\infty} A_k.
$$
\end{definition}

Let us show that every finite measure is $\sigma$-finite.

\begin{proposition} \label{43}
Let $(X,\mathsf{S})$ be a measurable space and let $\mu:\mathsf{S} \to \overline{\mathbb{R}}$ be a measure on $(X,\mathsf{S})$. If $\mu$ is finite, then $\mu$ is $\sigma$-finite.
\end{proposition}

\begin{proof}
Define the sequence $(A_k)$ of elements of $\mathsf{S}$ by setting $A_k:=X$ if $k$ is even and $A_k:=\varnothing$ if $k$ is odd. Then
$$
X=\bigcup_{k=1}^{\infty} A_{k}
$$
with $\mu(A_k)< +\infty$ for every $k \in \mathbb{N}$. Consequently, $\mu$ is $\sigma$-finite.
\end{proof}

The converse of the previous proposition is not true in general [Exercise \ref{E45}].

Any measure $\mu$ on $(X,\mathsf{S})$ satisfies the following properties, whose proofs follow easily from the definition.

\begin{proposition} \label{44}
Let $(X,\mathsf{S})$ be a measurable space and let $\mu:\mathsf{S} \to \overline{\mathbb{R}}$ be a measure on $(X,\mathsf{S})$.
\begin{itemize}
\item[(1)]{\rm {\scshape Additivity:}} $\mu(A \cup B)=\mu(A)+\mu(B)$ for every $A,B \in \mathsf{S}$ such that $A \cap B = \varnothing$.
\item[(2)]{\rm {\scshape Monotonicity:}} $\mu(A) \leq \mu(B)$ for every $A,B\in \mathsf{S}$ such that $A \subset B$.
\item[(3)]{\rm{\scshape Subtractivity:}}  $\mu(B \smallsetminus A) = \mu(B)-\mu(A)$ for every $A,B\in \mathsf{S}$ such that $A \subset B$ and $\mu(A)<+\infty$.
\item[(4)] $\mu(A \cup B) + \mu(A \cap B) = \mu(A) + \mu(B)$ for every $A,B\in \mathsf{S}$.
\end{itemize}
\end{proposition}

\begin{proof}
\textit{(1):} For every pair of disjoint sets $A,B \in \mathsf{S}$, define the disjoint sequence $(A_k)$ of elements of $\mathsf{S}$ by
$$
\begin{aligned}
A_1 &:=A, \\
A_2 &:=B, \\
A_k &:=\varnothing,\quad \forall k \geq 3.
\end{aligned}
$$

By the $\sigma$-additivity of $\mu$, it follows that
$$
\mu\left( A \cup B \right) = \mu \left( \bigcup_{k=1}^{\infty} A_k \right) = \sum_{k=1}^{\infty} \mu(A_k) =\mu(A) + \mu(B).
$$

\textit{(2):} Let $A,B \in \mathsf{S}$ be such that $A \subset B$. The set $B$ can be written as the disjoint union
$$
B=(B \smallsetminus A)\cup A.
$$

Thus, by the previous item and the nonnegativity of $\mu$, we conclude that
$$
\mu(A) \leq \mu(A) + \mu(B\smallsetminus A) = \mu(B).
$$

\begin{figure}[ht!]
\centering
\begin{tikzpicture}[xscale=1.5, yscale=1.5]
\filldraw[draw=black,fill=blue!15]
plot[smooth cycle] coordinates{
(-1,0)(-0.5,1)(-0.3,0.8)(1,1.25)(1.25,0)(0.75,-0.5)(-0.5,-1)};
\draw[draw=black,fill=red!15] (0,0) circle(0.5);
\draw(0,0)node{$_{A}$};
\draw(0.45,0.6)node{$_{B\smallsetminus A}$};
\draw(1.1,1.35)node{$_{B}$};
\end{tikzpicture}
\end{figure}

\textit{(3):} By the previous item,
$$
\mu(B)=\mu(B\smallsetminus A) + \mu(A)
$$
for every $A,B \in \mathsf{S}$ such that $A \subset B$. Subtracting $\mu(A)<+\infty$ from the previous identity, we obtain
$$
\mu(B\smallsetminus A)=\mu(B)-\mu(A).
$$

\textit{(4):} Let $A,B \in \mathsf{S}$. The sets $A$ and $B$ can be written as the following disjoint unions of elements of $\mathsf{S}$:
$$
A=(A \cap B) \cup (A\smallsetminus B) \quad\mbox{and}\quad B=(A \cap B)\cup(B\smallsetminus A).
$$

\begin{figure}[ht!]
\centering
\begin{tikzpicture}[xscale=1.85, yscale=1.5]
\filldraw[draw=black,fill=blue!15] (-1,-1)--(-1,1)--(2,1)--(2,-1)--(-1,-1);
\draw[fill=lime!15] (0,0) circle(0.75);
\draw[draw=black,fill=lime!15] (1,0) circle(0.75);
\draw[draw=black] (0,0) circle(0.75);

\draw(-0.2,0.25)node{$_{A}$};
\draw(1.2,0.25)node{$_{B}$};
\draw(0.5,0)node{$_{A\cap B}$};
\draw(1.8,0.8)node{$_{X}$};
\end{tikzpicture}
\end{figure}

Thus, by item \textit{(1)} we have
$$
\mu(A)=\mu(A \cap B) + \mu(A \smallsetminus B) \quad \mbox{and}\quad \mu(B)=\mu(A \cap B)+\mu(B\smallsetminus A),
$$
and therefore,
$$
\mu(A)+\mu(B)=\mu(A \cap B) + \mu(A\smallsetminus B) + \mu(A \cap B) + \mu(B\smallsetminus A) = \mu(A \cap B) + \mu(A \cup B)
$$
since
$$
A \cup B=(A \smallsetminus B)\cup (A \cap B) \cup (B \smallsetminus A)
$$
is a disjoint union of elements of $\mathsf{S}$.
\end{proof}

Let us now consider some examples.

\begin{example} \label{45}
Let $(X,\mathsf{S})$ be a measurable space and let $x_{0} \in X$ be fixed. Define $\delta_{x_{0}}:\mathsf{S} \to \{0,1\}$ as follows:
$$
\delta_{x_{0}}(A):=\left\{
\begin{array}{lcl}
1 & & \mbox{if}\,\,\, x_{0}\in A,\\
0 & & \mbox{if}\,\,\, x_{0}\not\in A.
\end{array}
\right.
$$

Then $\delta_{x_{0}}$ is a finite measure on $(X,\mathsf{S})$ called the \textbf{unit mass concentrated at $x_{0}$}. \index{measure!unit mass}
\end{example}

\begin{proof}
Properties (M1) and (M2) are immediate. Let $(A_k)$ be a disjoint sequence of elements of $\mathsf{S}$. If $x_{0} \not\in A_k$ for every $k \in \mathbb{N}$, then $\delta_{x_{0}}(A_k)=0$ for every $k \in \mathbb{N}$ and, therefore,
$$
\delta_{x_{0}}\left(\bigcup_{k=1}^{\infty} A_k\right)=0=\sum_{k=1}^{\infty} \delta_{x_{0}}(A_k).
$$

Otherwise, let $k_{0}$ be the unique natural number such that $x_{0} \in A_{k_{0}}$. Then $\delta_{x_{0}}(A_{k_{0}})=1$ and $\delta_{x_{0}}(A_k)=0$ for every $k \neq k_{0}$. Therefore,
$$
\delta_{x_{0}}\left(\bigcup_{k=1}^{\infty} 
A_k\right)=1=\sum_{k=1}^{\infty} \delta_{x_{0}}(A_k).
$$

Consequently, $\delta_{x_{0}}$ satisfies (M3), and it follows that it is a finite measure since clearly
$$
\delta_{x_{0}}(X)=1.
$$
\end{proof}

\begin{example}\label{46}
Let $X$ be nonempty and let $\mathsf{S}=\mathcal{P}(X)$. Define $\mu^{\sharp}:\mathsf{S} \to \overline{\mathbb{R}}$ as follows:
$$
\mu^{\sharp}(A)=\left\{
\begin{array}{lcl}
|A| & &\mbox{if}\,\,\, A\,\,\mbox{is finite},\\
+\infty & &\mbox{if}\,\,\, A\,\,\mbox{is infinite},
\end{array}
\right.
$$
where $|A|$ denotes the cardinality of the set $A$.

Then $\mu^{\sharp}$ is a measure on $(X,\mathsf{S})$ called the \textbf{counting measure}. \index{measure!counting}
\end{example}

\begin{proof}
Properties (M1) and (M2) are immediate. Let $(A_k)$ be a disjoint sequence of elements of $\mathsf{S}$. Consider the following two cases:

{\scshape Case 1.} There exists $k_{0} \in \mathbb{N}$ such that $A_{k_{0}}$ is infinite. The union $\bigcup_{k=1}^{\infty} A_{k}$ contains $A_{k_{0}}$ and is therefore an infinite set. Consequently,
$$
\mu^{\sharp}\left(\bigcup_{k=1}^{\infty} A_k \right)=+\infty=|A_{k_{0}}|+\sum_{k \neq k_{0}} |A_k|=\sum_{k=1}^{\infty} \mu^{\sharp}(A_k).
$$

{\scshape Case 2.} $A_{k}$ is a finite set for every $k \in \mathbb{N}$.

The result is immediate if $A_{k}=\varnothing$ for every $k \in \mathbb{N}$ or if there exists a finite collection of natural numbers $k_{1},\ldots,k_{N}$ such that $A_{k_{1}},A_{k_{2}},\ldots,A_{k_N}$ are nonempty. Thus, suppose that $A_{k} \neq \varnothing$ for every $k \in \mathbb{N}$. Since $1\leq |A_{k}|$ for every $k \in \mathbb{N}$, then
$$
\sum_{k=1}^{\infty}|A_k| =+\infty.
$$

Moreover, since the sequence $(A_{k})$ is disjoint, it follows that $\bigcup_{k=1}^{\infty} A_k$ is an infinite set, and consequently
$$
\mu^{\sharp}\left(\bigcup_{k=1}^{\infty}A_{k} \right)=+\infty=\sum_{k=1}^{\infty}\mu^{\sharp}(A_k).
$$

From both cases we conclude that $\mu^{\sharp}$ satisfies (M3), and therefore it is a measure.
\end{proof}

In this example it is easy to verify [Exercise \ref{E45}] that $\mu^{\sharp}$ is finite if and only if $X$ is finite, and that $\mu^{\sharp}$ is $\sigma$-finite if and only if $X$ is countable.

\begin{example} \label{47}
Let $X=\mathbb{N}$, $\mathsf{S}=\mathcal{P}(\mathbb{N})$, and let $\overline{x}=(x_j)$ be a sequence of nonnegative extended real numbers. Define $\mu_{\overline{x}}:\mathsf{S} \to \overline{\mathbb{R}}$ as follows:
$$
\mu_{\bar{x}}(A):=\left\{
\begin{array}{lcl}
0 & & \mbox{if}\,\,\,A=\varnothing, \\
\displaystyle\sum_{j \,\in \,A} x_{j} & &\mbox{if}\,\,\,A \neq \varnothing.
\end{array}
\right.
$$

Then $\mu_{\bar{x}}$ is a measure, which is finite if and only if $\sum_{j=1}^{\infty} x_{j} < +\infty,$ and $\sigma$-finite if and only if $x_j < +\infty$ for every $j \in \mathbb{N}$.
\end{example}

\begin{proof}
Properties (M1) and (M2) follow directly from the definition of $\mu_{\bar{x}}$. Let $(A_{k})$ be a disjoint sequence of elements of $\mathsf{S}$. Then
$$
\mu_{\bar{x}}\left(\bigcup_{k=1}^{\infty} A_{k} \right) = \sum_{j \,\in\, \cup_{k \geq 1} A_k } x_j = \sum_{k=1}^{\infty}\left( \sum_{j \,\in\,A_k} x_j \right) = \sum_{k=1}^{\infty} \mu_{\bar{x}}(A_{k}),
$$
that is, $\mu_{\bar{x}}$ satisfies (M3). Therefore, $\mu_{\bar{x}}$ is a measure.

Observe that $\mu_{\bar{x}}$ is finite if and only if $\mu_{\bar{x}}(\mathbb{N})<+\infty$ if and only if $\sum_{j=1}^{\infty} x_j < +\infty.$

Now suppose that $\mu_{\bar{x}}$ is $\sigma$-finite. Then there exists a sequence $(A_{k})$ of subsets of $\mathbb{N}$ such that
$$
X=\bigcup_{k=1}^{\infty} A_k\quad \mbox{with}\quad \mu_{\bar{x}}(A_{k})<+\infty\quad \mbox{for every}\,\, k\in \mathbb{N}.
$$

Suppose there exists $j_{0} \in \mathbb{N}$ such that $x_{j_0}=+\infty$. Since $j_0 \in \mathbb{N}$, then $j_{0} \in A_{k_0}$ for some $k_{0} \in \mathbb{N}$, so that
$$
\mu_{\bar{x}}(A_{k_0})=\sum_{j \in A_{k_0}} x_{j} \geq x_{j_0}=+\infty,
$$
which contradicts our assumption. Therefore, $x_{j} <+\infty$ for every $j \in \mathbb{N}.$

Conversely, if $x_{j} <  +\infty$ for every $j \in \mathbb{N}$, defining the sequence $(A_k)$ of elements of $\mathsf{S}$ by $A_{j}:=\{j\},$
we have
$$
X=\bigcup_{j=1}^{\infty} A_j
$$
and $\mu_{\bar{x}}(A_{j})=x_{j} < +\infty$ for every $j \in \mathbb{N}$. That is, $\mu_{\bar{x}}$ is $\sigma$-finite.
\end{proof}

\begin{example} \label{48}
Let $(X,\mathsf{S})$ be a measurable space, let $\mu:\mathsf{S} \to \overline{\mathbb{R}}$ be a measure, and fix $A \in \mathsf{S}$ with $\mu(A)>0$. Define $\mu^{A}:\mathsf{S} \to \overline{\mathbb{R}}$ as
$$
\mu^{A}(B):=\mu(A \cap B).
$$

Then $\mu^{A}$ is a measure on $(X,\mathsf{S})$ called the \textbf{contraction of $\mu$ to $A$}. \index{measure!contraction} 
\end{example}

\begin{proof}
Clearly $\mu^{A}(\varnothing)=\mu(A \cap \varnothing)=\mu(\varnothing)=0$ and $\mu^{A}(B)=\mu(A \cap B) \geq 0$ for every $B \in \mathsf{S}$ since $\mu$ is a measure. This proves (M1) and (M2).

If $(B_k)$ is a disjoint sequence of elements of $\mathsf{S}$, then
$$
\mu^{A} \left(\bigcup_{k=1}^{\infty} B_k \right)=\mu\left( \bigcup_{k=1}^{\infty} (A \cap B_k) \right)=\sum_{k=1}^{\infty} \mu(A \cap B_k) = \sum_{k=1}^{\infty} \mu^{A}(B_k)
$$
since $\mu$ is a measure and $(A \cap B_k)$ is a disjoint sequence of elements of $\mathsf{S}$. Consequently, $\mu^{A}$ satisfies (M3).

Therefore, $\mu^{A}$ is a measure on $(X,\mathsf{S})$.
\end{proof}

\begin{example} \label{49}
Let $(X,\mathsf{S})$ be a measurable space, let $\mu:\mathsf{S} \to \overline{\mathbb{R}}$ be a measure, and fix $A \in \mathsf{S}$ with $\mu(A)>0$. Define $\mu_{A}:A\cap \mathsf{S} \to \overline{\mathbb{R}}$ as $\mu_{A}(B):=\mu(B)$. Then $\mu_{A}$ is a measure on $A$ called the \textbf{restriction of $\mu$ to $A$}. \index{measure!restriction}
\end{example}

The proof is straightforward and is left as an exercise [Exercise \ref{E48}].

\begin{example} \label{410}
Let $(X,\mathsf{S})$ be a measurable space, let $\mu:\mathsf{S}  \to \overline{\mathbb{R}}$ be a measure on $(X,\mathsf{S})$, and let $f:X \to Y$ be a function. Consider the measurable space $(Y,\mathsf{T})$, where $\mathsf{T}$ is the $\sigma$-algebra of subsets of $Y$ given by
$$
\mathsf{T}:=\{F \subset Y\,:\,f^{-1}(F) \in \mathsf{S} \}.
$$

The function $\nu:\mathsf{T} \to \overline{\mathbb{R}}$ defined by $\nu(F):=\mu(f^{-1}(F))$ is a measure on $(Y,\mathsf{T})$.
\end{example}

\begin{proof}
Since $\mu(\varnothing)=0$ and $\mu(f^{-1}(F)) \geq 0$ for every $f^{-1}(F) \in \mathsf{S}$, it follows that $\nu(\varnothing)=0$ and $\nu(F) \geq 0 $ for every $F \in \mathsf{T}$. That is, $\nu$ satisfies (M1) and (M2).

Let $(F_k)$ be an arbitrary disjoint sequence of elements of $\mathsf{T}$. The sequence $(f^{-1}(F_k))$ is a disjoint sequence of elements of $\mathsf{S}$ and, therefore,
$$
\mu\left(\bigcup_{k=1}^{\infty} f^{-1}(F_k) \right)=\sum_{k=1}^{\infty} \mu(f^{-1}(F_k)).
$$

Consequently,
$$
\begin{aligned}
\nu\left( \bigcup_{k=1}^{\infty} F_k \right)=\mu \left(  f^{-1}\left(\bigcup_{k=1}^{\infty} F_k \right)\right)=\mu\left(\bigcup_{k=1}^{\infty} f^{-1}(F_k) \right)=\sum_{n=1}^{\infty} \mu(f^{-1}(F_k))=\sum_{k=1}^{\infty} \nu(F_k),
\end{aligned}
$$
and it follows that $\nu$ satisfies (M3).
\end{proof}

\begin{example} \label{411} \index{measure!Lebesgue measure on $\mathbb{R}$}
Let $X=\mathbb{R}$ and $\mathsf{S}=\mathcal{B}(\mathbb{R})$. There exists a \textbf{unique} $\sigma$-finite measure defined on a $\sigma$-algebra containing $\mathcal{B}(\mathbb{R})$ and assigning to each interval its length. We denote this measure by $\lambda$ and call it the \textbf{Lebesgue measure\footnote{Henri Léon Lebesgue (1875--1941) was a French mathematician. He was born in Beauvais, Oise, Picardie, France. He studied at the École Normale Supérieure and, during the period 1899--1902, taught at the Lycée of Nancy. In 1910 he obtained a chair at the University of the Sorbonne. Lebesgue is mainly known for his contributions to measure and integration theory. He defined the Lebesgue integral, which generalizes the notion of the Riemann integral by extending the concept of area under a curve to include discontinuous functions. This is one of the achievements of modern analysis that expanded the scope of Fourier analysis. After 1910 he no longer concentrated on the area of study that he had initiated, since his work was a generalization, and he was wary of further abstractions.} on $\mathbb{R}$}. In the following chapters we will focus on the construction and study of this measure.
\end{example}

\begin{figure}[!ht]
\centering
\includegraphics[scale=0.45]{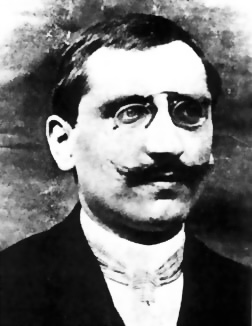} 
\begin{center}
H. Lebesgue (1875-1941)
\end{center}
\end{figure}

\begin{example} \label{412} \index{measure!Lebesgue-Stieltjes measure on $\mathbb{R}$}
Let $X=\mathbb{R}$, $\mathsf{S}=\mathcal{B}(\mathbb{R})$, and let $F:\mathbb{R} \to \mathbb{R}$ be a nondecreasing and right-continuous function. Then there exists a \textbf{unique} measure defined on a $\sigma$-algebra containing $\mathcal{B}(\mathbb{R})$, which we denote by $\overline{\lambda_{F}}$. This measure satisfies that if $E=(a,b)$, then $\overline{\lambda_{F}}(E)=F(b)-F(a)$, and it is called the \textbf{Lebesgue--Stieltjes measure\footnote{Thomas Joannes Stieltjes (1856--1894) was a Dutch mathematician of the nineteenth century who worked on numerous topics, including Gaussian quadrature, orthogonal polynomials, and continued fractions, among others. He is known for the Riemann--Stieltjes integral and the Lebesgue--Stieltjes integral. A graduate of the École Normale Supérieure in Paris in 1886, he was a student of the French mathematicians Charles Hermite and Gaston Darboux.} generated by $F$}. Its construction will be studied later on.
\end{example}

\begin{figure}[ht!]
\centering
\includegraphics[scale=0.225]{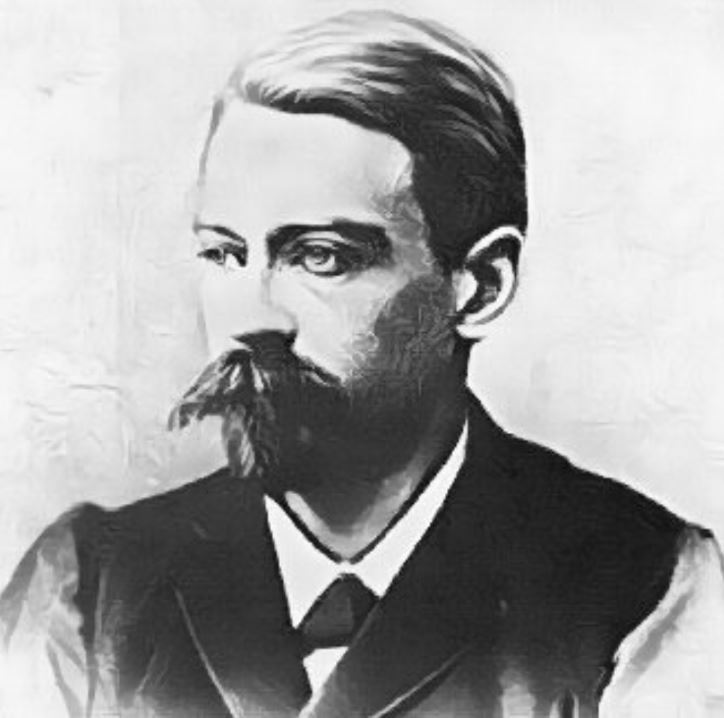} 
\begin{center}
T. Stieltjes (1856-1894)
\end{center}
\end{figure}

If there exists a finite collection of measures on a measurable space, then every nonnegative linear combination of them is again a measure.

\begin{proposition} \label{413}
Let $(X,\mathsf{S})$ be a measurable space, let $\mu_{1},\mu_{2},\ldots,\mu_{n} : \mathsf{S} \to \overline{\mathbb{R}}$ be measures on $(X,\mathsf{S})$, and let $c_{1},c_{2},\ldots,c_{n} \geq 0$ be given. The function $\mu:\mathsf{S} \to \overline{\mathbb{R}}$ defined by
$$
\mu(A):=\sum_{i=1}^{n} c_{i}\,\mu_{i}(A)
$$
is a measure on $(X,\mathsf{S})$.
\end{proposition}

\begin{proof}
It is immediate to verify that $\mu$ satisfies properties (M1) and (M2). Let $(A_k)$ be a disjoint sequence of elements of $\mathsf{S}$ whose union is equal to $A$, that is, $A=\bigcup_{k=1}^{\infty}A_{k}$. By $\sigma$-additivity, for each $i=1,\ldots,n$, we have $\mu_{i}(A)=\sum_{k=1}^{\infty} \mu_{i}(A_k)$. Multiplying both sides of the equality by $c_{i}\geq 0$, it follows that, for each $i=1,\ldots,n$, $c_{i}\mu_{i}(A)=\sum_{k=1}^{\infty} c_{i}\mu_{i}(A_k)$. Therefore,
$$
\mu(A)=\sum_{i=1}^{n} c_{i}\mu_{i}(A)=\sum_{i=1}^{n} \sum_{k=1}^{\infty} c_{i}\mu_{i}(A_k) =\sum_{k=1}^{\infty} \sum_{i=1}^{n} c_{i}\mu_{i}(A_k) = \sum_{k=1}^{\infty} \mu(A_k)
$$
since all terms are nonnegative. Consequently, $\mu$ satisfies (M3).
\end{proof}

\begin{theorem} \label{414}
Let $(X,\mathsf{S})$ be a measurable space and let $\mu:\mathsf{S} \to \overline{\mathbb{R}}$ be a measure on it. 
\begin{itemize}
\item[(a)] If $(A_k)$ is an increasing sequence of elements of $\mathsf{S}$, then
$$
\mu\left( \bigcup_{k=1}^{\infty} A_k \right)=\lim_{k\,\to\,\infty}\mu(A_k).
$$
\item[(b)] If $(B_k)$ is a decreasing sequence of elements of $\mathsf{S}$, then
$$
\mu\left(\bigcap_{k=1}^{\infty} B_k \right) \leq \lim_{k\to\infty} \mu(B_k).
$$

If moreover $\mu(B_k)<+\infty$ for some $k \in \mathbb{N}$, then equality holds.
\end{itemize}
\end{theorem}

\begin{proof}
\textit{(a):} If $\mu(A_k)=+\infty$ for some $k \in \mathbb{N}$, the equality holds trivially since both expressions are equal to $+\infty$, so we may assume that $\mu(A_k)<+\infty$ for every $k \in \mathbb{N}$. 

Let $k \in \mathbb{N}$. Since the sequence $(A_k)$ is increasing, we have $A_{k-1} \subset A_k$ and, therefore, $A_k=(A_k\smallsetminus A_{k-1}) \cup A_{k-1}$. Next, the set $A_{k-2}$ satisfies $A_{k-2} \subset A_{k-1}$ and, therefore, $A_{k-1}=(A_{k-1}\smallsetminus A_{k-2})\cup A_{k-2}$. Continuing in this way, we obtain that, for each $j \in \{1,2,\ldots,k \}$,
$$
A_{j}=(A_j \smallsetminus A_{j-1})\cup A_{j-1}
$$
with $A_{0}:=\varnothing$. Consequently,
$$
A_k = (A_k\smallsetminus A_{k-1})\,\cup (A_{k-1}\smallsetminus A_{k-2}) \,\cup\,\cdots \,\cup (A_2\smallsetminus A_1)\,\cup\,(A_1\smallsetminus A_0).
$$

\begin{figure}[ht!]
\centering
\begin{tikzpicture}[xscale=0.6, yscale=0.6]
\draw[fill=blue!15] (0,0) circle (3);
\draw[fill=red!12] (0,0) circle (1.5);

\draw (0,0) node{${A_{k-1}}$};
\draw (0,3.4) node{${A_{k}}$};
\draw (0,-2.2) node{$A_{k}\smallsetminus A_{k-1}$};
\end{tikzpicture}
\end{figure}

Defining $D_{k}:=A_{k}\smallsetminus A_{k-1}$, the sequence $(D_k)$ is a disjoint sequence of elements of $\mathsf{S}$ such that $A_{k}=D_{1}\,\cup\,D_{2}\,\cup\,\cdots\,\cup\,D_{k}$ for every $k \in \mathbb{N}$ and, therefore,
$$
\bigcup_{k=1}^{\infty} A_{k} = \bigcup_{k=1}^{\infty} D_{k}.
$$

Since $A_{k} \subset A_{k+1}$ for every $k \in \mathbb{N}$, we have $\mu(A_k) \leq \mu(A_{k+1})$ for every $k \in \mathbb{N}$ by the monotonicity of $\mu$. Thus, the sequence $(\mu(A_k))$ is nondecreasing and, consequently, converges in $\overline{\mathbb{R}}$. Applying the $\sigma$-additivity and subtractivity of the measure $\mu$, we obtain
$$
\begin{aligned}
\mu\left(\bigcup_{k=1}^{\infty} A_{k}  \right)&=\mu\left(\bigcup_{k=1}^{\infty} D_{k}  \right)=\sum_{k=1}^{\infty}\mu(D_k)=\sum_{k=1}^{\infty} \mu(A_k\smallsetminus A_{k-1})=\sum_{k=1}^{\infty} (\mu(A_k)-\mu(A_{k-1}))\\
&=\lim_{n \to \infty} \sum_{k=1}^{n} (\mu(A_k)-\mu(A_{k-1}))=\lim_{n \to\infty} \mu(A_n) -\mu(A_0)=\lim_{n\to\infty} \mu(A_n),
\end{aligned}
$$
as stated.

\textit{(b):} Since $\displaystyle\bigcap_{k=1}^{\infty} B_{k} \subset B_{k}$ for every $k \in \mathbb{N}$, then
$$
\mu\left( \bigcap_{k=1}^{\infty} B_{k} \right) \leq \mu(B_k)\quad \mbox{for every}\,\, k \in \mathbb{N}
$$
and, since $(\mu(B_k))$ is decreasing, then
$$
\mu\left( \bigcap_{k=1}^{\infty} B_{k} \right) \leq \lim_{k\to\infty}\mu(B_k).
$$

Suppose that $\mu(B_k)<+\infty$ for some $k \in \mathbb{N}$. Let $k_{0}$ be the first natural number satisfying this property. Since $(B_k)$ is decreasing, $B_{k_{0}+(k+1)} \subset B_{k_{0}+k} \subset B_{k_{0}}$ for every $k \in \mathbb{N}$ and, therefore, $B_{k_0}\smallsetminus B_{k_{0}+k} \subset B_{k_0}\smallsetminus B_{k_{0}+(k+1)}$ for every $k \in \mathbb{N}$. Moreover, for each $j \in \{1,2,\ldots,k_{0} \}$, $B_{k_{0}}\smallsetminus B_{j}=\varnothing$, and therefore
$$
\begin{aligned}
&\,B_{k_0}\smallsetminus \bigcap_{k=1}^{\infty} B_{k}=B_{k_{0}} \cap \left(X\smallsetminus \bigcap_{k=1}^{\infty} B_k \right)=B_{k_{0}} \cap \left(\bigcup_{k=1}^{\infty} (X \smallsetminus B_k) \right)=\bigcup_{k=1}^{\infty} (B_{k_0}\cap (X\smallsetminus B_k)) \\
&=\bigcup_{k=1}^{\infty} (B_{k_0}\smallsetminus B_k)=\bigcup_{k=k_{0}+1}^{\infty} (B_{k_0}\smallsetminus B_k) = \bigcup_{k=k_0+1}^{\infty} (B_{k_{0}}\smallsetminus B_{k_{0}+(k-k_{0})}) = \bigcup_{n=1}^{\infty} (B_{k_{0}}\smallsetminus B_{k_{0}+n}).
\end{aligned}
$$

Define $A_k:=B_{k_{0}}\smallsetminus B_{k_{0}+k}$ for every $k \in \mathbb{N}$. Then $(A_k)$ is an increasing sequence of elements of $\mathsf{S}$ such that
$$
\bigcup_{k=1}^{\infty} A_k= B_{k_0}\smallsetminus \bigcap_{k=1}^{\infty} A_k.
$$

From part \textit{(a)} and the subtractivity of $\mu$, we obtain
$$
\begin{aligned}
\mu(B_{k_0})-\mu\left(\bigcap_{k=1}^{\infty} B_k \right) = \mu \left( \bigcup_{k=1}^{\infty}  A_k\right) = \lim_{k\to \infty} \mu(A_k) &= \mu(B_{k_{0}})-\lim_{k\to\infty}\mu(B_{k_0+k}).
\end{aligned}
$$

Subtracting $\mu(B_{k_{0}})<+\infty$ and multiplying by $-1$, the result follows.
\end{proof}

\begin{corollary} \label{415}
Let $(X,\mathsf{S})$ be a measure space, let $\mu:\mathsf{S} \to \overline{\mathbb{R}}$, and let $(A_k)$ be a sequence of elements of $\mathsf{S}$. Then,
\begin{itemize}
\item[(a)] $\mu\left( \displaystyle\liminf_{k \to \infty} A_k \right) \leq \displaystyle\liminf_{k \to \infty} \mu(A_k)$.
\item[(b)] $ \displaystyle\limsup_{k \to \infty} \mu(A_k) \leq \mu\left( \displaystyle\limsup_{k \to \infty} A_k \right)$ whenever $\mu \left(\displaystyle\bigcup_{k=1}^{\infty} A_k \right) < +\infty.$
\end{itemize}
\end{corollary}

\begin{proof}
\textit{(a):} The sequence $(B_{k})$ of elements of $\mathsf{S}$ defined by
$$
B_{k}:=\bigcap_{j=k}^{\infty} A_{j}
$$
is increasing and, therefore,
$$
\lim_{k \to \infty}B_k=\bigcup_{k=1}^{\infty} B_k = \bigcup_{k=1}^{\infty}\bigcap_{j=k}^{\infty} A_j = \liminf_{k \to \infty} A_k
$$
by Proposition \ref{242}. Theorem \ref{414} ensures that
$$
\begin{aligned}
\mu\left( \liminf_{k \to \infty}A_{k} \right)=\mu\left(\bigcup_{k=1}^{\infty} B_k \right)=\lim_{k \to\infty} \mu(B_k).
\end{aligned}
$$

Since $B_{k} \subset A_{k}$ for every $k \in \mathbb{N}$, the monotonicity of the measure and of the lower limit in the inequality above imply that
$$
\begin{aligned}
\mu\left( \liminf_{k \to \infty}A_{k} \right)= \liminf_{k \to \infty} \mu(B_k) \leq \liminf_{k \to \infty} \mu (A_k),
\end{aligned}
$$
as stated.

\textit{(b):} The sequence $(D_k)$ of elements of $\mathsf{S}$ given by 
$$
D_{k}:=\bigcup_{j=k}^{\infty}A_{j},
$$
is decreasing and, consequently, Proposition \ref{242} ensures that
$$
\lim_{k \to \infty} D_k=\bigcap_{k=1}^{\infty}D_k=\bigcap_{k=1}^{\infty}\bigcup_{j=k}^{\infty} A_{j} = \limsup_{k \to \infty} A_{k}.
$$

Since $\mu(D_{1})<+\infty$, Theorem \ref{414} ensures that
$$
\lim_{k \to \infty}\mu(D_{k}) =\mu\left( \bigcap_{k=1}^{\infty} D_k\right)=\mu\left( \limsup_{k \to\infty}A_{k} \right).
$$

Since $A_k\subset D_k$ for every $k \in \mathbb{N}$, applying the monotonicity of the measure and of the upper limit, we obtain
$$
\limsup_{k \to \infty}\mu(A_k) \leq \limsup_{k \to \infty} \mu(D_k)=\mu\left( \limsup_{k \to \infty} A_{k} \right).
$$

This proves the result.
\end{proof}

\begin{corollary} \label{416}
Let $(X,\mathsf{S})$ be a measure space, let $\mu:\mathsf{S} \to \overline{\mathbb{R}}$, and let $(A_k)$ be a sequence of elements of $\mathsf{S}$ converging to the set $A \in \mathsf{S}$. Then,
$$
\lim_{k \to \infty} \mu(A_k) = \mu(A)
$$
whenever $\mu \left(\displaystyle\bigcup_{k=1}^{\infty} A_k \right) < +\infty.$
\end{corollary}

\begin{proof}
Since the sequence $(A_k)$ converges to the measurable set $A$, the upper and lower limits of $(A_k)$ are both equal to $A$. Corollary \ref{315} ensures that
$$
\limsup_{k \to \infty} \mu(A_k) \leq \mu \left( \limsup_{k \to \infty} A_k \right) = \mu (A) = \mu \left( \liminf_{k \to \infty} A_k \right) \leq \liminf_{k \to \infty} \mu(A_k) \leq \limsup_{k \to \infty}\mu(A_{k})
$$
and, consequently, $\displaystyle\lim_{k \to \infty} \mu(A_k) = \mu(A)$.
\end{proof}

Let us now consider some examples related to the previously established results.

\begin{example} \label{417}
Let $X=\mathbb{N}$, $\mathsf{S}=\mathcal{P}(X)$, and let $\mu=\mu^{\sharp}$ be the counting measure. For each $k \in \mathbb{N}$ define
$$
A_k :=\{k,k+1,\ldots \}.
$$
The sequence $(A_k)$ is a decreasing sequence of elements of $\mathsf{S}$ with $\mu^{\sharp}(A_k)=+\infty$ for every $k \in \mathbb{N}$ since $A_k$ is countably infinite. Note that
$$
X\smallsetminus A_k=\{1,2,\ldots,k-1 \}
$$
for every $k >1$ and $X\smallsetminus A_1 = \varnothing$, so that
$$
\bigcup_{k=1}^{\infty} (X\smallsetminus A_k) = X = \mathbb{N}.
$$
Consequently,
$$
X\smallsetminus \bigcup_{k=1}^{\infty} (X\smallsetminus A_k) = \bigcap_{k=1}^{\infty} A_k = \varnothing.
$$

Therefore,
$$
\mu^{\sharp} \left(  \bigcap_{k=1}^{\infty} A_k \right) = 0 < \lim_{k\to\infty} \mu^{\sharp}(A_k) = +\infty.
$$

This example shows that the inequality in part (b) of {\rm Theorem \ref{414}} may be strict.
\end{example}

\begin{example} \label{418}
Let $X=\mathbb{N}$, $\mathsf{S}=\mathcal{P}(X)$, and let $\mu=\mu^{\sharp}$ be the counting measure. 

Fix $j \in \mathbb{N}$. Define $A_k:=\{ k,k+1,\ldots \}$ for each $k \leq j$ and $A_k:=\varnothing$ for $k>j$. Then $(A_{k})$ is a decreasing sequence of elements of $\mathsf{S}$. Note that the numerical sequence $(\mu^{\sharp}(A_k))$ is eventually constant equal to zero, and therefore $\displaystyle\lim_{k \to \infty} \mu^{\sharp}(A_k)=0$.

Moreover,
$$
\bigcap_{k=1}^{\infty} A_k = \varnothing \quad \mbox{so that}\quad \mu^{\sharp}\left( \bigcap_{k=1}^{\infty} A_k\right)=0.
$$

Consequently,
$$
\mu^{\sharp}\left( \displaystyle\bigcap_{k=1}^{\infty} A_k\right)=\displaystyle\lim_{k \to \infty} \mu^{\sharp}(A_k).
$$
\end{example}

\begin{example} \label{419}
Let $X=[0,1]$, $\mathsf{S}=\mathcal{B}([0,1])$, and let $\mu=\lambda$ be the Lebesgue measure on $\mathsf{S}$. Let $A=[1/4,3/4] \in \mathsf{S}$ and define the sequence $(A_k)$ as follows:
$$
A_k:=\left\{\begin{array}{lcl}
X\smallsetminus A & & \mbox{if }\,\, k \,\,\,\mbox{is even},\\
A & & \mbox{if }\,\, k \,\,\,\mbox{is odd}.
\end{array}
\right.
$$

Since $\lambda(A)=\lambda([1/4,3/4])=1/2$ and
$$
\lambda(X\smallsetminus A)=\lambda([0,1/4)\cup (3/4,1])=1/2,
$$
the numerical sequence $(\lambda(A_k))$ is constant equal to $1/2$ and, therefore,
$$
\displaystyle\limsup_{k \to \infty} \lambda(A_k)= 1/2 = \displaystyle\liminf_{k \to \infty} \lambda(A_k).
$$

Furthermore,
$$
\liminf_{k \to \infty} A_k =\varnothing\quad \mbox{and}\quad  \limsup_{k \to \infty} A_k =X,
$$
so that
$$
\lambda\left(\liminf_{k \to \infty} A_k \right) =0\quad \mbox{and}\quad \lambda\left(\limsup_{k \to \infty} A_k \right)=1.
$$

Therefore,
$$
\lambda\left(\liminf_{k \to \infty} A_k \right)  < \liminf_{k \to \infty} \lambda(A_k) \quad \mbox{and} \quad \limsup_{k \to \infty} \lambda(A_k) < \lambda\left(\limsup_{k \to \infty} A_k \right).
$$

That is, the inequalities given in {\rm Corollary \ref{415}} may be strict.
\end{example}

\begin{example} \label{420}
Let $X=\mathbb{N}$, $\mathsf{S}=\mathcal{P}(X)$, and let $\mu=\mu^{\sharp}$ be the counting measure. Define the sequence $(A_k)$ by
$$
A_k:=\{k,k+1,\ldots \}
$$
for every $k \in \mathbb{N}$. The sequence $(A_k)$ is decreasing and satisfies
$$
\bigcup_{k=1}^{\infty} A_k = A_1 =\mathbb{N}. 
$$

Consequently,
$$
\displaystyle\mu^{\sharp}\left( \bigcup_{k=1}^{\infty} A_k \right) = +\infty.
$$

Now, since the numerical sequence $(\mu^{\sharp}(A_k))$ is constantly equal to $+\infty$, then
$$
\displaystyle\limsup_{k \to \infty} \mu^{\sharp}(A_k) = +\infty.
$$
On the other hand,
$$
\limsup_{k \to \infty} A_k = \bigcap_{k=1}^{\infty} A_k = \varnothing\quad \mbox{so that}\quad \mu^{\sharp}\left( \limsup_{k \to \infty} A_k \right)=0,
$$
and therefore part (b) of {\rm Corollary \ref{415}} does not hold.
\end{example}

The following property of measures is of fundamental importance in their study and is known as \textbf{$\sigma$-subadditivity}. \index{sigma@$\sigma$!subadditivity}

\begin{proposition} \label{421}
Let $(X,\mathsf{S})$ be a measurable space and let $\mu:\mathsf{S}\to \overline{\mathbb{R}}$ be a measure. If $(A_k)$ is a sequence of elements of $\mathsf{S}$, then
$$
\mu \left(  \bigcup_{k=1}^{\infty} A_k \right) \leq \sum_{k=1}^{\infty} \mu (A_k).
$$
\end{proposition}

\begin{proof}
By Proposition \ref{236}, there exists a disjoint sequence $(E_k)$ of elements of $\mathsf{S}$ such that $E_k \subset A_k$ for every $k \in \mathbb{N}$ and
$$
\bigcup_{k=1}^{\infty} A_k = \bigcup_{k=1}^{\infty}E_k.
$$

Therefore,
$$
\mu\left( \bigcup_{k=1}^{\infty} A_k  \right) = \mu \left(  \bigcup_{k=1}^{\infty} E_k  \right) = \sum_{k=1}^{\infty} \mu(E_k) \leq \sum_{k=1}^{\infty} \mu(A_k),
$$
as stated.
\end{proof}

We now state and prove the famous Borel--Cantelli lemma\footnote{Francesco Paolo Cantelli (1875--1966) was an Italian mathematician. He was born in Palermo. He received his doctorate in mathematics in 1899 from the University of Palermo with a thesis on celestial mechanics and continued his interest in astronomy by working until 1903 at the osservatorio astronomico cittadino (National Astronomical Observatory) of Palermo, under the direction of Annibale Riccò. Cantelli's early works dealt with problems in astronomy and celestial mechanics. During 1903--1923 he carried out research on the mathematics of financial theory and actuarial science, as well as on probability theory, for which he became famous. Cantelli's later work was entirely devoted to probability, and it is in this field that his name appears in the Borel--Cantelli lemma and the Glivenko--Cantelli theorem. In 1916--1917 he contributed to the theory of stochastic convergence. In 1923 he resigned from his actuarial position when he was appointed professor of actuarial mathematics at the University of Catania. From there he moved to the University of Naples as professor and later, in 1931, to Sapienza University of Rome, where he remained until his retirement in 1951. He died in Rome.} in its easy direction, which we will apply later on.

\begin{figure}[ht!]
\centering
\includegraphics[scale=0.25]{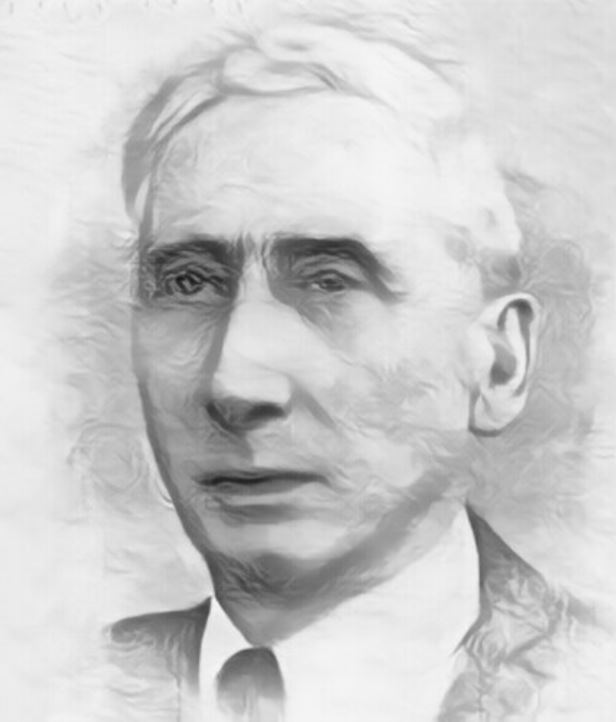} 
\begin{center}
F. Cantelli (1875-1966)
\end{center}
\end{figure}

\begin{theorem}[Borel--Cantelli Lemma] \label{422} \index{lemma!Borel--Cantelli}
Let $(X,\mathsf{S})$ be a measurable space and let $\mu:\mathsf{S} \to \overline{\mathbb{R}}$ be a measure. If $(A_k)$ is a sequence of elements of $\mathsf{S}$ such that $\sum_{k=1}^{\infty} \mu(A_k)<+ \infty,$ then $\mu\left(\displaystyle\limsup_{k \to \infty} A_k \right)=0.$
\end{theorem}

\begin{proof}
The sequence $(B_k)$ of elements of $\mathsf{S}$ given by $B_{k}:=\bigcup_{j=k}^{\infty} A_{j}$
is decreasing and therefore
$$
\lim_{k \to \infty}B_{k}=\bigcap_{k=1}^{\infty}B_{k} =\bigcap_{k=1}^{\infty} \bigcup_{j=k}^{\infty} A_j=\limsup_{k \to \infty} A_k.
$$

From Proposition \ref{421} it follows that
\begin{equation}\label{F41}
    0\leq \mu \left( B_{k} \right)= \mu \left( \bigcup_{j=k}^{\infty} A_j \right) \leq \sum_{j=k}^{\infty} \mu(A_j),\qquad \forall \, k \in \mathbb{N}.
\end{equation}

Since the series $\sum_{k=1}^{\infty} \mu(A_k)$ converges in $\mathbb{R}$, the Cauchy criterion for series ensures that $\sum_{j=k}^{\infty}\mu(A_j) \to 0$ as $k \to \infty$. Taking limits as $k \to \infty$ in (\ref{F41}), we obtain
$$
\lim_{k \to \infty} \mu(B_{k})=0.
$$

Now, from inequality (\ref{F41}) we obtain $\mu(B_{1}) \leq \sum_{k=1}^{\infty} \mu(A_{k})<+\infty,$ and therefore
$$
\mu\left(\displaystyle\limsup_{k \to \infty} A_k \right)=\lim_{k \to \infty}\mu(B_{k})=0
$$
by Theorem \ref{414}. This concludes the proof.
\end{proof}

The Borel--Cantelli lemma also ensures that if $(A_k)$ is a sequence of elements of $\mathsf{S}$ such that $\sum_{k=1}^{\infty}\mu(A_k)<+\infty,$
then $\mu\left( \displaystyle\liminf_{k \to \infty} A_k\right)=0$ since $\displaystyle\liminf_{k \to \infty} A_k \subset \displaystyle\limsup_{k \to \infty} A_k.$

The following example shows that the converse of the Borel--Cantelli lemma is not true in general.

\begin{example} \label{423}
Let $X=(0,1)$, $\mathsf{S}=\mathcal{B}(0,1)$, and let $\mu=\lambda$ be the Lebesgue measure on $\mathsf{S}$.

Let $(a_k)$ be a sequence in $X$ such that $a_{k} \to 0$ as $k \to \infty$. If $A_{k}:=(0,a_{k})$, then $\limsup_{k \to \infty}A_{k}=\varnothing$. Indeed, let $x \in (0,1)$. Since $a_{k} \to 0$ and $x>0$, there exists $k(x) \in \mathbb{N}$ such that $0<a_{k} < x$ for every $k \geq k(x)$. Consequently, $x \in X \smallsetminus \limsup_{k \to \infty}A_{k}$ {\rm [Exercise \ref{E230}]} and, therefore, $\lambda\left(\limsup_{k \to \infty}A_{k} \right)=0$.

However, if $a_{k} \geq \frac{1}{k+1}$ for every $k \in \mathbb{N}$, then $\sum_{k=1}^{\infty}\lambda(A_k)=\sum_{k=1}^{\infty}a_{k}=+\infty$.
\end{example}

The reader may find in \cite{Rinconmono} a very interesting application of the Borel--Cantelli lemma.

\section{Measure spaces}

\begin{definition} \label{424}
A \textbf{measure space} is a triple $(X,\mathsf{S},\mu)$, where $(X,\mathsf{S})$ is a measurable space and $\mu:\mathsf{S}\to \overline{\mathbb{R}}$ is a measure. \index{space!measure space}
\end{definition}

\begin{definition} \label{425} \index{set!null}
Let $(X,\mathsf{S},\mu)$ be a measure space. A set $N \in \mathsf{S}$ is called \textbf{$\mu$-null} if $\mu(N)=0$.

We denote by $\mathcal{N}(\mu)$ the class of $\mu$-null sets, that is,
$$
\mathcal{N}(\mu):=\{N \in \mathsf{S}\,:\,\mu(N)=0 \}.
$$
\end{definition}

By the $\sigma$-subadditivity of $\mu$, the class of $\mu$-null sets is a $\sigma$-ring of subsets of $X$, as shown in the following proposition.

\begin{proposition} \label{426}
Let $(X,\mathsf{S},\mu)$ be a measure space. The class $\mathcal{N}(\mu)$ is a $\sigma$-ring of subsets of $X$.
\end{proposition}

\begin{proof}
$\mathcal{N}(\mu)$ is nonempty since $\varnothing \in \mathcal{N}(\mu)$. Let $N,M \in \mathcal{N}(\mu)$ be arbitrary. Since $N\smallsetminus M \subset N$, we have $0 \leq \mu(N\smallsetminus M) \leq \mu(N) =0,$ and therefore $N\smallsetminus M \in \mathcal{N}(\mu)$. This proves (SR1).

Let $(N_k)$ be a sequence of elements of $\mathcal{N}(\mu)$. Proposition \ref{421} ensures that
$$
0\leq \mu \left(\bigcup_{k=1}^{\infty} N_k \right) \leq \sum_{k=1}^{\infty} \mu(N_k ) =0.
$$

Consequently,
$$
\mu \left(\displaystyle\bigcup_{k=1}^{\infty} N_k \right)=0,
$$
that is, $\mathcal{N}(\mu)$ satisfies (SR2).
\end{proof}

For example, in the measure space $(\mathbb{N},\mathcal{P}(\mathbb{N}),\mu^{\sharp})$, the empty set is the only subset of measure zero. Indeed, the well-ordering principle ensures that every nonempty subset $A$ of $\mathbb{N}$ has a minimum element and, therefore,
$$
\mu^{\sharp}(A) \geq 1.
$$

\begin{definition} \label{427} \index{almost everywhere}
Let $(X,\mathsf{S},\mu)$ be a measure space. We say that a property $\wp(x)$ holds \textbf{almost everywhere relative to $\mu$ ($\mu$ a.e.) on $X$} if there exists $N \in \mathcal{N}(\mu)$ such that $\wp(x)$ is true for every $x \in X\smallsetminus N$.
\end{definition}

Let us now look at some simple examples.

\begin{example} \label{428}
Let $(X,\mathsf{S},\mu)$ be a measure space and let $f,g,h:X \to \overline{\mathbb{R}}$ be functions. By writing $f \leq g$ $\mu$ a.e. we mean that there exists $N \in \mathcal{N}(\mu)$ such that
$$
f(x) \leq g(x)\quad \mbox{for every}\,\, x \in X\smallsetminus N.
$$

Thus $f=g$ $\mu$ a.e. is equivalent to $f \leq g$ $\mu$ a.e. and $g\leq f$ $\mu$ a.e.. Indeed, if $f \leq g$ $\mu$ a.e. and $g\leq f$ $\mu$ a.e., then there exist $N_1,N_2 \in \mathcal{N}(\mu)$ such that
$$
f(x) \leq g(x)\quad \mbox{for every}\,\, x \in X\smallsetminus N_1
$$
$$
g(x) \leq f(x)\quad \mbox{for every}\,\, x \in X\smallsetminus N_2.
$$

Defining $N:=N_1 \cup N_2$, it is clear that $N \in \mathcal{N}(\mu)$ since it is a $\sigma$-ring of subsets of $X$ and, therefore, $f(x)=g(x)$ for every $x \in X\smallsetminus N$. The converse implication follows directly.

Moreover, if $f \leq g$ $\mu$ a.e. and $g \leq h$ $\mu$ a.e., then $f \leq h$ $\mu$ a.e.. Again, taking $M:=M_1 \cup M_2 \in \mathcal{N}(\mu)$ where $M_1$ and $M_2 \in \mathcal{N}(\mu)$ are such that
$$
f(x) \leq g(x)\quad \mbox{for every}\,\, x \in X\smallsetminus M_1
$$
$$
g(x) \leq h(x)\quad \mbox{for every}\,\, x \in X\smallsetminus M_2
$$
then $f(x) \leq h(x)$ for every $x \in X\smallsetminus M$.
\end{example}

\begin{example} \label{429}
Let $(X,\mathsf{S},\mu)$ be a measure space and let $(f_k)$ be a sequence of $\mathsf{S}$-measurable functions. The set $A=\{x \in X\,:\, (f_k(x))\,\,\mbox{converges in }\,\mathbb{R} \}$ satisfies {\rm [Exercise \ref{E315}]}
$$
A= \bigcap_{n=1}^{\infty} \bigcup_{k=1}^{\infty} \bigcap_{p=1}^{\infty} \left\{ x \in X\,:\, |f_k(x)-f_{k+p}(x)| < \frac{1}{n}\right\} \in \mathsf{S}.
$$

And, if $N:=X\smallsetminus A$ is a set of measure zero, then $(f_k)$ converges $\mu$ a.e. to a function $f:X \to \mathbb{R}$ where $f$ may be chosen $\mathsf{S}$-measurable if it is defined by
$$
f(x):=\left\{
\begin{array}{ccl}
\displaystyle\lim_{k \to \infty} f_k(x) & & \mbox{if}\,\,x \in A,\\
0 & & \mbox{if}\,\,x\in N. 
\end{array}
\right.
$$
\end{example}

That a property $\wp(x)$ holds $\mu$ a.e. should not be confused with the statement that the set
$$
\{ x\in X\,:\, \wp(x)\,\,\text{is false}\}
$$
is necessarily $\mu$-null, since this set may, in principle, fail to be measurable. For example, let $X=\{1,2,3\}$, $\mathsf{S}=\sigma(\{1\})$, and $\mu=\delta_{1}$. The property $\wp(x): x \neq 2$ holds $\delta_{1}$ a.e. since it is true for every $x \in X \smallsetminus \{2,3\}$ and $\delta_{1}(\{2,3\})=0$. However,
$$
\{x \in X \,:\, \wp(x)\,\,\text{is false}\}=\{2\}
$$
and $\{2\} \notin \mathsf{S}$.

In the next section we shall construct, from a given measure space $(X,\mathsf{S},\mu)$, a new measure space in which every subset of a $\mu$-null set becomes measurable.

\section{Complete measure spaces}

\begin{definition}\label{430} \index{space!complete measure}
A measure space $(X,\mathsf{S},\mu)$ is said to be \textbf{complete} if, for every set $N\in\mathcal{N}(\mu)$ and every subset $M\subset N$, one has $M\in \mathcal{N}(\mu)$.
\end{definition}

For example, the measure space $(X,\mathcal{P}(X),\mu)$ is complete for every measure $\mu$ defined on $(X,\mathcal{P}(X))$.

Next we give a characterization of complete measure spaces, which is more intuitive in terms of completeness.

\begin{proposition} \label{431}
A measure space $X=(X,\mathsf{S},\mu)$ is complete if and only if $A \in \mathsf{S}$ whenever there exist subsets $N_1,N_2 \in \mathsf{S}$ such that $N_1 \subset A \subset N_2$ and $\mu(N_2\smallsetminus N_1)=0$.
\end{proposition}

\begin{figure}[ht!]
\centering
\begin{tikzpicture}[xscale=1.5, yscale=1.5]
\filldraw[draw=black, fill=gray!5]
plot[smooth cycle] coordinates{
(-0.85,0)(-0.5,0.8)(-0.3,0.8)(1,0.8)(0.785,0)(0.75,-0.5)(-0.5,-0.9)};
\draw[draw=black] (0,0) circle(0.5);
\draw[draw=black] (0,0) circle(0.75);
\draw(0,0)node{$_{N_1}$};
\draw(0.325,0.55)node{$_{A}$};
\draw(0.85,0.65)node{$_{N_2}$};
\end{tikzpicture}
\end{figure}

\begin{proof}
$\Rightarrow):$ Let $A \subset X$ be such that there exist $N_1,N_2 \in \mathsf{S}$ with $N_1 \subset A \subset N_2$ and $\mu(N_2\smallsetminus N_1)=0$. Since $A \smallsetminus N_1 \subset N_2\smallsetminus N_1$, then $A\smallsetminus N_1 \in \mathsf{S}$ because $X$ is complete. Consequently, $A=(A\smallsetminus N_1) \cup N_1 \in \mathsf{S}$.

$\Leftarrow):$ Let $M,N \subset X$ be such that $N \in \mathcal{N}(\mu)$ and $M \subset N$. Since $\varnothing \subset M \subset N$ and $\mu(N\smallsetminus\varnothing)=\mu(N)=0$, then $M \in \mathsf{S}$. That is, $X$ is complete.
\end{proof}

If $(X,\mathsf{S},\mu)$ is a non-complete measure space, it is possible to construct a new complete measure space by enlarging its $\sigma$-algebra. Indeed, one defines a new $\sigma$-algebra on $X$ generated by $\mathsf{S}$ and by all subsets of $\mu$-null sets, and extends the measure $\mu$ to a measure $\bar{\mu}$ defined on this $\sigma$-algebra. The resulting measure space is complete. This procedure is called the \textbf{completion of a measure space}, and the resulting space is known as the \textbf{$\mu$-completion} of $(X,\mathsf{S},\mu)$.

\begin{theorem}[Completion of a Measure Space] \label{432} 
Let $(X,\mathsf{S},\mu)$ be a measure space and let $\mathcal{N}(\mu)$ be the $\sigma$-ring of $\mu$-null sets. Define
$$
\overline{\mathsf{S}}:=\left\{ F \subset X\,:\, F=A \cup M_0\,\,\text{with}\,\,A\in \mathsf{S}\,\,\text{and}\,\,M_0\subset N_0\,\,\text{for some}\,\,N_0 \in \mathcal{N}(\mu) \right\}.
$$

Then,
\begin{itemize}
\item[(a)] $\overline{\mathsf{S}}$ is a $\sigma$-algebra of subsets of $X$ such that $\mathsf{S} \subset \overline{\mathsf{S}}$.
\item[(b)] The function $\overline{\mu}:\overline{\mathsf{S}} \to \overline{\mathbb{R}}$ given by $\overline{\mu} (A \cup M_0):=\mu(A)$ is well defined and is a measure on $(X,\overline{\mathsf{S}})$ with $\overline{\mu}_{|_{\mathsf{S}}}=\mu$.
\end{itemize}
\end{theorem}

\begin{proof}
\textit{(a):} Any element $A$ of $\mathsf{S}$ can be written as $A=A\cup \varnothing$, which implies that $\mathsf{S} \subset \overline{\mathsf{S}}$ and, consequently, $X \in \overline{\mathsf{S}}$. 

Let $F \in \overline{\mathsf{S}}$. There exist $A \in \mathsf{S}$ and $M_{0} \subset N_{0}$ with $N_{0} \in \mathcal{N}(\mu)$ such that $F=A \cup M_{0}$. Then,
$$
\begin{aligned}
X\smallsetminus F&=X\smallsetminus (A\cup M_0) = (X\smallsetminus A) \cap (X\smallsetminus M_{0}) \\
&= \left[ (X\smallsetminus A) \cap (X\smallsetminus M_{0}) \right] \cap X\\
&= \left[ (X\smallsetminus A) \cap (X\smallsetminus M_{0}) \right] \cap \left[ (X\smallsetminus N_0)\cup N_0\right]\\
&=\left[(X\smallsetminus A) \cap (X\smallsetminus M_{0}) \cap (X\smallsetminus N_0) \right]\cup \left[ (X\smallsetminus A) \cap (X\smallsetminus M_{0}) \cap N_0\right]\\
&=\left[(X\smallsetminus A) \cap (X\smallsetminus N_0) \right]\cup \left[ (X\smallsetminus A) \cap (X\smallsetminus M_{0}) \cap N_0\right]
\end{aligned}
$$
with $(X\smallsetminus A) \cap (X\smallsetminus N_0) \in \mathsf{S}$ and $(X\smallsetminus A) \cap (X\smallsetminus M_{0}) \cap N_0 \subset N_0$. Therefore, $X\smallsetminus F \in \overline{\mathsf{S}}$. 

Let $(F_k)$ be a sequence of elements of $\overline{\mathsf{S}}$. There exist sequences $(A_k)$, $(M_k)$, and $(N_k)$ of subsets of $X$ with $A_k \in \mathsf{S}$ and $M_k \subset N_k$ for $N_k \in \mathcal{N}(\mu)$ such that $F_k=A_k \cup M_k$ for every $k \in \mathbb{N}$. Then,
$$
\bigcup_{k=1}^{\infty} F_k = \left( \bigcup_{k=1}^{\infty} A_k \right) \cup \left( \bigcup_{k=1}^{\infty} M_k\right)
$$
where $\bigcup_{k=1}^{\infty} A_k \in \mathsf{S}$ and $\bigcup_{k=1}^{\infty} M_k \subset \bigcup_{k=1}^{\infty} N_k$ with $\bigcup_{k=1}^{\infty} N_k \in \mathcal{N}(\mu)$ by Proposition \ref{426}.

Therefore, $\bigcup_{k=1}^{\infty} F_k \in \overline{\mathsf{S}}$ and, consequently, $\overline{\mathsf{S}}$ is a $\sigma$-algebra (see Proposition \ref{28}).

\textit{(b):} Let $F \in \overline{\mathsf{S}}$ be arbitrary. Then $F=A_1 \cup M_1$ with $A_1 \in \mathsf{S}$ and $M_1 \subset N_1$ for some $N_1 \in \mathcal{N}(\mu)$. Consider another representation of $F$, namely, $F=A_2 \cup M_2$, where $A_2 \in \mathsf{S}$, $M_2 \subset N_2$, and $N_2 \in \mathcal{N}(\mu)$. 

Therefore, $A_1 \cup M_1 = A_2 \cup M_2$ and, since $M_1 \subset N_1$ and $M_2 \subset N_2$, it follows that
$$
\begin{aligned}
A_1 \subset A_1 \cup M_1 \cup N_1 &= A_2 \cup M_2 \cup N_1 \subset A_2 \cup N_2 \cup N_1,\\
A_2 \subset A_2 \cup M_2 \cup N_2 &= A_1 \cup M_1 \cup N_2 \subset A_1 \cup N_1 \cup N_2.\\
\end{aligned}
$$

Thus,
$$
\begin{aligned}
\mu(A_1) &\leq \mu(A_2 \cup N_2 \cup N_1) \leq \mu(A_2) + \mu(N_2) + \mu(N_1) = \mu(A_2),\\
\mu(A_2) &\leq \mu(A_1 \cup N_1 \cup N_2) \leq \mu(A_1) + \mu(N_1) + \mu(N_2) = \mu(A_1)
\end{aligned}
$$
and, therefore, $\mu(A_1)=\mu(A_2)$, which shows that $\overline{\mu}$ is well defined.

It is immediate to verify that $\overline{\mu}:\overline{\mathsf{S}} \to \overline{\mathbb{R}}$ satisfies (M1) and (M2). 

Let $(F_k)$ be a disjoint sequence of elements of $\overline{\mathsf{S}}$. By the previous part we have
$$
\bigcup_{k=1}^{\infty} F_k = \left( \bigcup_{k=1}^{\infty} A_k \right) \cup \left( \bigcup_{k=1}^{\infty} M_k\right)
$$
where $(A_k)$ is a sequence of elements of $\mathsf{S}$ and $(M_k)$ is a sequence such that $M_k \subset N_k$ for every $k \in \mathbb{N}$ with $N_k \in \mathcal{N}(\mu)$. Note that if $k,j \in \mathbb{N}$ are such that $k \neq j$, then
$$
A_k \cap A_j \subset (A_k \cup M_k) \cap (A_j \cup M_j) = F_k \cap F_j = \varnothing
$$
which means that $(A_k)$ is a disjoint sequence of elements of $\mathsf{S}$ and, therefore,
$$
\overline{\mu} \left( \bigcup_{k=1}^{\infty} F_k  \right) = \mu \left( \bigcup_{k=1}^{\infty} A_k  \right)=\sum_{k=1}^{\infty} \mu(A_k) = \sum_{k=1}^{\infty} \overline{\mu}(F_k).
$$

That is, $\overline{\mu}$ satisfies (M3).

Finally, since $A=A\cup\varnothing$ for every $A \in \mathsf{S} \subset \overline{\mathsf{S}}$, then $\overline{\mu}(A)=\overline{\mu}(A \cup \varnothing) = \mu(A)$.

Consequently, $\overline{\mu}$ is a measure on $(X,\overline{\mathsf{S}})$ such that $\overline{\mu}_{|_{\mathsf{S}}}=\mu$.
\end{proof}

We now proceed to verify that the measure space constructed in the previous theorem is indeed complete.

\begin{theorem}
Let $(X,\mathsf{S},\mu)$ be a measure space and let $(X,\overline{\mathsf{S}},\bar{\mu})$ be its $\mu$-completion. Then $(X,\overline{\mathsf{S}},\bar{\mu})$ is a complete measure space.
\end{theorem}

\begin{proof}
Let $F \in \overline{\mathsf{S}}$ be such that $\overline{\mu}(F)=0$ and let $G \subset F$ be arbitrary. There exist $A \in \mathsf{S}$ and $M_0 \subset N_0$ for some $N_0 \in \mathcal{N}(\mu)$ such that $F=A \cup M_0$. Since $\overline{\mu}(F)=\mu(A)$, it follows that $A \in \mathcal{N}(\mu)$ and, consequently, $A \cup N_{0} \in \mathcal{N}(\mu)$. Since $G \subset F\subset A \cup N_{0}$, by writing $G=\varnothing \cup G$, we conclude that $G \in \overline{\mathsf{S}}$. That is, $(X,\overline{\mathsf{S}},\overline{\mu})$ is complete.
\end{proof}

In Proposition \ref{431} a characterization of complete measure spaces was established that clarifies the notion of completeness. It is then easy to see that $F \in \overline{\mathsf{S}}$ if and only if there exist $N_1,N_2 \in \mathsf{S}$ such that $N_1 \subset F \subset N_2$ and $\mu(N_2\smallsetminus N_1)=0$ [Exercise \ref{E430}].

Let us now look at an example where we complete a measure space.

\begin{example} \label{434}
Let $(X,\mathsf{S},\delta_1)$ be the measure space given by $X=\{1,2,3\}$, the $\sigma$-algebra $\mathsf{S}=\{\varnothing,\{1\},\{2,3\},X \}$, and $\delta_1$ the unit measure concentrated at $1$.

It is clear that the class of $\mu$-null sets is
$$
\mathcal{N}(\delta_1)=\{N \in \mathsf{S}\,:\, \delta_1(N)=0 \}=\{N\in \mathsf{S}\,:\, 1 \not \in N \}=\{ \varnothing, \{2,3 \}\}.
$$

The space $(X,\mathsf{S},\delta_1)$ is not complete since $\{2\} \subset \{2,3 \}$ with $\{2,3\} \in \mathcal{N}(\delta_1)$ and $\{2\} \not \in \mathsf{S}$.

The subsets of the elements of $\mathcal{N}(\delta_1)$ are $\varnothing$, $\{2\}$, and $\{3\}$ and, by adjoining them to the elements of $\mathsf{S}$, we obtain
$$
\overline{\mathsf{S}}=\{\varnothing,\{1\},\{2\},\{3\},\{1,2\},\{1,3\},\{2,3\},X \}=\mathcal{P}(X)
$$
and
$$
\bar{\delta_1}(F)=\left\{
\begin{array}{lcl}
1 & &si\,\,\, 1\in F,\\
0 & &si\,\,\, 1\not\in F.
\end{array}
\right.
$$

Therefore, $(X,\mathcal{P}(X),\overline{\delta_1})$ is the $\delta_{1}$-completion of $(X,\mathsf{S},\delta_1)$.
\end{example}

We mentioned in the previous section that for a proposition $\wp(x)$ in a measure space $(X,\mathsf{S},\mu)$ which holds $\mu$ a.e., in general it is not true that the set $\{x \in X\,:\,\wp(x)\,\,\text{is false}\}$ belongs to $\mathcal{N}(\mu)$. However, if $(X,\mathsf{S},\mu)$ is complete, then $\{x \in X\,:\,\wp(x)\,\,\text{is false}\} \in \mathcal{N}(\mu)$ is equivalent to saying that $\wp(x)$ holds $\mu$ a.e. [Exercise \ref{E427}].

We conclude this section with a very important result that we will use in the following chapters.

\begin{proposition} \label{435}
Let $(X,\mathsf{S},\mu)$ be a complete measure space, let $f:X\to \overline{\mathbb{R}}$ be an $\mathsf{S}$-measurable function, and let $g:X \to \overline{\mathbb{R}}$ be a function. If $g=f$ $\mu$ a.e., then $g$ is $\mathsf{S}$-measurable.
\end{proposition}

\begin{figure}[ht!]
\begin{minipage}[r]{0.5\textwidth}
\begin{center}
\begin{tikzpicture}[xscale=0.8, yscale=0.8]
\draw[->,gray] (0,0) -- (7,0);
\draw [->,gray] (0,-1.2) -- (0,1.2);
\draw[ultra thick, blue] plot[smooth] coordinates
{(0.00,0.00)(0.79,0.71)(1.57,1.00)(2.36,0.71)
(3.14,0.00)(3.93,-0.71)(4.71,-1.00)(5.50,-0.71)(6.28,-0.00)};

\draw (3.25,1) node{$_{f}$};

\draw[dotted] (0,0.7)--(7,0.7);
\draw(0,0.7)node[left]{$_{c}$};
\end{tikzpicture}
\end{center}
\end{minipage} \hfill 
 \begin{minipage}[l]{0.5\textwidth}
\begin{center}
\begin{tikzpicture}[xscale=0.8,yscale=0.8]
\draw[->,gray] (0,0) -- (7,0);
\draw [->,gray] (0,-1.2) -- (0,1.2);
\draw[ultra thick] plot[smooth] coordinates
{(0.00,0.00)(0.79,0.71)(1.57,1.00)(2.36,0.71)
(3.14,0.00)(3.93,-0.71)(4.71,-1.00)(5.50,-0.71)(6.28,-0.00)};

\draw (3.25,1) node{$_{g}$};
\draw[dotted] (4.71,-1.00)--(4.71,0.00);
\draw[dotted] (1.57,1.00)--(1.57,0.00);

\draw(4.71,0.00) node{$_{\bullet}$};
\draw(1.57,0.00)node{$_{\bullet}$};
\draw(4.71,-1.00) node{$_{\circ}$};
\draw(1.57,1.00)node{$_{\circ}$};

\draw[dotted] (0,0.7)--(7,0.7);
\draw(0,0.7)node[left]{$_{c}$};
\end{tikzpicture}
\end{center}
\end{minipage}
\end{figure}

\begin{proof}
Let $N:=\{x \in X\,:\,g(x) \neq f(x) \}$. Then $N \in \mathcal{N}(\mu)$ since $(X,\mathsf{S},\mu)$ is complete [Exercise \ref{E427}]. Thus, for every $ c \in \mathbb{R}$ we have
$$
\{x \in X\,:\, g(x)>c \}=\left( \{x \in X\,:\,f(x)>c \} \cup \{ x\in N\,:\,g(x)>c\} \right) \smallsetminus \{x \in N\,:\,g(x) \leq c \}
$$ 
where $\{ x \in X\,:\,f(x)>c\} \in \mathsf{S}$ because $f$ is an $\mathsf{S}$-measurable function, and the sets $\{x \in N\,:\,g(x)>c \}$ and $\{x\in N\,:\,g(x) \leq c\}$ are subsets of $N$. The completeness of the measure space implies that $g^{-1}((c,+\infty]) =\{x \in X\,:\,g(x)>c\} \in \mathsf{S}$ for every $c \in \mathbb{R}$. Consequently, $g$ is $\mathsf{S}$-measurable.
\end{proof}

It is an easy exercise to show that the previous proposition may fail if the measure space $(X,\mathsf{S},\mu)$ is not complete [Exercise \ref{E427}]. However, it is possible to consider the completion of the measure space in order to obtain a result analogous to Proposition \ref{435}. We propose this here as an exercise [Exercise \ref{E431}].

\section{Uniqueness of measures}

The following theorems provide necessary conditions for two measures defined on a measurable space $(X,\mathsf{S})$ to coincide on it. These results are similar to the Hahn uniqueness theorem that we will study later and are a direct application of Dynkin's theorem.

\begin{theorem} \label{436}
Let $(X,\mathsf{S})$ be a measurable space and let $\mathcal{C} \subset \mathcal{P}(X)$ be a $\pi$-system such that $\sigma(\mathcal{C})=\mathsf{S}$. If $\mu,\nu:\mathsf{S} \to \mathbb{R}$ are two finite measures such that $\mu(X)=\nu(X)$ and $\mu(K)=\nu(K)$ for every $K \in \mathcal{C}$, then $\mu = \nu$ on $\mathsf{S}$.
\end{theorem}

\begin{proof}
Let $\mathcal{L}_{0}=\left\{A \in \mathsf{S}\,:\, \mu(A)=\nu(A) \right\}$. By hypothesis it is clear that $X \in \mathcal{L}_{0}$. Let us show that $\mathcal{L}_{0}$ is in fact a $\lambda$-system. Let $A,B \in \mathcal{L}_{0}$ be such that $B \subset A$. From Proposition \ref{44} it follows that
$$
\mu(A \smallsetminus B)=\mu(A)-\mu(B)=\nu(A)-\nu(B)=\nu(A \smallsetminus B)
$$
and, therefore, $A \smallsetminus B \in \mathcal{L}_{0}$. This proves (L2). 

Let $(A_k)$ be an increasing sequence of elements of $\mathcal{L}_{0}$. Applying Theorem \ref{414} we obtain
$$
\mu\left( \bigcup_{k=1}^{\infty} A_k\right)=\lim_{k\to\infty}\mu(A_k)=\lim_{k\to\infty}\nu(A_k)=\nu\left( \bigcup_{k=1}^{\infty} A_k\right)
$$ 
and it follows that $\bigcup_{k=1}^{\infty}A_k \in \mathcal{L}_{0}$. That is, $\mathcal{L}_{0}$ satisfies (L3).

Consequently, $\mathcal{C}$ is a $\pi$-system such that $\mathcal{C} \subset \mathcal{L}_{0}$. Theorem \ref{255} ensures that $\sigma(\mathcal{C})=\mathsf{S} \subset \mathcal{L}_{0}$, which concludes the proof.
\end{proof}

The following example shows that the conclusion of the previous theorem may fail if the hypothesis that the class $\mathcal{C}$ is a $\pi$-system is omitted.

\begin{example} \label{437}
Let $X=\{1,2,3,4\}$, $\mathsf{S}=\mathcal{P}(X)$, and $\mathcal{C}=\{\{1,2\},\{2,3\},\{3,4\}\}$. It is immediate that $\sigma(\mathcal{C})=\mathcal{P}(X)$ and that $\mathcal{C}$ is not a $\pi$-system.

Consider the measures $\mu,\nu:\mathcal{P}(X) \to \mathbb{R}$ given by $\mu(\{1\})=\mu(\{2\})=\mu(\{3\})=\mu(\{4\})=\frac{1}{4}$ and $\nu(\{1\})=\nu(\{3\})=\frac{1}{2}$ and $\nu(\{2\})=\nu(\{4\})=0$. It is clear that $\mu(X)=1=\nu(X)$ and $\mu(K)=\nu(K)$ for every $K \in \mathcal{C}$, but $\nu \neq \mu$ on $\mathcal{P}(X)$.
\end{example}

The following result now addresses the case of two $\sigma$-finite measures and is obtained as a direct consequence of the previous theorem.

\begin{theorem} \label{438}
Let $(X,\mathsf{S})$ be a measurable space and let $\mathcal{C} \subset \mathcal{P}(X)$ be a $\pi$-system such that $\sigma(\mathcal{C})=\mathsf{S}$. If $\mu, \nu: \mathsf{S} \to \overline{\mathbb{R}}$ are two measures such that $\mu(K)=\nu(K)$ for every $K \in \mathcal{C}$ and if there exists an increasing sequence $(K_j)$ of elements of $\mathcal{C}$ that converges to $X$ and is such that $\mu(K_{j})=\nu(K_{j})<+\infty$ for every $j \in \mathbb{N}$, then $\mu=\nu$ on $\mathsf{S}$.
\end{theorem}

\begin{proof}
Denote by $\mu_{j}$ and $\nu_{j}:K_{j}\cap \mathsf{S}  \to \mathbb{R}$ the restrictions of $\mu$ and $\nu$ determined by the set $K_{j}$, which are finite measures.

Let $j \in \mathbb{N}$. The class $K_{j} \cap \mathcal{C}$ is a $\pi$-system of subsets of $K_{j}$ such that $\sigma(K_{j} \cap \mathcal{C})=K_{j} \cap \sigma(\mathcal{C})=K_{j} \cap \mathsf{S}$ by Theorem \ref{221}. Thus, by hypothesis, $\mu_{j}(K_j)=\nu_{j}(K_{j})$ and $\mu_{j}(K)=\nu_{j}(K)$ for every $K \in K_{j} \cap \mathcal{C}$. Theorem \ref{436} ensures that $\mu_{j}=\nu_{j}$ on $K_{j} \cap \mathsf{S}$. 

Let $A \in \mathsf{S}$. Since $(K_{j})$ is an increasing sequence of elements of $\mathsf{S}$ such that $K_{j} \to X$, then $(K_{j} \cap A)$ is an increasing sequence of elements of $\mathsf{S}$ such that $K_{j} \cap A \to A$. From Theorem \ref{414} and what was previously shown, it follows that 
$$
\mu(A)=\lim_{j \to \infty} \mu(K_{j} \cap A)=\lim_{j \to \infty}\nu(K_{j} \cap A)=\nu(A). 
$$

Consequently, $\mu=\nu$ on $\mathsf{S}$.
\end{proof}

\section{Exercises}

\begin{exercise} \label{E41}
Let $X$ be an uncountable set and let
$$
\mathsf{S}=\{A \subset X\,:\, A\,\,\mbox{or}\,\,X\smallsetminus A\,\,\mbox{is finite or countable}\}.
$$
Define $\mu:\mathsf{S} \to \{0,1\}$ as follows
$$
\mu(A)=\left\{
\begin{array}{lcl}
0 & & \text{if}\,\,A\,\,\text{ is finite or countable},\\
1 & & \text{if}\,\,X\smallsetminus A\,\,\text{ is finite or countable}.
\end{array}
\right.
$$

Show that $\mu$ is a finite measure on $(X,\mathsf{S})$.
\end{exercise}

\begin{exercise} \label{E42}
Let $X$ be an uncountable set and let $\mathsf{S}=\mathcal{P}(X)$. Define $\mu:\mathsf{S} \to \overline{\mathbb{R}}$ as follows
$$
\mu(A)=\left\{
\begin{array}{lcl}
0 & & \text{if}\,\,A\,\,\text{ is finite or countable},\\
+\infty & & \text{if}\,\,A\,\,\text{ is uncountable}.
\end{array}
\right.
$$

Show that $\mu$ is a measure on $(X,\mathsf{S})$.
\end{exercise}

\begin{exercise} \label{E43}
Let $X$ be an infinite set and let $\mathsf{S}=\mathcal{P}(X)$. Define $\mu,\nu:\mathsf{S} \to \overline{\mathbb{R}}$ as follows
$$
\mu(A)=\left\{
\begin{array}{lcl}
0 & & \mbox{if}\,\,A\,\,\mbox{is finite},\\
+\infty & & \mbox{if}\,\,A\,\,\mbox{is infinite}
\end{array}
\right. \,\,\,\,\,\,\,\,\,\, \nu(E)=\left\{
\begin{array}{lcl}
0 & & \mbox{if}\,\,A\,\,\mbox{is finite},\\
1 & & \mbox{if}\,\,A\,\,\mbox{is infinite}.
\end{array}
\right.
$$

Which properties of {\rm Definition \ref{41}} do $\mu$ and $\nu$ satisfy? Is it possible to conclude that $\mu$ and $\nu$ are measures?
\end{exercise}

\begin{exercise} \label{E44}
Let $(X,\mathsf{S})$ be an arbitrary measurable space and let $x,y \in X$ be fixed. Show that $\delta_{x}$ and $\delta_{y}$, the unit measures concentrated at $x$ and $y$ respectively, are equal if and only if $x$ and $y$ belong exactly to the same sets in $\mathsf{S}$.
\end{exercise}

\begin{exercise} \label{E45}
Let $X$ be nonempty, let $\mathsf{S}=\mathcal{P}(X)$, and let $\mu=\mu^{\sharp}$ be the counting measure. Show the following statements:
\begin{itemize}
\item[(a)] $\mu^{\sharp}$ is finite if and only if $X$ is finite.
\item[(b)] $\mu^{\sharp}$ is $\sigma$-finite if and only if $X$ is countable.
\end{itemize}
\end{exercise}

\begin{exercise} \label{E46}
Let $X=\mathbb{N}$ and $\mathsf{S}=\mathcal{P}(\mathbb{N})$. Show that every measure $\mu$ on $(X,\mathsf{S})$ is obtained from a unique sequence of nonnegative extended real numbers $(x_{j})$ in the following way:
$$
\mu(A)=\left\{
\begin{array}{lcl}
0 & & \mbox{if}\,\,\,A=\varnothing, \\
\displaystyle\sum_{j \,\in \,A} x_{j} &&\mbox{if}\,\,\,A \neq \varnothing.
\end{array}
\right.
$$
\end{exercise}

{\setlength{\parindent}{0pt}
\begin{exercise} \label{E47}
Let $A \subset \{1,2,\ldots,99\}$ be a subset consisting of exactly 10 elements. Show that there exist disjoint nonempty subsets $I,J$ of $A$ such that
$$
\sum_{a\,\in\,I} a = \sum_{a\,\in\,J} a.
$$
\end{exercise}

(Hint: Show that the function $\varphi:\mathcal{P}(A) \to \mathbb{R}$ given by $\varphi(A)=\displaystyle\sum_{a\in A} a$ cannot be injective).}

\begin{exercise} \label{E48}
Let $(X,\mathsf{S})$ be a measurable space, let $\mu:\mathsf{S} \to \overline{\mathbb{R}}$ be a measure, and fix $A \in \mathsf{S}$ with $\mu(A)>0$. Define $\mu_{A}:A \cap \mathsf{S} \to \overline{\mathbb{R}}$ by $\mu_{A}(B):=\mu(B)$. Show that $\mu_{A}$ is a measure on $(A,A\cap \mathsf{S})$ called the \textbf{restriction of $\mu$ to $A$}, and determine the conditions under which $\mu_{A}$ is finite and $\sigma$-finite. 

What happens if $A \in \mathsf{S}$ is such that $\mu(A)=0$? Justify your answer.
\end{exercise}

\begin{exercise} \label{E49}
Let $(Y,\mathsf{T},\nu)$ be a measure space and let $f:X \to Y$ be a bijective function. Consider the measurable space $(X,\mathsf{S})$ where $\mathsf{S}$ is the $\sigma$-algebra of subsets of $X$ given by
$$
\mathsf{S}:=\{f^{-1}(F) \,:\,F \in \mathsf{T} \}.
$$

Define $\mu:\mathsf{S} \to \overline{\mathbb{R}}$ by $\mu(\mathfrak{E}):=\nu(f(\mathfrak{E}))$. Show that $\mu$ is a measure on $(X,\mathsf{S})$.

What happens if the hypothesis that $f$ is bijective is omitted? Justify your answer.
\end{exercise}

\begin{exercise} \label{E410}
Let $(X,\mathsf{S})$ be a measurable space. Show the following statements:
\begin{itemize}
\item[(a)]  If $(\mu_k)$ is a nondecreasing sequence of measures on $(X,\mathsf{S})$, then the formula $\mu(A):=\displaystyle\lim_{k \to \infty}\mu_k(A)$ with $A \in \mathsf{S}$ defines a measure on $(X,\mathsf{S})$.
\item[(b)] If $(\mu_k)$ is an arbitrary sequence of measures on $(X,\mathsf{S})$, then the formula $\mu(A):=\sum_{k=1}^{\infty}\mu_k(A)$ with $A \in \mathsf{S}$ defines a measure on $(X,\mathsf{S})$.
\end{itemize}
\end{exercise}

{\setlength{\parindent}{0pt}
\begin{exercise} \label{E411}
Let $(\mathbb{R},\mathcal{B}(\mathbb{R}))$ be a measurable space and let $(x_j)$ be a sequence of real numbers. Show the following statements:
\begin{itemize}
\item[(a)] The function $\mu:\mathcal{B}(\mathbb{R}) \to \overline{\mathbb{R}}$ given by $\mu(B):=\sum_{j=1}^{\infty} \delta_{x_{j}}(B)$ is a measure.
\item[(b)] The measure $\mu$ from the previous item assigns finite values to bounded intervals of $\mathbb{R}$ if and only if $\displaystyle\lim_{j \to \infty}|x_j|=+\infty$.
\item[(c)] Solve: for what type of sequences $(x_{j})$ in $\mathbb{R}$ is the measure $\mu$ from item (a) $\sigma$-finite? Justify your answer.
\end{itemize} 
\end{exercise}

(Hint: Use Exercise \ref{E410}).}

\begin{exercise} \label{E412}
Consider the measurable space $(\mathbb{R},\mathcal{B}(\mathbb{R}))$. Define $\mu:\mathcal{B}(\mathbb{R}) \to \overline{\mathbb{R}}$ by letting $\mu(A)=$ the number of rational numbers in $A$, with $\mu(A):=+\infty$ if there are infinitely many rational numbers in $A$. Show that $\mu$ is a $\sigma$-finite measure on $(\mathbb{R},\mathcal{B}(\mathbb{R}))$ under which every open interval of $\mathbb{R}$ has measure equal to $+\infty$.
\end{exercise}

\begin{exercise} \label{E413}
Let $(X,\mathsf{S})$ be a measurable space.
\begin{itemize}
\item[(1)] Let $\mu:\mathsf{S} \to \overline{\mathbb{R}}$ be a nonnegative and $\sigma$-additive function. Show that if $\mu(A)<+\infty$ for some $A \in \mathsf{S}$, then $\mu(\varnothing)=0$. Conclude that $\mu$ is a measure on $(X,\mathsf{S})$.
\item[(2)] Give an example showing that, in general, the condition $\mu(\varnothing)=0$ is not a consequence of properties (M2) and (M3) from {\rm Definition \ref{41}}.
\end{itemize}
\end{exercise}

\begin{exercise} \label{E414}
Let $(X,\mathsf{S},\mu)$ be a measure space, and let $(A_k)$ and $(B_k)$ be two sequences of elements of $\mathsf{S}$. Show the following statements:
\begin{itemize}
\item[(a)] If $\mu(A_j \cap A_k)=0$ for every $j,k$ with $j \neq k$, then
$$
\mu \left( \bigcup_{k=1}^{\infty} A_k \right) = \sum_{k=1}^{\infty} \mu(A_k)\qquad (almost\,\,\sigma\text{-additivity}).
$$
\item[(b)] If $\mu(A_k \bigtriangleup B_k)=0$ for every $k \in \mathbb{N}$, then 
$$
\mu \left( \bigcup_{k=1}^{\infty} A_k\,\bigtriangleup\,\bigcup_{k=1}^{\infty} B_k \right),\quad \mu\left( \bigcap_{k=1}^{\infty} A_k\,\bigtriangleup\,\bigcap_{k=1}^{\infty} B_k \right)
$$
$$
\mu\left( \limsup_{k\to\infty} A_k\,\bigtriangleup\,\limsup_{k\to\infty}B_k \right),\quad \mu\left( \liminf_{k\to\infty} A_k\,\bigtriangleup\,\liminf_{k\to\infty}B_k \right)
$$
are all equal to zero.
\end{itemize}
\end{exercise}

\begin{exercise} \label{E415}
Let $(X,\mathsf{S},\mu)$ be a measure space. Is it true that if $A,B\in \mathsf{S}$ are such that $\mu(A \cap B)=0$, then $A\cap B=\varnothing$? Justify your answer.
\end{exercise}

\begin{exercise} \label{E416}
Let $(X,\mathsf{S})$ be a measurable space and let $\mu:\mathsf{S} \to \mathbb{R}$ be a finite, nonnegative, and additive function. Show that the following conditions are equivalent for sequences $(A_k)$ of elements of $\mathsf{S}$:
\begin{itemize}
\item[(a)] $\mu\left(\displaystyle\bigcup_{k=1}^{\infty} A_k \right)=\sum_{k=1}^{\infty} \mu(A_k)$ whenever $A_k \cap A_{\ell} = \varnothing$ for $k \neq \ell$ ($\sigma$-additivity).
\item[(b)] If $(A_k)$ is decreasing and $A=\displaystyle\bigcap_{k=1}^{\infty} A_k$, then $\mu(A)=\displaystyle\lim_{k \to \infty}\mu(A_k)$ (continuity from above). 
\item[(c)] If $(A_k)$ is increasing and $A=\displaystyle\bigcup_{k=1}^{\infty} A_k$, then $\mu(A)=\displaystyle\lim_{k\to \infty}\mu(A_k)$ (continuity from below). 
\item[(d)] $\mu\left(\displaystyle\lim_{k\to\infty} A_k \right)=\displaystyle\lim_{k\to\infty} \mu(A_k)$ (continuity).
\end{itemize}
\end{exercise}

\begin{exercise}[Inclusion--Exclusion Formula] \label{E417}
Let $\mu$ be a measure on $(X,\mathsf{S})$ and let $A_{1},\ldots,A_{n} \in \mathsf{S}$ be such that $\mu(A_i)<+\infty$ for each $i =1,\ldots, n$. Prove the inclusion--exclusion formula, also known as the Poincaré formula, that is, show that
$$
\begin{aligned}
\mu\left( \bigcup_{i=1}^{n} A_i \right)&=\sum_{i=1}^{n} \mu(A_i)-\sum_{1\leq i < j \leq n} \mu(A_i \cap A_j)+\sum_{1\leq i < j <k\leq n} \mu(A_i \cap A_j \cap A_k)+ \ldots\\
&\quad +\ldots+(-1)^{n+1} \mu \left( \bigcap_{i=1}^{n} A_i \right).
\end{aligned}
$$
\end{exercise}

\begin{exercise} \label{E418}
Let $(X,\mathsf{S},\mu)$ be a measure space and let $A_1,\ldots,A_{N} \in  \mathsf{S}$. For each fixed $\ell \in \{ 1,2,\ldots,N\}$, let
$$
\mathfrak{C}_{\ell}:=\{x \in X\,:\,x\in A_j\,\,\text{for exactly}\,\,\ell\,\,\text{indices}\,\,j \in \{ 1,2,\ldots,N\}  \}.
$$
Show the following statements:
\begin{itemize}
\item[(a)] $\mathfrak{C}_{\ell} \in \mathsf{S}$.
\item[(b)] $\sum_{\ell=1}^{N}\mu(A_\ell)=\sum_{\ell=1}^{k} \ell\,\mathfrak{C}_{\ell}$.
\item[(c)] If
$$
\mathfrak{D}_{\ell}=\{x \in X\,:\,x\in E_{j}\,\,\text{for at most}\,\,\ell\,\,\mbox{indices}\,\,j \in \{ 1,2,\ldots,N\}\},
$$
then
$$
\sum_{\ell=1}^{N} \mu(E_\ell) \leq \sum_{\ell=1}^{N} \ell\,\mu(\mathfrak{D}_\ell).
$$
\item[(d)] Assume that $\mu(X)=1$ and that the sets $A_1,\ldots,A_N \in  \mathsf{S}$ are such that each $x \in X$ belongs to at least $r$ of the sets $A_{1},\ldots,A_{N}$. Conclude that there exists at least one index $i$ such that $\mu(A_i) \geq r/k$.
\end{itemize}
\end{exercise}

{\setlength{\parindent}{0pt}
\begin{exercise} \label{E419}
Let $(X,\mathsf{S},\mu)$ be a $\sigma$-finite measure space and let $\mathfrak{D} \subset \mathsf{S}$ be a family of pairwise disjoint sets. Let $A \in \mathsf{S}$ be fixed with $\mu(A)>0$. Show that the family
$$
\mathfrak{D}_{A}=\left\{ D \in \mathfrak{D}\,:\, \mu(A \cap D)>0\right\}
$$
is at most countable.
\end{exercise}

(Hint: First consider the case $\mu(A)<+\infty$).}

\begin{exercise} \label{E420}
Let $\{ (X_{j}, \mathsf{S}_{j},\mu_{j})\,:\, j \in \mathbb{N}\}$ be a countable collection of measure spaces such that $X_{i} \cap X_k = \varnothing$ for every $i,k \in \mathbb{N}$ with $i \neq k$. Define the following:
\begin{itemize}
\item[(1)] $X:=\bigcup_{j\,\in\,\mathbb{N}}\, X_{j}$.
\item[(2)] $ \mathsf{S}:=\{A \subset X\,:\,A \cap X_{j} \in  \mathsf{S}_j\,\,\,\text{for every}\,\,j \in \mathbb{N} \}$.
\item[(3)] $\mu:\mathsf{S} \to \overline{\mathbb{R}}$ by
$$
\mu(A):=\sum_{j \in \mathbb{N}} \,\mu_{j}(A \cap X_{j}).
$$
\end{itemize}

Show the following statements:
\begin{itemize}
\item[(a)] $(X, \mathsf{S},\mu)$ is a measure space (called the \textbf{direct sum} of $\{ (X_{j},\mathsf{S}_{j},\mu_{j})\}_{j\, \in \,\mathbb{N}}$). 
\item[(b)] $\mu$ is $\sigma$-finite if and only if $\mu_{j}$ is $\sigma$-finite for every $j \in \mathbb{N}$.
\end{itemize}
\end{exercise}

{\setlength{\parindent}{0pt}
\begin{exercise}[Nikodým Metric Space] \index{space!metric!Nikodým} \label{E421}
Let $(X,\mathsf{S},\mu)$ be a measure space and let
$$
\mathfrak{F}:=\{A \in \mathsf{S}\,:\,\mu(A)<+\infty \}.
$$
Let $A,B \in \mathfrak{F}$. We say that $A$ and $B$ are equivalent, $A \sim B$, if and only if $\mu(A \bigtriangleup B)=0$.

Show the following statements:
\begin{itemize}
\item[(a)] $\sim$ is an equivalence relation on $\mathfrak{F}$.
\item[(b)] Denote by $[A]$ the equivalence class of $A \in \mathfrak{F}$ and define $\widetilde{\mathfrak{F}}:=\{[A]\,:\,A \in \mathfrak{F} \}.$
The function $d:\widetilde{\mathfrak{F}} \times \widetilde{\mathfrak{F}} \to \mathbb{R}$ given by $d([A],[B]):=\mu(A \bigtriangleup B)$ is well defined, that is, it does not depend on the chosen representatives of $A$ and $B$, and it is a metric on $\widetilde{\mathfrak{F}}$.
\item[(c)] $(\tilde{\mathfrak{F}}, d)$ is a complete metric space.
\end{itemize}
\end{exercise}

(Hint: If $([A_k])$ is a $d$-Cauchy sequence in $\widetilde{\mathfrak{F}}$, find a subsequence $([A_{k_j}])$ such that $d([A_{k_{j}}],[A_{k_{j+1}}])<\frac{1}{2^j}$ for every $j \in \mathbb{N}$. Define $A_{\ast}:=\displaystyle\liminf_{j \to \infty} A_{k_{j}}$ and prove that $A_{\ast} \in \mathfrak{F}$ and that $A_{\ast} \bigtriangleup A_{k_{j}} \subset \bigcup_{\ell = j}^{\infty}\,(A_{k_{\ell}} \bigtriangleup A_{k_{\ell+1}}).$ Conclude that $d([A_{k_{j}}],[A_{\ast}]) \to 0$ and, therefore, that $d([A_k],[A_{\ast}])\to 0$).}

\begin{exercise} \label{E422}
Assume the hypotheses and notation of the previous exercise. Show that the functions $\phi:\widetilde{\mathfrak{F}} \times \widetilde{\mathfrak{F}} \to \widetilde{\mathfrak{F}}$ with values $[A \bigtriangleup B], \qquad [A \cup B], \qquad [A \cap B], \qquad [A-B]$ are $d$-uniformly continuous if $\mu(X)<+\infty$. Prove the analogous statement for the function $\varphi:\widetilde{\mathfrak{F}} \to \widetilde{\mathfrak{F}}$ given by
$\varphi([A]):=[X\smallsetminus A].$
\end{exercise}

\begin{exercise} \label{E423}
\begin{itemize}
\item[(1)] Give an example where the strict inequality in item (a) of {\rm Corollary \ref{415}} holds, and another one where equality holds.
\item[(2)] Give an example where the strict inequality in item (b) of {\rm Corollary \ref{415}} holds, and another one where equality holds.
\item[(3)] Prove by means of an example that item (b) of {\rm Corollary \ref{415}} may fail if
$$
\mu\left(\bigcup_{k=1}^{\infty}\, A_k\right)=+\infty.
$$
\end{itemize}
\end{exercise}

{\setlength{\parindent}{0pt}
\begin{exercise} \label{E424}
Consider the measurable space $(\mathbb{R},\mathcal{B}(\mathbb{R}),\lambda)$ and let $(p_k)$ be a sequence of polynomials converging uniformly on $[a,b]$ to a function $f$. Define $H_{k}:=\{x \in [a,b]\,:\,p_{k}(x)=f(x) \}$ and show that if $f$ is not a polynomial on $[a,b]$, then $\lambda(H_k) \to 0$ as $k \to \infty$.
\end{exercise}

(Hint: Show that if $p_j \neq p_k$ then $\lambda(H_j \cap H_k)=0$ and conclude that $\sum_{k=1}^{\infty}\lambda(H_k)<+\infty.$).}

\begin{exercise} \label{E425}
Let $(X,\mathsf{S},\mu)$ be a $\sigma$-finite measure space and let $f:X \to \overline{\mathbb{R}}$ be an $\mathsf{S}$-measurable function. Show that there exists a sequence of $\mathsf{S}$-simple functions $(s_k)$ such that
$$
s_k(x) \to f(x)
$$
for every $x \in X$ and
$$
\mu(\{x \in X\,:\,s_k(x) \neq 0 \})<+\infty
$$
for every $k \in \mathbb{N}$.
\end{exercise}

{\setlength{\parindent}{0pt}
\begin{exercise} \label{E426}
Let $(X,\mathsf{S},\mu)$ be a measure space and let $A,B \in \mathsf{S}$. Show that $\chi_{A}=\chi_{B}$ $\mu$ a.e. if and only if
$$
A \bigtriangleup B \subset N
$$
for some $N \in \mathcal{N}(\mu)$.
\end{exercise}

(Hint: Show that $\chi_{A}(x) \neq \chi_{B}(x)$ if and only if $x \in A \bigtriangleup B$).}

{\setlength{\parindent}{0pt}
\begin{exercise} \label{E427}
Let $(X,\mathsf{S},\mu)$ be a measure space.
\begin{itemize}
    \item[(a)] Prove, by means of an example, that {\rm Proposition \ref{435}} may fail if $(X,\mathsf{S},\mu)$ is not complete.
    \item[(b)] Show that if $\wp(x)$ is a property on $X$ and
$$
\{x \in X\,:\,\wp(x)\,\,\text{is false} \} \in \mathcal{N}(\mu),
$$
then $\wp(x)$ holds $\mu$ a.e.
    \item[(c)] Show that if $(X,\mathsf{S},\mu)$ is complete and $\wp(x)$ is a property on $X$, then
$$
\{x \in X\,:\,\wp(x)\,\,\text{is false} \} \in \mathcal{N}(\mu)
$$
if and only if $\wp(x)$ holds $\mu$ a.e.
\end{itemize}
\end{exercise}

(Hint: For item \textit{(a)}, {\rm Exercise \ref{E36}} is useful.)}

\begin{exercise} \label{E428}
Let $(X,\mathsf{S},\mu)$ be a complete measure space. Show that if $A \in \mathsf{S}$ and $\mu(A \bigtriangleup B)=0$ for every $B \subset X$, then $B \in \mathsf{S}$ and $\mu(B)=\mu(A)$. 
\end{exercise}

\begin{exercise} \label{E429}
Let $(X,\mathsf{S},\mu)$ be a measure space and let $(X,\overline{\mathsf{S}},\overline{\mu})$ be its $\mu$-completion. Show that if $\left(X,\overline{\overline{\mathsf{S}}},\overline{\overline{\mu}}\right)$ is another $\mu$-completion of $(X,\mathsf{S},\mu)$, then $\left(X,\overline{\overline{\mathsf{S}}},\overline{\overline{\mu}}\right)$ is also a completion of $(X,\overline{\mathsf{S}},\overline{\mu})$.
\end{exercise}

\begin{exercise} \label{E430}
Let $(X,\mathsf{S},\mu)$ be a measure space and let $(X,\bar{\mathsf{S}},\bar{\mu})$ be its $\mu$-completion. Show that $F \in \overline{\mathsf{S}}$ if and only if there exist $N_1,N_2 \in \mathsf{S}$ such that $N_1 \subset F \subset N_2$ and $\mu(N_2\smallsetminus N_1)=0.$
\end{exercise}

\begin{exercise} \label{E431}
Let $(X,\mathsf{S},\mu)$ be a measure space and let $(X,\overline{\mathsf{S}},\overline{\mu})$ be its $\mu$-completion. Show that if $\overline{f}:X\to \overline{\mathbb{R}}$ is an $\overline{\mathsf{S}}$-measurable function, then there exists an $\mathsf{S}$-measurable function $f:X \to \overline{\mathbb{R}}$ such that the set $\overline{N}:=\{x \in X\,:\, f(x)\neq\bar{f}(x) \}$ belongs to $\mathcal{N}(\overline{\mu})$.
\end{exercise}

\begin{exercise} \label{E432}
Let $(X,\mathsf{S})$ be a measurable space and define $\mu:\mathsf{S} \to \overline{\mathbb{R}}$ as follows
$$
\mu(A)=\left\{
\begin{array}{lcl}
0 & &\text{if}\,\,\,A=\varnothing,\\
+\infty & &\text{if}\,\,\,A \neq \varnothing.
\end{array}
\right.
$$

Show that $\mu$ is a measure on $(X,\mathsf{S})$ and determine its $\mu$-completion.
\end{exercise}

\begin{exercise} \label{E433}
Show that if $X=(X,\mathsf{S},\mu)$ is a complete measure space, then $X$ is equal to its $\mu$-completion.
\end{exercise}

\begin{exercise} \label{E434}
Let $(X,\mathsf{S},\mu)$ and $(Y,\mathsf{T},\nu)$ be the measure spaces determined in {\rm Example \ref{410}}. Let $(X,\overline{\mathsf{S}},\overline{\mu})$ be the $\mu$-completion of $(X,\mathsf{S},\mu)$ and let $(Y,\overline{\mathsf{T}},\overline{\nu})$ be the $\nu$-completion of $(Y,\mathsf{S},\nu)$. Show that
$
\overline{\nu}(A)=\overline{\mu}(f^{-1}(A))
$
for every $A \in \overline{\mathsf{T}}$ and deduce that $(Y,\mathsf{T},\nu)$ is complete whenever $(X,\mathsf{S},\mu)$ is complete.
\end{exercise}

\begin{exercise} \label{E435}
Let $(X,\mathsf{S},\mu)$ be a measure space. Define $\mu^{\bullet}:\mathsf{S} \to \overline{\mathbb{R}}$ as follows
$$
\mu^{\bullet}(A):=\sup\left\{\mu(B)\,:\, B \subset A,\,\,B \in \mathsf{S}\,\,\text{and}\,\,\mu(B)<+\infty \right\}.
$$

Show the following statements:
\begin{itemize}
\item[(a)] $\mu^{\bullet}$ is a measure on $(X,\mathsf{S})$.
\item[(b)] If $\mu$ is $\sigma$-finite, then $\mu=\mu^{\bullet}$.
\item[(c)] Determine $\mu^{\bullet}$ when $(X,\mathsf{S},\mu)$ is the measure space given in {\rm Exercise \ref{432}}.
\end{itemize}
\end{exercise}

\begin{exercise} \label{E436}
Let $(X,\mathsf{S},\mu)$ be a measure space. A set $A \subset X$ is said to be \textbf{locally in $\mathsf{S}$} (or \textbf{locally measurable}) if $A \cap B \in \mathsf{S}$ for every $B \in \mathsf{S}$ with $\mu(B)<+\infty$. Define
$$
\widetilde{\mathsf{S}}:=\{A \subset X\,:\, A\,\,\text{is locally in }\,\mathsf{S} \}.
$$
\index{set!locally measurable}

Show the following statements:
\begin{itemize}
\item[(a)] $\widetilde{\mathsf{S}}$ is a $\sigma$-algebra of subsets of $X$ such that $\mathsf{S} \subset \widetilde{\mathsf{S}}$. 
\item[(b)] We say that the measure space $(X,\mathsf{S},\mu)$ is \textbf{saturated} if $\mathsf{S}=\widetilde{\mathsf{S}}$. \index{space!measure!saturated}

If $\mu$ is $\sigma$-finite, then $(X,\mathsf{S},\mu)$ is a saturated measure space.
\item[(c)] The function $\widetilde{\mu}:\widetilde{\mathsf{S}} \to \overline{\mathbb{R}}$ given by
$$
\widetilde{\mu}(A)=\left\{
\begin{array}{lcl}
\mu(A) &&\text{if}\,\,\,A \in \mathsf{S},\\
+\infty & &\text{if}\,\,\,A \in \widetilde{\mathsf{S}}\smallsetminus\mathsf{S}
\end{array}
\right.
$$
is a measure on $(X,\widetilde{\mathsf{S}})$. 

The measure space $(X,\widetilde{\mathsf{S}},\widetilde{\mu})$ is called the \textbf{$\mu$-saturation} of $(X,\mathsf{S},\mu)$.
\item[(d)] If $(X,\mathsf{S},\mu)$ is complete, then $(X,\widetilde{\mathsf{S}},\widetilde{\mu})$ is complete.
\item[(e)] $(X,\widetilde{\mathsf{S}},\widetilde{\mu})$ is a saturated measure space.
\end{itemize}
\end{exercise}

\begin{exercise} \label{E437}
Give an example of a measurable space $(X,\mathsf{S})$ with two finite measures $\mu$ and $\nu$ such that $\mu(X)=\nu(X)$ and such that the set
$$
\left\{A \in \mathsf{S}\,:\,\mu(A)=\nu(A) \right\}
$$
is not a $\sigma$-algebra.
\end{exercise}

\begin{exercise} \label{E438} \index{measure!semifinite}
Let $X=(X,\mathsf{S})$ be a measurable space. We say that a measure $\mu:\mathsf{S} \to \overline{\mathbb{R}}$ is \textbf{semifinite} if, for every $A \in \mathsf{S}$ with $\mu(A)=+\infty$, there exists $B \in \mathsf{S}$ such that $B \subset A$ and $0 < \mu(B)<+\infty$.

Show the following statements:
\begin{itemize}
    \item[(a)] If $\mu$ is $\sigma$-finite, then $\mu$ is semifinite. Give an example showing that the converse is not true in general.
    \item[(b)] For every $A \in \mathsf{S}$ with $\mu(A)=+\infty$ and every $\varepsilon>0$, there exists $B(\varepsilon) \in \mathsf{S}$ such that $B(\varepsilon) \subset A$ and $\varepsilon < \mu(B(\varepsilon))<+\infty$. {\rm (Hint: Show that the set $\{\mu(B)\,:\,B \in \mathsf{S},\,\,\,B \subset A\,\,\,\text{and}\,\,\,\mu(B)<+\infty\}$ is not bounded above.)}
\end{itemize}
\end{exercise}

\begin{exercise} \label{E439}
State and prove an analogue of the uniqueness theorem for measures (see {\rm Theorem \ref{436}}) using the monotone class lemma.
\end{exercise}

\begin{exercise}\label{E440}
Let $(\mathbb{R},\mathcal{B}(\mathbb{R}),\lambda)$ be a measure space. Show the following:
\begin{itemize}
    \item[(a)] Every countable subset $N$ of $\mathbb{R}$ is Borel measurable and has Lebesgue measure equal to $0$.
    \item[(b)] Every compact subset $K$ of $\mathbb{R}$ is Borel measurable and has finite Lebesgue measure.
\end{itemize}
\end{exercise}
\chapter{The Lebesgue Integral of a nonnegative measurable function} \label{Capitulo5}
\markboth{{\scriptsize 5. THE INTEGRAL OF NONNEGATIVE MEASURABLE FUNCTIONS}}{ {\scriptsize 5. THE INTEGRAL OF NONNEGATIVE MEASURABLE FUNCTIONS}}

The modern theory of integration received a decisive impulse around 1904 through the work of Henri Lebesgue. His construction, known as the Lebesgue integral, became a fundamental concept in measure theory; moreover, remarkably interesting applications were found in areas such as probability, differential equations, mathematical analysis, and harmonic analysis, as well as in applied disciplines such as physics, finance, and engineering.

The Lebesgue integral generalizes the Riemann integral, studied in introductory calculus courses, by allowing its definition for a much broader class of real-valued functions and on domains that are not necessarily Euclidean. This extension provides an appropriate framework for problems where the Riemann integral is insufficient.

There are several ways to introduce the Lebesgue integral. In this chapter, we begin with the study of nonnegative simple functions defined on an arbitrary measure space and, from them, extend the definition to any nonnegative measurable function. We will analyze its fundamental properties and the procedures that allow the integration of more general functions.

One of the essential virtues of the Lebesgue integral is that it provides precise conditions under which limits and integration may be interchanged, unlike what occurs in the Riemann theory. This capability is crucial in the analysis of series and of the Fourier transform, among other topics. In this context, we shall present and prove two basic results: the monotone convergence theorem and Fatou's lemma, which will be used repeatedly in subsequent chapters.

Finally, we will be able to establish a comparison between the Riemann and Lebesgue integrals: we shall prove that the Dirichlet function, which we already saw is not integrable in the Riemann sense, is nevertheless Lebesgue integrable, thereby resolving the problem posed at the beginning of this book. The chapter concludes with Lebesgue's theorem for nonnegative, bounded, Borel-measurable functions, based on the analysis of the structure of the set of discontinuities of such functions.

\section{The integral of nonnegative simple functions}

Let $(X,\mathsf{S},\mu)$ be a measure space. We denote by $\overline{\mathbb{M}}(X,\mathsf{S})$ the set of all $\mathsf{S}$-measurable functions $f:X \to \overline{\mathbb{R}}$, and by $\overline{\mathbb{M}}^{+}(X,\mathsf{S})$ the subset consisting of the nonnegative functions.

We shall define the Lebesgue integral for increasingly general measurable functions, starting from the notion of the Lebesgue integral for nonnegative $\mathsf{S}$-simple functions.

We denote by $\mathbb{S}^{+}(X,\mathsf{S})$ the set of all $\mathsf{S}$-simple nonnegative functions $s:X \to \mathbb{R}$.

\begin{definition} \label{51} \index{integral!Lebesgue} \index{integral!Lebesgue!of a nonnegative simple function}
Let $s:X \to \mathbb{R}$ be a nonnegative $\mathsf{S}$-simple function and let
$$
s=\sum_{j=1}^{N} \alpha_{j}\,\chi_{A_j}
$$
be its canonical representation. We define the \textbf{(Lebesgue) integral of $s$ on $X$ with respect to the measure $\mu$}, denoted by $\int_{X} s\,d\mu$, as the nonnegative extended real number given by
$$
\int_{X} s\,d\mu :=\sum_{j=1}^{N}\alpha_{j}\,\mu(A_j).
$$
\end{definition}

For example, if we denote by $0$ the constant function with value $0$ on $X$, it is clear that it is a nonnegative $\mathsf{S}$-simple function and has canonical representation $0=0\cdot\chi_{X} +1\cdot\chi_{\varnothing}$,
so that its Lebesgue integral is equal to $0\cdot\mu(X)=0.$

We now prove that this integral is linear.

\begin{proposition}[Linearity of the Lebesgue integral for simple functions] \label{52}
For every $s,t \in \mathbb{S}^{+}(X,\mathsf{S})$ and every $c \geq 0$, the following hold:
\begin{itemize}
\item[(a)] $\int_X cs\,d\mu=c\int_X s\,d\mu$.
\item[(b)] $\int_X (s+t)\,d\mu = \int_X s\,d\mu + \int_X t\,d\mu$.
\end{itemize}
\end{proposition}

\begin{proof}
We write $s$ and $t$ in canonical form
$$
s=\sum_{j=1}^{N} \alpha_{j}\,\chi_{A_j}\quad\mbox{and}\quad  t=\sum_{k=1}^{M} \beta_{k}\,\chi_{B_k}.
$$

\textit{(a):} If $c=0$, then
$
cs=0=0\cdot\chi_{X}
$
and therefore
$
\int_X cs\,d\mu=\int_X 0\cdot\chi_{X}\,d\mu=0\cdot\mu(X)=0=0\int_X s\,d\mu=c\int_X s\,d\mu.
$
Assume now that $c>0$. Then $cs=\sum_{j=1}^{N} c\alpha_{j}\,\chi_{A_j}$ is the canonical representation of the nonnegative $\mathsf{S}$-simple function $cs$ and, consequently,
$$
\int_X cs\,d\mu =\sum_{j=1}^{N} c\alpha_{j}\,\mu(A_j)=c\sum_{j=1}^{N}\alpha_{j}\,\mu(A_j)=c\int_X s\,d\mu.
$$

\textit{(b):} The function $s+t$ admits the representation
$$
s+t=\sum_{j=1}^{N}\sum_{k=1}^{N}(\alpha_{j}+\beta_{k})\chi_{A_j \cap B_k},
$$
which is not necessarily canonical.

\begin{figure}[ht!]
\begin{minipage}[r]{0.45\textwidth}
\begin{center}
\begin{tikzpicture}[xscale=1,yscale=0.75]
	\draw[->,gray] (0,-2) -- (5.5,-2); \draw [->,gray] (0,-2) -- (0,4); 
\draw (5.8,-2) node{$X$}; \draw (0,4.1) node[right] {$\mathbb{R}$}; 

\draw (0,2.8) node{$-$}; \draw (0,2.8) node[left] {$_{\alpha_{i}}$};
\draw [ thick] (1.2,2.8)--(2.22,2.8);
\draw (1.2,2.8) node{$_{\bullet}$}; 
\draw (2.3,2.8) node{$_{\circ}$};

\draw (0,-1.3) node{$-$}; \draw (0,-1.3) node[left] {$_{\alpha_{n}}$};
\draw [ thick] (2.3,-1.3)--(2.72,-1.3);
\draw (2.3,-1.3) node{$_{\bullet}$}; 
\draw (2.8,-1.3) node{$_{\circ}$};

\draw (0,1) node{$-$}; \draw (0,1) node[left] {$_{\alpha_{j}}$};
\draw [ thick] (0.3,1)--(1.12,1);
\draw (0.3,1) node{$_{\bullet}$}; 
\draw (1.2,1) node{$_{\circ}$};

\draw (0,-0.75) node{$-$}; \draw (0,-0.75) node[left] {$_{\beta_{r}=\alpha_{\ell}}$};
\draw [ thick] (2.8,-0.75)--(4,-0.75);
\draw (2.8,-0.75) node{$_{\bullet}$}; 
\draw (4.1,-0.75) node{$_{\circ}$};

\draw (0,1.8) node{$-$}; \draw (0,1.8) node[left] {$_{\beta_{k}}$};
\draw [ thick, red] (0.3,1.8)--(1.63,1.8);
\draw (0.3,1.8) node{$_{\bullet}$}; 
\draw (1.7,1.8) node{$_{\circ}$};

\draw (0,3.4) node{$-$}; \draw (0,3.4) node[left] {$_{\beta_{m}}$};
\draw [ thick, red] (1.7,3.4)--(3.33,3.4);
\draw (1.7,3.4) node{$_{\bullet}$}; 
\draw (3.4,3.4) node{$_{\circ}$};

\draw [ thick, red] (3.4,-0.75)--(5.5,-0.75);
\draw (3.4,-0.75) node{$_{\bullet}$}; 

\draw [-,dotted ] (0.3,1.8)--(0.3,-2);
\draw [-,dotted ] (1.2,-2)--(1.2,2.8);
\draw [-,dotted ] (1.7,-2)--(1.7,3.4);
\draw [-,dotted ] (2.3,-2)--(2.3,2.8);
\draw [-,dotted ] (2.8,-2)--(2.8,3.4);
\draw [-,dotted ] (3.4,-2)--(3.4,3.4);
\draw [-,dotted ] (4.1,-0.75)--(4.1,-2);
\draw [-,dotted ] (4.7,-0.75)--(4.7,-2);

\draw (0.7,-2.5) node{$_{A_{j}}$};
\draw (1.7,-2.5) node{$_{A_{i}}$};
\draw (2.57,-2.5) node{$_{A_{n}}$};
\draw (3.37,-2.5) node{$_{A_{\ell}}$};

\draw (0.3,-2) node{$_{|}$};
\draw (1.2,-2) node{$_{|}$};
\draw (1.7,-2) node{$_{|}$};
\draw (2.3,-2) node{$_{|}$};
\draw (2.8,-2) node{$_{|}$};
\draw (3.4,-2) node{$_{|}$};
\draw (4.1,-2) node{$_{|}$};

\draw (1.2,-3) node{$_{B_{k}}$};
\draw (2.3,-3) node{$_{B_{m}}$};
\draw (4.8,-2.5) node{$_{B_{r}}$};
\end{tikzpicture}
\begin{center}
Canonical representation of $s$ and $t$
\end{center}
\end{center}
\end{minipage} 
\hfill 
 \begin{minipage}[l]{0.55\textwidth}
\begin{center}
\begin{tikzpicture}[xscale=1,yscale=0.75]
	\draw[->,gray] (0,-2) -- (5.5,-2); \draw [->,gray] (0,-2) -- (0,4); 
\draw (5.8,-2) node{$X$}; \draw (0,4.1) node[right] {$\mathbb{R}$};

\draw (0,1.5) node{$-$}; \draw (0,1.5) node[left] {$_{\alpha_{j}+\beta_{k}=\alpha_{\ell}+\beta_{m}}$};
\draw [ thick, blue] (0.3,1.5)--(1.18,1.5);
\draw (0.3,1.5) node{$_{\bullet}$}; 
\draw (1.2,1.5) node{$_{\circ}$};

\draw (0,2.5) node{$-$}; \draw (0,2.5) node[left] {$_{\alpha_{i}+\beta_{k}}$};
\draw [ thick, blue] (1.2,2.5)--(1.65,2.5);
\draw (1.2,2.5) node{$_{\bullet}$}; 
\draw (1.7,2.5) node{$_{\circ}$};

\draw (0,3.5) node{$-$}; \draw (0,3.5) node[left] {$_{\alpha_{i}+\beta_{m}}$};
\draw [ thick, blue] (1.7,3.5)--(2.25,3.5);
\draw (1.7,3.5) node{$_{\bullet}$}; 
\draw (2.3,3.5) node{$_{\circ}$};

\draw (0,0.5) node{$-$}; \draw (0,0.5) node[left] {$_{\alpha_{n}+\beta_{m}}$};
\draw [ thick, blue] (2.3,0.5)--(2.75,0.5);
\draw (2.3,0.5) node{$_{\bullet}$}; 
\draw (2.8,0.5) node{$_{\circ}$};

\draw (0,1.5) node{$-$}; 
\draw [ thick, blue] (2.8,1.5)--(3.37,1.5);
\draw (2.8,1.5) node{$_{\bullet}$}; 
\draw (3.4,1.5) node{$_{\circ}$};

\draw (0,-0.25) node{$-$}; \draw (0,-0.25) node[left] {$_{\alpha_{\ell}+\beta_{r}}$};
\draw [ thick, blue] (3.4,-0.25)--(4.08,-0.25);
\draw (3.4,-0.25) node{$_{\bullet}$}; 
\draw (4.1,-0.25) node{$_{\circ}$};

\draw (0,-0.85) node{$-$}; \draw (0,-0.85) node[left] {$_{\beta_{r}}$};
\draw [ thick, blue] (4.1,-0.85)--(5.5,-0.85);
\draw (4.1,-0.85) node{$_{\bullet}$}; 

\draw [-,dotted ] (0.3,1.5)--(0.3,-2);
\draw [-,dotted ] (1.2,2.5)--(1.2,-2);
\draw [-,dotted ] (1.7,3.5)--(1.7,-2);
\draw [-,dotted ] (2.3,3.5)--(2.3,-2);
\draw [-,dotted ] (2.8,1.5)--(2.8,-2);
\draw [-,dotted ] (3.4,1.5)--(3.4,-2);
\draw [-,dotted ] (4.1,-0.25)--(4.1,-2);

\draw (0.7,-2.5) node{$_{A_{j}\cap B_k}$};
\draw (1.43,-3) node{$_{A_{i}\cap B_k}$};
\draw (2.1,-2.5) node{$_{A_{i}\cap B_{m}}$};
\draw (2.57,-3) node{$_{A_{n}\cap B_{m}}$};
\draw (3.1,-2.5) node{$_{A_{\ell} \cap B_{m}}$};
\draw (3.8,-3) node{$_{A_{\ell}\cap B_{r}}$};
\draw (4.8,-2.5) node{$_{B_{r}}$};

\draw (0.3,-2) node{$_{|}$};
\draw (1.2,-2) node{$_{|}$};
\draw (1.7,-2) node{$_{|}$};
\draw (2.3,-2) node{$_{|}$};
\draw (2.8,-2) node{$_{|}$};
\draw (3.4,-2) node{$_{|}$};
\draw (4.1,-2) node{$_{|}$};

\end{tikzpicture}
\begin{center}
Graph of $s+t$
\end{center}
\end{center}
\end{minipage}
\end{figure}

Let $c_{1},\ldots,c_{p}$ be the distinct real numbers of the set $\{\alpha_{j}+\beta_{k}\,:\,j=1,\ldots,N\,\,\mbox{and}\,\,k=1,\ldots,M \}$. For each $\ell=1,\ldots,p$ define the set:
$$
G_{\ell}=\bigcup\left\{A_{j}\cap B_{k}\,:\,\alpha_{j}+\beta_{k}=c_{\ell} \right\}.
$$

Then, the canonical representation of the function $s+t$ is given by
$$
s+t=\sum_{\ell=1}^{p}\,c_{\ell}\,\chi_{G_{\ell}}.
$$

\begin{figure}[ht!]
\centering
\begin{tikzpicture}[xscale=1,yscale=0.75]
	\draw[->,gray] (0,-2) -- (5.5,-2); \draw [->,gray] (0,-2) -- (0,4); 
\draw (5.8,-2) node{$X$}; \draw (0,4.1) node[right] {$\mathbb{R}$};

\draw (0,1.5) node{$-$}; \draw (0,1.5) node[left] {$_{c_{\ell}}$};
\draw [ thick, blue] (0.3,1.5)--(1.18,1.5);
\draw (0.3,1.5) node{$_{\bullet}$}; 
\draw (1.2,1.5) node{$_{\circ}$};

\draw (0,2.5) node{$-$}; \draw (0,2.5) node[left] {$_{c_{q}}$};
\draw [ thick, blue] (1.2,2.5)--(1.65,2.5);
\draw (1.2,2.5) node{$_{\bullet}$}; 
\draw (1.7,2.5) node{$_{\circ}$};

\draw (0,3.5) node{$-$}; \draw (0,3.5) node[left] {$_{c_{o}}$};
\draw [ thick, blue] (1.7,3.5)--(2.25,3.5);
\draw (1.7,3.5) node{$_{\bullet}$}; 
\draw (2.3,3.5) node{$_{\circ}$};

\draw (0,0.5) node{$-$}; \draw (0,0.5) node[left] {$_{c_{m}}$};
\draw [ thick, blue] (2.3,0.5)--(2.75,0.5);
\draw (2.3,0.5) node{$_{\bullet}$}; 
\draw (2.8,0.5) node{$_{\circ}$};

\draw (0,1.5) node{$-$}; 
\draw [ thick, blue] (2.8,1.5)--(3.37,1.5);
\draw (2.8,1.5) node{$_{\bullet}$}; 
\draw (3.4,1.5) node{$_{\circ}$};

\draw (0,-0.25) node{$-$}; \draw (0,-0.25) node[left] {$_{c_{n}}$};
\draw [ thick, blue] (3.4,-0.25)--(4.08,-0.25);
\draw (3.4,-0.25) node{$_{\bullet}$}; 
\draw (4.1,-0.25) node{$_{\circ}$};

\draw (0,-0.85) node{$-$}; \draw (0,-0.85) node[left] {$_{c_{p}}$};
\draw [ thick, blue] (4.1,-0.85)--(5.5,-0.85);
\draw (4.1,-0.85) node{$_{\bullet}$}; 

\draw [-,dotted ] (0.3,1.5)--(0.3,-2);
\draw [-,dotted ] (1.2,2.5)--(1.2,-2);
\draw [-,dotted ] (1.7,3.5)--(1.7,-2);
\draw [-,dotted ] (2.3,3.5)--(2.3,-2);
\draw [-,dotted ] (2.8,1.5)--(2.8,-2);
\draw [-,dotted ] (3.4,1.5)--(3.4,-2);
\draw [-,dotted ] (4.1,-0.25)--(4.1,-2);

\draw (0.7,-2.5) node{$_{G_{\ell}}$};
\draw (1.43,-3) node{$_{G_{q}}$};
\draw (2.1,-2.5) node{$_{G_{o}}$};
\draw (2.57,-3) node{$_{G_{m}}$};
\draw (3.1,-2.5) node{$_{G_{\ell} }$};
\draw (3.8,-3) node{$_{G_{n}}$};
\draw (4.8,-2.5) node{$_{G_{p}}$};

\draw (0.3,-2) node{$_{|}$};
\draw (1.2,-2) node{$_{|}$};
\draw (1.7,-2) node{$_{|}$};
\draw (2.3,-2) node{$_{|}$};
\draw (2.8,-2) node{$_{|}$};
\draw (3.4,-2) node{$_{|}$};
\draw (4.1,-2) node{$_{|}$};

\end{tikzpicture}
\begin{center}
Canonical representation of $s+t$
\end{center}
\end{figure}

Let us denote by $\sum_{(\ell)}$ the sum taken over the pairs $(j,k)$ such that $\alpha_{j}+\beta_{k}=c_{\ell}$. It is then clear that $\mu(G_{\ell})=\sum_{(\ell)}\,\mu(A_j \cap B_k)$ and, therefore,

$$
\begin{aligned}
\int_X(s+t)\,d\mu &= \sum_{\ell=1}^{p} c_{\ell}\mu(G_{\ell}) =\sum_{\ell=1}^{p}\sum_{(\ell)} c_{\ell}\mu(A_j \cap B_k)=\sum_{\ell=1}^{p}\sum_{(\ell)} (\alpha_j + \beta_k)\mu(A_j \cap B_k)\\
&=\sum_{j=1}^{N}\sum_{k=1}^{M} (\alpha_j + \beta_k)\mu(A_j \cap B_k)\\
&=\sum_{j=1}^{N}\alpha_{j}\sum_{k=1}^{M}\mu(A_j \cap B_k) + \sum_{k=1}^{M}\beta_k \sum_{j=1}^{N}\mu(A_j \cap B_k).\\
\end{aligned}
$$

Since $X=\bigcup_{j=1}^{N} A_j = \bigcup_{k=1}^{M} B_k$, then $A_j=\bigcup_{k=1}^{M}(A_j \cap B_k)$ and $B_k=\bigcup_{j=1}^{N}(A_j \cap B_k)$ for each $j=1,\ldots,N$ and $k=1,\ldots,M$. Applying the additivity of $\mu$ we obtain
$$
\begin{aligned}
\int_X(s+t)\,d\mu &=\sum_{j=1}^{N}\alpha_{j}\sum_{k=1}^{M}\mu(A_j \cap B_k) + \sum_{k=1}^{M}\beta_k \sum_{j=1}^{N}\mu(A_j \cap B_k)\\
&=\sum_{j=1}^{N} \alpha_j \mu(A_j)  + \sum_{k=1}^{M} \beta_k \mu(B_k) = \int_X s\,d\mu + \int_X t\,d\mu,
\end{aligned}
$$
as stated.
\end{proof}

Arguing by induction and using Proposition \ref{52}, it is possible to prove that for any finite number $s_1,\ldots,s_n$ of functions in $\mathbb{S}^{+}(X,\mathsf{S})$ and nonnegative real values $c_1,\ldots,c_n$, $\int_X \sum_{j=1}^{n}c_js_j\,d\mu = \sum_{j=1}^{n} c_j\int_X s_j$ [Exercise \ref{E57}].

We shall now see that the definition of the integral does not depend on the canonical representation of the nonnegative $\mathsf{S}$-simple function. 

\begin{theorem} \label{53}
Let $s\in \mathbb{S}^{+}(X,\mathsf{S})$ be such that $s=\sum_{j=1}^{M}\,\beta_{j}\,\chi_{B_j}$ with the $\beta_j> 0$ not necessarily distinct and the $B_j \in \mathsf{S}$ not necessarily disjoint. Then,
$$
\int_X s\,d\mu=\sum_{j=1}^{M}\beta_j\,\mu(B_j).
$$
\end{theorem}

\begin{proof}
For each $j=1,\ldots,M$, define the function $s_j:X \to \mathbb{R}$ by
$$s_j(x) :=\left\{
\begin{array}{l c l}
\beta_j & &\mbox{if}\,\, x\in B_j\\
0       & &\mbox{if}\,\, x\in X\smallsetminus B_j. 
\end{array}
\right.
$$

It is possible to write, for each $j=1,\ldots,M$, the function $s_j$ as $s_j=\beta_j\,\chi_{B_j} + 0\,\chi_{X\smallsetminus B_j}$ where $\{B_j,X\smallsetminus B_j\}$ is a partition of $X$ with elements in $\mathsf{S}$ and $\beta_{j} \neq 0$. Therefore, $s_j \in \mathbb{S}^{+}(X,\mathsf{S})$ and $\int_{X}s_{j}\,d\mu=\beta_{j}\mu(B_{j})$ for each $j=1,\ldots,M$.

It is clear that $s=\sum_{j=1}^{M}s_j$ and from the previous proposition and induction [Exercise \ref{E53}] it follows that
$$
\int_X s\,d\mu = \sum_{j=1}^{M}\int_X s_j\,d\mu = \sum_{j=1}^{M} \beta_j\,\mu(B_j),
$$
as stated.
\end{proof}

The following result establishes the monotonicity of the Lebesgue integral for simple functions as an immediate consequence of the previous result.

\begin{proposition}[Monotonicity of the Lebesgue integral for simple functions] \label{54}
For every $s,t \in \mathbb{S}^{+}(X,\mathsf{S})$ such that $s \leq t$, one has $\int_X s\,d\mu \leq \int_X t\,d\mu$.
\end{proposition}

\begin{proof}
Let $s=\sum_{j=1}^{N} \alpha_j \chi_{A_j}$ and $t=\sum_{k=1}^{M} \beta_k \chi_{B_k}$ be the canonical representations of the functions $s,t \in \mathbb{S}^{+}(X,\mathsf{S})$.

\begin{figure}[ht!]
\begin{minipage}[r]{0.45\textwidth}
\begin{center}
\begin{tikzpicture}[xscale=1.22,yscale=0.85]
	\draw[->,gray] (0,0) -- (4.5,0); \draw [->, gray] (0,0) -- (0,4); 
\draw (4.8,0.25) node[below] {$X$}; \draw (0,4) node[right] {$\mathbb{R}$}; 

\draw (0,0.5) node{$-$}; \draw (0,0.5) node[left] {$_{\alpha_{i}}$};
\draw [thick] (0.2,0.5)--(1.45,0.5);
\draw (0.2,0.5) node{$_{\bullet}$}; 
\draw (1.5,0.5) node{$_{\circ}$};

\draw (0,1) node{$-$}; \draw (0,1) node[left] {$_{\alpha_{n}}$};
\draw [ thick] (1.5,1)--(2.45,1);
\draw (1.5,1) node{$_{\bullet}$}; 
\draw (2.5,1) node{$_{\circ}$};

\draw (0,1.5) node{$-$}; \draw (0,1.5) node[left] {$_{\alpha_{j}}$};
\draw [ thick] (2.5,1.5)--(3.43,1.5);
\draw (2.5,1.5) node{$_{\bullet}$}; 
\draw (3.5,1.5) node{$_{\circ}$};

\draw (0,2) node{$-$}; \draw (0,2) node[left] {$_{\alpha_{\ell}}$};
\draw [ thick] (3.5,2)--(4.5,2);
\draw (3.5,2) node{$_{\bullet}$}; 

\draw (0,2.5) node{$-$}; \draw (0,2.5) node[left] {$_{\beta_{k}}$};
\draw [ thick, red] (0.2,2.5)--(1.92,2.5);
\draw (0.2,2.5) node{$_{\bullet}$}; 
\draw (2,2.5) node{$_{\circ}$};

\draw (0,3) node{$-$}; \draw (0,3) node[left] {$_{\beta_{m}}$};
\draw [ thick, red] (2,3)--(3.19,3);
\draw (2,3) node{$_{\bullet}$}; 
\draw (3.25,3) node{$_{\circ}$};

\draw (0,3.5) node{$-$}; \draw (0,3.5) node[left] {$_{\beta_{r}}$};
\draw [ thick, red] (3.25,3.5)--(4.5,3.5);
\draw (3.25,3.5) node{$_{\bullet}$}; 

\draw [dotted] (0.2,2.5)--(0.2,0);
\draw [dotted] (1.5,1)--(1.5,0);
\draw [dotted] (2.5,1.5)--(2.5,0);
\draw [dotted] (3.5,2)--(3.5,0);

\draw [dotted] (2,3)--(2,0);
\draw [dotted] (3.25,3.5)--(3.25,0);

\draw (0.2,0) node{$_{|}$};
\draw (1.5,0) node{$_{|}$};
\draw (2,0) node{$_{|}$};
\draw (2.5,0) node{$_{|}$};
\draw (3.25,0) node{$_{|}$};
\draw (3.5,0) node{$_{|}$};

\draw (1,-0.5) node{$_{A_{i}}$};
\draw (2,-0.5) node{$_{A_{n}}$};
\draw (3,-0.5) node{$_{A_{j}}$};
\draw (4,-0.5) node{$_{A_{\ell}}$};

\draw (1.2,-1) node{$_{B_{k}}$};
\draw (2.5,-1) node{$_{B_{m}}$};
\draw (4.3,-1) node{$_{B_{r}}$};

\end{tikzpicture}
\end{center}
\end{minipage} \hfill 
 \begin{minipage}[l]{0.55\textwidth}
\begin{center}

\begin{tikzpicture}[xscale=1.22,yscale=0.85]
	\draw[->,gray] (0,0) -- (4.5,0); \draw[->,gray] (0,0) -- (0,4); 
	
\draw (4.8,0.25) node[below] {$X$}; \draw (0,4) node[right] {$\mathbb{R}$}; 

\draw (0,0.5) node{$-$}; \draw (0,0.5) node[left] {$_{\alpha_{i}}$};
\draw [thick] (0.2,0.5)--(1.45,0.5);
\draw (0.2,0.5) node{$_{\bullet}$}; 
\draw (1.5,0.5) node{$_{\circ}$};

\draw (0,1) node{$-$}; \draw (0,1) node[left] {$_{\alpha_{n}}$};
\draw [ thick] (1.5,1)--(2.45,1);
\draw (1.5,1) node{$_{\bullet}$}; 
\draw (2.5,1) node{$_{\circ}$};

\draw (0,1.5) node{$-$}; \draw (0,1.5) node[left] {$_{\alpha_{j}}$};
\draw [ thick] (2.5,1.5)--(3.43,1.5);
\draw (2.5,1.5) node{$_{\bullet}$}; 
\draw (3.5,1.5) node{$_{\circ}$};

\draw (0,2) node{$-$}; \draw (0,2) node[left] {$_{\alpha_{\ell}}$};
\draw [ thick] (3.5,2)--(4.5,2);
\draw (3.5,2) node{$_{\bullet}$}; 

\draw (0,2.5) node{$-$}; \draw (0,2.5) node[left] {$_{\beta_{k}}$};
\draw [ thick, red] (0.2,2.5)--(1.92,2.5);
\draw (0.2,2.5) node{$_{\bullet}$}; 
\draw (2,2.5) node{$_{\circ}$};

\draw (0,3) node{$-$}; \draw (0,3) node[left] {$_{\beta_{m}}$};
\draw [ thick, red] (2,3)--(3.19,3);
\draw (2,3) node{$_{\bullet}$}; 
\draw (3.25,3) node{$_{\circ}$};

\draw (0,3.5) node{$-$}; \draw (0,3.5) node[left] {$_{\beta_{r}}$};
\draw [ thick, red] (3.25,3.5)--(4.5,3.5);
\draw (3.25,3.5) node{$_{\bullet}$}; 

\draw [dotted] (0.2,2.5)--(0.2,0);
\draw [dotted] (1.5,2.5)--(1.5,0);
\draw [dotted] (2.5,3)--(2.5,0);
\draw [dotted] (3.5,3.5)--(3.5,0);

\draw [dotted] (2,3)--(2,0);
\draw [dotted] (3.25,3.5)--(3.25,0);

\draw (0.2,0) node{$_{|}$};
\draw (1.5,0) node{$_{|}$};
\draw (2,0) node{$_{|}$};
\draw (2.5,0) node{$_{|}$};
\draw (3.25,0) node{$_{|}$};
\draw (3.5,0) node{$_{|}$};

\draw (1,-0.5) node{$_{A_{i}\cap B_{k}}$};
\draw (1.7,-1) node{$_{A_{n}\cap B_{k}}$};
\draw (2.2,-0.5) node{$_{A_{n} \cap B_{m}}$};
\draw (2.85,-1) node{$_{A_{j} \cap B_{m}}$};
\draw (3.35,-0.5) node{$_{A_{j} \cap B_{r}}$};
\draw (4,-1) node{$_{A_{\ell} \cap B_{r}}$};

\draw (1.5,2.5) node{$_{\odot}$};
\draw (2,1) node{$_{\odot}$};
\draw (2.5,3) node{$_{\odot}$};
\draw (3.25,1.5) node{$_{\odot}$};
\draw (3.5,3.5) node{$_{\odot}$};
\end{tikzpicture}
\end{center}
\end{minipage}
\end{figure}

Since $s(x) \leq t(x)$ for every $x \in X$, it follows that $\alpha_j \leq \beta_k$ whenever $A_j \cap B_k \neq \varnothing$. Consequently,
$$
\alpha_j \mu(A_j \cap B_k) \leq \beta_k \mu(A_j \cap B_k),
$$
for $j=1,\ldots,N$ and for $k=1,\ldots,M$. Thus, we write $s$ and $t$ through the following representations
$$
s=\sum_{j,k}\alpha_j \chi_{A_j \cap B_k},\qquad t=\sum_{j,k} \beta_k \chi_{A_j \cap B_k},
$$
where $\sum_{j,k}$ denotes that the sum extends over the pairs $(j,k)$ such that $\alpha_j \leq \beta_k$. 

Applying Theorem \ref{53} we obtain that 
$$
\int_X s\,d\mu = \sum_{j,k} \alpha_j \mu(A_j \cap B_k) \leq \sum_{j,k} \beta_k \mu(A_j \cap B_k) = \int_X t\,d\mu,
$$
as stated.
\end{proof}

Let us look at some examples.

\begin{example} \label{55}
Consider the measure space given by $X=\mathbb{N}$, $\mathsf{S}=\mathcal{P}(\mathbb{N})$, and $\mu=\mu^{\sharp}$ the counting measure. Let $(x_k)$ be a sequence of nonnegative real numbers.

Then, for a fixed natural number $N \in \mathbb{N}$, define $s_N :\mathbb{N} \to \mathbb{R}$ by
$$
s_N(k):=\left\{
\begin{array}{l c l}
x_k & &\mbox{if}\,\, 1 \leq k \leq N,\\
0  & & \mbox{if}\,\, N<k.
\end{array}
\right.
$$

\begin{figure}[ht!]
\centering
\begin{tikzpicture}[xscale=0.8,yscale=0.8]
	\draw[->, gray] (0,0) -- (7.5,0); \draw [->, gray] (0,0) -- (0,4); 
\draw (7.5,0) node[right] {$\mathbb{N}$}; \draw (0,4) node[right] {$\mathbb{R}$}; 

\draw (0,3) node{$_{-}$}; \draw (0,3) node[left]{$_{x_{N}}$};

\draw (0,2) node{$_{-}$}; \draw (0,2) node[left]{$_{x_{N-1}}$};

\draw (0,0.5) node{$_{-}$}; \draw (0,0.5) node[left]{$_{x_{1}}$};

\draw (1,0) node{$_{|}$}; \draw (1,-0.1) node[below]{$_{1}$}; \draw (2,-0.1) node[below]{$_{\ldots}$};

\draw (3,0) node{$_{|}$}; \draw (3,-0.1) node[below]{$_{N-1}$};

\draw (4,0) node{$_{|}$}; \draw (4,-0.1) node[below]{$_{N}$};

\draw (5,0) node{$_{|}$}; \draw (5,-0.1) node[below]{$_{N+1}$};

\draw (6,0) node{$_{|}$}; \draw (6,-0.1) node[below]{$_{N+2}$}; \draw (7,-0.1) node[below]{$_{\ldots}$};

\draw (1,0.5) node{$_{\bullet}$};
\draw (3,2) node{$_{\bullet}$};
\draw (4,3) node{$_{\bullet}$};

\draw (5,0) node{$_{\bullet}$};
\draw (6,0) node{$_{\bullet}$};

\draw (0,-0.1) node[below]{$_{0}$};
\end{tikzpicture}
\begin{center}
$s_{N}(x)$
\end{center}
\end{figure}

It is clear that $s_N \in \mathbb{S}^{+}(\mathbb{N},\mathcal{P}(\mathbb{N}))$ and we may write it as
$$
s_N=\sum_{k=1}^{N} x_{k}\,\chi_{\{k\}} + 0\chi_{\{N+1,N+2,\ldots \}}.
$$

Thus, by {\rm Theorem \ref{53}} we have
$$
\int_{\mathbb{N}} s_N\,d\mu^{\sharp} = \sum_{k=1}^{N} x_{k}\,\mu^{\sharp}(\{k\}) + 0\,\mu^{\sharp}(\{N+1,N+2,\ldots \}) = \sum_{k=1}^{N} x_{k}.
$$

This integral corresponds to the $N$-th partial sum of the series $\sum_{k=1}^{\infty} x_{k}$.
\end{example}

\begin{example} \label{56}
Let $X=[0,1]$, $\mathsf{S}=\mathcal{B}([0,1])$, and $\mu=\lambda$ the Lebesgue measure on $\mathsf{S}$. For fixed $k \in \mathbb{N}$, define the function $s_k:X \to \mathbb{R}$ as follows
$$
s_k(x):=\left\{
\begin{array}{l c l}
k & &\mbox{if}\quad x \in [0,\frac{1}{k}),\\
1  & & \mbox{if}\quad x \in [\frac{1}{k},1-\frac{1}{k}),\\
\frac{1}{k} & & \mbox{if}\quad x \in [1-\frac{1}{k},1].\\
\end{array}
\right.
$$

\begin{figure}[ht!]
\centering
\begin{tikzpicture}[xscale=1,yscale=0.8]
	\draw[->,gray] (0,0) -- (3,0); \draw [->,gray] (0,0) -- (0,4);
\draw (0,4) node[right] {$\mathbb{R}$}; 

\draw (0,3) node{$_{-}$}; \draw (0,3) node[left]{$_{k}$};
\draw (0,2) node{$_{-}$}; \draw (0,2) node[left]{$_{1}$};
\draw (0,1) node{$_{-}$}; \draw (0,1) node[left]{$_{\frac{1}{k}}$};

\draw (3,-0.1) node[below]{$_{1}$};
\draw (2,0) node{$_{|}$}; \draw (2,-0.1) node[below]{$_{\frac{k-1}{k}}$};
\draw (1,0) node{$_{|}$}; \draw (1,-0.1) node[below]{$_{\frac{1}{k}}$};
\draw (0,-0.1) node[below]{$_{0}$};
 
\draw[thick] (0,3)--(0.99,3); \draw (0,3) node{$_{\bullet}$}; \draw (1,3) node{$_{\circ}$};
\draw[thick] (1,2)--(1.99,2); \draw (1,2) node{$_{\bullet}$}; \draw (2,2) node{$_{\circ}$};
\draw[thick] (2,1)--(3,1); \draw (2,1) node{$_{\bullet}$};

\draw[dotted] (1,3)--(1,0);
\draw[dotted] (2,2)--(2,0);
\end{tikzpicture}
\begin{center}
$s_{k}(x)$
\end{center}
\end{figure}

It is clear that $s_k \in \mathbb{S}^{+}([0,1],\mathcal{B}([0,1]))$ and has the following canonical representation
$$
s_k = k\,\chi_{[0,\frac{1}{k})} + 1\,\chi_{[\frac{1}{k},1-\frac{1}{k})} + \frac{1}{k}\,\chi_{[1-\frac{1}{k},1]}.
$$

Therefore,
$$
\int_{[0,1]} s_k\,d\lambda = k\,\lambda\left([0,1/k) \right) + \lambda\left([1/k,1-1/k) \right) +\frac{1}{k}\,\lambda\left([1-1/k,1] \right)=2-\frac{2}{k}+\frac{1}{k^2}.
$$
\end{example}

\begin{definition} \label{57}
Let $(X,\mathsf{S},\mu)$ be a measure space and let $s:X\to \mathbb{R}$ be a nonnegative simple function with canonical representation $s=\sum_{j=1}^{N} \alpha_{j}\,\chi_{A_j}$. For each $A \in \mathsf{S}$, we define the \textbf{(Lebesgue) integral of $s$ on $A$ with respect to the measure $\mu$}, denoted by $\int_{A}\,s\,d\mu$, as the extended real number given by
$$
\int_A\, s\,d\mu := \int_X (s\cdot\chi_{A})\,d\mu :=\sum_{j=1}^{N} \alpha_{j}\,\mu(A_j \cap A).
$$
\end{definition}

The following result shows a well-known interpretation from our differential and integral calculus courses of the concept of integral.

\begin{theorem} \label{58}
Let $(X,\mathsf{S},\mu)$ be a measure space and let $s \in \mathbb{S}^{+}(X,\mathsf{S})$ be fixed. The function $\mu_{s}:\mathsf{S} \to \overline{\mathbb{R}}$ given by $\mu_{s}(A):=\int_{A} s\,d\mu$ defines a measure on $(X,\mathsf{S})$.
\end{theorem}

\begin{proof}
Let $s=\sum_{j=1}^{N} \alpha_{j}\,\chi_{A_j}$ be the canonical representation of $s$. From Definition \ref{57} it follows that
$$
\mu_{s}(A):=\int_{A} \,s\,d\mu =\sum_{j=1}^{N} \alpha_{j}\,\mu(A_j \cap A) = \sum_{j=1}^{N}\alpha_{j}\mu^{A_j}(A)\quad \mbox{for every}\,\,\,A \in \mathsf{S}.
$$
where $\mu^{A_j}$ is the contraction measure of $\mu$ to $A_j$ (see Example \ref{48}). Consequently, $\mu_{s}$ is a nonnegative linear combination of measures on $(X,\mathsf{S})$. Applying Proposition \ref{413}, we conclude that $\mu_{s}$ is a measure on $(X,\mathsf{S})$.
\end{proof}

\section{The integral of nonnegative measurable functions}

Let $f \in \overline{\mathbb{M}}^{+}(X,\mathsf{S})$. We denote by $\underline{\mathcal{S}}(f)$ the set of all functions $s:X \to \mathbb{R}$ in $\mathbb{S}^{+}(X,\mathsf{S})$ such that $s \leq f$. Observe that $\underline{\mathcal{S}}(f)$ is nonempty since the constant function equal to $0$ belongs to $\mathbb{S}^{+}(X,\mathsf{S})$ and $0 \leq f$ on $X$.

The above allows us to give the following definition.

\begin{definition} \label{59} \index{integral!Lebesgue!of a nonnegative measurable function}
Let $f \in \overline{\mathbb{M}}^{+}(X,\mathsf{S})$. We define the \textbf{(Lebesgue) integral of $f$ on $X$ with respect to $\mu$} as the extended real number given by
$$
\int_X f\,d\mu :=\sup\left\{ \int_X s\,d\mu \,:\,s \in \underline{\mathcal{S}}(f) \right\}.
$$
\end{definition}

It is easy to prove that if $f,g\in \overline{\mathbb{M}}^{+}(X,\mathsf{S})$ are such that $f \leq g$, then $\int_X f\,d\mu \leq \int_X g\,d\mu$ [Exercise \ref{E58}]. In particular, $0 \leq \int_{X}f\,d\mu$ for every $f \in \overline{\mathbb{M}}^{+}(X,\mathsf{S})$.

We now prove that the previous definition extends the notion of integral given in \ref{51} to a broader family of measurable functions.

\begin{proposition} \label{510}
Let $\phi \in \mathbb{S}^{+}(X,\mathsf{S})$ with canonical representation $\phi:=\sum_{j=1}^{N}\,\alpha_{j}\,\chi_{A_j}$. Then,
$$
\sup\left\{ \int_X s\,d\mu \,:\,s \in \underline{\mathcal{S}}(\phi) \right\}=\sum_{j=1}^{N}\,\alpha_{j}\,\mu(A_j).
$$
\end{proposition}

\begin{proof}
Since $\phi$ is a nonnegative $\mathsf{S}$-simple function, it is clear that $\phi \in \underline{\mathcal{S}}(\phi)$ and, therefore,
$$
\sum_{j=1}^{N}\,\alpha_{j}\,\mu(A_j)\leq \sup\left\{ \int_X s\,d\mu \,:\,s \in \underline{\mathcal{S}}(\phi) \right\}.
$$

Conversely: let $s \in \underline{\mathcal{S}}(\phi)$ be arbitrary. Since $s \leq \phi$, from Proposition \ref{54} we obtain that $\int_X s\,d\mu \leq \sum_{j=1}^{N}\,\alpha_{j}\,\mu(A_j)$. Consequently,
$$
\sup\left\{ \int_X s\,d\mu \,:\,s \in \underline{\mathcal{S}}(\phi) \right\} \leq \sum_{j=1}^{N}\,\alpha_{j}\,\mu(A_j).
$$

This concludes the proof.
\end{proof}

Given a bounded function $f \in \mathbb{M}^{+}(X,\mathsf{S})$, we denote by $\overline{\mathcal{S}}(f)$ the set of all functions $t:X \to \mathbb{R}$ in $\mathbb{S}^{+}(X,\mathsf{S})$ such that $f \leq t$. Observe again that $\overline{\mathcal{S}}(f)$ is nonempty since, as $f$ is bounded, there exists $M>0$ such that $f \leq M$ on $X$. Then, it is clear that $M = M\chi_{X} + 0\chi_{\varnothing}$ so that $M \in\mathbb{S}^{+}(X,\mathsf{S})$.

The following result provides necessary and sufficient conditions to characterize some measurable functions and their integral.

\begin{theorem} \label{511}
Let $(X,\mathsf{S},\mu)$ be a finite measure space and let $f:X \to \mathbb{R}$ be nonnegative and bounded.
\begin{itemize}
\item[(a)] If $f \in \mathbb{M}^{+}(X,\mathsf{S})$, then
\begin{eqnarray} \label{F51}
\inf\left\{ \int_X t\,d\mu\,:\, t \in \overline{\mathcal{S}}(f) \right\}=\sup \left\{ \int_X s\,d\mu\,:\, s \in \underline{\mathcal{S}}(f)\right\}.
\end{eqnarray}
\item[(b)] If the previous equality holds for $f$ and $(X,\mathsf{S},\mu)$ is complete, then $f \in \mathbb{M}^{+}(X,\mathsf{S})$.
\end{itemize}
\end{theorem}

\begin{proof}
\textit{(a):} For arbitrary $s \in \underline{\mathcal{S}}(f)$ and $t \in \overline{\mathcal{S}}(f)$ we have $s \leq t$. From Proposition \ref{54} it follows that
$$
\int_X s\,d\mu \leq \int_X t\,d\mu.
$$

Consequently,
$$
\sup \left\{ \int_X s\,d\mu\,:\, s \in \underline{\mathcal{S}}(f)\right\} \leq \inf\left\{ \int_X t\,d\mu\,:\, t \in \overline{\mathcal{S}}(f) \right\}.
$$

Let $M>0$ be such that $0 \leq f(x) < M$ for every $x \in X$. Fix $k \in \mathbb{N}$. For each $j \in \{1,\ldots,k \}$ define the following sets
$$
A_{j}(k):=\left\{ x \in X\,:\, \frac{(j-1)M}{k} \leq f(x) < \frac{jM}{k}\right\},
$$
which belong to $\mathsf{S}$ and form a partition of $X$.

Therefore, the functions
$$
s_{k}:=\frac{M}{k}\sum_{j=1}^{k} (j-1)\chi_{A_{j}(k)}\quad \mbox{and}\quad t_{k}:=\frac{M}{k}\sum_{j=1}^{k} j\,\chi_{A_{j}(k)}
$$
are nonnegative $\mathsf{S}$-simple functions such that $s_{k} \in \underline{\mathcal{S}}(f)$ and $t_{k} \in \overline{\mathcal{S}}(f)$.

\begin{figure}[ht!]
\centering
\begin{tikzpicture}[yscale=1.5]
\draw[->,gray] (-0.5,0) -- (6,0);
\draw [->,gray] (0,-0.2) -- (0,2);
\draw[thick] plot[smooth] coordinates
{(0.00,0.00)(1,1)(3,0.51)(4,1.5)(5,1.5)};

\draw (6,0) node [right]{$X$}; \draw (0,2) node[right]{$\mathbb{R}$};
\draw[dotted] (2.8,0.5)--(0,0.5);
\draw[dotted] (3.75,0)--(3.75,1.25);
\draw[dotted] (0,1.25)--(3.75,1.25);
\draw[dotted] (2.8,0)--(2.8,0.5);

\draw (0,1.25) node{$_{-}$};
\draw (0,0.5) node{$_{-}$};
\draw (0,1.75) node{$_{-}$}; \draw (0,1.75) node[left]{$_{M}$};

\draw (0,1.25) node[left]{$_{\frac{jM}{k}}$};
\draw (0,0.5) node[left]{$_{\frac{(j-1)M}{k}}$};

\draw (2.8,0) node{$_{|}$};
\draw (3.75,0) node{$_{|}$};

\draw[thick,red] (2.8,0.5)--(3.74,0.5);
\draw[thick,blue] (2.8,1.25)--(3.74,1.25);
\draw (2.8,0.5) node{$_{\bullet}$};
\draw (2.8,1.25) node{$_{\bullet}$};

\draw (3.75,0.5) node{$_{\circ}$};
\draw (3.75,1.25) node{$_{\circ}$};


\draw (3.25,0) node[below]{$_{A_{j}(k)}$};
\end{tikzpicture}
\end{figure}

Then,
$$
\begin{aligned}
0 &\leq \inf\left\{ \int_X t\,d\mu\,:\, t \in \overline{\mathcal{S}}(f) \right\} - \sup\left\{ \int_X s\,d\mu \,:\,s \in \underline{\mathcal{S}}(f) \right\}\\
 &\leq \int_X t_k\,d\mu - \int_X s_k\,d\mu \\
 &=\frac{M}{k}\left(\sum_{j=1}^{k} j\,\mu(A_j(k))-\sum_{j=1}^{k} (j-1)\,\mu(A_j(k))  \right)\\
 &=\frac{M}{k}\sum_{j=1}^{k} \mu(A_{j}(k))=\frac{M}{k}\,\mu(X).
\end{aligned}
$$

Taking the limit as $k \to \infty$ in the previous inequality we conclude that 
$$
0 \leq \inf\left\{ \int_X t\,d\mu\,:\, t \in \overline{\mathcal{S}}(f) \right\} - \sup\left\{ \int_X s\,d\mu \,:\,s \in \underline{\mathcal{S}}(f) \right\} \leq 0,
$$
which is the desired equality.

\textit{(b):} Since equality (\ref{F51}) holds, then for each $k \in \mathbb{N}$ there exist $s_{k} \in \underline{\mathcal{S}}(f)$ and $t_{k}\in \overline{\mathcal{S}}(f)$ such that
$$
\int_X t_k\,d\mu - \int_X s_k\,d\,\mu < \frac{1}{k}.
$$

Let $t_{\ast}:=\displaystyle\inf_{k \geq 1} t_{k}$ and $s^{\ast}:=\displaystyle\sup_{k \geq 1} s_{k}$. We may assume that $(t_k(x))$ is a bounded sequence for every $x \in X$ since $f$ is bounded. Theorem \ref{321} ensures that $s^{\ast}, t_{\ast} \in \mathbb{M}^{+}(X,\mathsf{S})$ and it is also clear that $s^{\ast} \leq f \leq t_{\ast}$.

For each $j \in \mathbb{N}$, define the following sets in $\mathsf{S}$
$$
\begin{aligned}
N_{j}&:=\left\{ x \in X \,:\,t_{\ast}(x)-s^{\ast}(x) \geq \tfrac{1}{j}\right\}, \\
N &:=\left\{x \in X \,:\,t_{\ast}(x)-s^{\ast}(x) > 0\right\}.
\end{aligned}
$$

In fact, $N=\bigcup_{j=1}^{\infty}\,N_{j}$. From the definition of $t_{\ast}$ and $s^{\ast}$ it is clear that, for each $j \in \mathbb{N}$, the set $N_{j}$ is contained in the measurable set
$$
N_{j}(k):=\left\{ x \in X \,:\,t_{k}(x)-s_{k}(x) \geq \tfrac{1}{j}\right\}\quad \mbox{for every}\,\,k \in \mathbb{N}.
$$

On the other hand, $1 \leq j (t_k(x)-s_k(x))$ if and only if $x \in  N_{j}(k)$ for every $k \in \mathbb{N}$. Then,
$$
\chi_{N_{j}(k)} \leq {j}\,\chi_{N_{j}(k)}(t_k-s_k)\quad \text{for every}\,\,k \in \mathbb{N}.
$$

Thus,
$$
\mu(N_{j}(k)) = \int_X \chi_{N_{j}(k)}\,d\mu \leq {j} \int_X \chi_{N_{j}(k)}(t_k-s_k)\,d\mu \leq {j}\int_X (t_k-s_k)\,d\mu 
$$
for every $k \in \mathbb{N}$. Since $t_k=(t_k-s_k)+s_k$ for every $k \in \mathbb{N}$, applying the linearity of the integral for simple functions, we obtain 
$$
\mu(N_{j}(k)) \leq j\int_X (t_k-s_k)\,d\mu = j\left(\int_X t_k\,d\mu  -\int_X s_k\,d\mu \right) < \frac{j}{k}\quad \text{for every}\,\,k \in \mathbb{N}.
$$

Consequently, $\mu(N_{j}) < \frac{j}{k}$ for every $k \in \mathbb{N}$ and, therefore, $\mu(N_{j})=0$ for every $j \in \mathbb{N}$, which establishes that $\mu(N)=0$.

This implies that $t^{\ast}=s^{\ast}$ $\mu$ a.e. and, consequently, $t^{\ast}=f$ $\mu$ a.e. Since $(X,\mathsf{S},\mu)$ is complete, Proposition \ref{435} ensures that $f$ is $\mathsf{S}$-measurable.
\end{proof}

In the previous theorem, if $(X,\mathsf{S},\mu)$ is not complete, then $f \in \mathbb{M}^{+}(X,\overline{\mathsf{S}})$ where $\overline{\mathsf{S}}$ is the $\mu$-completion of $\mathsf{S}$.

\begin{definition} \label{512}
Let $f \in \overline{\mathbb{M}}^{+}(X,\mathsf{S})$. For each $A \in \mathsf{S}$, we define \textbf{the (Lebesgue) integral of $f$ on $A$ with respect to $\mu$} as the extended real number given by
$$
\int_{A} f\,d\mu :=\int_X (f\cdot\chi_{A})\,d\mu.
$$
\end{definition}

We now prove some important properties of the integral over arbitrary measurable sets.

\begin{proposition} \label{513}
Let $f \in \overline{\mathbb{M}}^{+}(X,\mathsf{S})$.
\begin{itemize}
\item[(a)] If $A,B \in \mathsf{S}$ are such that $A \subset B$, then $\int_{A}\,f\,d\mu \leq \int_{B}\, f\,d\mu$.
\item[(b)] If $\mu(A)=0$, then $\int_{A} f\,d\mu=0$.
\end{itemize}
\end{proposition}

\begin{proof}
\textit{(a):} It follows immediately from the monotonicity of the integral and the fact that $f\cdot\chi_{A} \leq f\cdot\chi_{B}$ on $X$ for any $A,B \in \mathsf{S}$ with $A \subset B$.

\textit{(b):} Let $s \in \underline{\mathcal{S}}(f\cdot \chi_{A})$ be arbitrary with canonical representation equal to $\sum_{j=1}^{N} \alpha_{j}\,\chi_{A_j}$. Since $0 \leq s \leq f\cdot \chi_{A}$ on $X$, then $s=s\cdot \chi_{A}$ on $X$. Thus,
$$
0 \leq \int_{X} s\,d\mu = \int_{A} s\,d\mu = \sum_{j=1}^{N} \alpha_{j}\,\mu(A_j \cap A) \leq \sum_{j=1}^{N} \alpha_{j}\,\mu(A) =0
$$
so necessarily $\int_X s\,d\mu=0$. Consequently, $\int_{A}\,f\,d\mu =0$.
\end{proof}

We conclude this section with some examples.

\begin{example} \label{514}
Let $(X,\mathsf{S})$ be a measurable space and let $x_0 \in X$ be fixed. Let $\delta_{x_{0}}$ be the unit measure concentrated at $x_{0}$. Then, for each $f \in \mathbb{M}^{+}(X,\mathsf{S})$ we have
$$
\int_X f\,d\delta_{x_{0}} = f(x_{0}).
$$
\end{example}

\begin{proof}
Consider the following two cases.

{\scshape Case 1:}\quad $f \in \mathbb{S}^{+}(X,\mathsf{S})$.

If $\alpha_{1},\ldots,\alpha_{N}$ are the distinct values taken by $f$, Proposition \ref{315} ensures that the sets $A_j:=f^{-1}(\{\alpha_j\})$ for each $j=1,\ldots,N$ belong to $\mathsf{S}$ and form a partition of $X$. 

Since $f=\sum_{j=1}^{N}\alpha_j\,\chi_{A_j}$ is the canonical representation and $X$ is the disjoint union of $A_1,\ldots,A_N$, let $j_{\ast}$ be the unique natural number in $\{1,\ldots,N\}$ such that $x_{0} \in A_{j_{\ast}}$. Then,
$$
\int_X f\,d\delta_{x_{0}}=\sum_{j=1}^{N} \alpha_j\,\delta_{x_{0}}(A_j) = \alpha_{j_{\ast}}\,\delta_{x_{0}}(A_{j_{\ast}}) = \alpha_{j_{\ast}} =f(x_0)
$$
since $\delta_{x_{0}}(A_{j_{\ast}}) = \delta_{x_{0}}(f^{-1}(\{ \alpha_{j_{\ast}}\}))=1$ if and only if $x_{0} \in f^{-1}(\{ \alpha_{j_{\ast}} \})$ if and only if $f(x_0)=\alpha_{j_{\ast}}$.

{\scshape Case 2:}\quad $f \in \mathbb{M}^{+}(X,\mathsf{S})$.

The previous case ensures that $\int_{X} s\,d\delta_{x_0}=s(x_0)$ for every $s \in \underline{\mathcal{S}}(f)$. Since $0\leq s(x_0) \leq f(x_0)$ for every $s \in \underline{\mathcal{S}}(f)$, then
$$
\sup\left\{ \int_{X} s\,d\delta_{x_0}\,:\, s \in \underline{\mathcal{S}}(f)\right\}\leq f(x_0).
$$

Conversely: define the set $A_0:=f^{-1}(\{ f(x_0) \})$. Since $\{ f(x_0) \} \in \mathcal{B}(\mathbb{R})$ and $f$ is $\mathsf{S}$-measurable, then $A_0 := f^{-1}(\{ f(x_0)\}) \in \mathsf{S}$. Proposition \ref{315} ensures that the function $s_0 :X \to \mathbb{R}$ given by
$$
s_0:=f(x_0)\cdot\chi_{A_0} + 0\cdot \chi_{X\smallsetminus A_0}
$$
is nonnegative $\mathsf{S}$-simple and satisfies $s_0 \leq f$ on $X$. Consequently, $s_0 \in \underline{\mathcal{S}}(f)$ and we have
$$
\int_X s_0\,d\delta_{x_{0}} = f(x_0)\,\delta_{x_{0}}(A_0) = f(x_0). 
$$

\begin{figure}[ht!]
\begin{minipage}[l]{0.5\textwidth}
\begin{center}
\begin{tikzpicture}[xscale=1.5, yscale=1.5]
\draw[->,gray] (0,0) -- (3,0);
\draw [->,gray] (0,0) -- (0,2);
\draw[-,thick] (0,0)--(1,1);
\draw[-,thick] (1,1)--(2,1);
\draw[-,thick] (2,1)--(3,2);
\draw[dotted] (0,1)--(1,1);
\draw[dotted] (1.5,0)--(1.5,1);
\draw(1.5,-0.05)node[below]{$_{x_0}$};
\draw(1.5,0)node{$_{|}$};

\draw(0,1)node[left]{$_{f(x_0)}$};
\draw(0,1)node{$_{-}$};
\draw(1.5,1.25)node{$_{f:X \to \mathbb{R}}$};
\end{tikzpicture}
\end{center}
\end{minipage} \hfill 
 \begin{minipage}[r]{0.5\textwidth}
\begin{center}
\begin{tikzpicture}[xscale=1.5, yscale=1.5]
\draw[->,gray] (0,0) -- (3,0);
\draw [->,gray] (0,0) -- (0,2);
\draw[-,thick] (0,0)--(1,0);
\draw[-,thick] (1,1)--(2,1);
\draw[-,thick] (2,0)--(3,0);
\draw[dotted] (1,0)--(1,1);
\draw[dotted] (2,0)--(2,1);
\draw(1.5,-0.05)node[below]{$_{x_0}$};
\draw(1.5,0)node{$_{|}$};

\draw(0,1)node[left]{$_{f(x_0)}$};
\draw(0,1)node{$_{-}$};
\draw(1.5,1.25)node{$_{s_0:X \to \mathbb{R}}$};
\end{tikzpicture}
\end{center}
\end{minipage}
\end{figure}

Therefore,
$$
f(x_0) \leq \sup\left\{ \int_{X} s\,d\delta_{x_0}\,:\, s \in \underline{\mathcal{S}}(f)\right\}.
$$

This proves that $\int_X f\,d\delta_{x_{0}} = f(x_0)$.
\end{proof}

From the previous example it is immediate to conclude that the Lebesgue integral of a function $f \in \mathbb{M}^{+}(X,\mathsf{S})$, over any subset $A \in \mathsf{S}$, with respect to the unit measure concentrated at $x_{0}$, $\delta_{x_{0}}$, is equal to $f(x_{0})\cdot\chi_{A}(x_{0})$.

\begin{example} \label{515}
Let $(\mathbb{N},\mathcal{P}(\mathbb{N}),\mu^{\sharp})$. For every nonnegative function $f:\mathbb{N} \to \mathbb{R}$ we have
$$
\int_{A} f\,d\mu^{\sharp} =\left\{
\begin{array}{lcl}
0 & & \mbox{if}\,\,\,A=\varnothing,\\
\sum_{k \in A} f(k) & & \mbox{if}\,\,\,A\neq \varnothing,\\
\end{array}
\right.
$$
\end{example}

\begin{proof}
Let $f \in \mathbb{M}^{+}(\mathbb{N},\mathcal{P}(\mathbb{N}))$ and $A \in \mathcal{P}(\mathbb{N})$ be arbitrary. 

Consider the following cases:

{\scshape Case 1:}\quad $A=\varnothing$.

The function $\chi_{A}$ is the constant function equal to $0$, so that
$$
\int_{A} f\,d\mu^{\sharp} = \int_{\mathbb{N}} (f \cdot \chi_{A})\,d\mu^{\sharp} = \int_{\mathbb{N}} 0\,d\mu^{\sharp} = 0.
$$

{\scshape Case 2:}\quad $A$ is a finite set.

Let $x_1,\ldots,x_N$ be the distinct values of the set $A$. It is then clear that the function $f\cdot \chi_{A}$ is $\mathcal{P}(\mathbb{N})$-simple and nonnegative. Moreover, its canonical representation is given by the expression
$$
(f\cdot \chi_{A} )(k)=\sum_{j=1}^{N} f(x_j)\cdot\chi_{\{x_j\}} + 0\chi_{\mathbb{N}\smallsetminus A}.
$$

\begin{figure}[ht!]
\begin{minipage}[l]{0.5\textwidth}
\begin{center}
\begin{tikzpicture}[xscale=1.75, yscale=1.5]
\draw[->,gray] (0,0) -- (2.5,0); \draw(2.5,0)node[right]{$_{\mathbb{N}}$};
\draw [->,gray] (0,0) -- (0,2.25); \draw(0,2.25)node[left]{$_{\mathbb{R}}$};

\draw(0.25,2)node{$_{\bullet}$};
\draw(0.5,1.75)node{$_{\bullet}$};
\draw(0.75,1.5)node{$_{\bullet}$};
\draw(1,1.25)node{$_{\bullet}$};
\draw(1.25,1)node{$_{\bullet}$};
\draw(1.5,.75)node{$_{\bullet}$};
\draw(1.75,.5)node{$_{\bullet}$};
\draw(2,.5)node{$_{\bullet}$};
\draw(2.25,.5)node{$_{\bullet}$};

\draw(0,2)node{$_{-}$};
\draw(0,1.75)node{$_{-}$};
\draw(0,1.5)node{$_{-}$};
\draw(0,1.25)node{$_{-}$};
\draw(0,1)node{$_{-}$};
\draw(0,.75)node{$_{-}$};
\draw(0,.5)node{$_{-}$};
\draw(0,0.5)node{$_{-}$};
\draw(0,0.5)node{$_{-}$};

\draw[dotted](0.25,2)--(0.25,0);
\draw[dotted](0.5,1.75)--(0.5,0);
\draw[dotted](0.75,1.5)--(0.75,0);
\draw[dotted](1,1.25)--(1,0);
\draw[dotted](1.25,1)--(1.25,0);
\draw[dotted](1.5,.75)--(1.5,0);
\draw[dotted](1.75,.5)--(1.75,0);
\draw[dotted](2,.5)--(2,0);
\draw[dotted](2.25,.5)--(2.25,0);

\draw(0.75,-0.05)node[below]{$_{x_1}$};
\draw(0.75,0)node{$_{|}$};

\draw(1,-0.05)node[below]{$_{x_2}$};
\draw(1,0)node{$_{|}$};

\draw(1.5,-0.05)node[below]{$_{x_3}$};
\draw(1.5,0)node{$_{|}$};

\draw(1.75,-0.05)node[below]{$_{x_4}$};
\draw(1.75,0)node{$_{|}$};

\end{tikzpicture}
\end{center}
\end{minipage} \hfill 
 \begin{minipage}[r]{0.5\textwidth}
\begin{center}
\begin{tikzpicture}[xscale=1.75, yscale=1.5]
\draw[->,gray] (0,0) -- (2.5,0); \draw(2.5,0)node[right]{$_{\mathbb{N}}$};
\draw [->,gray] (0,0) -- (0,2.25); \draw(0,2.25)node[left]{$_{\mathbb{R}}$};

\draw(0.25,0)node{$_{\bullet}$};
\draw(0.5,0)node{$_{\bullet}$};
\draw(0.75,1.5)node{$_{\bullet}$};
\draw(1,1.25)node{$_{\bullet}$};
\draw(1.25,0)node{$_{\bullet}$};
\draw(1.5,.75)node{$_{\bullet}$};
\draw(1.75,.5)node{$_{\bullet}$};
\draw(2,0)node{$_{\bullet}$};
\draw(2.25,0)node{$_{\bullet}$};

\draw(0,2)node{$_{-}$};
\draw(0,1.75)node{$_{-}$};
\draw(0,1.5)node{$_{-}$};
\draw(0,1.25)node{$_{-}$};
\draw(0,1)node{$_{-}$};
\draw(0,.75)node{$_{-}$};
\draw(0,.5)node{$_{-}$};
\draw(0,0.5)node{$_{-}$};
\draw(0,0.5)node{$_{-}$};

\draw[dotted](0.75,1.5)--(0.75,0);
\draw[dotted](1,1.25)--(1,0);
\draw[dotted](1.5,.75)--(1.5,0);
\draw[dotted](1.75,.5)--(1.75,0);

\draw(0.75,-0.05)node[below]{$_{x_1}$};
\draw(0.75,0)node{$_{|}$};

\draw(1,-0.05)node[below]{$_{x_2}$};
\draw(1,0)node{$_{|}$};

\draw(1.5,-0.05)node[below]{$_{x_3}$};
\draw(1.5,0)node{$_{|}$};

\draw(1.75,-0.05)node[below]{$_{x_4}$};
\draw(1.75,0)node{$_{|}$};
\end{tikzpicture}
\end{center}
\end{minipage}
\end{figure}

Therefore,
$$
\int_{A} f\,d\mu^{\sharp}= \int_{\mathbb{N}} (f\cdot \chi_{A})\,d\mu^{\sharp} = \sum_{j=1}^{N} f(x_j).
$$

{\scshape Case 3:}\quad $A$ is an infinite set.

For each $j \in \mathbb{N}$, define the function $s_j:\mathbb{N} \to \mathbb{R}$ by
$$
s_j(k):=\left\{
\begin{array}{lcl}
 f(k) & & \mbox{if}\,\,\, k \in A \cap \{1,\ldots,j\},\\
 0  & & \mbox{if}\,\,\, A\cap \{1,\ldots,j\}=\varnothing,\\
\end{array}
\right.
$$

It is clear that $s_j \in \mathbb{S}^{+}(\mathbb{N},\mathcal{P}(\mathbb{N}))$ and $0\leq s_j \leq f$ on $\mathbb{N}$, for every $j \in \mathbb{N}$. Consequently, $s_j \in \underline{\mathcal{S}}(f \cdot \chi_{A})$ for every $j \in \mathbb{N}$ and it follows from the previous case that
$$
\int_{\mathbb{N}} s_j\,d\mu^{\sharp} =\sum_{k \in A \cap \{1,\ldots,j \}} f(k)\leq \int_{A} f\,d\mu^{\sharp}\quad \text{for every }\, j \in \mathbb{N}.
$$

Therefore,
$$
\sum_{k \in A} f(k) \leq \int_{A} f\,d\mu^{\sharp}.
$$

For the reverse inequality it suffices to consider only the case when $\sum_{k\,\in\,A} f(k)<+\infty$. 

Let $s \in \underline{\mathcal{S}}(f\cdot \chi_{A})$ be arbitrary. Since $s$ takes only finitely many values, $0 \leq s \leq (f \cdot \chi_{A})$, and we are assuming that the set $A$ is infinite, then at least one of the values of $s$ is attained infinitely many times. That is, there exist $\alpha \geq 0$ and an infinite subset $B \subset A$ such that $s(k)=\alpha$ for every $k \in B$. Therefore,
$$
0 \leq \alpha\,|B| = \sum_{k\,\in\,B} \alpha =\sum_{k\,\in\,B} s(k) \leq \sum_{k\,\in\,A} f(k) < +\infty
$$
which implies that $\alpha=0$. Consequently, every positive value of $s$ is attained only finitely many times, so there exists $j_0 \in \mathbb{N}$ such that $s(k)=0$ for every $k >j_0$. Thus $0\leq s \leq s_{j_0}$ and this implies that
$$
\int_{\mathbb{N}} s\,d\mu^{\sharp} \leq \int_{\mathbb{N}} s_{j_0}\,d\mu^{\sharp} = \sum_{k \in A \cap \{1,\ldots,j_{0} \}} f(k) \leq \sum_{k\,\in\,A} f(k).
$$

Since $s \in \underline{\mathcal{S}}(f \cdot \chi_{A})$ was arbitrary, we conclude that $\int_{A}f\,d\mu^{\sharp} \leq \sum_{k\,\in\,A} f(k).$

This proves the result.
\end{proof}

The previous examples show how difficult it can be to compute integrals by directly applying Definition \ref{59}. In the next section we shall see some important results that will allow us, in certain cases, to compute integrals efficiently.

\section{Convergence theorems}

We begin this section by first considering the following example.

\begin{example} \label{516}
Let $X=[0,1)$, $\mathsf{S}=\mathcal{B}([0,1))$, and let $\mu=\lambda$ be the Lebesgue measure on $\mathsf{S}$. Define the function $f_k:X \to \mathbb{R}$ by $f_k(x):=\chi_{[\frac{1}{k},1)}(x)$.

\begin{figure}[ht!]
\begin{minipage}[c]{0.5\textwidth}
\begin{center}
\begin{tikzpicture}[xscale=0.75,yscale=0.75]
	\draw[->,gray] (0,0) -- (4,0); \draw [->,gray] (0,0) -- (0,4); 
	\draw[-,gray] (-0.2,0)--(0,0); \draw[-,gray] (0,0)--(0,-0.23);

\draw (0,0) node{$_{\bullet}$};
\draw (0,3) node{$_{-}$}; \draw (0,3) node[left]{$_{1}$};
\draw (1.5,0) node{$_{|}$}; \draw (1.5,-0.1) node[below]{$_{\frac{1}{k}}$};
\draw (3,0) node{$_{|}$}; \draw (3,-0.1) node[below]{$_{1}$};
\draw[ultra thick] (1.5,3)--(2.92,3); \draw (1.5,3) node{$_{\bullet}$}; \draw (3,3) node{$_{\circ}$};
\draw[ultra thick] (0,0)--(1.5,0); 
\draw (0,-0.1) node[below]{$_{0}$};
\draw [dotted] (1.5,0)--(1.5,3);
\end{tikzpicture}
\begin{center}
$f_{k}(x)=\chi_{[\frac{1}{k},1)}(x)$
\end{center}
\end{center}    
\end{minipage} \hfill \begin{minipage}[c]{0.5\textwidth}
\begin{center}
\begin{tikzpicture}[xscale=0.75,yscale=0.75]
	\draw[->,gray] (0,0) -- (4,0); \draw [->,gray] (0,0) -- (0,4); 
	\draw[-,gray] (-0.2,0)--(0,0); \draw[-,gray] (0,0)--(0,-0.23);

\draw (0,0) node{$_{\bullet}$};
\draw (0,3) node{$_{-}$}; \draw (0,3) node[left]{$_{1}$};
\draw (1,0) node{$_{|}$}; \draw (1,-0.1) node[below]{$_{\frac{1}{k+1}}$};
\draw (3,0) node{$_{|}$}; \draw (3,-0.1) node[below]{$_{1}$};
\draw[ultra thick] (1,3)--(2.92,3); \draw (1,3) node{$_{\bullet}$}; \draw (3,3) node{$_{\circ}$};
\draw[ultra thick] (0,0)--(1,0); 
\draw (0,-0.1) node[below]{$_{0}$};
\draw [dotted] (1,0)--(1,3);
\end{tikzpicture}
\begin{center}
$f_{k+1}(x)=\chi_{[\frac{1}{k+1},1)}(x)$
\end{center}
\end{center}   
\end{minipage}
\end{figure}

Then $(f_k)$ is a sequence of nonnegative $\mathsf{S}$-measurable functions that converges pointwise on $[0,1)$ to the function
$$
f(x)=\left\{ 
\begin{array}{lcl}
1 & &\mbox{if}\,\,\, x \in (0,1),\\
0 & &\mbox{if}\,\,\, x \not \in (0,1).
\end{array}
\right.
$$

Observe that, for each $k \in \mathbb{N}$,
$$
\int_{[0,1)} f_k\,d\lambda= \int_{[0,1)} \chi_{[\frac{1}{k},1)}\,d\lambda= 1-\frac{1}{k}
$$
and, consequently,
$$
\lim_{k \to \infty} \int_{[0,1)} f_k\,d\lambda=\lim_{k \to\infty}\left(1-\frac{1}{k} \right) = 1.
$$

On the other hand,
$$
\int_{[0,1)} \lim_{k \to \infty}\,f_k\,d\lambda=\int_{[0,1)} f\,d\lambda = \int_{[0,1)} \chi_{(0,1)}\,d\lambda = 1.
$$
\end{example}

The previous example naturally suggests the question of whether it is possible, in general, to interchange the limit with the integral for sequences of nonnegative measurable functions. However, this property is not valid without additional hypotheses. Indeed, even when a sequence of nonnegative measurable functions converges pointwise to a measurable function, it does not necessarily follow that the integral of the limit coincides with the limit of the integrals. Under what conditions can such an interchange be guaranteed? It was the mathematician Beppo Levi\footnote{Beppo Levi (1875--1961) was an Italian mathematician who later became an Argentine citizen. He made contributions to the study of curves on algebraic surfaces, the Lebesgue integral, and measure theory. He introduced the spaces of square-integrable functions whose derivatives are also square-integrable, now known as ``Beppo Levi spaces''. In his honor, that theorem is called the ``Beppo Levi theorem''. In set theory he anticipated what would later be known as the ``Zermelo postulate''.} who first established the hypotheses under which this conclusion is valid.

\begin{figure}[ht!]
\centering
\includegraphics[scale=0.3]{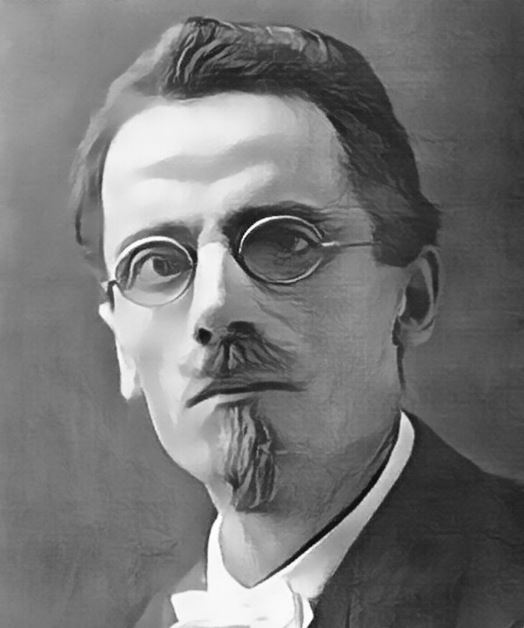} 
\begin{center}
B. Levi (1875-1961)
\end{center}
\end{figure}

\begin{theorem}[Monotone convergence theorem] \label{517} \index{theorem!of monotone convergence}
Let $(f_k)$ be a nondecreasing sequence of functions in $\overline{\mathbb{M}}^{+}(X,\mathsf{S})$ such that $f_k(x) \to f(x)$ for every $x \in X$. Then,
$$
\int_{A} f\,d\mu = \lim_{k\,\to\,\infty} \int_{A} f_k\,d\mu \quad\text{for every }\,\,A \in \mathsf{S}.
$$
\end{theorem}

\begin{proof}
By hypothesis, the sets $A_\pm(f)\in\mathsf{S}$, so Theorem \ref{331} ensures that $f \in \overline{\mathbb{M}}^{+}(X,\mathsf{S})$ and, consequently, $\int_{A} f\,d\mu$ is well defined in $\overline{\mathbb{R}}$.

Let $A \in \mathsf{S}$ be arbitrary but fixed and let $\varepsilon \in (0,1)$. For an arbitrary but fixed $s \in \underline{\mathcal{S}}(f)$, define the following sequence of sets in $\mathsf{S}$,
$$
A_{k}:=\left\{ x \in A\,:\, f_k(x) \geq (1-\varepsilon)\,s(x) \right\}.
$$

Note that since $f_k\leq f_{k+1}$, then $A_{k}\subset A_{k+1}$ for every $k \in \mathbb{N}$. Hence $(A_k)$ is a nondecreasing sequence in $\mathsf{S}$ such that $\lim_{k \to\infty} A_{k}=\bigcup_{k=1}^{\infty}A_k=A$. Indeed, if $x \in A$, since $f_{k}(x) \to f(x)$, there exists $k=k(\varepsilon,x) \in \mathbb{N}$ such that $f_{k}(x)\geq (1-\varepsilon)f(x) \geq (1-\varepsilon)s(x)$.

From the definition of the sets $A_{k}$ and Proposition \ref{513} we obtain
$$
\int_{A_{k}} (1-\varepsilon)s\,d\mu =(1-\varepsilon)\int_{A_{k}} s\,d\mu \leq \int_{A_k} f_{k}\,d\mu \leq \int_{A} f_{k}\,d\mu \qquad \forall k \in \mathbb{N}.
$$

Since $(1-\varepsilon)s$ is a nonnegative $\mathsf{S}$-simple function, Theorem \ref{58} together with Theorem \ref{414} ensure that
$$
\lim_{k \to \infty}\int_{A_k} (1-\varepsilon)s\,d\mu = \int_{A} (1-\varepsilon)s\,d\mu
$$
and, therefore,
$$
(1-\varepsilon)\int_{A} s\,d\mu \leq \lim_{k \to \infty} \int_{A}f_{k}\,d\mu.
$$

Taking the limit as $\varepsilon \to 0$ we conclude that
$$
\int_{A} s\,d\mu \leq \lim_{k \to \infty}\int_{A} f_k\,d\mu.
$$

Consequently,
$$
\int_{A} f\,d\mu = \sup\left\{ \int_{A} s\,d\mu\,:\, s \in \underline{\mathcal{S}}(f)\right\} \leq \lim_{k \to \infty}\int_{A} f_k\,d\mu.
$$

On the other hand, for every $k \in \mathbb{N}$, we have $0 \leq f_k \leq f_{k+1} \leq f$ on $X$ and, consequently, $0\leq (f_k \cdot \chi_{A}) \leq (f_{k+1} \cdot \chi_{A}) \leq (f \cdot \chi_{A})$ on $X$. Applying the monotonicity of the integral we obtain
$$
\int_{A} f_{k}\,d\mu \leq \int_{A} f_{k+1}\,d\mu \leq \int_{A} f\,d\mu
$$
and, therefore,
$$
\lim_{k \to \infty}\int_{A}f_{k}\,d\mu \leq \int_{A} f\,d\mu.
$$

This concludes the proof.
\end{proof}

In Example \ref{516} the sequence $(f_k)$ satisfies the hypotheses of Theorem \ref{517}, and thus it is confirmed that
$$
\int \lim_{k\to \infty} f_k\,d\mu = \lim_{k \to \infty} \int_X f_k\,d\mu = \lim_{k \to \infty} \left(1-\frac{1}{k} \right)=1.
$$

Let $f \in \overline{\mathbb{M}}^{+}(X,\mathsf{S})$. Theorem \ref{318} ensures that there exists a nondecreasing sequence $(s_k)$ of functions in $\mathbb{S}^{+}(X,\mathsf{S})$ that converges pointwise to $f$ on $X$. Applying Theorem \ref{517} we obtain
$$
\int_{X} f\,d\mu =\lim_{k\,\to\,\infty} \int_{X} s_k\,d\mu .
$$

It is possible to prove that the value of $\int_{X} f\,d\mu$ does not depend on the approximating sequence of simple functions [Exercise \ref{E525}]. This result provides a characterization of Definition \ref{59}.

\begin{figure}[ht!]
\centering   
\begin{tikzpicture}[xscale=2.85, yscale=2.85]


\draw[dashed, fill=red!10] (0.253,-0.01)--(0.253,0.2)--(1.52,0.2)--(1.52,-0.01);

\draw[dashed, fill=red!10] (0.354,0.2)--(0.354,0.4)--(1.4,0.4)--(1.4,0.2);

\draw[dashed, fill=red!10] (0.454,0.4)--(0.454,0.6)--(1.27,0.6)--(1.27,0.4);

\draw[dashed, fill=red!10] (0.554,0.6)--(0.554,0.8)--(1.17,0.8)--(1.17,0.6);

\draw[dashed, fill=red!10] (0.654,0.8)--(0.654,1)--(1.07,1)--(1.07,0.8);


\draw[-,red,thick] (-0.25,0)--(0.253,-0.01);
\draw[-,red,thick] (0.253,0.2)--(0.354,0.2);
\draw[-,red,thick] (0.354,0.4)--(0.454,0.4);
\draw[-,red,thick] (0.454,0.6)--(0.554,0.6);
\draw[-,red,thick] (0.554,0.8)--(0.654,0.8);

\draw[-,red,thick] (0.654,1)--(1.07,1);
\draw[-,red,thick] (1.07,0.8)--(1.17,0.8);
\draw[-,red,thick] (1.17,0.6)--(1.27,0.6);
\draw[-,red,thick] (1.27,0.4)--(1.4,0.4);
\draw[-,red,thick] (1.4,0.2)--(1.52,0.2);
\draw[-,red,thick] (1.52,-0.01)--(2.25,-0.01);


\draw (0,-0.01) node{$_{|}$}; 
\draw (2,-0.01) node{$_{|}$}; 

\draw[->,gray] (-0.25,-0.01)--(2.25,-0.01); 
\draw[->,gray] (0.875,-0.2)--(0.875,1.5); 

\draw[domain= 0:2,thick, smooth] plot(\x,{{3*\x^2*exp(-\x^3)}} );

\end{tikzpicture}
\begin{center}
$\int_{X}f\,d\mu=\lim_{k \to \infty}\int_{X}s_{k}\,d\mu$
\end{center}
\end{figure}

Finding a nondecreasing sequence of nonnegative simple functions that converges pointwise to a given nonnegative measurable function $f$ is, in some cases, an efficient procedure. This will allow us to compute integrals in a less tedious way than by directly applying Definition \ref{59}.

For example, consider the measure space $(\mathbb{N},\mathcal{P}(\mathbb{N}),\mu^{\sharp})$ and any nonnegative function $f:\mathbb{N} \to \mathbb{R}$. Defining, for each $n \in \mathbb{N}$, the function $s_n:\mathbb{N} \to \mathbb{R}$ by
$$
s_n:=\sum_{k=1}^{n} f(k)\chi_{\{k\}} + 0\chi_{\mathbb{N}\smallsetminus \{1,\ldots,n\} },
$$ 
it is clear that $s_n \in \mathbb{S}^{+}(\mathbb{N},\mathcal{P}(\mathbb{N}))$. Moreover, $(s_n)$ is nondecreasing and $s_n \to f$ pointwise on $\mathbb{N}$. Indeed, let $\varepsilon >0$ and $k \in \mathbb{N}$ be arbitrary. Taking $n_0:=k \in \mathbb{N}$, we have $s_n(k)=f(k)$ for every $n \geq n_0$ and, therefore, $|s_n(k)-f(k)|< \varepsilon$ for every $n \geq n_0$. Applying Theorem \ref{517} and Example \ref{55} we conclude that
$$
\int_{\mathbb{N}} f\,d\mu^{\sharp} = \lim_{n \to \infty}\int_{\mathbb{N}} s_n\,d\mu^{\sharp} = \lim_{n \to \infty} \sum_{k=1}^{n} f(k) = \sum_{k=1}^{\infty} f(k).
$$

We now prove some important properties of the integral of nonnegative measurable functions that follow immediately from Theorem \ref{517}.

\begin{corollary} \label{518}
Let $f,g \in \overline{\mathbb{M}}^{+}(X,\mathsf{S})$ and let $c \geq 0$. Then,
\begin{itemize}
\item[(a)] $\int_{A} cf\,d\mu = c\int_{A} f\,d\mu$ for every $A \in \mathsf{S}$.
\item[(b)] $\int_{A} (f+g)\,d\mu = \int_{A} f\,d\mu + \int_{A} g\,d\mu$ for every $A \in \mathsf{S}$.
\item[(c)] $\int_{A\cup B} f\,d\mu = \int_{A}f\,d\mu + \int_{B}f\,d\mu$ for every $A,B \in \mathsf{S}$ with $A \cap B =\varnothing$.
\end{itemize}
\end{corollary}

\begin{proof}
Let $A \in \mathsf{S}$. Since $f$ and $g$ are nonnegative measurable functions, Lemma \ref{318} ensures that there exist nondecreasing sequences $(s_k)$ and $(t_k)$ of elements of $\mathbb{S}^{+}(X,\mathsf{S})$ such that $s_k \to f$ and $t_k \to g$ on $X$. Then, $(\widetilde{s}_k):=(s_k \cdot \chi_{A})$ and $(\widetilde{t}_k):=(t_k \cdot \chi_{A})$ are nondecreasing sequences of elements of $\mathbb{S}^{+}(X,\mathsf{S})$ such that $\widetilde{s}_k \to f\cdot \chi_{A}$ and $\widetilde{t}_k \to g \cdot \chi_{A}$ pointwise on $X$.

\textit{(a):} Let $c \geq 0$. Proposition \ref{316} ensures that $c \widetilde{s}_k \in \mathbb{S}^{+}(X,S)$ for every $k \in \mathbb{N}$. Moreover, $c\widetilde{s}_k \to c(f\cdot \chi_{A})$ pointwise on $X$. Applying Proposition \ref{52} and Theorem \ref{517} we obtain
$$
\int_{A} c\,f\,d\mu = \lim_{k\to \infty}\int_{X} c\widetilde{s}_k\,d\mu = c\lim_{k \to \infty} \int_{X} \widetilde{s}_k  =c\int_{A} f\,d\mu
$$

\textit{(b):} Proposition \ref{52} and Theorem \ref{517} ensure that
$$
\int_{A} (f+g)\,d\mu = \lim_{k \to \infty} \int_{X} (\widetilde{s}_k  + \widetilde{t}_k )\,d\mu = \lim_{k \to \infty}\int_{X} \widetilde{s}_k \,d\mu + \lim_{k \to \infty} \int_{X} \widetilde{t}_k \,d\mu =\int_{A} f\,d\mu + \int_{B} g\,d\mu.
$$

\textit{(c):} Let $B\in \mathsf{S}$ be such that $A \cap B = \varnothing$. Then, $\chi_{A \cup B} = \chi_{A} + \chi_{B}$ [Exercise \ref{E33}] and, therefore, $f \cdot \chi_{A \cup B} = f\chi_{A} + f \chi_{B}$. Applying part \textit{(b)} we obtain the result.
\end{proof}

The analogue of Theorem \ref{517} for nonincreasing sequences of nonnegative measurable functions $(f_{k})$ is not true in general [Exercise \ref{E527}]. However, if $\int_{X} f_{1} \,d\mu < +\infty$ then the result is also valid.

\begin{corollary} \label{519}
Let $(f_k)$ be a sequence of functions in $\overline{\mathbb{M}}^{+}(X,\mathsf{S})$ that converges pointwise on $X$. If $f_{k+1} \leq f_{k} \leq \ldots \leq f_1$ for every $k \in \mathbb{N}$ and $\int_{A} f_1\,d\mu <+\infty$ for every $A \in \mathsf{S}$, then
$$
\int_{A}\, \lim_{k \to \infty} f_k\,d\mu =\lim_{k \to \infty} \int_{A} f_k\,d\mu \quad \mbox{for every}\,\,A \in \mathsf{S}.
$$
\end{corollary}

\begin{proof}
The result follows by applying Theorem \ref{517} to the sequence of functions $(f_1-f_k)$.
\end{proof}

The following result, due to Beppo Levi, generalizes the idea of Theorem \ref{517} to series of functions in $\overline{\mathbb{M}}^{+}(X,\mathsf{S})$.

\begin{theorem}[Beppo Levi] \label{520} \index{theorem!Beppo Levi theorem}
Let $(f_k)$ be a sequence of functions in $\overline{\mathbb{M}}^{+}(X,\mathsf{S})$ and define $f(x):=\sum_{k=1}^{\infty} f_k(x)$ for every $x \in X$. Then $f \in \overline{\mathbb{M}}^{+}(X,\mathsf{S})$ and, for every $A \in \mathsf{S}$, one has
$$
\int_{A} f\,d\mu = \sum_{k=1}^{\infty} \int_{A} f_k\,d\mu.
$$
\end{theorem}

\begin{proof}
Denote by
$$
g_{j}:=\sum_{k=1}^{j} f_k \quad \text{ and }\quad g:=\sum_{k=1}^{\infty} f_{k}.
$$

The sequence $(g_{j})$ is a nondecreasing sequence of elements of $\overline{\mathbb{M}}^{+}(X,\mathsf{S})$ that converges pointwise to $g$ on $X$. Therefore, $g \in \overline{\mathbb{M}}^{+}(X,\mathsf{S})$.

Applying Theorem \ref{517} to the sequence of partial sums $(g_{j})$ we obtain
$$
\int_X g\,d\mu = \lim_{j \to \infty} \int_{X} g_j\,d\mu =\lim_{j \to \infty} \int_X \sum_{k=1}^{j} f_k\,d\mu =\lim_{j \to \infty} \sum_{k=1}^{j} \int_X f_k\,d\mu = \sum_{k=1}^{\infty} \int_X f_k\,d\mu.
$$

Now, for each $A \in \mathsf{S}$, we apply the previous argument to the sequence of measurable functions $(f_k \cdot \chi_{A})$ and to the sequence of partial sums $(g_{j} \cdot \chi_{A})$, which converges pointwise to the function $(g \cdot \chi_{A})$ on $X$. This concludes the proof.
\end{proof}

The previous theorem has an application to double series of nonnegative extended real numbers by taking $X=\mathbb{N}$, $\mathsf{S}=\mathcal{P}(\mathbb{N})$, and $\mu=\mu^{\sharp}$ to be the counting measure.

\begin{corollary} \label{521}
Let $X=\mathbb{N}$, $\mathsf{S}=\mathcal{P}(\mathbb{N})$, and let $\mu^{\sharp}$ be the counting measure. For every function $f:\mathbb{N} \times \mathbb{N} \to [0,+\infty]$ given by $f(i,j)=x_{ij}$, one has
$$
\sum_{i=1}^{\infty} \left(\sum_{j=1}^{\infty} x_{ij} \right) =\sum_{j=1}^{\infty} \left(\sum_{i=1}^{\infty} x_{ij} \right).
$$
\end{corollary}

\begin{proof}
This is a direct consequence of Theorem 5.20 and Example \ref{515}.
\end{proof}

The following result establishes that every nonnegative measurable function $f:X \to \overline{\mathbb{R}}$ induces a measure on $(X,\mathsf{S})$.

\begin{theorem} \label{522}
Let $f \in \overline{\mathbb{M}}^{+}(X,\mathsf{S})$ be fixed. The function $\mu_{f}:\mathsf{S} \to \overline{\mathbb{R}}$ given by $\mu_{f}(A):=\int_{A} f\,d\mu$ defines a measure on $(X,\mathsf{S})$.
\end{theorem}

The proof is left as an exercise [Exercise \ref{E515}].

Let $(X,\mathsf{S},\mu)$ be a measure space. We say that a measure $\nu$ on $(X,\mathsf{S})$ is defined by a \textbf{density} if there exists a function $f \in \overline{\mathbb{M}}^{+}(X,\mathsf{S})$ such that $\nu=\mu_{f}$ on $\mathsf{S}$. \index{density} \index{function!density}
 
In the following example we compute the Lebesgue integral on $X$ with respect to the measure given in Theorem \ref{522}.

\begin{example} \label{523}
Let $f \in \overline{\mathbb{M}}^{+}(X,\mathsf{S})$ be fixed. Then, for every $g \in \overline{\mathbb{M}}^{+}(X,\mathsf{S})$, one has
$$
\int_{A} g\,d\mu_{f} = \int_{A} fg\,d\mu\quad \text{ for every }\,\,A \in \mathsf{S}.
$$
\end{example}

\begin{proof}
We prove only the case when $A=X$, since the general case follows easily from it. We proceed by cases.

{\scshape Case 1:}\quad $g \in \mathbb{S}^{+}(X,\mathsf{S})$.

Let $g=\sum_{j=1}^{N} \alpha_{j}\,\chi_{A_j}$ be its canonical representation. By the linearity of the integral (see Corollary \ref{519}) we have
$$
\int_X g\,d\mu_{f} = \sum_{j=1}^{N} \alpha_{j}\,\mu_{f}(A_j) = \sum_{j=1}^{N}  
  \alpha_{j} \int_X (f \cdot \chi_{A_j})\,d\mu = \int_X \sum_{j=1}^{N} (f \cdot \alpha_{j}\,\chi_{A_{j}})\,d\mu =\int_X fg\,d\mu.
$$

{\scshape Case 2:}\quad $g \in \overline{\mathbb{M}}^{+}(X,\mathsf{S})$.

Theorem \ref{318} ensures that there exists a nondecreasing sequence $(s_k)$ of functions in $\mathbb{S}^{+}(X,\mathsf{S})$ such that $s_k(x)\,\to\,g(x)$ for every $x \in X$. From the previous case and Theorem \ref{517} it follows that
$$
\int_X g\,d\mu_{f} = \lim_{k\,\to\,\infty} \int_X s_k\,d\mu_{f} = \lim_{k\,\to\,\infty} \int_X f\,s_k\,d\mu =\int_X fg\,d\mu,
$$
as stated.
\end{proof}

One of the disadvantages of Theorem \ref{517} is that it only works when the sequence of measurable functions $(f_k)$ is nondecreasing and convergent. It is then natural to ask what happens when we have an arbitrary sequence of nonnegative measurable functions. Let us look at the following example.

\begin{example} \label{524}
Let $X=\mathbb{N}$, $\mathsf{S}=\mathcal{P}(\mathbb{N})$, and let $\mu^{\sharp}$ be the counting measure. Define $f_k:\mathbb{N} \to \mathbb{R}$ by
$$
f_k:=\left\{ 
\begin{array}{lcl}
\chi_{P} & &\mbox{if}\,\,\,\,k\,\,\mbox{is even},\\
\chi_{X\smallsetminus P}& &\mbox{if}\,\,\,\,k\,\,\mbox{is odd},\\
\end{array}
\right.
$$
where $P:=\left\{2,4,6,\ldots \right\}$ is the set of even numbers in $\mathbb{N}$.

The sequence $(f_k)$ does not converge pointwise on $\mathbb{N}$ since
$$
\liminf_{k \to \infty} f_k = 0 \quad\text{ and }\quad \limsup_{k \to \infty}f_k =1.
$$

Consequently,
$$
\int_{\{1,2,3\}} \liminf_{k \to \infty} f_k \,d\mu^{\sharp} = 0\quad \text{and}\quad \int_{\{1,2,3\}} \limsup_{k \to \infty} f_k \,d\mu^{\sharp}=3.
$$

On the other hand, $f_k \cdot \chi_{\{1,2,3 \}}$ is given by
$$
(f_k \cdot \chi_{\{1,2,3 \}})=\left\{ 
\begin{array}{lcl}
\chi_{\{ 2\}} & &\mbox{if}\,\,\,\,k\,\,\mbox{is even},\\
\chi_{\{1,3\}}& &\mbox{if}\,\,\,\,k\,\,\mbox{is odd}.\\
\end{array}
\right.
$$

Hence,
$$
\liminf_{k \to \infty} \int_{\{1,2,3\}} f_k\,d\mu^{\sharp} = 1\quad \text{and}\quad \limsup_{k \to \infty} \int_{\{1,2,3\}} f_k\,d\mu^{\sharp} = 2.
$$
\end{example}

The following result, due to Pierre Fatou\footnote{Pierre Joseph Louis Fatou (1878-1929) was a French mathematician and astronomer who worked in the field of complex dynamics. He entered the École Normale Supérieure in Paris in 1898 to study mathematics and graduated in 1901. After graduating, he obtained a position as an astronomer at the Paris Observatory.}, establishes for sequences of functions in $\mathbb{M}^{+}(X,\mathsf{S})$, not necessarily convergent, a relation between the integral and the lower limit. This principle is, in general, weaker than Theorem \ref{517}, but it is especially useful in the study of arbitrary sequences. According to \cite{Grabinsky}, this result was proved by Fatou in his 1906 doctoral thesis, in the context of the study of trigonometric series.

\begin{figure}[ht!]
\centering
\includegraphics[scale=0.275]{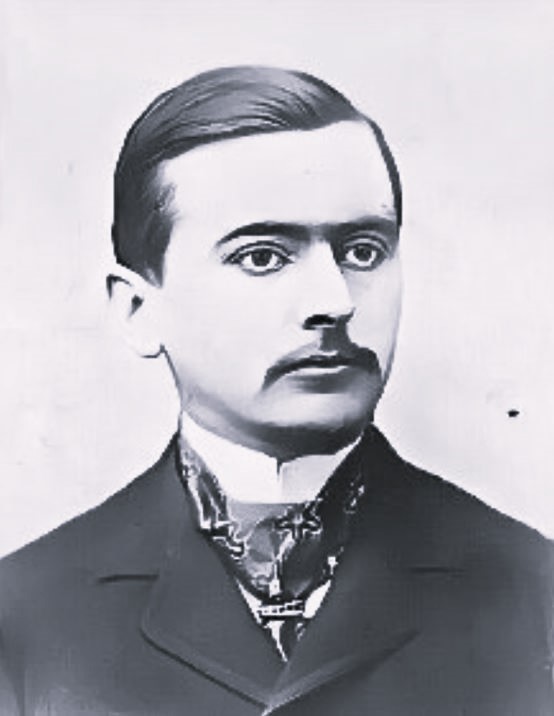} 
\begin{center}
P. Fatou (1878-1929)
\end{center}
\end{figure}

\begin{theorem}[Fatou's Lemma] \label{525} \index{lemma!Fatou's lemma}
Let $(f_k)$ be a sequence of functions in $\overline{\mathbb{M}}^{+}(X,\mathsf{S})$. Then,
$$
\int_{A} \liminf_{k \to \infty} f_k\,d\mu \leq \liminf_{k \to \infty} \int_{A} f_k\,d\mu \quad \mbox{for every}\,\,A \in \mathsf{S}.
$$
\end{theorem}

\begin{proof}
Let $A \in \mathsf{S}$ be arbitrary. For $k \in \mathbb{N}$, define $g_k:=\inf_{j \geq k} (f_j)$. Then $g_k$ is a nonnegative $\mathsf{S}$-measurable function such that $g_k \leq f_j$ for every $k \leq j$. By the monotonicity of the integral we have
$$
\int_{A} g_k\,d\mu \leq \int_{A} f_j\,d\mu \quad \text{ for every }\,\,\, k \leq j
$$
and, therefore, $\int_{A} g_k\,d\mu \leq \displaystyle\liminf_{k \to \infty}\int_{A} f_k\,d\mu$ for every $k \in \mathbb{N}$.

On the other hand, the sequence $(g_k)$ satisfies $0\leq g_k \leq g_{k+1}$ for every $k \in \mathbb{N}$ and
$$
\displaystyle\lim_{k \to\infty} g_k =\liminf_{k \to \infty} f_k.
$$
Applying Theorem \ref{517} and using the previous inequality, we obtain
$$
\int_{A} \liminf_{k \to \infty} f_k\,d\mu = \int_{A} \lim_{k \to\infty} g_k\,d\mu = \lim_{k \to\infty} \int_{A} g_k\,d\mu \leq \liminf_{k \to \infty} \int_{A} f_k\,d\mu,
$$
which is the desired identity.
\end{proof}

The sequence $(f_k)$ from Example \ref{524} satisfies the hypothesis of Theorem \ref{525}, and thus the inequality
$$
\int_{\{1,2,3\}} \liminf_{k \to \infty} f_k\,d\mu^{\sharp} = 0 < \liminf_{k \to \infty} \int_{\{1,2,3 \}} f_k\,d\mu^{\sharp} = 1
$$
is confirmed.

The following example shows that, in general, equality does not hold in Fatou's lemma even when the sequence $(f_k)$ converges pointwise.
 
\begin{example} \label{526} 
Let $X=\mathbb{N}$, $\mathsf{S}=\mathcal{P}(\mathbb{N})$, and let $\mu=\mu^{\sharp}$ be the counting measure. Let $f_k:\mathbb{N} \to \mathbb{R}$ be the function given by $f_k(n):=\frac{1}{k}\,\chi_{\{k,\ldots,2k-1 \}}(n).$

\begin{figure}[ht!]
\centering

\begin{minipage}{0.45\textwidth}
\centering
\begin{tikzpicture}[xscale=0.65,yscale=0.8]
	\draw[->, gray] (0,0) -- (8.5,0); \draw [->,gray] (0,0) -- (0,4); 
	\draw[-,gray] (-0.2,0)--(0,0); \draw[-,gray] (0,0)--(0,-0.23);
\draw (8.5,0) node[right] {$\mathbb{N}$}; \draw (0,4) node[right] {$\mathbb{R}$}; 

\draw (0,3) node{$_{-}$}; \draw (0,3) node[left]{$_{\frac{1}{k}}$};

\draw (1,0) node{$_{|}$}; \draw (1,-0.1) node[below]{$_{1}$}; \draw (2,-0.1) node[below]{$_{\ldots}$};

\draw (3,0) node{$_{|}$}; \draw (3,-0.1) node[below]{$_{k}$};

\draw (4,0) node{$_{|}$}; \draw (4,-0.1) node[below]{$_{k+1}$};

\draw (7,0) node{$_{|}$}; \draw (5,-0.1) node[below]{$_{\ldots}$};

\draw (6,0) node{$_{|}$}; \draw (6,-0.1) node[below]{$_{2k-1}$};  \draw (7,-0.1) node[below]{$_{2k}$}; 

\draw (8,-0.1) node[below]{$_{\ldots}$};

\draw (3,3) node{$_{\bullet}$};
\draw (4,3) node{$_{\bullet}$};
\draw (6,3) node{$_{\bullet}$};

\draw (1,0) node{$_{\bullet}$};
\draw (7,0) node{$_{\bullet}$};

\draw[dotted] (3,3)--(3,0);
\draw[dotted] (4,3)--(4,0);
\draw[dotted] (6,3)--(6,0);
\draw (0,-0.1) node[below]{$_{0}$};
\end{tikzpicture}

\vspace{0.2cm}
$f_{k}(n)=\frac{1}{k}\chi_{\{k,\ldots,2k-1\}}(n)$
\end{minipage}
\hfill
\begin{minipage}{0.45\textwidth}
\centering
\begin{tikzpicture}[xscale=0.65,yscale=0.8]
	\draw[->, gray] (0,0) -- (8.5,0); \draw [->,gray] (0,0) -- (0,4); 
	\draw[-,gray] (-0.2,0)--(0,0); \draw[-,gray] (0,0)--(0,-0.23);
\draw (8.5,0) node[right] {$\mathbb{N}$}; \draw (0,4) node[right] {$\mathbb{R}$}; 

\draw (0,2) node{$_{-}$}; \draw (0,2) node[left]{$_{\frac{1}{k+1}}$};

\draw (1,0) node{$_{|}$}; \draw (1,-0.1) node[below]{$_{1}$}; \draw (2,-0.1) node[below]{$_{\ldots}$};

\draw (3,0) node{$_{|}$}; \draw (3,-0.1) node[below]{$_{k}$};

\draw (4,0) node{$_{|}$}; \draw (4,-0.1) node[below]{$_{k+1}$};

\draw (7,0) node{$_{|}$}; \draw (5,-0.1) node[below]{$_{\ldots}$};

\draw (6,0) node{$_{|}$}; \draw (6,-0.1) node[below]{$_{2k-1}$};  \draw (7,-0.1) node[below]{$_{2k}$}; 

\draw (8,-0.1) node[below]{$_{\ldots}$};

\draw (3,0) node{$_{\bullet}$};
\draw (4,2) node{$_{\bullet}$};
\draw (6,2) node{$_{\bullet}$};

\draw (1,0) node{$_{\bullet}$};
\draw (7,2) node{$_{\bullet}$};

\draw[dotted] (4,2)--(4,0);
\draw[dotted] (6,2)--(6,0);
\draw[dotted] (7,2)--(7,0);
\draw (0,-0.1) node[below]{$_{0}$};
\end{tikzpicture}
\begin{center}
$f_{k+1}(n)=\frac{1}{k+1}\chi_{\{k+1,\ldots,2k\}}(n)$
\end{center}
\end{minipage}
\end{figure}

Then $f_k \geq 0$ and
$$
\int_{\mathbb{N}} f_k\,d\mu^{\sharp} = \frac{1}{k}\,\mu^{\sharp}(\{k,\ldots,2k-1\}) = 1 \quad\forall k \in \mathbb{N}.
$$

Since $\displaystyle\lim_{k \to \infty} f_k (n)=0$ for every $n \in \mathbb{N}$, we have
$$
\int_{\mathbb{N}} \lim_{k \to \infty} f_k\,d\mu^{\sharp} = 0 < 1 = \lim_{k \to \infty}\int_{\mathbb{N}} f_k\,d\mu^{\sharp}.
$$
\end{example}

Fatou's lemma can be extended in the following result.

\begin{theorem} \label{527}
Let $(f_k)$ be a sequence of elements of $\overline{\mathbb{M}}(X,\mathsf{S})$.
\begin{itemize}
\item[(a)] If there exists $g \in \mathbb{M}^{+}(X,\mathsf{S})$ with $\int_X g\,d\mu < +\infty$ such that $f_k \geq -g$ $\mu$ a.e. on $A \in \mathsf{S}$, then
$$
\int_{A} \liminf_{k \to \infty} f_k\,d\mu \leq \liminf_{k \to \infty} \int_{A} f_k\,d\mu.
$$
\item[(b)] If there exists $g \in \mathbb{M}^{+}(X,\mathsf{S})$ with $\int_X g\,d\mu < +\infty$ such that $f_k \leq -g$ $\mu$ a.e. on $A \in \mathsf{S}$, then
$$
 \limsup_{k \to \infty} \int_{A} f_k \,d\mu \leq \int_{A} \limsup_{k \to \infty} f_k\,d\mu .
$$
\end{itemize}
\end{theorem}

The proof is left as an exercise [Exercise \ref{E532}].

An interesting fact is that, assuming Fatou's lemma to be valid, one can deduce the monotone convergence theorem [Exercise \ref{E531}]. To do this, it is necessary to have a proof of Fatou's lemma that does not depend on the monotone convergence theorem, that is, an alternative and independent proof.

\vspace{0.5cm}

\begin{proof}[Proof of Theorem \ref{525}]
Let $A \in \mathsf{S}$ be arbitrary. Denote by
$$
f_{\ast}:=\displaystyle\liminf_{k \to \infty} f_k. 
$$

Let $s_{\ast} \in \underline{\mathcal{S}}(f_{\ast})$ and $\varepsilon \in (0,1)$ be arbitrary. For each $k \in \mathbb{N}$ define the set
$$
A_{k}:=\left\{x \in A\,:\,  f_j(x) \geq (1-\varepsilon)\,s_{\ast}(x)\,\,\,\mbox{for every}\,\,j \geq k \right\},
$$
which is clearly an element of $\mathsf{S}$.

Then $(A_k)$ is an increasing sequence of elements of $\mathsf{S}$ such that $A=\bigcup_{k=1}^{\infty}A_k$ since $(1-\varepsilon)s_{\ast} \leq f_{\ast}$. Therefore,
$$
\int_{A_k} (1-\varepsilon)s_{\ast}\,d\mu =(1-\varepsilon)\int_{A_k} s_{\ast}\,d\mu \leq \int_{A_k} f_k\,d\mu \leq \int_{A} f_k\,d\mu.
$$

Since $s_{\ast}$ is a nonnegative $\mathsf{S}$-simple function, applying Theorem \ref{58} and Corollary \ref{415} to the previous inequality, we obtain
$$
(1-\varepsilon)\int_{A} s_{\ast}\,d\mu = (1-\varepsilon)\liminf_{k \to \infty} \int_{A_k} s_{\ast}\,d\mu \leq \liminf_{k \to \infty} \int_{A} f_k\,d\mu.
$$

Since $\varepsilon \in (0,1)$ is arbitrary, we conclude that
$$
\int_{A} s_{\ast}\,d\mu \leq \liminf_{k \to \infty} \int_{A} f_k\,d\mu.
$$

Taking the supremum over all $s_{\ast} \in \underline{\mathcal{S}}(f_{\ast})$, it follows that
$$
\int_{A} f_{\ast}\,d\mu=\int_{A} \liminf_{k \to \infty} f_k\,d\mu \leq \liminf_{k \to \infty} \int_{A} f_k\,d\mu,
$$
which concludes the proof.
\end{proof}

We conclude this section with the following result.

\begin{theorem} \label{528}
Let $f\in \overline{\mathbb{M}}^{+}(X,\mathsf{S})$. Then,
$$
\int_{X} f\,d\mu =0 \quad \text{if and only if}\quad \mu\left(\left\{ x\in X\,:\, f(x)>0 \right\}  \right)=0.
$$
\end{theorem}

\begin{proof}
$\Rightarrow):$ Define the following sequence of sets
$$
A_k:=\left\{ x \in X\,:\, f(x) > \frac{1}{k} \right\}.
$$

Since $f$ is $\mathsf{S}$-measurable, $(A_k)$ is an increasing sequence of elements of $\mathsf{S}$ such that
$$
\bigcup_{k=1}^{\infty}\,A_k=\left\{ x\in X\,:\, f(x)>0 \right\}.
$$
By Proposition \ref{426}, it is enough to prove that $\mu(A_k) =0$ for every $k \in \mathbb{N}$.

From the definition of the sets $A_{k}$ it follows that
$$
0 \leq \left(\frac{1}{k} \cdot \chi_{A_k}\right) \leq ( f \cdot \chi_{A_k})
$$
on $X$ for every $k \in \mathbb{N}$, and applying the monotonicity of the integral we obtain
$$
0 \leq  \frac{1}{k}\,\mu(A_k)=\int_{A_k} \frac{1}{k}\,d\mu\leq \int_{A_{k}}f\,d\mu \leq \int_{X} f\,d\mu =0\quad \forall\,k \in \mathbb{N}.
$$

Consequently, $\mu(A_k)=0$ for every $k \in \mathbb{N}$.

$\Leftarrow):$ Define the sets
$$
A:=\left\{ x\in X\,:\, f(x)>0 \right\}
\quad \text{and} \quad
B:=\left\{ x \in X\,:\, f(x)=0\right\}.
$$

Clearly $A$ and $B$ belong to $\mathsf{S}$, are disjoint, and satisfy $X= A\cup B$. Using Corollary \ref{518} together with Proposition \ref{513}, we conclude that
$$
\int_{X} f\,d\mu = \int_{A} f\,d\mu + \int_{B} f\,d\mu =0.
$$

This proves the assertion.
\end{proof}

\section{Relation with the Riemann integral}

Let $a,b \in \mathbb{R}$ be such that $a < b$. Consider on $[a,b]$ the Borel $\sigma$-algebra $\mathcal{B}([a,b])$ and let $\lambda$ denote the Lebesgue measure on $[a,b]$.

Consider the following examples.

\begin{example} \label{529}
If $f:[0,1] \to \mathbb{R}$ is given by $f(x)=x$, then
$$
\int_{[0,1]} f\,d\lambda = \frac{1}{2}.
$$
\end{example}

\begin{proof}
 Since $f$ is continuous on $[0,1]$, it follows that $f$ is Borel-measurable and nonnegative on $[0,1]$.

For each $k \in \mathbb{N}$, consider the simple nonnegative function $s_k :  [0,1] \to \mathbb{R}$ given by
$$
s_{k}(x):=\sum_{j=0}^{k2^k} \frac{j}{2^k}\,\chi_{A_{k}(j)}(x),
$$
where
$$
\begin{aligned}
A_{k}(j)&:=\left\{ x \in [0,1]\,:\, \frac{j}{2^{k}} \leq f(x) < \frac{j+1}{2^{k}}  \right\},\quad j \in \{0, 1, 2, \ldots, k2^{k}-1\},\\
A_{k}(k2^{k})&:=\left\{ x \in  [0,1]\,:\, k \leq f(x)  \right\}.
\end{aligned}
$$

Lemma \ref{318} ensures that $(s_k)$ is an increasing sequence such that $s_k \to f$ pointwise on $ [0,1]$.

\begin{figure}[ht!]
\begin{minipage}[l]{0.3\textwidth}
\begin{center}
\begin{tikzpicture}[xscale=1.35, yscale=1.35]

\draw[->,gray] (0,0)--(3,0); 
\draw[->,gray] (0,0)--(0,3); 

\draw[-,thick] (0,0)--(2,2);

\draw (0,0) node[left]{$_{0/2}$}; 
\draw (0,1) node[left]{$_{1/2}$}; 
\draw (0,2) node[left]{$_{2/2}$};

\draw (0,0) node{$_{-}$}; 
\draw (0,1) node{$_{-}$}; 
\draw (0,2) node{$_{-}$};

\draw[dotted] (1,1)--(1,0);
\draw[dotted] (2,2)--(2,0);

\draw[-,ultra thick] (0,0)--(1,0);
\draw[-,ultra thick] (1,1)--(2,1);
\draw (2,2) node{$_{\bullet}$};
\draw (2,0) node[below]{$_{1}$};
\draw[dotted] (0,2)--(2,2);
\end{tikzpicture}
\begin{center}
$s_1(x)$
\end{center}

\end{center}
\end{minipage} \hfill 
 \begin{minipage}[c]{0.3\textwidth}
\begin{center}

\begin{tikzpicture}[xscale=1.35, yscale=1.35]

\draw[->,gray] (0,0)--(3,0); 
\draw[->,gray] (0,0)--(0,3); 

\draw[-,thick] (0,0)--(2,2);

\draw (0,0) node[left]{$_{0/4}$}; 
\draw (0,0.5) node[left]{$_{1/4}$}; 
\draw (0,1) node[left]{$_{2/4}$};
\draw (0,1.5) node[left]{$_{3/4}$};
\draw(0,2) node[left]{$_{4/4}$};
\draw (0,2.45) node[left]{$_{5/4}$};
\draw (0,2.8) node[left]{$_{\vdots}$};

\draw (0,0) node{$_{-}$}; 
\draw (0,0.5) node{$_{-}$}; 
\draw (0,1) node{$_{-}$}; 
\draw (0,1.5) node{$_{-}$}; 
\draw (0,2) node{$_{-}$};
\draw (0,2.45) node{$_{-}$};

\draw[-,ultra thick] (0,0)--(0.5,0);
\draw[-,ultra thick] (0.5,0.5)--(1,0.5);
\draw[-,ultra thick] (1,1)--(1.5,1);
\draw[-,ultra thick] (1.5,1.5)--(2,1.5);
\draw (2,2) node{$_{\bullet}$};

\draw[dotted] (0.5,0.5)--(0.5,0);
\draw[dotted] (1,1)--(1,0);
\draw[dotted] (1.5,1.5)--(1.5,0);
\draw[dotted] (2,2)--(2,0);
\draw[dotted] (0,2)--(2,2);

\draw (2,0) node[below]{$_{1}$};
\end{tikzpicture}
\begin{center}
$s_2(x)$
\end{center}

\end{center}
\end{minipage} \hfill 
 \begin{minipage}[r]{0.3\textwidth}
\begin{center}

\begin{tikzpicture}[xscale=1.35, yscale=1.35]

\draw[->,gray] (0,0)--(3,0); 
\draw[->,gray] (0,0)--(0,3); 

\draw[-,thick] (0,0)--(2,2);

\draw (0,0) node[left]{$_{0/8}$}; 
\draw (0,0.25) node[left]{$_{1/8}$}; 
\draw (0,0.5) node[left]{$_{2/8}$}; 
\draw (0,0.75) node[left]{$_{3/8}$}; 
\draw (0,1) node[left]{$_{4/8}$};
\draw (0,1.25) node[left]{$_{5/8}$}; 
\draw (0,1.5) node[left]{$_{6/8}$}; 
\draw (0,1.75) node[left]{$_{7/8}$}; 
\draw(0,2) node[left]{$_{8/8}$};
\draw (0,2.25) node[left]{$_{9/8}$};
\draw (0,2.75) node[left]{$_{\vdots}$};

\draw (0,0) node{$_{-}$}; 
\draw (0,0.25) node{$_{-}$}; 
\draw (0,0.5) node{$_{-}$}; 
\draw (0,0.75) node{$_{-}$}; 
\draw (0,1) node{$_{-}$};
\draw (0,1.25) node{$_{-}$}; 
\draw (0,1.5) node{$_{-}$}; 
\draw (0,1.75) node{$_{-}$};  
\draw (0,2) node{$_{-}$};
\draw (0,2.25) node{$_{-}$};

\draw[-,ultra thick] (0,0)--(0.25,0);
\draw[-,ultra thick] (0.25,0.25)--(0.5,0.25);
\draw[-,ultra thick] (0.5,0.5)--(0.75,0.5);
\draw[-,ultra thick] (0.75,0.75)--(1,0.75);
\draw[-,ultra thick] (1,1)--(1.25,1);
\draw[-,ultra thick] (1.25,1.25)--(1.5,1.25);
\draw[-,ultra thick] (1.5,1.5)--(1.75,1.5);
\draw[-,ultra thick] (1.75,1.75)--(2,1.75);
\draw (2,2) node{$_{\bullet}$};

\draw[dotted] (0.25,0.25)--(0.25,0);
\draw[dotted] (0.5,0.5)--(0.5,0);
\draw[dotted] (0.75,0.75)--(0.75,0);
\draw[dotted] (1,1)--(1,0);
\draw[dotted] (1.25,1.25)--(1.25,0);
\draw[dotted] (1.5,1.5)--(1.5,0);
\draw[dotted] (1.75,1.75)--(1.75,0);
\draw[dotted] (2,2)--(2,0);
\draw[dotted] (0,2)--(2,2);

\draw (2,0) node[below]{$_{1}$};
\end{tikzpicture}
\begin{center}
$s_3(x)$
\end{center}

\end{center}
\end{minipage}
\end{figure}

Observe that, for each $k \in \mathbb{N}$, $A_{k}(j)=\varnothing$ for every $j \in \{2^{k}+1,\ldots,k2^{k}-1,k2^{k} \}$, so that
$$
s_{k}(x):=\sum_{j=0}^{2^k} \frac{j}{2^k}\,\chi_{A_{k}(j)}(x),
$$
where
$$
\begin{aligned}
A_{k}(j)&:=\left\{ x \in  [0,1]\,:\, \frac{j}{2^{k}} \leq f(x) < \frac{j+1}{2^{k}}  \right\}=\left[\frac{j}{2^{k}},\frac{j+1}{2^{k}} \right)\quad j \in \{0, 1, 2, \ldots, 2^{k}-1\},\\
A_{k}(2^k)&:=\{1\}.\\
\end{aligned}
$$

Consequently,
$$
\begin{aligned}
\int_{[0,1]} s_k\,d\lambda &= \sum_{j=0}^{2^k} \frac{j}{2^k}\, \lambda(A_j) = \sum_{j=0}^{2^k-1} \frac{j}{2^k} \,\lambda\left( \left[ \frac{j}{2^k}, \frac{j+1}{2^k}\right) \right)  = \sum_{j=0}^{2^k-1} \frac{j}{4^k} = \frac{1}{4^k} \sum_{j=0}^{2^k-1} j \\
&=\frac{(2^k-1)2^{k-1}}{4^k} =\frac{1}{2}\left( 1 - \frac{1}{2^k}\right),\quad \mbox{for every}\,\,k \in \mathbb{N}.
\end{aligned}
$$

Applying the monotone convergence theorem, we obtain
$$
\int_{[0,1]} f\,d\lambda = \lim_{k \to \infty} \int_{[0,1]} s_k\,d\lambda = \lim_{k \to \infty} \frac{1}{2}\left( 1 - \frac{1}{2^k}\right)=\frac{1}{2}
$$
as stated.

\end{proof}

In general, it can be proved [Exercise \ref{E535}] that, for every $n \in \mathbb{N}$, the equality
$$
\int_{[0,1]} x^n\,d\lambda = \frac{1}{n+1}
$$
holds, since the function $x^n$ is continuous, Borel-measurable, and nonnegative on $[0,1]$.

Since $x^n$ is continuous on $[0,1]$, Theorem \ref{16} ensures that $x^n$ is Riemann integrable on $[0,1]$ and, in fact,
$$
R\int_{0}^{1} x^n\,dx = \frac{1}{n+1}.
$$

\begin{example} \label{530}
On $[0,1]$ consider $\widetilde{\mathbb{Q}}$, the set of all rational numbers in the interval $[0,1]$, and let $\{q_0,q_1,q_2,\ldots\}$ be an enumeration of it.

For any fixed $N \in \mathbb{N}$, the function given by
$$
f_N:=\chi_{\{q_0,q_1,\ldots,q_N\}}
$$
is Borel-measurable since $\{q_0,q_1,\ldots,q_N\} \in \mathcal{B}([0,1])$. Consequently,
$$
\int_{[0,1]} f_N\,d\lambda = \lambda(\{q_0,q_1,\ldots,q_N\}) = \sum_{j=0}^{N} \lambda(\{q_j\}) =0.
$$

Moreover, $f_N:[0,1]\to \mathbb{R}$ is Riemann integrable on $[0,1]$ and
$$
R\int_{0}^{1} f_N(x)\,dx=0.
$$
\end{example}

The previous examples naturally lead us to ask whether for every bounded, nonnegative, Borel-measurable function $f:[a,b] \to \mathbb{R}$ one has
$$
\int_{[a,b]} f\,d\lambda = R\int_{a}^{b} f(x)\,dx. 
$$

The answer is not true in general, as shown by the following example.

\begin{example} \label{531} \index{function!Dirichlet function}
On $[0,1]$ consider $\widetilde{\mathbb{Q}}$, the set of all rational numbers in the interval $[0,1]$. The function given by
$$
f:=\chi_{\widetilde{\mathbb{Q}}}
$$
is Borel-measurable since $\widetilde{\mathbb{Q}}\in \mathcal{B}([0,1])$ {\rm [Exercise \ref{E440}]}. Consequently,
$$
\int_{[0,1]} f\,d\lambda = \lambda(\widetilde{\mathbb{Q}}) =0.
$$

However, $f:[0,1]\to \mathbb{R}$ is not Riemann integrable on $[0,1]$ (see {\rm Example \ref{19}}).
\end{example}

We have exhibited a function that is not integrable in the Riemann sense, but whose Lebesgue integral exists and is equal to $0$. This example highlights the fact that the Lebesgue integral constitutes an extension and reformulation of the concept of the Riemann integral, since it encompasses a broader class of functions and thus resolves the first problem posed at the beginning of the text. This idea is fully formalized in the following theorem, due to Henri Lebesgue, which moreover provides a complete characterization of bounded functions $f:[a,b]\to\mathbb{R}$ that are Riemann integrable, in terms of the structure of their set of discontinuities. In what follows, we shall state a weaker version of this result for the case in which $f:[a,b]\to\mathbb{R}$ is a Borel-measurable function on $[a,b]$; the general version will be presented later as an exercise. 

\begin{theorem}[Henri Lebesgue] \label{532} \index{theorem! of Henri Lebesgue}
Let $f:[a,b] \to \mathbb{R}$ be a bounded nonnegative function. Then,
\begin{itemize}
\item[(a)] $f$ is Riemann-integrable if and only if $f$ is continuous $\lambda$ a.e.
\item[(b)] If $f$ is Borel-measurable and Riemann-integrable, then
$$
R\int_{a}^{b} f(x)\,dx = \int_{[a,b]} f\,d\lambda.
$$
\end{itemize}
\end{theorem}

To prove H. Lebesgue's theorem we require the following concepts and the following three lemmas.

\begin{definition} \label{533}
Let $f:[a,b] \to \mathbb{R}$ be a bounded function and let $\delta >0$. Let
$g_{\delta}(x):=\inf\{ f(t)\,:\, t\in(x-\delta,x+\delta) \cap [a,b] \}$ and
$h_{\delta}(x):=\sup\left\{f(t)\,:\, t\in(x-\delta,x+\delta) \cap [a,b] \right\}$.
We define the \textbf{lower envelope of $f$} and the \textbf{upper envelope of $f$} as follows: \index{lower!envelope} \index{upper!envelope}
$$
g(x):=\lim_{\delta\to0} g_{\delta}(x):=\sup_{\delta >0} g_{\delta}(x)\quad \mbox{and}\quad h(x):=\lim_{\delta\to 0} h_{\delta}(x):=\inf_{\delta >0}h_{\delta}(x)\quad \mbox{respectively}.
$$
\end{definition}

\begin{figure}[ht!]
\centering
\begin{tikzpicture}[yscale=1.655]
\draw[->, gray] (-0.5,0) -- (6,0);
\draw [->, gray] (0,-0.2) -- (0,2);
\draw[thick] plot[smooth] coordinates
{(0.5,0.5)(1,1)(3,0.51)(4,1.5)(5,1.5)};

\draw (6,0) node [right]{$_{[a,b]}$}; \draw (0,2) node[right]{$\mathbb{R}$};
\draw[dotted] (3.74,0)--(3.74,1.25);
\draw[dotted] (0.5,0)--(0.5,0.5);
\draw[dotted] (5,0)--(5,1.5);
\draw[dotted] (2.8,0)--(2.8,1.25);



\draw (2.8,0) node{$_{|}$}; \draw (2.8,-0.1) node[below]{$_{x-\delta}$};
\draw (3.75,0) node{$_{|}$}; \draw (3.75,-0.1) node[below]{$_{x+\delta}$};

\draw (0.5,0) node{$_{|}$}; \draw (0.5,-0.1) node[below]{$_{a}$};
\draw (5,0) node{$_{|}$}; \draw (5,-0.1) node[below]{$_{b}$};

\draw[thick,red] (2.8,0.5)--(3.74,0.5);
\draw[thick,blue] (2.8,1.25)--(3.74,1.25);
\draw (2.8,0.5) node{$_{\circ}$};
\draw (2.8,1.25) node{$_{\circ}$};

\draw (3.75,0.5) node{$_{\circ}$};
\draw (3.75,1.25) node{$_{\circ}$};

\draw (3.75,0.5) node[right]{$_{g_{\delta}(x)}$};
\draw (3.75,1.25) node[right]{$_{h_{\delta}(x)}$};

\draw (3.25,0) node{$_{|}$};\draw (3.25,-0.1) node[below]{$_{x}$};
\end{tikzpicture}
\begin{center}
$g_{\delta}(x)$ and $h_{\delta}(x)$
\end{center}
\end{figure}

\begin{lemma} \label{534}
Let $f:[a,b] \to \mathbb{R}$ be a bounded function. The lower envelope of $f$ is lower semicontinuous and the upper envelope of $f$ is upper semicontinuous.
\end{lemma}

\begin{proof}
Let $x_{0} \in [a,b]$. Since $f$ is bounded, then $g(x_0) \in \mathbb{R}$. For every arbitrary $\varepsilon >0$, we can find $\delta >0$ such that
$$
g(x_0) \leq g_{\delta}(x_0) + \varepsilon :=\inf\{ f(t)\,:\, t\in(x_{0}-\delta,x_{0}+\delta) \cap [a,b] \} + \varepsilon.
$$

Then,
$$
\begin{aligned}
g(x_0) &\leq \inf\{ f(t)\,:\, t\in(x_{0}-\delta,x_{0}+\delta) \cap [a,b] \} + \varepsilon \\
&\leq  \inf\{ f(t)\,:\, t\in\left(x-\frac{\delta}{2},x+\frac{\delta}{2}\right) \cap [a,b] \} + \varepsilon \\
&= g_{\frac{\delta}{2}}(x) + \varepsilon \\
&\leq g(x) + \varepsilon\quad\,\,\,\mbox{if}\,\,\,x \in \left(x_0-\frac{\delta}{2},x_0+\frac{\delta}{2} \right),
\end{aligned}
$$
since $\left(x-\frac{\delta}{2},x+\frac{\delta}{2}\right) \subset (x_{0}-\delta,x_{0}+\delta)$. Consequently, $g$ is lower semicontinuous.

The proof that $h$ is upper semicontinuous is analogous and is left as an exercise [Exercise \ref{E537}].
\end{proof}

The previous result, besides asserting that the lower envelope $g$ and the upper envelope $h$ of a bounded function $f$ are lower semicontinuous and upper semicontinuous, respectively, also establishes that they are Borel-measurable functions by virtue of Examples \ref{37} and \ref{38}.

\begin{lemma} \label{535}
Let $f:[a,b] \to \mathbb{R}$ be a bounded function, $g$ its lower envelope, and $h$ its upper envelope. Then,
\begin{itemize}
\item[(a)] $g \leq f \leq h$.
\item[(b)] $g(x_0)=f(x_0)$ if and only if $f$ is lower semicontinuous at $x_0 \in [a,b]$.
\item[(c)] $f(x_0)=h(x_0)$ if and only if $f$ is upper semicontinuous at $x_0 \in [a,b]$.
\item[(d)] $g(x_0)=f(x_0)=h(x_0)$ if and only if $f$ is continuous at $x_0 \in [a,b]$.
\end{itemize} 
\end{lemma}

\begin{proof}
\textit{(a):} It is an immediate consequence of the definitions of $g$ and $h$.

\textit{(b):}  Suppose that $g(x_0)=f(x_0)$ with $x_0 \in [a,b]$. Let $\varepsilon >0$ be arbitrary. We find $\delta >0$ such that 
$$
g(x_0)-\varepsilon \leq g_{\delta}(x_0):=\inf\{ f(x)\,:\, x\in(x_{0}-\delta,x_{0}+\delta) \cap [a,b] \} \leq f(x)
$$
for every $x \in (x_0-\delta,x_0+\delta) \cap [a,b]$. 

Therefore, $f(x_0)-\varepsilon \leq f(x)$ for every $x \in (x_0-\delta,x_0+\delta) \cap [a,b]$, that is, $f$ is lower semicontinuous at $x_0$.

Conversely, suppose that $f$ is lower semicontinuous at $x_0 \in [a,b]$. By part \textit{(a)}, it is enough to prove that $f(x_0) \leq g(x_0)$. Given $\varepsilon>0$, we find $\delta >0$ such that
$$
f(x_0)-\varepsilon \leq f(x)
$$
for every $x \in (x_0-\delta,x_0+\delta) \cap [a,b]$. Then,
$$
f(x_0)-\varepsilon \leq \inf\{ f(x)\,:\, x\in(x_{0}-\delta,x_{0}+\delta) \cap [a,b] \} :=g_{\delta}(x_0) \leq g(x_0)
$$
and, since $\varepsilon >0$ is arbitrary, we conclude that $f(x_0) \leq g(x_0)$.

\textit{(c):} Its proof is analogous to the previous part and is left as an exercise [Exercise \ref{E537}].

\textit{(d):} From the previous parts, $f$ is continuous at $x_0$ if and only if $f$ is lower semicontinuous and upper semicontinuous at $x_0$, if and only if $g(x_0)=f(x_0)=h(x_0)$.
\end{proof}

\begin{lemma} \label{536}
If $f:[a,b] \to \mathbb{R}$ is bounded and nonnegative, then:
$$
\underline{R}\int_{a}^{b} f(x)\,dx = \int_{[a,b]} g\,d\lambda \quad\mbox{and}\quad \overline{R}\int_{a}^{b} f(x)d\,dx=\int_{[a,b]} h\,d\lambda
$$
where $g$ and $h$ denote the lower and upper envelopes of $f$, respectively.
\end{lemma}

\begin{proof}
Let $(P_k)$ be a sequence of partitions of $[a,b]$ such that $P_{k} \subset P_{k+1}$ and $\Vert \,P_k\, \Vert \to 0$. If
$P_k=\left\{a=x^{k}_{0}<x^{k}_{1}<\ldots < x^{k}_{n_k}=b\right\}$, define
$$
\xi_k(x):=\left\{
\begin{array}{lcl}
\inf \left\{f(\zeta) \,:\, \zeta \in (x_{i-1}^{k},x_{i}^{k})\right\} & &\mbox{if}\,\,\,\,x \in (x_{i-1}^{k},x_{i}^{k}),\\
g(x) & &\mbox{if}\,\,\,\, x=x^{k}_{i}\,\,\,\mbox{for some}\,\,\,i=1,\ldots,n_{k}
\end{array}
\right.
$$

$$
\Xi_{k}(x):=\left\{
\begin{array}{lcl}
\sup \left\{f(\zeta) \,:\, \zeta \in (x_{i-1}^{k},x_{i}^{k})\right\}& &\mbox{if}\,\,\,\,x \in (x_{i-1}^{k},x_{i}^{k}),\\
h(x) & &\mbox{if}\,\,\,\, x=x^{k}_{i}\,\,\,\mbox{for some}\,\,\,i=1,\ldots,n_{k}
\end{array}
\right.
$$
where $g$ and $h$ denote the lower and upper envelopes of $f$, respectively. For each $k \in \mathbb{N}$, we may rewrite $\xi_{k}$ and $\Xi_{k}$ on $[a,b]$ as follows
$$
\begin{aligned}
\xi_{k}&=\sum_{i=1}^{n_{k}} \inf\left\{ f(\zeta)\,:\, \zeta \in (x_{i-1}^{k},x_{i}^{k})\right\}\,\chi_{(x_{i-1}^{k},x_{i}^{k})} + g(x)\chi_{\{x=x^{k}_{i}\,:\,\mbox{for some}\,\,i=1,\ldots,n_{k}\}},\\
\Xi_{k}&=\sum_{i=1}^{n_{k}} \sup\left\{ f(\zeta)\,:\, \zeta \in (x_{i-1}^{k},x_{i}^{k})\right\}\,\chi_{(x_{i-1}^{k},x_{i}^{k})} + h(x)\chi_{\{x=x^{k}_{i}\,:\,\mbox{for some}\,\,i=1,\ldots,n_{k}\}}.
\end{aligned}
$$

Consequently, $\xi_{k}$ and $\Xi_{k}$ are nonnegative Borel-simple functions and
$$
\begin{aligned}
\int_{[a,b]} \xi_{k}\,d\lambda &= \sum_{i=1}^{n_{k}} \inf\left\{ f(\zeta)\,:\, \zeta \in (x_{i-1}^{k},x_{i}^{k})\right\}\lambda((x_{i-1},x_{i})) = \underline{\mathcal{S}}(f,P_{k}),\\
\int_{[a,b]} \Xi_{k}\,d\lambda &=\sum_{i=1}^{n_{k}} \sup\left\{ f(\zeta)\,:\, \zeta \in (x_{i-1}^{k},x_{i}^{k})\right\}\lambda((x_{i-1},x_{i})) =\overline{\mathcal{S}}(f,P_{k}).
\end{aligned}
$$

On the other hand, for every $k \leq k+1$ we have
$\xi_{k} \leq \xi_{k+1} \leq f \leq \Xi_{k+1} \leq \Xi_{k}$
since $P_{k}\subset P_{k+1}$. Therefore, $(\xi_k)$ is an increasing sequence of functions in $\mathbb{S}^{+}([a,b],\mathcal{B}([a,b]))$ and $(\Xi_k)$ is a decreasing sequence of functions in $\mathbb{S}^{+}([a,b],\mathcal{B}([a,b]))$ such that $\int_{[a,b]} \Xi_{1}\,d\lambda < +\infty$ since $f$ is bounded. Applying Theorem \ref{517} and Corollary \ref{519} to $(\xi_k)$ and $(\Xi_k)$ respectively, we obtain
\begin{eqnarray} \label{F52}
\int_{[a,b]} \lim_{k\,\to\,\infty} \xi_k\,d\lambda =\lim_{k\,\to\,\infty} \int_{[a,b]} \xi_{k}\,d\lambda = \lim_{k\,\to\,\infty} \underline{\mathcal{S}}(f,P_{k}) =\sup_{k\,\in\,\mathbb{N}}\underline{\mathcal{S}}(f,P_{k}) =\underline{R}\int_{a}^{b} f(x)\,dx
\end{eqnarray}
\begin{eqnarray}  \label{F53}
\int_{[a,b]} \lim_{k\,\to\,\infty} \Xi_k\,d\lambda =\lim_{k\,\to\,\infty} \int_{[a,b]} \Xi_{k}\,d\lambda = \lim_{k\,\to\,\infty} \overline{\mathcal{S}}(f,P_{k}) =\inf_{k\,\in\,\mathbb{N}}\overline{\mathcal{S}}(f,P_{k}) =\overline{R}\int_{a}^{b} f(x)\,dx
\end{eqnarray}
where the last equalities follow from the fact that $\Vert \,P_{k} \,\Vert \rightarrow 0$ [Exercise \ref{E538}]. Now, let us see that $\xi_{k}\,\to\,g$ and $\Xi_{k}\,\to\,h$ on $[a,b]$, for which we consider the following two cases.

{\scshape Case 1.}\quad $x \not\in \bigcup_{k=1}^{\infty} P_k$.

In this case, we have $x \not\in P_k$ for every $k \in \mathbb{N}$. Consequently,
$$
\xi_k(x) \leq g(x) \leq f(x) \leq h(x) \leq \Xi_k(x) \quad\,\,\forall k \in \mathbb{N}.
$$

Taking the limit as $k \to \infty$, we conclude that
\begin{eqnarray} \label{F54}
\lim_{k \to \infty} \xi_k(x) \leq g(x) \leq f(x) \leq h(x) \leq \lim_{k \to \infty} \Xi_k(x).
\end{eqnarray}

On the other hand, given $\delta >0$, we can find $k \in \mathbb{N}$ sufficiently large such that
$(x_{i-1}^{k},x_{i}^{k}) \subset (x-\delta,x+\delta)$. Consequently,
$$
g_{\delta}(x) =\inf\left\{ f(t)\,:\,t\in (x-\delta,x+\delta)\cap [a,b] \right\} \leq \inf\left\{f(\zeta)\,:\,\zeta \in (x_{i-1}^{k},x_{i}^{k}) \right\}=\xi_{k}(x),
$$
$$
\Xi_{k}(t)=\sup\left\{f(\zeta)\,:\,\zeta \in (x_{i-1}^{k},x_{i}^{k}) \right\} \leq \sup\left\{ f(t)\,:\,t\in (x-\delta,x+\delta)\cap [a,b] \right\}=h_{\delta}(x).
$$

Taking once again the limit as $k \to \infty$, we obtain
$$
g_{\delta}(x)\leq \lim_{k \to \infty} \xi_{k}(x)\quad\mbox{and}\quad \lim_{k \to \infty} \Xi_{k}(x) \leq h_{\delta}(x),
$$
and therefore,
\begin{eqnarray} \label{F55}
g(x) \leq \lim_{k \to \infty} \xi_{k}(x)\quad\mbox{and}\quad \lim_{k \to \infty} \Xi_{k}(x) \leq h(x).
\end{eqnarray}

From inequalities (\ref{F54}) and (\ref{F55}) it follows that
$$
g(x) = \lim_{k \to \infty} \xi_{k}(x)\quad\mbox{and}\quad \lim_{k \to \infty} \Xi_{k}(x) = h(x), \qquad \forall x \not\in \bigcup_{k=1}^{\infty} P_k.
$$

{\scshape Case 2.}\quad $x \in \bigcup_{k=1}^{\infty} P_k$.

Let $k_0$ be the first natural number such that $x \in P_{k_0}$. Since $(P_k)$ is increasing, then $x \in P_k$ for every $k \geq k_0$. Consequently, $\xi_{k}(x)=g(x)$ and $\Xi_{k}(x)=h(x)$ for every $k \geq k_0$, so necessarily
$\displaystyle\lim_{k \to \infty} \xi_{k}(x)=g(x)$ and
$\displaystyle\lim_{k \to \infty} \Xi_{k}(x)=h(x)$.

From both cases, we conclude that $\xi_{k}(t)\,\to\,g(t)$ and $\Xi_{k}(t)\,\to\,h(t)$ for every $x \in [a,b]$. Thus, substituting these equalities into (\ref{F52}) and (\ref{F53}), we obtain
$$
\underline{R}\int_{a}^{b} f(x)\,dx = \int_{[a,b]} g\,d\lambda \quad\mbox{and}\quad \overline{R}\int_{a}^{b} f(x)\,dx=\int_{[a,b]} h\,d\lambda,
$$
as stated.
\end{proof}

\begin{proof}[Proof of Theorem \ref{532}]
\textit{(a):} If $f$ is Riemann-integrable on $[a,b]$, then
$$
\underline{R}\int_{a}^{b} f(x)\,dx = \overline{R}\int_{a}^{b} f(x)\,dx.
$$

Lemma \ref{536} then implies that
$$
\int_{[a,b]} g\,d\lambda =\int_{[a,b]} h\,d\lambda 
$$
where $g$ and $h$ denote the lower and upper envelopes of $f$, respectively. That is,
$\int_{[a,b]} (g-h)\,d\lambda =0$
and, consequently, $g=h$ $\lambda$ a.e. Lemma \ref{535} ensures that
$g=f=h$ $\lambda$ a.e. and, therefore, that $f$ is continuous $\lambda$ a.e.

Conversely: If $f$ is continuous $\lambda$ a.e., then
$g=f=h$ $\lambda$ a.e. and, therefore,
$g-h=0$ $\lambda$ a.e. Theorem \ref{528} ensures that
$\int_{[a,b]} g\,d\lambda = \int_{[a,b]} h\,d\lambda$
and from Lemma \ref{536} we conclude that
$\underline{R}\int_{a}^{b} f(x)\,dx=\overline{R}\int_{a}^{b} f(x)\,dx$.
That is, $f$ is Riemann-integrable on $[a,b]$.  

\textit{(b):} If $f$ is Borel-measurable and Riemann-integrable on $[a,b]$, Theorem \ref{528} together with Lemma \ref{536} ensure that
$g=h$ $\lambda$ a.e. Since $g \leq f \leq h$ and
$g=h$ $\lambda$ a.e., then
$f=g$ $\lambda$ a.e.

Applying Theorem \ref{528} and Lemma \ref{536}, we obtain
$$
\int_{[a,b]} f\,d\lambda = \int_{[a,b]} g\,d\lambda = \underline{R}\int_{a}^{b} f(x)\,dx = R\int_{a}^{b} f(x)\,dx,
$$
as stated.
\end{proof}

We conclude this chapter with some examples.

\begin{example}  \label{537}
Let $\alpha \geq 1$ and $\beta >0$, and let $f:[0,1] \to \mathbb{R}$ be the continuous function given by
$$
f(x)=\frac{x^{\alpha-1}}{1+x^{\beta}}.
$$

Then,
$$
R\int_{0}^{1} f(x)\,dx=\sum_{k=0}^{\infty} \frac{(-1)^{k}}{1+k\beta}.
$$
\end{example} 

\begin{proof}
We may decompose the function $f$ as follows
$$
f(x)=\frac{x^{\alpha-1}}{1+x^{\beta}}=x^{\alpha-1}\,\frac{1-x^{\beta}}{1-x^{2\beta}}=(1-x^{\beta})\,x^{\alpha-1}\,\frac{1}{1-x^{2\beta}},
$$
where
$$
\frac{1}{1-x^{2\beta}}=\sum_{k=0}^{\infty} x^{2\beta k} \quad \mbox{if}\,\,x\in [0,1).
$$

Thus,
$$
f(x)=\sum_{k=0}^{\infty} f_k(x) = \sum_{k=0}^{\infty}(1-x^{\beta})\,x^{\alpha-1+2\beta k}\quad\mbox{if}\,\,x\in [0,1),
$$
with
$f_k(x):=(1-x^{\beta})\,x^{\alpha-1+2\beta k} \geq 0$.

Applying the Beppo-Levi theorem (see Theorem \ref{520}), we obtain
\begin{eqnarray} \label{F56}
\int_{[0,1]}f\,d\lambda =\sum_{k=0}^{\infty} \int_{[0,1]}f_k\,d\lambda.
\end{eqnarray}

Since $f$ and $f_k$ are Riemann-integrable on $[0,1]$, Theorem \ref{532} allows us to rewrite expression (\ref{F56}) as
$$
R\int_{0}^{1} f(x)\,dx=\sum_{k=0}^{\infty}R\int_{0}^{1} f_k(x)\,dx.
$$

Consequently,
$$
R\int_{0}^{1} f(x)\,dx=\sum_{k=0}^{\infty}R\int_{0}^{1} f_k(x)\,dx =\sum_{k=0}^{\infty}\left(\frac{1}{\alpha+2k\beta}-\frac{1}{\alpha+(2k+1)\beta} \right)=\sum_{k=0}^{\infty} \frac{(-1)^{k}}{1+k\beta}
$$
as stated.
\end{proof}

Observe that in the previous example, if $\alpha=1=\beta$, then
$$
R\int_{0}^{1} \frac{1}{1+x}\,dx = \mbox{log}(2) = \sum_{k=0}^{\infty} \frac{(-1)^{k}}{1+k},
$$
and if $\alpha=1$ and $\beta=2$, then
$$
R\int_{0}^{1} \frac{1}{1+x^2}\,dx = \frac{\pi}{4}= \sum_{k=0}^{\infty} \frac{(-1)^{k}}{1+2k}.
$$

\begin{figure}[ht!]
\begin{minipage}[l]{0.5\textwidth}
\begin{center}
\begin{tikzpicture}[xscale=4, yscale=2]
\draw[domain= 0:1,thick] plot(\x,{1/(1+\x)} );

\draw[->,gray] (0,0)--(1.2,0); 
\draw[->,gray] (0,0)--(0,1.2); 

\draw(0,1)node[left]{$_{1}$};
\draw(0,0)node[left]{$_{0}$};
\draw(1,0)node[below]{$_{1}$};
\draw[dotted] (1,0)--(1,0.5);
\end{tikzpicture}
\begin{center}
$f(x)$ with $\alpha=1=\beta$
\end{center}

\end{center}
\end{minipage} \hfill 
 \begin{minipage}[r]{0.5\textwidth}
\begin{center}
\begin{tikzpicture}[xscale=4, yscale=2]
\draw[domain= 0:1,thick] plot(\x,{1/(1+\x^2)} );

\draw[->,gray] (0,0)--(1.2,0); 
\draw[->,gray] (0,0)--(0,1.2); 

\draw(0,1)node[left]{$_{1}$};
\draw(0,0)node[left]{$_{0}$};
\draw(1,0)node[below]{$_{1}$};
\draw[dotted] (1,0)--(1,0.5);
\end{tikzpicture}
\begin{center}
$f(x)$ with $\alpha=1$ and $\beta=2$
\end{center}
\end{center}
\end{minipage}
\end{figure}

The following function is a modification of Dirichlet's function (see Example \ref{537}) and is known as Thomae's function in honor of Carl Johannes Thomae\footnote{Carl Johannes Thomae (1840-1921) was a German mathematician. His research dealt with function theory and with what German-speaking mathematicians often call ``Epsilontik'', the rigorous development of analysis, differential geometry, and topology using epsilon neighborhoods in the style of Weierstrass.}.

\begin{example}[Thomae's function] \label{538} \index{function!Thomae function}
Thomae's function $f:[0,1]\to \mathbb{R}$ given by:
$$
f(x)=\left\{
\begin{array}{lcl}
\frac{1}{q} & &\mbox{if}\,\,\,\,x\in \widetilde{\mathbb{Q}},\,\,x=\frac{p}{q}\,\,\mbox{with}\,\,p,q \in \mathbb{Z}^{+}\,\,\mbox{and}\,\,\mbox{\rm gcd}(p,q)=1,\\
\\
0 & & \mbox{if}\,\,\,\,x \not\in \widetilde{\mathbb{Q}}\,\,\mbox{or}\,\,x=0\,,
\end{array}
\right.
$$
is Riemann-integrable on $[0,1]$.
\end{example}

\begin{figure}[ht!]
\begin{minipage}[c]{0.5\textwidth}
\begin{center}
\includegraphics[scale=1]{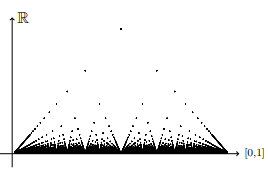} 
\end{center}
\end{minipage} \hfill \begin{minipage}[c]{0.5\textwidth}
\begin{center}
\includegraphics[scale=0.275]{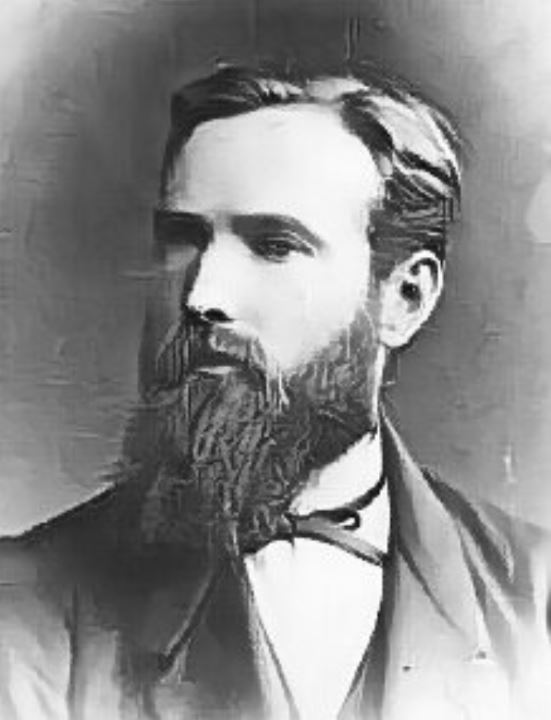} 
\begin{center}
C. Thomae (1840-1921)
\end{center}
\end{center}
\end{minipage}
\end{figure}

\begin{proof}
We shall prove that $f$ is continuous on the set $([0,1]\smallsetminus \widetilde{\mathbb{Q}}) \cup \{0\}$ but discontinuous on $\widetilde{\mathbb{Q}}$. The case when $x_0=0$ is immediate. Thus, assume that $x_{0} \in [0,1]\smallsetminus \widetilde{\mathbb{Q}}$. 

Let $\varepsilon >0$ be arbitrary. Choose $k_{0} \in \mathbb{N}$ such that $\frac{1}{k_0}<\varepsilon$. Denote by $Q_{0}$ the set of rational numbers $r$ in $\widetilde{\mathbb{Q}}$ such that $r=\frac{p}{q}$ with $p,q \in \mathbb{Z}^{+}$, $\mbox{gcd}(p,q)=1$, and $q \leq k_{0}$. This set is finite since the function $\iota: Q_0 \to \left\{1,\ldots,k_0\right\} \times \left\{1,\ldots,k_0\right\}$ given by $\iota\left( \frac{p}{q} \right)=(p,q)$ is injective. Hence, the set $\left\{ |r-x_0|\,:\,r \in Q_0 \right\}$ is finite and therefore has a minimum. Let
$m_{\ast}:=\inf\left\{ |r-x_0|\,:\,r \in Q_0 \right\}$.
Observe that $m_{\ast}>0$ since $x_0 \not \in Q_0$. Define $\delta:=\frac{m_{\ast}}{2}>0$ and let $x \in [0,1]$ satisfy $|x-x_0|<\delta$.

{\scshape Case 1:}\quad If $x \not \in \widetilde{\mathbb{Q}}$, then $f(x)=0$ and the result follows trivially.

{\scshape Case 2:}\quad Suppose that $x \in \widetilde{\mathbb{Q}}$. It is clear that $x \not \in Q_0$ and, consequently, $x=\frac{p}{q}$ with $\mbox{gcd}(p,q)=1$ and $q > k_{0}$. Therefore,
$|f(x)-f(x_0)|=\frac{1}{q}<\frac{1}{k_0}<\varepsilon$.

This proves that $f$ is continuous on $([0,1]\smallsetminus \widetilde{\mathbb{Q}}) \cup \{0\}$.

On the other hand, let $x_0 \in \widetilde{\mathbb{Q}}$. Then there exist $p_0,q_0 \in \mathbb{N}$ with $\mbox{gcd}(p_0,q_0)=1$ such that $x_0=\frac{p_0}{q_0}$. Define $\varepsilon_0:=\frac{1}{2q_0}>0$ and let $\delta >0$. Since the set of irrational numbers is dense in $\mathbb{R}$, we may choose
$x \in (x_0-\delta,x_0+\delta)\cap [0,1]$
irrational. Hence,
$|x-x_0|<\delta$
and
$|f(x)-f(x_0)|=\frac{1}{q_0} > \frac{1}{2q_0}=\varepsilon_0$. 

Therefore, $f$ is not continuous on $\widetilde{\mathbb{Q}}$.

Since $\widetilde{\mathbb{Q}}$ is countable, we have
$\lambda(\widetilde{\mathbb{Q}})=0$
[Exercise \ref{E440}], which shows that $f$ is continuous $\lambda$ a.e. Theorem \ref{532} then ensures that $f$ is Riemann-integrable on $[0,1]$.
\end{proof}

\section{Exercises}

\begin{exercise} \label{E51}
Let $(X,\mathsf{S},\mu)$ be a measure space. Prove that for every $A \in \mathsf{S}$, $\int_{X} \chi_{A}\,d\mu = \mu(A).$
\end{exercise}

\begin{exercise} \label{E52}
Let $(X,\mathsf{S},\mu)$ be a measure space. Prove that for every constant $c > 0$, the function defined by $f(x):=c$ is $\mathsf{S}$-simple and satisfies
$\int_{X} f\,d\mu = c\mu(X).$
\end{exercise}

\begin{exercise} \label{E53}
Let $(\mathbb{R},\mathcal{B}(\mathbb{R}))$ be a measurable space. Consider $\delta_0$ the unit measure concentrated at $0$, $\mu^{\sharp}$ the counting measure, and $\lambda$ the Lebesgue measure on $\mathcal{B}(\mathbb{R})$. For each item, find a function
$f \in \mathbb{S}^{+}(\mathbb{R},\mathcal{B}(\mathbb{R}))$
such that
\begin{itemize}
    \item[(a)] $\int_{\mathbb{R}}f\,d\delta_{0}=0$, $\int_{\mathbb{R}}f\,d\mu^{\sharp}=+\infty$ and $\int_{\mathbb{R}}f\,d\lambda=1$.
    \item[(b)] $\int_{\mathbb{R}}f\,d\delta_{0}=1$, $\int_{\mathbb{R}}f\,d\mu^{\sharp}=+\infty$ and $\int_{\mathbb{R}}f\,d\lambda=0$.
    \item[(c)] $\int_{\mathbb{R}}f\,d\delta_{0}=1$, $\int_{\mathbb{R}}f\,d\mu^{\sharp}=+\infty$ and $\int_{\mathbb{R}}f\,d\lambda=+\infty$.
     \item[(d)] $\int_{\mathbb{R}}f\,d\delta_{0}=0$, $\int_{\mathbb{R}}f\,d\mu^{\sharp}=1$ and $\int_{\mathbb{R}}f\,d\lambda=0$.
     \item[(e)] $\int_{\mathbb{R}}f\,d\delta_{0}=1$, $\int_{\mathbb{R}}f\,d\mu^{\sharp}=+\infty$ and $\int_{\mathbb{R}}f\,d\lambda=1$.
\end{itemize}
\end{exercise}

\begin{exercise} \label{E54}
Let $X=\mathbb{R}$,
$$
\mathsf{S}=\{A \subset \mathbb{R}\,:\,A\,\,\text{or}\,\,\mathbb{R} \smallsetminus A\,\,\text{is finite or countable}\},
$$
and let $\mu:\mathsf{S} \to \mathbb{R}$ be the measure given by
$$
\mu(A)=\left\{
\begin{array}{lcl}
0 &  & \text{if} \quad A \,\,\text{is finite or countable},  \\
1 &  & \text{if} \quad \mathbb{R} \smallsetminus A \,\,\text{is finite or countable}.
\end{array}
\right.
$$
\begin{itemize}
    \item[(a)] Describe the functions $f:X \to \overline{\mathbb{R}}$ that are $\mathsf{S}$-measurable.
    \item[(b)] Let $f \in \overline{\mathbb{M}^{+}}(\mathbb{R},\mathsf{S})$. Give a description of the integral
$\int_{\mathbb{R}} f\,d\mu$
in terms of the structure of $f$.
\end{itemize}
\end{exercise}

\begin{exercise} \label{E55} \index{measure!multiplicative}
Let $(X,\mathsf{S})$ be a measurable space. A measure
$\mu:\mathsf{S} \to \mathbb{R}$
on $(X,\mathsf{S})$ is called a \textbf{multiplicative measure} if
$\mu(A \cap B)=\mu(A)\mu(B)$
for every $A,B \in \mathsf{S}$. Prove the following:
\begin{itemize}
    \item[(a)] A measure $\mu$ on $(X,\mathsf{S})$ is multiplicative if and only if for every
$f,g \in \overline{\mathbb{M}}^{+}(X,\mathsf{S})$,
$$
\int_{X}fg\,d\mu = \left(\int_{X} f\,d\mu \right)\left(  \int_{X} g\,d\mu\right).
$$
\item[(b)] Determine all multiplicative measures on the measurable space
$(\mathbb{R},\mathcal{B}(\mathbb{R}))$.
\end{itemize}
\end{exercise}

\begin{exercise} \label{E56}
Let $(X,\mathsf{S})$ be a measurable space and let
$\mu,\nu :\mathsf{S} \to \overline{\mathbb{R}}$
be measures such that
$\mu(A) \leq \nu(A)$
for every $A \in \mathsf{S}$. Prove that for every
$s \in \mathbb{S}^{+}(X,\mathsf{S})$
and every
$A \in \mathsf{S}$,
one has
$$
\int_{A} s\,d\mu \leq \int_{A} s\,d\nu.
$$
\end{exercise}

\begin{exercise} \label{E57}
Let $(X,\mathsf{S},\mu)$ be a measure space and let $N \in \mathbb{N}$ be fixed. Prove that for every
$s_1,\ldots,s_N \in \mathbb{S}^{+}(X,\mathsf{S})$
and every
$c_1,\ldots,c_N \geq 0$,
the following equality holds:
$$
\int_{A} \left(\sum_{j=1}^{N} c_js_j  \right)\,d\mu = \sum_{j=1}^{N} c_j\int_{A} s_j\,d\mu,\quad\text{for every }\,\,A \in \mathsf{S}.
$$
\end{exercise}

\begin{exercise} \label{E58}
Let $(X,\mathsf{S},\mu)$ be a measure space. Prove that for every
$f,g \in \overline{\mathbb{M}}^{+}(X,\mathsf{S})$
such that $f\leq g$, one has
$$
\int_{X} f\,d\mu \leq \int_{X} g\,d\mu.
$$
\end{exercise}

\begin{exercise} \label{E59}
Let $(X,\mathsf{S},\mu)$ be a measure space, and let
$f \in \overline{\mathbb{M}}^{+}(X,\mathsf{S})$
and
$c >0$
be fixed. Prove that the function
$\iota: \underline{\mathcal{S}}(f) \to \underline{\mathcal{S}}(cf)$
defined by
$\iota(s):=cs$
is injective.
\end{exercise}

\begin{exercise} \label{E510}
Let $(X,\mathsf{S},\mu)$ be a measure space,
$f,g \in \overline{\mathbb{M}}^{+}(X,\mathsf{S})$,
and
$s \in \mathbb{S}^{+}(X,\mathsf{S})$
such that
$s \leq f+g$.
For every
$k,j \in \mathbb{N}$
with
$j \leq k$,
define the following functions:
$$
\begin{aligned}
\varphi_{k}(x)&:=\sup\left\{ \left(\frac{j}{k} \right)s(x)\,:\, \left(\frac{j}{k} \right)s(x) \leq f(x)\right\},\\
\psi_{k}(x)&:=\sup\left\{ \left(1-\frac{1}{k} \right)s(x)-\varphi(x),0 \right\}.
\end{aligned}
$$

Prove that
$\left(1-\frac{1}{k} \right)s \leq \varphi_k + \psi_k$,
$\varphi_k \leq f$,
and
$\psi_{k} \leq g$.
\end{exercise}

{\setlength{\parindent}{0pt}
\begin{exercise} \label{E511}
Let $(X,\mathsf{S},\mu)$ be a measure space,
$f,g \in \overline{\mathbb{M}}^{+}(X,\mathsf{S})$,
and
$c \geq 0$.
Without using the monotone convergence theorem, prove that for every
$A \in \mathsf{S}$:
\begin{itemize}
\item[(a)] $\int_{A} cf\,d\mu = c\int_{A} f\,d\mu$.
\item[(b)] $\int_{A}(f+g)\,d\mu = \int_{A}f\,d\mu + \int_{A}g\,d\mu$.
\end{itemize}
\end{exercise}

(Hint: Use Exercise \ref{E58} for part (a) and Exercise \ref{E59} for part (b)).}

\begin{exercise} \label{E512}
Let $(X,\mathsf{S},\mu)$ be a measure space. Prove that
$$
\int_X f\,d\mu =\sup\left\{ \int_{A} f\,d\mu\,:\, A\in \mathsf{S},\,\,\,\mu(A)<+\infty \right\},
$$
for every
$f \in \mathbb{M}^{+}(X,\mathsf{S})$
such that
$\int_X f \,d\mu <+\infty$.
\end{exercise}

\begin{exercise} \label{E513}
Let $(X,\mathsf{S},\mu)$ be a measure space and $f \in \overline{\mathbb{M}}^{+}(X,\mathsf{S})$. Prove that, if $N \in \mathcal{N}(\mu)$, then $\int_{A} f\,d\mu = \int_{X} f\,d\mu$ with $A=X\smallsetminus N$. 
\end{exercise}

{\setlength{\parindent}{0pt}
\begin{exercise} \label{E514}
Let $(X,\mathsf{S},\mu)$ be a measure space and $f \in \overline{\mathbb{M}}^{+}(X,\mathsf{S})$. Prove that, if $\omega:\mathsf{S} \to \overline{\mathbb{R}}$ is the measure on $(X,\mathsf{S})$ given by $\omega(A):=c\mu(A)$ for some $c \geq 0$, then
$$
\int_{X} f\,d\omega = c\int_{X} f\,d\mu.
$$

If now $\nu$ is another measure on $(X,\mathsf{S})$ and $\omega$ is the measure given by $\omega(A):=\mu(A)+\nu(A),$
then
$$
\int_{X} f\,d\omega = \int_{X} f\,d\mu + \int_{X}f\,d\nu.
$$
\end{exercise}

(Hint: First prove the case when $f$ is simple and use the monotone convergence theorem for the general case).}

{\setlength{\parindent}{0pt}
\begin{exercise} \label{E515}
Let $(X,\mathsf{S},\mu)$ be a measure space and let $f \in \overline{\mathbb{M}}^{+}(X,\mathsf{S})$ be given. Define the function $\mu_{f}:\mathsf{S} \to \overline{\mathbb{R}}$ by $\mu_{f}(A):=\int_{A} f\,d\mu$. Prove that $\mu_{f}$ is a measure on $(X,\mathsf{S})$.
\end{exercise}

(Hint: Use the monotone convergence theorem).}

\begin{exercise} \label{E516}
Let $(X,\mathsf{S},\mu)$ be a $\sigma$-finite measure space and let $f:X \to \overline{\mathbb{R}}$ be an $\mathsf{S}$-measurable function. Prove that there exists a sequence of $\mathsf{S}$-simple functions $(s_k)$ such that $s_k(x) \to f(x)$ for every $x \in X$ and
$$
\mu \left( \left\{ x \in X\,:\, s_k(x) \neq 0\right\}\right)<+\infty
$$
for every $k \in \mathbb{N}$.
\end{exercise}

\begin{exercise} \label{E517} \index{Lebesgue!sum}
Let $(X,\mathsf{S},\mu)$ be a measure space and $f \in \mathbb{M}^{+}(X,\mathsf{S})$. Let $P=(0=t_0<t_1<\ldots) $ be a partition of $[0,+\infty)$. Prove that, if $\int_X f\,d\mu <+\infty$, then
$$
\int_X f\,d\mu =\lim_{\Vert P \Vert \to 0} \sum_{k=0}^{\infty} \varsigma_{k}\,\mu\left( \left\{ x\in X\,:\, t_k \leq f(x) < t_{k+1} \right\} \right)
$$
with $(\varsigma_k)$ any sequence of points satisfying $\varsigma_{0}=0$ and $\varsigma \in [t_k,t_{k+1})$ if $k \geq 1$. (Each of the above expressions is known as a Lebesgue sum).
\end{exercise}

{\setlength{\parindent}{0pt}
\begin{exercise} \label{E518}
Let $(X,\mathsf{S},\mu)$ be a measure space and let $f \in \overline{\mathbb{M}}^{+}(X,\mathsf{S})$ be given such that $\int_X f\,d\mu <+\infty$. Prove the following statements:
\begin{itemize}
\item[(a)] $\mu \left( \left\{ x \in X\,:\,f(x)\geq \varepsilon \right\} \right)<+\infty$ for every $\varepsilon >0$.
\item[(b)] $\mu \left( \left\{ x \in X\,:\,f(x)=+\infty \right\} \right)=0$.
\item[(c)] $\left\{ x \in X\,:\, f(x) >0 \right\}$ is $\sigma$-finite.
\end{itemize}
\end{exercise}

(Hint: Define the set $F_{\varepsilon}:=\left\{ x \in X\,:\, \varepsilon \leq f(x)\right\}$ for every $\varepsilon>0$. For part (b) use Theorem \ref{414} considering the set $F_{k}$ for each $k \in \mathbb{N}$ and in part (c) consider the set $F_{1/k}$ for each $k \in \mathbb{N}$).}

\begin{exercise} \label{E519}
Let $(X,\mathsf{S},\mu)$ be a finite measure space and let $f \in \mathbb{M}^{+}(X,\mathsf{S})$. Prove that
$$
\int_X f\,d\mu < +\infty\quad\mbox{if and only if}\quad \sum_{k=0}^{\infty} 2^{k}\mu\left( \left\{ x \in X\,:\, f(x) \geq 2^{k}\right\} \right)<+\infty.
$$
\end{exercise}

\begin{exercise} \label{E520}
Let $(X,\mathsf{S},\mu)$ be a finite measure space and let $f \in \mathbb{M}^{+}(X,\mathsf{S})$. Prove that, if $f$ is bounded on $X$, then
$$
\int_X f\,d\mu < +\infty\quad\mbox{if and only if}\quad \sum_{k=0}^{\infty} \frac{1}{2^{k}}\mu\left( \left\{ x \in X\,:\, f(x) \geq 2^{-k}\right\} \right)<+\infty.
$$
\end{exercise}

{\setlength{\parindent}{0pt}
\begin{exercise} \label{E521}
Let $(X,\mathsf{S},\mu)$ be a measure space and let $f \in \overline{\mathbb{M}}^{+}(X,\mathsf{S})$ be given. Prove that if $\int_X f\,d\mu <+\infty$ then, for every $\varepsilon >0$, there exists a set $A_{\varepsilon} \in \mathsf{S}$ such that $\mu(A_{\varepsilon})<+\infty$ and
$$
\int_X f\,d\mu \leq \int_{A_{\varepsilon}} f\,d\mu + \varepsilon.
$$
\end{exercise}

(Hint: Let $A=\left\{x \in X\,:\, f(x) \neq 0 \right\}$. Define the sequence $(A_k)$ of elements of $\mathsf{S}$ by
$A_k:=\left\{x \in X\,:\, \frac{1}{k} \leq f(x) \right\}$.
Prove that $\int_X f\,d\mu = \int_{A} f\,d\mu$ and that the sequence defined above is nondecreasing such that $A_k\,\to\, A$ with $\mu(A_k)=k\int_X f\,d\mu$. Use Exercise \ref{E518} and Theorem \ref{517}).}

\begin{exercise} \label{E522}
Let $(X,\mathsf{S},\mu)$ be a measure space and let $f \in \overline{\mathbb{M}}^{+}(X,\mathsf{S})$ be such that
$$
0 < \iota :=\int_X f\,d\mu < +\infty.
$$

Prove the following:
$$
\lim_{k \to \infty} \int_X k\log\left(1+\frac{f^{\alpha}}{k^{\alpha}} \right)\,d\mu=\left\{
\begin{array}{lcl}
+\infty & & \mbox{if}\,\,\,\,\alpha < 1,\\
\iota & &\mbox{if}\,\,\,\, \alpha=1,\,\,\,\,\mbox{for}\,\,\,0 < \alpha < +\infty,\\
0 & &\mbox{if}\,\,\,\, \alpha >1.
\end{array}
\right.
$$
\end{exercise}

\begin{exercise} \label{E523}
Let $(X,\mathsf{S},\mu)$ be a measure space and let $f\in \overline{\mathbb{M}}^{+}(X,\mathsf{S})$ with $\int_X f\,d\mu <+\infty$. For each $0<\delta<1$ define the function $f_{\delta} : X \to [0,+\infty]$ as follows:
$$
f_{\delta}(x):=\left\{
\begin{array}{lcl}
0 & &\mbox{if}\,\,\,\, f(x) \leq \delta,\\
f(x) & &\mbox{if}\,\,\,\, \delta < f(x) \leq \delta^{-1},\\
\delta^{-1} & &\mbox{if}\,\,\,\, f(x)>\delta^{-1}.
\end{array}
\right.
$$

Prove that $f_{\delta}$ is $\mathsf{S}$-measurable for each $0<\delta<1$ and that
$\displaystyle\lim_{\delta \to 0} \int_X f_{\delta}\,d\mu = \int_X f\,d\mu.$
\end{exercise}

\begin{exercise} \label{E524}
Let $X=\mathbb{R}$, $\mathsf{S}=\mathcal{B}(\mathbb{R})$ and $\mu=\lambda$ be the Lebesgue measure on $\mathsf{S}$. Define the sequence of functions $f_k : \mathbb{R} \to \mathbb{R}$ by
$$
f_k(x)=c_k\,\left(1-\frac{x}{k} \right)^{k}\,\chi_{[0,k]}(x).
$$

\begin{itemize}
\item[(a)] Find $c_k \in \mathbb{R}$ such that $\int_{\mathbb{R}} f_k\,d\lambda=1$.
\item[(b)] Prove that the sequence $(f_k(x))$ converges for every $x  \in \mathbb{R}$.
\item[(c)] Let $f(x):=\displaystyle\lim_{k \to \infty} f_k(x)$. For each $A \in \mathcal{B}(\mathbb{R})$ define
$$
\nu_{k}(A):=\int_{A} f_k\,d\lambda \quad\mbox{and}\quad \nu(A):=\int_{A} f\,d\lambda. 
$$

Prove that $\nu_k \to \nu$ uniformly on $\mathcal{B}(\mathbb{R})$.
\end{itemize}
\end{exercise}

\begin{exercise} \label{E525}
Let $(X,\mathsf{S},\mu)$ be a measure space and let $f\in \overline{\mathbb{M}}^{+}(X,\mathsf{S})$. Prove that, if $(\varphi_k)$ and $(\phi_k)$ are two nondecreasing sequences of elements of $\mathbb{S}^{+}(X,\mathsf{S})$ such that they converge pointwise to $f$ on $X$, then
$$
\lim_{k \to \infty}\int_{X} \varphi_k\,d\mu = \lim_{k \to \infty} \int_{X} \phi_k\,d\mu.
$$
\end{exercise}

\begin{exercise} \label{E526}
Let $(\mathbb{R},\mathcal{B}(\mathbb{R}),\lambda)$ be a measure space. Define the sequence of functions $f_k:\mathbb{R} \to \mathbb{R}$ by $f_k(x):=\chi_{[0,k]}(x)$. Prove that $f_k \to \chi_{[0,+\infty)}$ pointwise on $\mathbb{R}$. Moreover, compute $\displaystyle\lim_{k \to \infty}$ $\int_{\mathbb{R}} f_k\,d\mu$ and $\int_{\mathbb{R}} f\,d\mu$. What can you conclude from these integrals? Is it possible to apply {\rm Theorem \ref{517}}? Justify your answers.
\end{exercise}

\begin{exercise} \label{E527}
Let $(\mathbb{R},\mathcal{B}(\mathbb{R}),\lambda)$ be a measure space. Define the sequence of functions $f_k:\mathbb{R} \to \mathbb{R}$ by $f_k(x):=\frac{1}{k}\,\chi_{[k,+\infty)}(x)$. Prove that $f_k \to 0$ pointwise on $\mathbb{R}$. Moreover, compute $\displaystyle\lim_{k \to \infty}$ $\int_{\mathbb{R}} f_k\,d\mu$ and $\int_{\mathbb{R}} f\,d\mu$. What can you conclude from these integrals? Is it possible to apply {\rm Theorem \ref{517}}? Does this contradict {\rm Corollary \ref{519}}? Justify your answers.
\end{exercise}

\begin{exercise} \label{E528}
Let $(\mathbb{R},\mathcal{B}(\mathbb{R}),\lambda)$ be a measure space. Define the sequence of functions $f_k:\mathbb{R} \to \mathbb{R}$ by $f_k(x):=\frac{1}{k}\,\chi_{[0,k]}(x)$. Prove that $(f_k)$ converges pointwise to $f=0$ on $\mathbb{R}$. Moreover, compute $\displaystyle\lim_{k \to \infty}$ $\int_{\mathbb{R}} f_k\,d\mu$ and $\int_{\mathbb{R}} f\,d\mu$. What can you conclude from these integrals? Is it possible to apply {\rm Theorem \ref{517}}? Is it possible to apply {\rm Theorem \ref{525}}? Justify your answers.
\end{exercise}

\begin{exercise} \label{E529}
Let $(\mathbb{R},\mathcal{B}(\mathbb{R}),\lambda)$ be a measure space. Define the sequence of functions $f_k:\mathbb{R} \to \mathbb{R}$ by $f_k(x):=k\,\chi_{[1/k,2/k]}(x)$ and let $f=0$. Does the sequence $(f_k)$ converge pointwise to $f$ on $\mathbb{R}$? Is it possible to apply {\rm Theorem \ref{517}}? Is it possible to apply {\rm Theorem \ref{525}}? Justify your answers.
\end{exercise}

\begin{exercise} \label{E530}
Let $X=[0,+\infty)$, $\mathsf{S}=\mathcal{B}([0,+\infty))$ and let $\lambda$ be the Lebesgue measure on $\mathsf{S}$. Define the sequence of functions $f_k:\mathbb{R} \to \mathbb{R}$ by $f_k(x)=-\frac{1}{k}\,\chi_{[0,k]}(x)$. Prove that the sequence $(f_k)$ converges pointwise to $f=0$ on $[0,+\infty)$. Moreover, compute $\displaystyle\liminf_{k \to \infty}$ $\int_{\mathbb{R}} f_k\,d\mu$ and $\int_{\mathbb{R}} f\,d\mu$. What can you conclude from these integrals? Is it possible to apply {\rm Theorem \ref{525}}? Justify your answers.
\end{exercise}

\begin{exercise} \label{E531}
Assume Fatou's lemma and prove the monotone convergence theorem.
\end{exercise}

\begin{exercise} \label{E532}
Let $(X,\mathsf{S},\mu)$ be a measure space and let $(f_k)$ be a sequence in $\overline{\mathbb{M}}(X,\mathsf{S})$. Prove the following statements:
\begin{itemize}
\item[(a)] If there exists $g \in \mathbb{M}^{+}(X,\mathsf{S})$ with $\int_X g\,d\mu < +\infty$ such that $f_k \geq -g$ $\mu$ a.e. on $A \in \mathsf{S}$, then
$$
\int_{A} \liminf_{k \to \infty} f_k\,d\mu \leq \liminf_{k \to \infty} \int_{A} f_k\,d\mu.
$$
\item[(b)] If there exists $g \in \mathbb{M}^{+}(X,\mathsf{S})$ with $\int_{X} g\,d\mu < +\infty$ such that $f_k \leq -g$ $\mu$ a.e. on $A \in \mathsf{S}$, then
$$
 \limsup_{k \to \infty} \int_{A} f_k\,d\mu \leq \int_{A} \limsup_{k \to \infty} f_k\,d\mu .
$$
\end{itemize}

Can part (a) be applied to {\rm Exercise \ref{E530}}? Justify your answer.
\end{exercise}

\begin{exercise} \label{E533}
Let $(X,\mathsf{S},\mu)$ be a finite measure space. Let $(f_k)$ be a sequence of elements of $\mathbb{M}^{+}(X,\mathsf{S})$ such that $f_k \to f$ uniformly on $X$. Prove that
$\int_X f\,d\mu =\displaystyle\lim_{k \to \infty}$ $\int_X f_k\,d\mu$.

Is it possible to omit the hypothesis $\mu(X)<+\infty$? Justify your answer.
\end{exercise}

\begin{exercise} \label{E534}
Let $(X,\mathsf{S},\mu)$ be a measure space and let $(f_k)$ be a sequence of elements of $\mathbb{M}^{+}(X,\mathsf{S})$ such that $f_k \to f$ on $X$ with
$\displaystyle\lim_{k \to \infty}$ $\int_X f_k \,d\mu = \int_X f\,d\mu < +\infty$.

Prove that, for every $A \in \mathsf{S}$,
$$\displaystyle\lim_{k \to \infty} \int_{A} f_k \,d\mu = \int_{A} f\,d\mu .$$

Also prove that the previous conclusion may fail if
$\displaystyle\lim_{k \to \infty}$ $\int_X f_k \,d\mu = \int_X f\,d\mu < +\infty$
does not hold.
\end{exercise}

\begin{exercise} \label{E535}
Let $[0,1]$ be endowed with the Borel $\sigma$-algebra $\mathcal{B}([0,1])$ and let $\lambda$ be the Lebesgue measure on $[0,1]$. Let $n \in \mathbb{N}$ be arbitrary but fixed.
\begin{itemize}
    \item[(a)] Construct a nondecreasing sequence $(s_k)$ of elements of $\mathbb{S}^{+}([0,1],\mathcal{B}([0,1]))$ such that $s_k \to x^n$ pointwise on $[0,1]$.
    \item[(b)] Using part \textit{(a)} and the monotone convergence theorem, prove that
    $$
    \int_{[0,1]} x^n\,d\lambda = \frac{1}{n+1}.
    $$
 \end{itemize}
\end{exercise}

\begin{exercise} \label{E536}
Let $[0,1]$ be endowed with the Borel $\sigma$-algebra $\mathcal{B}([0,1])$ and let $\lambda$ be the Lebesgue measure on $[0,1]$. Define $\varphi:[0,1] \to \mathbb{R}$ by
$$
\varphi(x):=\left\{
\begin{array}{lcl}
2x  & & \quad \mbox{if}\quad 0 \leq x \leq \displaystyle\tfrac{1}{2},\\
2(1-x) & & \quad \mbox{if}\quad \tfrac{1}{2} \leq x \leq 1.
\end{array}
\right.
$$
\begin{itemize}
    \item[(a)] Construct a nondecreasing sequence $(s_k)$ of elements of $\mathbb{S}^{+}([0,1],\mathcal{B}([0,1]))$ such that $s_k \to \varphi$ pointwise on $[0,1]$. 
    \item[(b)] Using part \textit{(a)} and the monotone convergence theorem, prove that
    $$
    \int_{[0,1]} \varphi\,d\lambda = \frac{1}{2}.
    $$
 \end{itemize}
\end{exercise}

\begin{exercise} \label{E537}
Let $f:[a,b] \to [0,\infty)$ be a bounded function and let $h$ be its upper envelope. Prove the following statements:
\begin{itemize}
\item[(a)] $h$ is upper semicontinuous.
\item[(b)] $f(t_0)=h(t_0)$ if and only if $f$ is upper semicontinuous at $t_0 \in [a,b]$.
\end{itemize}
\end{exercise}

{\setlength{\parindent}{0pt}
\begin{exercise} \label{E538}
Let $a,b \in \mathbb{R}$ be such that $a < b$. Consider on $[a,b]$ the Borel $\sigma$-algebra $\mathcal{B}([a,b])$ and let $\lambda$ be the Lebesgue measure on $[a,b]$. 

Let $f:[a,b] \to [0,\infty)$ be bounded and let $(P_k)$ be a sequence of partitions of $[a,b]$ such that $\Vert P_k \Vert \to 0$. Prove the following statements:
\begin{itemize}
\item[(a)] $\displaystyle\lim_{k\,\to\,\infty}$ $\underline{\mathcal{S}}(f,P_{k}) =\sup_{k\,\in\,\mathbb{N}}\underline{\mathcal{S}}(f,P_{k}) =\underline{R}\int_{a}^{b} f(t)\,dt$.
\item[(b)] $\displaystyle\lim_{k\,\to\,\infty}$ $\overline{\mathcal{S}}(f,P_{k}) =\inf_{k\,\in\,\mathbb{N}}\overline{\mathcal{S}}(f,P_{k}) =\overline{R}\int_{a}^{b} f(t)\,dt$.
\end{itemize}
\end{exercise}

(Hint: For $\varepsilon >0$, fix a partition $P=\{t_0=a<\ldots<t_n=b \}$ such that
$\overline{R}\int_{a}^{b} f(t)\,dt < \overline{\mathcal{S}}(f,P) < \overline{R}\int_{a}^{b} f(t)\,dt + \varepsilon$.
Let $\delta:=\displaystyle\min_{ 0 \leq i \leq {n-1}}$ $(t_{i+1}-t_{i})$.
Choose $k \in \mathbb{N}$ sufficiently large such that $\Vert P_k \Vert <\delta$. Prove that
$\overline{\mathcal{S}}(f,P_k) < \overline{\mathcal{S}}(f,P) + nB\Vert P_k \Vert$
where $B=\sup\left\{ |f(t)|\,:\, a \leq t \leq b\right\}$.
Conclude that
$\displaystyle\limsup_{k \to \infty}$ $\overline{\mathcal{S}}(f,P_k) \leq \overline{R}\int_{a}^{b} f(t)\,dt + \varepsilon$).}

{\setlength{\parindent}{0pt}
\begin{exercise} \label{E539}
Prove with an example that a Riemann-integrable function is not necessarily Borel-measurable.
\end{exercise}

(Hint: Use the fact that there exists $D \subset \mathbb{R}$ which is Lebesgue-measurable but not Borel-measurable).}

\begin{exercise}[Improper Riemann Integral] \label{E540} \index{integral!improper Riemann}
Let $a,b \in \mathbb{R}$ be such that $-\infty \leq a < b \leq \infty$. Consider on $(a,b)$ the Borel $\sigma$-algebra $\mathcal{B}(a,b)$ and let $\lambda$ be the Lebesgue measure on $(a,b)$. 

Let $f:(a,b) \to \mathbb{R}$ be a nonnegative function. If $f$ is Riemann-integrable on every closed and bounded interval $[c,d] \subset (a,b)$ and $\displaystyle\lim_{c \to a,\,d\to b}$ $ R \int_{c}^{d} f(x)\,dx$ exists, then we define the \textbf{improper Riemann integral} on the interval $(a,b)$ by
$R\int_{a}^{b}f(x)\,dx:=\displaystyle\lim_{c\to a,\,d \to b}$ $ R \int_{c}^{d} f(x)\,dx$. \index{integral!improper Riemann}
\begin{itemize}
    \item[(a)] On $(0,1)$, prove that $R\int_{0}^{1} \frac{1}{x}\sin\left(\frac{1}{x} \right)\,dx$ exists as an improper Riemann integral.
    \item[(b)] On the other hand, prove that $\int_{(0,1)} \left|\frac{1}{x}\sin\left(\frac{1}{x} \right)\right|\,d\lambda \nless +\infty$.
    \item[(c)] Does this contradict {\rm Theorem \ref{E532}}? Justify your answer.
\end{itemize}
\end{exercise}
\chapter{Lebesgue Integrable Functions} \label{Capitulo6}
\markboth{{\scriptsize 6. LEBESGUE INTEGRABLE FUNCTIONS}}{ {\scriptsize 6. LEBESGUE INTEGRABLE FUNCTIONS}}

Let us now return to one of the problems with which we began the first part of these notes:

\begin{center}
\textit{Is the pointwise limit of a sequence of integrable functions integrable?}
\end{center}

In Chapter \ref{Capitulo1} we saw that, for the Riemann integral, the answer to this question is negative in general. This motivated the need to extend the concept of integration to a broader class of functions. In the previous chapter we defined the Lebesgue integral for nonnegative measurable functions and presented an example showing that this integral constitutes an extension and reformulation of the Riemann integral. It is therefore natural to ask whether this definition can be extended to arbitrary measurable functions. The first objective of this chapter will be to establish such a criterion.

Let $(X,\mathsf{S},\mu)$ be a measure space and let $f \in \overline{\mathbb{M}}(X,\mathsf{S})$. By Proposition \ref{317}, the functions $f^{+}$ and $f^{-}$, defined by $f = f^{+} - f^{-}$, belong to the set $\overline{\mathbb{M}}^{+}(X,\mathsf{S})$, and therefore their Lebesgue integrals are well defined. A first attempt to define the integral of $f$ would be the expression
\begin{equation} \label{F61}
\int_X f\,d\mu = \int_X f^{+}\,d\mu - \int_X f^{-}\,d\mu.
\end{equation}

However, expression \eqref{F61} presents a difficulty: if $\int_X f^{+}\,d\mu = +\infty$ and $\int_X f^{-}\,d\mu = +\infty$, then we obtain an expression of the form $(+\infty) - (+\infty)$, which is not defined in $\overline{\mathbb{R}}$. On the other hand, if exactly one of the two integrals is equal to $+\infty$, then the difference is well defined and hence expression \eqref{F61} makes sense.

Consequently, if $f \in \overline{\mathbb{M}}(X,\mathsf{S})$ satisfies $\int_X f^{+}\,d\mu = +\infty$ or $\int_X f^{-}\,d\mu = +\infty$, but not both simultaneously, its \textbf{Lebesgue integral on $X$ with respect to $\mu$} is defined as the possibly extended real number given by \eqref{F61}. In this chapter, we shall first study those measurable functions whose Lebesgue integral is a finite real number.

The second objective will be to answer the initial problem in the context of the Lebesgue integral. To this end, we shall state one of the fundamental theorems of integration theory: \textit{the dominated convergence theorem}. This result establishes conditions under which the pointwise limit of a sequence of integrable functions remains integrable and, moreover, when it is possible to interchange the limit and the integral sign. Among its many applications is differentiation under the integral sign, of which we shall present several examples.

Finally, we shall conclude the chapter with the study of the integrability of measurable functions with complex values.

\section{Integrable Functions}

\begin{definition} \label{61} \index{function!Lebesgue-integrable} \index{integral!Lebesgue!of a measurable function}
Let $(X,\mathsf{S},\mu)$ be a measure space and let $f \in \overline{\mathbb{M}}(X,\mathsf{S})$ be given. We say that \textbf{$f$ is (Lebesgue) integrable with respect to $\mu$} if and only if $\int_X f^{+}\,d\mu <+\infty$ and $\int_X f^{-}\,d\mu <+\infty$. In this case, we define the \textbf{(Lebesgue) integral of $f$ on $X$ with respect to $\mu$} by the expression
$$
\int_X f\,d\mu :=\int_X f^{+}\,d\mu - \int_X f^{-}\,d\mu.
$$
\end{definition}

We shall denote by $\overline{\mathcal{L}}(X,\mathsf{S},\mu)$ the set consisting of all functions in $\overline{\mathbb{M}}(X,\mathsf{S})$ that are integrable with respect to $\mu$. That is,
$$
\overline{\mathcal{L}}(X,\mathsf{S},\mu):=\left\{ f \in \overline{\mathbb{M}}(X,\mathsf{S})\,:\, f\,\,\,\mbox{is integrable with respect to}\,\,\mu \right\}.
$$

Let $A \in \mathsf{S}$ and $f \in \overline{\mathcal{L}}(X,\mathsf{S},\mu)$ be arbitrary. Since $A \subset X$, we have $\chi_{A} \leq \chi_{X}$ and therefore
$$
f^{+}\cdot\chi_{A} \leq f^{+}\quad \mbox{and}\quad f^{-}\cdot\chi_{A} \leq f^{-}.
$$

Applying the monotonicity of the integral for nonnegative functions, we obtain
$$
\int_{A} f^{+}\,d\mu \leq \int_X f^{+}\,d\mu < +\infty \quad \mbox{and}\quad \int_{A} f^{-}\,d\mu \leq \int_X f^{-}\,d\mu < +\infty,
$$
which shows that $(f\cdot \chi_{A}) \in \overline{\mathcal{L}}(X,\mathsf{S},\mu)$. Moreover, we write
$$
\int_{A} f\,d\mu =\int_{X} (f\cdot\chi_{A})\,d\mu=\int_{X} (f\cdot\chi_{A})^{+}\,d\mu-\int_{X} (f\cdot\chi_{A})^{-}\,d\mu = \int_{A} f^{+}\,d\mu - \int_{A} f^{-}\,d\mu.
$$

We shall denote by $\mathcal{L}(X,\mathsf{S},\mu)$ the subclass of integrable functions taking values in $\mathbb{R}$. That is,
$$
\mathcal{L}(X,\mathsf{S},\mu):=\left\{ f \in \mathbb{M}(X,\mathsf{S})\,:\,f\,\,\,\mbox{is integrable with respect to}\,\,\mu \, \right\}.
$$

Let us now examine some important properties..

\begin{theorem} \label{62}
Let $(X,\mathsf{S},\mu)$ be a measure space.
\begin{itemize}
\item[(a)] $f \in \overline{\mathcal{L}}(X,\mathsf{S},\mu)$ if and only if $|f| \in \overline{\mathcal{L}}(X,\mathsf{S},\mu)$. Moreover,
$$
\left|\int_X f\,d\mu\right| \leq \int_X |f|\,d\mu.
$$
\item[(b)] If $f \in \overline{\mathcal{L}}(X,\mathsf{S},\mu)$ and $g \in \overline{\mathbb{M}}(X,\mathsf{S})$ is such that $|g|\leq |f|$ on $X$, then $g \in \overline{\mathcal{L}}(X,\mathsf{S},\mu)$.
\end{itemize}
\end{theorem}

\begin{proof}
\textit{(a):} Suppose that $f \in \overline{\mathcal{L}}(X,\mathsf{S},\mu)$. By Proposition \ref{317}, it follows that $f^{+},f^{-} \in \overline{\mathbb{M}}^{+}(X,\mathsf{S})$ and satisfy $|f(x)|=f^{+}(x)+f^{-}(x)$ for every $x \in X$.

Thus, $|f|^{+}=|f|$ on $X$ and therefore $|f|^{-}=0$ on $X$. Consequently, $\int_X |f|^{+}\,d\mu <  +\infty$ and $\int_X |f|^{-}\,d\mu=0$, which means that $|f| \in \overline{\mathcal{L}}(X,\mathsf{S},\mu)$.

Conversely, suppose that $|f| \in \overline{\mathcal{L}}(X,\mathsf{S},\mu)$. Since $f^{+},f^{-}$ are nonnegative and satisfy $f^{+},f^{-} \leq |f|$ on $X$, applying the monotonicity of the integral for nonnegative functions, we obtain
$$
\int_X f^{+}\,d\mu,\,\int_X f^{-}\,d\mu \leq \int_X |f|\,d\mu <+\infty.
$$

Therefore, $f \in \overline{\mathcal{L}}(X,\mathsf{S},\mu)$.

On the other hand, if $f \in \overline{\mathcal{L}}(X,\mathsf{S},\mu)$, then
$$
\left| \int_X f\,d\mu \right|=\left| \int_X f^{+}\,d\mu - \int_X f^{-}\,d\mu \right| \leq \int_X f^{+}\,d\mu  + \int_X f^{-}\,d\mu = \int_X |f|\,d\mu,
$$
as claimed.

\textit{(b):} Let $f \in \overline{\mathcal{L}}(X,\mathsf{S},\mu)$ and let $g \in \overline{\mathbb{M}}(X,\mathsf{S})$ satisfy $|g|\leq |f|$ on $X$. Then $g^{+},g^{-} \leq |g| \leq |f|$ on $X$, and applying the monotonicity of the integral for nonnegative functions, we obtain $\int_X g^{+}\,d\mu <+\infty$ and $\int_X g^{-}\,d\mu < +\infty$, since the previous part ensures that $\int_X |f|\,d\mu < +\infty$. Consequently, $g \in  \overline{\mathcal{L}}(X,\mathsf{S},\mu)$.
\end{proof}

The equality case in part \textit{(a)} of Theorem \ref{62} is established in the following proposition.

\begin{proposition} \label{63}
Let $f \in \overline{\mathcal{L}}(X,\mathsf{S},\mu)$. Then,
$$
\left|\int_{X} f\,d\mu\right| = \int_{X} |f|\,d\mu
$$
if and only if
$$
\mu\left( \left\{ x \in X\,:\, f(x)>0\right\}\right)=0\quad \text{or}\quad \mu \left( \left\{ x \in X\,:\, f(x)<0\right\}\right)=0.
$$

That is,
$$
\left|\int_{X} f\,d\mu\right| = \int_{X} |f|\,d\mu
$$
if and only if $f$ is nonnegative $\mu$-a.e. or if $f$ is nonpositive $\mu$-a.e.
\end{proposition}

The proof is left as an exercise [Exercise \ref{E62}].

One drawback of Definition \ref{61} is that it requires knowing the positive and negative parts of the function. We shall now prove that, if there exists another representation of a function as the difference of nonnegative measurable functions with finite Lebesgue integral, then its integral coincides with the difference of the integrals of those functions.

\begin{proposition} \label{64}
Let $f \in \overline{\mathcal{L}}(X,\mathsf{S},\mu)$ be given. If $f=f_1-f_2$ $\mu$-a.e., where $f_{1},f_{2} \in \overline{\mathbb{M}}^{+}(X,\mathsf{S})$, $\int_{X} f_{1}\,d\mu <+\infty$, and $\int_{X} f_{2}\,d\mu <+\infty$, then
$$
\int_{A} f\,d\mu =\int_{A} f_{1}\,d\mu - \int_{A} f_{2}\,d\mu \quad \text{for every}\,\,A \in \mathsf{S}.
$$
\end{proposition}

\begin{proof}
Since $f=f^{+}-f^{-}$ on $X$ and $f=f_{1}-f_{2}$ $\mu$-a.e., it follows that $f^{+}+f_{2}=f_{1}+f^{-}$ $\mu$-a.e. Applying Corollary \ref{518}, we obtain
$$
\int_{A} f^{+}\,d\mu + \int_{A} f_{2}\,d\mu = \int_{A} f_{1}\,d\mu + \int_{A} f^{-} \,d\mu \quad \text{for every}\,\,A \in \mathsf{S}.
$$

Since all the integrals are finite, we have
$$
\int_{A} f\,d\mu=\int_{A} f^{+}\,d\mu - \int_{A} f^{-}\,d\mu = \int_{A} f_{1}\,d\mu - \int_{A} f_{2} \,d\mu \quad \text{for every}\,\,A \in \mathsf{S},
$$
as claimed.
\end{proof}

Let us now look at some examples.

\begin{example} \label{65}
If $(X,\mathsf{S},\mu)$ is a measure space and $A \in \mathsf{S}$, then
$$
\int_{X} \chi_{A}\,d\mu = \mu(A).
$$

Therefore, $\chi_{A} \in \mathcal{L}(X,\mathsf{S},\mu)$ if and only if $\mu(A)<+\infty$.

Now, let $c \in \mathbb{R} \smallsetminus\{0\}$ be given. Consider the function $f(x):=c$ for every $x \in X$. It is not difficult to see that, if $c>0$, then $f^{+}(x)=\max\{c,0\}=c$ and $f^{-}(x)=\max\{-c,0\}=0$. Consequently,
$$
\int_{X} f\,d\mu = \int_{X} f^{+}\,d\mu = c\mu(X).
$$

Now, if $c<0$, then $f^{+}(x)=\max\{c,0\}=0$ and $f^{-}(x)=\max\{-c,0\}=-c$, so that
$$
\int_{X} f\,d\mu = -\int_{X} f^{-}\,d\mu = c\mu(X).
$$

Therefore, $f \in \mathcal{L}(X,\mathsf{S},\mu)$ if and only if $\mu(X)<+\infty$.

If $f$ is the constant function equal to $0$, it is clear that $f\in \mathcal{L}(X,\mathsf{S},\mu)$ independently of the measure $\mu$ on $(X,\mathsf{S})$.
\end{example}

\begin{example} \label{66}
Consider the measure space given by $X=[0,+\infty)$, $\mathsf{S}=\mathcal{B}([0,+\infty))$, and $\lambda$ the Lebesgue measure. Define the function
$$
f(x):=\sum_{k=1}^{\infty} \frac{(-1)^k}{k}\chi_{[k-1,k)}(x).
$$

It is not difficult to see that
$$
f^{+}(x)=\sum_{k \in \wp} \frac{1}{k}\chi_{[k-1,k)}(x),
$$
where $\wp:=\{2,4,\ldots,\}$ is the set of even numbers in $\mathbb{N}$. The Beppo Levi theorem implies that
$$
\int_{[0,+\infty)} f^{+}\,d\lambda = \sum_{k\in \wp} \int_{[0,+\infty)}  \frac{1}{k}\chi_{[k-1,k)}\,d\lambda =\sum_{k \in \wp} \frac{1}{k} = +\infty,
$$
and therefore $f$ is not Lebesgue-integrable on $X$ with respect to $\lambda$.
\end{example}

From now on, whenever no confusion arises, we shall simply say that a function $f$ is integrable instead of saying that it is Lebesgue-integrable with respect to the measure $\mu$.

The following result establishes some well-known properties of the integral for functions in $\overline{\mathcal{L}}(X,\mathsf{S},\mu)$.

\begin{theorem}[Linearity of the Lebesgue integral] \label{67}
The Lebesgue integral is a linear operator, that is:

\begin{itemize}
    \item[(a)] If $f \in \overline{\mathcal{L}}(X,\mathsf{S},\mu)$ and $\gamma\in \mathbb{R}$, then $\gamma \,f \in \overline{\mathcal{L}}(X,\mathsf{S},\mu)$ and 
    $$
\int_{A}\gamma\,f\,d\mu = \gamma\int_{A}f\,d\mu  \quad \text{for every}\,\,A\in \mathsf{S}.
    $$
    \item[(b)] If $f,g \in \overline{\mathcal{L}}(X,\mathsf{S},\mu)$, then $f+g \in \overline{\mathcal{L}}(X,\mathsf{S},\mu)$ and 
    $$
\int_{A}(f+g)\,d\mu = \int_{A}f\,d\mu + \int_{A}g\,d\mu \quad \text{for every}\,\,A\in \mathsf{S}.
    $$
\end{itemize}
\end{theorem}

\begin{proof}
Let $f,g \in \overline{\mathcal{L}}(X,\mathsf{S},\mu)$ and let $\gamma \in \mathbb{R}$. We shall prove the result on the set $X$, since the case of an arbitrary subset $A \in \mathsf{S}$ follows easily from it.

\textit{(a):} If $\gamma =0$, the result is immediate. Thus, suppose that $\gamma \neq 0$. Observe that
$$
(\gamma f)^{+}=\left\{
\begin{array}{lcl}
\gamma f^{+} & &\mbox{if}\quad \gamma >0,\\
-\gamma f^{-} & &\mbox{if}\quad \gamma <0,
\end{array}
\right. \qquad\text{and}\qquad
(\gamma f)^{-}=\left\{
\begin{array}{lcl}
\gamma f^{-} & &\mbox{if}\quad \gamma >0,\\
-\gamma f^{+} & &\mbox{if}\quad \gamma <0.
\end{array}
\right.
$$

By Corollary \ref{518}, we have
$$
\begin{aligned}
&\int_{X} (\gamma f)^{+}\,d\mu=\left\{
\begin{array}{lcl}
\gamma \int_{X} f^{+}\,d\mu & &\mbox{if}\quad \gamma >0,\\
-\gamma \int_{X} f^{-}\,d\mu & &\mbox{if}\quad \gamma <0,
\end{array}
\right.\\
&\int_{X}(\gamma f)^{-}\,d\mu=\left\{
\begin{array}{lcl}
\gamma \int_{X}f^{-}\,d\mu & &\mbox{if}\quad \gamma >0,\\
-\gamma \int_{X}f^{+}\,d\mu & &\mbox{if}\quad \gamma <0.
\end{array}
\right.
\end{aligned}
$$
and therefore $\int_{X} (\gamma f)^{+}\,d\mu < +\infty$ and $\int_{X} (\gamma f)^{-} \,d\mu <+\infty$. That is, $\gamma\, f \in \overline{\mathcal{L}}(X,\mathsf{S},\mu)$.

On the other hand,
$$
\gamma f = (\gamma f)^{+} - (\gamma f)^{-}=\left\{
\begin{array}{lcl}
\gamma f^{+} -\gamma f^{-}=\gamma(f^{+}-f^{-})  & &\mbox{if}\quad \gamma >0,\\
-\gamma f^{-} + \gamma f^{+}=-\gamma(f^{-}-f^{+}) & &\mbox{if}\quad \gamma<0,
\end{array}
\right.
$$
and consequently,
$$
\int_{X} \gamma\,f\,d\mu = \int_{X}(\gamma\,f)^{+}\,d\mu - \int_{X}(\gamma\,f)^{-}\,d\mu =\gamma\left( \int_{X}f^{+}\,d\mu-\int_{X}f^{-}\,d\mu\right)=\gamma\int_{X}f\,d\mu.
$$

\textit{(b):} By the triangle inequality in $\mathbb{R}$, we have $|f+g|\leq |f|+|g|$ on $X$, and applying the monotonicity of the integral for nonnegative functions, we obtain
$$
\int_X |f+g|\,d\mu \leq \int_X |f|\,d\mu + \int_X|g|\,d\mu  <+\infty.
$$

Consequently, $f+g \in \overline{\mathcal{L}}(X,\mathsf{S},\mu)$.

Now observe that
$$
f+g=(f^{+}-f^{-})+(g^{+}-g^{-})=(f^{+}+g^{+})-(f^{-}+g^{-})\quad \text{on}\,X,
$$
where the sum is as in Theorem \ref{334}. Proposition \ref{64} and the previous part then ensure that
$$
\begin{aligned}
\int_{X} (f+g)\,d\mu &= \int_{X} (f^{+}+g^{+})\,d\mu- \int_{X}(f^{-}+g^{-})\,d\mu \\
&=\left(\int_{X} f^{+}\,d\mu - \int_{X} f^{-}\,d\mu \right) + \left( \int_{X} g^{+}\,d\mu - \int_{X} g^{-}\,d\mu \right)=\int_{X} f\,d\mu + \int_{X} g\,d\mu.
\end{aligned}
$$

This proves linearity.
\end{proof}

\begin{corollary} \label{68}
Let $f,g \in \overline{\mathcal{L}}(X,\mathsf{S},\mu)$.
\begin{itemize}
\item[(a)] $\int_{A \cup B} f\,d\mu = \int_{A} f\,d\mu + \int_{B} f\,d\mu$ for every $A,B \in \mathsf{S}$ such that $A \cap B=\varnothing$.
\item[(b)] $\int_{A} f\,d\mu \leq \int_{A} g\,d\mu$ for every $A \in \mathsf{S}$ if and only if $f \leq g$ $\mu$ a.e.
\item[(c)] $\int_{A} f\,d\mu = \int_{A} g\,d\mu$ for every $A \in \mathsf{S}$ if and only if $f=g$ $\mu$ a.e.
\end{itemize}
\end{corollary}

\begin{proof}
\textit{(a):} Let $A,B \in \mathsf{S}$ be disjoint. Since $\chi_{A \cup B}=\chi_{A} + \chi_{B}$ on $X$, Theorem \ref{67} ensures that
$$
\int_{A \cup B} f\,d\mu = \int_X (f \cdot \chi_{A} + f \cdot \chi_{B})\,d\mu = \int_X f\cdot\chi_{A}\,d\mu + \int_X f\cdot \chi_{B}\,d\mu = \int_{A} f\,d\mu + \int_{B} f\,d\mu.
$$

\textit{(b):} $\Rightarrow):$ By Theorem \ref{67}, we have $\int_{A}(g-f)\,d\mu \geq 0$ for every $A \in \mathsf{S}$.

For each $k \in \mathbb{N}$, define the set
$$
N_{k}:=\left\{ x \in X\,:\, f(x)-g(x) > \frac{1}{k} \right\}. 
$$

Then $(N_{k})$ is an increasing sequence of elements of $\mathsf{S}$ such that $N_{k} \to N$, where
$$
N:=\left\{ x \in X\,:\, g(x)<f(x)  \right\}.
$$

From the inequality $(f-g)\cdot\chi_{N_{k}} \geq \frac{1}{k} \cdot \chi_{N_{k}} \geq 0$ and the initial assumption, we obtain
$$
0 \leq \int_{N_k} (g-f)\,d\mu \leq -\frac{1}{k}\,\mu(N_k),\quad \text{for every } k \in \mathbb{N}.
$$

Therefore $\mu(N_k)=0$ for every $k \in \mathbb{N}$ and consequently $\mu(N)=0$. That is, $f \leq g$ $\mu$ a.e.

$\Leftarrow):$ Let $A \in \mathsf{S}$ be arbitrary and let $N \in \mathcal{N}(\mu)$ be such that $f(x) \leq g(x)$ for every $x \in X\smallsetminus N$. 

Since $A=(A \cap N) \cup (A \smallsetminus N)$ and $\mu(N)=0$, then $\mu(A \cap N)=0$ and
$$
\int_{A \cap N} f\,d\mu =0 =\int_{A\cap N} g\,d\mu.
$$

Using part \textit{(a)} and Theorem \ref{65}, we conclude that
$$
\begin{aligned}
\int_{A} g\,d\mu &=\int_{A \smallsetminus N} g\,d\mu=\int_{A \smallsetminus N} f\,d\mu + \int_{A \smallsetminus N} (g-f)\,d\mu \geq \int_{A \smallsetminus N}f\,d\mu=\int_{A}f\,d\mu
\end{aligned}
$$
since $g-f \geq 0$ on $A \smallsetminus N$.

\textit{(c):} Observe that $f=g$ $\mu$ a.e. is equivalent to $f \leq g$ $\mu$ a.e. and to $g \leq f$ $\mu$ a.e. (see Example \ref{428}). The result follows directly by applying part \textit{(b)} to these inequalities.
\end{proof}

An important particular case of part \textit{(c)} of the previous result is the following:
$$
\int_{A} f\,d\mu =0 \quad \text{for every}\,\, A \in \mathsf{S}\quad \text{if and only if}\quad f=0 \quad \mu\,\,\text{a.e.}
$$

We leave it as a simple exercise to prove this result without using Corollary \ref{66} [Exercise \ref{E63}].

In Definition \ref{61}, the criterion for integrability was established by requiring that the positive and negative parts of a measurable function have finite Lebesgue integral. On the other hand, Theorem \ref{62} guarantees that a function $f$ is integrable with respect to a measure if and only if $|f|$ is. The usefulness of this equivalence lies in the fact that it allows us to directly apply many of the results developed in the previous chapter and, moreover, it coincides with the original formulation proposed by H. Lebesgue.

Let $f \in \overline{\mathbb{M}}(X,\mathsf{S})$. Exercise \ref{E518} ensures that, if $f$ is integrable with respect to $\mu$, then
$$
\mu \left( \left\{ x \in X\,:\, |f(x)|=+\infty\right\} \right)=0.
$$

That is, $|f|<+\infty$ $\mu$ a.e. This statement is important for the following result.

\begin{proposition} \label{69}
Let $(X,\mathsf{S},\mu)$ be a measure space and let $f \in \overline{\mathbb{M}}(X,\mathsf{S})$ be given. Then, $f \in \overline{\mathcal{L}}(X,\mathsf{S},\mu)$ if and only if there exists a function $f_{0} \in \mathcal{L}(X,\mathsf{S},\mu)$ such that $f=f_{0}$ $\mu$ a.e.
\end{proposition}

\begin{proof}
$\Leftarrow):$ If there exists $f_0 \in\mathcal{L}(X,\mathsf{S},\mu)$ such that $f_0=f$ $\mu$ a.e., then $\int_{X} f_0\,d\mu =\int_X f\,d\mu$ [Exercise \ref{E64}] and, since $\int_X f_0\,d\mu$ is finite, then $f \in \overline{\mathcal{L}}(X,\mathsf{S},\mu)$.

$\Rightarrow):$ Consider the measurable set $A=\left\{ x \in X\,:\, |f(x)|=+\infty\right\},$ whose measure is equal to $0$ [Exercise \ref{E518}]. Defining $f_0: X \to \mathbb{R}$ by $f_0:=f\cdot \chi_{X\smallsetminus A}$, it is clear that $f_0 \in \mathbb{M}(X,\mathsf{S})$ and $f_0 = f$ $\mu$ a.e. Therefore, $\int_X f_0\,d\mu = \int_X f\,d\mu$ [Exercise \ref{E64}], which proves that $f_0 \in \mathcal{L}(X,\mathsf{S},\mu)$.
\end{proof}

\section{Lebesgue dominated convergence theorem}

Let us consider the following examples.

\begin{example} \label{610}
Let $(\mathbb{R},\mathcal{B}(\mathbb{R}),\lambda)$ be a measure space. Define the function $f_k:\mathbb{R} \to \mathbb{R}$ by
$$
f_k:=\chi_{(-k,k)},
$$
that is, the characteristic function of the open ball centered at $0$ with radius $k \in \mathbb{N}$ in $\mathbb{R}$.

\begin{figure}[ht!]
\centering
\begin{tikzpicture}[xscale=0.8,yscale=0.8]
	\draw[->,gray] (0,0) -- (5,0); \draw [->,gray] (0,0) -- (0,4); 
	\draw[<-,gray] (-5,0)--(0,0); \draw[-,gray] (0,0)--(0,-0.23);
\draw (5,0) node[right] {$\mathbb{R}$}; \draw (0,4) node[right] {$\mathbb{R}$}; 

\draw (0,3) node{$_{-}$}; \draw (0,3.2) node[left]{$_{1}$};
\draw (3,0) node{$_{|}$}; \draw (3,-0.1) node[below]{$_{k}$};
\draw (-3,0) node{$_{|}$}; \draw (-3,-0.1) node[below]{$_{-k}$};
\draw[ultra thick] (-2.91,3)--(2.91,3); \draw (3,3) node{$_{\circ}$};
\draw[ultra thick] (5,0)--(3,0); \draw (-3,3) node{$_{\circ}$};
\draw[ultra thick] (-5,0)--(-3,0); 
\draw (0,-0.1) node[below]{$_{0}$};
\draw [dotted] (3,0)--(3,3);
\draw [dotted] (-3,0)--(-3,3);
\end{tikzpicture}
\begin{center}
$f_k(x)=\chi_{(-k,k)}(x)$
\end{center}
\end{figure}

Then $(f_k)$ is a nondecreasing sequence of nonnegative Borel-measurable functions that converges pointwise to the function $f:=\chi_{\mathbb{R}}$ on $\mathbb{R}$. Moreover, for every $k \in \mathbb{N}$, we have
$$
\int_{\mathbb{R}} f_k \,d\lambda = \lambda((-k,k))=2k,
$$
which means that $f_k \in \mathcal{L}(\mathbb{R},\mathcal{B}(\mathbb{R}),\lambda)$ for every $k \in \mathbb{N}$. Furthermore,
$$
\lim_{k \to \infty} \int_{\mathbb{R}} f_k\,d\lambda = +\infty.
$$

On the other hand,
$$
\int_{\mathbb{R}}\lim_{k \to \infty} f_k\,d\lambda = \lambda(\mathbb{R})=+\infty,
$$
and we conclude that
$$
\lim_{k \to \infty} \int_{\mathbb{R}} f_k\,d\lambda =\int_{\mathbb{R}} \lim_{k \to \infty} f_k\,d\lambda.
$$
\end{example}

The previous example exhibits a sequence of integrable functions $(f_k)$ in $(\mathbb{R},\mathcal{B}(\mathbb{R}),\lambda)$ whose pointwise limit is not integrable.

\begin{example} \label{611}
Let $(\mathbb{R},\mathcal{B}(\mathbb{R}),\lambda)$ be a measure space. Define the function $f_k:\mathbb{R} \to \mathbb{R}$ by
$$
f_k:=\chi_{[\frac{1}{k},1)}.
$$

\begin{figure}[ht!]
\begin{minipage}[c]{0.5\textwidth}
\begin{center}
\begin{tikzpicture}[xscale=0.75,yscale=0.75]
	\draw[->,gray] (0,0) -- (4,0); \draw [->,gray] (0,0) -- (0,4); 
	\draw[<-,gray] (-1,0)--(0,0); \draw[-,gray] (0,0)--(0,-0.23);

\draw (0,0) node{$_{\bullet}$};
\draw (0,3) node{$_{-}$}; \draw (0,3) node[left]{$_{1}$};
\draw (1.5,0) node{$_{|}$}; \draw (1.5,-0.1) node[below]{$_{\frac{1}{k}}$};
\draw (3,0) node{$_{|}$}; \draw (3,-0.1) node[below]{$_{1}$};
\draw[ultra thick] (1.5,3)--(2.92,3); \draw (1.5,3) node{$_{\bullet}$}; \draw (3,3) node{$_{\circ}$};
\draw[ultra thick] (0,0)--(1.5,0); 
\draw (0,-0.1) node[below]{$_{0}$};
\draw [dotted] (1.5,0)--(1.5,3);
\draw[ultra thick] (0,0)--(-1,0);
\draw[ultra thick] (3,0)--(4,0);
\draw[dotted] (3,3)--(3,0);
\end{tikzpicture}   
\begin{center}
$f_{k}(x)=\chi_{[\frac{1}{k},1)}(x)$
\end{center}
\end{center}    
\end{minipage} \hfill \begin{minipage}[c]{0.5\textwidth}
\begin{center}
\begin{tikzpicture}[xscale=0.75,yscale=0.75]
	\draw[->,gray] (0,0) -- (4,0); \draw [->,gray] (0,0) -- (0,4); 
	\draw[<-,gray] (-1,0)--(0,0); \draw[-,gray] (0,0)--(0,-0.23);

\draw (0,0) node{$_{\bullet}$};
\draw (0,3) node{$_{-}$}; \draw (0,3) node[left]{$_{1}$};
\draw (1,0) node{$_{|}$}; \draw (1,-0.1) node[below]{$_{\frac{1}{k+1}}$};
\draw (3,0) node{$_{|}$}; \draw (3,-0.1) node[below]{$_{1}$};
\draw[ultra thick] (1,3)--(2.92,3); \draw (1,3) node{$_{\bullet}$}; \draw (3,3) node{$_{\circ}$};
\draw[ultra thick] (0,0)--(1,0); 
\draw (0,-0.1) node[below]{$_{0}$};
\draw [dotted] (1,0)--(1,3);
\draw[ultra thick] (0,0)--(-1,0); 
\draw[ultra thick] (3,0)--(4,0);
\draw[dotted] (3,3)--(3,0);
\end{tikzpicture}
\begin{center}
$f_{k+1}(x)=\chi_{[\frac{1}{k+1},1)}(x)$
\end{center}
\end{center}   
\end{minipage}
\end{figure}

Then $(f_k)$ is a sequence of Borel-measurable functions whose pointwise limit is given by
$$
\lim_{k \to \infty} f_k(x)=\left\{ 
\begin{array}{lcl}
1 & &\mbox{if}\quad x \in (0,1),\\
0 & &\mbox{if}\quad x \not \in (0,1).
\end{array}
\right.
$$

Observe that $f_k \in \mathcal{L}(\mathbb{R},\mathcal{B}(\mathbb{R}),\lambda)$ for every $k \in \mathbb{N}$ since
$$
\int_{\mathbb{R}} f_k\,d\lambda= \int_{\mathbb{R}} \chi_{[\frac{1}{k},1)}\,d\lambda= 1-\frac{1}{k}.
$$

Consequently,
$$
\lim_{k \to \infty} \int_{\mathbb{R}} f_k\,d\lambda = 1.
$$

On the other hand,
$$
\int_{\mathbb{R}} \lim_{k \to \infty}\,f_k\,d\lambda = \int_{\mathbb{R}} \chi_{(0,1)}\,d\lambda = 1
$$
and we conclude that
$$
\lim_{k \to \infty}\int_{\mathbb{R}} f_{k}\,d\lambda = \int_{\mathbb{R}} \lim_{k \to\infty} f_{k}\,d\lambda.
$$
\end{example}

In this example we observe that both the functions $f_{k}$ and their limit $\displaystyle\lim_{k \to \infty} f_k$ are integrable in the measure space $(\mathbb{R},\mathcal{B}(\mathbb{R}),\lambda)$ and, moreover, that it is possible to interchange the limit and the integral sign.

\begin{example} \label{612}
Let $(\mathbb{R},\mathcal{B}(\mathbb{R}),\lambda)$ be a measure space. Define the function $f_k:\mathbb{R} \to \mathbb{R}$ by
$$
f_k:=k\chi_{[0,\frac{2}{k}]}.
$$

\begin{figure}[ht!]
\begin{minipage}[c]{0.5\textwidth}
\begin{center}
\begin{tikzpicture}[xscale=0.8,yscale=0.8]
	\draw[->,gray] (0,0) -- (4,0); \draw [->,gray] (0,0) -- (0,4); 
	\draw[<-,gray] (-2.2,0)--(0,0); \draw[-,gray] (0,0)--(0,-0.23);
\draw (4,0) node[right] {$\mathbb{R}$}; \draw (0,4) node[right] {$\mathbb{R}$}; 

\draw (0,3) node{$_{-}$}; \draw (0,3) node[left]{$_{k}$};
\draw (1.5,0) node{$_{|}$}; \draw (1.5,-0.1) node[below]{$_{\frac{2}{k}}$};
\draw (3,0) node{$_{|}$}; \draw (3,-0.1) node[below]{$_{2}$};
\draw[ultra thick] (0,3)--(1.47,3); \draw (0,3) node{$_{\bullet}$}; \draw (1.5,3) node{$_{\circ}$};
\draw[ultra thick] (0,0)--(0,0); 
\draw[ultra thick] (3,0)--(4,0); 
\draw[ultra thick] (-2.2,0)--(0,0); 
\draw[ultra thick] (1.5,0)--(3,0); 
\draw (0,-0.1) node[below]{$_{0}$};
\draw [dotted] (1.5,0)--(1.5,3);
\end{tikzpicture}
\begin{center}
$f_{k}(x)=k\chi_{[0,\frac{2}{k}]}(x)$
\end{center}
\end{center}
\end{minipage} \hfill \begin{minipage}[c]{0.5\textwidth}
\begin{center} 
\begin{tikzpicture}[xscale=0.8,yscale=0.8]
	\draw[->,gray] (0,0) -- (4,0); \draw [->,gray] (0,0) -- (0,4); 
	\draw[<-,gray] (-2.2,0)--(0,0); \draw[-,gray] (0,0)--(0,-0.23);
\draw (4,0) node[right] {$\mathbb{R}$}; \draw (0,4) node[right] {$\mathbb{R}$}; 

\draw (0,3.5) node{$_{-}$}; \draw (0,3.5) node[left]{$_{k+1}$};
\draw (1,0) node{$_{|}$}; \draw (1,-0.1) node[below]{$_{\frac{2}{k+1}}$};
\draw (3,0) node{$_{|}$}; \draw (3,-0.1) node[below]{$_{2}$};
\draw[ultra thick] (0,3.5)--(1,3.5); \draw (0,3.5) node{$_{\bullet}$}; \draw (1,3.5) node{$_{\circ}$};
\draw[ultra thick] (0,0)--(0,0); 
\draw[ultra thick] (3,0)--(4,0); 
\draw[ultra thick] (-2.2,0)--(0,0); 
\draw[ultra thick] (1,0)--(3.5,0); 
\draw (0,-0.1) node[below]{$_{0}$};
\draw [dotted] (1,0)--(1,3.5);
\end{tikzpicture}
\begin{center}
$f_{k+1}(x)=(k+1)\chi_{[0,\frac{2}{k+1}]}(x)$
\end{center}
\end{center}
\end{minipage}
\end{figure}

It is clear that the sequence $(f_k)$ is formed by nonnegative Borel-measurable functions. In this context, it is natural to think that, when $k \to \infty$, the sequence $(f_k)$ approaches the constant zero function on $\mathbb{R}$. Let us now analyze whether this intuition is correct.

Let $\varepsilon > 0$ and $x \in \mathbb{R}$. Consider the following cases:

{\rm{\scshape Case 1:}}\quad If $x <0$, then clearly $x<\frac{2}{k}$ for every $k \in \mathbb{N}$, so that $f_k(x)=0$ for every $k \in \mathbb{N}$ and, therefore, $|f_k(x)|<\varepsilon$ for every $k \in \mathbb{N}$.

{\rm{\scshape Case 2:}}\quad If $0< x \leq 2$, then there exists $k_{0}=k_{0}(x) \in \mathbb{N}$ such that $\frac{2}{k_0}<x$. Consequently, $f_k(x)=0$ for every $k \geq k_{0}$ and, therefore, $|f_k(x)|<\varepsilon$ for every $k \geq k_{0}$.

{\rm{\scshape Case 3:}}\quad If $2 <x$, then clearly $f_k(x)=0$ for every $k \in \mathbb{N}$ and, therefore, $|f_k(x)|<\varepsilon$ for every $k \in \mathbb{N}$.

{\rm{\scshape Case 4:}}\quad Finally, if $x=0$, then $x \in [0,\frac{2}{k}]$ for every $k \in \mathbb{N}$ and we have $f_k(x)=k$ for every $k \in \mathbb{N}$. In this case, if $\varepsilon \in (0,1)$, then $|f_k(0)|=k> \varepsilon$ for every $k \in \mathbb{N}$. That is, $(f_k(0))$ does not converge to $0$ in $\mathbb{R}$.

From the four previous cases, we conclude that $f_{k}(x) \to f(x)$ for every $x \in \mathbb{R}\smallsetminus\{0\}$, where $f$ denotes the constant zero function on $\mathbb{R}$. Observe moreover that $\lambda(\{0\})=0$, and therefore we may conclude that $f_k \to f$ $\lambda$ a.e.

On the other hand, $f_k \in \mathcal{L}(\mathbb{R},\mathcal{B}(\mathbb{R}),\lambda)$ since
$$
\int_{\mathbb{R}} f_k\,d\lambda =k \lambda([0,\tfrac{2}{k}])=2,
$$
for every $k \in \mathbb{N}$. Consequently,
$$
\lim_{k \to \infty} \int_{\mathbb{R}} f_k\,d\lambda = 2.
$$

We also have $f \in\mathcal{L}(\mathbb{R},\mathcal{B}(\mathbb{R}),\lambda)$ and $\int_{\mathbb{R}} f\,d\lambda =0.$
\end{example}

In this example we have a sequence $(f_k)$ of integrable functions in $(\mathbb{R},\mathcal{B}(\mathbb{R}),\lambda)$ that converges a.e. to a function $f$ which is integrable, but for which the relation
$$
\lim_{k \to \infty} \int_{X} f_k\,d\mu = \int_{X} \lim_{k \to \infty} f_k\,d\mu
$$
does not hold.

From the previous examples, we are interested in determining under which conditions the limit of a sequence of functions $(f_k)$ in $\overline{\mathcal{L}}(X,\mathsf{S},\mu)$, defined on a measure space $(X,\mathsf{S},\mu)$, again belongs to $\overline{\mathcal{L}}(X,\mathsf{S},\mu)$ and, moreover, satisfies the equality
$$
\lim_{k \to \infty} \int_{X} f_k\,d\mu = \int_{X} \lim_{k \to \infty} f_k\,d\mu.
$$

This formulation essentially captures one of the central problems with which we began this text (see Problem \ref{17}).

The conditions that answer the previous question are established by one of the most important theorems in the modern theory of integration. Lebesgue's dominated convergence theorem is fundamental because of the mildness of its hypotheses, whose conclusions guarantee the integrability of the pointwise limit of a sequence of integrable functions, as well as the possibility of interchanging the limit and the integral sign. We now state a version of this theorem and propose as an exercise a generalization of it [Exercise \ref{E616}]. Both results are due to H. Lebesgue.

\begin{theorem}[Lebesgue dominated convergence theorem] \label{613} \index{theorem!Lebesgue dominated convergence}
Let $(X,\mathsf{S},\mu)$ be a measure space and let $(f_k)$ be a sequence of functions in $\overline{\mathcal{L}}(X,\mathsf{S},\mu)$ such that $f_k \to f$ $\mu$ a.e. for some $f \in \overline{\mathbb{M}}(X,\mathsf{S})$. Suppose that there exists $g \in \overline{\mathcal{L}}(X,\mathsf{S},\mu)$ such that, for every $k \in \mathbb{N}$,
$$
|f_{k}|\leq g \quad \mu \text{ a.e.}
$$
Then, $f \in \overline{\mathcal{L}}(X,\mathsf{S},\mu)$ and
$$
\int_{A} f\,d\mu =\lim_{k \to \infty} \int_{A} f_k\,d\mu\quad \text{for every}\,\, A\in\mathsf{S}.
$$
\end{theorem}

\begin{proof}
We shall prove the result only for the case $A=X$, since the general case is easily obtained from this one by considering the sequence of functions $(f_{k}\cdot \chi_{A})$ and the function $f\cdot \chi_{A}$ for every $A \in \mathsf{S}$.

The sets
$$
\begin{aligned}
N_{-1}&:=\left\{ x \in X\,:\, g(x)=\pm \infty\right\},\\
N_{0}&:=\left\{ x \in X\,:\, f_{k}(x)\,\,\mbox{does not converge to}\,\,f(x)\right\},\\
N_{k}&:=\left\{ x \in X\,:\,|f_{k}(x)|>g(x)\right\},
\end{aligned}
$$
are $\mu$-null. Therefore, $N:=\bigcup_{k=-1}^{\infty} N_k$ is $\mu$-null. Observe that $|f(x)|\leq g(x)<+\infty$ for every $x \in X\smallsetminus N$. Replacing $f_k(x)$, $f(x)$, and $g(x)$ by $0$ for every $x \in N$, we may assume, without loss of generality, that
$$
\lim_{k \to \infty} f_k(x)=f(x)\quad \mbox{and}\quad |f_k(x)|\leq g(x)\quad \forall k \in \mathbb{N},\quad\forall x \in X.
$$ 

Since $|f_{k}(x)|\leq g(x)$ for every $k \in \mathbb{N}$ and every $x \in X$, then $|f(x)| \leq g(x)$ for every $x \in X$ and, therefore, $f \in \overline{\mathcal{L}}(X,\mathsf{S},\mu)$.

By hypothesis, $- g \leq f_k \leq g$ on $X$ for every $k \in \mathbb{N}$. Fatou's lemma applied to the sequence of functions $(f_k+g)$ in $\overline{\mathbb{M}}^{+}(X,\mathsf{S})$ yields
$$
\begin{aligned}
\int_{X} (f+g)\,d\mu =\int_{X} \liminf_{k \to \infty}  (f_k+g)\,d\mu \leq \liminf_{k \to \infty} \int_{X} (f_k+g)\,d\mu = \liminf_{k \to \infty}\int_{X} f_k\,d\mu + \int_{X} g\,d\mu.
\end{aligned}
$$

Since $\int_{X} g\,d\mu <+\infty$, the previous inequality implies
\begin{eqnarray} \label{F62}
\int_{X} f\,d\mu \leq \liminf_{k \to \infty} \int_{X} f_k\,d\mu.
\end{eqnarray}

Applying Fatou's lemma once again to the sequence $(g-f_k)$ in $\overline{\mathbb{M}}^{+}(X,\mathsf{S})$, we obtain
$$
\begin{aligned}
\int_{X} (g-f)\,d\mu =\int_{X} \liminf_{k \to \infty} (g-f_k)\,d\mu \leq \liminf_{k \to \infty} \int_{X} (g-f_k)\,d\mu = \int_{X}g\,d\mu -\limsup_{k \to \infty} \int_{X} f_k\,d\mu.
\end{aligned}
$$
and therefore
\begin{eqnarray} \label{F63}
\limsup_{k \to \infty}\int_{X}f_k\,d\mu \leq \int_{X} f\,d\mu.
\end{eqnarray}

From \eqref{F62} and \eqref{F63} it follows that
$$
\int_{X} f\,d\mu \leq \liminf_{k \to \infty} \int_{X} f_k\,d\mu \leq \limsup_{k \to \infty} \int_{X} f_k\,d\mu \leq \int_{X} f\,d\mu.
$$

Consequently,
$$
\int_{X} f\,d\mu = \lim_{k \to \infty}\int_{X} f_k\,d\mu,
$$
as stated.
\end{proof}

In the Lebesgue dominated convergence theorem, the hypothesis of the existence of a function dominating the sequence is essential, even if the space has finite measure [Exercise \ref{E626}]. Observe that the sequences of functions in Examples \ref{610} and \ref{612} are not dominated by any integrable function $g$. We shall prove this fact only for Example \ref{612}; the other one is left as a simple exercise.

\begin{example} \label{614}
Let $(\mathbb{R},\mathcal{B}(\mathbb{R}),\lambda)$ be a measure space. Define the function $f_k:\mathbb{R} \to \mathbb{R}$ by $f_k:=k\chi_{[0,\frac{2}{k}]}$. For the sequence of functions $(f_k)$ there does not exist any integrable function $g$ that dominates the sequence.
\end{example}

\begin{proof}
Define the function $\xi:\mathbb{R} \to \overline{\mathbb{R}}$ by
$$
\xi(x):=\sup_{k \in \mathbb{N}} |f_k(x)|.
$$

If we prove that the function $\xi$ is not integrable on $\mathbb{R}$, then it follows immediately that any function $g$ satisfying $|f_k| \leq g$ on $X$ for every $k \in \mathbb{N}$ is also not integrable on $\mathbb{R}$, since $\xi \leq g$ on $\mathbb{R}$.

Let $j \in \mathbb{N}$. Define the function $s_j: \mathbb{R} \to \mathbb{R}$ by
$$
s_j:=\sum_{k=1}^{j} k\chi_{\left(\tfrac{2}{k+1},\tfrac{2}{k}\right]}.
$$

\begin{figure}[ht!]
\centering
\begin{tikzpicture}[xscale=1,yscale=0.9]
	\draw[->,gray] (0,0) -- (4.75,0); \draw [->,gray] (0,0) -- (0,4); 
	\draw[<-,gray] (-2.2,0)--(0,0); \draw[-,gray] (0,0)--(0,-0.23);
\draw (4.75,0) node[right] {$\mathbb{R}$}; \draw (0,4) node[right] {$\mathbb{R}$}; 

\draw (0,0.75) node{$_{-}$}; \draw (0,0.75) node[left]{$_{1}$};
 \draw (0,1.2) node[left]{$_{\vdots}$};
\draw (0,1.5) node{$_{-}$}; \draw (0,1.5) node[left]{$_{k}$};
\draw (0,2) node{$_{-}$}; \draw (0,2) node[left]{$_{k+1}$};
\draw (0,3) node{$_{-}$}; \draw (0,3) node[left]{$_{j}$}; \draw (0,2.65) node[left]{$_{\vdots}$};

\draw (1.25,2.75) node{$_{\ddots}$};
\draw (1.625,2.425) node{$_{\ddots}$};
\draw (3.2,1.25) node{$_{\ddots}$};

\draw (2.5,0) node{$_{|}$}; \draw (2.5,-0.1) node[below]{$_{\frac{2}{k}}$};
\draw (2,0) node{$_{|}$}; \draw (2,-0.1) node[below]{$_{\frac{2}{k+1}}$};
\draw (0.5,0) node{$_{|}$}; 
\draw (1.5,-0.1) node[below]{$_{\ldots}$};
\draw (1,0) node{$_{|}$}; \draw (1,-0.1) node[below]{$_{\frac{2}{j}}$};
\draw (4,0) node{$_{|}$}; \draw (4,-0.1) node[below]{$_{2}$};
\draw (3.5,0) node{$_{|}$}; \draw (3.5,-0.1) node[below]{$_{\frac{2}{2}}$};
\draw (3,-0.1) node[below]{$_{\ldots}$};

\draw[ultra thick] (0.5,3)--(1,3);
\draw[ultra thick] (2,2)--(2.5,2);
\draw[ultra thick] (3,1.5)--(2.5,1.5);
\draw[ultra thick] (4,0.75)--(3.5,0.75);

\draw[ultra thick] (4,0)--(4.5,0);
\draw[ultra thick] (-2,0)--(0.5,0);

\draw (0,-0.1) node[below]{$_{0}$};

\end{tikzpicture}
\begin{center}
$s_j(x)$
\end{center}
\end{figure}

Then $s_j=f_k$ on $\left(\frac{2}{k+1},\frac{2}{k} \right]$ for each $1 \leq k \leq j$, and therefore we may assert that $0 \leq s_j \leq \xi$ for every $j \in \mathbb{N}$. By the monotonicity of the integral, we have
\begin{eqnarray} \label{F64}
\int_{\mathbb{R}} s_j \,d\lambda \leq \int_{\mathbb{R}} \xi \,d\lambda \quad \text{for every}\,\, j \in \mathbb{N}
\end{eqnarray}
where
$$
\begin{aligned}
\int_{\mathbb{R}} s_j\,d\lambda &= \sum_{k=1}^{j} k\,\lambda\left( \left( \frac{2}{k+1},\frac{2}{k} \right] \right)=\sum_{k=1}^{j} k \left( \frac{2}{k}-\frac{2}{k+1} \right)=\sum_{k=1}^{j} \frac{2}{k+1}.
\end{aligned}
$$

Taking the limit as $j \to \infty$ in \eqref{F64}, we conclude that
$$
\int_{\mathbb{R}} \xi \,d\lambda = +\infty,
$$
that is, $\xi \not \in \mathcal{L}(\mathbb{R},\mathcal{B}(\mathbb{R}),\lambda)$.

Consequently, there does not exist any $g \in \mathcal{L}(\mathbb{R},\mathcal{B}(\mathbb{R}),\lambda)$ that dominates the sequence $(f_k)$.
\end{proof}

The previous example confirms why it was not possible to interchange the limit of the sequence $(f_k)$ from Example \ref{612} with the integral sign, even though the sequence converges pointwise almost everywhere.

We now establish some results of interest whose proofs follow directly as consequences of the dominated convergence theorem. Various applications of this principle will be presented in later sections.

\begin{corollary}[Bounded convergence theorem] \label{615} \index{theorem!bounded convergence}
Let $(X,\mathsf{S},\mu)$ be a finite measure space and let $(f_k)$ be a sequence of functions in $\mathbb{M}(X,\mathsf{S})$ such that $f_k \to f$ $\mu$ a.e. for some $f \in \mathbb{M}(X,\mathsf{S})$, and suppose that there exists $\gamma>0$ such that, for every $k \in \mathbb{N}$,
$|f_k| \leq \gamma$ $\mu$ a.e. Then $f_k,f \in \mathcal{L}(X,\mathsf{S},\mu)$ for every $k \in \mathbb{N}$ and
$$
\int_{A} f \,d\mu = \lim_{k \to \infty} \int_{A} f_k\,d\mu \quad \text{for every}\,\,A \in \mathsf{S}.
$$
\end{corollary}

\begin{proof}
Defining the function $g:=\gamma$, we have that $g \in \mathbb{M}^{+}(X,\mathsf{S})$ and, therefore,
$$
\int_X g \,d\mu = \int_ X {\gamma}\,d\mu = \gamma \,\mu(X) < +\infty.
$$

Consequently, $g \in \mathcal{L}(X,\mathsf{S},\mu)$. Since $|f_k| \leq g$ $\mu$ a.e. for every $k \in \mathbb{N}$, it follows that $(f_k)$ is a sequence in $\mathcal{L}(X,\mathsf{S},\mu)$. Applying Theorem \ref{613}, we obtain the result.
\end{proof}

\begin{corollary} \label{616}
Let $(X,\mathsf{S},\mu)$ be a measure space and let $(f_j)$ be a sequence of functions in $\mathcal{L}(X,\mathsf{S},\mu)$ such that
$$
\sum_{j=1}^{\infty} \int_X |f_j|\,d\mu < +\infty.
$$

Then the series $\sum_{j=1}^{\infty} f_k$ converges $\mu$ a.e. to a function $f \in \mathcal{L}(X,\mathsf{S},\mu)$,
$$
\sum_{j=1}^{\infty} \int_X f_j\,d\mu = \int_X f\,d\mu\quad \text{and}\quad  \lim_{k \to \infty}\int_X \left| f-\sum_{j=1}^{k} f_j\right|\,d\mu=0.
$$
\end{corollary}

\begin{proof}
Let
$$
g_k(x):=\sum_{j=1}^{k}|f_j(x)| \quad \text{and}\quad g(x):=\sum_{j=1}^{\infty}|f_j(x)|.
$$

It is clear that $g_k,g \in \overline{\mathbb{M}}^{+}(X,\mathsf{S})$ for every $k \in \mathbb{N}$. The Beppo Levi theorem (see Theorem \ref{520}) ensures that
$$
\int_{X} g\,d\mu=\lim_{k \to \infty} \int_X g_k\,d\mu =\sum_{j=1}^{\infty} \int_X |f_j|\,d\mu  < +\infty,
$$
that is, $g \in \mathcal{L}(X,\mathsf{S},\mu)$. There exists a $\mu$-null subset $N$ of $X$ such that (see Exercise \ref{E518})
$$
g(x):=\sum_{j=1}^{\infty} |f_j(x)| <+\infty \quad \forall x \in X\smallsetminus N.
$$

Thus, the series $\sum_{j=1}^{\infty} f_j(x)$ converges in $\mathbb{R}$ for every $x \in X\smallsetminus N$.

Define $f:X \to \mathbb{R}$ by
$$
f(x):=\left\{
\begin{array}{lcl}
\sum_{j=1}^{\infty} f_j(x) & & \mbox{if}\quad x \in X\smallsetminus N,\\
0 & &\mbox{if}\quad x \in N.\\
\end{array}
\right.
$$

Then $\sum_{j=1}^{k} f_j \to f$ $\mu$ a.e. and
$$
\left| \sum_{j=1}^{k} f_j(x) \right| \leq g(x) \quad \forall x \in X,\quad \forall k \in \mathbb{N}.
$$

The Lebesgue dominated convergence theorem ensures that $f \in \mathcal{L}(X,\mathsf{S},\mu)$ and that
$$
\int_X f\,d\mu = \sum_{j=1}^{\infty} \int_X f_j\,d\mu .
$$

On the other hand, since $f-\sum_{j=1}^{k} f_j \in \mathcal{L}(X,\mathsf{S},\mu)$ and
$$
\left| f(x)- \sum_{j=1}^{k} f_j(x)\right| \leq |f(x)| + \left| \sum_{j=1}^{\infty} f_j(x)\right| \leq |f(x)| + g(x) \quad \forall x \in X\smallsetminus N\quad \forall k \in \mathbb{N},
$$
the Lebesgue dominated convergence theorem again implies that
$$
\lim_{k \to \infty} \int_X \left|f-\sum_{j=1}^{k} f_j \right|\,d\mu=0,
$$
as stated.
\end{proof}

\begin{corollary} \label{617}
Let $(X,\mathsf{S},\mu)$ be a measure space and let $a,b$ be fixed such that $-\infty \leq a < b \leq \infty$. If $f:X \times (a,b) \to \mathbb{R}$ is a function such that $f^{t}:X \to \mathbb{R}$ given by $f^{t}(x):=f(x,t)$ is $\mathsf{S}$-measurable for every $t \in (a,b)$, and $f_{x}:(a,b) \to \mathbb{R}$ given by $f_{x}(t):=f(x,t)$ is continuous for every $x \in X$, and moreover there exists $g \in \mathcal{L}(X,\mathsf{S},\mu)$ such that $|f(x,t)| \leq g(x)$ for every $x \in X$ and every $t \in (a,b)$, then the function $F:(a,b) \to \mathbb{R}$ defined by
$$
F(t):=\int_X f^{t}\,d\mu,
$$
\noindent is continuous on $(a,b)$.
\end{corollary}

\begin{proof}
Let $t_0 \in (a,b)$ be arbitrary and let $(t_k)$ be a sequence of elements of $(a,b)$ such that $t_k \to t_0$. Define the sequence of functions $f_k: X \to \mathbb{R}$ by $f_k:=f^{t_k}$. By hypothesis, $(f_k)$ is a sequence of elements of $\mathbb{M}(X,\mathsf{S})$ and $|f_k(x)|:=|f^{t_k}(x)|=|f(x,t_k)| \leq g(x)$ for every $x \in X$. Thus, $f_k \in \mathcal{L}(X,\mathsf{S},\mu)$ for every $k \in \mathbb{N}$. Moreover, by continuity, $f(x,t_k) \to f(x,t_0)$ for every $x \in X$, that is, $f_k(x) \to f^{t_0}(x)$ for every $x \in X$.

The Lebesgue dominated convergence theorem ensures that
$$
F(t_0)=\int_X f^{t_0}\,d\mu = \lim_{k \to \infty}\int_X f_{k}\,d\mu = \lim_{k \to \infty} F(t_k).
$$

This proves that $F$ is continuous on $(a,b)$.
\end{proof}

\begin{corollary}[Differentiation under the integral sign] \label{618} 
Consider the hypotheses and notation as in the previous corollary. Suppose moreover that there exists $t_0 \in (a,b)$ such that $f^{t_0} \in \mathcal{L}(X,\mathsf{S},\mu)$, that $\frac{\partial f}{\partial t}$ exists on $X \times (a,b)$, and that there exists $g \in  \mathcal{L}(X,\mathsf{S},\mu)$ such that
$$
\left| \frac{\partial f}{\partial t} (x,t)\right| \leq g(x)\quad \forall x\in X,\ \forall t\in(a,b).
$$

Then $F$ is differentiable on $(a,b)$ and
$$
F'(t)=\int_X \frac{\partial f}{\partial t}(x,t)\,d\mu
$$
for every $t \in (a,b)$.
\end{corollary}

\begin{proof}
Let $t \in (a,b)$ be arbitrary such that $t \neq t_0$. Given $x \in X$, applying the mean value theorem, we obtain that there exists $\xi \in (a,b)$ such that
$$
|f(x,t)-f(x,t_0)|=|t-t_0|\left| \frac{ \partial f}{\partial t}(x,\xi)\right|.
$$ 

Consequently,
$$
|f(x,t)| \leq |f(x,t_0)| +(b-a)g(x).
$$

Hence, $f^{t} \in \mathcal{L}(X,\mathsf{S},\mu)$ for every $t \in (a,b)$. 

Now let $ t \in (a,b)$ be arbitrary and let $(t_k)$ be a sequence of elements of $(a,b)$ such that $t_k \neq t$ for every $k \in \mathbb{N}$ and $t_k \to t$ in $\mathbb{R}$. Define the sequence of functions $f_k:X \to \mathbb{R}$ by 
$$
f_k(x):=\frac{f^{t_k}(x)-f^{t}(x)}{t_k-t}.
$$

Then
$$
\lim_{k \to \infty} f_k(x)=\frac{\partial f}{\partial t}(x,t)
$$
for every $x \in X$ and, therefore,
$$
\frac{\partial f}{\partial t}(\cdot,t):X \to \mathbb{R}
$$
is an $\mathsf{S}$-measurable function. Moreover, applying an argument similar to the previous one, we have that for every $k \in \mathbb{N}$ there exists $\xi_k \in (a,b)$ such that
$$
|f_k(x)| =\left|\frac{f^{t_k}(x)-f^{t}(x)}{t_k-t} \right| \leq \left|\frac{ \partial f}{\partial t}(x,\xi_k)\right| \leq g(x) \quad \text{for every }\,x \in X.
$$

Applying the Lebesgue dominated convergence theorem, we obtain
$$
\lim_{k \to \infty} \frac{F(t_k)-F(t)}{t_k-t}=\lim_{k \to \infty} \int_X f_k(x)\,d\mu = \int_X \frac{\partial f}{\partial t}(x,t)\,d\mu.
$$

That is, $F$ is differentiable on $(a,b)$ and
$$
F'(t)=\int_X \frac{\partial f}{\partial t}(x,t)\,d\mu
$$
for every $t \in (a,b)$.
\end{proof}

In the following example we shall show an application of the previous result to compute improper Riemann integrals [Exercise \ref{E634}].

\begin{example} \label{619}
On $(0,1]$ consider the Borel $\sigma$-algebra and let $\lambda$ be the Lebesgue measure. Then, for every $t>0$, we have
$$
R\int_{0}^{1} x^{t}\log(x)\,dx = -\frac{1}{(t+1)^2}.
$$
\end{example}

\begin{proof}
Define $f:(0,1]\times (0,+\infty) \to \mathbb{R}$ by $f(x,t):=x^{t}$.

\begin{figure}[ht!]
\begin{minipage}[r]{0.2\textwidth}
\begin{center}
\begin{tikzpicture}[xscale=2, yscale=2]
\draw[domain= 0:1,thick] plot(\x,{{\x^(1/2)}} );

\draw[->,gray] (-0.1,0)--(1.15,0); 
\draw[->,gray] (0,-0.1)--(0,1); 

\end{tikzpicture}
\begin{center}
$f(x,\frac{1}{2})=x^{1/2}$
\end{center}
\end{center}
\end{minipage} \hfill \begin{minipage}[c]{0.2\textwidth}
\begin{center}
\begin{tikzpicture}[xscale=2, yscale=2]
\draw[domain= 0:1,thick] plot(\x,{{\x}} );

\draw[->,gray] (-0.1,0)--(1.15,0); 
\draw[->,gray] (0,-0.1)--(0,1); 

\end{tikzpicture}
\begin{center}
$f(x,1)=x$
\end{center}
\end{center}
\end{minipage} \hfill \begin{minipage}[l]{0.2\textwidth}
\begin{center}

\begin{tikzpicture}[xscale=2, yscale=2]
\draw[domain= 0:1,thick] plot(\x,{{\x^(2)}} );

\draw[->,gray] (-0.1,0)--(1.15,0); 
\draw[->,gray] (0,-0.1)--(0,1); 
\end{tikzpicture}
\begin{center}
$f(x,2)=x^{2}$
\end{center}
\end{center} 
\end{minipage}  
\end{figure}

Let $t >0$. It is clear that the function $x^t$ is continuous in $x$ on every subinterval $[c,d]$ of $(0,1]$ and, therefore, Riemann-integrable on every subinterval $[c,d]$ of $(0,1]$. Hence, $f(x,t)$ is Lebesgue-integrable with respect to $\lambda$ for every $x\in (0,1]$.

On the other hand,
$$
\left| \frac{\partial}{\partial t}f(x,t) \right| = \left| x^t\log(x) \right| \leq |\log(x)| 
$$
for every $0< x \leq 1$ and $t >0$. Since $|\log(x)|$ is continuous and nonnegative on every closed subinterval of $(0,1]$, it follows that $|\log(x)|$ is Riemann-integrable on every closed subinterval of $(0,1]$. Thus, $|\log(x)|$ has an improper Riemann integral on $(0,1]$ and, consequently, it is Lebesgue-integrable with respect to $\lambda$ on $(0,1]$.

Corollary \ref{619} ensures that the function $F:(0,+\infty) \to \mathbb{R}$
$$
F(t):=R\int_{0}^{1}  x^t\,dx 
$$
is differentiable and satisfies
$$
F'(t)=R\int_{0}^{1} \frac{\partial}{\partial t}f(x,t)\,dx = R\int_{0}^{1} x^t\log(x)\,dx.
$$

Since $F(t)=\frac{1}{t+1}$ for every $t >0$, then $F'(t)=-\frac{1}{(t+1)^2}$ for every $t >0$. Therefore,
$$
R\int_{0}^{1} x^t\log(x)\,dx=F'(t)=-\frac{1}{(t+1)^2},
$$
for every $t >0$, as claimed.
\end{proof}

We conclude this section with an example that turns out to be quite challenging using the Riemann integral. However, by using the properties of the Lebesgue integral, it can be obtained without major difficulty.

\begin{example} \label{620}
On $[0,1]$ consider the Borel $\sigma$-algebra and let $\lambda$ be the Lebesgue measure. The sequence of functions $f_k:[0,1] \to \mathbb{R}$, $f_k(x)=kx^2e^{-kx^2}$, satisfies
$$
\lim_{k\to\infty} R\int_{0}^{1} kx^2e^{-kx^2}\,dx = 0.
$$
\end{example}

\begin{figure}[ht!]
\begin{minipage}[c]{0.2\textwidth}
\begin{center}
\begin{tikzpicture}[xscale=2.5, yscale=3]
\draw[domain= 0:1,thick] plot(\x,{{4*\x^2*exp(-4*\x^2)}} );

\draw[->,gray] (0,0)--(1.15,0); 
\draw[->,gray] (0,0)--(0,0.55); 

\end{tikzpicture}
\begin{center}
$y=4x^2e^{-4x^2}$
\end{center}
\end{center}
\end{minipage} \hfill \begin{minipage}[c]{0.2\textwidth}
\begin{center}
\begin{tikzpicture}[xscale=2.5, yscale=3]
\draw[domain= 0:1,thick] plot(\x,{{6*\x^2*exp(-6*\x^2)}} );

\draw[->,gray] (0,0)--(1.15,0); 
\draw[->,gray] (0,0)--(0,0.55); 

\end{tikzpicture}
\begin{center}
$y=6x^2e^{-6x^2}$
\end{center}
\end{center}
\end{minipage} \hfill \begin{minipage}[c]{0.2\textwidth}
\begin{center}

\begin{tikzpicture}[xscale=2.5, yscale=3]
\draw[domain= 0:1,thick] plot(\x,{{8*\x^2*exp(-8*\x^2)}} );

\draw[->,gray] (0,0)--(1.15,0); 
\draw[->,gray] (0,0)--(0,0.55); 

\end{tikzpicture}
\begin{center}
$y=8x^2e^{-8x^2}$
\end{center}
\end{center} 
\end{minipage}  \hfill \begin{minipage}[c]{0.2\textwidth}
\begin{center}
\begin{tikzpicture}[xscale=2.5, yscale=3]
\draw[domain= 0:1,thick] plot(\x,{{10*\x^2*exp(-10*\x^2)}} );

\draw[->,gray] (0,0)--(1.15,0); 
\draw[->,gray] (0,0)--(0,0.55); 

\end{tikzpicture}
\begin{center}
$y=10x^2e^{-10x^2}$
\end{center}
\end{center}
\end{minipage}
\end{figure}

\begin{proof}
Since $f_k$ is continuous on $[0,1]$ for every $k \in \mathbb{N}$, it follows that $f_k$ is Borel-measurable and Riemann-integrable on $[0,1]$ for every $k \in \mathbb{N}$. Henri Lebesgue's theorem (see Theorem \ref{532}) ensures that
$$
\int_{[0,1]} f_k\,d\lambda = R\int_{0}^{1} f_k(x)\,dx < + \infty,\quad\forall\, k \in \mathbb{N}.
$$

Observe that $0 \leq f_k(x)  \leq e^{-1}$ for every $x \in[0,1]$ and every $k \in \mathbb{N}$. Moreover, $\lim_{ k \to \infty} f_k(x)=0$ for every $x \in [0,1]$ since $\lim_{k \to \infty} e^{-kx^2}=0$ for every $x \in [0,1]$. Applying the Lebesgue dominated convergence theorem, we obtain
$$
\lim_{k \to \infty} R\int_{0}^{1} kx^2e^{-kx^2}\,dx = R\int_{0}^{1} \lim_{k \to \infty} (kx^2e^{-kx^2})\,dx =0,
$$
as claimed.
\end{proof}

\section{Complex integrable functions}

In this final section we shall study the integrability criterion for a measurable function taking values in the set of complex numbers $\mathbb{C}$, which reduces to examining its real and imaginary parts through the results established in the previous sections.

Let $(X,\mathsf{S},\mu)$ be a measure space. We denote by $\mathbb{M}_{\mathbb{C}}(X,\mathsf{S})$ the space of functions $f:X \to \mathbb{C}$ that are $\mathsf{S}$-measurable. Recall that, by Theorem \ref{336}, a function $f:X \to \mathbb{C}$ is $\mathsf{S}$-measurable if and only if its real part $\Re (f): X \to \mathbb{R}$ and its imaginary part $\Im(f): X \to \mathbb{R}$ are $\mathsf{S}$-measurable functions.

\begin{definition} \label{621} \index{integral!Lebesgue!of a complex measurable function}
A function $f \in \mathbb{M}_{\mathbb{C}}(X,\mathsf{S})$ is said to be \textbf{(Lebesgue) integrable on $X$ with respect to the measure $\mu$} if and only if $\Re(f)$ and $\Im(f):X \to \mathbb{R}$ are integrable on $X$ with respect to $\mu$.

In this case, the \textbf{integral of $f$ on $X$} is defined by
$$
\int_{X} f\,d\mu:=\int_{X} \Re(f)\,d\mu +i \int_{X} \Im(f)\,d\mu .
$$
\end{definition}

We denote by $\mathcal{L}_{\mathbb{C}}(X,\mathsf{S},\mu)$ the space of complex-valued functions that are integrable with respect to $\mu$. That is,
$$
\mathcal{L}_{\mathbb{C}}(X,\mathsf{S},\mu):=\left\{ f \in \mathbb{M}_{\mathbb{C}}(X,\mathsf{S})\,:\, \Re(f), \Im(f) \in \mathcal{L}(X,\mathsf{S},\mu) \right\}.
$$

For every $f \in \mathcal{L}_{\mathbb{C}}(X,\mathsf{S},\mu)$ and every $A \in \mathsf{S}$ we write
$$
\int_{A} f\,d\mu = \int_{A} \Re(f)\,d\mu + i \int_{A} \Im(f)\,d\mu. 
$$

\begin{theorem}[Linearity] \label{622}
If $f,g \in \mathcal{L}_{\mathbb{C}}(X,\mathsf{S},\mu)$ and $\gamma \in \mathbb{C}$, then $\gamma f +g \in \mathcal{L}_{\mathbb{C}}(X,\mathsf{S},\mu)$. Moreover,
$$
\int_{X} (\gamma f +g)\,d\mu = \gamma\int_{X} f \,d\mu + \int_{X} g\,d\mu.
$$
\end{theorem}

\begin{proof}
If we write $\gamma=\gamma_{1}+i\gamma_{2}$ with $\gamma_{1},\gamma_{2} \in \mathbb{R}$, $f=\Re(f)+i\Im(f)$, and $g=\Re(g)+i\Im(g)$, then
$$
\begin{aligned}
\gamma f +g &= (\gamma_{1}+i\gamma_{2})(\Re(f)+i\Im(f))+(\Re(g)+i\Im(g)) \\
&=(\gamma_{1}\Re(f)-\gamma_{2}\Im(f)) + i (\gamma_{1}\Im(f)+\gamma_{2}\Re(f)) + (\Re(g)+i\Im(g))\\
&=(\gamma_{1}\Re(f)-\gamma_{2}\Im(f) + \Re(g)) + i (\gamma_{1}\Im(f)+\gamma_{2}\Re(f) + \Im(g)).
\end{aligned}
$$

That is, $\Re(\gamma f + g)=\gamma_{1}\Re(f)-\gamma_{2}\Im(f) + \Re(g)$ and $\Im(\gamma f + g)=\gamma_{1}\Im(f)+\gamma_{2}\Re(f) + \Im(g)$.

Theorem \ref{67} then ensures that $\Re(\gamma f + g), \Im(\gamma f + g) \in \mathcal{L}(X,\mathsf{S},\mu)$ and, therefore, $\gamma f + g \in \mathcal{L}_{\mathbb{C}}(X,\mathsf{S},\mu)$. Moreover,
$$
\begin{aligned}
\int_{X} (\gamma f + g)\,d\mu &= \int_{X} \Re(\gamma f + g)\,d\mu + i \int_{X}\Im(\gamma f +g)\,d\mu\\
&=\gamma_{1}\int_{X}\Re(f)\,d\mu -\gamma_{2}\int_{X}\Im(f)\,d\mu + \int_{X} \Re(g)\,d\mu\\
&\quad +i\gamma_{1}\int_{X} \Im(f)\,d\mu  + i\gamma_{2}\int_{X}\Re(f)\,d\mu + i \int_{X} \Im(g)\,d\mu\\
&=\gamma\int_{X}f\,d\mu + \int_{X}g\,d\mu,
\end{aligned}
$$
as claimed.
\end{proof}

Recall that, if $f:X \to \mathbb{C}$ is a function, then its complex modulus is defined as the function $|f|:X \to \mathbb{R}$ given by $|f(x)|:=\left(|\Re(f)(x)|^{2} + |\Im(f)(x)|^{2} \right)^{1/2}$.

\begin{theorem}\label{623}
\begin{itemize}
\item[(a)] $f \in \mathcal{L}_{\mathbb{C}}(X,\mathsf{S},\mu)$ if and only if $|f| \in \mathcal{L}(X,\mathsf{S},\mu)$. In this case,
$$
\left|\int_{X} f\,d\mu \right| \leq \int_{X} |f|\,d\mu.
$$
\item[(b)] If $f \in \mathcal{L}_{\mathbb{C}}(X,\mathsf{S},\mu)$ and $g \in \mathbb{M}_{\mathbb{C}}(X,\mathsf{S})$ satisfies $|g|\leq |f|$, then $g \in \mathcal{L}_{\mathbb{C}}(X,\mathsf{S},\mu)$.
\end{itemize}
\end{theorem}

\begin{proof}
\textit{(a):} $\Rightarrow):$ Since $f \in \mathcal{L}_{\mathbb{C}}(X,\mathsf{S},\mu)$, we have $\Re(f), \Im(f) \in \mathcal{L}(X,\mathsf{S},\mu)$. From the following chain of inequalities valid for every complex number
$$
|\Re(f)|,|\Im(f)| \leq |f| \leq |\Re(f)| + |\Im(f)|
$$
we obtain that $|f| \in \mathcal{L}(X,\mathsf{S},\mu)$ since $|\Re(f)| + |\Im(f)| \in  \mathcal{L}(X,\mathsf{S},\mu)$ by Theorem \ref{67}.

$\Leftarrow):$ Again, from the inequality
$$
|\Re(f)|,|\Im(f)| \leq |f| \leq |\Re(f)| + |\Im(f)|
$$
we obtain that $|\Re(f)|,|\Im(f)| \in \mathcal{L}(X,\mathsf{S},\mu)$ since $|f| \in \mathcal{L}(X,\mathsf{S},\mu)$. This proves that $\Re(f),\Im(f) \in \mathcal{L}(X,\mathsf{S},\mu)$ and, consequently, $f \in \mathcal{L}_{\mathbb{C}}(X,\mathsf{S},\mu)$.

On the other hand, suppose that $\int_{X} f\,d\mu \neq 0$. Let $re^{i\theta}=\int_{X} f\,d\mu$ be its polar representation. Then
$$
r=e^{-i\theta}\int_{X} f\,d\mu = \int_{X} e^{-i\theta}f\,d\mu.
$$

Thus,
$$
\left| \int_{X} f\,d\mu \right|=|re^{i\theta}|=|r||e^{i\theta}|=r = \int_{X} e^{-i\theta}f\,d\mu
$$
and, moreover,
$$
r=\Re(r)=\Re\left( \int_{X} e^{-i\theta}f\,d\mu \right)=\int_{X} \Re\left( e^{-i\theta}f\right)\,d\mu.
$$

Hence,
$$
\left| \int_{X} f\,d\mu \right|=\int_X \Re\left( e^{-i\theta}f\right)\,d\mu \leq \int_{X} |e^{-i\theta}f|\,d\mu=\int_{X} |f|\,d\mu,
$$
as claimed. The case $\int_X f\,d\mu =0$ follows trivially.

\textit{(b):} Since $|\Re(g)|,|\Im(g)| \leq |g| \leq |f|$ and $|f| \in \mathcal{L}(X,\mathsf{S},\mu)$ by the previous part, it follows that $\Re(g),\Im(g) \in \mathcal{L}(X,\mathsf{S},\mu)$ and, therefore, $g \in \mathcal{L}_{\mathbb{C}}(X,\mathsf{S},\mu)$.
\end{proof}

\begin{proposition} \label{624}
If $f\in \mathcal{L}_{\mathbb{C}}(X,\mathsf{S},\mu)$, then
$$
\left| \int_{X} f\,d\mu \right|=\int_{X} |f|\,d\mu \quad \Longleftrightarrow \quad \text{there exists}\,\,\theta \in \mathbb{R}\,\,\text{such that}\,\,\mu\left(\left\{ x\in X\,:\,f(x)\neq e^{i\theta}|f(x)| \right\} \right)=0.
$$
\end{proposition}

\begin{proof}
$\Leftarrow):$ Clearly,
$$
\left| \int_{X} f\,d\mu \right|=\left| \int_X e^{-i\theta}|f|\,d\mu \right|=\int_{X}|f|\,d\mu.
$$

$\Rightarrow):$ Suppose that $\int_{X} f\,d\mu \neq 0$. Write $\int_{X}f\,d\mu=re^{i\theta}$. As in the previous proof, we have
$$
\int_{X} \Re(e^{-i\theta}f)\,d\mu = \left| \int_X f\,d\mu \right|=\int_{X}|f|\,d\mu.
$$

Thus, $\Re(e^{-i\theta}f) = |f|$ $\mu$ a.e., that is,
$$
\mu\left(\left\{ x\in X\,:\,|f(x)|\neq \Re(e^{-i\theta}f(x)) \right\} \right)=0.
$$

Since $\Re(e^{-i\theta}f) \leq |f|$ on $X$, we conclude that
$$
\mu\left(\left\{ x\in X\,:\,|f(x)|\neq e^{-i\theta}f(x) \right\} \right)=0
$$
where
$$
\theta=\mbox{arg}\left( \int_{X} f\,d\mu \right).
$$

If $\int_{X} f\,d\mu=0$, taking $\theta=0$ yields the result.
\end{proof}

\begin{proposition} \label{625}
Let $f,g \in \mathcal{L}_{\mathbb{C}}(X,\mathsf{S},\mu)$. Then,
\begin{itemize}
\item[(a)] $\int_{A \cup B} f\,d\mu = \int_{A} f\,d\mu + \int_{B} f\,d\mu$ for every $A,B \in \mathsf{S}$ such that $A \cap B = \varnothing$.
\item[(b)] $\int_{X} f\,d\mu = \int_{X} g\,d\mu$ if and only if $f=g$ $\mu$ a.e.
\end{itemize}
\end{proposition}

\begin{proof}
The proof is straightforward since it reduces to examining the real and imaginary parts [Exercise \ref{E639}].
\end{proof}

We now present an extension of the dominated convergence theorem for complex integrable functions.

\begin{theorem}[Dominated convergence theorem in $\mathcal{L}_{\mathbb{C}}(X)$] \label{626} \index{theorem!dominated convergence}
Let $(X,\mathsf{S},\mu)$ be a measure space and let $f_k:X \to \mathbb{C}$ be a sequence of functions in $\mathcal{L}_{\mathbb{C}}(X,\mathsf{S},\mu)$ such that $f_k\to f$ $\mu$ a.e. for some $f \in \mathbb{M}_{\mathbb{C}}(X,\mathsf{S})$. Suppose that there exists a function $g:X \to \mathbb{R}$ in $\mathcal{L}(X,\mathsf{S},\mu)$ such that, for every $k \in \mathbb{N}$, we have $|f_k|\leq g$ $\mu$ a.e. Then, $f \in \mathcal{L}_{\mathbb{C}}(X,\mathsf{S},\mu)$ and
$$
\int_{X} f\,d\mu =\lim_{k \to \infty} \int_{X} f_k\,d\mu.
$$
\end{theorem}

\begin{proof}
The proof is straightforward since it reduces to examining the real and imaginary parts [Exercise \ref{E639}].
\end{proof}

\section{Exercises}

\begin{exercise} \label{E61}
Let $(X,\mathsf{S},\mu)$ be a measure space. Prove that, if $f,g \in \overline{\mathcal{L}}(X,\mathsf{S},\mu)$, then $\max\{f,g\} \in \overline{\mathcal{L}}(X,\mathsf{S},\mu)$ and $\min\{f,g\} \in \overline{\mathcal{L}}(X,\mathsf{S},\mu)$.
\end{exercise}

\begin{exercise} \label{E62}
Let $(X,\mathsf{S},\mu)$ be a measure space, and let $f \in \overline{\mathcal{L}}(X,\mathsf{S},\mu)$ and $A \in \mathsf{S}$ be given. Prove that
$$
\begin{aligned}
\left|\int_{A} f\,d\mu\right| = \int_{A} |f|\,d\mu \,\,\,\mbox{if and only if}\,\,\, &\mu\left( \left\{ x \in A\,:\, f(x)>0\right\}\right)=0\,\,\mbox{or}\,\,\\
&\mu \left( \left\{ x \in A\,:\, f(x)<0\right\}\right)=0.
\end{aligned}
$$
\end{exercise}

\begin{exercise} \label{E63}
Let $(X,\mathsf{S},\mu)$ be a measure space and let $f \in \overline{\mathcal{L}}(X,\mathsf{S},\mu)$. Prove, without using {\rm Corollary \ref{66}}, that the following statements are equivalent:
\begin{itemize}
\item[(a)] $f=0$ $\mu$ a.e.
\item[(b)] $\int_{A} f\,d\mu=0$ for every $A \in \mathsf{S}$.
\item[(c)] $\int_X |f|\,d\mu=0$.
\end{itemize}
\end{exercise}

{\setlength{\parindent}{0pt}
\begin{exercise} \label{E64}
Let $(X,\mathsf{S},\mu)$ be a measure space and let $f,g \in \overline{\mathbb{M}}(X,\mathsf{S})$ be such that $f=g$ $\mu$ a.e. Prove that, if $f \in \overline{\mathcal{L}}(X,\mathsf{S},\mu)$ or $g \in \overline{\mathcal{L}}(X,\mathsf{S},\mu)$, then $f,g \in \overline{\mathcal{L}}(X,\mathsf{S},\mu)$ and $\int_X f\,d\mu = \int_X g\,d\mu$.
\end{exercise}

(Hint: First consider the case where $f$ and $g$ are nonnegative. Define the measurable set $A:=\{ x \in X\,:\, f(x) \neq g(x) \}$ and the function $h$ as follows: $h(x):=+\infty$ if $x \in A$ and $h(x):=0$ if $x \not\in A$. Prove that $\int_X h\,d\mu =0$ using the monotone convergence theorem).}

\begin{exercise} \label{E65}
Let $(X,\mathsf{S},\mu)$ be a measure space, and let $f \in \overline{\mathcal{L}}(X,\mathsf{S},\mu)$ and $\gamma >0$ be given. Prove that the set $\left\{ x \in X\,:\, |f(x)|>\gamma \right\}$ has finite measure.
\end{exercise}

\begin{exercise} \label{E66}
Let $(X,\mathsf{S},\mu)$ be a measure space and let $f \in \overline{\mathcal{L}}(X,\mathsf{S},\mu)$ be given. Prove that the set $\left\{ x \in X\,:\, f(x) \neq 0 \right\}$ is $\sigma$-finite.
\end{exercise}

\begin{exercise} \label{E67}
Let $(X,\mathsf{S},\mu)$ be a measure space, and let $f \in \overline{\mathcal{L}}(X,\mathsf{S},\mu)$ and $\varepsilon >0$ be given. Prove that there exists an $\mathsf{S}$-simple integrable function $\zeta$ such that
$$
\int_{X} |f-\zeta|\,d\mu< \varepsilon.
$$
\end{exercise}

{\setlength{\parindent}{0pt}
\begin{exercise} \label{E68}
Let $(X,\mathsf{S},\mu)$ be a measure space and let $f \in \overline{\mathcal{L}}(X,\mathsf{S},\mu)$ be given. Prove the following statements:
\begin{itemize}
\item[(a)] For every $\varepsilon >0$, there exists $\delta>0$ such that, for every $A\in \mathsf{S}$,
$$
\left|\int_{A} f\,d\mu\right|<\varepsilon\quad \mbox{whenever}\quad \mu(A)<\delta.
$$ 
\item[(b)] For every $\varepsilon >0$, there exists $A \in \mathsf{S}$ such that, for every $B \in \mathsf{S}$,
$$
\left|\int_{X} f\,d\mu - \int_{B} f\,d\mu \right|<\varepsilon\quad \text{whenever}\quad A \subset B.
$$
\end{itemize}
\end{exercise}

(Hint: In part \textit{(a)} argue by contradiction and use the Borel--Cantelli lemma.)}

\begin{exercise} \label{E69}
Let $(X,\mathsf{S},\mu)$ be a measure space and let $f,g,h\in \overline{\mathbb{M}}(X,\mathsf{S})$. Prove that, if $f,g \in \overline{\mathcal{L}}(X,\mathsf{S},\mu)$ and $h$ satisfies $h=f+g$ $\mu$ a.e., then $h \in \overline{\mathcal{L}}(X,\mathsf{S},\mu)$ and $\int_X h\,d\mu = \int_X f\,d\mu + \int_X g\,d\mu$.
\end{exercise}

\begin{exercise} \label{E610}
Let $(\mathbb{R},\mathcal{B}(\mathbb{R}),\lambda)$ be a measure space. Define the sign function $\mbox{\rm sgn}:\mathbb{R} \to \mathbb{R}$ by
$$
\mbox{\rm sgn}(x):=\left\{
\begin{array}{rcl}
 1    & & \mbox{if}\,\,\, x>0, \\
 0    & & \mbox{if}\,\,\, x=0, \\
 -1   & & \mbox{if}\,\,\, x<0.\\
\end{array}
\right.
$$

Is it possible to compute the Lebesgue integral of $\mbox{\rm sgn}$ on $\mathbb{R}$? Justify your answer.
\end{exercise}

\begin{exercise} \label{E611}
Give an example of a measurable function $f$ such that $f \in \overline{\mathcal{L}}(X,\mathsf{S},\mu)$ but $f^{2} \not\in \overline{\mathcal{L}}(X,\mathsf{S},\mu)$.
\end{exercise}

\begin{exercise} \label{E612}
Let $X=\mathbb{N}$, $\mathsf{S}=\mathcal{P}(\mathbb{N})$, and $\mu=\mu^{\sharp}$ be the counting measure. Prove that $f \in \mathcal{L}(X,\mathsf{S},\mu)$ if and only if the series $\sum_{k=1}^{\infty} f(k)$ converges absolutely.
\end{exercise}

\begin{exercise}[P. Chebyshev's inequality] \label{E613} \index{inequality!of Chebyshev}
Let $(X,\mathsf{S},\mu)$ be a measure space and let $\varphi:[0,+\infty) \to [0,+\infty)$ be a nondecreasing function which is positive on $(0,+\infty)$. Prove that, if $f \in \mathbb{M}(X,\mathsf{S})$ is such that $\varphi \circ |f| \in \mathcal{L}(X,\mathsf{S},\mu)$, then for every $\varepsilon>0$
$$
\mu(\{x \in X\,:\, |f(x)| \geq \varepsilon \}) \leq \frac{1}{\varphi(\varepsilon)} \int_X \varphi \circ |f| \,d\mu.
$$
\end{exercise}

{\setlength{\parindent}{0pt}
\begin{exercise}[Mean value lemma] \label{E614} \index{lemma!mean value}
Let $(X,\mathsf{S},\mu)$ be a measure space and let $f \in \mathcal{L}(X,\mathsf{S},\mu)$ be such that $\left| \frac{1}{\mu(A)}\int_{A} f\,d\mu \right| \leq \gamma$ for some constant $\gamma \geq 0$, for every $A \in \mathsf{S}$ with $0<\mu(A)<+\infty$. Prove that $|f| \leq \gamma$ $\mu$ a.e.
\end{exercise}

(Hint: Use Exercise \ref{E518}).}

{\setlength{\parindent}{0pt}
\begin{exercise} \label{E615}
Let $(X,\mathsf{S},\mu)$ and $(X,\widetilde{\mathsf{S}},\widetilde{\mu})$ be two measure spaces such that $\mathsf{S} \subset \widetilde{\mathsf{S}}$ and $\widetilde{\mu}_{|_{\mathsf{S}}}=\mu$. Prove that $\mathcal{L}(X,\mathsf{S},\mu) \subset \mathcal{L}(X,\widetilde{\mathsf{S}},\widetilde{\mu})$ and $\int_X f\,d\mu = \int_X f\,d\widetilde{\mu}$ for every $f \in \mathcal{L}(X,\mathsf{S},\mu)$.
\end{exercise}

(Hint: Let $0<\gamma<+\infty$ and let $f:X \to [0,\gamma]$ be an $\mathsf{S}$-measurable function. Assume that
$$
\mu\left( \left\{ x \in X\,:\,f(x)>0\right\}\right)<+\infty.
$$

Prove that $\int_X f\,d\mu = \int_X f\,d\widetilde{\mu}$ and similarly for the function $\gamma-f$. Also use Exercise \ref{E523}).}

\begin{exercise}[A generalization of the dominated convergence theorem] \label{E616} \index{theorem!dominated convergence}
Let $(X,\mathsf{S},\mu)$ be a measure space and let $(g_k)$ be a sequence of elements of $\mathcal{L}^{+}(X,\mathsf{S},\mu)$ such that $g_k \to g$ $\mu$ a.e. for some function $g \in \mathcal{L}^{+}(X,\mathsf{S},\mu)$ and $\int_X g_k \,d\mu \to \int_X g\,d\mu$.

Prove that, if $(f_k)$ is a sequence in $\mathbb{M}(X,\mathsf{S})$ such that $|f_k|\leq g_k$ $\mu$ a.e. for every $k \in \mathbb{N}$ and $f_k \to f$ for some function $f \in \mathbb{M}(X,\mathsf{S})$ $\mu$ a.e., then
\begin{itemize}
\item[(a)] $f \in \mathcal{L}(X,\mathsf{S},\mu)$, and
\item[(b)] $\int_{A} f\,d\mu=$ $\displaystyle\lim_{k \to \infty}$ $\int_{A} f_k\,d\mu$ for every $A \in \mathsf{S}$.
\end{itemize}
\end{exercise}

\begin{exercise} \label{E617}
Let $(X,\mathsf{S},\mu)$ be a finite measure space and let $(f_k)$ be a sequence of elements of $\mathcal{L}(X,\mathsf{S},\mu)$ that converges uniformly on $X$ to a function $f \in \mathbb{M}(X,\mathsf{S})$. Prove that
$$
\lim_{k \to \infty} \int_X f_k\,d\mu = \int_X f\,d\mu.
$$

What happens if $\mu(X)=+\infty$? Justify your answer.
\end{exercise}

{\setlength{\parindent}{0pt}
\begin{exercise} \label{E618}
Consider the measurable space $(\mathbb{R},\mathcal{B}(\mathbb{R}))$. Prove that if $\mu,\nu$ are two finite measures on $(\mathbb{R},\mathcal{B}(\mathbb{R}))$ such that $\int_{\mathbb{R}} f\,d\nu = \int_{\mathbb{R}} f\,d\mu$ for every bounded continuous function $f:\mathbb{R} \to \mathbb{R}$, then $\nu=\mu$. 
\end{exercise}

(Hint: Consult the uniqueness theorem for measures).}

\begin{exercise}[Lieb's theorem] \label{E619} \index{theorem!of Lieb}
Let $(X,\mathsf{S},\mu)$ be a measure space and let $(f_k)$ be a sequence of elements of $\mathcal{L}(X,\mathsf{S},\mu)$ converging to a function $f \in \mathcal{L}(X,\mathsf{S},\mu)$ $\mu$ a.e. Prove that
$$
\lim_{k \to \infty}\left\vert \int_{X}|f_k|\,d\mu-\int_{X}|f|\,d\mu-\int_{X}|f_k-f|\,d\mu \right\vert=\lim_{k \to\infty}\int_{X}||f_k|-|f|-|f_k-f||\,d\mu=0.
$$
\end{exercise}

\begin{exercise}[Scheffé's lemma] \label{E620} \index{lemma!of Scheffé}
Let $(X,\mathsf{S},\mu)$ be a measure space and let $(f_k)$ be a sequence of elements of $\mathcal{L}(X,\mathsf{S},\mu)$ converging to a function $f \in \mathcal{L}(X,\mathsf{S},\mu)$ $\mu$ a.e. Prove that $\lim_{k \to \infty} \int_X |f_k| \,d\mu = \int_X |f|\,d\mu$ if and only if $\lim_{k \to \infty} \int_{X} |f_k-f|\,d\mu=0$.
\end{exercise}

\begin{exercise} \label{E621}
Let $(X,\mathsf{S},\mu)$ be a measure space and let $f \in \mathcal{L}(X,\mathsf{S},\mu)$ be such that $\int_A f\,d\mu \leq M < \infty$ whenever $\mu(A)<+\infty$. Prove that $\int_X f\,d\mu \leq M$.
\end{exercise}

\begin{exercise} \label{E622}
Let $(X,\mathsf{S},\mu)$ be a measure space and let $f \in \mathcal{L}(X,\mathsf{S},\mu)$. Suppose that for some $A \in \mathsf{S}$ there exists a disjoint sequence $(A_k)$ of elements of $\mathsf{S}$ such that $\bigcup_{k=1}^{\infty}A_k = A$. Prove that
$$
\int_{A} f\,d\mu = \sum_{k=1}^{\infty} \int_{A_k} f\,d\mu.
$$
\end{exercise}

\begin{exercise} \label{E623}
Let $(X,\mathsf{S},\mu)$ be a measure space and let $f \in \mathcal{L}(X,\mathsf{S},\mu)$. Prove the following statements:
\begin{itemize}
\item[(a)] If $(A_k)$ is an increasing sequence of elements of $\mathsf{S}$, then
$$
\int_{\bigcup_{k=1}^{\infty}\,A_k} f\,d\mu = \lim_{k \to \infty} \int_{A_k} f\,d\mu.
$$
\item[(b)] If $(A_k)$ is a decreasing sequence of elements of $\mathsf{S}$, then
$$
\int_{\cap_{k=1}^{\infty}\,A_k} f\,d\mu = \lim_{k \to \infty} \int_{A_k} f\,d\mu.
$$
\end{itemize}
\end{exercise}

\begin{exercise} \label{E624}
Let $(X,\mathsf{S})$ be a measurable space and let $(\mu_k)$ be an arbitrary sequence of measures on it. Consider the measure $\mu:=\sum_{k=1}^{\infty} \mu_k$ on $(X,\mathsf{S})$ {\rm [Exercise \ref{E410}]}. Prove that, if $f \in \overline{\mathbb{M}}^{+}(X,\mathsf{S})$, then $\int_X f\,d\mu <+\infty$ if and only if
$$
\sum_{k=1}^{\infty} \int_X f\,d\mu_k <+\infty.
$$

Moreover, prove that, if $\int_X f\,d\mu <+\infty$, then $\int_X f\,d\mu = \sum_{k=1}^{\infty} \int_X f\,d\mu_k$.
\end{exercise}

\begin{exercise} \label{E625}
Let $(X,\mathsf{S},\mu)$ be a measure space and let $(f_k)$ be a sequence of elements of $\mathcal{L}(X,\mathsf{S},\mu)$ such that $0 \leq f_k(x) \leq 1$ for every $x \in X$ and such that $f_k \to 1$ on $X$. Suppose that for some $A \in \mathsf{S}$ with $\mu(A)<\infty$ one has $f_k(x)=1$ for every $k \in \mathbb{N}$ and every $x \in X \smallsetminus A$. Prove that $\lim_{k \to \infty} \int_X (1-f_k)\,d\mu =0$.
\end{exercise}

\begin{exercise} \label{E626}
Let $(\mathbb{N},\mathcal{P}(\mathbb{N}))$ be a measurable space and let $\mu:\mathcal{P}(\mathbb{N}) \to \mathbb{R}$ be the measure given by $\mu(A):=\sum_{i \in A} \frac{1}{2^{i}}$ or $0$ if $A=\varnothing$. Define the sequence of functions $\eta_k : \mathbb{N} \to \mathbb{R}$ by $\eta_k:=2^{k}\,\chi_{\{k+1,k+2,\ldots \}}$. Prove that $\eta_k \to 0$ pointwise on $\mathbb{N}$ and that there does not exist any integrable function $\xi$ dominating the sequence.
\end{exercise}

{\setlength{\parindent}{0pt}
\begin{exercise} \label{E627}
Let $(\mathbb{R},\mathcal{B}(\mathbb{R}),\lambda)$ be a measure space. In the following items explain whether it is possible to apply the Lebesgue dominated convergence theorem to the sequence of functions $(f_k)$.
\begin{itemize}
\item[(a)] $f_k: \mathbb{R} \to \mathbb{R}$ given by $f_k:=k\chi_{[0,\frac{1}{k}]}$.
\item[(b)] $f_k: \mathbb{R} \to \mathbb{R}$ given by $f_k:=\chi_{[k,k+1)}$.
\end{itemize}
\end{exercise}

(Hint: Consider the function $\xi(x):=\displaystyle\sup_{k \in \mathbb{N}}|f_k(x)|$ in both items).}

\begin{exercise} \label{E628}
Let $(X,\mathsf{S},\mu)$ be a measure space and let $f \in \mathbb{M}(X,\mathsf{S})$. For each $k \in \mathbb{N}$, let $f_k:X \to \mathbb{R}$ be the $k$-truncation of $f$ {\rm [Exercise \ref{E322}]}. Prove the following statements:
\begin{itemize}
\item[(a)] If $f \in \mathcal{L}(X,\mathsf{S},\mu)$, then $\int_X f\,d\mu =$ $\displaystyle\lim_{k \to \infty}$ $\int_X f_k\,d\mu$.
\item[(b)] If $\displaystyle\sup_{k \in \mathbb{N}} \int_X |f_k|\,d\mu < \infty$, then $f \in \mathcal{L}(X,\mathsf{S},\mu)$.
\end{itemize}
\end{exercise}

\begin{exercise} \label{E629}
Let $([0,1],\mathcal{B}([0,1]),\lambda)$ be a measure space. Define the sequence of functions $f_k:=kxe^{-kx}\chi_{[0,1]}$. Prove the following statements:
\begin{itemize}
    \item[(a)] $f_k$ is Borel-measurable for every $k \in \mathbb{N}$.
    \item[(b)] $f_k \to 0$ pointwise on $[0,1]$.
    \item[(c)] There exists an integrable function $g$ dominating the sequence $(f_k)$. 
    \item[(d)] Compute the following limit
    $$
    \lim_{k\to \infty} R\int_{0}^{1} kxe^{-kx}\,dx.
    $$
\end{itemize}
\end{exercise}

\begin{exercise} \label{E630}
Let $([0,1],\mathcal{B}([0,1]),\lambda)$ be a measure space and let $(\xi_k)$ be a sequence of real numbers. Define $f_k:[0,1] \to \mathbb{R}$ by
$$
f_k(x):=\left\{
\begin{array}{lcl}
2\xi_k kx  & & \,\,\,\,\mbox{if}\,\,\, 0 \leq x \leq \displaystyle\tfrac{1}{2k},\\
2\xi_k(1-kx) & & \,\,\,\,\mbox{if}\,\,\, \tfrac{1}{2k} \leq x \leq \tfrac{1}{k},\\
0 & & \,\,\,\,\mbox{if}\,\,\, \tfrac{1}{k} \leq x \leq 1.\\
\end{array}
\right.
$$

 Prove the following statements:
\begin{itemize}
    \item[(a)] $f_k$ is Borel-measurable for every $k \in \mathbb{N}$.
    \item[(b)] $f_k \to 0$ pointwise on $[0,1]$.
    \item[(c)] For every $k \in \mathbb{N}$ one has
    $$
    \int_{[0,1]} f_k\,d\lambda = \frac{\xi_k}{2k}.
    $$
    \item[(d)] What conditions must the sequence $(\xi_k)$ satisfy in order to apply the Lebesgue dominated convergence theorem? Justify your answer.
\end{itemize}
\end{exercise}

\begin{exercise} \label{E631}
Let $([0,1],\mathcal{B}([0,1]),\lambda)$ be a measure space and let $(\xi_k)$ be a sequence of real numbers. Define $f_k :[0,1] \to \mathbb{R}$ by
$$
f_k(x):=\left\{
\begin{array}{lcl}
0 & & \,\,\,\,\mbox{if}\,\,\,  x\in \left[0,\frac{1}{k}\right],\\
\frac{1}{\sqrt{x}}  & & \,\,\,\,\mbox{if}\,\,\, x \in \left(\frac{1}{k},1 \right].
\end{array}
\right.
$$

\begin{itemize}
    \item[(a)] Prove that $f_k$ is Borel-measurable for every $k \in \mathbb{N}$.
    \item[(b)] For each $k \in \mathbb{N}$, compute the integral
    $$
    \int_{[0,1]} f_k\,d\lambda .
    $$
    \item[(c)] Show that
    $$
\lim_{k \to \infty} f_k(x):=\left\{
\begin{array}{lcl}
\frac{1}{\sqrt{x}}  & & \,\,\,\,\mbox{if}\,\,\, x \in \left(0,1 \right],\\
0 & & \,\,\,\,\mbox{if}\,\,\,  x=0.\\
\end{array}
\right.
$$
   \item[(d)] Is $\lim_{k \to \infty} f_k$ integrable on $[0,1]$ with respect to $\lambda$? If so, compute its integral.
\end{itemize}
\end{exercise}

\begin{exercise} \label{E632}
Let $([0,1],\mathcal{B}([0,1]),\lambda)$ be a measure space and let $f:[0,1] \to \mathbb{R}$ be a continuous function. Prove that
$$
\lim_{k \to \infty} R\int_{0}^{1} kf(x)e^{-kx}\,dx = f(0).
$$
\end{exercise}

{\setlength{\parindent}{0pt}
\begin{exercise} \label{E633}
Let $X=\mathbb{R}$, $\mathsf{S}=\mathcal{B}(\mathbb{R})$, and $\mu=\lambda$ the Lebesgue measure on $\mathsf{S}$. Let $f:\mathbb{R} \to \mathbb{R}$ be a Borel-measurable function such that its restriction to $[a,b]$ is Riemann-integrable whenever $a<b \in \mathbb{R}$. Suppose that $\displaystyle\lim_{k \to \infty} $ $R\int_{-k}^{k} |f(x)|\,dx < +\infty$. Prove that $f \in \mathcal{L}(\mathbb{R},\mathcal{B}(\mathbb{R}),\lambda)$ and that
$$
\int_{\mathbb{R}} f\,d\lambda = R\int_{-\infty}^{\infty} f(x)\,dx.
$$
\end{exercise}

(Hint: Use the dominated convergence theorem).}

\begin{exercise} \label{E634}
Let $a \in (-\infty,+\infty)$ be fixed. On $(a,+\infty)$ consider the Borel $\sigma$-algebra and the Lebesgue measure. Prove that, if $f:(a,+\infty)\to \mathbb{R}$ is a Borel-measurable function and Riemann-integrable on every closed interval $[c,d]$ with $d > c$ in $(a,+\infty)$, then $f$ is Lebesgue-integrable if and only if $|f|$ is improperly Riemann-integrable on $(a,+\infty)$. Moreover, prove that
$$
\int_{(a,+\infty)} f\,d\lambda = R\int_{a}^{\infty} f(x)\,dx.
$$

Finally, give an example where such a relation does not hold.
\end{exercise}

\begin{exercise} \label{E635}
Prove that, if $\gamma >0$, then
$$
R\int_{0}^{\infty} e^{-\gamma \xi}\,d\xi = \frac{1}{\gamma}.
$$
\end{exercise}

\begin{exercise} \label{E636}
Prove the following identity:
$$
R\int_{0}^{\infty} e^{-\xi^2}\,d\xi = \frac{\sqrt{\pi}}{2}.
$$
\end{exercise}

\begin{exercise} \label{E637}
Let $n \in \mathbb{N}$ be fixed. Prove that
$$
R\int_{0}^{\infty} \xi^n e^{-\xi}\,d\xi = n!.
$$
\end{exercise}

\begin{exercise} \label{E638}
Let $(X,\mathsf{S},\mu)$ be a measure space and let $a,b$ be fixed real numbers such that $a<b$. Suppose that $f:X \times (a,b) \to \mathbb{R}$ is a function such that $f^{t}:X \to \mathbb{R}$ given by $f^{t}(x):=f(x,t)$ is $\mathsf{S}$-measurable for every $t \in (a,b)$ and that there exist $t_{0},t_{1} \in [a,b]$ such that $f^{t_0} \in \mathcal{L}(X,\mathsf{S},\mu)$; that $\frac{\partial f^{t_{1}}}{\partial t}$ exists and that there exists an integrable function $g:X \to \mathbb{R}$ such that:
$$
\left| \frac{f^{t}(x)-f^{t_1}(x)}{t-t_1}\right| \leq g(x) \quad \forall x \in X,\quad \forall t \in [a,b]\smallsetminus\{t_1\}.
$$

Prove that
$$
\left(\frac{\partial}{\partial t} \int_X f^{t}\,d\mu\right)_{t=t_1} = \int_X \frac{\partial f^{t_1}}{\partial t}\,d\mu.
$$
\end{exercise}

\begin{exercise} \label{E639}
Let $(X,\mathsf{S},\mu)$ be a measure space. Let $A,B \in \mathsf{S}$ be disjoint and let $f,g \in \mathcal{L}_{\mathbb{C}}(X,\mathsf{S},\mu)$ be given. 
\begin{itemize}
\item[(a)] Prove that $\int_{A \cup B} f\,d\mu = \int_{A} f\,d\mu + \int_{B} f\,d\mu$.
\item[(b)] Prove that $\int_{X} f\,d\mu = \int_{X} g\,d\mu$ if and only if $f=g$ $\mu$ almost everywhere.
\item[(c)] State and prove the Lebesgue dominated convergence theorem for complex-valued functions.
\end{itemize}
\end{exercise}

Let $\mathcal{I}:=(a,b)$ be an interval of $\mathbb{R}$ with $-\infty \leq a < b \leq \infty$. A function $\varphi:\mathcal{I} \to \mathbb{R}$ is said to be \textbf{convex} if for every $x_0,x_1 \in \mathcal{I}$ and $t \in [0,1]$ one has:
$$
\varphi(tx_0+(1-t)x_1) \leq t\,\varphi(x_0)+(1-t)\,\varphi(x_1).
$$

\begin{exercise}[Jensen's inequality] \label{E640} \index{inequality!Jensen's}
Solve the following items:

\begin{itemize}
\item[(a)] Let $\mathcal{I}:=(a,b)$ be an interval of $\mathbb{R}$ with $-\infty \leq a < b \leq \infty$ and let $\varphi:\mathcal{I} \to \mathbb{R}$ be a convex function. Prove the following:
          \begin{itemize}
          \item[(a.1)] For every point $x_0 \in \mbox{\rm int}(\mathcal{I})$ (interior of $\mathcal{I}$), there exists a line passing through the point $(x_0,\varphi(x_0))$, say with equation $y=ax+b$, such that $ax+b \leq \varphi(x)$ for every $x \in \mathcal{I}$.
          \item[(a.2)] $\varphi$ is continuous on $\mbox{\rm int}(\mathcal{I})$. 
          \end{itemize}
\item[(b)] Let $(X,\mathsf{S},\mu)$ be a measure space such that $\mu(X)=1$ and let $f \in \mathcal{L}(X,\mathsf{S},\mu)$. If $\mathcal{I}:=(a,b)$ is an interval in $\mathbb{R}$ ($-\infty \leq a < b \leq \infty$) such that $f \in \mathcal{I}$ $\mu$ almost everywhere, prove that $\int_{X} f\,d\mu \in \mathcal{I}$. 
\item[(c)] Let $(X,\mathsf{S},\mu)$ be a measure space such that $\mu(X)=1$, let $f \in \mathcal{L}(X,\mathsf{S},\mu)$, and let $\mathcal{I}:=(a,b)$ be an interval in $\mathbb{R}$ ($-\infty \leq a < b \leq \infty$). Prove that, if $\varphi:\mathcal{I} \to \mathbb{R}$ is a convex function and $f \in \mathcal{I}$ $\mu$ almost everywhere, then $\int_X f \,d\mu \in \mathcal{I}$, and the inequality
$$
\varphi\left(\int_X f\,d\mu \right) \leq \int_{X} (\varphi\circ f)\,d\mu
$$
holds.
\end{itemize}
\end{exercise}

\section{Project: The change of variables theorem}

Let $(X,\mathsf{S},\mu)$ be a measure space. In general, the computation of the integral
$$
\int_X f\,d\mu
$$
of a function $f\in \mathcal{L}(X,\mathsf{S},\mu)$ may be a complex and highly abstract process, since the measure space does not, in principle, possess any additional structure that facilitates the explicit computation of integrals.

Even in apparently simple situations conceptual difficulties arise. For example, if $X=[a,b]$ with $0\leq a<b$, $\mathsf{S}=\mathcal{B}([a,b])$, and $\mu=\lambda$ is the Lebesgue measure, the computation of the integral
$$
\int_{[a,b]} x\,d\lambda
$$
is not immediate if one dispenses with Henri Lebesgue's theorem, even though it is known that the identity function can be approximated pointwise by a nondecreasing sequence of simple functions.

This observation highlights the need for theoretical tools that make it possible to transform complicated integrals into more manageable expressions. With this objective in mind, we propose to the reader a series of exercises aimed at establishing fundamental results that facilitate the computation of integrals through suitable transformations of the measure space. We shall refer to these results as \textbf{change of variables theorems}.

\subsection*{6.6.1.\quad Objective}

The objective of this project is to present and prove the change of variables theorems, beginning with the particular case of the measure space $(\mathbb{R},\mathcal{B}(\mathbb{R}),\lambda)$ and subsequently extending the result to the general case of a measure space $(X,\mathsf{S},\mu)$.

\subsection*{6.6.2.\quad Procedure}

\begin{itemize}
    \item[1.] Let $(\mathbb{R},\mathcal{B}(\mathbb{R}),\lambda)$ be a measure space. We say that a function $\phi:\mathbb{R} \to \mathbb{R}$ is a \textbf{homeomorphism} if $\phi$ is continuous and bijective and its inverse $\phi^{-1}:\mathbb{R} \to \mathbb{R}$ is continuous. \index{homeomorphism}

    Prove that if $\phi:\mathbb{R} \to \mathbb{R}$ is a homeomorphism, then $\phi(\mathcal{B}(\mathbb{R}))=\mathcal{B}(\mathbb{R})$.

    \item[2.] We say that a function $\tau:\mathbb{R} \to \mathbb{R}$ is an \textbf{affine transformation} if it is of the form $\tau(x):=ax+\zeta$ with $a,\zeta \in \mathbb{R}$ and $a\neq 0$. \index{function!affine}

    Prove that if $\tau:\mathbb{R} \to \mathbb{R}$ is an affine transformation, then $\tau(\mathcal{B}(\mathbb{R}))=\mathcal{B}(\mathbb{R})$.

    \item[3.] Prove that if $\tau:\mathbb{R} \to \mathbb{R}$ is an affine transformation with $\tau(x)=ax+\zeta$, then $\lambda(\tau(B))=|a|\lambda(B)$ for every $B \in \mathcal{B}(\mathbb{R})$.

    \item[4.] Let $a \in \mathbb{R}$ with $a \neq 0$ and $\zeta \in \mathbb{R}$. Prove that, if $f \in \mathcal{L}(\mathbb{R},\mathcal{B}(\mathbb{R}),\lambda)$ and $\tau(x)=ax+\zeta$ is an affine transformation, then the function $f\circ \tau$ belongs to $\mathcal{L}(\mathbb{R},\mathcal{B}(\mathbb{R}),\lambda)$ and the following identity holds:
$$
\int_{B} (f\circ \tau)\,d\lambda = |a|^{-1} \int_{\tau(B)} f\,d\lambda\quad \mbox{for every}\,\,\,B \in \mathcal{B}(\mathbb{R}).
$$

We call this result the \textbf{linear change of variables theorem on $\mathbb{R}$}. \index{theorem!linear change of variables} 

\item[5.] Let $a,b\in \mathbb{R}$ be such that $0\leq a < b$. Prove that if $f(x)=x$ on $[a,b]$, then
$$
\int_{[a,b]} f\,d\lambda = \frac{1}{2}(b^2-a^2).
$$ 

\item[6.] Let $(X,\mathsf{S},\mu)$ be a measure space and let $(Y,\mathsf{T})$ be a measurable space. Given a function $T:X \to Y$ measurable with respect to $\mathsf{S}$ and $\mathsf{T}$, define $\mu_{T}:\mathsf{T} \to \overline{\mathbb{R}}$ by
\begin{eqnarray}
\mu_{T}(F):=\mu(T^{-1}(F))    
\end{eqnarray}
and call it the \textbf{measure induced by $T$ on $Y$}.

Prove that if $f \in \mathcal{L}(Y,\mathsf{T},\mu_{T})$, then $f \circ T \in \mathcal{L}(X,\mathsf{S},\mu)$ and the following identity holds:
$$
\int_{F} f\,d\mu_{T} = \int_{T^{-1}(F)} (f \circ T)\,d\mu\quad \text{for every}\,\, F \in \mathsf{T}.
$$

We call this result the \textbf{change of variables theorem on abstract measure spaces}. \index{theorem!change of variables}
\end{itemize}

\chapter{Lebesgue Spaces} \label{Capitulo7}
\markboth{{\scriptsize 7. LEBESGUE SPACES}}{ {\scriptsize 7. LEBESGUE SPACES}}

The Lebesgue spaces that we shall present in this chapter constitute one of the most relevant families of normed spaces in mathematics, with fundamental applications in functional analysis, partial differential equations, probability, and in disciplines such as finance and engineering. 

We begin by introducing the spaces $\mathcal{L}^p(X)$, formed by measurable functions $f: X \to \mathbb{R}$ such that $|f|^p$ is integrable for $p \in [1,\infty)$, and such that $f$ is essentially bounded for $p = \infty$. These sets constitute seminormed vector spaces in which we shall study key results from measure theory and analysis, such as the Hölder--Riesz and Minkowski inequalities.

One difficulty with the spaces $\mathcal{L}^p(X)$ is that the seminorm defining them does not satisfy all the properties of a norm. To overcome this issue, we shall introduce an equivalence relation identifying functions that differ only on a set of measure zero. In this way, we define the spaces $L^p(X)$ as the collection of equivalence classes of functions in $\mathcal{L}^p(X)$, thus obtaining normed vector spaces. Finally, we shall state and prove the Riesz--Fischer theorem, which guarantees that Lebesgue spaces are Banach spaces, that is, normed spaces in which every Cauchy sequence converges.

\section{Definitions and basic properties}

Let $(X,\mathsf{S},\mu)$ be a measure space. If $p \in (0,\infty)$, we define
$$
\mathcal{L}^{p}(X,\mathsf{S},\mu):=\{ f \in \mathbb{M}(X,\mathsf{S})\,:\,|f|^{p}\in \mathcal{L}(X,\mathsf{S},\mu) \}.
$$

Observe that if $f \in \mathbb{M}(X,\mathsf{S})$ is such that $|f|$ is integrable on $X$, then $f$ is integrable on $X$ by Theorem \ref{62}. Therefore, $\mathcal{L}^{1}(X,\mathsf{S},\mu)$ is simply the collection of measurable and integrable functions $\mathcal{L}(X,\mathsf{S},\mu)$ studied in the previous chapter. Hence, it is clear that $f \in \mathcal{L}^{p}(X,\mathsf{S},\mu)$ if and only if $|f|^{p} \in \mathcal{L}^{1}(X,\mathsf{S},\mu)$.

For $f \in \mathcal{L}^{p}(X,\mathsf{S},\mu)$, we denote by
\begin{eqnarray} \label{F71} \index{space!Lp@$\mathcal{L}^{p}(X,\mathsf{S},\mu)$}
|f|_{p}:=\left(\int_X |f|^p\,d\mu \right)^{1/p}.
\end{eqnarray}

Now, we define \index{space!Lpinfty@$\mathcal{L}^{\infty}(X,\mathsf{S},\mu)$}
$$
\mathcal{L}^{\infty}(X,\mathsf{S},\mu):=\{ f \in \mathbb{M}(X,\mathsf{S})\,:\,f\,\,\text{is essentially bounded} \}.
$$

A function $f \in \mathcal{L}^{\infty}(X,\mathsf{S},\mu)$ is called an \textbf{essentially bounded function}. For $f \in \mathcal{L}^{\infty}(X,\mathsf{S},\mu)$, we denote by \index{function!essentially bounded}
\begin{eqnarray} \label{F72}
|f|_{\infty}:=\inf\left\{ c \geq 0\,:\,\mu\left( \{x \in X\,:\, |f(x)|> c \} \right)=0 \right\}.
\end{eqnarray}

It is clear from the definition of an essentially bounded function that the set in (\ref{F72}) over which the infimum is taken is nonempty. Let us establish some properties.

\begin{theorem} \label{71}
Let $f \in \mathcal{L}^{\infty}(X,\mathsf{S},\mu)$ be fixed. Then,
\begin{itemize}
\item[(a)] $|f|_{\infty} \geq 0$.
\item[(b)] $|f| \leq |f|_{\infty}$ $\mu$ a.e.
\item[(c)] If $|f|_{\infty}>0$, then $\mathfrak{I}_{f}:=\{c>0\,:\,\mu(\{x \in X\,:\,|f(x)|>c \})=0 \}$ is equal to the interval $\left[|f|_{\infty},+\infty \right)$.
\item[(d)] If $f$ is bounded on $X$, then $|f|_{\infty} \leq \sup\left\{|f(x)|\,:\,x \in X \right\}$.
\item[(e)] $|f|_{\infty}=\sup\left\{c \geq 0 \,:\,\mu(\{ x \in X\,:\,|f(x)| \geq c \})>0 \right\}$.
\end{itemize}
\end{theorem}

\begin{proof}
\textit{(a):} It follows immediately from the definition. 

\textit{(b):} By definition of infimum, for each $k \in \mathbb{N}$, there exists $c_{k} \geq 0$ such that $|f|_{\infty} \leq c_{k} < |f|_{\infty}+\frac{1}{k}$ and $\mu\left(\{x \in X\,:\, |f(x)|>c_{k}\} \right)=0$.

Since
$$
\left\{ x\in X\,:\,|f(x)|>|f|_{\infty}\right\}=\bigcup_{k=1}^{\infty}\left\{ x\in X\,:\,|f(x)|>|f|_{\infty}+\frac{1}{k}\right\},
$$
then
$$
\begin{aligned}
\mu\left(\left\{x \in X\,:\,\vert f(x)\vert > \vert f \vert_\infty \right\} \right) &\leq\sum_{k=1}^{\infty} \mu\left(\left\{x \in X\,:\,\vert f(x)\vert > \vert f \vert_\infty + \frac{1}{k}\right\} \right)\\
&\leq\sum_{k=1}^{\infty} \mu\left(\left\{x \in X\,:\,\vert f(x)\vert > c_k \right\} \right)=0.
\end{aligned}
$$

Therefore, $|f| \leq \vert f \vert_{\infty}$ $\mu$ a.e.

\textit{(c):} Suppose that $\vert f \vert_{\infty}>0$. It is clear that $\vert f \vert_{\infty} \leq c$ for every $c \in \mathfrak{I}_{f}$, so $\mathfrak{I}_{f} \subset [\vert f \vert_{\infty},+\infty)$.

Conversely: let $c \geq \vert f \vert_{\infty}$ be arbitrary. Since
$$
\left\{x \in X\,:\,|f(x)|>c \right\} = \bigcup_{k=1}^{\infty} \left\{x \in X\,:\,|f(x)|>c+\frac{1}{k} \right\}
$$
and $c+\frac{1}{k} >\vert f \vert_{\infty}$ for every $k \in \mathbb{N}$, then
$$
\begin{aligned}
\mu\left(\left\{x \in X\,:\,|f(x)|>c \right\} \right) &\leq \sum_{k=1}^{\infty} \mu\left(\left\{x \in X\,:\,|f(x)|>c+\frac{1}{k} \right\} \right)\\
& \leq \sum_{k=1}^{\infty} \mu\left(\left\{x \in X\,:\,|f(x)|>\vert f \vert_{\infty}\right\} \right) =0.
\end{aligned}
$$

Consequently, $c \in \mathfrak{I}_{f}$.

\textit{(d):} If $f$ is bounded on $X$, then $b:=\sup\left\{|f(x)|\,:\,x \in X \right\} \in \mathbb{R}$ so that the set $\{x \in X\,:\,|f(x)|>b\}$ is empty and, therefore,
$$
\mu\left(\{x \in X\,:\,|f(x)|>b\} \right)=0.
$$

Consequently, $\vert f\vert_{\infty} \leq b$.

\textit{(e):} Let $\mathfrak{s}:=\sup\left\{c \geq 0 \,:\,\mu(\{ x \in X\,:\,|f(x)| \geq c \})>0 \right\}$. If $\mathfrak{s}=0$, then 
$$
\mu\left(\{ x \in X\,:\,|f(x)|>c \}\right)=0\quad \forall c > 0
$$
and therefore necessarily $\vert f \vert_{\infty}=0$. Now, if $\vert f \vert_{\infty}=0$, part \textit{(b)} implies that $f=0$ $\mu$ a.e. and, consequently, $\mathfrak{s}=0$. Suppose then that $\vert f \vert_{\infty} >0$. For each $c \in [0,\vert f \vert_{\infty})$, we have $\mu\left(\{x \in X\,:\,|f(x)| \geq c \} \right) >0$ and, consequently, $c \leq \mathfrak{s}$ for every $c \in [0,\vert f \vert_{\infty})$. Thus, $\vert f \vert_{\infty} \leq \mathfrak{s}$. 

If $\vert f \vert_{\infty} <\mathfrak{s}$, then
$$
\mu\left(\{ x \in X\,:\,|f(x)| \geq c\} \right)=0\quad \forall c \in (\vert f \vert_{\infty},\mathfrak{s})
$$
by definition of $\vert f \vert_{\infty}$. Likewise, by definition of $\mathfrak{s}$, we have
$$
\mu\left(\{ x \in X\,:\,|f(x)| \geq c\} \right)>0\quad \forall c \in (\vert f \vert_{\infty},\mathfrak{s}).
$$

Consequently, it is impossible that $\vert f \vert_{\infty} < \mathfrak{s}$ and we conclude that $\vert f \vert_{\infty} =\mathfrak{s}$.
\end{proof}

Next, we shall prove that the set $\mathcal{L}^{p}(X,\mathsf{S},\mu)$ is a vector space over $\mathbb{R}$ for every $p\in(0,\infty]$ with the operations given by:
$$
(f+g)(x):=f(x)+g(x)\quad (\gamma\,f)(x):=\gamma\,f(x)
$$
where $f,g \in \mathcal{L}^{p}(X,\mathsf{S},\mu)$ and $\gamma \in \mathbb{R}$. To do this, we first establish the following inequality.

\begin{lemma}  \label{72}
For every $a,b \in \mathbb{R}$ and $p \in (0,\infty)$, we have
\begin{eqnarray}
|a+b|^{p} \leq 2^{p}(|a|^{p}+|b|^{p}).
\end{eqnarray}
\end{lemma}

\begin{proof}
By the triangle inequality in $\mathbb{R}$, we have
$$
|a+b|^{p} \leq (|a|+|b|)^{p} \leq (2\max\{|a|,|b| \})^{p}=2^{p}\max\{|a|^{p},|b|^{p} \} \leq 2^{p}(|a|^{p}+|b|^{p})
$$
as stated.
\end{proof}

\begin{theorem} \label{73}
$\mathcal{L}^{p}(X,\mathsf{S},\mu)$ is a vector space over $\mathbb{R}$ for every $p \in (0,\infty]$.
\end{theorem}

\begin{proof}
Consider the following two cases:

{\scshape Case 1.}\quad $p \in (0,\infty)$.

First observe that if $f$ and $g$ are measurable, then $|f+g|^{p}$ is measurable for every $p \in (0,\infty)$. Now, if $f,g \in \mathcal{L}^{p}(X,\mathsf{S},\mu)$, Lemma \ref{72} ensures that $|f(x) +g(x)|^{p} \leq 2^{p} (|f(x)|^{p} + |g(x)|^{p})$ for every $x \in X$, and Theorem \ref{62} implies that $f+g \in \mathcal{L}^{p}(X,\mathsf{S},\mu)$.

If $f \in \mathcal{L}^{p}(X,\mathsf{S},\mu)$ and $\gamma \in \mathbb{R}$, then $|\gamma\,f|^{p}=|\gamma|^{p}|f|^{p}$ is integrable on $X$ by Theorem \ref{67}. Consequently, $\gamma\,f \in\mathcal{L}^{p}(X,\mathsf{S},\mu)$.

{\scshape Case 2.}\quad $p =\infty$.

If $f,g \in \mathcal{L}^{\infty}(X,\mathsf{S},\mu)$, Theorem \ref{71} ensures that $|f+g| \leq |f|+|g| \leq |f|_{\infty} + |g|_{\infty}$ $\mu$ a.e. and, therefore, $f+g \in \mathcal{L}^{\infty}(X,\mathsf{S},\mu)$.

Now, let $\gamma \in \mathbb{R}$. Observe that $|\gamma\,f| \leq |\gamma|\,c$ $\mu$ a.e. if $|f|\leq c$ $\mu$ a.e. for some $c \geq 0$. Therefore, $\gamma\,f \in \mathcal{L}^{\infty}(X,\mathsf{S},\mu)$ whenever $f \in \mathcal{L}^{\infty}(X,\mathsf{S},\mu)$. 
\end{proof}

A particularly important property for the remainder of the text is the following.

\begin{proposition} \label{74}
Let $p \in (0,\infty]$ and $f \in \mathcal{L}^{p}(X,\mathsf{S},\mu)$. Then, $|f|_{p}=0$ if and only if $f=0$ $\mu$ a.e.
\end{proposition}

\begin{proof}
Again, we divide the proof into two cases.

{\scshape Case 1.}\quad $p \in (0,\infty)$.

It follows immediately from Corollary \ref{68} that
$$
|f|_{p}^{p}=\int_X |f|^p\,d\mu =0 \,\,\,\,\Longleftrightarrow \,\,\,\,|f|^p=0\,\,\,\mu\text{ a.e.}
$$
and therefore the assertion follows directly.

{\scshape Case 2.} $p =\infty$.

By part \textit{(b)} of Theorem \ref{71}, we have $|f| \leq \vert f \vert_{\infty}$ $\mu$ a.e. If $\vert f \vert_{\infty}=0$, then $0 \leq |f| \leq \vert f \vert_{\infty}=0$ $\mu$ a.e., and therefore necessarily $f=0$ $\mu$ a.e. 

Conversely: if $f=0$ $\mu$ a.e., then $|f|=0$ $\mu$ a.e. Thus, $|f|_{\infty} \leq 0$ and, therefore, $\vert f \vert_{\infty}=0$.

This completes the proof.
\end{proof}

\section{H\"older-Riesz and Minkowski inequalities}

The goal of this section is to establish some properties of the vector spaces $\mathcal{L}^{p}(X,\mathsf{S},\mu)$ when $p \in [1,\infty]$ that will allow us to equip them with a function providing the structure of seminormed spaces. To this end, we require the following results.

The following result, due to Otto H\"older\footnote{Otto H\"older (1857--1937) was a mathematician born in Stuttgart, Germany. He is famous for many theorems including H\"older's inequality, the Jordan-H\"older theorem, the theorem stating that every linearly ordered group satisfying an Archimedean property is isomorphic to a subgroup of the additive group of real numbers, the classification of simple groups up to order 200, and H\"older's theorem implying that the Gamma function does not satisfy any algebraic differential equation. Another important concept related to H\"older is the H\"older condition, which is used in many areas of mathematical analysis, including the theory of partial differential equations and function spaces.}, is fundamental in the study of the structure of the vector spaces $\mathcal{L}^{p}(X,\mathsf{S},\mu)$ as well as in others, for example, spaces of continuous functions, series, and Euclidean spaces. This statement is a generalization of the well-known Cauchy-Buniakowski-Schwarz inequality. According to \cite{Grabinsky}, it was originally proved by O. H\"older for the usual inner product of vectors in $\mathbb{R}^{n}$ and generalized by the Hungarian mathematician F. Riesz\footnote{Frigyes Riesz (1880--1956) was a Hungarian mathematician who made fundamental contributions to functional analysis. He carried out part of the foundational work in the development of functional analysis, and his work has had many important applications in physics. He established the spectral theory for bounded symmetric operators in a form very close to what is now regarded as standard. He also made many contributions to other areas, including ergodic theory and topology, and gave an elementary proof of the mean ergodic theorem.} for functions in $\mathcal{L}^{p}(X,\mathsf{S},\mu)$. 

To prove the H\"older-Riesz inequality, we shall use the well-known Young inequality\footnote{William Henry Young (1863--1942) was an English mathematician. He worked in measure theory, Fourier series, and differential calculus, among other fields, and made contributions to the study of functions of several complex variables. He was the husband of Grace Chisholm Young, with whom he authored or coauthored 214 articles and 4 books. Young's theorem is named in his honor. He was elected Fellow of the Royal Society on May 2, 1907, and served as president of the London Mathematical Society from 1922 to 1924.}.

\begin{lemma}[Young's inequality] \label{75} \index{inequality!Young's}
Let $p,q \in (1,\infty)$ be such that $\tfrac{1}{p}+\tfrac{1}{q} =1$. Then, for every pair of real numbers $a,b \geq 0$, we have
$$
ab \leq \frac{1}{p}\,a^{p} + \frac{1}{q}\,b^{q},
$$
with equality if and only if $a^{p}=b^{q}$.
\end{lemma}

\begin{figure}[ht!]
\centering
\includegraphics[scale=0.275]{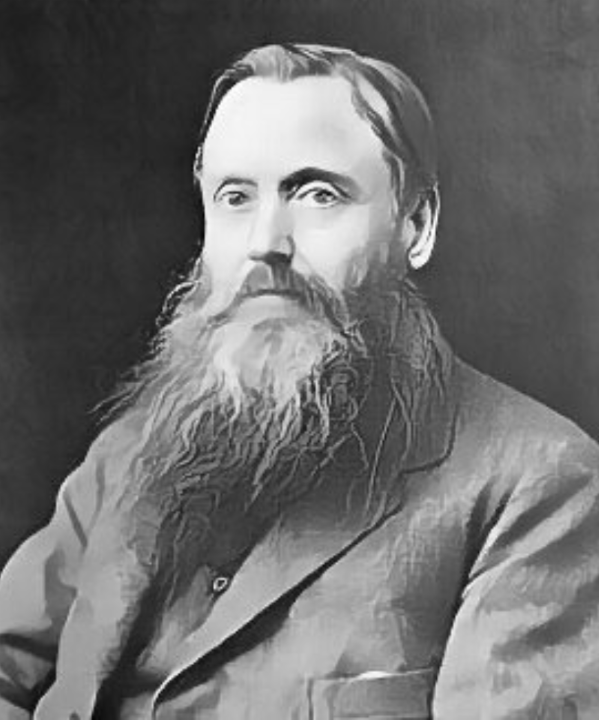} 
\begin{center}
William Henry Young (1863-1942)
\end{center}
\end{figure}

\begin{proof}
The case when $a=0$ or $b=0$ is trivial. Thus, suppose that $a,b \neq 0$. The exponential function is convex, that is, for every $x_0,x_1 \in \mathbb{R}$ and $t \in [0,1]$, we have
$$
(1-t)e^{x_0} + te^{x_1} \geq e^{[(1-t)x_0+tx_1]}.
$$

\begin{figure}[ht!]
\centering
\begin{tikzpicture}[xscale=1.5, yscale=1.55]
\draw[domain= -2.5:0.5, thick] plot(\x,{{exp(\x)}} );
\draw[-, thick] (-2,0.1353)--(0,1);

\draw[->,gray] (-2.8,-0.2)--(1.3,-0.2); 
\draw[->,gray] (-2.7,-0.3)--(-2.7,2); 

\draw (-2,-0.23) node[below]{$_{x_0}$};
\draw (0,-0.23) node[below]{$_{x_1}$};
\draw (-1,-0.22) node[below]{$_{(1-t)x_0+ tx_1}$};

\draw (-2.7,0.56) node[left]{$_{(1-t)\varphi(x_0)+t\varphi(x_1)}$};
\draw (-2.7,0.3678) node[left]{$_{\varphi((1-t)x_0+tx_1)}$};

\draw (-1,0.56) node{$_{\bullet}$};
\draw (-1,0.3678) node{$_{\bullet}$};

\draw (-2,-0.2) node{$_{|}$};
\draw (0,-0.2) node{$_{|}$};
\draw (-1,-0.2) node{$_{|}$};
\draw (-2.7,0.56) node{$-$};
\draw (-2.7,0.3678) node{$-$};
\draw[-, dotted] (0,-0.2)--(0,1);
\draw[-, dotted] (-2,-0.2)--(-2,0.1353);
\draw[-, dotted] (-1,-0.2)--(-1,0.56);

\draw[-, dotted] (-2.7,1)--(0,1);
\draw[-, dotted] (-2.7,0.1353)--(-2,0.1353);
\draw[-, dotted] (-2.7,0.56)--(-1,0.56);
\draw[-, dotted] (-2.7,0.3678)--(-1,0.3678);


\end{tikzpicture}
\begin{center}
$\varphi(x) = e^{x}$
\end{center}
\end{figure}

Taking $x_0=\ln(a^{p})$, $x_1=\ln(b^{q})$, and $t=\tfrac{1}{q}$, we obtain
$$
\frac{1}{p}a^{p} + \frac{1}{q}b^{q} \geq e^{(\tfrac{1}{p}\ln(a^{p}) + \tfrac{1}{q}\ln(b^{q}))}=ab,
$$
which is the desired inequality.

Verifying the equality case in Young's inequality is left as an exercise [Exercise \ref{E71}].
\end{proof}

\begin{theorem}[H\"older-Riesz inequality] \label{76} \index{inequality!H\"older-Riesz}
(a) Let $p,q \in (1,\infty)$ be such that $\frac{1}{p}+\frac{1}{q}=1$. If $f \in \mathcal{L}^{p}(X,\mathsf{S},\mu)$ and $g \in \mathcal{L}^{q}(X,\mathsf{S},\mu)$, then $fg \in \mathcal{L}^{1}(X,\mathsf{S},\mu)$ and
\begin{eqnarray} \label{F74}
\vert fg \vert_1 \leq \vert f \vert_{p}\,\vert g \vert_{q}.
\end{eqnarray}

(b) If $f \in \mathcal{L}^{1}(X,\mathsf{S},\mu)$ and $g \in \mathcal{L}^{\infty}(X,\mathsf{S},\mu)$, then $fg \in \mathcal{L}^{1}(X,\mathsf{S},\mu)$ and
\begin{eqnarray} \label{F75}
\vert fg \vert_1 \leq \vert f \vert_{1}\,\vert g \vert_{\infty}.
\end{eqnarray}
\end{theorem}

\begin{figure}[ht!]
\begin{minipage}[c]{0.5\textwidth}
\begin{center}
\includegraphics[scale=0.525]{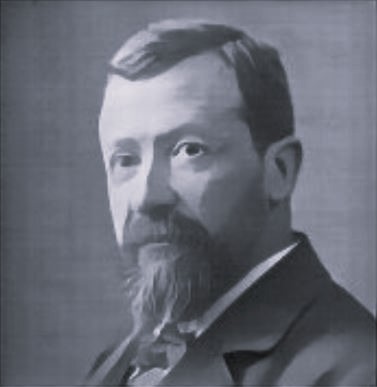} 
\begin{center}
Otto H\"older (1857-1937)
\end{center}
\end{center}
\end{minipage} \hfill \begin{minipage}[c]{0.5\textwidth}
\begin{center}
\includegraphics[scale=0.28]{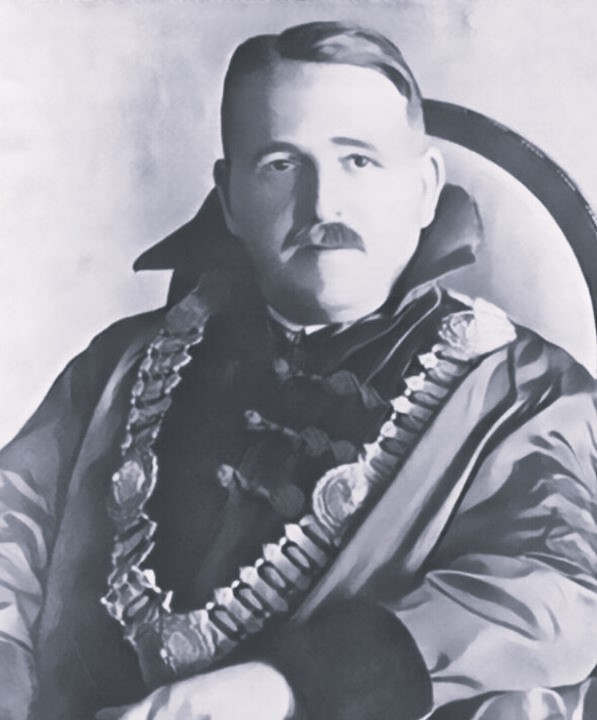} 
\begin{center}
Frigyes Riesz (1880-1956)
\end{center}
\end{center}
\end{minipage}
\end{figure}

\begin{proof}
\textit{(a):} Let $f \in \mathcal{L}^{p}(X,\mathsf{S},\mu)$ and $g \in \mathcal{L}^{q}(X,\mathsf{S},\mu)$. If $f=0$ or $g=0$ the statement is trivially satisfied. Suppose then that both functions are nonzero. By Proposition \ref{74} it follows that $|f|_p^p=\int_X |f|^{p}\,d\mu$ and $|g|_q^q=\int_X |g|^{q}\,d\mu$ are nonzero. 

For each $x \in X$, define the following positive real numbers
$$
a_x:=\frac{|f(x)|}{|f|_p}\quad \mbox{and}\quad  b_x:=\frac{|g(x)|}{|g|_q}\,.
$$

Applying Young's inequality to the pair of real numbers $a_x$ and $b_x$ for each $x \in X$, we obtain
\begin{eqnarray} \label{F76}
\frac{|f(x)g(x)|}{|f|_p|g|_q} \leq \frac{1}{p}\,\frac{|f(x)|^p}{|f|_p^p} + \frac{1}{q}\,\frac{|g(x)|^q}{|g|_q^q}.
\end{eqnarray}

Since the function on the right-hand side is integrable on $X$, inequality \eqref{F76} implies that $|fg|$ is integrable on $X$. Hence, $fg \in \mathcal{L}^{1}(X,\mathsf{S},\mu)$. Finally, applying the monotonicity of the integral to \eqref{F76}, we obtain
$$
\begin{aligned}
\frac{1}{|f|_p|g|_q}\int_X |fg|\,d\mu &\leq \frac{1}{p}\,\frac{1}{|f|_p^p}\int_X |f|^{p}\,d\mu + \frac{1}{q}\,\frac{1}{|g|_q^q}\int_X |g|^{q}\,d\mu = \frac{1}{p} + \frac{1}{q} = 1
\end{aligned}
$$
and, consequently,
$$
|fg|_1 \leq |f|_p|g|_q
$$
as claimed.

\textit{(b):} Since $f \in \mathcal{L}^{1}(X,\mathsf{S},\mu)$ and $g \in \mathcal{L}^{\infty}(X,\mathsf{S},\mu)$, the function $|f||g|_{\infty}$ is integrable on $X$ by Theorem \ref{67}. Using part \textit{(b)} of Theorem \ref{71}, we have $|fg| \leq |f||g|_{\infty}$ $\mu$ a.e. and, integrating this inequality, we obtain
\begin{eqnarray} \label{F77}
\int_{X} |fg| \,d\mu \leq \int_{X} |f||g|_{\infty}\,d\mu = |f|_{1}|g|_{\infty}.
\end{eqnarray}

Therefore, $fg \in \mathcal{L}^{1}(X,\mathsf{S},\mu)$.
\end{proof}

Observe that equality in (\ref{F74}) holds if and only if $\frac{1}{|f|_p|g|_q}\int_X |fg|\,d\mu=1$ if and only if $\left(\int_{X} |g|^{q}\,d\mu \right)|f|^{p} = \left(\int_{X} |f|^{p}\,d\mu \right)|g|^{q}$ $\mu$ a.e. 

On the other hand, in (\ref{F75}), $\int_{X} |fg|\,d\mu = |f|_{1}|g|_{\infty}$ if and only if $ |f|_{1}|g|_{\infty} - \int_{X} |fg|\,d\mu =0$ if and only if $\int_{X} |f|(|g|_{\infty}-|g|)\,d\mu =0$ if and only if $|f|(|g|_{\infty} - |g|)=0$ $\mu$ a.e. if and only if 
$$
\mu(\left\{x \in X\,:\, f(x) \neq 0\quad\mbox{and}\quad|g(x)|< |g|_{\infty} \right\})=0.
$$

There exists a generalization of the H\"older-Riesz inequality which we propose here as an interesting exercise [Exercise \ref{E712}]. 

From now on we shall adopt the convention that $\frac{1}{\infty}:=0$ so that the expression $\frac{1}{p}+\frac{1}{q}=1$ makes sense for every $p \in [1,\infty]$. The values $p$ and $q$ in $[1,\infty]$ satisfying the previous relation are called \textbf{conjugate exponents}. \index{conjugate exponents}

The inequality that we shall prove next, due to H. Minkowski\footnote{Hermann Minkowski (1864--1909) was a German mathematician who developed the geometric theory of numbers. His most important work was carried out in the areas of number theory, mathematical physics, and the theory of relativity. Minkowski investigated the arithmetic of quadratic forms in $n$ variables. His research in this field led him to consider the geometric properties of $n$-dimensional spaces. In 1896 he introduced his geometry of numbers, a geometric method for solving problems in number theory. In 1902 he joined the mathematics department at the University of Göttingen, where he collaborated closely with David Hilbert.} in 1896, provides us with the final ingredient needed to achieve the objective of this section.

\begin{figure}[ht!]
\centering
\includegraphics[scale=0.3]{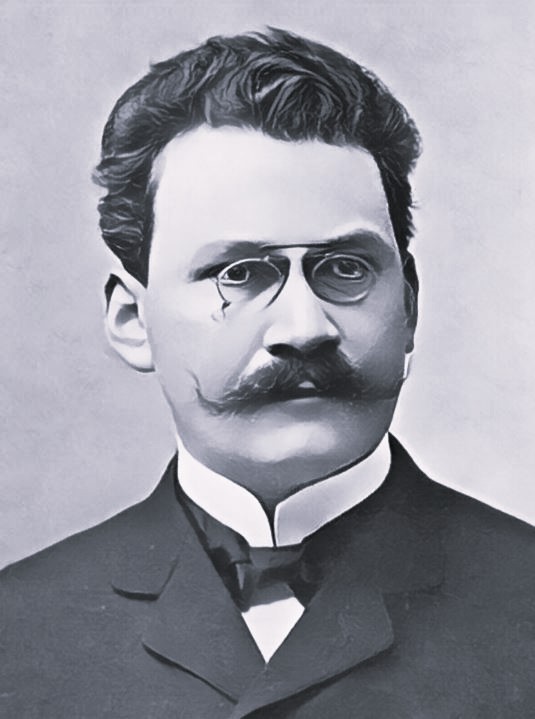} 
\begin{center}
Hermann Minkowski (1864-1909)
\end{center}
\end{figure}

\begin{theorem}[Minkowski's inequality] \label{77} \index{inequality!Minkowski's}
Let $p \in [1,\infty]$. If $f,g \in \mathcal{L}^{p}(X,\mathsf{S},\mu)$, then
\begin{eqnarray} \label{F78}
|f+g|_{p} \leq |f|_{p} + |g|_{p}.
\end{eqnarray}
\end{theorem}

\begin{proof}
We proceed by cases.

{\scshape Case 1.}\quad $p =1$.

Theorem \ref{73} ensures that $f+g \in \mathcal{L}^{1}(X,\mathsf{S},\mu)$. Applying the monotonicity of the integral to the inequality $|f+g| \leq |f| + |g|$, we obtain
$$
\int_X |f+g|\,d\mu \leq \int_X |f|\,d\mu + \int_X |g|\,d\mu,
$$
as claimed.

{\scshape Case 2.}\quad $p \in (1,\infty)$.

Theorem \ref{73} ensures that $f+g \in \mathcal{L}^{p}(X,\mathsf{S},\mu)$. Define $q:=\frac{p}{p-1} \in (1,\infty)$. Clearly, $p$ and $q$ are conjugate exponents and
$$
\int_X \left(|f+g|^{(p-1)} \right)^{q}\,d\mu = \int_X |f+g|^{p}\,d\mu  < +\infty.
$$ 

That is, $|f+g|^{p-1} \in \mathcal{L}^{q}(X,\mathsf{S},\mu)$ and $||f+g|^{p-1}|_q=|f+g|_p^{p-1}$.

Since
$$
|f+g|^{p}=|f+g||f+g|^{p-1} \leq \big(|f|+|g| \big)|f+g|^{p-1}  \leq |f||f+g|^{p-1} + |g||f+g|^{p-1},
$$
the monotonicity of the integral yields
\begin{eqnarray} \label{F79}
\int_X |f+g|^{p}\,d\mu \leq \int_X |f||f+g|^{p-1}\,d\mu + \int_X |g||f+g|^{p-1}\,d\mu. 
\end{eqnarray}

Therefore, applying the H\"older-Riesz inequality to the functions $f,g \in \mathcal{L}^{p}(X,\mathsf{S},\mu)$ and $|f+g|^{p-1} \in \mathcal{L}^{q}(X,\mathsf{S},\mu)$, we obtain from inequality (\ref{F79}) that
$$
\begin{aligned}
|f+g|_{p}^{p}&=\int_X |f+g|^{p}\,d\mu \leq \int_X |f||f+g|^{p-1}\,d\mu + \int_X |g||f+g|^{p-1}\,d\mu \\
& \leq |f|_p||f+g|^{p-1}|_q + |g|_p||f+g|^{p-1}|_q   = (|f|_p+|g|_p )\,|f+g|_p^{p-1}.
\end{aligned}
$$

If $|f+g|_{p}=0$, then inequality (\ref{F78}) is trivially satisfied. Hence, if $|f+g|_{p} >0$, dividing the previous inequality by $|f+g|_p^{p-1}>0$, we obtain
$$
|f+g|_p = |f+g|^{p}_p |f+g|^{1-p}_p \leq |f|_p + |g|_p,
$$
which is the desired inequality.

{\scshape Case 3.}\quad $p =\infty$.

From Theorem \ref{71} we conclude that $|f+g| \leq |f|_{\infty} + |g|_{\infty}$ $\mu$ a.e. and, therefore, $|f +g|_{\infty} \leq |f|_{\infty} + |g|_{\infty}$, which proves (\ref{F78}). 
\end{proof}

Before proving the result that fulfills the main objective of this section, we give the following definition.

\begin{definition} \label{78} \index{space!seminormed} \index{seminorm}
Let $V$ be a vector space over $\mathbb{R}$. A \textbf{seminorm} on $V$ is a function $\rho:V \to \mathbb{R}$ satisfying the following properties:
\begin{itemize}
\item[\rm (SN1)] $\rho(v) \geq 0$ for every $v \in V$.
\item[\rm (SN2)] $\rho(\alpha\,v)=|\alpha|\,\rho(v)$ for all $v \in V$, $\alpha \in \mathbb{R}$.
\item[\rm (SN3)] $\rho(v+w) \leq \rho(v) + \rho(w)$ for all $v,w \in V$.
\end{itemize}

A \textbf{seminormed space} is a vector space $V$ endowed with a seminorm $\rho(\cdot)$.
\end{definition}

\begin{theorem} \label{79}
$\mathcal{L}^{p}(X,\mathsf{S},\mu)=(\mathcal{L}^{p}(X,\mathsf{S},\mu), \vert \cdot \vert_p)$ is a seminormed space for every $p \in [1,\infty]$.
\end{theorem}

\begin{proof}
Theorem \ref{73} ensures that $\mathcal{L}^{p}(X,\mathsf{S},\mu)$ is a vector space over $\mathbb{R}$ for every $p \in [1,\infty]$, so it only remains to prove that the function $| \cdot |_p$ is a seminorm for every $p \in [1,\infty]$. It follows immediately from Theorems \ref{65} and \ref{71} that $| \cdot |_p$ satisfies property (SN1) for every $p \in [1,\infty]$.

Now let $\gamma \in \mathbb{R}$ and $f \in \mathcal{L}^{p}(X,\mathsf{S},\mu)$. If $p \in [1,\infty)$, it is clear that $|\gamma\,f|_{p}^{p}=|\gamma|^{p}|f|_{p}^{p}$ by Theorem \ref{67}.

If $p=\infty$, since $|f|\leq |f|_{\infty}$ $\mu$ a.e., then $|\gamma\, f| \leq |\gamma||f|_{\infty}$ $\mu$ a.e. and, therefore, $|\gamma\,f|_{\infty} \leq  |\gamma||f|_{\infty}$. 

Conversely, suppose that $\gamma \neq 0$ and $f \neq 0$. Since $|\gamma\,f|  <  |\gamma|\,c$ $\mu$ a.e. for every $c \geq 0$ such that $\mu\left(\{x\in X\,:\, |f(x)|>c \}\right)=0$, then $|\gamma||f| \leq |\gamma\,f|_{\infty}$ $\mu$ a.e. by part \textit{(b)} of Theorem \ref{71}. Thus, $|f| \leq \frac{|\gamma\,f|_{\infty}}{|\gamma|}$ $\mu$ a.e. and it follows that $|\gamma||f|_{\infty} \leq |\gamma\,f|_{\infty}$. Consequently, $|\gamma f|_{\infty} = |\gamma||f|_{\infty}$. The case when $\gamma=0$ or $f=0$ is trivially satisfied.

This last statement together with Theorem \ref{77} ensures that $|\cdot|_{p}$ is a seminorm for every $p \in [1,\infty]$.
\end{proof}

In the next section we shall study the difficulty that arises when trying to extend $| \cdot |_p$ to a norm. For the moment, let us look at some interesting examples of spaces $\mathcal{L}^{p}(X,\mathsf{S},\mu)$ and some relations between them.

Consider the measure space given by $X=\mathbb{N}$, $\mathsf{S}=\mathcal{P}(\mathbb{N})$, and $\mu=\mu^{\sharp}$ the counting measure. Clearly, in this space every function $f:\mathbb{N} \to \mathbb{R}$ is $\mathcal{P}(\mathbb{N})$-measurable. 

Let $f:\mathbb{N} \to \mathbb{R}$ be a function. From Example \ref{515} and Exercise \ref{E612}, we conclude that if $p \in [1,\infty)$, then
$$
f \in \mathcal{L}^{p}(\mathbb{N},\mathcal{P}(\mathbb{N}),\mu^{\sharp}) \quad \Longleftrightarrow \quad \sum_{k=1}^{\infty} |f(k)|^{p} < +\infty.
$$

Thus, for each $p \in [1,\infty)$,
$$
\mathcal{L}^{p}(\mathbb{N},\mathcal{P}(\mathbb{N}),\mu^{\sharp})=\left\{ f:\mathbb{N} \to \mathbb{R}\,:\, \sum_{k=1}^{\infty}|f(k)|^{p} < +\infty \right\}.
$$

The advantage of the previous result and observations is that they will allow us to derive a wide variety of results concerning series of real numbers from those already proved for integrals. 

Let us now consider the case when $p=\infty$. First observe that the Well-Ordering Principle ensures that the only $\mu^{\sharp}$-null subset of $\mathbb{N}$ is the empty set. Thus, if $f: \mathbb{N} \to \mathbb{R}$ is bounded $\mu^{\sharp}$ a.e., then $f:\mathbb{N} \to \mathbb{R}$ is bounded, and conversely. Therefore,
$$
\mathcal{L}^{\infty}(\mathbb{N},\mathcal{P}(\mathbb{N}),\mu^{\sharp})=\left\{ f:\mathbb{N} \to \mathbb{R}\,:\, \mbox{there exists}\,\,c \in \mathbb{R}^+\,\,\mbox{such that}\,\,|f(k)|< c\,\,\mbox{for every}\,\,k \in \mathbb{N}\right\}.
$$

Next, we present the definition of the previous spaces in the classical context of a basic course in mathematical analysis, where they are known in many texts as \textit{discrete Lebesgue spaces}. This terminology is justified by the previous observations. \index{space!discrete Lebesgue}

\begin{definition}\label{710} \index{space!ellp@$\ell^{p}(\mathbb{R})$} \index{space!ellinf@$\ell^{\infty}(\mathbb{R})$}
If $p \in [1,\infty)$, the set of all sequences of real numbers $(x_k)$ such that the series 
$$
\sum_{k=1}^{\infty} |x_k|^{p}
$$
converges is denoted by $\ell^{p}(\mathbb{R})$. This set is precisely the space $\mathcal{L}^{p}(\mathbb{N},\mathcal{P}(\mathbb{N}),\mu^{\sharp})$.

Similarly, if $p=\infty$, the set of all bounded sequences of real numbers $(x_k)$ is denoted by $\ell^{\infty}(\mathbb{R})$, which coincides with the space $\mathcal{L}^{\infty}(\mathbb{N},\mathcal{P}(\mathbb{N}),\mu^{\sharp})$.
\end{definition}

Let us now look at some examples.

\begin{example} \label{711}
Let $X=\mathbb{N}$, $\mathsf{S}=\mathcal{P}(\mathbb{N})$ and $\mu=\mu^{\sharp}$ be the counting measure. Define the function $f:\mathbb{N} \to \mathbb{R}$ by $f(k):=\frac{1}{k}$. It is immediate that $f$ is bounded since $|f(k)| \leq 1$ for all $k \in \mathbb{N}$, and therefore $f \in \mathcal{L}^{\infty}(\mathbb{N},\mathcal{P}(\mathbb{N}),\mu^{\sharp})$. 

On the other hand, notice that
$$
\int_{\mathbb{N}} |f|\,d\mu^{\sharp} = \sum_{k=1}^{\infty} \frac{1}{k} = +\infty
$$
and therefore $f \not\in \mathcal{L}^{1}(\mathbb{N},\mathcal{P}(\mathbb{N}),\mu^{\sharp})$. However, $f \in \mathcal{L}^{2}(\mathbb{N},\mathcal{P}(\mathbb{N}),\mu^{\sharp})$ since
$$
\int_{\mathbb{N}} |f|^2\,d\mu^{\sharp} = \sum_{k=1}^{\infty} \frac{1}{k^2}=\frac{\pi^2}{6} < +\infty.
$$
\end{example}

\begin{example} \label{712}
Let $(X,\mathsf{S},\mu)$ be a measure space such that $\mu(X)=+\infty$. It is clear that the characteristic function $\chi_{X}$ is bounded and, therefore, $\chi_{X} \in \mathcal{L}^{\infty}(X,\mathsf{S},\mu)$. However, $\chi_{X} \not \in \mathcal{L}^{p}(X,\mathsf{S},\mu)$ for every $p \in [1,\infty)$ since
$$
\int_{X} |\chi_{X}|^{p}\,d\mu =\mu(X) = +\infty.
$$
\end{example}

A common feature of the spaces appearing in the previous examples is that they have infinite measure. The following proposition establishes a sufficient condition for an inclusion relation between the spaces $\mathcal{L}^{p}$.

\begin{theorem} \label{713}
Let $(X,\mathsf{S},\mu)$ be a finite measure space and let $1 \leq p < q \leq \infty$. Then $\mathcal{L}^{q}(X,\mathsf{S},\mu) \subset \mathcal{L}^{p}(X,\mathsf{S},\mu)$. Moreover, for every $f \in \mathcal{L}^{q}(X,\mathsf{S},\mu)$,
\begin{itemize}
\item[(a)] $\vert f \vert_{p} \leq \mu(X)^{(q-p)/qp}\vert f \vert_{q}$ if $q \in [1,\infty)$.
\item[(b)] $\vert f \vert_{p} \leq \mu(X)^{1/p}\vert f \vert_{\infty}$ if $q=\infty$.
\end{itemize}
\end{theorem}

\begin{proof}
Since $\mu(X) < +\infty$, it follows that $\chi_{X} \in \mathcal{L}^{s}(X,\mathsf{S},\mu)$ for every $s \in [1,\infty)$ and
$$
\vert \chi_{X} \vert_{s} = \left(\int_{X} |\chi_{X}|^s\,d\mu \right)^{1/s} = \mu(X)^{1/s}.
$$

\textit{(a):} Let $1 \leq p < q < \infty$ and let $f \in \mathcal{L}^{q}(X,\mathsf{S},\mu)$. Notice that
$$
\int_{X} |f|^{q}\,d\mu =\int_{X} |f|^{p(q/p)}\,d\mu= \int_{X} ||f|^{p}|^{q/p}\,d\mu
$$ 
and therefore $|f|^{p} \in \mathcal{L}^{q/p}(X,\mathsf{S},\mu)$ and
$$
\vert |f|^{p} \vert_{q/p} = \left(\int_{X} |f|^{p(q/p)}\,d\mu \right)^{p/q} = \left(\int_{X} |f|^{q}\,d\mu \right)^{p/q} =\vert f \vert_{q}^{p}.
$$

On the other hand, the function $\chi_{X} \in\mathcal{L}^{q/(q-p)}(X,\mathsf{S},\mu)$. By the H\"older-Riesz inequality it follows that $\chi_{X}|f|^{p}  = |f|^{p}$ is integrable and
$$
\int_{X} |f|^{p} \,d\mu \leq \vert \chi_{X} \vert_{q/(q-p)} \vert |f|^{p} \vert_{q/p}  = \mu(X)^{(q-p)/q}  \vert f \vert_{q}^{p}.
$$

Consequently, $f \in \mathcal{L}^{p}(X,\mathsf{S},\mu)$ and
\begin{eqnarray} \label{F710}
\vert f \vert_{p} \leq \mu(X)^{(q-p)/pq}  \vert f \vert_{q}.
\end{eqnarray}

\textit{(b):} Let $1 \leq p < \infty$ and let $f \in \mathcal{L}^{\infty}(X,\mathsf{S},\mu)$. Then $|f|^{p} \leq \chi_{X} \vert f \vert_{\infty}^{p}$ $\mu$ a.e. Since $\mu(X)<+\infty$, the function $\chi_{X} \vert f \vert_{\infty}^{p}$ is integrable and, therefore, $f \in \mathcal{L}^{p}(X,\mathsf{S},\mu)$. Now, integrating the previous inequality we obtain $\vert f \vert_{p}^{p} \leq \mu(X) \vert f \vert_{\infty}^{p}$, that is,
\begin{eqnarray} \label{F711}
\vert f \vert_{p} \leq \mu(X)^{1/p}  \vert f \vert_{\infty}.
\end{eqnarray}

This completes the proof.
\end{proof}

Some additional inclusion relations between Lebesgue spaces can be found in the exercises section [Exercises \ref{E77}, \ref{E78} and \ref{E79}].

\section{The spaces \texorpdfstring{$L^{p}$}{}}

Let us begin by recalling the following concept from previous courses.

\begin{definition} \label{714} \index{space!normed} \index{norm}
Let $V$ be a vector space over $\mathbb{R}$. A \textbf{norm} on $V$ is a function $\Vert \cdot \Vert: V \to \mathbb{R}$ satisfying the following properties:
\begin{itemize}
\item[\rm (N1)] $\Vert v \Vert=0$ if and only if $v=0_V$.
\item[\rm (N2)] $\Vert \alpha\,v \Vert=|\alpha|\,\Vert v \Vert$ for every $v \in V$ and every $\alpha \in \mathbb{R}$.
\item[\rm (N3)] $\Vert v + w \Vert \leq \Vert v \Vert + \Vert w \Vert$ for every $v,w \in V$.
\end{itemize}

A \textbf{normed space} is a vector space $V$ together with a given norm $\Vert \cdot \Vert$.
\end{definition}

Let $(X,\mathsf{S},\mu)$ be a measure space. Theorem \ref{79} ensures that the functions defined in (\ref{F71}) and (\ref{F72}) satisfy properties (N2) and (N3) of Definition \ref{714} on the spaces $\mathcal{L}^{p}(X,\mathsf{S},\mu)$ for every $p \in [1,\infty]$. Moreover, Proposition \ref{74} ensures that $\vert f \vert_{p}=0$ if and only if $f=0$ $\mu$ a.e. for every $p \in [1,\infty]$. Let us examine a concrete example.

\begin{example}\label{715}
Let $(X,\mathsf{S},\mu)$ be a measure space such that $\mu(X)< +\infty$.

For every subset $A \in \mathsf{S}$ it is clear that $\chi_{A} \in \mathcal{L}^{\infty}(X,\mathsf{S},\mu)$ and
$$
\{x \in X \,:\, \chi_{A}(x)>c \}=\left\{
\begin{array}{c c l}
\varnothing & & \mbox{if}\,\, c \geq 1,\\
A & & \mbox{if}\,\, 0 \leq c <1,\\
X & & \mbox{if}\,\, c<0.
\end{array}
\right.
$$

Thus, if $A \in \mathsf{S}$ is such that $\mu(A)=0$, then $\vert \chi_{A} \vert_{\infty}=0$. On the other hand, if $\mu(A)>0$, then $\vert \chi_{A} \vert_{\infty}=1$.

Now, $\chi_{A} \in \mathcal{L}^{p}(X,\mathsf{S},\mu)$ for every $p \in [1,\infty)$ and every $A \in \mathsf{S}$. Moreover,
$$
\vert \chi_{A} \vert_{p}^{p} =\int_{X} |\chi_{A}|^{p}\,d\mu = \mu(A)\quad \text{for every}\,\,A \in \mathsf{S}.
$$

Choose a nonempty set $A \in \mathsf{S}$ such that $\mu(A)=0$. In this case we have $\vert \chi_{A} \vert_{p}=0$ for every $p \in [1,\infty]$, but $\chi_{A}(x) \neq 0 $ for every $x \in A$.
\end{example}

The previous discussion shows that the functions defined in (\ref{F71}) and (\ref{F72}) cannot define a norm on the spaces $\mathcal{L}^{p}(X,\mathsf{S},\mu)$ for $p \in [1,\infty]$ since property (N1) is not satisfied. Thus, the objective of this section is to endow the expression $\vert \cdot \vert_{p}$ with the structure of a norm. To achieve this, we proceed as follows:

We define on $\mathcal{L}^{p}(X,\mathsf{S},\mu)$ for $p \in [1,\infty]$ the following equivalence relation $\sim_{\mu}$ [Exercise \ref{E74}]:
\begin{eqnarray} \label{F712}
f \sim_{\mu} g \quad \Longleftrightarrow\quad f=g \quad\mu\text{ a.e.}. 
\end{eqnarray}

Let $[f]_{\mu}:=\left\{ g \in \mathcal{L}^{p}(X,\mathsf{S},\mu) \,:\, f \sim_{\mu} g \right\}$ denote the equivalence class of $f \in \mathcal{L}^{p}(X,\mathsf{S},\mu)$ relative to $\mu$. The set of equivalence classes on $\mathcal{L}^{p}(X,\mathsf{S},\mu)$ relative to $\sim_{\mu}$ is denoted by
\begin{eqnarray} \label{F713} \index{space!Lpa@$L^{p}(X,\mathsf{S},\mu)$} \index{space!Lpa1@$L^{\infty}(X,\mathsf{S},\mu)$}
L^{p}(X,\mathsf{S},\mu):= \mathcal{L}^{p}(X,\mathsf{S},\mu)\,\diagup\,\sim_{\mu}.
\end{eqnarray}

Notice that if $f_{i}=g_{i}$ $\mu$ a.e. on $X$, $i=1,2$, then $\gamma_{1} f_{1} + \gamma_{2} f_{2} = \gamma_{1} g_{1} + \gamma_{2} g_{2}$ $\mu$ a.e. on $X$ for every $\gamma_{1},\gamma_{2} \in \mathbb{R}$. Consequently, the vector space structure of $\mathcal{L}^{p}(X,\mathsf{S},\mu)$ induces a vector space structure on $L^{p}(X,\mathsf{S},\mu)$ through the operations
$$
[\gamma f]_{\mu}:=\gamma [f]_{\mu}\quad \mbox{and}\quad [f+g]_{\mu}:=[f]_{\mu}+[g]_{\mu},
$$
which are well defined by Proposition \ref{74}.

On the space $L^{p}(X,\mathsf{S},\mu)$ we define $\Vert \cdot \Vert_{p}:L^p(X,\mathsf{S},\mu)\to \mathbb{R}$ by \index{norm!on Lp@on $L^{p}(X,\mathsf{S},\mu)$} \index{norm!on Linfinity@on $L^{\infty}(X,\mathsf{S},\mu)$}
\begin{eqnarray} \label{F714}
\Vert [f]_{\mu} \Vert_{p}:=\left( \int_{X} |f|^{p}\,d\mu \right)^{1/p},\qquad p \in [1,\infty);
\end{eqnarray}
\begin{eqnarray} \label{F715} 
\Vert [f]_{\mu} \Vert_{\infty}:=\inf\left\{ c \geq 0\,:\,\mu\left( \{x \in X\,:\, |f(x)|> c \} \right)=0 \right\}.
\end{eqnarray}
whose values are independent of the representative of the equivalence class [Exercise \ref{E75}]. Once again, Proposition \ref{74} together with Theorem \ref{79} ensures that the function $\Vert \cdot \Vert_{p}$ is a norm on the vector space $L^{p}(X,\mathsf{S},\mu)$ for every $p \in [1,\infty]$. Therefore, $(L^{p}(X,\mathsf{S},\mu), \Vert \cdot \Vert_{p})$ is a normed space for every $p \in [1,\infty]$.

\begin{remark} \label{716}
In what follows, we shall identify each function $f$ with its equivalence class modulo $\sim_\mu$, and we shall denote it simply by $f:X\to\mathbb{R}$. Under this convention, the space $L^{p}(X,\mathsf S,\mu)$, for $1\leq p\leq \infty$, may be understood as a space of functions in which two functions are regarded as equal whenever they differ only on a set of measure zero.

Likewise, every function defined on $X$ except on a set of measure zero will be regarded as defined on all of $X$ by assigning, for instance, the value $0$ on the set where it was not originally defined.
\end{remark}

Next, we prove that the normed spaces $(L^{p}(X,\mathsf{S},\mu), \Vert \cdot \Vert_{p})$ are Banach spaces for every $p \in [1,\infty]$, that is, spaces in which every Cauchy sequence converges. The proof of this fact is due to the mathematicians Frigyes Riesz and Ernst Fischer\footnote{Ernst Sigismund Fischer (1875--1954) was an Austrian mathematician who worked alongside Franz Mertens and Hermann Minkowski at the Universities of Vienna and Zürich. He later became a professor at the University of Erlangen, where he collaborated with Emmy Noether. His research area was mathematical analysis, where he became interested in orthonormal sequences of functions, a topic that laid the foundations for the theory of Hilbert spaces. In 1907, independently of Frigyes Riesz, he proved the Riesz--Fischer theorem, which still bears their names.}, and we require the following results.

\begin{lemma} \label{717}
Let $(V,\Vert \cdot \Vert)$ be a normed space and $(v_k)$ a Cauchy sequence in $V$. If some subsequence of $(v_k)$ converges to $v$ in $V$, then $(v_k)$ converges to $v$ in $V$.
\end{lemma}

The proof is left as an exercise [Exercise \ref{E73}].

\begin{lemma}[Cauchy Condition] \label{718} \index{condition!Cauchy}
Let $(V,\Vert \cdot \Vert)$ be a normed space, $(v_k)$ a Cauchy sequence in $V$, and $\delta>0$ fixed. Then there exists a subsequence $(v_{k_j})$ of $(v_k)$ such that:
$$
\Vert v_{k_{j+1}}-v_{k_j} \Vert < \delta^{j},\qquad \forall\,j\in \mathbb{N}.
$$
\end{lemma}

\begin{proof}
Since $(v_k)$ is a Cauchy sequence, there exists $k_{1}:=k_{\delta} \in \mathbb{N}$ such that $\Vert v_{k}-v_{i} \Vert < \delta$ for all $k,i \geq k_{1}$. Similarly, there exists $k_{\delta,2} \in \mathbb{N}$ such that $\Vert v_{k}-v_{i} \Vert < \delta^{2}$ for all $k,i \geq k_{\delta,2}$. Defining $k_{2}:=\max\{k_{1},k_{\delta,2}\}+1 \in \mathbb{N}$ it follows that $\Vert v_{k}-v_{i} \Vert < \delta^{2}$ for all $k,i \geq k_{2}$. Proceeding in this way, for each $j \in \mathbb{N}$ we obtain a natural number $k_{j}>k_{j-1}$ such that $\Vert v_{k}-v_{i} \Vert \leq \delta^{j}$ for all $k,i\geq k_{j}$. Thus, by induction we construct a subsequence satisfying the required property.
\end{proof}

\begin{theorem}[Dominated Convergence Theorem in $L^{p}$] \label{719} \index{theorem!dominated convergence Lp}
Let $p \in [1,\infty)$ and $(f_k)$ be a sequence in $L^{p}(X,\mathsf{S},\mu)$ such that $f_k\to f$ $\mu$ a.e. If there exists $g \in L^{p}(X,\mathsf{S},\mu)$ such that $|f_k| \leq g$ $\mu$ a.e. for all $k \in \mathbb{N}$, then $f \in L^{p}(X,\mathsf{S},\mu)$ and $\lim_{k \to \infty}\Vert f_k-f\Vert_{p}=0$.
\end{theorem}

\begin{proof}
Note that $|f_k|^{p} \in L^{1}(X,\mathsf{S},\mu)$ and $|f_k|^{p} \leq |g|^{p}$ for all $k \in \mathbb{N}$ $\mu$ a.e. Since $|f_k|^{p} \to |f|^{p}$ $\mu$ a.e., by the Lebesgue dominated convergence theorem we obtain that $|f|^{p} \in L^{1}(X,\mathsf{S},\mu)$ and, hence, $f \in L^{p}(X,\mathsf{S},\mu)$. Now, since $|f_k-f|^{p} \leq (|f_k|+|f|)^{p} \leq (2g)^{p}$ for all $k \in \mathbb{N}$ $\mu$ a.e., it follows that $|f_k-f|^{p} \in L^{1}(X,\mathsf{S},\mu)$, and by the Lebesgue dominated convergence theorem we conclude that $\lim_{k \to \infty}\Vert f_k-f\Vert_{p}=0$.
\end{proof}

\begin{theorem}[Riesz--Fischer] \label{720} \index{theorem!Riesz-Fischer}
$\left(L^{p}(X,\mathsf{S},\mu), \Vert \cdot \Vert_{p}\right)$ is a Banach space for all $p \in [1,\infty]$.
\end{theorem}

\begin{figure}[ht!]
\begin{minipage}[c]{0.5\textwidth}
\begin{center}
\includegraphics[scale=0.3]{riez.jpeg} 
\begin{center}
Frigyes Riesz (1880-1956)
\end{center}
\end{center}
\end{minipage} \hfill \begin{minipage}[c]{0.5\textwidth}
\begin{center}
\includegraphics[scale=0.3]{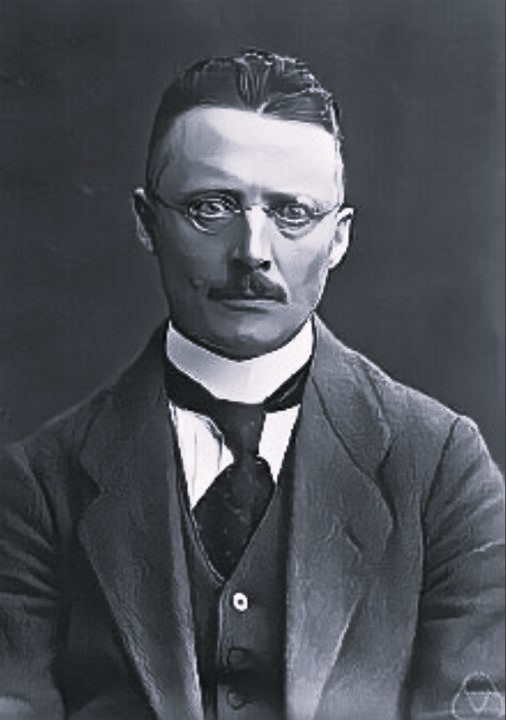} 
\begin{center}
Ernst Fischer (1875-1954)
\end{center}
\end{center}
\end{minipage}
\end{figure}

\begin{proof}
We proceed by cases.

{\scshape Case 1.}\quad $p \in [1,\infty)$.

Let $(f_k)$ be a Cauchy sequence in $(L^{p}(X,\mathsf{S},\mu), \Vert \cdot \Vert_{p})$. By Lemma \ref{718}, for $\delta:=\frac{1}{2}>0$, there exists a subsequence $(f_{k_j})$ of $(f_k)$ such that
$$
\left\Vert f_{k_{j+1}}-f_{k_{j}} \right\Vert_{p} \leq \frac{1}{2^j} \quad\forall\,j \in \mathbb{N}.
$$

Then,
$$
\sum_{j=1}^{\infty} \left\Vert f_{k_{j+1}}-f_{k_{j}} \right\Vert_{p} \leq \sum_{j=1}^{\infty} \frac{1}{2^j}=1.
$$

For each $m \in \mathbb{N}$, define the functions $g_{m},g \in \mathbb{M}^{+}(X, \mathsf{S})$ by
$$
\begin{aligned}
g_{m}&:=|f_{k_{1}}| + \sum_{j=1}^{m-1} |f_{k_{j+1}}-f_{k_{j}}|\\
g&:=|f_{k_{1}}|+ \sum_{j=1}^{\infty} |f_{k_{j+1}}-f_{k_{j}}|=\sup_{m \in \mathbb{N}}g_m(x).
\end{aligned}
$$

By Minkowski's inequality we have
$$
\Vert g_{m} \Vert_{p} \leq \Vert f_{k_{1}} \Vert_{p} + \sum_{j=1}^{m-1} \Vert f_{k_{j+1}} - f_{k_{j}} \Vert_{p} \leq \Vert f_{k_{1}} \Vert_{p} + \sum_{j=1}^{m-1} \frac{1}{2^{j}} \leq \Vert f_{k_{1}} \Vert_{p} + 1 \quad \forall m \in \mathbb{N}
$$
and therefore $g_{m} \in L^{p}(X,\mathsf{S},\mu)$ for all $m \in \mathbb{N}$.

By Fatou's lemma, we obtain
$$
\left(\int_{X} |g|^{p}\,d\mu\right)^{\frac{1}{p}} \leq \left( \liminf_{m \to \infty} \int_{X} |g_{m}|^{p}\,d\mu \right)^{\frac{1}{p}}
\leq \liminf_{m \to \infty} \left( \int_{X} |g_{m}|^{p}\,d\mu\right)^{\frac{1}{p}} \leq \Vert f_{k_1}\Vert_{p} +1 < +\infty.
$$

That is, $g \in L^{p}(X,\mathsf{S},\mu)$. There exists a null set $N \subset X$ such that [Exercise \ref{E518}]
$$
g(x)=|f_{k_{1}}(x)|+ \sum_{j=1}^{\infty} |f_{k_{j+1}}(x)-f_{k_{j}}(x)|<+\infty \quad\forall \,x \in X\smallsetminus N.
$$

Hence, the series
$$
f_{k_{1}}(x) + \sum_{j=1}^{\infty}(f_{k_{j+1}}(x)-f_{k_{j}}(x))
$$
converges in $\mathbb{R}$ for every $x \in X\smallsetminus N$. We define $f:X \to \mathbb{R}$ by
$$
f(x):=\left\{
\begin{array}{ll}
f_{k_{1}}(x) + \sum_{j=1}^{\infty} (f_{k_{j+1}}(x)-f_{k_{j}}(x)) & \text{if } x \in X\smallsetminus N\\
0 & \text{if } x \in N.
\end{array}
\right.
$$

Since
$$
f_{k_{1}}+\sum_{j=1}^{m-1}\big(f_{k_{j+1}}-f_{k_{j}}\big)=f_{k_m}\qquad \forall m \in \mathbb{N},
$$
it follows that
$$
 f_{k_{m}}\to f \quad \mu\text{ a.e.} \quad \text{and}\quad |f_{k_{m}}| \leq g_{m}\leq g \quad \forall m \in \mathbb{N}, \quad \mu\text{ a.e.}
$$

By Theorem \ref{719} it follows that $f \in L^{p}(X,\mathsf{S},\mu)$ and $\lim_{j \to \infty} \Vert f-f_{k_j}\Vert_{p}=0$. That is, the subsequence $(f_{k_j})$ converges to $f$ in $L^{p}(X,\mathsf{S},\mu)$. Since $(f_k)$ is Cauchy, we conclude from Lemma \ref{717} that $f_k \to f$ in $L^{p}(X,\mathsf{S},\mu)$.

{\scshape Case 2.}\quad $p=\infty$.

Let $(f_k)$ be a Cauchy sequence in $\left(L^{\infty}(X,\mathsf{S},\mu), \Vert \cdot \Vert_{\infty}\right)$. For each $m \in \mathbb{N}$, there exists $k_m \in \mathbb{N}$ such that
$$
\Vert f_j - f_k\Vert_{\infty} \leq \frac{1}{m} \qquad \forall\,j,k \geq k_m.
$$

Let $N_{j,k}$ be a $\mu$-null subset of $X$ such that (see Theorem \ref{71})
$$
|f_j(x)-f_k(x)|\leq \Vert f_j - f_k\Vert_{\infty}\qquad \forall x \in X\smallsetminus N_{j,k}.
$$

Then the set $N:=\bigcup_{j,k=1}^{\infty} N_{j,k}$ is $\mu$-null and, for each $x \in X\smallsetminus N$, we have
\begin{eqnarray} \label{F716}
|f_{j}(x)-f_{k}(x)|\leq \frac{1}{m}\qquad \forall\,j,k \geq k_m.
\end{eqnarray}

Consequently, for each $x \in X\smallsetminus N$, $(f_k(x))$ is a Cauchy sequence in $\mathbb{R}$ and therefore, for each $x \in X\smallsetminus N$, there exists $f(x) \in \mathbb{R}$ such that $f_k(x)\to f(x)$ in $\mathbb{R}$. We define $f$ on $X$ as follows
$$
f(x):=\left\{ 
\begin{array}{ll}
f(x) & \text{if } x \in X\smallsetminus N,\\
0 & \text{if } x \in N. 
\end{array}
\right.
$$

Fix $k \in \mathbb{N}$. In (\ref{F716}), taking the limit as $j \to \infty$ we obtain
$$
|f(x)-f_k(x)| \leq \frac{1}{m}\quad\forall k \geq k_m,\quad\forall x \in X\smallsetminus N.
$$

Hence, $f-f_k \in L^{\infty}(X,\mathsf{S},\mu)$ and
\begin{eqnarray} \label{F717}
\Vert f-f_k \Vert_{\infty} \leq \frac{1}{m} \qquad \forall\,k \geq k_m.
\end{eqnarray}

Since $f=(f-f_{k_1})+f_{k_1}$, it follows that $f \in L^{\infty}(X,\mathsf{S},\mu)$ and from (\ref{F717}) we conclude that $f_k \to f$ in $L^{\infty}(X,\mathsf{S},\mu)$.

This completes the proof.
\end{proof}

Within the proof of the Riesz--Fischer theorem, we have implicitly established the following result.

\begin{corollary} \label{721}
Let $(f_k)$ be a sequence in $L^{p}(X,\mathsf{S},\mu)$ such that $f_k \to f$ in $L^{p}(X,\mathsf{S},\mu)$.
\begin{itemize}
\item[(a)] If $p\in [1,\infty)$, there exists a subsequence $(f_{k_j})$ of $(f_k)$ and a function $g \in L^{p}(X,\mathsf{S},\mu)$ such that
$$
\begin{aligned}
&f_{k_j}\to f \quad \mu \text{ a.e.,}\\
&|f_{k_j}| \leq g \quad \mu \text{ a.e.},\quad \forall j \in \mathbb{N}.
\end{aligned}
$$
\item[(b)] If $p=\infty$, then $f_k(x) \to f(x)$ $\mu$ a.e. and there exists $c \in \mathbb{R}$ such that $|f_k|\leq c$ $\mu$ a.e., for all $k \in \mathbb{N}$.
\end{itemize}
\end{corollary}

\begin{proof}
\textit{(a):} Since $(f_k)$ is a Cauchy sequence in $L^{p}(X,\mathsf{S},\mu)$, following the proof of Theorem \ref{720} we conclude that there exist a subsequence $(f_{k_j})$ of $(f_k)$ and functions $\tilde{f},g \in L^{p}(X,\mathsf{S},\mu)$ such that
$$
\begin{aligned}
&f_{k_j} \to \tilde{f} \quad \mu \text{ a.e.,}\\
&|f_{k_j}| \leq g \quad \mu \text{ a.e.},\quad\forall j \in \mathbb{N}.\\
\end{aligned}
$$

Thus, $f_k \to \tilde{f}$ in $L^{p}(X,\mathsf{S},\mu)$ and, since $f_k \to f$ in $L^{p}(X,\mathsf{S},\mu)$, it follows that $f=\tilde{f}$.

\textit{(b):} Theorem \ref{720} ensures that for each $k \in \mathbb{N}$ there exists a $\mu$-null set $N_k$ such that
$$
|f_k(x)| \leq \Vert f_k \Vert_{\infty}\quad \text{and}\quad |f(x)-f_k(x)| \leq \Vert f-f_k\Vert_{\infty} \qquad\forall\,x \in X\smallsetminus N_k.
$$

Let $N:=\bigcup_{k=1}^{\infty} N_k$. Then
$$
|f_k(x)| \leq \Vert f_k \Vert_{\infty}\quad \text{and}\quad |f(x)-f_k(x)| \leq \Vert f-f_k\Vert_{\infty}\quad \forall\,x \in X\smallsetminus N,\quad \forall\, k \in \mathbb{N}.
$$

Since $f_k \to f$ in $L^{\infty}(X,\mathsf{S},\mu)$, we conclude that $f_k(x) \to f(x)$ for all $x \in X\smallsetminus N$. Moreover, since $(f_k)$ is bounded in $L^{\infty}(X,\mathsf{S},\mu)$, there exists $c \in \mathbb{R}$ such that $|f_k(x)|\leq \Vert f_k\Vert_{\infty} \leq c$ for all $x \in X\smallsetminus N$ and all $k \in \mathbb{N}$, as stated.
\end{proof}

The following result follows directly from the Riesz--Fischer theorem and Example \ref{515}.

\begin{corollary} \label{722}
The space of real sequences $\ell^{p}(\mathbb{R})$ is a Banach space for all $p \in [1,\infty]$.
\end{corollary}

We conclude this section by constructing the spaces $L^{p}_{\mathbb{C}}(X,\mathsf{S},\mu)$ for $p \in [1,\infty]$ using the material developed in the previous chapter and in the preceding sections.

\begin{definition} \label{723}
Let $(X,\mathsf{S},\mu)$ be a measure space and let $p \in [1,\infty]$ be given. We define
$$
\begin{aligned}
\mathcal{L}^{p}_{\mathbb{C}}(X,\mathsf{S},\mu)&:=\left\{f \in \mathbb{M}_{\mathbb{C}}(X,\mathsf{S})\,:\,|f|^{p} \in \mathcal{L}(X,\mathsf{S},\mu) \right\}\qquad \text{if}\quad p \in [1,\infty),\\
\mathcal{L}^{\infty}_{\mathbb{C}}(X,\mathsf{S},\mu)&:=\left\{f \in \mathbb{M}_{\mathbb{C}}(X,\mathsf{S})\,:\,f\,\,\text{is $\mu$ a.e. bounded} \right\}\qquad \text{if}\quad p = \infty.\\
\end{aligned}
$$

Hence, $L^{p}_{\mathbb{C}}(X,\mathsf{S},\mu)$ is the set obtained by identifying all functions that are equal $\mu$ a.e. in $\mathcal{L}^{p}_{\mathbb{C}}(X,\mathsf{S},\mu)$ for every $p \in [1,\infty]$.
\end{definition}

\begin{remark} \label{724}
Everything stated, without exception, in the previous sections for $L^{p}(X,\mathsf{S},\mu)$ remains valid for $L^{p}_{\mathbb{C}}(X,\mathsf{S},\mu)$, except for minor modifications when considering scalars $\gamma \in \mathbb{C}$ and functions in $\mathbb{M}_{\mathbb{C}}(X,\mathsf{S})$. Therefore, $(L^{p}_{\mathbb{C}}(X,\mathsf{S},\mu), \Vert \cdot \Vert_{p})$ is a normed vector space and is complete for all $p \in [1,\infty]$.
\end{remark}

\subsection*{7.3.1.\quad The case when $p \in (0,1)$}

Let $(X,\mathsf{S},\mu)$ be a measure space and let $p \in (0,1)$. Recall the following definition
$$
\mathcal{L}^{p}(X,\mathsf{S},\mu):=\{ f \in \mathbb{M}(X,\mathsf{S})\,:\,|f|^{p}\in \mathcal{L}(X,\mathsf{S},\mu) \}.
$$

We denote by $L^{p}(X,\mathsf{S},\mu)$ the set of equivalence classes under the relation given in (\ref{F712}), which is a vector space with the same operations as in the case $p \in [1,\infty]$, thanks to Theorem \ref{73}. As before, we identify two functions on $X$ if and only if they are equal $\mu$ a.e. Thus, we simply denote the equivalence class $[f]_{\mu}$ by $f:X \to \mathbb{R}$.

We now ask whether the expression
$$
\Vert f \Vert_{p} :=\left( \int_{X} |f|^p\,d\mu \right)^{1/p}
$$
also defines a norm on $L^{p}(X,\mathsf{S},\mu)$ for $p \in (0,1)$, as it did in the previous section. The answer is that this is not true in general, since Minkowski's inequality does not extend to values of $p$ in $(0,1)$, as the following example shows.

\begin{example} \label{725}
Let $p:=1/2$ and consider the measure space $(X,\mathsf{S},\mu)$ given by $X=\mathbb{R}$, $\mathsf{S}=\mathcal{B}(\mathbb{R})$ and $\mu=\lambda$, the Lebesgue measure on $\mathsf{S}$.

It is clear that the functions $f_1:=\chi_{[0,1]}$ and $f_2:=\chi_{[1/2,3/2]}$ belong to the space $L^{1/2}(\mathbb{R},\mathcal{B}(\mathbb{R}),\lambda)$ and satisfy $\Vert f_1 \Vert_{1/2} = 1 =\Vert f_2 \Vert_{1/2}$.

On the other hand, we have
$$
f_1+f_2(x)=\left\{
\begin{array}{lcl}
1 & & \mbox{if}\,\,\,x \in [0,1/2)\cup (1,3/2],\\
2 & & \mbox{if}\,\,\,x \in [1/2,1],\\
0 & & \mbox{otherwise},
\end{array}
\right.
$$
and, consequently,
$$
\begin{aligned}
\Vert f_1 + f_2 \Vert_{1/2} &= \left(\int_{\mathbb{R}} |f_1+f_2|^{1/2}\,d\lambda \right)^{2}=\left(1 + \frac{\sqrt{2}}{2} \right)^2 = \frac{3}{2} + \sqrt{2}.
\end{aligned}
$$

Therefore, $\Vert f_1+f_2 \Vert_{1/2} > \Vert f_1 \Vert_{1/2} + \Vert f_2 \Vert_{1/2}$.
\end{example}

No coincidence is it that the triangle inequality is reversed for values of $p$ in $(0,1)$, and the proof of this fact is left as an exercise [Exercise \ref{E723}]. Similarly, there is a version of Hölder–Riesz inequality for the case when $p \in (0,1)$, and in that setting its conjugate exponent $q:=\frac{p}{p-1}$ takes values in the interval $(-\infty,0)$ [Exercise \ref{E724}].

Thus, the goal of this section is to endow the space $L^{p}(X,\mathsf{S},\mu)$ with a metric structure when $p \in (0,1)$. To this end, we briefly recall the definition of a metric.

\begin{definition} \label{726} \index{metric} 
Let $X$ be a non-empty set. A metric or distance on $X$ is a function $d:X \times X \to \mathbb{R}$ satisfying the following three properties:
\begin{itemize}
\item[\rm(D1)] $d(x,y)=0$ if and only if $x=y$.
\item[\rm (D2)] $d(x,y)=d(y,x)$ for all $x,y \in X$.
\item[\rm (D3)] $d(x,z) \leq d(x,y) + d(y,z)$ for all $x,y,z \in X$.
\end{itemize}

A \textbf{metric space} is a set $X$ equipped with a given metric $d$. \index{metric space}
\end{definition}

The following inequality is a key ingredient for achieving the objective.

\begin{lemma} \label{727}
Let $p \in (0,1)$ and let $a,b$ be two non-negative real numbers. Then, the following inequality holds:
\begin{eqnarray}
(a+b)^{p} \leq a^{p} + b^{p}.
\end{eqnarray}
\end{lemma}

\begin{proof}
For each fixed $a \geq 0$, define the function $\varsigma_a:[0,\infty) \to \mathbb{R}$ as
$$
\varsigma_a(x):=x^{p}+a^{p}-(a+x)^{p}.
$$

This function satisfies $\varsigma_{a}(0)=0$, is differentiable, and for each $x \geq 0$,
$$
\varsigma'_{a}(x)=p\left(x^{p-1}-(a+x)^{p-1} \right).
$$

\begin{figure}[ht!]
\centering
\begin{tikzpicture}[xscale=2.6, yscale=1]
\draw[domain= 0:2, thick] plot(\x,{\x^(1/2)+12-(144+\x)^(1/2)} );

\draw[->,gray] (-0.2,0)--(2,0); 
\draw[<-,gray] (0,2.35)--(0,-0.2); 







\end{tikzpicture}
\begin{center}
$\varsigma_{12}$ with $p=1/2$.
\end{center}
\end{figure}

Since $p-1<0$, we have $\varsigma'_{a}(x)=p\left(x^{p-1}-(a+x)^{p-1} \right)>0$, and this implies that the derivative $\varsigma_{a}'$ is positive and therefore $\varsigma_{a}$ is increasing on $[0,\infty)$. Thus, for all $b \geq 0$, we have $\varsigma_{a}(0)=0 \leq \varsigma_a(b)=b^{p}+a^{p}-(a+b)^{p},$
which is the desired inequality.
\end{proof}

\begin{theorem} \label{728}
For each $p \in (0,1)$, the function $d_{p}:L^{p}(X,\mathsf{S},\mu)\times L^{p}(X,\mathsf{S},\mu) \to \mathbb{R}$ defined by
\begin{eqnarray} \label{F719}
d_{p}(f,g):=\int_{X} |f-g|^{p}\,d\mu
\end{eqnarray}
defines a metric on $L^{p}(X,\mathsf{S},\mu)$.
\end{theorem}

\begin{proof}
Let $p\in (0,1)$. Proposition \ref{74} ensures that $d_{p}(f,g)=0$ if and only if $f=g$ in $L^{p}(X,\mathsf{S},\mu)$. Moreover, it is immediate from definition (\ref{F719}) that $d_{p}(f,g)=d_{p}(g,f)$ for all $f,g \in L^{p}(X,\mathsf{S},\mu)$.

Finally, Lemma \ref{727} ensures that
$$
|a+b|^{p} \leq (|a|+|b|)^{p} \leq |a|^{p} + |b|^{p}
$$
for all real numbers $a$ and $b$. Hence, for all $f,g,h \in L^{p}(X,\mathsf{S},\mu)$, we have
$$
|f-g|^{p}=|(f-h)+(h-g)|^{p} \leq |f-h|^{p} + |h-g|^{p}.
$$

Applying the monotonicity of the integral we obtain
$$
d_{p}(f,g) \leq d_{p}(f,h) + d_{p}(h,g),
$$
which proves the result.
\end{proof}

It is possible to show that the metric defined in (\ref{F719}) is complete on the space $L^{p}(X,\mathsf{S},\mu)$ for $p \in (0,1)$, following, with suitable modifications, the steps given in the proof of the Riesz–Fischer theorem [Exercise \ref{E727}].

We conclude that the most general topological structure shared by all Lebesgue spaces $L^{p}(X,\mathsf{S},\mu)$ with $p \in (0,\infty]$ is that of a metric space.

\section{Exercises}

\begin{exercise}[Generalized Young inequality] \label{E71} \index{Young inequality!generalized}
Let $a_{1},\ldots,a_{N}$ and $p_{1},\ldots,p_{N}$ be real numbers such that $a_{j} \geq 0$ and $p_{j}>0$ for each $j=1,\ldots,N$. Prove that if $\sum_{j=1}^{N}p_{j}=1$, then the following inequality holds
$$
\prod_{j=1}^{N} a_{j} \leq \sum_{j=1}^{N} p_{j}a_{j}^{1/p_{j}};
$$
and equality holds if and only if $a_{i}^{1/p_{i}}=a_{j}^{1/p_{j}}$ for all $i,j = 1,\ldots,N$.
\end{exercise}

\begin{exercise} \label{E72}
We denote by $\mathcal{C}([0,1])$ the set of all continuous functions $f:[0,1] \to \mathbb{R}$. The usual pointwise addition and scalar multiplication give $\mathcal{C}([0,1])$ the structure of a vector space. 
\begin{itemize}
    \item[(a)] Prove that the function $N_{1}:\mathcal{C}([0,1)] \to \mathbb{R}$ given by
$$
N_{1}(f):=R\int_{0}^{1} |f(x)|\,dx
$$
defines a norm on the space $\mathcal{C}([0,1])$.
    \item[(b)] Prove that the sequence of functions $f_k:[0,1] \to \mathbb{R}$ given by
$$
f_k(x):=\left\{
\begin{array}{lcl}
0 & & \mbox{if}\,\,\,\, 0 \leq x \leq \frac{1}{2}\left(1-\tfrac{1}{k}\right),\\
2kx-(k-1) & & \mbox{if}\,\,\,\,  \tfrac{1}{2}\left(1-\tfrac{1}{k}\right) \leq x \leq \tfrac{1}{2},\\
1 & & \mbox{if}\,\,\,\,  \tfrac{1}{2} \leq x \leq 1,
\end{array}
\right.
$$
is Cauchy in $(\mathcal{C}([0,1]),N_{1})$.

\item[(c)] Is $(\mathcal{C}([0,1]),N_{1})$ a Banach space? Justify your answer.
\end{itemize}
\end{exercise}

\begin{exercise} \label{E73}
Let $(x_k)$ be a Cauchy sequence in a metric space $X=(X,d)$.
\begin{itemize}
\item[(a)] Prove that $(x_k)$ is bounded in $X$.
\item[(b)] Prove that, if some subsequence of $(x_k)$ converges to $x$ in $X$, then $(x_k)$ converges to $x$ in $X$.
\end{itemize}
\end{exercise}

\begin{exercise} \label{E74}
Prove that the relation $\sim_{\mu}$ defined in {\rm (\ref{F712})} is an equivalence relation on $\mathcal{L}^{p}(X,\mathsf{S},\mu)$.
\end{exercise}

\begin{exercise} \label{E75}
Prove that the expressions given in {\rm(\ref{F714})} and {\rm(\ref{F715})} are well defined, i.e., they do not depend on the representative of the equivalence class.
\end{exercise}

\begin{exercise} \label{E76}
Let $a,b$ be real numbers such that $a < b$. Prove that the set
$$
\mathcal{C}([a,b]):=\{f:[a,b] \to \mathbb{R}\,:\, f\,\,\mbox{is continuous}\}
$$ 
is a closed vector subspace of $L^{\infty}([a,b],\mathcal{B}([a,b]),\lambda)$ and that
$$
\Vert f \Vert_{\infty}=\sup_{x \in [a,b]}|f(x)| \qquad \forall\,f \in \mathcal{C}([a,b]).
$$
\end{exercise}

\begin{exercise} \label{E77}
Let $(X,\mathsf{S},\mu)$ be a measure space and let $1 \leq p < q < r < \infty$. Prove that if $f \in L^{p}(X,\mathsf{S},\mu) \cap L^{r}(X,\mathsf{S},\mu)$, then $f \in L^{q}(X,\mathsf{S},\mu)$ and
$$
\Vert f \Vert_{q}^{q} \leq \Vert f \Vert_{p}^{\theta p} \Vert f \Vert_{r}^{(1-\theta)r},
$$
where $\theta \in (0,1)$ satisfies $q=\theta p + (1-\theta)r$.
\end{exercise}

\begin{exercise}[Interpolation inequality] \label{E78} \index{inequality!interpolation}
Let $(X,\mathsf{S},\mu)$ be a measure space and let $1 \leq p < q < r \leq \infty$. Prove that if $f \in L^{p}(X,\mathsf{S},\mu) \cap L^{r}(X,\mathsf{S},\mu)$, then $f \in L^{q}(X,\mathsf{S},\mu)$ and
$$
\Vert f \Vert_{q} \leq \Vert f \Vert_{p}^{\theta} \Vert f \Vert_{r}^{(1-\theta)},
$$
where $\theta \in (0,1)$ satisfies $\tfrac{1}{q} = \tfrac{\theta}{p} + \tfrac{1-\theta}{r}$ if $r<\infty$ and $\theta:=\tfrac{p}{q}$ if $r=\infty$.
\end{exercise}

{\setlength{\parindent}{0pt}
\begin{exercise} \label{E79}
Let $(X,\mathsf{S},\mu)$ be a measure space and let $f \in L^{1}(X,\mathsf{S},\mu) \cap L^{\infty}(X,\mathsf{S},\mu)$. Prove that $f \in L^{p}(X,\mathsf{S},\mu)$ for all $p \in (1,\infty)$ and that
$$
\Vert f \Vert_{p} \leq \Vert f \Vert_{1} + \Vert f \Vert_{\infty}.
$$
\end{exercise}

(Hint: Use the interpolation inequality.)}

{\setlength{\parindent}{0pt}
\begin{exercise} \label{E710}
Give an example of an infinite measure space $(X,\mathsf{S},\mu)$ for which $L^{r}(X,\mathsf{S},\mu) \subset L^{s}(X,\mathsf{S},\mu)$ holds for $1 \leq r < s < \infty$.
\end{exercise}

(Hint: Consider $\mathbb{N}$ with the discrete $\sigma$-algebra and counting measure.)}

{\setlength{\parindent}{0pt}
\begin{exercise} \label{E711}
Let $(X,\mathsf{S},\mu)$ be a finite measure space and let $f \in L^{\infty}(X,\mathsf{S},\mu)$. Prove that $f \in L^{p}(X,\mathsf{S},\mu)$ for each $p \in [1,\infty)$ and that
$$
\lim_{p \to \infty} \Vert f \Vert_{p} = \Vert f \Vert_{\infty}.
$$
\end{exercise}

(Hint: First use Theorem \ref{713}. Then consider $\varepsilon > 0$ with $0<\varepsilon<\Vert f \Vert_{\infty}$ and define $X_{\varepsilon}:=\{x \in X \,:\, |f(x)| \geq \Vert f \Vert_{\infty}-\varepsilon \}$. Prove that $\mu(X_{\varepsilon})<+\infty$.)}

\begin{exercise}[Generalized Hölder–Riesz inequality] \label{E712} \index{inequality!generalized Hölder–Riesz}
Let $(X,\mathsf{S},\mu)$ be a measure space. Let $r,p_1,\ldots,p_m \in [1,\infty)$ satisfy $\frac{1}{p_1}+\cdots+\frac{1}{p_m}=\frac{1}{r}$. Prove that, for any $f_j \in L^{p_j}(X,\mathsf{S},\mu)$, $1\leq j \leq m$, it holds that $\prod_{j=1}^{m} f_j \in L^{r}(X,\mathsf{S},\mu)$ and
$$
\left\Vert \prod_{j=1}^{m} f_{j} \right\Vert_{r} \leq \prod_{j=1}^{m}\Vert f_j\Vert_{p_j}.
$$

When is it possible to conclude equality? Justify your answer in detail.
\end{exercise}

\begin{exercise} \label{E713}
Let $(X,\mathsf{S},\mu)$ be a measure space and let $p,q \in (1,\infty)$ satisfy $\frac{1}{p}+\frac{1}{q}=1$. Given $g \in L^{q}(X,\mathsf{S},\mu)$, define $\eta_{g}:L^{p}(X,\mathsf{S},\mu) \to \mathbb{R}$ by
$$
\eta_{g}(f):=\int_{X} fg\,d\mu.
$$ 

Prove that $\eta_{g}$ is a linear and continuous function.
\end{exercise}

\begin{exercise} \label{E714}
Let $(X,\mathsf{S},\mu)$ be a measure space, let $f \in L^{p}(X,\mathsf{S},\mu)$ with $p \in [1,\infty)$ and let $\varepsilon >0$ be given. Prove that there exists an $\mathsf{S}$-simple function $\varsigma \in L^{p}(X,\mathsf{S},\mu)$ such that $\Vert f - \varsigma \Vert_{p} < \varepsilon$.
\end{exercise}

{\setlength{\parindent}{0pt}
\begin{exercise} \label{E715}
Let $(X,\mathsf{S},\mu)$ be a finite measure space and let $p\in (1,\infty)$. Prove that, for a given $\vartheta \in L^{1}(X,\mathsf{S},\mu)$, there exists a sequence $(\vartheta_k)$ of elements of $L^{p}(X,\mathsf{S},\mu)$ such that $\lim_{k \to \infty}\Vert \vartheta - \vartheta_{k} \Vert_{1} = 0$.
\end{exercise}

(Hint: Define the measurable sets $A_{k}:=\{x \in X \,:\, |\vartheta(x)| < k \}$ for all $k \in \mathbb{N}$ and the sequence $\vartheta_k:=(\vartheta\cdot \chi_{A_k})$. Use the dominated convergence theorem.)}

{\setlength{\parindent}{0pt}
\begin{exercise} \label{E716}
Let $X=\mathbb{R}$, $\mathsf{S}=\mathcal{B}(\mathbb{R})$ and $\mu=\lambda$ be the Lebesgue measure on $\mathsf{S}$. Let $\varphi:\mathbb{R} \to \mathbb{R}$ be defined by
$$
\varphi(x)=\left\{
\begin{array}{lcl}
x^{-\frac{1}{2}} & & \mbox{if}\,\,\,x\in (0,1),\\
0 & & \mbox{if}\,\,\,x\notin (0,1).\\
\end{array}
\right.
$$

Let $\sum_{k=1}^{\infty} \varsigma_{k}$ be a convergent series of positive real numbers and let $\mathbb{Q}=\{q_{1},q_{2},\ldots \}$ be an enumeration. Consider the function $f:\mathbb{R} \to [0,+\infty)$ defined by $f(x):=\sum_{k=1}^{\infty}\varsigma_{k}\,\varphi(x-q_{k})$. Prove that $f \in L^{1}(\mathbb{R},\mathcal{B}(\mathbb{R}),\lambda).$ (In fact, $\int_{\mathbb{R}} f\,d\lambda =2\sum_{k=1}^{\infty}\varsigma_{k}$, so that $f(x)<\infty$ $\lambda$ a.e.).
\end{exercise}

(Hint: Use the dominated convergence theorem and the Lebesgue integral theorem.)}

\begin{exercise} \label{E717}
Let $(X,\mathsf{S},\mu)$ be a measure space and let $1 \leq p < \infty$ and $1 \leq q \leq \infty$. Prove that the set
$$
\mathcal{H}:=\left\{f \in L^p(X,\mathsf{S},\mu)\cap L^{q}(X,\mathsf{S},\mu)\,:\, \Vert f \Vert_{q} \leq 1 \right\}
$$
is closed in $L^p(X,\mathsf{S},\mu)$.
\end{exercise}

\begin{exercise} \label{E718}
Let $(X,\mathsf{S},\mu)$ be a measure space. Define the function $\Vert \cdot \Vert_{1\cap \infty}:L^1(X,\mathsf{S},\mu)\cap L^{\infty}(X,\mathsf{S},\mu) \to \mathbb{R}$ by
$$
\Vert f \Vert_{1\cap \infty}:=\max\{\Vert f \Vert_1, \Vert f \Vert_{\infty} \}.
$$

Prove that $\Vert \cdot \Vert_{1\cap \infty}$ is a norm. Is $L^1(X,\mathsf{S},\mu)\cap L^{\infty}(X,\mathsf{S},\mu)$ a Banach space with this norm if $\mu(X)<+\infty$? Justify your answer.
\end{exercise}

\begin{exercise} \label{E719}
Let $(X,\mathsf{S},\mu)$ be a measure space and let $s \in [1,\infty]$. For given $f_{1},\ldots,f_{N} \in L^{s}(X,\mathsf{S},\mu)$, define $f:X \to \mathbb{R}^{N}$ by $f(x):=(f_{1}(x),\ldots,f_{N}(x))$ and denote
$$
\Vert f \Vert_{s}:=\left(\Vert f_{1} \Vert_{s}^{s} + \ldots + \Vert f_{N} \Vert_{s}^{s} \right)^{1/s}.
$$

Prove that if $\mu(X)< +\infty$, $1 \leq p < s \leq \infty$ and $f_{1},\ldots,f_{N} \in L^{s}(X,\mathsf{S},\mu)$, then
\begin{itemize}
\item[(a)] $\Vert f \Vert_{p} \leq (N\mu(X))^{(s-p)/sp}\Vert f \Vert_{s}$ if $s \in [1,\infty)$.
\item[(b)] $\Vert f \Vert_{p} \leq (N\mu(X))^{1/p}\Vert f \Vert_{s}$ if $s=\infty$.
\end{itemize}
\end{exercise}

\begin{exercise} \label{E720}
Let $(X,\mathsf{S},\mu)$ be a finite measure space, let $f \in \mathbb{M}(X,\mathsf{S})$, and let $ p \in [1,\infty)$. Prove that $f \in L^{p}(X,\mathsf{S},\mu)$ if and only if
$$
\sum_{k=1}^{\infty} k^{p}\,\mu\left(\left\{ x \in X\,:\, (k-1) \leq |f(x)| < k \right\} \right) < +\infty.
$$
\end{exercise}

\begin{exercise} \label{E721}
Let $(X,\mathsf{S},\mu)$ be a measure space and let $f \in L^{1}(X,\mathsf{S},\mu)$ be such that $\mu(A_f)<+\infty$, where
$$
A_{f}:=\left\{ x \in X\,:\,f(x) \neq 0 \right\}.
$$

Prove that
$$
\lim_{k \to \infty}\int_{X}|f|^{1/k}\,d\mu  = \mu(A_f).
$$

Is the result still true if we do not assume $\mu(A_{f})<+\infty$? Justify your answer.
\end{exercise}

{\setlength{\parindent}{0pt}
\begin{exercise} \label{E722}
Let $X=(0,1)$, $\mathsf{S}=\mathcal{B}(0,1)$ and let $\lambda$ be the Lebesgue measure on $\mathsf{S}$. Prove that
$$
L^{p}((0,1),\mathcal{B}(0,1),\lambda) \neq \bigcap_{0<r<p} L^{r}((0,1),\mathcal{B}(0,1),\lambda)\quad \text{ for every fixed }\,p\geq 1.
$$
\end{exercise}

(Hint: If $p=1$ consider $f_{1}:=\frac{1}{x}$ and if $p \neq1$ modify it appropriately.)}

{\setlength{\parindent}{0pt}
\begin{exercise}[Inversion of the Hölder–Riesz inequality] \label{E723} \index{inversion!of Hölder–Riesz inequality}
Let $(X,\mathsf{S},\mu)$ be a measure space and let $p\in (1,\infty)$. Set $r:=\frac{1}{p} \in (0,1)$ and $s:=\frac{r}{r-1} \in (-\infty,0)$. Let $f,g:X \to \mathbb{R}$ be such that $f^{r}, g^{s}, fg \in L^{1}(X,\mathsf{S},\mu)$ and $g \neq 0$. Prove that
$$
\left(\int_{X} |f|^{r}\,d\mu \right)^{1/r}\left(\int_{X} |g|^{s}\,d\mu \right)^{1/s} \leq \int_{X} |fg|\,d\mu.
$$
\end{exercise}

(Hint: Justify that $p$ and $-\frac{s}{r}$ are conjugate exponents. Show that $|fg|^{r} \in L^{p}$ and $|g|^{-r} \in L^{-s/r}$, and use the usual Hölder–Riesz inequality.)}

\begin{exercise}[Inversion of Minkowski’s inequality] \label{E724}  \index{inversion!of Minkowski's inequality}
Let $(X,\mathsf{S},\mu)$ be a measure space, let $p\in (1,\infty)$, and let $r:=\frac{1}{p} \in (0,1)$. Prove that if $f,g:X \to \mathbb{R}$ are positive functions such that $f^{r},g^{r} \in L^{1}(X,\mathsf{S},\mu)$, then
$$
\left(\int_{X} |f+g|^{r}\,d\mu \right)^{1/r} \geq \left( \int_{X} |f|^{r}\,d\mu \right)^{1/r} + \left( \int_{X} |g|^{r}\,d\mu \right)^{1/r}.
$$
\end{exercise}

(Hint: Use the inversion of the Hölder–Riesz inequality)

\begin{exercise} \label{E725}
Let $X=[0,1]$, $\mathsf{S}=\mathcal{B}[0,1]$ and let $\mu=\lambda$ be the Lebesgue measure on $\mathsf{S}$. Consider $p=1/2$ and the functions $f:=\chi_{[0,1]}$ and $g:=-\tfrac{1}{2}\chi_{[0,1]}$. Verify that $f,g \in L^{1/2}$ and that the following inequality holds:
$$
\Vert f + g \Vert_{1/2} \leq \Vert f \Vert_{1/2} + \Vert g \Vert_{1/2}.
$$

Is this a contradiction to the inversion of Minkowski’s inequality? Justify your answer.
\end{exercise}

\begin{exercise} \label{E726}
Let $(X,\mathsf{S},\mu)$ be a measure space and let $p \in (0,1)$. Prove that for any $f,g \in L^{p}(X,\mathsf{S},\mu)$ the following inequality holds:
$$
\Vert f +g \Vert_{p} \leq 2^{\frac{1}{p}-1}\left(\Vert f \Vert_{p} + \Vert g \Vert_{p} \right).
$$
\end{exercise}

\begin{exercise} \label{E727}
Let $(X,\mathsf{S},\mu)$ be a measure space. Prove that the space $(L^{p}(X,\mathsf{S},\mu),d_p)$ is a complete metric space for all $p \in (0,1)$.
\end{exercise}

\begin{exercise} \label{E728}
Let $(X,\mathsf{S},\mu)$ be a measure space and let $p \in (0,\infty]$. Is it true that if $(f_k)$ is a sequence in $L^{p}(X,\mathsf{S},\mu)$ such that $\lim_{k \to \infty}\Vert f_{k}\Vert_{p}=\Vert f \Vert_{p}$, then $f_{k} \to f$ in $L^{p}(X,\mathsf{S},\mu)$? Justify your answer in detail.
\end{exercise}

{\setlength{\parindent}{0pt}
\begin{exercise}  \label{E729}
Let $(X,\mathsf{S},\mu)$ be a measure space and let $0 < p < q < r  \leq \infty$. Prove that
$$
L^{p}(X,\mathsf{S},\mu) \cap L^{r}(X,\mathsf{S},\mu) \subset L^{q}(X,\mathsf{S},\mu) \subset L^{p}(X,\mathsf{S},\mu) + L^{r}(X,\mathsf{S},\mu),
$$
where
$$
 L^{p}(X,\mathsf{S},\mu) + L^{r}(X,\mathsf{S},\mu):=\left\{f+g\,:\,f \in L^{p}(X,\mathsf{S},\mu),\, g \in L^{r}(X,\mathsf{S},\mu) \right\}.
$$
\end{exercise}

(Hint: Consider the set $A:=\{x \in X\,:\,|f(x)| \leq 1\}$ and use Exercise \ref{E613} to show that $\mu(X\smallsetminus A)<+\infty$).}

\begin{exercise} \label{E730}
Let $(X,\mathsf{S},\mu)$ be a measure space and let $0 < r < s < \infty$. Prove the following:
\begin{itemize}
    \item[(a)] If $f \in L^{r}(X,\mathsf{S},\mu)\cap L^{s}(X,\mathsf{S},\mu)$, then $f \in L^{p}(X,\mathsf{S},\mu)$ for every $p \in [r,s]$.
    \item[(b)] The function $\Phi:[r,s] \to \mathbb{R}$ defined by
    $$
\Phi(p):=\log\left( \int_{X}|f|^{p}\,d\mu \right)
    $$
    is convex.
    \item[(c)] $\Vert f \Vert_{p} \leq \max\left\{\Vert f \Vert_{r},\Vert f \Vert_{s} \right\}$ for every $p \in [r,s]$.  
\end{itemize}
\end{exercise}

\section{Project: Constructing an isometric isomorphism}

\begin{definition} \label{729}  \index{isometry}
Let $V=(V,\Vert \cdot \Vert_V)$ and $W=(W,\Vert \cdot \Vert_W)$ be normed vector spaces. A function $\iota : V \to W$ is called an \textbf{isometry} if
$$
\Vert \iota(u)-\iota(v) \Vert_W = \Vert u-v \Vert_V \qquad \forall \,u,v \in V.
$$

If, in addition, $\iota$ is linear and bijective, we call it an \textbf{isometric isomorphism}. \index{isometric isomorphism}
\end{definition}

From the previous definition we immediately observe that every isometry is injective. Indeed, if $v_1,v_2 \in V$ satisfy $\iota(v_1)=\iota(v_2)$, then
$$
\Vert v_1-v_2\Vert_V=\Vert \iota(v_1) - \iota(v_2) \Vert_W =0,
$$
and hence $v_1=v_2$.

From a geometric point of view, two normed spaces are considered the same if there exists a linear bijection between them which is an isometry. In that case we say that they are isometrically isomorphic.

\subsection*{7.5.1.\quad Objective}

Let $(X,\mathsf{S},\mu)$ be an arbitrary measure space. We denote by \index{space!L1ast@$(L^{1}(X,\mathsf{S},\mu))^{\ast}$}
$$ 
(L^{1}(X,\mathsf{S},\mu))^{\ast}:=\{T:L^{1}(X,\mathsf{S},\mu) \to \mathbb{R}\,:\,T\text{ is linear and continuous}\},
$$
which is a vector space over $\mathbb{R}$ with the usual pointwise addition and scalar multiplication of functions.

For $T \in (L^{1}(X,\mathsf{S},\mu))^{\ast}$ we define
$$
\Vert T \Vert:=\sup\left\{ |T(g)|\,:\,\Vert g \Vert_{1} \leq 1 \right\}.
$$

Prove the following theorem.

\begin{theorem}
Let $(X,\mathsf{S},\mu)$ be a measure space.
\begin{itemize}
    \item[(a)] $\Vert \cdot \Vert$ is a norm on the vector space $(L^{1}(X,\mathsf{S},\mu))^{\ast}$.
    \item[(b)] Let $f \in L^{\infty}(X,\mathsf{S},\mu)$ be fixed. The map $\Psi_{f}:L^{1}(X,\mathsf{S},\mu) \to \mathbb{R}$ defined by $\Psi_{f}(g):=\int_{X} fg\,d\mu$ is linear and continuous on $L^{1}(X,\mathsf{S},\mu)$. That is, $\Psi_{f} \in (L^{1}(X,\mathsf{S},\mu))^{\ast}$.
\end{itemize}    
\end{theorem}

Let $(X,\mathsf{S},\mu)$ be a $\sigma$-finite measure space. The goal of this project is to construct a finite measure $\nu$ on $(X,\mathsf{S})$ such that, for $p\in\{1,\infty\}$, the spaces $L^{p}(X,\mathsf{S},\mu)$ and $L^{p}(X,\mathsf{S},\nu)$ are isometrically isomorphic. 

From this identification, we aim to establish an explicit relationship between the spaces $(L^{1}(X,\mathsf{S},\mu))^{\ast}$ and $(L^{1}(X,\mathsf{S},\nu))^{\ast}$, and the spaces $L^{\infty}(X,\mathsf{S},\mu)$ and $L^{\infty}(X,\mathsf{S},\nu)$.

\subsection*{7.5.2.\quad Procedure}

Let $(X,\mathsf{S},\mu)$ be a $\sigma$-finite measure space.

\begin{itemize}
    \item[1.] Prove that there exists a disjoint sequence $(A_k)$ of elements of $\mathsf{S}$ such that $X=\bigcup_{k=1}^{\infty} A_k$ and such that $0< \mu(A_k)<+\infty$ for all $k \in \mathbb{N}$.

    \item[2.] Prove that the map $\Psi:L^{\infty}(X,\mathsf{S},\mu) \to (L^{1}(X,\mathsf{S},\mu))^{\ast}$ defined by $\Psi(f):=\Psi_{f}$ is a linear isometry.

    (Hint: If $f \neq 0$, consider $\varepsilon >0$ such that $0<\varepsilon<\Vert f \Vert_{\infty}$ and $A \in \mathsf{S}$ with $0< \mu(A) <+\infty$ such that $|f(x)|\geq \Vert f \Vert_{\infty}-\varepsilon$ for all $x \in A$. Apply $\Psi$ to the function $\widetilde{f}:=\frac{|f|}{f}\chi_{A}$.)
    
    \item[3.] Prove that there exists a finite measure $\nu:\mathsf{S} \to \mathbb{R}$ such that if $\mu(A)=0$ then $\nu(A)=0$ for all $A \in \mathsf{S}$. 
    
    (Hint: Consider the formula $\nu(A):=\sum_{k=1}^{\infty} \frac{1}{2^k\,\mu(A_k)}\,\mu(A \cap A_k).$)
    
    \item[4.] Let $f:X \to \mathbb{R}$ be a non-negative measurable function. Prove that there exist a non-negative measurable function $\vartheta:X \to \mathbb{R}$ and a finite measure $\nu:\mathsf{S} \to \mathbb{R}$ such that
$$
\int_{X} f\,d\mu = \int_{X} f\vartheta\,d\nu.
$$ 

(Hint: For each $j \in \mathbb{N}$ define $\vartheta_j:=2^{j}\mu(A_j)\,\chi_{A_j}$ and $\vartheta:=\sum_{j=1}^{\infty}\vartheta_j$.)

    \item[5.] Consider $\vartheta:X \to \mathbb{R}$ and the finite measure $\nu:\mathsf{S} \to \mathbb{R}$ from the previous item. Prove that the map $\eta:L^{1}(X,\mathsf{S},\mu) \to  L^{1}(X,\mathsf{S},\nu)$ defined by $\eta(u):=u\vartheta$ is an isometric isomorphism.

    \item[6.] As in the previous item, consider $\vartheta:X \to \mathbb{R}$ and the finite measure $\nu:\mathsf{S} \to \mathbb{R}$ from item 3. Prove that the identity map $\mathrm{id}:L^{\infty}(X,\mathsf{S},\mu) \to  L^{\infty}(X,\mathsf{S},\nu)$ is an isometric isomorphism.
    
    \item[7.] Let $\eta:L^{1}(X,\mathsf{S},\mu) \to L^{1}(X,\mathsf{S},\nu)$ be the isometric isomorphism from item 5. Define $\eta^{\ast}:(L^{1}(X,\mathsf{S},\nu))^{\ast} \to (L^{1}(X,\mathsf{S},\mu))^{\ast}$ by
    $$
\eta^{\ast}(T):=(T \circ \eta).
    $$

   Prove that $\eta^{\ast}$ is an isometric isomorphism. We call $\eta^{\ast}$ the \textbf{adjoint operator of $\eta$}. 

    \item[8.] Let $(X,\mathsf{S},\mu)$ be a $\sigma$-finite measure space, let $\nu$ be the finite measure from item 3, and let $\eta:L^{1}(X,\mathsf{S},\mu) \to L^{1}(X,\mathsf{S},\nu)$ be the isometric isomorphism from item 5. Prove that the following diagram commutes:
    \begin{figure}[ht!]
        \centering
\begin{tikzpicture}[>=stealth, baseline=(current bounding box.center),scale=0.5]
  \node (A) {$L^{\infty}(X,\mathsf{S},\mu)$};
  \node (B) [right=2cm of A] {$(L^{1}(X,\mathsf{S},\mu))^{\ast}$};
  \node (C) [below=1cm of A] {$L^{\infty}(X,\mathsf{S},\nu)$};
  \node (D) [right=2cm of C] {$(L^{1}(X,\mathsf{S},\nu))^{\ast}$};

  \draw[->] (A) -- node[above] {$\Psi$} (B);
  \draw[->] (C) -- node[below] {$\Psi$} (D);

  \draw[->] (A) -- node[left] {$\text{id}$} (C);
  \draw[<-] (B) -- node[right] {$\eta^{\ast}$} (D);
\end{tikzpicture}
    \end{figure}
\end{itemize}

This result shows that the $L^{1}$ and $L^{\infty}$ spaces associated with equivalent measures are isometrically isomorphic. In particular, constructing an equivalent finite measure allows us to work with more manageable function spaces without altering their essential structure. Under either measure, functions and linear functionals can be described in an equivalent way, which justifies the use of change-of-measure techniques as a fundamental tool in analysis and probability to transfer properties and results between different contexts.

\chapter{Modes of convergence} \label{Capitulo8}
\markboth{{\scriptsize 8. MODES OF CONVERGENCE }}{ {\scriptsize 8. MODES OF CONVERGENCE}}

In an introductory course on mathematical analysis, two modes of convergence of sequences of functions are studied: pointwise convergence and uniform convergence. The latter is important because, under it, the limit of a sequence of continuous functions is itself a continuous function. Moreover, it is stronger than pointwise convergence in the sense that uniform convergence implies pointwise convergence.

In this text we focus on measurable functions, a broader class than that of continuous functions. We prove that the pointwise limit of a sequence of measurable functions is still measurable, which is why in previous chapters we focused on this notion of convergence.

In this chapter we introduce different modes of convergence for sequences of measurable functions, in which the concepts of measure and integration play an essential role. We study the relationships between these notions in the general case, that is, when the measure space is arbitrary, and we also discuss additional implications that arise when the space has finite measure.

Many of these modes of convergence have relevant applications in several areas related to mathematics, such as probability theory and ergodic theory.

\section{Definitions and examples}

Let $(X,\mathsf{S},\mu)$ be a measure space. 

\begin{definition} \label{81} \index{convergence!almost everywhere}
A sequence of $\mathsf{S}$-measurable functions $f_k:X \to \mathbb{R}$ converges \textbf{$\mu$ almost everywhere} (\textbf{$\mu$ a.e.}) to an $\mathsf{S}$-measurable function $f:X \to \mathbb{R}$ if there exists a $\mu$-null subset $N$ of $X$ such that $f_k(x) \to f(x)$ in $\mathbb{R}$ for every $x \in X\smallsetminus N$. That is, $(f_k)$ converges $\mu$ almost everywhere to $f \in \mathbb{M}(X,\mathsf{S})$ if there exists $N \in \mathcal{N}(\mu)$ such that for every $\varepsilon >0$ and every $x \in X \smallsetminus N$, there exists $k_0 \in \mathbb{N}$ (which depends on $\varepsilon$ and $x$) such that
$$
|f_k(x)-f(x)|<\varepsilon\quad \forall\,k\geq k_0.
$$
\end{definition}

It is clear that if a sequence of measurable functions $(f_k)$ converges pointwise on $X$ to a measurable function $f$, then it converges $\mu$ almost everywhere to $f$. The converse is not true, since, for example, the sequence $f_k:=k\chi_{[0,\frac{2}{k}]}$ from Example \ref{612} does not converge pointwise to $0$ on $[0,1]$.

We now consider another example.

\begin{example} \label{82}
Let $X=[0,1]$, $\mathsf{S}=\mathcal{B}([0,1])$ and let $\mu=\lambda$ be the Lebesgue measure on $\mathsf{S}$. Let $f_k:[0,1]\to \mathbb{R}$ be the function defined by
$$
f_k(x):=\max\left\{1-k\left|x-\tfrac{1}{k} \right|,0 \right\}.
$$

\begin{figure}[ht!]
\centering
\begin{tikzpicture}[xscale=1,yscale=0.85]
	\draw[->,gray] (0,0) -- (4,0); \draw [->,gray] (0,0) -- (0,4); 
	\draw[<-,gray] (-0.45,0)--(0,0); \draw[-,gray] (0,0)--(0,-0.23);
\draw (4,0) node[right] {$[0,1]$}; 

\draw (0,3.5) node{$_{-}$}; \draw (0,3.5) node[left]{$_{1}$};
\draw (1.5,0) node{$_{|}$}; \draw (1.5,-0.1) node[below]{$_{\frac{2}{k}}$};
\draw (0.8,0) node{$_{|}$}; \draw (0.8,-0.1) node[below]{$_{\frac{1}{k}}$};
\draw (3,0) node{$_{|}$}; \draw (3,-0.1) node[below]{$_{1}$};

\draw[ultra thick] (0,0)--(0.8,3.5); 
\draw[ultra thick] (0.8,3.5)--(1.5,0); 
\draw[ultra thick] (1.5,0)--(4,0); 
\draw (0,-0.2) node[below]{$_{0}$};
\draw [dotted] (0.8,0)--(0.8,3);
\end{tikzpicture}
\begin{center}
    $f_k(x)$
\end{center}  
\end{figure}

The sequence $(f_k)$ converges to $0$ $\lambda$ almost everywhere.
\end{example}

\begin{proof}
It is clear that $(f_k)$ is a sequence of Borel-measurable functions on $[0,1]$, since $f_k$ is continuous on $[0,1]$ for every $k \in \mathbb{N}$. We now prove that $(f_k)$ converges pointwise to $0$ on $[0,1]$. The case $x=0$ is immediate, so assume that $x \in (0,1]$. For each $\varepsilon >0$, choose $k_0 \in \mathbb{N}$ such that $\frac{2}{k_0} \leq x$. Consequently, $f_k(x)=0$ for all $k \geq k_0$, and therefore $|f_k(x)|<\varepsilon$ for all $k \geq k_0$.

This completes the proof.
\end{proof}

An important remark regarding Definition \ref{81} is that, in general, if a sequence of functions $(f_k)$ in $\mathbb{M}(X,\mathsf{S})$ converges $\mu$ almost everywhere to a function $f:X \to \mathbb{R}$, then the latter need not be $\mathsf{S}$-measurable. Indeed, if $N \subset X$ is a $\mu$-null set, $M \subset N$, and we consider the constant sequence $f_k \equiv 0$ on $X$, then we have $f_k \to \chi_{M}$ $\mu$ almost everywhere. However, $\chi_M$ may fail to be $\mathsf{S}$-measurable, since $M$ need not belong to $\mathsf{S}$ when the measure space $(X,\mathsf{S},\mu)$ is not complete [Exercise \ref{E826}]. For this reason, throughout this text we will work with sequences of $\mathsf{S}$-measurable functions that converge $\mu$ almost everywhere to a function which is also $\mathsf{S}$-measurable.

\begin{definition} \label{83} \index{convergence!in measure}
A sequence of $\mathsf{S}$-measurable functions $f_k:X \to \mathbb{R}$ converges \textbf{in measure $\mu$} to an $\mathsf{S}$-measurable function $f:X \to \mathbb{R}$ if, for every $\xi >0$ and every $\varepsilon >0$, there exists $k_0 \in \mathbb{N}$ (which may depend on $\varepsilon$) such that
$$
\mu\left( \left\{ x \in X\,:\,|f_k(x)-f(x)|\geq \xi \right\} \right) < \varepsilon \quad\forall k \geq k_0.
$$
We denote this type of convergence by $f_{k} \xrightarrow[\mu]{} f$.

Equivalently, $f_{k} \xrightarrow[\mu]{} f$ if, for every $\xi >0$,
$$
\lim_{k \to \infty}\,\mu\left(\left\{x \in X\,:\,|f_k(x)-f(x)|\geq \xi \right\} \right)=0.
$$
\end{definition}

We now consider an example.

\begin{example} \label{84}
Let $X=[0,1]$, $\mathsf{S}=\mathcal{B}([0,1])$ and let $\mu=\lambda$ be the Lebesgue measure on $\mathsf{S}$. The sequence of functions $f_k:=k\chi_{[0,\frac{1}{k}]}$ converges to $0$ in measure.
\end{example}

\begin{figure}[ht!]
\begin{minipage}[c]{0.5\textwidth}
\begin{center}
\begin{tikzpicture}[xscale=0.85,yscale=0.8]
	\draw[->,gray] (0,0) -- (4,0); \draw [->,gray] (0,0) -- (0,4); 
	\draw[<-,gray] (-0.45,0)--(0,0); \draw[-,gray] (0,0)--(0,-0.23);

\draw (0,2) node{$_{-}$}; \draw (0,2) node[left]{$_{k}$};
\draw (1.5,0) node{$_{|}$}; \draw (1.5,-0.1) node[below]{$_{\frac{1}{k}}$};
\draw (3,0) node{$_{|}$}; \draw (3,-0.1) node[below]{$_{1}$};
\draw[ultra thick] (0,2)--(1.47,2); \draw (0,2) node{$_{\bullet}$}; \draw (1.5,2) node{$_{\circ}$}; 
\draw[ultra thick] (1.5,0)--(3,0); 
\draw (0,-0.1) node[below]{$_{0}$};
\draw [dotted] (1.5,0)--(1.5,2);
\end{tikzpicture}
\begin{center}
    $k\chi_{[0,\frac{1}{k}]}(x)$
\end{center}  
\end{center} 
\end{minipage} \hfill \begin{minipage}[c]{0.5\textwidth}
\begin{center}
\begin{tikzpicture}[xscale=0.85,yscale=0.8]
	\draw[->,gray] (0,0) -- (4,0); \draw [->,gray] (0,0) -- (0,4); 
	\draw[<-,gray] (-0.45,0)--(0,0); \draw[-,gray] (0,0)--(0,-0.23);

\draw (0,2.5) node{$_{-}$}; \draw (0,2.5) node[left]{$_{k+1}$};
\draw (1,0) node{$_{|}$}; \draw (1,-0.1) node[below]{$_{\frac{1}{k+1}}$};
\draw (3,0) node{$_{|}$}; \draw (3,-0.1) node[below]{$_{1}$};
\draw[ultra thick] (0,2.5)--(1,2.5); \draw (0,2.5) node{$_{\bullet}$}; \draw (1,2.5) node{$_{\circ}$}; 
\draw[ultra thick] (1,0)--(3,0); 
\draw (0,-0.1) node[below]{$_{0}$};
\draw [dotted] (1,0)--(1,2.5);
\end{tikzpicture}
\begin{center}
    $(k+1)\chi_{[0,\frac{1}{k+1}]}(x)$
\end{center}    
\end{center} 
\end{minipage}
\end{figure}

\begin{proof}
Let $\xi >0$ and let $k \in \mathbb{N}$. Observe that if $\xi\leq k$, then
$\{ x \in [0,1]\,:\, |f_k(x)| \geq \xi\}=\left[0,\frac{1}{k}\right]$
and, consequently,
$$
\lambda \left( \left\{ x \in [0,1]\,:\,|f_k(x)| \geq \xi \right\} \right)
=\lambda \left(\left[0,\frac{1}{k} \right]\right)=\frac{1}{k}.
$$

On the other hand, if $\xi > k$, then the set $\left\{ x \in [0,1]\,:\,|f_k(x)| \geq \xi \right\}$ is empty and therefore
$$
\lambda \left( \left\{ x \in [0,1]\,:\,|f_k(x)| \geq \xi \right\} \right)=0.
$$

Thus, in any case we have
$$
0\leq \lambda \left( \left\{ x \in [0,1]\,:\,|f_k(x)| \geq \xi \right\} \right) \leq \frac{1}{k}, \qquad \forall k \in \mathbb{N}.
$$

Taking the limit as $k \to \infty$ we conclude that
$$
\lim_{k \to \infty}\lambda \left( \left\{ x \in X\,:\,|f_k(x)| \geq \xi \right\} \right) \leq \lim_{k \to \infty}\frac{1}{k}=0.
$$

This shows that $f_{k} \xrightarrow[\lambda]{} 0$.
\end{proof}

\begin{definition} \label{85} \index{convergence!almost uniform}
A sequence of $\mathsf{S}$-measurable functions $f_k:X \to \mathbb{R}$ converges \textbf{$\mu$ almost uniformly} ($\mu$ \textbf{a.u.}) to an $\mathsf{S}$-measurable function $f:X \to \mathbb{R}$ if, for every $\delta >0$, there exists a measurable set $F \in \mathsf{S}$ with $\mu(F)<\delta$ such that $f_k \to f$ uniformly on $X\smallsetminus F$.

That is, $(f_k)$ converges $\mu$ almost uniformly to $f \in \mathbb{M}(X,\mathsf{S})$ if, for every $\delta >0$, there exists a measurable set $F \in \mathsf{S}$ with $\mu(F)<\delta$ such that for every $\varepsilon >0$, there exists $k_0 \in \mathbb{N}$ such that
$$
|f_k(x)-f(x)|<\varepsilon \quad \forall k \geq k_0,\quad \forall x \in X \smallsetminus F.
$$

We denote this type of convergence by $f_k \to f$ $\mu$ a.u.
\end{definition}

It is immediate that if $(f_k)$ is a sequence of measurable functions that converges uniformly to $f$ on $X$, then $(f_k)$ converges to $f$ $\mu$ almost uniformly. The converse is not true in general, as shown by the following example.
 
\begin{example} \label{86}
Let $X=[0,1]$, $\mathsf{S}=\mathcal{B}([0,1])$ and let $\mu=\lambda$ be the Lebesgue measure on $\mathsf{S}$. Let $f_k:[0,1]\to \mathbb{R}$ be the function defined by $f_k:=k\chi_{\left[\frac{1}{k},\frac{2}{k} \right]}$.

\begin{figure}[ht!]
    \centering
	\begin{tikzpicture}[xscale=1,yscale=0.8]
	\draw[->,gray] (0,0) -- (4,0); \draw [->,gray] (0,0) -- (0,4); 
	\draw[-,gray] (4,0)--(0,0); \draw[-,gray] (0,0)--(0,4);
	\draw[-,gray] (-0.4,0)--(0,0); \draw[-,gray] (0,0)--(0,-0.45);

\draw (0,2.8) node{$-$}; \draw (0,2.8) node[left] {$_{k}$};
\draw (1,0) node{$|$}; \draw (1,-0.27) node[below] {$_{\frac{1}{k}}$};
\draw (2.3,0) node{$|$}; \draw (2.3,-0.2) node[below] {$_{\frac{2}{k}}$};
\draw [ultra thick] (1,2.8)--(2.3,2.8);
\draw [ultra thick] (1,0)--(-0.8,0); \draw [ultra thick] (2.3,0)--(3.8,0);
\draw (1,2.8) node{$_{\bullet}$}; \draw (2.3,2.8) node{$_{\bullet}$};
\draw [dotted] (1,0)--(1,2.8);
\draw [dotted] (2.3,0)--(2.3,2.8);
\end{tikzpicture}
\begin{center}
    $k\chi_{[\frac{1}{k},\frac{2}{k}]}(x)$
\end{center}  
\end{figure}

The sequence of functions $(f_k)$ converges to $0$ $\lambda$ almost uniformly, but it does not converge uniformly to $0$ on $[0,1]$.
\end{example}

\begin{proof}
Let $\delta >0$. Choose $k_{\delta} \in \mathbb{N}$ such that $\frac{2}{k_{\delta}} < \delta$. Define the Borel measurable set $F:=[0,\frac{2}{k_{\delta}}]$; clearly $\lambda(F)<\delta$. Thus, for every $\varepsilon >0$, defining $k_0:=k_{\delta} \in \mathbb{N}$, we have
$$
\left|f_k(x) \right|<\varepsilon \qquad \forall k \geq k_0,\quad\forall x \in [0,1]\smallsetminus F = (\tfrac{2}{k_{0}},1].
$$

\begin{figure}[ht!]
\begin{minipage}[c]{0.3\textwidth}
\begin{center}
\begin{tikzpicture}[xscale=0.85,yscale=0.85]
	\draw [fill=red!15,  dotted]  (2.5,0)--(2.5,4)--(0,4)--(0,0);
	\draw[->,gray] (0,0) -- (4,0); \draw [->,gray] (0,0) -- (0,4); 
	\draw[-,gray] (4,0)--(0,0); \draw[-,gray] (0,0)--(0,4);
	\draw[-,gray] (-0.5,0)--(0,0); \draw[-,gray] (0,0)--(0,-0.5);
\draw (0,2.25) node{$-$}; \draw (0,2.25) node[left] {$_{k_\delta}$}; \draw (0,-0.2) node[left] {$_{0}$};
\draw (2.5,0) node{$|$}; \draw (2.5,-0.27) node[below] {$_{\frac{2}{k_{\delta}}}$};
\draw (1.8,0) node{$|$}; \draw (1.8,-0.27) node[below] {$_{\frac{1}{k_{\delta}}}$};
\draw (3.5,0) node{$|$}; \draw (3.5,-0.2) node[below] {$_{1}$};
\draw [ultra thick] (1.8,2.25)--(2.5,2.25);
\draw [ultra thick] (0,0)--(1.8,0);
\draw [ultra thick] (2.5,0)--(3.5,0);
\draw [-, dotted] (1.8,2.25)--(1.8,0);
\draw [-, dotted] (2.5,2.25)--(2.5,0);

\end{tikzpicture}
\end{center}
\end{minipage} \hfill 
 \begin{minipage}[c]{0.3\textwidth}
\begin{center}
\begin{tikzpicture}[xscale=0.85,yscale=0.85]
	\draw [fill=red!15,  dotted]  (2.5,0)--(2.5,4)--(0,4)--(0,0);
	\draw[->,gray] (0,0) -- (4,0); \draw [->,gray] (0,0) -- (0,4); 
	\draw[-,gray] (4,0)--(0,0); \draw[-,gray] (0,0)--(0,4);
	\draw[-,gray] (-0.5,0)--(0,0); \draw[-,gray] (0,0)--(0,-0.5);

\draw (0,2.5) node{$-$}; \draw (0,2.5) node[left] {$_{k_\delta +1}$}; \draw (0,-0.2) node[left] {$_{0}$};
\draw (1.75,0) node{$|$}; \draw (1.75,-0.27) node[below] {$_{\frac{2}{k_{\delta}+1}}$};
\draw (1,0) node{$|$}; \draw (1,-0.27) node[below] {$_{\frac{1}{k_{\delta}+1}}$};
\draw (3.5,0) node{$|$}; \draw (3.5,-0.2) node[below] {$_{1}$};
\draw [ultra thick] (1,2.5)--(1.75,2.5);
\draw [-, dotted] (1,2.5)--(1,0);
\draw [-, dotted] (1.75,2.5)--(1.75,0);
\draw [ultra thick] (0,0)--(1,0);
\draw [ultra thick] (1.75,0)--(3.5,0);
\draw (2.5,0) node{$|$}; \draw (2.5,-0.27) node[below] {$_{\frac{2}{k_{\delta}}}$};

\end{tikzpicture}
\end{center}
\end{minipage} \hfill
\begin{minipage}[c]{0.3\textwidth}
\begin{center}
\begin{tikzpicture}[xscale=0.85,yscale=0.85]
	\draw [fill=red!15,  dotted]  (2.5,0)--(2.5,4)--(0,4)--(0,0);
	\draw[->,gray] (0,0) -- (4,0); \draw [->,gray] (0,0) -- (0,4); 
	\draw[-,gray] (4,0)--(0,0); \draw[-,gray] (0,0)--(0,4);
	\draw[-,gray] (-0.5,0)--(0,0); \draw[-,gray] (0,0)--(0,-0.5);

\draw (0,2.75) node{$-$}; \draw (0,2.75) node[left] {$_{k_\delta +2}$}; \draw (0,-0.2) node[left] {$_{0}$};
\draw (1.25,0) node{$|$}; \draw (1.25,-0.27) node[below] {$_{\frac{2}{k_{\delta}+2}}$};
\draw (0.5,0) node{$|$}; \draw (0.5,-0.27) node[below] {$_{\frac{1}{k_{\delta}+2}}$};
\draw (3.5,0) node{$|$}; \draw (3.5,-0.2) node[below] {$_{1}$};
\draw [ultra thick] (0.5,2.75)--(1.25,2.75);
\draw [-, dotted] (0.5,2.75)--(0.5,0);
\draw [-, dotted] (1.25,2.75)--(1.25,0);
\draw [ultra thick] (0,0)--(0.5,0);
\draw [ultra thick] (1.25,0)--(3.5,0);
\draw (2.5,0) node{$|$}; \draw (2.5,-0.27) node[below] {$_{\frac{2}{k_{\delta}}}$};
\end{tikzpicture}
\end{center}
\end{minipage}
\end{figure}

Consequently, $f_k \to 0$ $\lambda$ almost uniformly. On the other hand, $(f_k)$ does not converge uniformly to $0$, since if $\varepsilon_{0} \in (0,1)$, there is no $k \in \mathbb{N}$ such that $|f_k(x)|< \varepsilon_{0}$ for all $x \in [0,1]$. Indeed,
$$
\left| f_k\left( \frac{1}{k} \right) \right|=k  > \varepsilon_{0},\quad\forall\, k \in \mathbb{N}.
$$

This concludes the proof.
\end{proof}

The notion of almost uniform convergence should not be interpreted as being equivalent to the concept of uniform convergence almost everywhere. That is, there exists a sequence of measurable functions $(f_k)$ that converges $\mu$ almost uniformly, but there does not exist a $\mu$-null subset $N$ of $X$ such that $(f_k)$ converges uniformly on $X \smallsetminus N$ [Exercise \ref{E810}].

\begin{definition} \label{87} \index{convergence!in mean $p$}
Let $p \in [1,\infty)$. A sequence of $\mathsf{S}$-measurable functions $f_k:X \to \mathbb{R}$ in $L^{p}(X,\mathsf{S},\mu)$ converges \textbf{in mean $p$} to a function $f:X \to \mathbb{R}$ in $L^{p}(X,\mathsf{S},\mu)$ if, for every $\varepsilon >0$, there exists $k_0 \in \mathbb{N}$ (which may depend on $\varepsilon$) such that
$$
\Vert f_k - f \Vert_{p} < \varepsilon \quad\forall k \geq k_0.
$$

We denote this type of convergence by $f_{k} \xrightarrow[L^{p}]{} f$.
\end{definition}

We now consider a simple example.

\begin{example} \label{88}
Fix $p \in [1,\infty)$. Let $X=[0,1]$, $\mathsf{S}=\mathcal{B}([0,1])$ and let $\mu=\lambda$ be the Lebesgue measure on $\mathsf{S}$. Let $f_k:[0,1]\to \mathbb{R}$ be the function defined by $f_k:=\frac{1}{k}$.

Since $\lambda([0,1])=1$, it is clear that $f_k \in L^{p}([0,1],\mathcal{B}([0,1]),\lambda)$ for every $k \in \mathbb{N}$ (see Theorem \ref{713}). Consequently,
$$
\|f_k \|_{p}^{p}=\int_{[0,1]} \frac{1}{k^{p}}\,d\lambda = \frac{1}{k^{p}} \qquad \forall k \in \mathbb{N}.
$$

Taking limits as $k \to \infty$ in the previous identity, we conclude that $(f_k)$ converges in mean $p$ to $0$.
\end{example}

We now analyze the relationships between the different notions of convergence introduced.

\begin{theorem}[Mean $p$ convergence implies convergence in measure] \label{89}
Let $(X,\mathsf{S},\mu)$ be a measure space and fix $p \in [1,\infty)$. If $(f_k)$ is a sequence of $\mathsf{S}$-measurable functions that converges in mean $p$ to $f \in L^{p}(X,\mathsf{S},\mu)$, then $f_{k} \xrightarrow[\mu]{} f$.
\end{theorem}

\begin{proof}
Let $\xi>0$ be arbitrary. For each $k \in \mathbb{N}$ define the measurable set
$$
A_{k}(\xi):=\left\{x \in X\,:\,|f_k(x)-f(x)|\geq \xi \right\}.
$$

Then,
$$
\xi\chi_{A_{k}(\xi)} \leq |f_k-f|\chi_{A_{k}(\xi)} \leq |f_k-f|\chi_{X}
$$
and it follows that
$$
\xi^{p}\chi_{A_{k}(\xi)} \leq |f_k-f|^{p}\chi_{A_{k}(\xi)} \leq |f_k-f|^{p}\chi_{X},
$$
where $|f_k-f|^{p} \in L^{p}(X,\mathsf{S},\mu)$ by hypothesis.

By monotonicity of the integral we obtain
$$
\int_{A_{k}(\xi)}\xi^{p}\,d\mu \leq \int_{A_{k}(\xi)}|f_k-f|^{p}\,d\mu \leq \int_X |f_k-f|^{p}\,d\mu.
$$

Consequently,
$$
\mu(A_k(\xi))\xi^{p} \leq \int_{A_{k}(\xi)}|f_k-f|^{p}\,d\mu \leq  \Vert f_k - f\Vert_{p}^{p},
$$
which implies that
$$
0\leq \mu(A_k(\xi))\leq \frac{1}{\xi^{p}} \Vert f_k - f\Vert_{p}^{p}, \qquad\forall k \in \mathbb{N}.
$$ 
   
Taking limits as $k \to \infty$ in the previous inequality we conclude that
$$
0 \leq \lim_{k \to \infty}\mu(A_k(\xi)) \leq \lim_{k \to \infty}\frac{1}{\xi^{p}} \Vert f_k - f\Vert_{p}^{p} =0.
$$
   
Therefore, $f_{k} \xrightarrow[\mu]{} f$.
\end{proof}

The converse is not true in general, as shown by the following example.

\begin{example} \label{810}
Let $X=[0,1]$, $\mathsf{S}=\mathcal{B}([0,1])$ and let $\mu=\lambda$ be the Lebesgue measure on $\mathsf{S}$. The sequence of functions $f_k:[0,1] \to \mathbb{R}$ defined by $f_k:=k\chi_{[0,\frac{1}{k}]}$ converges to $0$ in measure, but does not converge in mean $p$ for any $p \in [1,\infty)$.
\end{example}

\begin{proof}
Example \ref{84} shows that $f_{k} \xrightarrow[\lambda]{} 0$.

On the other hand, let $p \in [1,\infty)$. Arguing by contradiction, suppose that $f_k \to f$ in $L^{p}([0,1],\mathcal{B}([0,1]),\lambda)$. Then $(f_k)$ is a Cauchy sequence in $(L^{p},\Vert \cdot \Vert_{p})$. However, for $k \in \mathbb{N}$, we have
$$
\Vert f_{2k}-f_{k}\Vert_{p}^{p}=\int_{[0,1]} |f_{2k}-f_{k}|^{p}\,d\lambda
$$
where
$$
|f_{2k}-f_{k}|^{p}(x)=\left\{
\begin{array}{lcl}
k^{p} & & \mbox{if }\,\,\,0 \leq x < \frac{1}{k},\\
0 & & \mbox{if }  \,\,\, \frac{1}{k} \leq x < 1.
\end{array}
\right.
$$

\begin{figure}[ht!]
\centering
\begin{tikzpicture}[xscale=1,yscale=0.85]
	\draw[->,gray] (0,0) -- (4,0); \draw [->,gray] (0,0) -- (0,4.5); 
	\draw[<-,gray] (-2.2,0)--(0,0); \draw[-,gray] (0,0)--(0,-0.23);
\draw (4,0) node[right] {$[0,1]$}; 

\draw (0,3) node{$_{-}$}; \draw (0,3) node[left]{$_{k}$};
\draw (0,4) node{$_{-}$}; \draw (0,4) node[left]{$_{2k}$};
\draw (1.5,0) node{$_{|}$}; \draw (1.5,-0.1) node[below]{$_{\frac{1}{k}}$};
\draw (0.8,0) node{$_{|}$}; \draw (0.8,-0.1) node[below]{$_{\frac{1}{2k}}$};
\draw (3,0) node{$_{|}$}; \draw (3,-0.1) node[below]{$_{1}$};
\draw[ultra thick] (0,3)--(1.47,3); \draw (0,3) node{$_{\bullet}$}; \draw (1.5,3) node{$_{\circ}$};
\draw[ultra thick] (0,4)--(0.75,4); \draw (0,4) node{$_{\bullet}$}; \draw (0.8,4) node{$_{\circ}$};
\draw[ultra thick] (0,0)--(0,0); 
\draw[ultra thick] (3,0)--(4,0); 
\draw[ultra thick] (-2.2,0)--(0,0); 
\draw[ultra thick] (1.5,0)--(3,0); 
\draw (0,-0.1) node[below]{$_{0}$};
\draw [dotted] (1.5,0)--(1.5,3);
\draw [dotted] (0.8,0)--(0.8,4);
\end{tikzpicture}
\begin{center}
 $|f_{2k}(x)-f_{k}(x)|$
\end{center}  
\end{figure}

Therefore,
$$
\Vert f_{2k}-f_{k}\Vert_{p}^{p}=\int_{[0,1]} |f_{2k}-f_{k}|^{p}\,d\lambda =k^{p} \lambda([0,\tfrac{1}{k}))=k^{p-1} \quad\forall k \in \mathbb{N}.
$$

Taking $\varepsilon_{0} \in (0,1)$ it is clear that
$$
\Vert f_{2k}-f_{k}\Vert_{p}^{p} \geq \varepsilon_{0}  \quad\forall k \in \mathbb{N}.
$$

This contradicts our assumption.
\end{proof}

In the previous example, the sequence of functions $(f_k)$ converges $\lambda$ almost everywhere to $0$ [Exercise \ref{E83}], so we conclude that almost everywhere convergence does not imply convergence in mean $p$ either. We recall, however, the dominated convergence theorem in $L^{p}$ (see Theorem \ref{719}).

\begin{theorem}[Dominated convergence in $L^{p}$] \label{811}
Let $p \in [1,\infty)$ and let $(f_k)$ be a sequence in $L^{p}(X,\mathsf{S},\mu)$ such that $f_k\to f$ $\mu$ almost everywhere on $X$. If there exists $g \in L^{p}(X,\mathsf{S},\mu)$ such that $|f_k| \leq g$ $\mu$ almost everywhere on $X$ for all $k \in \mathbb{N}$, then $f \in L^{p}(X,\mathsf{S},\mu)$ and
$$
\lim_{k \to \infty}\Vert f_k-f\Vert_{p}=0.
$$
\end{theorem}

Although convergence in mean $p$ does not, in general, imply almost everywhere convergence, the Riesz–Fischer theorem guarantees the existence of a subsequence that converges almost everywhere, as stated in Corollary \ref{721}. In the next section we present a result that generalizes this assertion to the case of convergence in measure.

\section{Riesz–Weyl Theorem}

Let $(X,\mathsf{S},\mu)$ be a measure space.

The result studied in this section clarifies fundamental aspects of the relationship between convergence in mean $p$ and almost everywhere convergence. Moreover, it will serve as a bridge to other notions of convergence, such as convergence in measure and almost uniform convergence. With this objective, we introduce the following definition.

\begin{definition} \label{812} \index{sequence!Cauchy in measure}
Let $(f_k)$ be a sequence of functions in $\mathbb{M}(X,\mathsf{S})$. We say that $(f_k)$ is \textbf{Cauchy in measure $\mu$} if, for every $\xi >0$ and every $\varepsilon >0$, there exists $k_{0} \in \mathbb{N}$ such that
$$
\mu\left( \left\{ x \in X\,:\,|f_k(x)-f_j(x)|\geq \xi \right\}\right)<\varepsilon \quad\forall k,j \geq k_{0}.
$$

Equivalently, for every $\xi >0$,
$$
\lim_{k,j \to \infty} \mu\left( \left\{ x \in X\,:\,|f_k(x)-f_j(x)|\geq \xi \right\}\right) =0.
$$
\end{definition}

\begin{theorem} \label{813}
Every sequence $(f_k)$ in $\mathbb{M}(X,\mathsf{S})$ that converges in measure is Cauchy in measure $\mu$.
\end{theorem}

\begin{proof}
Let $(f_k)$ be a sequence in $\mathbb{M}(X,\mathsf{S})$ such that $f_{k} \xrightarrow[\mu]{} f$. Then, for every $\xi >0$, we have
$$
\lim_{k\to \infty} \mu\left( \left\{ x \in X\,:\,|f_k(x)-f(x)|\geq \frac{\xi}{2} \right\}\right) =0.
$$

Let $A_{k}(\xi)$ denote the measurable set
$$
A_{k}(\xi):=\left\{ x \in X\,:\,|f_k(x)-f(x)|\geq \xi \right\}.
$$

Observe that if $x \in (X\smallsetminus A_{k}(\xi/2)) \cap (X\smallsetminus A_{j}(\xi/2))$, then
$$
\begin{aligned}
|f_k(x)-f(x)| & < \frac{\xi}{2},\\
|f_j(x)-f(x)| & < \frac{\xi}{2},
\end{aligned}
$$
and therefore
$$
|f_k(x)-f_j(x)| \leq |f_k(x)-f(x)|+|f(x)-f_j(x)| < \frac{\xi}{2} + \frac{\xi}{2} = \xi.
$$

Thus,
$$
x \in \{x \in X\,:\,|f_k(x)-f_{j}(x)| < \xi \} = X \smallsetminus A_{k,j}(\xi),
$$
where
$$
A_{k,j}(\xi):=\left\{ x \in X\,:\,|f_k(x)-f_j(x)|\geq \xi \right\}.
$$

Consequently,
$$
(X\smallsetminus A_{k}(\xi/2)) \cap (X\smallsetminus A_{j}(\xi/2)) \subset X\smallsetminus A_{k,j}(\xi),
$$
and taking complements we obtain
$$
A_{k,j}(\xi) \subset A_{k}(\xi/2) \cup A_{j}(\xi/2).
$$

By monotonicity of the measure,
$$
\mu (A_{k,j}(\xi)) \leq \mu (A_{k}(\xi/2)) + \mu (A_{j}(\xi/2)),
$$
and taking the limit as $k,j \to \infty$, we conclude that
$$
\lim_{k,j \to \infty} \mu (A_{k,j}(\xi)) \leq \lim_{k,j \to \infty}\left[ \mu (A_{k}(\xi/2)) + \mu (A_{j}(\xi/2)) \right] =0
\quad \forall \xi >0.
$$

That is, $(f_k)$ is Cauchy in measure $\mu$.
\end{proof}

We now consider an interesting example.

\begin{example} \label{814}
Let $X=[0,1)$, $\mathsf{S}=\mathcal{B}([0,1))$ and let $\mu=\lambda$ be the Lebesgue measure on $\mathsf{S}$. We define the following sequence of intervals:

$$
\begin{array}{lcll}
L_1 & & I_{1}:=[0,1) &\\
L_{2} & & I_{2}:=[0,1/2)\,\,\,\,I_{3}:=[1/2,1)&\\
L_{3} & & I_{4}:=[0,1/4)\,\,\,\,I_{5}:=[1/4,1/2)&I_{6}:=[1/2,3/4)\,\,\,\,I_{7}:=[3/4,1)\\
\vdots & & & \\
L_{k} & & I_{2^{k-1}}:=[0,1/2^{k-1}) \,\,\ldots \ldots & I_{2^k-1}:=[1-1/2^{k-1},1).
\end{array}
$$

We now define the following sequence of functions:
$$
\begin{array}{lcll}
L_1 & & f_{1}:=\chi_{I_1} &\\
L_{2} & & f_{2}:=\chi_{I_2}\,\,\,\,f_{3}:=\chi_{I_3}&\\
L_{3} & & f_{4}:=\chi_{I_4}\,\,\,\,f_{5}:=\chi_{I_5}&f_{6}:=\chi_{I_6}\,\,\,\,f_{7}:=\chi_{I_7}\\
\vdots & & & \\
L_{k} & & f_{2^{k-1}}:=\chi_{I_{2^{k-1}}} \,\,\ldots \ldots & f_{2^k-1}:=\chi_{I_{2^k-1}}
\end{array}
$$

\begin{figure}[ht!]
    \centering
	\begin{tikzpicture}[xscale=0.8,yscale=0.8]
	\draw [fill=lime!15,  dotted]  (3,0)--(3,2.8)--(0,2.8)--(0,0);
	\draw[->,gray] (0,0) -- (4,0); \draw [->,gray] (0,0) -- (0,4); 
	\draw[-,gray] (4,0)--(0,0); \draw[-,gray] (0,0)--(0,4);
	\draw[-,gray] (-0.5,0)--(0,0); \draw[-,gray] (0,0)--(0,-0.5);
\draw (4.3,0.25) node[below] {$x$}; \draw (0,4.1) node[right] {$f_{1}(x)$}; 
\draw (0,2.8) node{$-$}; \draw (0,2.8) node[left] {$_{1}$}; \draw (0,-0.2) node[left] {$_{0}$};
\draw (3,0) node{$|$}; \draw (3,-0.2) node[below] {$_{1}$};
\draw [ultra thick] (0,2.8)--(3,2.8);
\draw [-, dotted] (3,2.8)--(3,0);
\end{tikzpicture}
\end{figure}

\begin{figure}[ht!]
\begin{minipage}[c]{0.3\textwidth}
\begin{center}
\begin{tikzpicture}[xscale=0.8,yscale=0.8]
	\draw [fill=lime!15,  dotted]  (1.5,0)--(1.5,2.8)--(0,2.8)--(0,0);
	\draw[->,gray] (0,0) -- (4,0); \draw [->,gray] (0,0) -- (0,4); 
	\draw[-,gray] (4,0)--(0,0); \draw[-,gray] (0,0)--(0,4);
	\draw[-,gray] (-0.5,0)--(0,0); \draw[-,gray] (0,0)--(0,-0.5);
\draw (4.3,0.25) node[below] {$x$}; \draw (0,4.1) node[right] {$f_{2}(x)$}; 
\draw (0,2.8) node{$-$}; \draw (0,2.8) node[left] {$_{1}$}; \draw (0,-0.2) node[left] {$_{0}$};
\draw (1.5,0) node{$|$}; \draw (1.5,-0.27) node[below] {$_{\frac{1}{2}}$};
\draw (3,0) node{$|$}; \draw (3,-0.2) node[below] {$_{1}$};
\draw [ultra thick] (0,2.8)--(1.5,2.8);
\draw [ultra thick] (1.5,0)--(3,0);
\end{tikzpicture}
\end{center}
\end{minipage} \hfill 
 \begin{minipage}[c]{0.3\textwidth}
\begin{center}
\begin{tikzpicture}[xscale=0.8,yscale=0.8]
	\draw [fill=lime!15,  dotted]  (3,0)--(3,2.8)--(1.5,2.8)--(1.5,0);
	\draw[->,gray] (0,0) -- (4,0); \draw [->,gray] (0,0) -- (0,4); 
	\draw[-,gray] (4,0)--(0,0); \draw[-,gray] (0,0)--(0,4);
	\draw[-,gray] (-0.5,0)--(0,0); \draw[-,gray] (0,0)--(0,-0.5);
\draw (4.3,0.25) node[below] {$x$}; \draw (0,4.1) node[right] {$f_{3}(x)$}; 
\draw (0,2.8) node{$-$}; \draw (0,2.8) node[left] {$_{1}$}; \draw (0,-0.2) node[left] {$_{0}$};
\draw (1.5,0) node{$|$}; \draw (1.5,-0.27) node[below] {$_{\frac{1}{2}}$};
\draw (3,0) node{$|$}; \draw (3,-0.2) node[below] {$_{1}$};
\draw [ultra thick] (1.5,2.8)--(3,2.8);
\draw [ultra thick] (0,0)--(1.5,0);
\draw [-, dotted] (3,2.8)--(3,0);
\end{tikzpicture}
\end{center}
\end{minipage} \hfill 
\begin{minipage}[c]{0.3\textwidth}
\begin{center}
\begin{tikzpicture}[xscale=0.8,yscale=0.8]
	\draw [fill=lime!15,  dotted]  (1,0)--(1,2.8)--(0,2.8)--(0,0);
	\draw[->,gray] (0,0) -- (4.5,0); \draw [->,gray] (0,0) -- (0,4); 
	\draw[-,gray] (4.5,0)--(0,0); \draw[-,gray] (0,0)--(0,4);
	\draw[-,gray] (-0.5,0)--(0,0); \draw[-,gray] (0,0)--(0,-0.5);
\draw (4.8,0.25) node[below] {$x$}; \draw (0,4.1) node[right] {$f_{4}(x)$}; 
\draw (0,2.8) node{$-$}; \draw (0,2.8) node[left] {$_{1}$}; \draw (0,-0.2) node[left] {$_{0}$};
\draw (1,0) node{$|$}; \draw (1,-0.27) node[below] {$_{\frac{1}{4}}$};
\draw (2,0) node{$|$}; \draw (2,-0.27) node[below] {$_{\frac{2}{4}}$};
\draw (3,0) node{$|$}; \draw (3,-0.27) node[below] {$_{\frac{3}{4}}$};
\draw (4,0) node{$|$}; \draw (4,-0.2) node[below] {$_{1}$};
\draw [ultra thick] (0,2.8)--(1,2.8);
\draw [ultra thick] (1,0)--(4,0);
\end{tikzpicture}
\end{center}
\end{minipage}  
\end{figure}

\begin{figure}[ht!]
\begin{minipage}[c]{0.3\textwidth}
\begin{center}
\begin{tikzpicture}[xscale=0.8,yscale=0.8]
	\draw [fill=lime!15,  dotted]  (2,0)--(2,2.8)--(1,2.8)--(1,0);
	\draw[->,gray] (0,0) -- (4.5,0); \draw [->,gray] (0,0) -- (0,4); 
	\draw[-,gray] (4.5,0)--(0,0); \draw[-,gray] (0,0)--(0,4);
	\draw[-,gray] (-0.5,0)--(0,0); \draw[-,gray] (0,0)--(0,-0.5);
\draw (4.8,0.25) node[below] {$x$}; \draw (0,4.1) node[right] {$f_{5}(x)$}; 
\draw (0,2.8) node{$-$}; \draw (0,2.8) node[left] {$_{1}$}; \draw (0,-0.2) node[left] {$_{0}$};
\draw (1,0) node{$|$}; \draw (1,-0.27) node[below] {$_{\frac{1}{4}}$};
\draw (2,0) node{$|$}; \draw (2,-0.27) node[below] {$_{\frac{2}{4}}$};
\draw (3,0) node{$|$}; \draw (3,-0.27) node[below] {$_{\frac{3}{4}}$};
\draw (4,0) node{$|$}; \draw (4,-0.2) node[below] {$_{1}$};
\draw [ultra thick] (1,2.8)--(2,2.8);
\draw [ultra thick] (0,0)--(1,0);
\draw [ultra thick] (2,0)--(4,0);
\end{tikzpicture}
\end{center}
\end{minipage} \hfill \begin{minipage}[c]{0.3\textwidth}
\begin{center}
\begin{tikzpicture}[xscale=0.8,yscale=0.8]
	\draw [fill=lime!15,  dotted]  (2,0)--(2,2.8)--(3,2.8)--(3,0);
	\draw[->,gray] (0,0) -- (4.5,0); \draw [->,gray] (0,0) -- (0,4); 
	\draw[-,gray] (4.5,0)--(0,0); \draw[-,gray] (0,0)--(0,4);
	\draw[-,gray] (-0.5,0)--(0,0); \draw[-,gray] (0,0)--(0,-0.5);
\draw (4.8,0.25) node[below] {$x$}; \draw (0,4.1) node[right] {$f_{6}(x)$}; 
\draw (0,2.8) node{$-$}; \draw (0,2.8) node[left] {$1$}; \draw (0,-0.2) node[left] {$_0$};
\draw (1,0) node{$|$}; \draw (1,-0.27) node[below] {$_{1/4}$};
\draw (2,0) node{$|$}; \draw (2,-0.27) node[below] {$_{2/4}$};
\draw (3,0) node{$|$}; \draw (3,-0.27) node[below] {$_{3/4}$};
\draw (4,0) node{$|$}; \draw (4,-0.2) node[below] {$1$};
\draw [ultra thick] (2,2.8)--(3,2.8);
\draw [ultra thick] (0,0)--(2,0);
\draw [ultra thick] (3,0)--(4,0);
\end{tikzpicture}
\end{center}
\end{minipage} \hfill 
 \begin{minipage}[c]{0.3\textwidth}
\begin{center}
\begin{tikzpicture}[xscale=0.8,yscale=0.8]
	\draw [fill=lime!15,  dotted]  (3,0)--(3,2.8)--(4,2.8)--(4,0);
	\draw[->,gray] (0,0) -- (4.5,0); \draw [->,gray] (0,0) -- (0,4); 
	\draw[-,gray] (4.5,0)--(0,0); \draw[-,gray] (0,0)--(0,4);
	\draw[-,gray] (-0.5,0)--(0,0); \draw[-,gray] (0,0)--(0,-0.5);
\draw (4.8,0.25) node[below] {$x$}; \draw (0,4.1) node[right] {$f_{7}(x)$}; 
\draw (0,2.8) node{$-$}; \draw (0,2.8) node[left] {$1$}; \draw (0,-0.2) node[left] {$_0$};
\draw (1,0) node{$|$}; \draw (1,-0.27) node[below] {$_{1/4}$};
\draw (2,0) node{$|$}; \draw (2,-0.27) node[below] {$_{2/4}$};
\draw (3,0) node{$|$}; \draw (3,-0.27) node[below] {$_{3/4}$};
\draw (4,0) node{$|$}; \draw (4,-0.2) node[below] {$1$};
\draw [ultra thick] (3,2.8)--(4,2.8);
\draw [ultra thick] (0,0)--(3,0);
\end{tikzpicture}
\end{center}
\end{minipage}
\end{figure}

The sequence $(f_k)$ converges in measure $\lambda$ to $0$, but it does not converge to $0$ $\lambda$ almost everywhere.
\end{example}

\begin{proof}
Let $\xi >0$. We consider the following cases.

{\scshape Case 1.}\quad It is clear that if $\xi >1$, then the set
$\left\{x \in [0,1)\,:\, |f_k(x)| \geq \xi \right\}$
is empty for every $k \in \mathbb{N}$. Hence,
$$
\lambda\left(\left\{x \in [0,1)\,:\, |f_k(x)| \geq \xi \right\} \right)=0 \qquad\forall k \in \mathbb{N}.
$$

{\scshape Case 2.}\quad Suppose that $0 < \xi \leq 1$.

For each $j \in \mathbb{N}$ and each $k \in \mathbb{N}$ such that $I_{k} \in L_{j}$, we have $2^{j-1} \leq k \leq 2^{j}-1$ and therefore $\lambda(I_{k})=\frac{1}{2^{j-1}}$.

Consequently,
$$
\lambda\left(\left\{x \in [0,1)\,:\, |f_k(x)| \geq \xi \right\} \right)=\lambda(I_{k})=\frac{1}{2^{j-1}}\leq\frac{2}{k},
$$
since
$$
2^{j-1} \leq k \leq 2^{j}-1 = 2\cdot 2^{j-1}-1 \leq 2\cdot 2^{j-1},
$$
which implies
$$
\frac{1}{2^{j-1}} \leq \frac{2}{k} \leq \frac{2}{2^{j-1}}.
$$

Taking limits as $k \to \infty$, we conclude that
$$
0 \leq \lim_{k \to \infty}\lambda\left(\left\{x \in [0,1)\,:\, |f_k(x)| \geq \xi \right\} \right)\leq\lim_{k\to \infty}\frac{2}{k}=0.
$$

That is, $f_{k} \xrightarrow[\lambda]{} 0$.

However, note that $f_{k} \nrightarrow 0$ $\lambda$ almost everywhere. Indeed, let $\varepsilon_{0} \in (0,1)$. For each $k \in \mathbb{N}$,
$$
\{x \in [0,1)\,:\, |f_{k}(x) | > \varepsilon_{0}\}=I_{k}.
$$

Let $x \in [0,1)$. For each $k \in \mathbb{N}$, there exists a unique $j_{k} \in \mathbb{N}$ such that $x \in I_{j_{k}}$ and $2^{j_{k}-1} \leq j_{k} \leq 2^{j_{k}}-1$. Hence $\limsup_{k \to \infty} I_{k}=[0,1)$ and
$$
\lambda\left(\limsup_{k \to \infty} \{x \in [0,1)\,:\, |f_{k}(x) | > \varepsilon_{0}\}\right)
=\lambda\left( \limsup_{k \to\infty}I_{k}\right) = 1 .
$$

Exercise \ref{E86} implies that $f_{k} \nrightarrow 0$ $\lambda$ almost everywhere.
\end{proof}

From the previous example we conclude that convergence in measure does not imply almost everywhere convergence. It is a simple exercise to show that the sequence $(f_k)$ from Example \ref{814} converges in mean $p$ [Exercise \ref{E85}], so we also conclude that convergence in mean $p$ does not imply almost everywhere convergence.

The following result, due to F. Riesz and H. Weyl\footnote{Hermann Weyl (1885–1955) was a German mathematician and one of the most influential mathematicians of the 20th century. He made fundamental contributions to differential geometry, group theory, relativity, and quantum mechanics, and developed modern representation theory as well as a deep reflection on the foundations of mathematics. He was a student of David Hilbert and professor in Zürich, Göttingen, and Princeton.}, establishes a criterion relating sequences $(f_k)$ of measurable functions that converge in measure with almost everywhere convergence. This result also guarantees that every sequence that is Cauchy in measure is convergent in measure.

\begin{theorem}[Riesz–Weyl] \label{815} \index{theorem!Riesz Weyl}
Let $(X,\mathsf{S},\mu)$ be a measure space and let $(f_k)$ be a sequence of functions in $\mathbb{M}(X,\mathsf{S})$ that is Cauchy in measure $\mu$. Then there exist $f \in \mathbb{M}(X,\mathsf{S})$ and a subsequence $(f_{k_j})$ of $(f_k)$ such that
\begin{itemize}
\item[(a)] $f_{k_{j}} \to f$ $\mu$ almost everywhere,
\item[(b)] $f_{k} \xrightarrow[\mu]{} f$.
\end{itemize}
\end{theorem}

\begin{figure}[ht!]
\begin{minipage}[c]{0.5\textwidth}
\begin{center}
\includegraphics[scale=0.3]{riez.jpeg} 
\begin{center}
F. Riesz (1880-1956)
\end{center}
\end{center}
\end{minipage} \hfill \begin{minipage}[c]{0.5\textwidth}
\begin{center}
\includegraphics[scale=0.3]{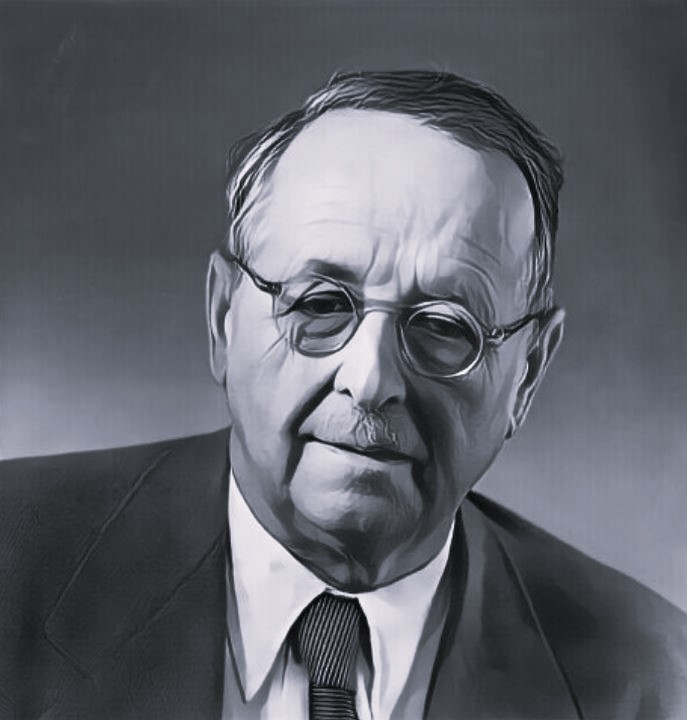} 
\begin{center}
H. Weyl (1885-1955)
\end{center}
\end{center}
\end{minipage}
\end{figure}

\begin{proof}
Since the proof is long, we divide it into several steps.

{\scshape Step 1:} Let $\xi_{1}=\varepsilon_{1}:=\frac{1}{2}>0$. Since $(f_k)$ is Cauchy in measure $\mu$, there exists $k_{1} \in \mathbb{N}$ such that
$$
\mu\left( \left\{ x \in X\,:\, |f_k(x)-f_{i}(x)| \geq \frac{1}{2} \right\}\right) < \frac{1}{2} \qquad\forall k,i \geq k_1.
$$

Similarly, for $\xi_{2}=\varepsilon_{2}:=\frac{1}{2^2}>0$, there exists $k_{2,\varepsilon_2} \in \mathbb{N}$ such that
$$
\mu\left( \left\{ x \in X\,:\, |f_k(x)-f_{i}(x)| \geq \frac{1}{2^2} \right\}\right) < \frac{1}{2^2} \qquad\forall k,i \geq k_{2,\varepsilon_2}.
$$

Defining $k_{2}:=\max\{k_{1},k_{2,\varepsilon_2} \}+1 \in \mathbb{N}$, it is clear that $k_{1} < k_{2}$ and that it satisfies
$$
\mu\left( \left\{ x \in X\,:\, |f_k(x)-f_{i}(x)| \geq \frac{1}{2^2} \right\}\right) < \frac{1}{2^2} \qquad\forall k,i \geq k_{2}.
$$

Continuing in this way, we obtain that for each $j \in \mathbb{N}$ there exists $k_{j} \in \mathbb{N}$ with $k_{j} >k_{j-1}>\cdots >k_{1}$ such that
$$
\mu\left(\left\{x \in X\,:\,|f_{k}(x)-f_{i}(x)| \geq \frac{1}{2^{j}} \right\} \right)<\frac{1}{2^{j}} \quad\forall k,i \geq k_{j}.
$$

Therefore, there exists a subsequence $(f_{k_j})$ of $(f_k)$ such that
$$
\sum_{j=1}^{\infty} \mu\left(\left\{x \in X\,:\,|f_{k_{j+1}}(x)-f_{k_{j}}(x)| \geq \frac{1}{2^{j}} \right\} \right) < \sum_{j=1}^{\infty} \frac{1}{2^j}=1.
$$

{\scshape Step 2:} For each $j \in \mathbb{N}$, define the measurable set
$$
A_j:=\left\{x \in X\,:\,|f_{k_{j+1}}(x)-f_{k_{j}}(x)| \geq \frac{1}{2^{j}} \right\}.
$$

The Borel–Cantelli lemma (see Theorem \ref{422}) implies that
$$
\mu\left( \limsup_{j \to \infty} A_{j} \right)=0
$$
since $\sum_{j=1}^{\infty} \mu(A_j) < +\infty$.

Notice that, for each $x \in X \smallsetminus \limsup_{j \to \infty} A_{j}$, there exists $j(x) \in \mathbb{N}$ (depending on $x$) such that $x \in X \smallsetminus A_{i}$ for all $i \geq j(x)$.

Let $\varepsilon >0$. Choose $j_0 \in \mathbb{N}$ such that $\frac{1}{2^{j_0}} < \varepsilon$. Thus, given $x \in X \smallsetminus \limsup_{j \to \infty} A_{j}$, set $j_{\ast}:=\max\{j_0,j(x)\} \in \mathbb{N}$. Then, for any $i > j > j_{\ast}$, we have
$$
\begin{aligned}
\left|f_{k_j}(x)-f_{k_i}(x) \right|
&\leq  \left|f_{k_j}(x)-f_{k_{j+1}}(x) \right|
+ \left|f_{k_{j+1}}(x)-f_{k_{j+2}}(x) \right|
+ \ldots
+ \left|f_{k_{i-1}}(x)-f_{k_{i}}(x) \right| \\
&< \frac{1}{2^j}+\frac{1}{2^{j+1}}+\ldots + \frac{1}{2^{i-1}}\\
& < \frac{1}{2^{j-1}} \leq \frac{1}{2^{j_{\ast}}} < \varepsilon.
\end{aligned}
$$

Consequently, $(f_{k_j}(x))$ is a Cauchy sequence in $\mathbb{R}$ for every $x \in X \smallsetminus \limsup_{j \to \infty} A_{j}$.

Defining $f:X \to \mathbb{R}$ by
$$
f(x):=\left\{
\begin{array}{lcl}
\lim_{j \to \infty} f_{k_j}(x) & & \mbox{if}\quad x \in X \smallsetminus \limsup_{j \to \infty} A_{j},\\
0 & & \mbox{if}\quad x \in \limsup_{j \to \infty} A_{j},
\end{array}
\right.
$$
it is clear that $f \in \mathbb{M}(X,\mathsf{S})$ and that $f_{k_j} \to f$ $\mu$ almost everywhere.

{\scshape Step 3:} Let $\xi >0$ and $\varepsilon >0$. Choose $j_{0} \in \mathbb{N}$ such that $\frac{1}{2^{j_0-1}} < \xi$ and $\frac{1}{2^{j_0-1}} < \varepsilon$. Consider $x \in X \smallsetminus A_{j}$ for all $j \geq j_0$. Then, for $i > j \geq j_{0}$, we have
\begin{eqnarray} \label{F81}
\left|f_{k_j}(x)-f_{k_i}(x) \right|< \frac{1}{2^{j_{0}-1}} < \xi.
\end{eqnarray}

Taking the limit as $i \to \infty$ in (\ref{F81}) yields
$$
\left|f_{k_{j}}(x)-f(x) \right| \leq \frac{1}{2^{j_0-1}} < \xi \qquad\forall j \geq j_{0}.
$$

Thus, $X \smallsetminus A_{j} \subset \{x \in X\,:\, |f_{k_{j}}(x)-f(x)| < \xi\}$ for all $j \geq j_0$, which implies
$$
\{x \in X\,:\, |f_{k_{j}}(x)-f(x)| \geq \xi\} \subset A_{j} \qquad \forall j \geq j_0.
$$

By monotonicity of the measure $\mu$, we conclude that
$$
\mu\left( \left\{ x \in X\,:\, |f_{k_{j}}(x)-f(x)| \geq \xi \right\} \right)
\leq \mu (A_j)
< \frac{1}{2^{j}} < \frac{1}{2^{j-1}} \leq \frac{1}{2^{j_{0}-1}} < \varepsilon
$$
for all $j \geq j_0$. That is, $f_{k_j} \xrightarrow[\mu]{} f$.

{\scshape Step 4:} Since $(f_k)$ is Cauchy in measure $\mu$ and $(f_{k_j})$ is a subsequence of $(f_k)$ that converges in measure to $f$, it follows that $f_k \xrightarrow[\mu]{} f$ [Exercise \ref{E813}].

This concludes the proof.
\end{proof}

Let us return to our previous example and observe the usefulness of Theorem \ref{815}.

\begin{example} \label{816}
Let $X=[0,1)$, $\mathsf{S}=\mathcal{B}([0,1))$ and let $\mu=\lambda$ be the Lebesgue measure on $\mathsf{S}$, with the sequence of functions $(f_k)$ given in Example \ref{814}, namely,
$$
\begin{array}{lcll}
L_{1} & & f_{1}:=\chi_{I_1} &\\
L_{2} & & f_{2}:=\chi_{I_2}\,\,\,\,f_{3}:=\chi_{I_3}&\\
L_{3} & & f_{4}:=\chi_{I_4}\,\,\,\,f_{5}:=\chi_{I_5}&f_{6}:=\chi_{I_6}\,\,\,\,f_{7}:=\chi_{I_7}\\
\vdots & & & \\
L_{k} & & f_{2^{k-1}}:=\chi_{I_{2^{k-1}}} \,\,\ldots \ldots & f_{2^k-1}:=\chi_{I_{2^k-1}}
\end{array}
$$

Theorem \ref{815} ensures that there exists a subsequence of $(f_k)$ which converges to $f=0$ $\lambda$ almost everywhere. We now construct this subsequence explicitly.

The function $k:\mathbb{N} \to \mathbb{N}$ defined by
$$
k(j):=\left\{
\begin{array}{lcl}
1 & & \mbox{if } j=1,\\
2^{j-1} & & \mbox{if } j>1,
\end{array}
\right.
$$
is strictly increasing on $\mathbb{N}$. We consider the subsequence $(f_{k_j})$ of $(f_k)$. That is, $(f_{k_j})$ is given by
$$
\begin{array}{c c c c c c}
f_{k_{1}}:=\chi_{I_1} & f_{k_{2}}:=\chi_{I_2}& f_{k_{3}}:=\chi_{I_4}&
\ldots & f_{k_{j}}:=\chi_{I_{2^{j-1}}}.
\end{array}
$$

\begin{figure}[ht!]
\begin{minipage}[c]{0.5\textwidth}
\begin{center}
\begin{tikzpicture}[xscale=0.8,yscale=0.8]
\draw [fill=lime!15,  dotted]  (3,0)--(3,2.8)--(0,2.8)--(0,0);
	\draw[->,gray] (0,0) -- (4,0); \draw [->,gray] (0,0) -- (0,4); 
	\draw[-,gray] (4,0)--(0,0); \draw[-,gray] (0,0)--(0,4);
	\draw[-,gray] (-0.5,0)--(0,0); \draw[-,gray] (0,0)--(0,-0.5);
\draw (4.3,0.25) node[below] {$x$}; \draw (0,4.1) node[right] {$f_{k_{1}}(x)$}; 
\draw (0,2.8) node{$-$}; \draw (0,2.8) node[left] {$_{1}$}; \draw (0,-0.2) node[left] {$_{0}$};
\draw (3,0) node{$|$}; \draw (3,-0.2) node[below] {$_{1}$};
\draw [ultra thick] (0,2.8)--(3,2.8);
\draw [-, dotted] (3,2.8)--(3,0);
\end{tikzpicture}
\end{center}
\end{minipage} \hfill 
 \begin{minipage}[c]{0.5\textwidth}
\begin{center}
\begin{tikzpicture}[xscale=0.8,yscale=0.8]
	\draw [fill=lime!15,  dotted]  (1.5,0)--(1.5,2.8)--(0,2.8)--(0,0);
	\draw[->,gray] (0,0) -- (4,0); \draw [->,gray] (0,0) -- (0,4); 
	\draw[-,gray] (4,0)--(0,0); \draw[-,gray] (0,0)--(0,4);
	\draw[-,gray] (-0.5,0)--(0,0); \draw[-,gray] (0,0)--(0,-0.5);
\draw (4.3,0.25) node[below] {$x$}; \draw (0,4.1) node[right] {$f_{k_{2}}(x)$}; 
\draw (0,2.8) node{$-$}; \draw (0,2.8) node[left] {$_{1}$}; \draw (0,-0.2) node[left] {$_{0}$};
\draw (1.5,0) node{$|$}; \draw (1.5,-0.27) node[below] {$_{\frac{1}{2}}$};
\draw (3,0) node{$|$}; \draw (3,-0.2) node[below] {$_{1}$};
\draw [ultra thick] (0,2.8)--(1.5,2.8);
\draw [ultra thick] (1.5,0)--(3,0);
\end{tikzpicture}
\end{center}
\end{minipage}
\end{figure}

\begin{figure}[ht!]
\begin{minipage}[c]{0.3\textwidth}
\begin{center}
\begin{tikzpicture}[xscale=0.8,yscale=0.75]
\draw [fill=lime!15,  dotted]  (1,0)--(1,2.8)--(0,2.8)--(0,0);
	\draw[->,gray] (0,0) -- (4.5,0); \draw [->,gray] (0,0) -- (0,4); 
	\draw[-,gray] (4.5,0)--(0,0); \draw[-,gray] (0,0)--(0,4);
	\draw[-,gray] (-0.5,0)--(0,0); \draw[-,gray] (0,0)--(0,-0.5);
\draw (4.8,0.25) node[below] {$x$}; \draw (0,4.1) node[right] {$f_{k_{3}}(x)$}; 
\draw (0,2.8) node{$-$}; \draw (0,2.8) node[left] {$_{1}$}; \draw (0,-0.2) node[left] {$_{0}$};
\draw (1,0) node{$|$}; \draw (1,-0.27) node[below] {$_{\frac{1}{4}}$};
\draw (4,0) node{$|$}; \draw (4,-0.2) node[below] {$_{1}$};
\draw [ultra thick] (0,2.8)--(1,2.8);
\draw [ultra thick] (1,0)--(4,0);
\end{tikzpicture}
\end{center}
\end{minipage} \hfill 
 \begin{minipage}[c]{0.3\textwidth}
\begin{center}
\begin{tikzpicture}[xscale=0.8,yscale=0.75]
	\draw [fill=lime!15,  dotted]  (0.5,0)--(0.5,2.8)--(0,2.8)--(0,0);
	\draw[->,gray] (0,0) -- (4.5,0); \draw [->,gray] (0,0) -- (0,4); 
	\draw[-,gray] (4.5,0)--(0,0); \draw[-,gray] (0,0)--(0,4);
	\draw[-,gray] (-0.5,0)--(0,0); \draw[-,gray] (0,0)--(0,-0.5);
\draw (4.8,0.25) node[below] {$x$}; \draw (0,4.1) node[right] {$f_{k_{4}}(x)$}; 
\draw (0,2.8) node{$-$}; \draw (0,2.8) node[left] {$_{1}$}; \draw (0,-0.2) node[left] {$_{0}$};
\draw (0.5,0) node{$|$}; \draw (0.5,-0.27) node[below] {$_{\frac{1}{8}}$};
\draw (4,0) node{$|$}; \draw (4,-0.2) node[below] {$_{1}$};
\draw [ultra thick] (0,2.8)--(0.5,2.8);
\draw [ultra thick] (0.5,0)--(4,0);
\end{tikzpicture}
\end{center}
\end{minipage} \hfill 
\begin{minipage}[c]{0.3\textwidth}
\begin{center}
\begin{tikzpicture}[xscale=0.8,yscale=0.75]
		\draw [fill=lime!15,  dotted]  (0,0)--(0,2.8)--(0.25,2.8)--(0.25,0);
	\draw[->,gray] (0,0) -- (4.5,0); \draw [->,gray] (0,0) -- (0,4); 
	\draw[-,gray] (4.5,0)--(0,0); \draw[-,gray] (0,0)--(0,4);
	\draw[-,gray] (-0.5,0)--(0,0); \draw[-,gray] (0,0)--(0,-0.5);
\draw (4.8,0.25) node[below] {$x$}; \draw (0,4.1) node[right] {$f_{k_{5}}(x)$}; 
\draw (0,2.8) node{$-$}; \draw (0,2.8) node[left] {$_{1}$}; \draw (0,-0.2) node[left] {$_{0}$};
\draw (0.25,0) node{$|$}; \draw (0.25,-0.27) node[below] {$_{\frac{1}{16}}$};
\draw (4,0) node{$|$}; \draw (4,-0.2) node[below] {$_{1}$};
\draw [ultra thick] (0,2.8)--(0.25,2.8);
\draw [ultra thick] (0.25,0)--(4,0);
\end{tikzpicture}
\end{center}
\end{minipage}  
\end{figure}

We claim that $f_{k_j}(x) \to 0$ $\lambda$ almost everywhere.
\end{example}

\begin{proof}
Let $\varepsilon >0$ and $x \in (0,1)$. There exists $j_{0} \in \mathbb{N}\smallsetminus\{1\}$ such that $\frac{1}{2^{j_0-1}}< x$. Therefore, $f_{k_j}(x)=0$ for all $j \geq j_{0}$ and hence
$$
|f_{k_j}(x)|<\varepsilon \qquad\forall j \geq j_{0}.
$$

If $x=0$, then $x \in [0,\tfrac{1}{2^{j-1}})$ for every $j \in \mathbb{N}$. Hence,
$$
f_{k_j}(x)=1 \qquad\forall j \in \mathbb{N},
$$
that is, $f_{k_j}(0) \nrightarrow 0$.

We conclude that $f_{k_j}(x) \to 0$ for every $x \in [0,1)\smallsetminus\{0\}$ and, since $\lambda(\{0\})=0$, it follows that $f_{k_j}(x) \to 0$ $\lambda$ almost everywhere, as claimed.
\end{proof}

We now study some properties of convergence in measure.

\begin{proposition} \label{817}
Let $(f_{k})$, $(g_{k})$ be sequences of functions in $\mathbb{M}(X,\mathsf{S})$, and let $f,g \in \mathbb{M}(X,\mathsf{S})$.

\item[\textit{(a)}] If $f_{k}\xrightarrow[\mu]{} f$ and $f_{k}\xrightarrow[\mu]{} g$, then $f=g$ $\mu$ a.e.

\item[\textit{(b)}] If $f_{k}\xrightarrow[\mu]{} f$ and $f_{k} \to g$ $\mu$ a.e., then $f=g$ $\mu$ a.e.

\item[\textit{(c)}] If $f_{k}\xrightarrow[\mu]{} f$ and $g_{k}=f_{k}$ $\mu$ a.e. for all $k \in \mathbb{N}$, then $g_{k}\xrightarrow[\mu]{} f$.
\end{proposition}

\begin{proof}
\textit{(a):} It suffices to prove that
$$
\mu \left( \left\{ x \in X\,:\, |f(x)-g(x)| \geq \varepsilon \right\} \right)=0
$$
for every $\varepsilon >0$. This follows from the following inclusion, valid for every $\varepsilon >0$:
$$
\begin{aligned}
\left\{ x \in X\,:\, |f(x)-g(x)| \geq \varepsilon \right\}
&\subset \left\{ x \in X\,:\, |f(x)-f_k(x)| \geq \frac{\varepsilon}{2} \right\} \\
&\quad \cup \left\{ x \in X\,:\, |f_k(x)-g(x)| \geq \frac{\varepsilon}{2} \right\}.
\end{aligned}
$$

By monotonicity of the measure, we obtain
$$
\begin{aligned}
0 \leq \mu\left(\left\{ x \in X\,:\, |f(x)-g(x)| \geq \varepsilon \right\} \right)
&\leq \mu\left(\left\{ x \in X\,:\, |f(x)-f_k(x)| \geq \frac{\varepsilon}{2} \right\} \right) \\
&\quad + \mu\left(\left\{ x \in X\,:\, |f_k(x)-g(x)| \geq \frac{\varepsilon}{2} \right\} \right).
\end{aligned}
$$

Thus, since $f_{k} \xrightarrow[\mu]{} f$ and $f_{k} \xrightarrow[\mu]{} g$, taking the limit as $k \to \infty$ yields the result.

\textit{(b):} Since $f_{k} \xrightarrow[\mu]{} f$, the sequence $(f_k)$ is Cauchy in measure $\mu$. By the Riesz–Weyl theorem, there exists $h \in \mathbb{M}(X,\mathsf{S})$ and a subsequence $(f_{k_j})$ of $(f_k)$ such that $f_{k_j} \to h$ $\mu$ a.e. and $f_k \xrightarrow[\mu]{} h$. Since $f_{k_j} \to g$ $\mu$ a.e., it follows that $h=g$ $\mu$ a.e., and by part (a), since also $f_k \xrightarrow[\mu]{} f$, we conclude that $f=h$ $\mu$ a.e. Hence $f=g$ $\mu$ a.e.

\textit{(c):} For each $k \in \mathbb{N}$, let $N_k \in \mathcal{N}(\mu)$ be such that $f_k(x)=g_k(x)$ for all $x \in X\smallsetminus N_k$. For $\xi >0$ arbitrary, setting $N:=\bigcup_{k=1}^{\infty} N_k$, we have
$$
\left\{x \in X\,:\,|g_k(x)-f(x)|\geq\xi \right\}
\subset N \cup \left\{x \in X\,:\,|f_k(x)-f(x)|\geq\xi \right\} \quad \forall k \in \mathbb{N}.
$$

By monotonicity of the measure,
$$
\begin{aligned}
\mu\left( \left\{x \in X\,:\,|g_k(x)-f(x)|\geq\xi \right\} \right)
\leq \mu(N) + \mu\left( \left\{x \in X\,:\,|f_k(x)-f(x)|\geq\xi \right\} \right)
\quad \forall k\in\mathbb{N}.
\end{aligned}
$$

Taking the limit as $k \to \infty$, we conclude that
$$
\lim_{k \to \infty}\mu\left( \left\{x \in X\,:\,|g_k(x)-f(x)|\geq\xi \right\} \right)=0,
$$
which shows that $g_{k} \xrightarrow[\mu]{} f$.
\end{proof}

In general, convergence in measure does not imply convergence in mean $p$ for $p \in [1,\infty)$. However, the Riesz–Weyl theorem allows us to obtain the following result.

\begin{theorem}[Dominated convergence in measure] \label{818} \index{theorem!dominated convergence in measure}
Let $(X,\mathsf{S},\mu)$ be a measure space, $p \in [1,\infty)$, and let $(f_k)$ be a sequence of functions in $L^{p}(X,\mathsf{S},\mu)$ such that $f_{k} \xrightarrow[\mu]{} f$ with $f \in \mathbb{M}(X,\mathsf{S})$. If there exists a non-negative function $g \in L^{1}(X,\mathsf{S},\mu)$ such that $|f_k|^{p} \leq g$ $\mu$ almost everywhere for all $k \in \mathbb{N}$, then $f \in L^{p}(X,\mathsf{S},\mu)$ and $f_{k} \xrightarrow[L^{p}]{} f$.
\end{theorem}

\begin{proof}
Arguing by contradiction, suppose that $f_{k} \not\to f$ in mean $p$. Then there exists a subsequence $(f_{k_j})$ of $(f_k)$ and $\varepsilon>0$ such that
$$
\Vert f_{k_j}-f\Vert_{p} \geq \varepsilon \quad\forall\, j \in \mathbb{N}.
$$

Moreover, $f_{k_j}\xrightarrow[\mu]{} f$ since $f_k \xrightarrow[\mu]{} f$ [Exercise \ref{E814}]. Hence, by the Riesz–Weyl theorem, there exists a subsequence $(f_{k_{j_\ell}})$ of $(f_{k_j})$ and a function $h \in \mathbb{M}(X,\mathsf{S})$ such that $f_{k_{j_{\ell}}} \to h$ $\mu$ a.e. and $f_{k_j} \xrightarrow[\mu]{} h$. Thus, $f=h$ $\mu$ a.e., and therefore $f_{k_{j_{\ell}}} \to f$ $\mu$ a.e. Moreover, $|f_{k_{j_\ell}}|^{p} \leq g$. Hence $|f_{k_{j_\ell}}-f|^{p} \to 0$ $\mu$ a.e., and
$$
|f_{k_{j_\ell}}-f|^{p}\leq (|f_{k_{j_\ell}}|+|f|)^{p} \leq (2g^{1/p})^{p}=2^{p}g,
$$
with $2^{p}g \in L^{1}(X,\mathsf{S},\mu)$. By the dominated convergence theorem of Lebesgue,
$$
\lim_{\ell \to \infty}\int_X |f_{k_{j_\ell}}-f|^{p}\,d\mu = 0,
$$
which contradicts our assumption.
\end{proof}

Next, we study some monotonicity properties that are preserved under convergence in measure.

\begin{theorem} \label{819}
Let $(f_k)$ be a sequence of functions in $\mathbb{M}(X,\mathsf{S})$, and let $f,g \in \mathbb{M}(X,\mathsf{S})$.
\begin{itemize}
\item[(a)] If $f_{k} \xrightarrow[\mu]{} f$ with $f_k \geq 0$ $\mu$ a.e., then $f \geq 0$ $\mu$ a.e.
\item[(b)] If $f_{k} \xrightarrow[\mu]{} f$ and $f_k \leq g$ $\mu$ a.e. for all $k \in \mathbb{N}$, then $f \leq g$ $\mu$ a.e.
\end{itemize}
\end{theorem}

\begin{proof}
\textit{(a):} Since $f_k \geq 0$ $\mu$ a.e., there exists $N_k \in \mathcal{N}(\mu)$ such that $f_k(x)\geq 0$ for all $x \in X\smallsetminus N_{k}$. Define, for each $k \in \mathbb{N}$,
$$
\tilde{f}_k(x):=\left\{
\begin{array}{lcl}
f_k(x) & & \mbox{if } x \in X\smallsetminus N_k,\\
0 & & \mbox{if } x \in N_k.
\end{array}
\right.
$$

Thus, $\tilde{f}_k=f_k$ $\mu$ a.e. for all $k \in \mathbb{N}$ and, consequently, $\tilde{f}_k \xrightarrow[\mu]{} f$.

To show that $f \geq 0$ $\mu$ a.e., it suffices to prove that
$$
\mu \left(\left\{ x \in X\,:\, f(x) \leq -\varepsilon\right\} \right)=0 \qquad\forall \varepsilon > 0.
$$

If $f(x) \leq -\varepsilon$, then since $f(x)=(f(x)-\tilde{f}_k(x))+\tilde{f}_k(x)$ and $\tilde{f}_k(x)\geq 0$, we obtain $f(x)-\tilde{f}_k(x) \leq -\varepsilon,$ and in particular $|f(x)-\tilde{f}_k(x)| \geq \varepsilon$.

Hence,
$$
\left\{ x \in X\,:\, f(x) \leq -\varepsilon \right\}
\subset \left\{ x \in X\,:\, |f(x)-\tilde{f}_k(x)| \geq \varepsilon \right\},
$$
and therefore
$$
\mu\left(\left\{ x \in X\,:\, f(x) \leq -\varepsilon \right\} \right)
\leq \mu \left( \left\{ x \in X\,:\, |f(x)-\tilde{f}_k(x)| \geq \varepsilon \right\}\right).
$$

Taking the limit as $k \to \infty$ we conclude that
$$
\mu\left(\left\{ x \in X\,:\, f(x) \leq -\varepsilon \right\} \right)=0 \quad\forall \varepsilon>0.
$$

\textit{(b):} Since $f_k \leq g$ $\mu$ a.e., we have $g-f_k \geq 0$ $\mu$ a.e., and moreover $g-f_k \xrightarrow[\mu]{} g-f$ [Exercise \ref{E815}]. By part \textit{(a)}, we conclude that $g-f \geq 0$ $\mu$ a.e., which proves the result.
\end{proof}

Next, we establish the relationships between almost uniform convergence, convergence in measure, and almost everywhere convergence. It is natural, after a first course in real analysis, to state that almost uniform convergence implies almost everywhere convergence.

\begin{theorem} \label{820}
Let $(f_k)$ be a sequence of functions in $\mathbb{M}(X,\mathsf{S})$ and let $f \in \mathbb{M}(X,\mathsf{S})$. If $f_{k} \to f$ $\mu$ a.u., then
\begin{itemize}
\item[(a)] $f_k \to f$ $\mu$ a.e.,
\item[(b)] $f_{k} \xrightarrow[\mu]{} f$.
\end{itemize}
\end{theorem}

\begin{proof}
\textit{(a):} For each $k \in \mathbb{N}$, there exists a set $F_k \in \mathsf{S}$ with $\mu(F_k) < \frac{1}{k}$ such that $f_k \to f$ uniformly on $X\smallsetminus F_k$. Define $F:=\bigcap_{k=1}^{\infty} F_k \in \mathsf{S}$. Then, for every $k \in \mathbb{N}$, $\mu(F)\leq \mu(F_k) < \frac{1}{k}$ and hence $\mu(F)=0$.

For $x \in X \smallsetminus F$, there exists $k(x) \in \mathbb{N}$ such that $x \in X\smallsetminus F_{k(x)}$. Since $f_k \to f$ uniformly on $X \smallsetminus F_{k(x)}$, Exercise \ref{Ej11} implies that $f_k(x) \to f(x)$ in $\mathbb{R}$. Therefore, $f_k(x) \to f(x)$ for all $x \in X \smallsetminus F$, which yields $f_k \to f$ $\mu$ a.e.

\textit{(b):} Let $\delta >0$ and $\varepsilon >0$ be given. There exists $F \in \mathsf{S}$ with $\mu(F)<\delta$ and $k_{\varepsilon} \in \mathbb{N}$ such that
$$
|f_k(x)-f(x)|<\varepsilon\quad\forall k \geq k_{\varepsilon},\quad\forall x\in X\smallsetminus F.
$$

Hence,
$$
\left\{x\in X\,:\, |f_k(x)-f(x)|\geq\varepsilon \right\}\subset F\quad\forall k \geq k_{\varepsilon},
$$
and by monotonicity of the measure,
$$
\mu\left( \left\{x\in X\,:\, |f_k(x)-f(x)|\geq\varepsilon \right\}\right) \leq \mu(F)<\delta\quad\forall k \geq k_{\varepsilon}.
$$

This shows that $f_{k} \xrightarrow[\mu]{} f$.
\end{proof}

We now present a final example showing that, in general, almost everywhere convergence does not imply almost uniform convergence.

\begin{example} \label{821}
Let $(\mathbb{R},\mathcal{B}(\mathbb{R}),\lambda)$ be a measure space. Define the sequence $f_k:\mathbb{R} \to \mathbb{R}$ by $f_{k}:=\chi_{[k,k+1]}$. Then the sequence $(f_k)$ converges to $0$ $\lambda$ almost everywhere, but it does not converge to $0$ $\lambda$ almost uniformly.
\end{example}

\begin{figure}[ht!]
    \centering
	\begin{tikzpicture}[xscale=1,yscale=0.7]
	\draw [fill=white,  dotted]  (2.3,0)--(2.3,2.8)--(1,2.8)--(1,0);
	\draw[->,gray] (0,0) -- (4,0); \draw [->,gray] (0,0) -- (0,4); 
	\draw[-,gray] (4,0)--(0,0); \draw[-,gray] (0,0)--(0,4);
	\draw[-,gray] (-1,0)--(0,0); \draw[-,gray] (0,0)--(0,-1);
\draw (0,2.8) node{$-$}; \draw (0,2.8) node[left] {$_{1}$};
\draw (1,0) node{$|$}; \draw (1,-0.27) node[below] {$_{k}$};
\draw (2.3,0) node{$|$}; \draw (2.3,-0.2) node[below] {$_{k+1}$};
\draw [ultra thick] (1,2.8)--(2.3,2.8);
\draw [ultra thick] (1,0)--(-0.8,0); \draw [ultra thick] (2.3,0)--(3.8,0);
\draw (1,2.8) node{$_{\bullet}$}; \draw (2.3,2.8) node{$_{\bullet}$};
\end{tikzpicture}
\begin{center}
$\chi_{[k,k+1]}(x)$
\end{center}
\end{figure}

\begin{proof}
It is a simple exercise to show that $f_k \to 0$ $\lambda$ almost everywhere [Exercise \ref{E88}].

Let $\delta_{0},\varepsilon_{0} \in (0,1)$ and let $F \in \mathcal{B}(\mathbb{R})$ be such that $\lambda(F)<\delta_{0}<1$. Then, for each $k \in \mathbb{N}$, there exists $x_{k} \in \mathbb{R}$ such that $x_{k} \in (\mathbb{R} \smallsetminus F) \cap [k,k+1]$. Hence,
$$
|f_{k}(x_{k})|=1 > \varepsilon_{0} \quad \forall k \in \mathbb{N},
$$
which implies that $f_{k} \not\to 0$ uniformly on $\mathbb{R} \smallsetminus F$. That is, $f_{k} \not\to 0$ $\lambda$ almost uniformly.
\end{proof}

From Example \ref{E814} and Theorem \ref{820}, we conclude that convergence in measure does not imply, in general, almost uniform convergence. However, by revisiting the proof of the Riesz–Weyl theorem, it is possible to obtain a partial relationship between these two notions of convergence.

\begin{theorem}[Riesz–Weyl] \label{822} \index{theorem!Riesz Weyl}
Let $(f_k)$ be a sequence of functions in $\mathbb{M}(X,\mathsf{S})$ that converges in measure to $f \in \mathbb{M}(X,\mathsf{S})$. Then there exists a subsequence $(f_{k_j})$ of $(f_k)$ such that $f_{k_j} \to f$ $\mu$ a.u.
\end{theorem}

\begin{proof}
Since $f_{k} \xrightarrow[\mu]{} f$, the sequence $(f_k)$ is Cauchy in measure $\mu$. Following the proof of Theorem \ref{815}, we conclude that there exists a subsequence $(f_{k_j})$ of $(f_k)$ such that
$$
\mu(A_{j})<\frac{1}{2^{j}} \quad \forall j \in \mathbb{N}
$$
where $A_{j}:=\{x \in X\,:\,|f_{k_{j+1}}(x)-f_{k_j}(x)|\geq 2^{-j}\}$.

Let $\delta>0$ and $\varepsilon >0$. Choose $j_{0} \in \mathbb{N}$ such that $\frac{1}{2^{j_0-1}}<\delta$ and $\frac{1}{2^{j_0-1}}<\varepsilon$. Define the set $F_{j_{0}}:=\bigcup_{\ell=j_{0}}^{\infty}A_{\ell} \in \mathsf{S}$. It is clear that $\mu(F_{j_{0}})< \frac{1}{2^{j_{0}-1}} < \delta$. Hence, for any $i > j > j_0$ and any $x \in X \smallsetminus F_{j_{0}}$, we have
\begin{eqnarray}
\left|f_{k_j}(x)-f_{k_i}(x) \right| < \frac{1}{2^{j_{0}-1}} < \varepsilon.
\end{eqnarray}

Therefore, $(f_{k_j})$ is uniformly Cauchy on $X \smallsetminus F_{j_{0}}$ and, by the uniform Cauchy convergence criterion [Exercise \ref{Ej13}], we conclude that $(f_{k_j})$ converges uniformly on $X \smallsetminus F_{j_0}$ to a function $\widetilde{f} \in \mathbb{M}(X,\mathsf{S})$. Thus, $f_{k_{j}} \to \widetilde{f}$ $\mu$ a.u. and, by Theorem \ref{820}, we conclude that $f_{k_j}\xrightarrow[\mu]{} \widetilde{f}$.

On the other hand, $f_{k_j}\xrightarrow[\mu]{} f$ [Exercise \ref{E814}], so $f=\widetilde{f}$ $\mu$ a.e. by Theorem \ref{817}. It follows that $f_{k_j} \to f$ $\mu$ a.u.
\end{proof}

The previous result motivates the following definition.

\begin{definition} \label{823} \index{sequence!Cauchy almost uniformly}
Let $(f_k)$ be a sequence of functions in $\mathbb{M}(X,\mathsf{S})$. We say that $(f_k)$ is \textbf{Cauchy $\mu$ almost uniformly} on $X$ if, for every $\delta>0$, there exists $F \in \mathsf{S}$ with $\mu(F)<\delta$ such that $(f_k)$ is uniformly Cauchy on $X\smallsetminus F$. That is, $(f_k)$ is Cauchy $\mu$ almost uniformly if, for every $\delta>0$, there exists $F \in \mathsf{S}$ with $\mu(F)<\delta$ such that for every $\varepsilon >0$, there exists $k_0 \in \mathbb{N}$ such that
$$
|f_k(x)-f_j(x)|<\varepsilon\quad\forall\, k,j \geq k_0,\quad\forall x \in X \smallsetminus F.
$$
\end{definition}

As one would expect, the Cauchy criterion holds for this notion of convergence, as stated in the following theorem.

\begin{theorem}[Cauchy criterion for almost uniform convergence] \label{824} \index{criterion!Cauchy criterion for almost uniform convergence}
Let $(f_{k})$ be a sequence of functions in $\mathbb{M}(X,\mathsf{S})$ and let $f \in \mathbb{M}(X,\mathsf{S})$. Then, $f_{k} \to f$ $\mu$ almost uniformly if and only if $(f_k)$ is Cauchy $\mu$ almost uniformly on $X$.
\end{theorem}

The proof is straightforward in view of the Riesz–Weyl theorem and is left as an exercise [Exercise \ref{E818}].

\section{Egorov's Theorem}

From a basic course in real analysis, we know that pointwise convergence of a sequence of functions does not imply uniform convergence [Exercise \ref{E81}]. In this section we state a result due to Dmitri Egorov\footnote{Dmitri Fedorovich Egorov (1869--1931) was a Russian mathematician known mainly for Egorov's theorem in measure theory and for his contributions to differential geometry and real analysis. He was a teacher of Nikolai Luzin and one of the founders of the Moscow School of Mathematics. In addition to his mathematical work, he suffered political and religious persecution in the Soviet Union and died in prison in 1931.} which, under suitable assumptions, asserts that a sequence that converges almost everywhere also converges almost uniformly and, by Theorem \ref{820}, converges in measure. Before proceeding to its proof, we give a characterization of the almost uniform convergence.

\begin{figure}[ht!]
\centering
\includegraphics[scale=0.3]{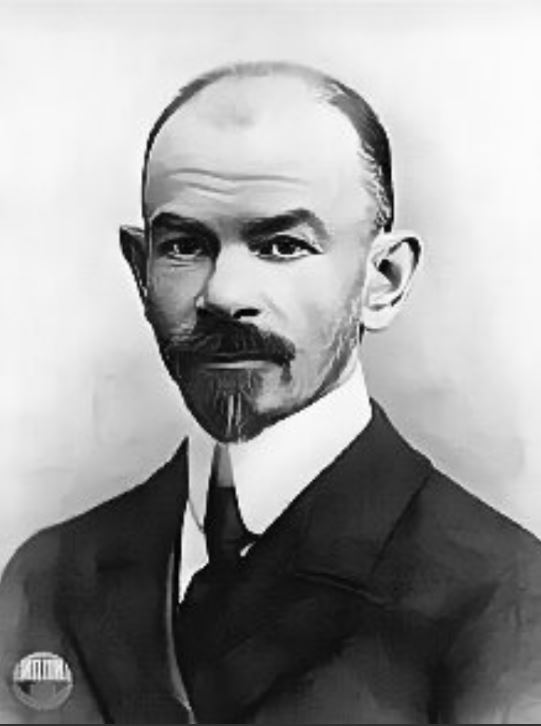} 
\begin{center}
D.F. Egorov (1869-1931)
\end{center}
\end{figure}

\begin{theorem} \label{825}
Let $(X,\mathsf{S},\mu)$ be a measure space, $(f_k)$ a sequence of functions in $\mathbb{M}(X,\mathsf{S})$, and $f \in \mathbb{M}(X,\mathsf{S})$. Then $f_k \to f$ $\mu$ almost uniformly if and only if for every fixed $\varepsilon > 0$ we have
$$
\lim_{k \to \infty} \mu (B_{k}(\varepsilon))=0,
$$
where
$$
B_{k}(\varepsilon)=\bigcup_{j=k}^{\infty}\left\{x \in X\,:\,|f_j(x)-f(x)|\geq \varepsilon \right\}.
$$
\end{theorem}

\begin{proof}
$\Rightarrow)$ Fix $\varepsilon>0$. Since $f_k \to f$ $\mu$ a.u., for every $\delta > 0$ there exists $F \in \mathsf{S}$ such that $\mu(F)<\delta$ and such that $f_k \to f$ uniformly on $X\smallsetminus F$. Hence, there exists $k_{0} \in \mathbb{N}$ such that
$$
|f_k(x)-f(x)|<\varepsilon \quad\forall k \geq k_{0},\quad\forall x \in X\smallsetminus F.
$$

Consequently,
$$
X\smallsetminus F \subset \bigcap_{k=k_{0}}^{\infty}\left\{ x \in X\,:\,|f_k(x)-f(x)|<\varepsilon \right\},
$$
and therefore
$$
\bigcup_{k=k_{0}}^{\infty}\left\{ x \in X\,:\,|f_k(x)-f(x)|\geq \varepsilon\right\}
=: B_{k_0}(\varepsilon)\subset F.
$$

It is clear that the sequence $(B_{k}(\varepsilon))$ is decreasing; hence, for all $k \geq k_{0}$,
$$
\mu(B_{k}(\varepsilon)) \leq \mu(B_{k_0}(\varepsilon)) \leq \mu(F) < \delta.
$$

Since $\delta>0$ was arbitrary, we conclude that $\mu(B_{k}(\varepsilon)) \to 0$.

$\Leftarrow)$ Let $\delta > 0$ be arbitrary and fix $\varepsilon>0$. Since $\mu(B_{k}(\varepsilon)) \to 0$, for each $j \in \mathbb{N}$ and $\delta_{j}:=\frac{\delta}{2^{j}} > 0$, there exists $k_{j} \in \mathbb{N}$ such that
$$
\mu(B_{k}(\varepsilon)) < \delta_{j}\quad\forall k \geq k_{j}.
$$

We may assume that $k_{1}<k_{2}<\cdots<k_{j}<k_{j+1}$ for all $j\in \mathbb{N}$. Define
$$
F:=\bigcup_{j=1}^{\infty} B_{k_j}(\varepsilon) \in \mathsf{S},
$$
then
$$
\mu(F) \leq \sum_{j=1}^{\infty} \mu(B_{k_j}(\varepsilon)) < \sum_{j=1}^{\infty} \frac{\delta}{2^{j}} = \delta.
$$

Moreover,
$$
X\smallsetminus F = \bigcap_{j=1}^{\infty} \bigcap_{k=k_j}^{\infty}\left\{x \in X\,:\,|f_k(x)-f(x)|<\varepsilon \right\}.
$$

Hence, for every $x \in X\smallsetminus F$ there exists $j$ such that $|f_k(x)-f(x)|<\varepsilon$ for all $k \geq k_j$, which shows that $f_k \to f$ uniformly on $X\smallsetminus F$.
\end{proof}

\begin{theorem}[Egorov] \label{826} \index{theorem!Egorov}
Let $(X,\mathsf{S},\mu)$ be a finite measure space, and let $(f_k)$ be a sequence of functions in $\mathbb{M}(X,\mathsf{S})$ and $f \in \mathbb{M}(X,\mathsf{S})$ such that $f_k \to f$ $\mu$ almost everywhere. Then $f_k \to f$ $\mu$ almost uniformly.
\end{theorem}

\begin{proof}
The set given by
$$
A:=\left\{x \in X\,:\,f_{k}(x)\to f(x) \right\},
$$
is an element of the $\sigma$-algebra $\mathsf{S}$ [Exercise \ref{E315}]. By hypothesis $\mu(X\smallsetminus A)=0$, and by properties of finite measures we have $0=\mu(X\smallsetminus A)=\mu(X)-\mu(A)$. Consequently, $\mu(X)=\mu(A)$.

Let $\varepsilon >0$. For each $j \in \mathbb{N}$, we define the measurable set
$$
C_{j}(\varepsilon):=\left\{ x \in X\,:\,|f_{j}(x)-f(x)|<\varepsilon \right\}.
$$

Let $x \in A$. Then there exists $k_{0} \in \mathbb{N}$ such that $|f_{k}(x)-f(x)|<\varepsilon$ for all $k \geq k_{0}$. That is, $x \in C_{k}(\varepsilon)$ for all $k \geq k_{0}$, which implies that
$$
x \in  \bigcap_{k=k_{0}}^{\infty} C_{k}(\varepsilon)\quad\text{for some }k_{0} \in \mathbb{N}.
$$

Therefore,
$$
A \subset \bigcup_{k=1}^{\infty}\bigcap_{j=k}^{\infty} C_{j}(\varepsilon)=\liminf_{k \to \infty} C_{k}(\varepsilon).
$$

Using properties of measures (see Corollary \ref{415}), we have
$$
\mu(X)=\mu(A) \leq \mu\left(\liminf_{k \to \infty} C_{k}(\varepsilon) \right)
= \lim_{k\,\to\,\infty} \mu\left( \bigcap_{j=k}^{\infty} C_{j}(\varepsilon)\right) \leq \mu(X),
$$
from which it follows that
$$
\mu(X)=\lim_{k\,\to\,\infty}\mu\left( \bigcap_{j=k}^{\infty} C_{j}(\varepsilon)\right).
$$

Then,
$$
\begin{aligned}
0=\lim_{k\,\to\,\infty}\left[  \mu\left( X \smallsetminus \bigcap_{j=k}^{\infty} C_{j}(\varepsilon)\right) \right]
=\lim_{k\,\to\,\infty}\left[  \mu\left( \bigcup_{j=k}^{\infty} (X\smallsetminus C_{j}(\varepsilon))\right) \right],
\end{aligned}
$$
where
$$
\bigcup_{j=k}^{\infty} (X\smallsetminus C_{j}(\varepsilon))
=\bigcup_{j=k}^{\infty}\left\{x \in X\,:\,|f_{j}(x)-f(x)|\geq \varepsilon \right\}.
$$

Therefore,
$$
\lim_{k\to\,\infty}\mu\left( \bigcup_{j=k}^{\infty}\left\{x \in X\,:\,|f_{j}(x)-f(x)|\geq \varepsilon \right\}  \right)=0.
$$

Theorem \ref{825} implies that $f_{k} \to f$ $\mu$ a.u.
\end{proof}

Example \ref{821} shows that the assumption $\mu(X)<+\infty$ is essential in Egorov’s theorem. Moreover, the theorem may fail if the limit function takes extended real values on a set of positive measure [Exercise \ref{E820}].

To conclude this section, we state Egorov’s almost uniform convergence theorem, which provides conditions under which a sequence of measurable functions that converges almost everywhere also converges almost uniformly, even when the underlying measure space is not necessarily finite. To prove this result, we will use the well-known Chebyshev inequality\footnote{Pafnuti Lvovich Chebyshev (1821–1894), also known as “Tchebychev” and in other transliterations of his surname, was a Russian mathematician who made fundamental contributions to number theory, probability, approximation theory, and mechanics. He is known for Bertrand–Chebyshev’s theorem, Chebyshev’s inequality, and Chebyshev polynomials. He also founded an important mathematical school in Saint Petersburg and had a major influence on the development of modern Russian mathematics.} [Exercise \ref{E613}].

\begin{figure}[ht!]
\centering
\includegraphics[scale=0.275]{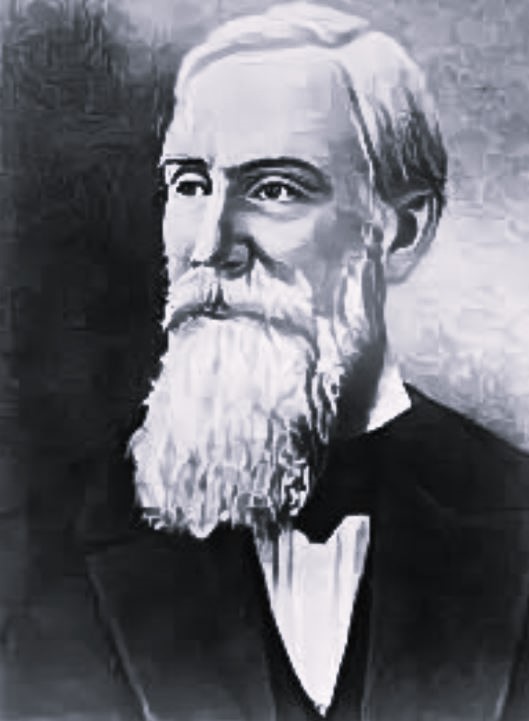} 
\begin{center}
P. Chebyshev (1821-1894)
\end{center}
\end{figure}

\begin{theorem}[Egorov's dominated convergence theorem] \label{827} \index{theorem!Egorov dominated convergence}
Let $(X,\mathsf{S},\mu)$ be a measure space and let $(f_k)$ be a sequence in $\mathbb{M}(X,\mathsf{S})$ such that $f_k \to f$ $\mu$ a.e., with $f \in \mathbb{M}(X,\mathsf{S})$. Assume that there exists $g \in L^{p}(X,\mathsf{S},\mu)$ with $p \in [1,\infty)$ such that $|f_k| \leq g$ $\mu$ a.e. for all $k \in \mathbb{N}$. Then $f_k \to f$ $\mu$ a.u.
\end{theorem}

\begin{proof}
Let $N$ be the $\mu$-null set given by
$$
N:=\left\{ x\in X\,:\,f_k(x) \nrightarrow f(x) \right\}.
$$

By hypothesis, $|f| \leq g$ on $X\smallsetminus N$ and therefore $|f_k-f| \leq |f_k|+|f| \leq 2g$ on $X\smallsetminus N$. For $\varepsilon>0$, define the measurable sets
$$
D_{j}(\varepsilon):=\left\{ x \in X\,:\,|f_j(x)-f(x)| \geq \varepsilon \right\}.
$$

From the above we conclude that
$$
D_{j}(\varepsilon)\smallsetminus N \subset \left\{x \in X\,:\,2g(x) \geq \varepsilon \right\}\quad\forall j \in \mathbb{N}.
$$

Thus,
$$
\bigcup_{j=1}^{\infty} D_{j}(\varepsilon)\smallsetminus N \subset \left\{x \in X\,:\,2g(x) \geq \varepsilon \right\},
$$
and applying P. Tchebyshev’s inequality [Exercise \ref{E613}] to the function $\varphi(t)=t^{p}$ we obtain
$$
\mu\left(\bigcup_{j=1}^{\infty} D_{j}(\varepsilon) \right) \leq \mu(\left\{x \in X\,:\,2g(x) \geq \varepsilon \right\}) \leq \frac{1}{\varepsilon^{p}}\int_X (2g)^{p}\,d\mu <+\infty.
$$

Since $f_k \to f$ $\mu$ a.e., we have
$$
\bigcap_{k=1}^{\infty}\bigcup_{j=k}^{\infty} D_{j}(\varepsilon) \in \mathcal{N}(\mu)
$$
and, since $\bigcup_{j=k}^{\infty} D_{j}(\varepsilon)$ has finite measure by the previous estimate, applying Theorem \ref{414} we conclude that
$$
\lim_{k\,\to\,\infty} \mu\left(\bigcup_{j=k}^{\infty} D_{j}(\varepsilon) \right)
=\mu\left(\bigcap_{k=1}^{\infty}\bigcup_{j=k}^{\infty} D_{j}(\varepsilon) \right)=0.
$$

Theorem \ref{825} implies that $f_k \to f$ $\mu$ a.u.
\end{proof}

\section{Relations between notions of convergence}

Let $X=(X,\mathsf{S},\mu)$ be a measure space.

In this brief section, as a conclusion and following the idea presented in \cite[Chapter 8]{Bartle}, we provide a graphical representation of the relationships between the different notions of convergence studied in this chapter.

In each figure, under the corresponding hypotheses, solid arrows between nodes represent direct implications between the different notions of convergence, while dashed arrows indicate relations that hold only after passing to a subsequence.

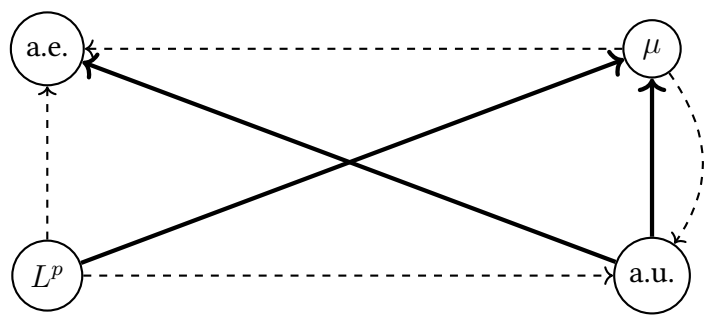
\begin{figure}[ht!]
    \centering
    \begin{tikzpicture}[xscale=2,yscale=1.5]

        \tikzstyle{vertex}=[circle, draw, thick]

        \node[vertex] (A) at (0,0) {$L^{p}$};
        \node[vertex] (B) at (4,0) {a.u.};
        \node[vertex] (C) at (4,2) {$\mu$};
        \node[vertex] (D) at (0,2) {a.e.};

        \draw[->, ultra thick] (A) -- (C);
        \draw[->, ultra thick] (B) -- (D);
        \draw[->, ultra thick] (B) -- (C);
        \draw[->, dashed, thick] (A) -- (D);
       \draw[->, dashed, thick] (C) -- (D);
       \draw[->, dashed, thick] (A) -- (B);
        \draw[->, dashed, thick]
    (C) to[out=-60,in=60] (B);
    \end{tikzpicture}
    \caption{The measure space is arbitrary.}
\end{figure}

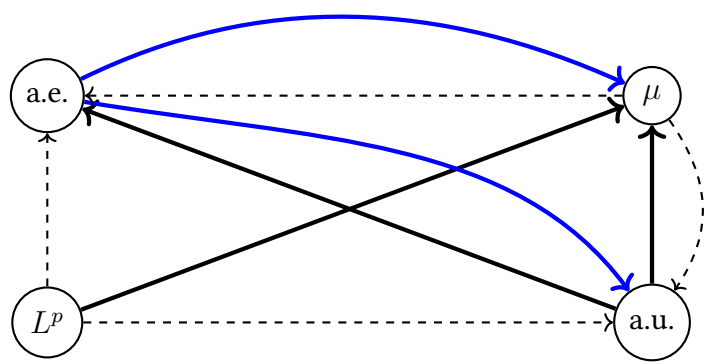
\begin{figure}[ht!]
    \centering
    \begin{tikzpicture}[xscale=2,yscale=1.5]

        \tikzstyle{vertex}=[circle, draw, thick]

        \node[vertex] (A) at (0,0) {$L^{p}$};
        \node[vertex] (B) at (4,0) {a.u.};
        \node[vertex] (C) at (4,2) {$\mu$};
        \node[vertex] (D) at (0,2) {a.e.};

        \draw[->, ultra thick] (A) -- (C);
        \draw[->, ultra thick] (B) -- (D);
        \draw[->, ultra thick] (B) -- (C);

        \draw[->, ultra thick, blue]
        (D) to[out=-12,in=120] (B);

        \draw[->, ultra thick, blue]
        (D) to[out=32,in=150] (C);
        \draw[->, dashed, thick] (A) -- (D);
       \draw[->, dashed, thick] (C) -- (D);
       \draw[->, dashed, thick] (A) -- (B);
        \draw[->, dashed, thick]
    (C) to[out=-60,in=60] (B);
    \end{tikzpicture}
    \caption{The measure space has finite measure.}
\end{figure}

\begin{figure}[ht!]
    \centering
    \begin{tikzpicture}[xscale=2,yscale=1.5]

        \tikzstyle{vertex}=[circle, draw, thick]

        \node[vertex] (A) at (0,0) {$L^{p}$};
        \node[vertex] (B) at (4,0) {a.u.};
        \node[vertex] (C) at (4,2) {$\mu$};
        \node[vertex] (D) at (0,2) {a.e.};

        \draw[->, ultra thick] (A) -- (C);
        \draw[->, ultra thick] (B) -- (D);
        \draw[->, ultra thick] (B) -- (C);

        \draw[->, ultra thick, red]
        (D) to[out=230,in=140] (A);

         \draw[->, ultra thick, red]
        (D) to[out=-12,in=120] (B);

        \draw[->, ultra thick, red]
        (C) to[out=150,in=90] (A);

        \draw[->, ultra thick, red]
        (D) to[out=32,in=150] (C);

        \draw[->, ultra thick, red]
        (B) to[out=200,in=-20] (A);
        \draw[->, dashed, thick] (A) -- (D);
       \draw[->, dashed, thick] (C) -- (D);
       \draw[->, dashed, thick] (A) -- (B);
        \draw[->, dashed, thick]
    (C) to[out=-30,in=30] (B);
    \end{tikzpicture}
    \caption{The measure space is arbitrary, and there exists a function $g \in L^{p}$ that dominates the sequence.}
\end{figure}
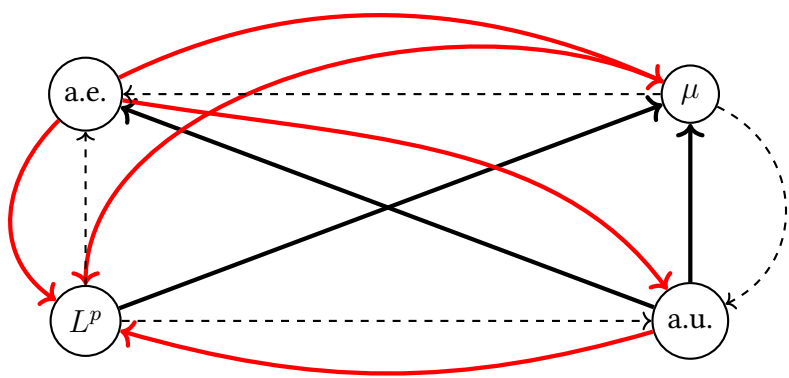

\section{Exercises}

\begin{exercise} \label{E81}
Let $X=[0,1]$, $\mathsf{S}=\mathcal{B}([0,1])$, and let $\mu=\lambda$ be the Lebesgue measure on $\mathsf{S}$. Prove the following statements:
\begin{itemize}
\item[(a)] The sequence of functions $f_k(x)=x^k$ converges pointwise on $[0,1]$.
\item[(b)] The sequence of functions $f_k(x)=x^k$ does not converge uniformly on $[0,1]$.
\item[(c)] Does the sequence of functions $f_k(x)=x^k$ converge almost uniformly on $[0,1]$?
\end{itemize}
\end{exercise}

\begin{exercise} \label{E82}
Let $X=\mathbb{R}$, $\mathsf{S}=\mathcal{B}(\mathbb{R})$, and let $\mu=\lambda$ be the Lebesgue measure on $\mathsf{S}$. Prove the following statements:
\begin{itemize}
\item[(a)] The sequence of functions $f_k(x)=\chi_{(k,k+\frac{1}{k}]}$ converges pointwise on $\mathbb{R}$.
\item[(b)] The sequence of functions $f_k(x)=\chi_{(k,k+\frac{1}{k}]}$ does not converge uniformly on $\mathbb{R}$.
\item[(c)] Does the sequence of functions $f_k(x)=\chi_{(k,k+\frac{1}{k}]}$ converge almost uniformly on $\mathbb{R}$?
\end{itemize}
\end{exercise}

\begin{exercise} \label{E83}
Let $X=[0,1]$, $\mathsf{S}=\mathcal{B}([0,1])$, and let $\mu=\lambda$ be the Lebesgue measure on $\mathsf{S}$. Prove the following:
\begin{itemize}
\item[(a)] The sequence of functions $f_k:=k\chi_{[0,\frac{1}{k}]}$ converges almost everywhere on $[0,1]$.
\item[(b)] Does the sequence of functions $f_k:=k\chi_{[0,\frac{1}{k}]}$ converge almost uniformly on $[0,1]$?
\end{itemize}
\end{exercise}

\begin{exercise} \label{E84}
Let $X=[0,1]$, $\mathsf{S}=\mathcal{B}([0,1])$, and let $\mu=\lambda$ be the Lebesgue measure on $\mathsf{S}$. Prove the following:
\begin{itemize}
\item[(a)] The sequence of functions $f_k(x)=\max\{k-k^2|x-\frac{1}{k}|,0 \}$ converges pointwise on $[0,1]$.
\item[(b)] The sequence of functions $f_k(x)=\max\{k-k^2|x-\frac{1}{k}|,0 \}$ does not converge uniformly on $[0,1]$.
\item[(c)] Does the sequence of functions $f_k(x)=\max\{k-k^2|x-\frac{1}{k}|,0 \}$ converge almost uniformly on $[0,1]$?
\end{itemize}
\end{exercise}

\begin{exercise} \label{E85}
Let $X=[0,1)$, $\mathsf{S}=\mathcal{B}([0,1))$, and let $\mu=\lambda$ be the Lebesgue measure on $\mathsf{S}$. Define the following sequence of intervals:

$$
\begin{array}{lcll}
L_1 & & I_{1}:=[0,1) &\\
L_{2} & & I_{2}:=[0,1/2)\,\,\,\,I_{3}:=[1/2,1)&\\
L_{3} & & I_{4}:=[0,1/4)\,\,\,\,I_{5}:=[1/4,1/2)&I_{6}:=[1/2,3/4)\,\,\,\,I_{7}:=[3/4,1)\\
\vdots & & & \\
L_{k} & & I_{2^{k-1}}:=[0,1/2^{k-1}) \,\,\ldots \ldots & I_{2^k-1}:=[1-1/2^{k-1},1)
\end{array}
$$

We now define the sequence of functions:
$$
\begin{array}{lcll}
L_1 & & f_{1}:=\chi_{I_1} &\\
L_{2} & & f_{2}:=\chi_{I_2}\,\,\,\,f_{3}:=\chi_{I_3}&\\
L_{3} & & f_{4}:=\chi_{I_4}\,\,\,\,f_{5}:=\chi_{I_5}&f_{6}:=\chi_{I_6}\,\,\,\,f_{7}:=\chi_{I_7}\\
\vdots & & & \\
L_{k} & & f_{2^{k-1}}:=\chi_{I_{2^{k-1}}} \,\,\ldots \ldots & f_{2^k-1}:=\chi_{I_{2^k-1}}
\end{array}
$$

\begin{itemize}
    \item[(a)] Prove that the sequence $(f_k)$ converges to $0$ in mean $p$ for every $p \in [1,\infty)$.
    \item[(b)] Prove that $\liminf_{k \to \infty}f_{k}(x)=0$ and $\limsup_{k \to \infty}f_k(x)=1$ for every $x \in [0,1]$.
\end{itemize}
\end{exercise}

\begin{exercise} \label{E86}
Let $(X,\mathsf{S},\mu)$ be a measure space, $(f_k)$ a sequence of functions in $\mathbb{M}(X,\mathsf{S})$, and $f \in \mathbb{M}(X,\mathsf{S})$. Prove the following:
\begin{itemize}
    \item[(a)] $f_k \to 0$ $\mu$ a.e. if and only if for every $\varepsilon >0$,
    $$
    \mu\left( \limsup_{k \to \infty}\{x \in X\,:\, |f_{k}(x)|>\varepsilon\} \right)=0.
    $$
    \item[(b)] From the previous item, conclude that $f_k \to f$ $\mu$ a.e. if and only if for every $\varepsilon >0$,
    $$
    \mu\left( \limsup_{k \to \infty}\{x \in X\,:\, |f_{k}(x)-f(x)|>\varepsilon\} \right)=0.
    $$
\end{itemize}
\end{exercise}

\begin{exercise} \label{E87}
Let $X=\mathbb{R}$, $\mathsf{S}=\mathcal{B}(\mathbb{R})$, and let $\mu=\lambda$ be the Lebesgue measure on $\mathsf{S}$. Show that the sequence of functions $f_k:=\frac{1}{k^{1/p}}\,\chi_{(0,k)}$ converges uniformly to $0$ but does not converge in mean $p$ for $p \in [1,\infty)$.
\end{exercise}

\begin{exercise} \label{E88}
Let $X=\mathbb{R}$, $\mathsf{S}=\mathcal{B}(\mathbb{R})$, and let $\mu=\lambda$ be the Lebesgue measure on $\mathsf{S}$. Show that the sequence of functions $f_k:=\chi_{[k,k+1]}$ converges almost everywhere to $0$ but does not converge in measure $\lambda$ to $0$.
\end{exercise}

\begin{exercise} \label{E89}
Let $X=[0,2]$, $\mathsf{S}=\mathcal{B}([0,2])$, and let $\mu=\lambda$ be the Lebesgue measure on $\mathsf{S}$. Show that the sequence of functions $f_k:=k^{1/p}\,\chi_{[\frac{1}{k},\frac{2}{k}]}$ converges pointwise to $0$ but does not converge in mean $p$ for $p \in [1,\infty)$.
\end{exercise}

{\setlength{\parindent}{0pt}
\begin{exercise} \label{E810}
Show that, in general, almost uniform convergence does not imply uniform convergence even on the complement of a set of measure zero.

\medskip
(Hint: Consider the sequence of functions $f_k:=k^{1/p}\,\chi_{[\frac{1}{k},\frac{2}{k}]}$ on $[0,2]$ endowed with its Borel $\sigma$-algebra and the Lebesgue measure $\lambda$.)
\end{exercise}}

\begin{exercise} \label{E811}
Let $X=(0,1)$, $\mathsf{S}=\mathcal{B}(0,1)$, and let $\mu=\lambda$ be the Lebesgue measure on $\mathsf{S}$. Show that the sequence of functions $f_k:=k\chi_{(0,\frac{1}{k^2})}$ converges to $0$ in mean $1$ but does not converge in mean $2$.
\end{exercise}

\begin{exercise} \label{E812}
Let $X=[0,1]$, $\mathsf{S}=\mathcal{B}([0,1])$, and let $\mu=\lambda$ be the Lebesgue measure on $\mathsf{S}$. Show that the sequence of functions $f_k:[0,1] \to \mathbb{R}$ given by
$$
f_k(x):=\left\{
\begin{array}{lcl}
1-kx & & \mbox{if } x \in [0,\frac{1}{k}],\\
0  & & \mbox{if } x \in [\frac{1}{k},1],\\
\end{array}
\right.
$$ 
converges almost everywhere on $[0,1]$.
\end{exercise}

\begin{exercise} \label{E813}
Let $(X,\mathsf{S},\mu)$ be a measure space. Show that if $(f_k)$ is a sequence of measurable functions which is Cauchy in measure $\mu$ and if some subsequence of $(f_k)$ converges in measure $\mu$ to a measurable function $f$, then $(f_k)$ converges in measure $\mu$ to $f$.
\end{exercise}

\begin{exercise} \label{E814}
Let $(X,\mathsf{S},\mu)$ be a measure space. The symbol $\circledast$ will denote any of the four notions of convergence studied. Show that if $(f_k)$ is a sequence of measurable functions converging $\circledast$ to a measurable function $f$, then any subsequence $(f_{k_j})$ also converges $\circledast$ to $f$.
\end{exercise}

\begin{exercise} \label{E815}
Let $(X,\mathsf{S},\mu)$ be a measure space. The symbol $\circledast$ will denote any of the four notions of convergence studied. Let $(f_k)$ and $(g_k)$ be sequences of measurable functions such that $f_k \to f$ $\circledast$ and $g_k \to g$ $\circledast$, where $f$ and $g$ are measurable functions. Prove the following statements:
\begin{itemize}
\item[(a)] $\gamma f_k \to \gamma f$ $\circledast$ for all $\gamma \in \mathbb{R}$.
\item[(b)] $f_k + g_k \to f + g$ $\circledast$.
\item[(c)] $|f_k| \to |f|$ $\circledast$.
\item[(d)] $\max\{f_k,g_k\} \to \max\{f,g\}$ $\circledast$ and $\min\{f_k,g_k\} \to \min\{f,g\}$ $\circledast$.
\item[(e)] $f_k^{+} \to f^{+}$ $\circledast$ and $f_k^{-} \to f^{-}$ $\circledast$.
\item[(f)] $\chi_{A} f_k \to \chi_{A} f$ $\circledast$ for all $A \in \mathsf{S}$.
\end{itemize}
\end{exercise}

{\setlength{\parindent}{0pt}
\begin{exercise} \label{E816}
Let $(X,\mathsf{S},\mu)$ be a measure space and let $(f_k)$, $(g_k)$ be two sequences of measurable functions. Show that if $f_k\xrightarrow[\mu]{} 0$ and $g_k \xrightarrow[\mu]{} 0$, then $f_k\,g_k \xrightarrow[\mu]{} 0$.
\end{exercise}

(Hint: Prove the following inclusions, valid for every $\varepsilon >0$ and every $k \in \mathbb{N}$,
$$
\begin{aligned}
\{ x \in X\,:\, |f_{k}(x)g_{k}(x)| \geq \varepsilon \}
&\subset \left\{x \in  X\,:\, \tfrac{|f_{k}(x)|^{2}}{2} + \tfrac{|g_{k}(x)|^{2}}{2} \geq \varepsilon \right\}\\
&\subset  \left\{x \in  X\,:\, \tfrac{|f_{k}(x)|^{2}}{2} \geq \tfrac{\varepsilon}{2} \right\}
\cup \left\{x \in  X\,:\, \tfrac{|g_{k}(x)|^{2}}{2} \geq \tfrac{\varepsilon}{2} \right\}.
\end{aligned}
$$)}

\begin{exercise} \label{E817}
Let $(X,\mathsf{S},\mu)$ be a measure space and let $\Psi:\mathbb{R} \to \mathbb{R}$ be a uniformly continuous function. Prove the following statements:
\begin{itemize}
\item[(a)] If $f_k:X \to \mathbb{R}$, $k \in \mathbb{N}$, is a sequence of measurable functions that converges uniformly to $f:X \to \mathbb{R}$, then the sequence $(\Psi \circ f_k)$ converges uniformly to $\Psi \circ f$.
\item[(b)] If $f_k:X \to \mathbb{R}$, $k \in \mathbb{N}$, is a sequence of measurable functions that converges almost uniformly to $f:X \to \mathbb{R}$, then the sequence $(\Psi \circ f_k)$ converges almost uniformly to $\Psi \circ f$.
\item[(c)] If $f_k:X \to \mathbb{R}$, $k \in \mathbb{N}$, is a sequence of measurable functions that converges in measure to $f:X \to \mathbb{R}$, then the sequence $(\Psi \circ f_k)$ converges in measure to $\Psi \circ f$.
\end{itemize}
\end{exercise}

\begin{exercise} \label{E818}
Let $(X,\mathsf{S},\mu)$ be a measure space and let $(f_k:X \to \mathbb{R})$ be a sequence of measurable functions. Show that $(f_k)$ converges almost uniformly on $X$ if and only if $(f_k)$ is Cauchy almost uniformly on $X$.
\end{exercise}

\begin{exercise} \label{E819} \index{Fatou's lemma!in measure} \index{monotone convergence theorem!in measure}
Let $(X,\mathsf{S},\mu)$ be a measure space, $(f_{k})$ a sequence of elements of $\mathbb{M}(X,\mathsf{S})$, and $f \in \mathbb{M}(X,\mathsf{S})$. Prove the following:
\begin{itemize}
    \item[(a)] \textbf{Fatou's lemma in measure:} If $f_{k} \geq 0$ $\mu$ a.e. for all $k \in \mathbb{N}$ and $f_{k} \xrightarrow[\mu]{} f$, then
    $$
    \int_{X}f\,d\mu \leq \liminf_{k \to\infty}\int_{X} f_k\,d\mu.
    $$
    \item[(b)] \textbf{Monotone convergence theorem in measure:} If $f_{k} \in L^{1}(X,\mathsf{S},\mu)$ and $f_{k} \leq f_{k+1}$ $\mu$ a.e. for all $k \in \mathbb{N}$, the sequence $(\int_{X}f_k\,d\mu)$ is bounded in $\mathbb{R}$, and $f_{k} \xrightarrow[\mu]{} f$, then $f \in L^{1}(X,\mathsf{S},\mu)$ and
    $$
    \int_{X}f \,d\mu = \lim_{k \to \infty}\int_{X}f_{k}\,d\mu.
    $$
\end{itemize}
\end{exercise}

\begin{exercise} \label{E820}
Prove the following:
\begin{itemize}
    \item[(a)] The conclusion of Egorov's theorem may fail if $\mu(X)=+\infty$, by considering the measure space $(\mathbb{R},\mathcal{B}(\mathbb{R}),\lambda)$ and the sequence of functions $f_{k}(x):=\chi_{(k,\infty)}(x)$.
    \item[(b)] The conclusion of Egorov's theorem may fail if the limit function takes extended real values on a set of positive measure, by considering the sequence $f_k:=k\chi_{[0,1]}$ on the measure space $([0,1], \mathcal{B}([0,1]),\lambda)$.
\end{itemize}
\end{exercise}

\begin{exercise} \label{E821}
Let $X=\mathbb{N}$, $\mathsf{S}=\mathcal{P}(\mathbb{N})$, and let $\mu=\mu^{\sharp}$ be the counting measure. Show that, given a sequence of measurable functions $(f_k)$, the sequence $(f_k)$ converges in measure to a measurable function $f$ if and only if $(f_k)$ converges uniformly to $f$.
\end{exercise}

{\setlength{\parindent}{0pt}
\begin{exercise}[The space $L^{0}$] \label{E822} \index{space!L0@$L^{0}(X,\mathsf{S},\mu)$}
Let $(X,\mathsf{S},\mu)$ be a finite measure space. Consider the equivalence relation on $\mathbb{M}(X,\mathsf{S})$ given by: $f\sim_\mu g$ if and only if $f=g$ $\mu$ almost everywhere.

We denote the set of equivalence classes by
$$
L^{0}(X,\mathsf{S},\mu):=\mathbb{M}(X,\mathsf{S})\diagup \sim_\mu.
$$

For simplicity, we identify measurable functions that are equal almost everywhere, so that $L^{0}(X,\mathsf{S},\mu)=\mathbb{M}(X,\mathsf{S})$.

Prove the following statements:
\begin{itemize}
\item[(a)] The function $d_\mu:L^{0}(X,\mathsf{S},\mu) \times L^0(X,\mathsf{S},\mu)\to \mathbb{R}$ given by
$$
d_\mu(f,g):=\int_X \frac{|f-g|}{1+|f-g|}\,d\mu
$$
is a metric.
\item[(b)] If $(f_k)$ is a sequence in $L^{0}(X,\mathsf{S},\mu)$, then $f_k \to f$ in measure if and only if $d_\mu(f_k,f) \to 0$.
\item[(c)] If $(f_k)$ is a Cauchy sequence in $(L^0,d_\mu)$, then there exists $f \in L^{0}(X,\mathsf{S},\mu)$ such that $d_\mu(f_k,f) \to 0$.
\end{itemize}
\end{exercise}

(Hint: Apply Chebyshev’s inequality to the function $\varphi(x):=\frac{x}{1+x}$.)}

\begin{exercise} \label{E823}
Define the function $d:L^{0}(X,\mathsf{S},\mu) \times L^0(X,\mathsf{S},\mu) \to \mathbb{R}$ by
$$
d(f,g):=\int_X \min\{|f-g|,1\}\,d\mu.
$$
Prove the following statements:
\begin{itemize}
\item[(a)] $d$ is a metric on $L^{0}(X,\mathsf{S},\mu)$.
\item[(b)] If $(f_k)$ is a sequence in $L^{0}(X,\mathsf{S},\mu)$, then $f_k \to f$ in measure if and only if $d(f_k,f) \to 0$.
\item[(c)] $(f_k)$ is a Cauchy sequence in measure if and only if $(f_k)$ is Cauchy in $(L^0,d)$.
\end{itemize}
\end{exercise}

\begin{exercise} \label{E824}
Let $(\mathbb{R},\mathcal{B}(\mathbb{R}))$ be a measurable space and let $\mu:\mathcal{B}(\mathbb{R}) \to [0,\infty)$ be a finite measure. Let $f:\mathbb{R} \to \mathbb{R}$ be a bounded function that is continuous at $0$. Show that if $(f_k)$ is a sequence of measurable functions such that $f_k \xrightarrow[\mu]{} 0$, then
$$
\lim_{k \to\infty}\int_{\mathbb{R}} (f \circ f_k)\,d\mu = \mu(\mathbb{R})\cdot f(0).
$$
\end{exercise}

\begin{exercise} \label{E825}
Let $(X,\mathsf{S},\mu)$ be a measure space. Determine whether the following statements are true or false:
\begin{itemize}
    \item[(a)] Uniform convergence implies convergence in measure.
    \item[(b)] Uniform convergence implies convergence in mean $p \in [1,\infty)$.
    \item[(c)] Convergence in measure implies uniform convergence.
    \item[(d)] Convergence in mean $p \in [1,\infty)$ implies uniform convergence.
\end{itemize}
\end{exercise}

\begin{exercise}\label{E826}
Let $(X,\mathsf{S},\mu)$ be a measure space and let $f:X \to \mathbb{R}$ be a function.
\begin{itemize}
    \item[(a)] Show that if $X$ is complete, then $f$ is $\mathsf{S}$-measurable if and only if there exists a sequence of functions $(f_k)$ in $\mathbb{M}(X,\mathsf{S})$ such that $f_k \to f$ $\mu$ almost everywhere.
    \item[(b)] Does the result in the previous item remain valid if the assumption that $X$ is complete is removed? Justify your answer.
\end{itemize}
\end{exercise}

Let $(X,\mathsf{S},\mu)$ be a measure space. A subset $\mathcal{H} \subset L^{1}(X,\mathsf{S},\mu)$ is said to be \textbf{uniformly integrable} \index{sequence!uniformly integrable} if
$$
\lim_{M \to \infty} \sup\left\{\int_{\{x \in X\,:\,|f(x)|\geq M\}}|f|\,d\mu\,:\, f\in \mathcal{H} \right\}=0.
$$

\begin{exercise}\label{E827}
Let $(X,\mathsf{S},\mu)$ be a measure space. Prove the following statements:
\begin{itemize}
    \item[(a)] Consider $X=(0,1)$, $\mathsf{S}=\mathcal{B}(0,1)$, and let $\mu=\lambda$ be the Lebesgue measure on $\mathsf{S}$. The set $\mathcal{H}:=\{x^{k}\,:\,k \in \mathbb{N}\}$ is uniformly integrable.
    \item[(b)] Let $(f_k)$ and $(g_k)$ be two sequences in $L^{1}(X,\mathsf{S},\mu)$ that are uniformly integrable. Then the sequence $(f_k+g_k)$ is uniformly integrable.
    \item[(c)] Let $\mathcal{H}$ be a subset of $L^{1}(X,\mathsf{S},\mu)$. If there exists $g \in L^{1}(X,\mathsf{S},\mu)$ such that $|f| \leq g$ $\mu$ a.e. for every $f \in \mathcal{H}$, then $\mathcal{H}$ is uniformly integrable.
    \item[(d)] Consider $X=[0,1]$, $\mathsf{S}=\mathcal{B}([0,1])$, and let $\mu=\lambda$ be the Lebesgue measure on $\mathsf{S}$. The set $\mathcal{H}:=\{f_k\,:\,k \geq 2\}$, where $f_k:=\frac{k}{\log(k)}\chi_{[0,\tfrac{1}{k}]}$, is uniformly integrable and there is no function $g \in {L}^{1}$ such that $|f_k| \leq g$ for all $k \geq 2$.
\end{itemize}
\end{exercise}

Let $(X,\mathsf{S},\mu)$ be a measure space. We say that a subset $\mathcal{H} \subset L^{1}(X,\mathsf{S},\mu)$ is \textbf{uniformly absolutely continuous} \index{sequence!uniformly absolutely continuous} if, for every $\varepsilon >0$, there exists $\delta >0$ (depending on $\varepsilon$) such that
$$
\sup\left\{\int_{A} f\,d\mu\,:\, f \in \mathcal{H} \right\}<\varepsilon \quad \text{whenever }\,\,\mu(A)<\delta.
$$

\begin{exercise}\label{E828}
Let $(X,\mathsf{S},\mu)$ be a measure space. Prove the following statements:    
\begin{itemize}
    \item[(a)] Consider $X=(0,1)$, $\mathsf{S}=\mathcal{B}(0,1)$, and let $\mu=\lambda$ be the Lebesgue measure on $\mathsf{S}$. The set $\mathcal{H}:=\{x^{k}\,:\,k \in \mathbb{N}\}$ is uniformly absolutely continuous.
    \item[(b)] Let $(f_k)$ be a bounded sequence in $L^{1}(X,\mathsf{S},\mu)$. If $(f_k)$ is uniformly integrable, then $(f_k)$ is uniformly absolutely continuous. Is the converse true? Justify your answer.
    \item[(c)] If $(f_k)$ is a bounded sequence in $L^{1}(X,\mathsf{S},\mu)$ and uniformly absolutely continuous, then $(f_k)$ is uniformly integrable.
    \item[(d)] Let $(f_k)$ be a bounded sequence in $L^{1}(X,\mathsf{S},\mu)$ which is uniformly integrable. If for every $\varepsilon >0$, there exists $A_{\varepsilon} \in \mathsf{S}$ with $\mu(A_{\varepsilon})<+\infty$ such that
    $$
\sup_{k \in \mathbb{N}}\left\{\int_{X \smallsetminus A_{\varepsilon}}|f_k|\,d\mu  \right\}<\varepsilon,
    $$
    then $(f_k)$ is bounded in $L^{1}(X,\mathsf{S},\mu)$.
    \item[(e)] Is it true that if $\mu(X)<+\infty$, then every uniformly integrable sequence in $L^{1}(X,\mathsf{S},\mu)$ satisfies the hypotheses of the previous item? Justify your answer.
\end{itemize}
\end{exercise}

\begin{exercise}\label{E829}
Let $(X,\mathsf{S},\mu)$ be a measure space and let $p \in [1,\infty)$ be fixed.

Prove the following statements:
\begin{itemize}
    \item[(a)] If $f \in L^{p}(X,\mathsf{S},\mu)$, then for every $\varepsilon >0$, there exists $A_{\varepsilon} \in \mathsf{S}$ with $\mu(A_{\varepsilon})<+\infty$ such that
    $$
\int_{X \smallsetminus A_{\varepsilon}}|f|^{p}\,d\mu <\varepsilon^{p}.
    $$
    \item[(b)] If $(f_k)$ is a sequence of elements of $L^{p}(X,\mathsf{S},\mu)$ such that $f_{k} \xrightarrow[L^{p}]{} f$, then the sequence $(|f_k|^{p})$ satisfies that for every $\varepsilon >0$, there exists $A_{\varepsilon} \in \mathsf{S}$ with $\mu(A_{\varepsilon})<+\infty$ such that
    $$
\sup_{k \in \mathbb{N}}\left\{\int_{X \smallsetminus A_{\varepsilon}}|f_k|^{p}\,d\mu  \right\}<\varepsilon^{p}.
    $$
    \item[(c)] If $(f_k)$ is a sequence of elements of $L^{p}(X,\mathsf{S},\mu)$ such that $f_{k} \xrightarrow[L^{p}]{} f$, then $(|f_k|^{p})$ is uniformly absolutely continuous.
    \item[(d)] If $(f_k)$ is a sequence of elements of $L^{p}(X,\mathsf{S},\mu)$ such that $f_{k} \xrightarrow[L^{p}]{} f$, then $(|f_k|^{p})$ is uniformly integrable.
    \end{itemize}
\end{exercise}

\begin{exercise}[Vitali convergence theorem] \label{E830} \index{theorem!Vitali convergence}
Let $(X,\mathsf{S},\mu)$ be a measure space, let $p \in [1,\infty)$ be fixed, and let $(f_k)$ be a sequence of functions in $L^{p}(X,\mathsf{S},\mu)$ satisfying the following three conditions:
\begin{itemize}
    \item[(1)] $f_{k} \xrightarrow[\mu]{} f$.
    \item[(2)] For every $\varepsilon >0$, there exists $A_{\varepsilon} \in \mathsf{S}$ with $\mu(A_{\varepsilon})<+\infty$ such that
    $$
\sup_{k \in \mathbb{N}}\left\{\int_{X \smallsetminus A_{\varepsilon}}|f_k|^{p}\,d\mu  \right\}<\varepsilon^{p}.
    $$ 
    \item[(3)] For every $\varepsilon >0$, there exists $\delta >0$ (depending on $\varepsilon$) such that 
    $$
\sup_{k \in \mathbb{N}}\left\{\int_{A} |f|^{p}\,d\mu\right\}<\varepsilon^{p} \quad \text{whenever }\,\,\mu(A)<\delta.
    $$
\end{itemize}

Show that the sequence $(f_k)$ converges to $f$ in mean $p$.
\end{exercise}

\section{Project: Almost everywhere uniform convergence}

Let $(X,\mathsf{S},\mu)$ be a measure space.

Prove the following theorem.

{\setlength{\parindent}{0pt}
\begin{theorem} \label{828} 
Let $(f_k)$ be a sequence of functions in $L^{\infty}(X,\mathsf{S},\mu)$ and let $f \in L^{\infty}(X,\mathsf{S},\mu)$. Then $\|f_k - f\|_{\infty} \to 0$ as $k \to \infty$ if and only if $f_k \to f$ uniformly $\mu$ almost everywhere. That is, $\|f_k - f\|_{\infty} \to 0$ if and only if there exists a $\mu$-null subset $N \subset X$ such that $f_k \to f$ uniformly on $X \smallsetminus N$.
\end{theorem}

(Hint: Review the proof of the Riesz–Fischer theorem).}

We give the following definition.

\begin{definition} \label{829} \index{convergence!uniformly almost everywhere}
We say that a sequence of functions $(f_k)$ in $L^{\infty}(X,\mathsf{S},\mu)$ converges \textbf{uniformly $\mu$ almost everywhere} to $f \in L^{\infty}(X,\mathsf{S},\mu)$ if $\|f_k - f\|_{\infty} \to 0$.

We denote this type of convergence by $\mu$ u.a.e.
\end{definition}

\subsection*{8.5.1.\quad Objective}

The objective of this project is to relate uniform convergence almost everywhere with the other notions of convergence studied.

\subsection*{8.5.2.\quad Procedure}

\begin{itemize}
    \item[1.] Show that if $(f_k)$ is a sequence of functions in $L^{\infty}(X,\mathsf{S},\mu)$ that converges $\mu$ u.a.e. to some function $f \in L^{\infty}(X,\mathsf{S},\mu)$, then $f_k \to f$ $\mu$ a.u.
    
    \item[2.] Deduce from the previous result that $\mu$ u.a.e. convergence implies $\mu$ a.e. convergence and convergence in measure.

    \item[3.] Show that if $\mu(X)<+\infty$ and $(f_k)$ converges $\mu$ u.a.e. to some function $f \in L^{\infty}(X,\mathsf{S},\mu)$, then $f_k \to f$ in mean $p$ for every $p \in [1,\infty)$.

    \item[4.] Show, through the following example, that the conclusion of the previous item may fail if $\mu(X)=+\infty$.

    Let $X=\mathbb{R}$, $\mathsf{S}=\mathcal{B}(\mathbb{R})$, and let $\mu=\lambda$ be the Lebesgue measure on $\mathsf{S}$. Show that the sequence of functions $f_k:\mathbb{R} \to \mathbb{R}$ given by $f_k:=\frac{1}{k}\chi_{[0,2^k]}$ converges $\lambda$ u.a.e. to $0$ but does not converge in mean $p$ for any $p \in [1,\infty)$.

    \item[5.] Show, however, that if the sequence $(f_k)$ is dominated in absolute value by some function $f \in L^{p}(X,\mathsf{S},\mu)$ for $p \in [1,\infty)$, then $\mu$ u.a.e. convergence implies convergence in mean $p$, even if $\mu(X)=+\infty$.

    \item[6.] Finally, show through the following example that $\mu$ u.a.e. convergence is a strictly more general notion than the other types of convergence, even if $\mu(X)<+\infty$.

    Let $X=[0,1]$, $\mathsf{S}=\mathcal{B}([0,1])$, and let $\mu=\lambda$ be the Lebesgue measure on $\mathsf{S}$. The sequence of functions $f_k:[0,1] \to \mathbb{R}$ given by $f_k:=k\chi_{[0,\frac{1}{2^k}]}$ converges to $0$ $\lambda$ a.e., $\lambda$ a.u., in measure $\lambda$, and in mean $p \in [1,\infty)$, but does not converge $\lambda$ u.a.e.
\end{itemize}

.

\setlength{\parskip}{1em} 
\lhead[{\scriptsize \thepage}]{ {\scriptsize  \rightmark }}
\rhead[ {\scriptsize   \leftmark}]{ {\scriptsize \thepage}}

\markboth{{\scshape Bibliografía}}{ {\scshape Bibliografía}}

\end{document}